    \documentclass[a4paper,12pt]{book}
    \counterwithin{figure}{chapter}

    \usepackage[utf8]{inputenc}
    \usepackage{anyfontsize}
    \usepackage{sectsty}
        \allsectionsfont{\raggedright}
    \usepackage{verbatim}
    \usepackage{amsfonts}
    \usepackage{amsmath}
    \usepackage{amsthm}
    \usepackage{amssymb} % for semidirect products \rtimes and \ltimes
    \usepackage{letltxmacro} % https://tex.stackexchange.com/questions/305174/proof-environment-produces-proof-proof-only-when-thmbox-is-used
        \LetLtxMacro\amsproof\proof
        \LetLtxMacro\amsendproof\endproof
    \usepackage{mathtools} % for injection arrows \xhookrightarrow{} and for \smallmatrix
    \usepackage{stmaryrd} % for double brackets
    \usepackage{calc}
    \usepackage{fancyhdr}
        \fancypagestyle{mainstyle}{%
            \fancyhead[RO]{\slshape\nouppercase{\leftmark}}
            \fancyhead[LE]{\slshape\nouppercase{\rightmark}}
            \fancyfoot[RO,LE]{\thepage}
            
            }
        \fancypagestyle{introstyle}{%
            \fancyhead{}
            \fancyfoot[RO,LE]{\thepage}
            
            }
        \fancypagestyle{plain}{%
            \fancyhf{} 
            \fancyfoot[RO,LE]{\thepage}
            
            }
    \usepackage[Lenny]{fncychap}
    \usepackage{enumitem} % for customizing lists
        \setlist[enumerate]{itemsep=1pt,topsep=2pt}
        \setlist[itemize]{itemsep=1pt,topsep=2pt}
        \setlist[description]{itemsep=1pt,topsep=2pt}
    \usepackage{rotating} % for sidewaystable
    \usepackage{multirow}
    \usepackage{tabularx} % for tabularx
        \newcolumntype{C}{>{\centering\arraybackslash}X}
        \newcolumntype{K}[1]{>{\centering\arraybackslash}m{#1}}
        \newcolumntype{5}{>{\hsize=1.5\hsize}X}
        \newcolumntype{2}{>{\hsize=.6\hsize}X}
        \newcolumntype{3}{>{\hsize=.9\hsize}X}
    \usepackage{ltablex}
    \usepackage[
        left=1.25in,
        right=1.25in,
        top=1.5in,
        bottom=1.75in,
        bindingoffset=0cm,
        heightrounded,
    ]{geometry}
    
    \usepackage[style=alphabetic,maxbibnames=99]{biblatex}
        \DeclareFieldFormat{title}{\mkbibquote{#1}}
    
    \usepackage{subcaption}
    \usepackage{color}
    \usepackage[dvipsnames]{xcolor}
    \usepackage{graphicx}\graphicspath{{./images/}}
    \usepackage{tikz-cd} % commutative diagrams: https://ctan.math.washington.edu/tex-archive/graphics/pgf/contrib/tikz-cd/tikz-cd-doc.pdf
    \usepackage{tikz}
        \usetikzlibrary{calc,arrows,arrows.meta,shapes}
        \usetikzlibrary{decorations.pathmorphing} % for curvy
        \tikzset{every picture/.style=thick}
        \tikzset{state/.style = {shape=ellipse,draw,inner sep=1pt, outer sep=1pt}}
        \tikzset{vertex/.style = {shape=circle,draw,inner sep=1pt, outer sep=1pt}}
        \tikzset{node/.style = {shape=circle,draw,fill,inner sep=1pt, outer sep=1pt}}
        \tikzset{edge/.style = {->, >={latex[scale=1.5]}}}
        \tikzset{curvy/.style={decorate, decoration=snake,segment length=.25cm}}

    \usepackage{thmtools}
    \usepackage{hyperref}
    \usepackage[noabbrev,nameinlink,capitalise]{cleveref}
        \crefname{claim}{Claim}{Claims}
        \crefname{subsection}{Subsection}{Subsections}
        \crefname{subsubsection}{Subsubsection}{Subsubsections}

\usepackage[mono=false]{libertine}
\usepackage{libertinust1math}

    \declaretheorem[numberwithin=chapter,thmbox=L]{theorem}
    \declaretheorem[numberlike=theorem,thmbox=S]{corollary}
    \declaretheorem[numberlike=theorem,thmbox=S]{lemma}
    \declaretheorem[numberlike=theorem,thmbox=S]{proposition}
    \declaretheorem[numberlike=theorem,thmbox=S]{claim}
    \declaretheorem[numberlike=theorem,shaded={rulecolor=black,rulewidth=.5pt,bgcolor=white},style=definition]{definition}
    \declaretheorem[numberlike=theorem,shaded={rulecolor=black,rulewidth=.5pt,bgcolor=white},style=definition]{example}
    \declaretheorem[numberlike=theorem,shaded={rulecolor=black,rulewidth=.5pt, bgcolor=white},style=definition]{assumption}
    \declaretheorem[numberlike=theorem,style=remark]{remark}
    \declaretheorem[numberlike=theorem,style=remark]{question}
    \newenvironment{sketch}{%
    \proof}{\endproof}
    \newenvironment{maintheorem*}{%
        \thmbox[L]{\textbf{Main Theorem}}%
        \hspace*{-1.5em}\slshape\ignorespaces%
        }
        {%
        \endthmbox\vspace*{.75ex}%
        }
    \newenvironment{corollary*}{%
        \thmbox[S]{\textbf{Corollary}}%
        \hspace*{-1.5em}\slshape\ignorespaces%
        }
        {%
        \endthmbox\vspace*{.75ex}%
        }

    \title{Rearrangement Groups of Fractals: Structure and Conjugacy}
    \author{Matteo Tarocchi}
    \date{}

    \newcommand{\credits}[1]{
        \vspace{-\baselineskip/2}
        \caption*{\textcolor{gray}{\footnotesize{Credits: {#1}}}}
        } %https://tex.stackexchange.com/questions/95029/add-source-to-figure-caption

\AtBeginDocument{%
  \LetLtxMacro\proof\amsproof
  \LetLtxMacro\endproof\amsendproof
}

\begin{document}

%\documentclass[12pt,a4paper,oneside]{book}

% \usepackage[utf8]{inputenc}
% \usepackage[english]{babel}
% \usepackage{multicol}
% \usepackage{graphicx}\graphicspath{{./images/}}
% \usepackage[
%     left=1.25in,
%     right=1.25in,
%     top=1.5in,
%     bottom=1.75in,
%     bindingoffset=0cm,
%     heightrounded,
% ]{geometry}
% \linespread{1.25}

%%%%%%%%%%%%%%%%%%%%%%%%%

%% FONTS 

% \usepackage[mono=false]{libertine}
% \usepackage{libertinust1math}

%%%%%%%%%%%%%%%%%%%%%%%%%

%\begin{document}
\thispagestyle{empty}
%
% UNIVERSITY
\begin{center}
{\Large \sc Università di Pavia e di Milano-Bicocca
\\
e Istituto Nazionale di Alta Matematica}
\vspace{.8cm}\\

\begin{figure}[!htb]
\minipage{0.25\textwidth}
\includegraphics[width=\textwidth]{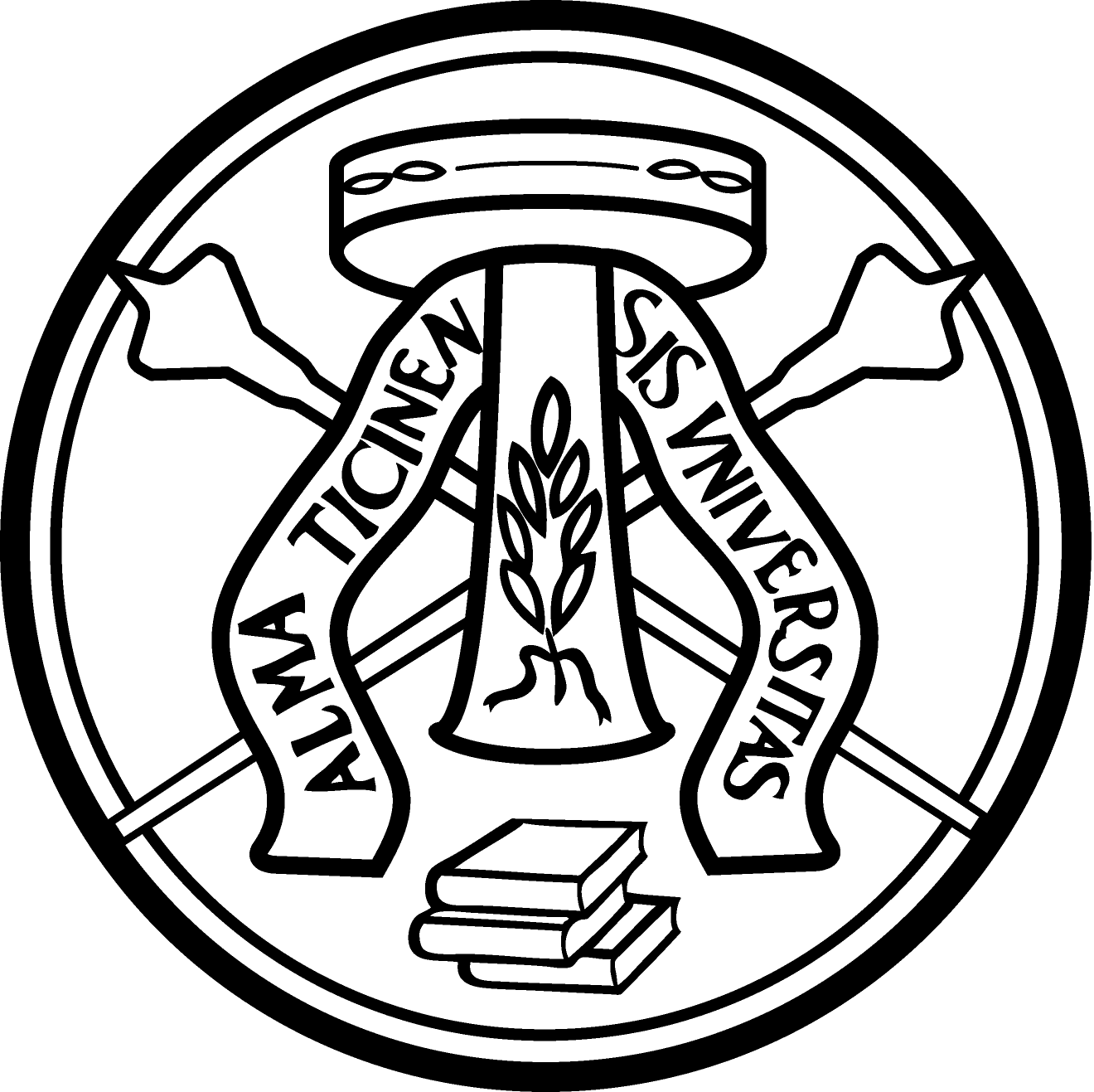}
\endminipage\hfill
\minipage{0.275\textwidth}
\includegraphics[width=\textwidth]{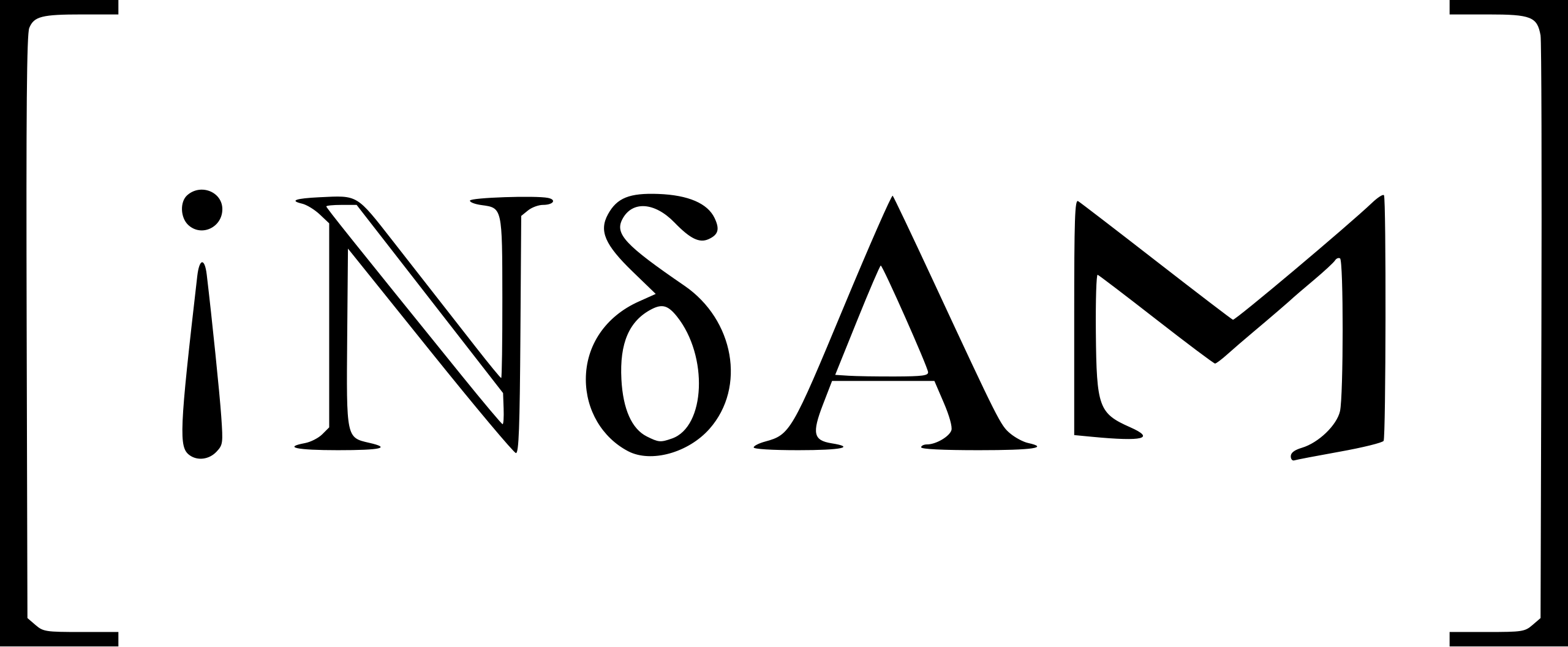}
\endminipage\hfill
\minipage{0.25\textwidth}%
\includegraphics[width=\textwidth]{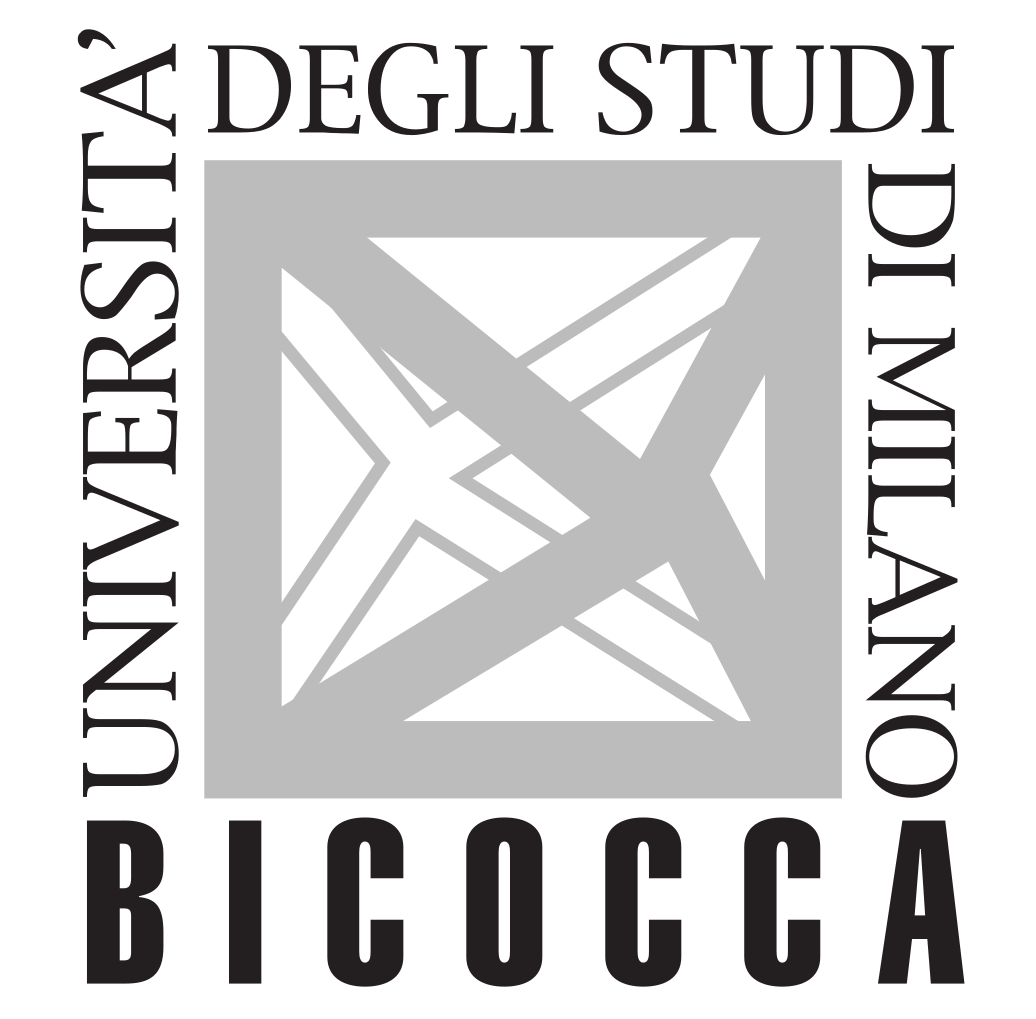}
\endminipage
\end{figure}
% \includegraphics[width=.25\textwidth]{UniPv.png}
% \hspace{1.25cm}
% \includegraphics[width=.275\textwidth]{indam.png}
% \hspace{1.25cm}
% \includegraphics[width=.25\textwidth]{UniMiB.png}
\vspace{.8cm}
%
% THESIS TITLE
\hrule height 1 pt
\vspace{.8cm}
\huge \sc Rearrangement Groups of Fractals:\\Structure and Conjugacy
\vspace{.8cm}
\hrule height 1 pt
\end{center}
\vspace{.4cm}
%
% AUTHOR AND SUPERVISORS
\begin{center}
%\begin{minipage}{.9\textwidth}
{\Large Author: \hfill Supervisor:}
\vspace{.125cm} \\
{\Large Matteo Tarocchi \hfill Francesco Matucci}
%\\ \\
%\phantom{x} \hfill \textit{Co-Supervisor:}\\
%\phantom{x} \hfill {Your Co-Supervisor's name} \\
%\end{minipage}
\vfill
{\large A thesis submitted in fulfillment of the requirements for the degree of} \\ \ \\
{\Large Doctor of Philosophy} \\ \ \\
{\large in the} \\ \ \\
{\Large Department of Mathematics and Applications}
\\ \ \\
XXXVII cycle
\end{center}

\newpage
\thispagestyle{empty}
\
%\end{document}

\pagenumbering{gobble}
\pagestyle{introstyle}

%%%%%%%%%%%%%%%%%%%%%%%%%

\section*{Abstract}

This dissertation is about \textbf{rearrangement groups}:
a class of (countable) groups of homeomorphisms of fractal topological spaces.
Introduced in 2019 by J. Belk and B. Forrest, this class generalizes the famous trio of Thompson groups $F$, $T$ and $V$ and includes some of their relatives and generalizations.

Each rearrangement group and the space on which it acts is determined by an \textbf{edge replacement system}, which is a specific type of graph rewriting systems.
The edge replacement system defines an edge shift and a \textbf{gluing relation} on it;
the fractal space $X$ is defined as the quotient under this relation, whenever it is an equivalence relation (for which there are natural sufficient conditions).
Rearrangements are those homeomorphisms of $X$ that descend from prefix-exchange homeomorphisms of the edge shifts that preserve the gluing relation.
They are represented by isomorphisms between graphs that approximate the fractal.
After an introduction to this topic, this dissertation branches into different aspects of rearrangement groups.

We first focus on a class of rearrangement groups of tree-like fractals known as \textbf{Wa\.zewski dendrites}.
We find finite generating sets for them and their commutator subgroups and we prove that the commutator subgroups are simple (with one possible exception).
We also show that these groups are dense in the groups of all homeomorphisms of dendrites.

We then extend our focus to the entire class of rearrangement groups and we provide a sufficient condition to solve the \textbf{conjugacy problem} (the decision problem of establishing whether two elements are conjugate) using the graphical tool of \textbf{strand diagram}.
This condition is not necessary, but understanding how to extend it is related to open problems in computer science.
However, this method solves the conjugacy problem in essentially all known rearrangement groups.

Next we study \textbf{invariable generation}:
the property of admitting a set of elements that generate the group even if each element is replaced by a conjugate.
Generalizing results of Gelander, Golan and Juschenko about Thompson groups,
with a dynamical approach we prove that transitive enough rearrangement groups are not invariably generated, which applies to most rearrangement group that has been considered so far.

Then we study the gluing relation that defines the fractal topological spaces on which rearrangement groups act.
We show that it is a \textbf{rational relation}, i.e., there exists a finite-state automaton which reads a pair of elements of the edge shift if and only if the two elements are related.
In fact, we provide an algorithmic construction of such finite-state automata.

Finally, we collect results that are based on constructions and modifications of edge replacement systems.
We first prove that finite groups and finitely generated abelian groups are rearrangement groups.
Then we show that the stabilizer of a finite set of rational points under the action of a rearrangement group is itself a rearrangement group.
Lastly, we prove that every rearrangement group embeds into Thompson's group $V$.

%%%%%%%%%%%%%%%%%%%%%%%%%

\cleardoublepage
\section*{Acknowledgements}

Let me start by thanking Francesco Matucci for everything:
for his mathematical supervision and his career advice, for his support when I felt lost, for his kindness and infinite patience during the last four years.
I am truly sorry that he had to endure so many of my stupid questions.
My thanks also go to Carlo Casolo and Orazio Puglisi, without whom I would not have never met Francesco.

I am genuinely fond of my time in Bicocca, and that is only thanks to all of the wondeful PhD students and postdocs who enlivened the U5 during my stay.
They have been way too important for me to omit their names:
my sincere thanks to Alberto, Andrea, Andrea, Anna, Bianca, Claudia, Claudio, Davide, Elena, Eleonora, Ettore, Fabio, Federico, Francesca, Gabriele, Giovanni, Giorgio, Giulia, Ilaria, Islam, Julian, Kirubell, Luca, Ludovico, Luigi, Marco, Marco, Marta, Nicola, Nicola, Ninni, Paolo, Sara, Simone, Simone, Stefano, Tommaso, Valentina.
An additional thank to Ninni for introducing me to the U5 when I first arrived.

I am especially grateful to Davide for so many things.
He made working so much more fun than it would have been by myself, taught me plenty of new maths and was always there for a chat on how badly everything was going.
He taught me that the complement of an open set is closed, for which I am indebted.
It should also be acknowledged that for three entire years he had to endure my awful sense of humor.

A huge thank to Rachel Skipper for her hospitality in Paris, for patiently explaining her work to me and for her wise words when I felt discouraged.
Her invitation to give a seminar talk was greatly impactful and I am also grateful to Yuri Neretin for the conversation that it sparked.
I would also like to thank Eduardo and the others at the ENS for their company during my stay there.

I want to express my gratitude to Jim Belk for his hospitality and for dedicating so much of his time and patience teaching me maths and helping me with my thesis.
His love for maths is deeply infectious and his cookies at the G\&T seminar are appreciated.
I would also like to give my heartfelt thanks to all the PhD students that I met in Glasgow for their outstanding friendliness.

I want to sincerely thank Bruno Duchesne for his amiability and hospitality at Orsay, for his patience listening to my bad ideas, for sharing his good ones with me and for believing in me.

I wish to thank Feyishayo Olukoya and Stefan Witzel for their useful comments that lead to improvements in the clarity and exposition of this thesis.

This section would be incomplete if I did not thank the plethora of wonderful people that I met at conferences.
Special thanks to Maarten for bearing me for two weeks non-stop and to Marco for his all-encompassing sense of humor.

Aside from those I met at work, I want to thank everyone who was patient enough to be around me for the last three years.
My first thank here goes to my sister's dog Zuma, always ready for a cuddle and for accompanying me on a short hike, joyful and energetic despite her age and health.
I am still terribly sorry for getting her tail stuck in the car door and I wish I could tell that to her.

Now that I mentioned my sister's dog, it is time that I thank my sister Ilaria too:
I am grateful to her for not despising me every time I was despicable, for her jokes, our chats about movies and all the memes she shared.
And, of course, I would like to thank my parents for their love and for supporting me throughout this long journey.
I would not be here without them (literally).
Let me also thank my uncle Stefano and my grandma Franca for their endless love.

It could go without saying that I have to thank Madalina for her intimate understanding of my struggles, for her tolerance and supporting me despite my difficult decisions, for her insights and our discussions and her humor, for her love.
It could go without saying all of this, but I want to mention it anyway, because it should never be taken for granted.

Without my friends, without their support even when we were far apart, these three years would have been impossible to survive.
I would like to thank them for saving my life plenty of times in roleplaying games and Giacomo in particular for saving my actual life.
I also want to thank Jacopo for all of our chats about life, the universe and everything.

Chances are, if you are reading this and are not mentioned here, I should probably thank you as well.

Thank you all.

%%%%%%%%%%%%%%%%%%%%%%%%%

\frontmatter
\renewcommand{\contentsname}{Table of contents}
\tableofcontents

%%%%%%%%%%%%%%%%%%%%%%%%%

\chapter*{List of Definitions and Symbols}
\addcontentsline{toc}{chapter}{List of Definitions and Symbols}
\markboth{List of definitions and symbols}{List of definitions and symbols}

The following table lists, in the order in which they are introduced, the notions that are defined throughout this dissertation.
For ease of reference, the table is subdivided into categories.

The second column features a commonly used symbol or an example of a symbol used for the notion, when available.
The third column contains (clickable) reference to where the notion is defined.
If a (sub)section is indicated as reference, it is understood that the notion is defined at the beginning or in some other obvious spot of the (sub)section.

\begingroup
\renewcommand{\arraystretch}{1.25}
\begin{tabularx}{\textwidth}{5|2|3}
    \textbf{Notion} & \textbf{Symbol} & \textbf{Reference}
    \\
    \hline
    \hline
    \multicolumn{3}{>{\hsize=\dimexpr3\hsize+4\tabcolsep+2\arrayrulewidth\relax}X}{\textbf{Graphs}}
    \rule{0pt}{20pt}
    \\
    \hline
    Graph & $(V,E,\iota,\tau)$ & \cref{def.graph}
    \\
    Adjacency and Incidence & & \cpageref{def.graph}
    \\
    Graph isomorphism & $\phi$ & \cref{def.graph.morphism}
    \\
    Subgraph & & \cref{def.subgraph}
    \\
    Spanned subgraph & & \cref{def.subgraph.spanned}
    \\
    Out-degree, In-degree, Degree & $\mathrm{out}(v), \mathrm{in}(v)$ & \cref{def.degree}
    \\
    Sink, Source & & \cref{def.sink.source}
    \\
    Loop & & \cref{def.loop}
    \\
    Isolated vertex & & \cref{def.isolated.vertex}
    \\
    Rotation system & & \cref{def.rotation.system}
    \\
    Walk, path or cycle & $(e_1, e_2, \dots)$ & \cref{def.walk.path.cycle}
    \\
    Inescapable cycle & & \cref{def.inescapable.cycle}
    \\
    Reachable vertex & & \cref{def.reachable}
    \\
    Biorientation & $\overline{\Gamma}$ & \cref{def.biorientation}
    \\
    Connected graph & & \cref{def.connected.graph}
    \\
    Connected component & & \cref{def.connected.graph}
    \\
    Reduced cycle & & \cref{def.reduced.cycle}
    \\
    Forest, tree & & \cref{def.forest.tree}
    \\
    Root & & \cpageref{txt.roots.leaves}
    \\
    Leaf & & \cpageref{txt.roots.leaves}
    \\
    Complete subforest & & \cref{def.complete.subforest}
    \\
    Labeled (and colored) graph & & \cref{def.labeled.graph}
    \\
    Morphism of colored graph & & \cref{def.colored.graph.morphism}
    \\
    \multicolumn{3}{>{\hsize=\dimexpr3\hsize+4\tabcolsep+2\arrayrulewidth\relax}X}{\textbf{Rewriting system}}
    \rule{0pt}{20pt}
    \\
    \hline
    Rewriting & $a \longrightarrow b$ & \cref{sub.rewriting.systems}
    \\
    Locally confluent rewriting system & & \cref{def.locally.confluent}
    \\
    Terminating rewriting system & & \cref{def.terminating}
    \\
    \multicolumn{3}{>{\hsize=\dimexpr3\hsize+4\tabcolsep+2\arrayrulewidth\relax}X}{\textbf{Automata}}
    \rule{0pt}{20pt}
    \\
    \hline
    Automaton & $(\Sigma, Q, \rightarrow, q_{0})$ & \cref{def.automaton}
    \\
    Transition & $q \xrightarrow{x} q'$ & \cref{def.automaton}
    \\
    Set of infinite sequences & $\Sigma^\omega$ & \cref{txt.infinite.sequences.symbol}
    \\
    Run & & \cref{def.run.recognize}
    \\
    Automaton recognizing a word or a sequence & & \cref{def.run.recognize}
    \\
    Word & $w = w_1 \dots w_n$ & \cpageref{txt.word.sequence}
    \\
    Sequence & $\alpha = x_1 x_2 \dots$ & \cpageref{txt.word.sequence}
    \\
    \multicolumn{3}{>{\hsize=\dimexpr3\hsize+4\tabcolsep+2\arrayrulewidth\relax}X}{\textbf{Edge shifts}}
    \rule{0pt}{20pt}
    \\
    \hline
    Cantor space & $\mathfrak{C}$ & \cref{sub:Cantor}
    \\
    Periodic sequence & $v \overline{w}$ & \cref{ex.interval}
    \\
    Edge shift & $\Omega(\Gamma,v_0)$ & \cref{def.edge.shift}
    \\
    Alphabet & $\mathbb{A}$ & \cref{def.alphabet.language}
    \\
    Language & $\mathbb{L}$ & \cref{def.alphabet.language}
    \\
    Cone & $C(x)$ & \cref{def.cones}
    \\
    \multicolumn{3}{>{\hsize=\dimexpr3\hsize+4\tabcolsep+2\arrayrulewidth\relax}X}{\textbf{edge replacement systems}}
    \rule{0pt}{20pt}
    \\
    \hline
    Edge replacement rule & $(R,C)$ & \cref{def.replacement.rules}
    \\
    $c$ replacement graph & $X_c$ & \cref{def.replacement.rules}
    \\
    Initial and terminal vertices of $X_c$ & $\iota_c$, $\tau_c$ & \cref{def.replacement.rules}
    \\
    Edge replacement system & $\mathcal{R} = (X_0,R,C)$ & \cref{def.replacement.system}
    \\
    Base graph & $X_0$ & \cref{def.replacement.system}
    \\
    Monochromatic or polychromatic edge replacement system & & \cpageref{txt.mono.poly.chromatic}
    \\
    Airplane edge replacement system & $\mathcal{A}$ & \cpageref{txt.airplane.replacement.system}
    \\
    Edge expansion & $\Gamma \triangleleft e$ & \cref{def.edge.expansion}
    \\
    Graph expansion & & \cref{def.edge.expansion}
    \\
    Full expansion sequence & $E_0, E_1, \dots$ & \cref{def.full.expansion}
    \\
    Color graph of $\mathcal{R}$ & $A_\mathcal{R}$ & \cref{def.color.graph}
    \\
    Alphabet of $\mathcal{R}$ & $\mathbb{A}_\mathcal{R}$ & \cref{prop.language.alphabet}
    \\
    Language of $\mathcal{R}$ & $\mathbb{L}_\mathcal{R}$ & \cref{prop.language.alphabet}
    \\
    Symbol space of $\mathcal{R}$ & $\Omega_\mathcal{R}$ & \cref{def.symbol.space}
    \\
    Expanding edge replacement system & & \cref{def.expanding}
    \\
    \\ %layout
    \multicolumn{3}{>{\hsize=\dimexpr3\hsize+4\tabcolsep+2\arrayrulewidth\relax}X}{\textbf{Limit spaces}}
    \rule{0pt}{20pt}
    \\
    \hline
    Gluing relation & $\sim$ & \cref{def.glue}
    \\
    Limit space & & \cref{def.limit.space}
    \\
    Point of a limit space & $\llbracket \alpha \rrbracket$ & \cpageref{txt.point.symbol}
    \\
    Cell & $\llbracket w \rrbracket$ & \cref{def.cell}
    \\
    Type of a cell & & \cref{def.cell.type}
    \\
    Topological interior of a cell & $\llparenthesis w \rrparenthesis$ & \cpageref{txt.topological.interior.symbol}
    \\
    Cellular partition & & \cref{def.cellular.partition}
    \\
    Cellular partition finer than another & & \cref{def.finer.partitions}
    \\
    Vertex of a limit space & & \cref{def.vertices}
    \\
    Boundary vertices of a cell & & \cpageref{txt.boundary.vertices}
    \\
    Rational, irrational points of a limit space & & \cref{def.rational.irrational.points}
    \\
    edge replacement system with finite branching & & \cref{def.finite.branching}
    \\
    Forest of graph expansions & $\mathbb{F}_\mathcal{R}$ & \cpageref{txt.forest.of.expansions}
    \\
    Forest expansion & $F \triangleleft e$ & \cref{def.forest.expansion}
    \\
    Base forest & $F_0$ & \cref{sub.forests.and.graphs}
    \\
    Replacement tree & $T_c$ & \cref{sub.forests.and.graphs}
    \\
    Leaf graph & & \cref{def.leaf.graph}
    \\
    \multicolumn{3}{>{\hsize=\dimexpr3\hsize+4\tabcolsep+2\arrayrulewidth\relax}X}{\textbf{Rearrangement groups}}
    \rule{0pt}{20pt}
    \\
    \hline
    Canonical homeomorphism & & \cref{def.canonical.homeomorphism}
    \\
    Rearrangement & $g$ & \cref{def.rearrangement}
    \\
    Graph pair diagram & $(D,\phi,R)$ & \cref{def.graph.pair.diagrams}
    \\
    Expansion of graph pair diagram & & \cref{sub.reduced.graph.pair.diagrams}
    \\
    Reduced graph pair diagram & & \cref{sub.reduced.graph.pair.diagrams}
    \\
    Composition of graph pair diagrams & & \cref{sub.composition.graph.pair.diagrams}
    \\
    Undirected color, undirected edge & & \cref{def.undirected.edges}
    \\
    Isolated color & & \cref{def.null.expanding.isolated}
    \\
    Null-expanding color & & \cref{def.null.expanding.isolated}
    \\
    Forest pair diagram & $(F_D, \phi, F_R)$ & \cref{def.forest.pair.diagram}
    \\
    Reduced forest pair diagrams & & \cpageref{txt.reduced.forest.pair.diagrams}
    \\
    Composition of forest pair diagrams & & \cpageref{txt.composition.forest.pair.diagrams}
    \\
    \multicolumn{3}{>{\hsize=\dimexpr3\hsize+4\tabcolsep+2\arrayrulewidth\relax}X}{\textbf{Examples of rearrangement groups}}
    \rule{0pt}{20pt}
    \\
    \hline
    Thompson's group $F$ & $F$ & \cref{sub.F}
    \\
    Thompson's group $T$ & $T$ & \cref{sub.T}
    \\
    Thompson's group $V$ & $V$ & \cref{sub.V}
    \\
    Higman-Thompson groups & $F_{n,r}$, $T_{n,r}$, $V_{n,r}$ & \cref{sub.higman.thompson.groups}
    \\
    Airplane rearrangement groups & $T_A$ & \cref{sub.airplane}
    \\
    Basilica edge replacement system & $\mathcal{B}$ & \cref{fig.basilica.replacement}
    \\
    Basilica rearrangement group & $T_B$ & \cref{sub.basilica}
    \\
    Topological full group of an edge shift & $V(\Gamma,v_0)$ & \cref{sub.topological.full.groups}
    \\
    Thompson-like groups $QF$, $QT$, $QV$ & $QF$, $QT$, $QV$ & \cref{sub.thompson.like}
    \\
    Houghton groups & $H_n$ & \cref{sub.Houghton}
    \\
    \multicolumn{3}{>{\hsize=\dimexpr3\hsize+4\tabcolsep+2\arrayrulewidth\relax}X}{\textbf{Dendrites}}
    \rule{0pt}{20pt}
    \\
    \hline
    Dendrite, subdendrite & & \cref{def:dendrite}
    \\
    Order of a point & & \cpageref{txt.order.of.points}
    \\
    Wa\.zewski dendrite & $D_n$ & \cref{def:wazewski:dendrite}
    \\
    Dendrite edge replacement system & $\mathcal{D}_n$ & \cref{sec:dendrite:rearrangement:groups}
    \\
    Dendrite rearrangement group & $G_n$ & \cref{sec:dendrite:rearrangement:groups}
    \\
    Branch point, branch & & \cref{def:branch}
    \\
    Endpoint & & \cref{def:ext}
    \\
    Set of branch points or endpoints in $X$ & $\mathrm{Br}(X)$, $\mathrm{En}(X)$ & \cref{sub:points:arcs}
    \\
    Rational endpoint & & \cpageref{txt.rationa.endpoints}
    \\
    Set of rational endpoints in $X$ & $\mathrm{REn}(X)$ & \cpageref{txt.rationa.endpoints}
    \\
    Arc & $[p,q]$ & \cref{def:arc}
    \\
    Rational arc & & \cref{def:arc:type}
    \\
    BB-arc, EE-arc, BE-arc & & \cref{def:arc:type}
    \\
    \multicolumn{3}{>{\hsize=\dimexpr3\hsize+4\tabcolsep+2\arrayrulewidth\relax}X}{\textbf{Dendrite rearrangement groups}}
    \rule{0pt}{20pt}
    \\
    \hline
    Group of local permutations at $p$ & $K_p$ & \cref{sub:perm}
    \\
    Main copy Thompson's $F$ & $H$ & \cref{sub:thomp}
    \\
    Set of subsets of $X$ with $k$ elements & $X^{(k)}$ & \cref{sub.transitive.properties.thompson.subgroups}
    \\
    Copy of Thompson's $F$ on the arc $A$ & $H_A$ & \cref{lem:thomps}
    \\
    Homeomorphism group of $D_n$ & $\mathbb{H}_n$ & \cref{sec:dns}
    \\
    Minimal subdendrite containing $\mathcal{F}$ & $[\mathcal{F}]$ & \cpageref{txt:trees}
    \\
    Tree of $\mathcal{F}$ & $T(\mathcal{F})$ & \cpageref{txt:trees}
    \\
    Oligomorphic group action & & \cpageref{txt.oligomorphic}
    \\
    Local parity map & $\pi_p$ & \cref{sub:parity}
    \\
    Parity map & $\Pi$ & \cpageref{txt.global.parity.map}
    \\
    Local endpoint derivative & $\partial_q$ & \cref{sub.endpoint.derivative}
    \\
    Endpoint derivative & $\Delta$ & \cpageref{txt.global.endpoint.derivative}
    \\
    Rigid stabilizer of $U$ in $G$ & $\mathrm{Rist}_G(U)$ & \cref{sub.simplicity.commutator}
    \\
    Oriented dendrite edge replacement system & $\mathcal{D}_n^+$ & \cref{sub:orientation}
    \\
    Orientation-preserving dendrite rearrangement group & $G_n^+$ & \cpageref{txt.oriented.dendrite.rearrangements}
    \\
    Generalized Wa\.zewski dendrite & $D_S$ & \cref{sub:gen:den}
    \\
    \multicolumn{3}{>{\hsize=\dimexpr3\hsize+4\tabcolsep+2\arrayrulewidth\relax}X}{\textbf{Strand diagrams}}
    \rule{0pt}{20pt}
    \\
    \hline
    Conjugacy problem & & \cref{cha.conjugacy}
    \\
    Strand diagram & & \cref{def.SDs}
    \\
    Split, merge & & \cref{def.SDs}
    \\
    Branching strand & & \cpageref{txt.branching.strands}
    \\
    Symbol generated by a split or merge & & \cpageref{txt.branching.strands}
    \\
    $R$-branhing strand diagram & & \cref{def.r.branching}
    \\
    Replacement groupoid & & \cref{def.replacement.groupoid}
    \\
    Reduction of a strand diagram & & \cref{def.SD.reductions}
    \\
    $\mathcal{X}$-strand diagram & & \cpageref{TXT X-SDs}
    \\
    Composition of strand diagrams & & \cref{sub.SDs.composition}
    \\
    Split, merge and permutation diagrams & & \cref{SUB groupoid generators}
    \\
    \multicolumn{3}{>{\hsize=\dimexpr3\hsize+4\tabcolsep+2\arrayrulewidth\relax}X}{\textbf{Closed strand diagrams}}
    \rule{0pt}{20pt}
    \\
    \hline
    $R$-branching closed strand diagrams & & \cref{SUB closed strand diagrams}
    \\
    Base points, base line, base graph & & \cpageref{txt.base.things}
    \\
    Closure of a strand diagram $f$ & $\llbracket f \rrbracket$ & \cpageref{txt.base.things}
    \\
    Open strand diagram of $\eta$ & $o(\eta)$ & \cpageref{txt.opening.closed.strand.diagrams}
    \\
    Similarities & & \cref{sub.similarities}
    \\
    Permutation of the base line & & \cref{sub.similarities}
    \\
    Shift of the base line & & \cpageref{TXT shifts}
    \\
    Expanding or reducing shifts & & \cpageref{TXT shifts}
    \\
    Type 1 and type 2 reductions & & \cref{sub.closed.strand.diagrams.reductions}
    \\
    Looping strand & & \cpageref{TXT 3 reductions}
    \\
    Type 3 reduction & & \cpageref{TXT 3 reductions}
    \\
    Reduction system & & \cref{SUB reduction systems}
    \\
    Reduction-confluence & & \cref{def.reduction.confluent}
    \\
    Infinite power of a strand diagram & $Z^\infty$ & \cref{SUB stable and vanishing}
    \\
    Main base line of $Z^\infty$ & & \cref{SUB stable and vanishing}
    \\
    Stable and vanishing symbols & & \cref{def.stable.vanishing}
    \\
    DPO graph rewriting system & & \cref{sub.DPO}
    \\
    \multicolumn{3}{>{\hsize=\dimexpr3\hsize+4\tabcolsep+2\arrayrulewidth\relax}X}{\textbf{Dynamics and invariable generation}}
    \rule{0pt}{20pt}
    \\
    \hline
    Invariable generation & IG & \cref{cha.IG}
    \\
    CO-transitive group action & & \cpageref{txt.CO.transitivity}
    \\
    Weak cell-transitivity & & \cref{def.wct}
    \\
    Minimal group action & & \cpageref{txt.minimal}
    \\
    Flexible group action & & \cref{prop.orbit.dense}
    \\
    Vigorous group action & & \cpageref{txt.vigorous}
    \\
    (Weakly) $g$-wandering set & & \cref{def.wandering}
    \\
    Caret & & \cpageref{txt.caret}
    \\
    Difference of forests & $F_D - F_R$ & \cpageref{txt.forest.difference}
    \\
    Domain imbalance & & \cref{def.domain.imbalance}
    \\
    Range imbalance & & \cref{rmk.imbalances}
    \\
    Components of $F_D - F_R$ & & \cpageref{txt.components}
    \\
    Expandable sequence & & \cpageref{txt.expandable.sequence}
    \\
    \multicolumn{3}{>{\hsize=\dimexpr3\hsize+4\tabcolsep+2\arrayrulewidth\relax}X}{\textbf{Rationality of the gluing relation}}
    \rule{0pt}{20pt}
    \\
    \hline
    Rational relation & & \cref{def.rational.relation}
    \\
    Gluing automaton of $\mathcal{R}$ & $\mathrm{Gl}_\mathcal{R}$ & \cref{sub.gluing.automaton}
    \\
    States of the gluing automaton & $q_0(i)$, $q_1 \big(\begin{smallmatrix} i & \gamma \\ j & \delta \end{smallmatrix}\big)$ & \cref{sub.states.gluing.automaton}
    \\
    Transitions of the gluing automaton & & \cref{sub.transitions.gluing.automaton}
    \\
    Type of adjacency & & \cref{tab:transition:type:1}
    \\
    \multicolumn{3}{>{\hsize=\dimexpr3\hsize+4\tabcolsep+2\arrayrulewidth\relax}X}{\textbf{Operations on edge replacement systems}}
    \rule{0pt}{20pt}
    \\
    \hline
    Product of $\mathcal{R}$ and $\mathcal{S}$ & $\mathcal{R} \oplus \mathcal{S}$ & \cref{sub.products}
    \\
    Edge replacement systems for stabilizers & $\mathcal{R}^S$, $\mathcal{R}^p$, $\mathcal{R}^q$ & \cref{sec.rational.stabilizers}
    \\
    Binary edge shift & & \cref{def.binary.edge.shift}
\end{tabularx}
\endgroup

%%%%%%%%%%%%%%%%%%%%%%%%%

\mainmatter
\pagestyle{mainstyle}
\setcounter{chapter}{-1}
\chapter{Introduction}
%\addcontentsline{toc}{chapter}{Introduction}
%\markboth{Introduction}{Introduction}

Highly symmetric topological spaces usually have rich groups of homeomorphisms and, more broadly, there are plenty of groups that act on such spaces by homeomorphisms with interesting dynamics and with fascinating algebraic and combinatorial or geometric properties.

\begin{figure}\centering
\includegraphics[width=.6\textwidth]{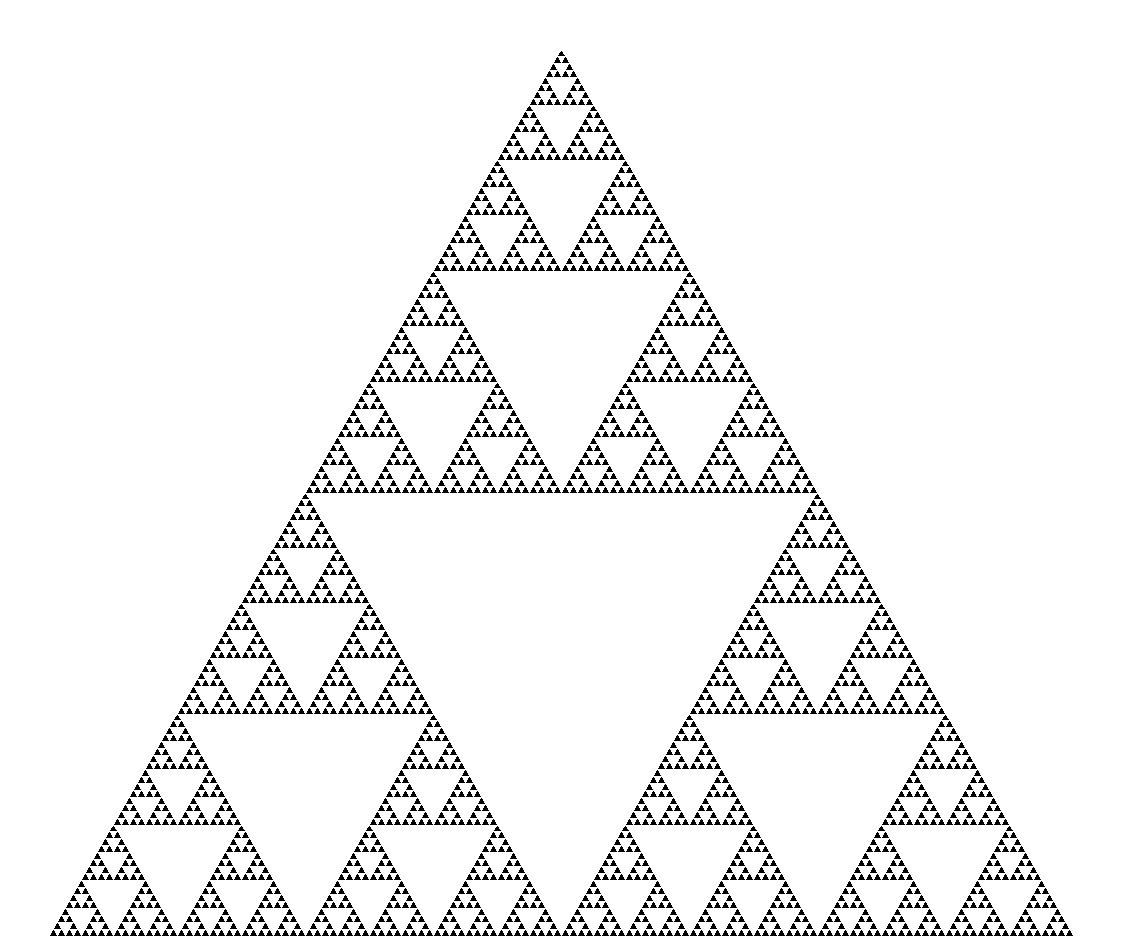}
\caption{The Sierpinski triangle.}
\label{fig.sierpinski.triangle}
\end{figure}

When looking for such spaces, it is natural to consider \textbf{self-similar fractals}, by which we informally mean a space with a structure of nested subspaces that have finitely many homeomorphism types.
For example, the \textbf{Cantor space} has a complicated uncountable group of homeomorphisms that has been studied under various topologies \cite{CantorTopologies}.
\textbf{Dendrites}, which we will discuss in \cref{cha.dendrites}, also have rich and huge (uncountable) groups of homeomorphisms and group actions by homeomorphisms on dendrites have been studied extensively in the literature \cite{UniversalDendrites,Shi12,AmenableDendrites,DM18,DM19,Kaleid,Duc20,FiniteOrbits}.
On the other hand, there exist self-similar fractals with tiny groups of homeomorphisms:
for instance, the group of homeomorphism of the famous Sierpinski triangle (depicted in \cref{fig.sierpinski.triangle}) consists of only six elements and is isomorphic to the dihedral group $D_3$, which is the group of symmetries of a regular triangle \cite{Louwsma}.
In general, it seems that groups acting on fractals can be expected to range from finite to uncountable, which makes their study even more compelling.

For uncountable groups of homeomorphisms of fractals, studying combinatorial, computational and algorithmic properties is evidently meaningless.
Fortunately, many of these spaces can be ``codified'' in a way that prompts the construction of natural discrete subgroups of the homeomorphism groups.
In general, the \hyperref[thm.Alexandroff.Hausdorff]{Alexandroff-Hausdorff Theorem} states that a compact metrizable space is a quotient of the Cantor space, meaning that any compact metrizable fractal can be built as a set of infinite words on finite alphabets modulo an equivalent relation that ``glues'' some of the words.
When such a ``gluing'' relation is particularly nice, it is natural to look at those homeomorphisms that descend from certain actions on the Cantor space of infinite words;
one simple such action that already produces promising results is that of \textbf{prefix-exchange transformations}, which is precisely the key idea behind rearrangement groups of fractals.

\subsection*{Graph Approximations of Fractals}

\begin{figure}\centering
\begin{subfigure}{.32\textwidth}
\centering
\begin{tikzpicture}[scale=1]
\draw[blue] (-2,0) -- (2,0);
\draw[red,fill=white] (0,0) circle (0.5);
\end{tikzpicture}
\end{subfigure}
\hfill
\begin{subfigure}{.32\textwidth}
\centering
\begin{tikzpicture}[scale=1]
\draw[blue] (-2,0) -- (2,0);
\draw[blue] (0,1) -- (0,-1);
\draw[red,fill=white] (0,0) circle (0.5);
\draw[red,fill=white] (-1.25,0) circle (0.25);
\draw[red,fill=white] (1.25,0) circle (0.25);
\end{tikzpicture}
\end{subfigure}
\hfill
\begin{subfigure}{.32\textwidth}
\centering
\begin{tikzpicture}[scale=1]
\draw[blue] (-2,0) -- (2,0);
\draw[blue] (0,1) -- (0,-1);
\draw[blue] (-1.25,0.5) -- (-1.25,-0.5);
\draw[blue] (1.25,0.5) -- (1.25,-0.5);
\draw[blue] (-0.6,-0.6) -- (0.6,0.6);
\draw[blue] (-0.6,0.6) -- (0.6,-0.6);
\draw[red,fill=white] (0,0) circle (0.5);
\draw[red,fill=white] (-1.25,0) circle (0.25);
\draw[red,fill=white] (1.25,0) circle (0.25);
\draw[red,fill=white] (0.75,0) circle (0.083);
\draw[red,fill=white] (1.75,0) circle (0.083);
\draw[red,fill=white] (-0.75,0) circle (0.083);
\draw[red,fill=white] (-1.75,0) circle (0.083);
\draw[red,fill=white] (0,0.75) circle (0.083);
\draw[red,fill=white] (0,-0.75) circle (0.083);
\end{tikzpicture}
\end{subfigure}
\caption{Graph approximations of the airplane Julia sets.}
\label{fig.airplane.approximation}
\end{figure}
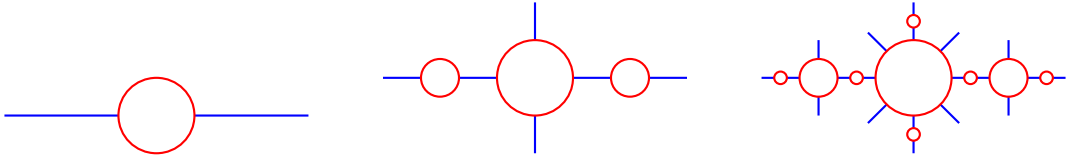

\begin{figure}\centering
\includegraphics[width=.6\textwidth]{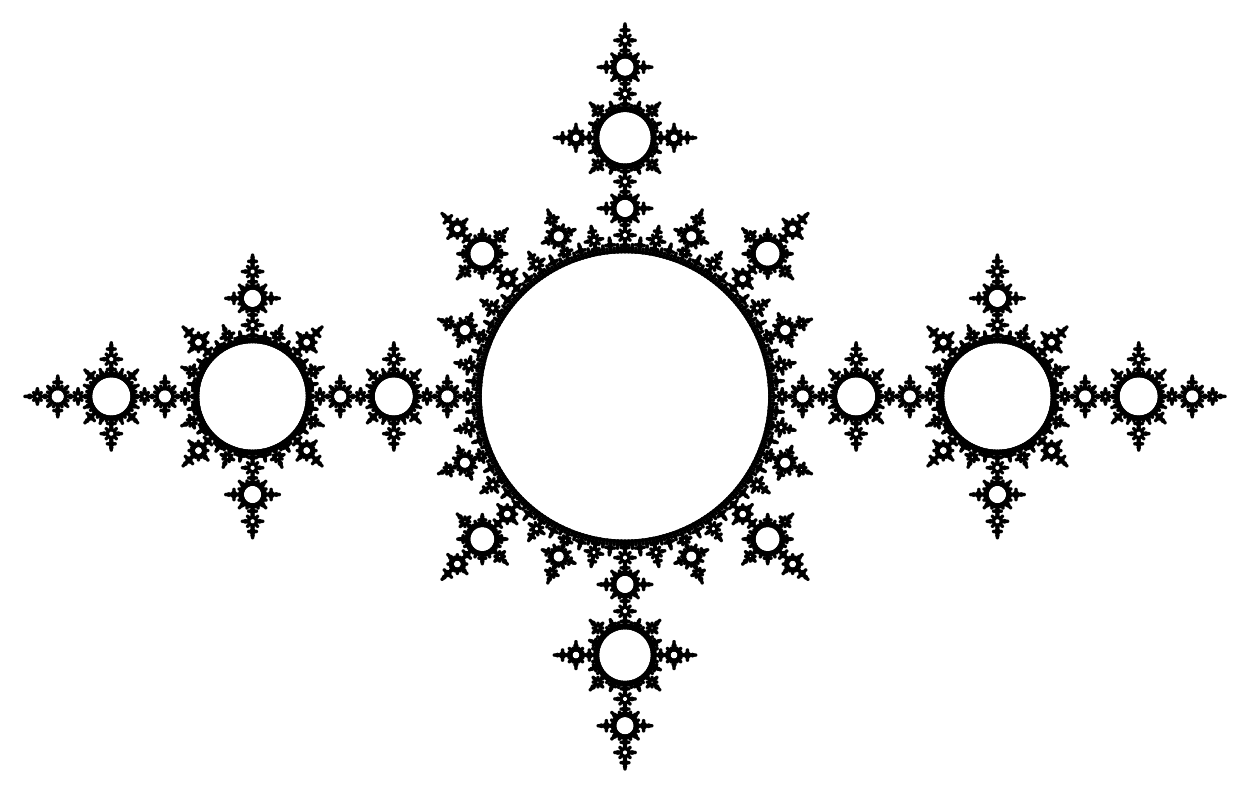}
\caption{The airplane limit space.}
\label{fig.airplane.intro}
\credits{J. Belk and B. Forrest, \cite[Figure 27]{BF19}}
\end{figure}

In 2019 James Belk and Bradley Forrest introduced \textbf{edge replacement systems}, which are special instances of graph rewriting systems \cite{BF19}.
Under mild conditions, an edge replacement system allows us to approximate a self-similar topological space with increasingly more complex graphs.
For example, \cref{fig.airplane.approximation} shows three such graphs that approximate the so-called airplane Julia set portrayed in \cref{fig.airplane.intro}.

The graph approximations produced by edge replacement systems do not just visually resemble a fractal space:
they actually define an equivalence ``gluing'' relation $\sim$ on a space $\Omega$ of infinite words, and the quotient $\Omega / \sim$ is homeomorphic to a fractal.
This equivalence relation is particularly well-behaved in the sense that, as will be proven in \cref{cha.rational.gluing}, it is \textbf{rational}:
there exists an automaton that reads precisely those pairs of infinite words that are glued together.

\subsection*{Prefix-exchange Transformations}

It turns out that, since the family of graph approximations produced by an edge replacement system defines the equivalence relation $\sim$, graph isomorphisms between them induce homeomorphisms of the fractal $\Omega / \sim$.
Edges of these graphs are finite prefixes of the infinite words of $\Omega$, so graph isomorphisms act on the set of infinite words $\Omega$ simply by changing finitely many prefixes and this action descends to an action on the fractal $\Omega / \sim$ by permutation of finitely many self-similar subspaces.
A homeomorphism of the fractal $\Omega / \sim$ that is induced by such a prefix-exchange transformation of $\Omega$ is what we call a \textbf{rearrangement} of that fractal.
Despite only being countable, groups of rearrangements of fractals are often dense in a natural overgroup of homeomorphisms of the fractal on which they act.
This will be discussed in \cref{cha.dendrites}, which shows how the rearrangement groups of dendrites can approximate homeomorphisms of dendrites arbitrarily well, and the same holds for the interval, the circle and the Cantor space.

Groups that act in this way include the famous three \textbf{Thompson groups} $F$, $T$ and $V$, which consist of prefix-exchange homeomorphisms of the unit interval, the unit circle and the Cantor space, respectively.
Introduced in the 60's by mathematician Richard Thompson for his studies in logic, they immediately caught the attention of group theorists such as Graham Higman for their finiteness properties:
all of them are finitely presented, which means that they can be entirely described by a finite amount of generators and relations between them.
Furthermore, $T$ and $V$ are simple, and they were the first known examples of infinite simple finitely presented simple groups.
Groups with these properties are still hard to produce, and generalizations of $F$, $T$ and $V$ are still the main source of such groups.

A key conjecture in combinatorial and geometric group theory is the \textbf{Boone-Higman conjecture}, which states that a finitely generated group has solvable word problem if and only if it embeds into a finitely presented simple group.
The recent progress and increasing interest in this conjecture (see for example \cite{BooneHigman}) highlights the importance of finding and studying finitely presented simple groups.
Moreover, generalizations of Thompson groups have produced examples of simple groups with arbitrary higher-dimensional finiteness properties \cite{SWZ19}.
This shows promise towards the study and classification of infinite simple groups, which is arguably an even harder endeavor than the already mammoth task of classifying finite simple groups.

Thompson groups are often studied with a dynamical approach.
It is natural then that similar methods have proved fruitful for the study of rearrangement groups of fractals as well, as is the case for the arguments used throughout \cref{cha.IG}.
Thompson groups also show good algorithmic behaviour, such as having solvable \textbf{conjugacy problem} \cite{BM14}.
In \cref{cha.conjugacy} we show that similar methods also apply to a vast class of rearrangement groups of fractals.
This also highlights that many rearrangement groups, despite often being large enough to approximate homeomorphisms with arbitrary precision, are still computationally manageable.

\subsection*{Structure of this dissertation}

The first three chapters are mostly introductory:
\cref{cha.background} provides a few basic notions and fixes some notation for the rest of the dissertation.
\cref{cha.fractals,cha.rearrangements} contain the necessary background on edge replacement systems and rearrangements.
Except for a few statements of \cref{cha.fractals} and for \cref{sec.null.expanding,sub.undirected.edges}, these chapters do not feature new results that are due to the author, but they offer a novel and unique treatment on many aspects of this topic.
For instance, \cref{sub.symbol.space} describes a fully explicit connection between the symbol space and edge shifts and \cref{sec.examples} exhibits a wide list of rearrangement groups of which some were not initially known to belong to this family.

\cref{cha.dendrites,cha.conjugacy,cha.IG,cha.rational.gluing} are adaptations of publications and preprints written by the author during his PhD (\cite{dendrite,conjugacy,IG,rationalgluing}, respectively), while \cref{cha.additional.results} exhibits other results that are unpublished at the time of writing this dissertation.
The content of \cref{cha.IG,cha.rational.gluing,sec.rational.stabilizers} is the fruit of a collaboration between the author and Davide Perego, while everything else from the fourth chapter onwards is due to the author alone.
More specifically:
\begin{description}
    \item[\cref{cha.dendrites}] is about rearrangement groups of dendrites and it offers an example of an extensive study of an infinite class of Thompson-like groups;
    \item[\cref{cha.conjugacy}] features a solution for the conjugacy problem in a vast family of rearrangement groups with the aid of strand diagrams that represent the elements of the groups and, when ``closed'', their conjugacy classes.
    \item [\cref{cha.IG}] employs dynamical arguments to show that rearrangement groups whose action is very transitive do not have the property of being invariably generated;
    \item [\cref{cha.rational.gluing}] exhibits an algorithmic construction of automata that recognize the gluing relation which defines the fractal spaces as quotients of spaces of infinite words.
\end{description}

%%%%%%%%%%%%%%%%%%%%%%%%%

\chapter{Background}
\label{cha.background}

The main objective of the upcoming \cref{cha.fractals} is to build fractal topological spaces as quotients of edge shifts under equivalence relations defined by graphs.
Before delving into that machinery, this chapter provides the necessary background on graphs and edge shifts.

\section{Directed Graphs}

We are going to frequently use the notion of directed graphs in this dissertation.
Those presented here are standard notions in the theory of directed graphs, possibly except for sinks and sources (\cref{def.sink.source}) and inescapable cycles (\cref{def.inescapable.cycle}).

\subsection{Basic Definitions}

There are plenty of different notions of graphs in the literature.
One of the most common is that of undirected graphs, where $E$ is a set (or a multiset, where multiple instances of the same element are allowed) of unordered pairs of elements of $V$.
We will not use undirected graphs, so from now on we will completely omit terms such as \textit{directed} or \textit{undirected} when talking about graphs, unless it is useful to remind the reader of this distinction.

We refer to \cite{digraphs} for a complete treatment on the topic of directed graphs (sometimes also called \textit{digraphs} for short).
We will not always use their notation, however.
For example, their definition of directed graphs does not allow the presence of parallel edges nor loops.
Instead, what \cite{digraphs} refers to as \textit{pseudographs} is going to be our definition of graphs, which is the following.

\begin{definition}
\label{def.graph}
A \textbf{graph} is a quadruple $(V, E, \iota, \tau)$, where:
\begin{itemize}
    \item $V$ is a nonempty set, whose elements are called \textbf{vertices};
    \item $E$ is a set, whose elements are called \textbf{edges};
    \item $\iota$ and $\tau$ are maps $E \to V$ and for each edge $e \in E$ the vertices $\iota(e)$ and $\tau(e)$ are called the \textbf{initial} and \textbf{terminal} vertices of $e$, respectively.
\end{itemize}
We say that a graph is \textbf{finite} if both $V$ and $E$ are finite sets.
\end{definition}

We say that an edge $e$ is \textbf{incident} on a vertex $v$ when $\iota(e)=v$ or $\tau(e)=v$.
When two edges are incident on a common vertex, we say that they are \textbf{adjacent}.

It is natural to define the standard notions of morphisms and substructures.

\begin{definition}
\label{def.graph.morphism}
Let $\Gamma_1=(V_1,E_1,\iota_1,\tau_1)$ and $\Gamma_2=(V_2,E_2,\iota_2,\tau_2)$ be two graphs.
A \textbf{graph morphism} $\phi \colon \Gamma_1 \to \Gamma_2$ consists of two maps
\[ \phi_V \colon V_1 \to V_2, \quad \phi_E \colon E_1 \to E_2, \]
such that $\iota_2(\phi_E(e))=\phi_V(\iota_1(e))$ and $\tau_2(\phi_E(e))=\phi_V(\tau_1(e))$ for all $e \in E_1$.

A \textbf{graph isomorphism} is a graph morphism $\phi \colon \Gamma_1 \to \Gamma_2$ such that the maps $\phi_E$ and $\phi_V$ are bijections and their inverse maps produce a graph morphism $\phi^{-1} \colon \Gamma_2 \to \Gamma_1$.
\end{definition}

\begin{definition}
\label{def.subgraph}
A \textbf{subgraph} of a graph $(V,E,\iota,\tau)$ is a graph $(V',E',\iota',\tau')$ such that $V' \subseteq V$, $E' \subseteq E$ and $\iota'$, $\tau'$ are the restrictions of $\iota$, $\tau$ to $E'$.
\end{definition}

There is a specific notion of subgraph that is often useful.

\begin{definition}
\label{def.subgraph.spanned}
Given a subset of vertices $V' \subseteq V$, the subgraph \textbf{spanned} by $V'$ is the unique subgraph $(V',E',\iota',\tau')$ with the property that $E'$ is maximal (with respect with the inclusion relation) among all subgraphs with vertex set $V'$.
\end{definition}

Differently from undirected graphs, there are two natural notions of vertex degrees in directed graphs.

\begin{definition}
\label{def.degree}
Given a vertex $v$, its \textbf{out-degree} and \textbf{in-degree} are:
\[ \mathrm{out}(v) \coloneq |\{ e \in E \mid \iota(e)=v \}| \, \text{ and } \, \mathrm{in}(v) \coloneq |\{ e \in E \mid \tau(e)=v \}|. \]
The \textbf{degree} of $v$ is the sum of its in- and out-degrees.
\end{definition}

These notions prompt the following definitions, which will be useful when discussing strand diagrams in \cref{sec.SDs}.

\begin{definition}
\label{def.sink.source}
A \textbf{sink} is a vertex whose out-degree is zero and a \textbf{source} is a vertex whose in-degree is zero.
\end{definition}

Furthermore, it is useful to give the two following definitions.

\begin{definition}
\label{def.loop}
A \textbf{loop} is an edge $e$ such that $\iota(e)=\tau(e)$.
\end{definition}

\begin{definition}
\label{def.isolated.vertex}
We say that a vertex $v$ is \textbf{isolated} if its in-degree and out-degree are zero, i.e., if there is no edge having $v$ as its initial or terminal vertex.
\end{definition}

\begin{assumption}
\label{ass.isolated.vertices}
From here on, we will always assume that graphs do not have isolated vertices.
\end{assumption}

This is a natural assumption for all of the graphs in this dissertation.
Indeed, isolated vertices ``disappear'' in the limit spaces of edge replacement systems (which will be introduced in \cref{cha.fractals});
in automata (\cref{sub.automata}), isolated vertices are inactive states without any use;
in rewriting systems (\cref{sub.rewriting.systems}), isolated vertices represent objects that are not involved in any rewriting.

Under this assumption, a graph isomorphism $\phi$ (\cref{def.graph.morphism}) is always entirely determined by its permutation of edges $\phi_E$.
Hence it will be common for us to think of the isomorphism simply as the single bijection $\phi_E$, completely forgetting about the permutation of vertices $\phi_V$ because it is entirely determined by $\phi_E$.

Finally, it is often useful to endow vertices with an ordering of their incident edges.

\begin{definition}[\cite{GraphsOnSurfaces}]
\label{def.rotation.system}
A \textbf{rotation system} on a graph is an assignment of a circular order to the edges incident on each vertex.
\end{definition}

Throughout this dissertation, rotation systems will be represented simply by the way in which drawing a graph locally embeds it into the plane around each vertex.

\subsection{Walks, Paths and Cycles}

It is natural and very useful to introduce the following sequences of edges.

\begin{definition}
\label{def.walk.path.cycle}
Let $S = (e_1, e_2, \dots)$ be a (finite or infinite, but nonempty) sequence of edges.
\begin{itemize}
    \item $S$ is a \textbf{walk} if $\tau(e_i) = \iota(e_{i+1})$ for every $i$.
    \item $S$ is a \textbf{path} if it is a walk devoid of repetitions of edges (i.e., $\nexists i \neq j$ such that $e_i = e_j$). 
    \item $S$ is a \textbf{cycle} if it is a finite path $(e_1, \dots, e_n)$ such that $\tau(e_n)=\iota(e_1)$.
\end{itemize}
If $S = (e_1, \dots, e_n)$, its \textbf{length} is $n$;
if $S$ is infinite, we say that its length is infinite.
Moreover, a finite subsequence $(e_i, e_{i+1}, \dots)$ or an infinite subsequence $(e_i, \dots, e_{i+k})$ of a walk (or path) $S = (e_1, e_2, \dots)$ is called a \textbf{subwalk} (or \textbf{subpath}) of $S$.
\end{definition}

Even if walks, paths and cycles are formally defined as sequences of edges, it is often natural to think of them as subsets of a graph.
Because of this, we say that a walk, path or cycle \textbf{contains a vertex} if that vertex is initial or terminal for some edge of the walk, path or cycle.

\begin{definition}
\label{def.inescapable.cycle}
An \textbf{inescapable cycle} is a cycle $C$ such that every walk that starts from a vertex of $C$ is entirely contained in $C$.
\end{definition}

Walks and paths motivate the following notions of connectivity in graphs.

\begin{definition}
\label{def.reachable}
Given two vertices $v$ and $w$, we say that $w$ is \textbf{reachable} from $v$ if there exists a finite walk (or, equivalently, a finite path) $(e_1, \dots, e_n)$ that goes from $v$ to $w$, i.e., such that $v = \iota(e_1)$ and $w = \tau(e_n)$.
\end{definition}

Reachability defines a binary relation that is transitive but generally not reflexive nor symmetric, and it is often used to define the notion of \textit{strong connectivity}.
We are instead going to need a weaker notion of connectivity that does not take into account the orientation of edges.

\begin{definition}
\label{def.biorientation}
Given a graph $\Gamma = (V,E,\iota,\tau)$, its \textbf{biorientation} $\overline{\Gamma}$ is the graph $(V,\overline{E}, \overline{\iota},\overline{\tau})$ where we define $\overline{E} = \{ e^{+}, e^{-} \mid e \in E \}$ and the maps $\overline{\iota}$ and $\overline{\tau}$ are defined as
\[ \overline{\iota}(e^+) = \iota(e),\, \overline{\tau}(e^+) = \tau(e),\, \overline{\iota}(e^-) = \tau(e),\, \overline{\tau}(e^-) = \tau(e). \]
\end{definition}

In essence, in the biorientation of a graph walks can travel in the opposite direction of edges.
This definition allows us to ignore the orientation of edges and treat them as if they did not have their direction from the initial to the terminal vertex.

\begin{definition}
\label{def.connected.graph}
Two vertices of a graph $\Gamma$ are \textbf{connected} if they are the same vertex or if there is a walk (or, equivalently, a path) between them in the biorientation $\overline{\Gamma}$.
\end{definition}

This last definition prompts a natural \textbf{connectivity relation} on the set of vertices.
It is reflexive by definition and it is easy to show that it is symmetric and transitive.
Hence it is an equivalence relation, which prompts the following definition.

\begin{definition}
\label{def.connected.component.graph}
A \textbf{connected component} of a graph $\Gamma$ is a subgraph spanned by an equivalence relation of the connectivity relation (\cref{def.subgraph.spanned}).
A graph is \textbf{connected} if it features a unique connected component.
\end{definition}

\subsection{Trees and Forests}

Trees have always inspired fruitful interactions with the world of group theory.
In this dissertation we will make frequent use of them.
Indeed, even if it is not stated formally here, we will often implicitly use the fact that prefix-exchange homeomorphisms of edge shifts (which make up the so-called topological full groups of edge shifts in which rearrangement groups embed, see \cref{sub.topological.full.groups}) are ``finitary'' quasi-automorphisms of forests.

\begin{definition}
\label{def.reduced.cycle}
Consider the biorientation $\overline{\Gamma}$ of a graph $\Gamma = (V,E,\iota,\tau)$.
A cycle $(e_1, \dots, e_n)$ is \textbf{reduced} if there is no $i \in \mathbb{Z}/n\mathbb{Z}$ such that $e_{i+1} = e_i^{-1}$.
\end{definition}

\begin{definition}
\label{def.forest.tree}
A \textbf{forest} is a graph whose biorientation does not have reduced cycles and a \textbf{tree} is a connected forest.
\end{definition}

\phantomsection\label{txt.roots.leaves}
In this dissertation every forest is going to be made of a finite number of trees and the edge orientation is always going to be pointing away from a unique source of each tree.
More precisely, our forests are going to be such that, for each source $r$ and each vertex $v$ that lie in the same connected component, there is a unique path from $r$ to $v$ (note that a path in the biorientation of a graph $\Gamma$ need not be a path in $\Gamma$).
In this context, we say that $r$ is a \textbf{root} and that the forest is \textbf{rooted}.
In a similar fashion, sinks are going to be called \textbf{leaves}.

\begin{definition}
\label{def.complete.subforest}
Given a rooted forest, a \textbf{complete subforest} is a subgraph that includes every root and such that
\begin{enumerate}
    \item if $v$ belongs to the subgraph, it also includes the unique path from a root to $v$;
    \item whenever the subgraph includes an edge starting from some vertex $v$, it includes every other such edge from $v$.
\end{enumerate}
\end{definition}

\section{Other Notions Related to Graphs}

This section contains the relevant background on three additional topics:
labelings and colorings of graphs, rewriting systems and automata.

\subsection{Colorings and Labelings of Graphs}

It is often useful to endow graphs with additional data from some set $X$.
This is usually done by some sort of labeling, which is simply a map from the set of edges or vertices to $X$.
We will only need edge-labelings in this dissertation, which are defined right below.

\begin{definition}
\label{def.labeled.graph}
A \textbf{labeled graph} is a graph $\Gamma = (V,E,\iota,\tau)$ together with a labeling map $l \colon E \to X$.
We say that $X$ is the set of \textbf{labels} and, if $l(e) = L$, we say that $e$ is \textbf{labeled} by $L$.
\end{definition}

Sometimes, especially when it is finite, the set of labels is identified with a set of colors and the labeling is conveniently depicted by coloring the edges.
In this case, we refer to the labeled graph as a \textbf{colored graph}.
It will be convenient to write $\mathrm{c} \colon E \to \mathrm{C}$ as the map that associates to each edge $e$ its color $\mathrm{c}(e)$.

Moreover, when the set of labels is a cartesian product $A \times B$, the labeling map $l$ can be naturally decomposed as $l_A \colon E \to A$ and $l_B \colon E \to B$, which allows us to think of $l$ as if it were two distinct labelings.
We are often going to tacitly use this decomposition to easily distinguish between a coloring and an additional symbolic labeling, for example starting from \cref{cha.fractals}, where graph expansions of edge replacement systems will be colored and also labeled by finite words, and later in \cref{cha.conjugacy}, where strand diagrams will be both colored and labeled by pairs of vertices and an additional symbol.

Even if there is no formal distinction between generic labelings and colorings, in this dissertation we are always going to need graph morphisms to preserve a given coloring but not necessarily a given labeling.
Thus, when dealing with colored graphs we will use the following definition.

\begin{definition}
\label{def.colored.graph.morphism}
Consider two graphs $\Gamma_1$ and $\Gamma_2$ that are colored by the same set of colors $\mathrm{C}$, and let $c_i \colon E_i \to \mathrm{C}$ be their coloring maps ($i = 1,2$).
In addition to the requirements described in \cref{def.graph.morphism} a \textbf{graph morphism} $\phi \colon \Gamma_1 \to \Gamma_2$ needs to satisfy $c_1(e) = c_2(\phi_E(e))$ for every $e \in E_1$.
\end{definition}

Finally, when the set of colors is a singleton (which is going to be the case for Thompson groups and dendrite rearrangement groups discussed in \cref{sub.Thompson} and \cref{cha.dendrites}, respectively), the coloring map is trivial and the coloration is automatically preserved by any morphism.
In this case, we will completely omit any mention of colors.

\subsection{Rewriting Systems}
\label{sub.rewriting.systems}

Graphs are often used to describe certain types of relations on the set of vertices.
This is the case for (abstract) rewriting systems, which consist of a set of objects and rewriting rules that allow us to transform an object to a (possibly) different one.

In general, an (abstract) \textbf{rewriting system} is simply a directed graph, where we say that the vertices are \textbf{objects} of the rewriting system and the edges are \textbf{rewritings}.
It is convenient to write $a \longrightarrow b$ to denote a rewriting from $a$ to $b$, which is simply an edge whose initial and terminal vertices are the objects $a$ and $b$, respectively.

\begin{example}
\label{ex.rewriting}
Let $V = \{0,1\}^* = \bigcup_{n \geq 0} \{0,1\}^n$ be the set of finite words in the alphabet $\{0,1\}$, included the empty word $\varepsilon$.
This is going to be our set of objects.
We introduce rewritings that remove the last digit of an object, i.e., $E = \{ w0 \overset{0}{\longrightarrow} w \mid w \in V \} \cup \{ w1 \overset{1}{\longrightarrow} w \mid w \in V \}$.
As a graph, this rewriting system is simply a tree where every vertex has out-degree $1$ and in-degree $2$, except for the empty word $\varepsilon$ which instead has out-degree $0$.
\end{example}

In practice, it is often convenient to refer to a graph as a rewriting system when there are clear rule that establish when
an object is rewritten to another.
There might be different rewriting rules, so it is often useful to write $a \overset{R}{\longrightarrow} b$ to specify that the rewriting follows a rule $R$.
This is essentially just a labeling of the edges of the graph with the set of rewriting rules.

\begin{definition}
\label{def.locally.confluent}
A rewriting system is \textbf{locally confluent} if, whenever there are reductions $s \longrightarrow a$ and $s \longrightarrow b$, there exists an object $t$ and walks (or, equivalently, paths) from $a$ to $t$ and from $b$ to $t$.
\end{definition}

Local confluence means that, when an object can be rewritten in different ways, there exist sequences of subsequent rewritings that allow us to reach a common object.

In the terminology of rewriting systems, sinks (vertices of out-degree $0$) are usually called \textbf{reduced objects}, since they cannot be rewritten.

\begin{definition}
\label{def.terminating}
A rewriting system $(V,E,\iota,\tau)$ is \textbf{terminating} if it does not feature any infinite walk.
\end{definition}

Termination essentially means that the rewritings rules allow us to eventually transform any object to a reduced one.
However, there might be multiple reduced objects that are reachable from the same object.

The following Lemma from \cite{NewmanDiamond} is a classical and fundamental result in the theory of rewriting systems and it is going to be of use multiple times in this dissertation.
It is often used to find normal forms and reduced elements in a variety of different settings.

\medskip %layout
\begin{lemma}[Newman's Diamond Lemma]
\label{lem.diamond}
If a rewriting system is locally confluent and terminating, then each of its connected components features one and only one reduced object.
\end{lemma}

The previous \cref{ex.rewriting} is both locally confluent (simply because there is at most one rewriting from each object) and terminating (because the length of each word is reduced to a minimum).
Thus, \hyperref[lem.diamond]{Newman's Diamond Lemme} tells us that there is a unique reduced object.
In fact, the empty word is clearly the unique reduced element of that rewriting system.

\subsection{Automata}
\label{sub.automata}

In this subsection we briefly introduce the notion of automata, which is going to be useful in \cref{cha.rational.gluing}.
There are many kinds of automata, but the so-called \textit{synchronous deterministic finite-state automaton} is arguably the most commonly used.
Since we will only be using this kind of automata in this dissertation, we will omit these words, the meaning of which is discussed later in \cref{rmk.automata.definition}.

\phantomsection\label{txt.infinite.sequences.symbol}
In order to give definition of automata, recall that a \textbf{partial function} is a binary relation between two sets that associates to each element of the first set at most one element of the second. 

\begin{definition}
\label{def.automaton}
Let $\Sigma$ be a finite set called \textit{alphabet}.
A (synchronous deterministic finite-state) \textbf{automaton} is a quadruple $(\Sigma, Q, \rightarrow, q_{0})$ with $Q$ a finite set whose elements are called \textbf{states}, $q_{0} \in Q$ called \textbf{initial state} and a partial function $\rightarrow$ from $Q \times \Sigma$ to $Q$ which is called \textbf{transition function}.
We write $q \xrightarrow{x} q'$ if the transition function features $(q,x) \mapsto q'$.
\end{definition}

From here onwards we will denote by $\Sigma^\omega$ the set of all infinite sequences in an alphabet $\Sigma$.

\begin{definition}
\label{def.run.recognize}
A \textbf{run} of the automaton is a (finite or infinite) sequence of transitions
\[ q_0 \xrightarrow{x_1} q_1 \xrightarrow{x_2} \dots.\]
An automaton \textbf{recognizes} a finite word or an infinite sequence $x_1 x_2 \dots$ when there exists such a run.
We say that an automaton \textbf{recognizes} a subset $\Omega \subseteq \Sigma^\omega$
when it recognizes a sequence if and only if it belongs to $\Omega$.
\end{definition}

Even if these definitions are not directly related to graphs, an automaton $\left( \Sigma, Q, \rightarrow, q_{0} \right)$ can be thought of as finite labeled graph with a distinguished vertex.
For this comparison, the set of vertices is $Q$, the distinguished vertex is $q_0$ and each transition $q \xrightarrow{x} q'$ corresponds to an edge from $q$ to $q'$ labeled by $x \in \Sigma$.
The requirement that $\rightarrow$ is a partial function is equivalent to assuming that, for every $x \in \Sigma$, each vertex is the origin of at most one edge labeled by $x$.
Under this point of view, a run is a (finite or infinite) walk starting at $q_0$.

With this graph-theoretic notion of automaton in mind, the reader may sense a strong similarity between recognizability and edge shifts (which will be defined shortly in \cref{def.edge.shift}).
Nevertheless, these two concepts do not coincide and one should instead resort to the more general notion of \textit{sofic shifts} (discussed for example in \cite[Chapter 3]{LindMarcus}) in order to give a symbolic dynamical definition of recognizability.

\begin{remark}
\label{rmk.automata.definition}
There is an extensive literature on a plethora of different kinds of automata, of which the synchronous deterministic finite-state ones only represent a specific family.
To understand what makes them special, we briefly expose the meaning of these words.
\begin{itemize}
    \item \textit{Synchronous} means that edges are labeled by elements of $\Sigma$, whereas finite words in the alphabet $\Sigma$ could be employed in more general notions of automata.
    \item \textit{Deterministic} means that the transition function must be a partial function instead of any binary relation.
    \item \textit{Finite-state} simply means that the set of states $Q$ is finite.
\end{itemize}
\end{remark}

\section{Edge Shifts}

In this section we introduce edge shifts.
These objects are central in symbolic dynamics, but here we will see them as ways to explicitely construct compact metrizable and totally disconnected topological spaces that are made of certain infinite sequences of symbols.
For more information on the symbolic dynamical aspects of edge shifts, we refer the reader to \cite{LindMarcus}.

\phantomsection\label{txt.word.sequence}
For the sake of clarity, throughout this dissertation we will generally use the term \textit{word} to refer to finite words and the term \textit{sequence} for infinite words.

\subsection{The Cantor space}
\label{sub:Cantor}

Traditionally, the Cantor space $\mathfrak{C}$ is obtained as a limit of the following iterative procedure.

\begin{enumerate}
    \item[0)] Start with the unit interval $\mathfrak{C}_\varepsilon = [0,1]$.
    \item[1)] Remove from $\mathfrak{C}_\varepsilon$ the middle third $(1/3, 2/3)$, which results in the union of two disjoint sets
    \[ \mathfrak{C}_{0} = [0, 1/3] \text{ and } \mathfrak{C}_{1} = [2/3, 1]. \]
    \item[2)] Remove the middle thirds of both $\mathfrak{C}_0$ and $\mathfrak{C}_1$, which results in the union of four disjoint sets
    \[ \mathfrak{C}_{00} = [0, 1/9],\, \mathfrak{C}_{01} = [2/9, 1/3],\, \mathfrak{C}_{10} = [2/3, 7/9],\, \mathfrak{C}_{11} = [8/9, 1]. \]
    \item[n)] Remove the middle thirds from each $\mathfrak{C}_w$, where $w$ spans among the $2^{n-1}$ words of length $n-1$ in the alphabet $\{0,1\}$, which results in the union of $2^n$ disjoint sets $\mathfrak{C}_{w'}$, where $w'$ spans among the $2^n$ words of length $n$ in the alphabet $\{0,1\}$.
\end{enumerate}
The resulting topological subspace $\mathfrak{C}$ of $[0,1]$ seems small at first glance, since it looks like we have removed almost everything from $[0,1]$.
However, it features exactly one point for each infinite sequence in the alphabet $\{0,1\}$, so it is actually an uncountable set.
The Cantor space was a surprising example of an uncountable subset of $\mathbb{R}$ whose (Lebesgue) measure equals zero.

As hinted behind the notation of this construction, it is not hard to see that the Cantor set is homeomorphic to the set $\{0,1\}^\omega$ of all infinite sequences in the alphabet $\{0,1\}$, if endowed with the product topology.
By Tychonoff's theorem, it is thus a compact space.
Being a subspace of the metric space $[0,1]$, the Cantor space is also a metric space.
It is easy to see that it is perfect (i.e., it does not contain isolated points) and totally disconnected (i.e., its connected components are the singletons).
In truth, a converse of the previous assertion holds by the following classical result (see for example \cite[Theorem 30.3]{Willard}):

\medskip %layout
\begin{theorem}[Brouwer's Theorem]
\label{thm.brouwer}
A nonempty topological space is homeomorphic to the Cantor space if and only if it is compact, metrizable, perfect and totally disconnected.
\end{theorem}

This allows us to say that \textit{a} Cantor space is any topological space with these properties, and although these are all homeomorphic, sometimes it is more convenient to use \textit{one} Cantor space in place of \textit{another}.
This will be clear, for example, when in \cref{sec.limit.spaces} we will build limit spaces of colored edge replacement systems as quotients of \textit{certain} specific Cantor spaces defined in \cref{sub.edge.shifts}.

\subsection{Codings of Compact Metrizable Spaces}

Our interest towards the Cantor space is mainly motivated by the following classical theorem (see for example \cite[Theorem 30.7]{Willard}):

\medskip %layout
\begin{theorem}[Alexandroff-Hausdorff Theorem]
\label{thm.Alexandroff.Hausdorff}
Every compact metrizable space is a quotient of the Cantor space.
\end{theorem}

The fractal spaces that we are interested in are in fact compact and metrizable, so the \hyperref[thm.Alexandroff.Hausdorff]{Alexandroff-Hausdorff Theorem} essentially allows us to codify them.
More precisely, given a compact metrizable space $X$, there will always be a quotient map $\pi \colon \mathfrak{C} = \{0,1\}^\omega \to X$ which continuously associates to each infinite binary sequence a point of $X$.
We will thus refer to such a map $\pi$ as a \textbf{coding} of the space $X$, and here we give the simplest classical example.

\begin{example}
\label{ex.interval}
Let $X = [0,1]$ be the unit interval.
Consider the map
\[ \pi \colon \mathfrak{C} = \{0,1\}^\omega \to X = [0,1] \]
that associates to each infinite binary sequence $\alpha = \alpha_1 \alpha_2 \dots$ the point $\sum_{i \geq 1} \alpha_i 2^{-i}$ of $[0,1]$.
It is easy to verify that this is the quotient map associated to the equivalence relation $w 1 \overline{0} \sim w 0 \overline{1}$ for each finite binary word $w$ (where by $\overline{v}$ we mean the infinite periodic sequence $v v v \ldots$ for any finite word $v$).
\end{example}

Note that the \hyperref[thm.Alexandroff.Hausdorff]{Alexandroff-Hausdorff Theorem} allows us to say that there exists such a map $\pi$, but it does not provide us with one.
However, the specific map $\pi$ from the previous \cref{ex.interval} is particularly nice in a sense that will be described in \cref{cha.rational.gluing}, but in general the equivalence relation that defines the quotient map $\pi$ could be much more unwieldy.
In \cref{sec.limit.spaces} we will define specific equivalence relations on Cantor spaces that are determined by edge replacement systems.

Before we do that, we need to describe a way to explicitly construct and handle a specific family of Cantor spaces that arises from the theory of symbolic dynamical systems.

\subsection{Definition of Edge Shifts}
\label{sub.edge.shifts}

See \cite[\S 2.2]{LindMarcus} for more information about the symbolic dynamical aspects of edge shifts.

\begin{definition}
\label{def.edge.shift}
Given a finite graph $\Gamma = (V,E,\iota,\tau)$ and a specific vertex $v_0 \in V$, its (one-sided initial) \textbf{edge shift} $\Omega(\Gamma, v_0)$ is the set of infinite walks on $\Gamma$ starting at $v_0$, i.e.,
\[ \Omega(\Gamma, v_0) = \{ e_1 e_2 \dots \mid \iota(e_1) = v_0 \text{ and } \tau(e_n) = \iota(e_{n+1}) \}. \]
\end{definition}

Observe that sinks are essentially ignored by the edge shift, since there cannot be infinite walks that travel through sinks.

\begin{remark}
More generally, as done in \cite{LindMarcus}, ``two-sided'' edge shifts are defined as bi-infinite walks on the graph $\Gamma$.
Moreover, one could consider one-sided ``non initial'' edge shifts without requiring that the infinite walks start at $v_0$.
However, in this dissertation we will only make use of one-sided initial edge shifts, so from now on we will refer to one-sided initial edge shifts simply as edge shifts, for the sake of brevity.
There is also a more general notion of \textit{subshifts of finite type}, but in a certain sense these can be recoded to edge shifts (see \cite[Theorem 2.3.2]{LindMarcus}).
\end{remark}

\begin{example}
\label{ex.full.shift}
Given a finite set $E$, the \textbf{full shift} on $E$ is the Cantor space
\[ E^\omega = \{ \alpha_1 \alpha_2 \dots \mid \alpha_i \in E \}. \]
The full shift is an edge shift for the graph $(V,E,\iota,\tau)$ whose vertex set is a singleton $V = \{v_0\}$ and such that $\iota(e) = \tau(e) = v_0$ for each $e \in E$.
\end{example}

In general, edge shifts are the set of sequences that satisfy the restriction of being readable on a finite graph $\Gamma$.
The edge shifts that we are going to consider are going to be called \textit{symbol spaces} and are introduced later in \cref{sub.symbol.space}.
Each edge replacement system (\cref{def.replacement.system}) will define such an edge shift, so additional examples will be unveiled in the next chapter.

Since edge shifts are sets of certain infinite sequences whose digits are edges of a graph, it is natural to use the following nomenclature.

\begin{definition}
\label{def.alphabet.language}
Given an edge shift $\Omega(\Gamma,v_0)$, its \textbf{alphabet} is the set of edges of $\Gamma$ and its \textbf{language} is the set of finite (non-empty) prefixes of elements of $\Omega(\Gamma,v_0)$.
\end{definition}

\subsection{The Topology of Edge Shifts}

Let $E$ be a finite set with the discrete topology and endow $E^\omega$ with the product topology.
It is known and straightforward to see that each full shift $E^\omega$ (\cref{ex.full.shift}) is a Cantor space.

Given a finite graph $\Gamma$ with set of edges $E$, any edge shift $\Omega(\Gamma, v_0)$ is a subset of $E^\omega$.
We endow each edge shift $\Omega(\Gamma, v_0)$ with the subspace topology induced from the Cantor space $E^\omega$.
The following topological properties of $\Omega(\Gamma, v_0)$ are easily proven by observing that edge shifts are closed subspaces of the Cantor space $E^\omega$.

\medskip %layout
\begin{proposition}
\label{prop.shift.is.nice}
Edge shifts are compact, metrizable and totally disconnected topological spaces.
\end{proposition}

The following Proposition expands on the previous one.
We are going to use the notion of \textit{inescapable cycles} of graphs (\cref{def.inescapable.cycle}).

\medskip %layout
\begin{proposition}
\label{prop.shift.is.Cantor}
An edge shift $\Omega(\Gamma,v_0)$ is a Cantor space if and only if $\Gamma$ features no inescapable cycle that is reachable from $v_0$.
In particular, if the out-degree of every vertex that is reachable from $v_0$ is at least two then the edge shift $\Omega(\Gamma,v_0)$ is a Cantor space.
\end{proposition}

\begin{proof}
In light of \cref{prop.shift.is.nice}, it suffices to prove that $\Omega(\Gamma,v_0)$ is perfect if and only if there is no inescapable cycle that is reachable from $v_0$.
To do this, it is convenient to consider the following standard metric on the Cantor space $E^\omega$:
\[ d(\alpha,\beta) = 2^{-\delta(\alpha,\beta)} \]
where $\delta(\alpha,\beta)$ denotes the length of the maximum common prefix of $\alpha$ and $\beta$ (possibly empty and thus of length zero, or infinite when $\alpha=\beta$).
The topology of $\Omega(\Gamma,v_0)$ described above (induced from the inclusion in the Cantor space $E^\omega$) is then induced by the restriction of the metric $d$.

If $\alpha$ and $\beta$ belong to the edge shift, from a graph-theoretic point of view $\delta(\alpha,\beta)$ is simply the length of the longest common subwalk (\cref{def.walk.path.cycle}) of $\alpha$ and $\beta$ in $\Gamma$ starting from $v_0$.

For a point $\alpha$ of $\Omega(\Gamma,v_0)$ to be isolated then means that there exists a $C>0$ such that the radius-$C$ ball of $\Omega(\Gamma,v_0)$ centered at $\alpha$ consists solely of $\alpha$.
This is equivalent to asserting that, for any infinite walk on $\Gamma$ starting from $v_0$, if the walk follows $\alpha$ for at least $-\mathrm{log}_2 (C)$ edges, then that walk must be precisely $\alpha$.
This means precisely that $\alpha$ ends in some inescapable cycle.
Indeed, since $\Gamma$ is finite, if $\alpha = a_1 a_2 \dots$ does not end in an inescapable cycle then there are arbitrarily large $n \in \mathbb{N}$ such that $a_n$ has out-degree at least $2$ and so there are distinct elements $a_1 \dots a_n a_{n+1}$ and $a_1 \dots a_n x$ of $\Omega(\Gamma,v_0)$, which implies the existence of elements of $\Omega(\Gamma,v_0)$ that are arbitrarily close to $\alpha$.
\end{proof}

Edge shifts are ``tassellated'' by special subsets commonly referred to as cones (sometimes also called cylinders in the literature), which are defined right below.

\begin{definition}
\label{def.cones}
Given an edge shift $\Omega(\Gamma,v_0)$ and an element $x$ of its language, the \textbf{cone} determined by $x$ is the subspace
\[ C(x) = \{ x \alpha \mid x \alpha \in \Omega(\Gamma,v_0) \}. \]
\end{definition}

It is famously known that the set of cones of a Cantor space $E^\omega$ (seen as a full shift as in \cref{ex.full.shift}) is a basis of clopen sets.
The same holds for every edge shift, since it inherits its topology from the Cantor space $E^\omega$ and a cone for the edge shift $\Omega$ is the intersection of $\Omega$ with a cone for $E^\omega$.

\medskip %layout
\begin{proposition}
\label{prop.basis.clopen}
The set of cones of $\Omega(\Gamma,v_0)$ is a basis of clopen sets.
\end{proposition}

If an edge shift is not a Cantor space, meaning that it contains an isolated point, then a cone is a Cantor space if and only if it does not contain an isolated point itself.
Moreover, if $\alpha$ is an isolated point then there exists some $N \in \mathbb{N}$ such that, for every prefix of $\alpha$ of length greater or equal to $N$, the corresponding cone is the singleton $\{\alpha\}$.

%%%%%%%%%%%%%%%%%%%%%%%%%

\chapter{From Graphs to Fractals}
\label{cha.fractals}

In this chapter we introduce edge replacement systems.
This machinery, defined in \cite{BF19}, allows us to build a plethora of fractal topological spaces as quotients of edge shifts under a gluing relation that is determined by sequences of graphs.
We will also see how to represent this sequence of graphs with so-called forest expansions.

\section{Edge replacement systems}
\label{sub.replacement.systems}

\begin{definition}
\label{def.replacement.rules}
A set of \textbf{edge replacement rules} is a pair $(R, \mathrm{C})$, with $R = \{X_c \mid c \in \mathrm{C} \}$, where
\begin{itemize}
    \item $\mathrm{C}$ is a finite set of colors;
    \item for each color $c \in \mathrm{C}$, $X_c \in R$ is a finite graph edge-colored by $\mathrm{C}$ and equipped with distinct vertices $\iota_c$ and $\tau_c$.
\end{itemize}
Each $X_c$ is called \textbf{$c$ replacement graph} (for example, \textit{red replacement graph}) and the vertices $\iota_c$ and $\tau_c$ are called the \textbf{initial} and \textbf{terminal vertices} of $X_c$, respectively.
\end{definition}

\begin{definition}
\label{def.replacement.system}
An \textbf{edge replacement system} $(X_0, R, \mathrm{C})$ consists of a set of edge replacement rules $(R, \mathrm{C})$ together with a finite graph $X_0$ that is edge-colored by $\mathrm{C}$.
The graph $X_0$ is called the \textbf{base graph} of the edge replacement system.
\end{definition}

\begin{remark}
One could define much more general replacement rules (systems) where entire subgraphs are replaced instead of single edges.
Miguel Del R\'io's PhD dissertation defines \textit{graph diagram groups}, which are rearrangement-like groups that arise in this way \cite{graphdiagramgroups}.
Those groups are more complicated and do not usually act on a fractal-like limit space.
\end{remark}

\phantomsection\label{txt.mono.poly.chromatic}
If the set of colors is a singleton, we say that the edge replacement system is \textbf{monochromatic} and we omit any mention of the coloring.
If there is more than one color, it is natural to say that the edge replacement system is instead \textbf{polychromatic}.

\begin{figure}
\centering
\begin{tikzpicture}
    \node at (0,1.6) {$X_0$};
    \node[vertex] (l) at (-.75,0) {};
    \node[vertex] (r) at (.75,0) {};
    \draw[edge,blue] (l) to node[above]{$s$} (r);
    \begin{scope}[xshift=3.5cm]
    \draw[edge,red,domain=5:175] plot ({.5*cos(\x)}, {.5*sin(\x)});
    \draw (90:.5) node[above,red] {$b_2$};
    \draw[edge,red,domain=185:355] plot ({.5*cos(\x)}, {.5*sin(\x)});
    \draw (270:.5) node[below,red] {$b_3$};
    \node at (0,1.6) {$X_{\text{\textcolor{blue}{b}}}$};
    \node[vertex] (l) at (-1.75,0) {}; \draw (-1.75,0) node[above]{$\iota_{\text{\textcolor{blue}{b}}}$};
    \node[vertex] (cl) at (-.5,0) {};
    \node[vertex] (cr) at (.5,0) {};
    \node[vertex] (r) at (1.75,0) {}; \draw (1.75,0) node[above]{$\tau_{\text{\textcolor{blue}{b}}}$};
    \draw[edge,blue] (cl) to node[above]{$b_1$} (l);
    \draw[edge,blue] (cr) to node[above]{$b_4$} (r);
    \end{scope}
    \begin{scope}[xshift=7.5cm]
    \node at (0,1.6) {$X_{\text{\textcolor{red}{r}}}$};
    \node[vertex] (l) at (-1.25,0) {}; \draw (-1.25,0) node[above]{$\iota_{\text{\textcolor{red}{r}}}$};
    \node[vertex] (r) at (1.25,0) {}; \draw (1.25,0) node[above]{$\tau_{\text{\textcolor{red}{r}}}$};
    \node[vertex] (c) at (0,0) {};
    \node[vertex] (ct) at (0,1) {};
    \draw[edge,red] (l) to node[below]{$r_1$} (c);
    \draw[edge,red] (c) to node[below]{$r_2$} (r);
    \draw[edge,blue] (c) to node[left]{$r_3$} (ct);
    \end{scope}
\end{tikzpicture}
\caption{The airplane edge replacement system.}
\label{fig.airplane.replacement}
\end{figure}
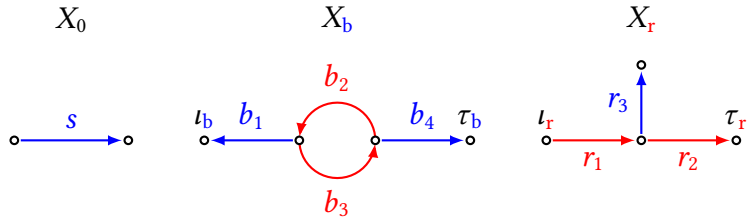

\phantomsection\label{txt.airplane.replacement.system}
\cref{fig.airplane.replacement} depicts the so-called \textbf{airplane edge replacement system}, which we are going to denote by $\mathcal{A}$ and use as a guiding example throughout this section, as it features enough complexity while still being easy to visualize.
It was introduced in \cite[Example 2.13]{BF19} and the corresponding rearrangement group was studied in \cite{airplane}.
\cref{sec.examples} features many more examples of edge replacement systems.

\subsection{Graph Expansions}

Fix an edge replacement system $\mathcal{R} = (X_0, R, \mathrm{C})$ and, without loss of generality, let $\mathrm{C} = \{ 1, 2, \dots, n \}$ and $R = \{ X_1, X_2, \dots, X_n \}$.

As the names in the previous definitions suggest, the idea behind edge replacement systems is that, starting from the base graph, one can replace each edge with the replacement graph of its color.

\begin{definition}
\label{def.edge.expansion}
Let $\Gamma$ be a graph that is edge-colored by $\mathrm{C}$.
Consider an edge $e$ of $\Gamma$ and let $c$ be its color.
The \textbf{edge expansion} of $\Gamma$ by $e$ is the graph $\Gamma \triangleleft e$ obtained by replacing the edge $e$ with the graph $X_c$, identifying the initial and terminal vertices $\iota_c$ and $\tau_c$ of the replacement graph $X_c$ with the initial and terminal vertices $\iota(e)$ and $\tau(e)$ of the edge $e$ of $\Gamma$, respectively.
The graphs that are obtained from finite sequences of edge expansions of the base graph $X_0$ are called \textbf{graph expansions} of the edge replacement system $(X_0, R, \mathrm{C})$.
\end{definition}

Each edge replacement system then defines a rewriting system of graphs whose set of objects is the set of all (isomorphism classes of) graphs that can be obtained as graph expansions of the base graph and whose rewritings are the edge expansions $X \overset{\mathrm{c}(e)}{\longrightarrow} X \triangleleft e$ for each edge $e$ of $X$, distinguished by types that correspond to colors.

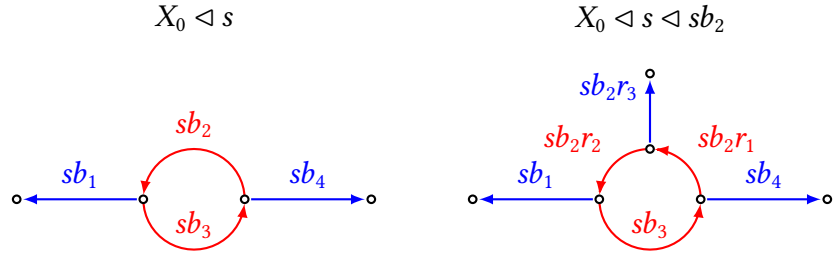
\begin{figure}
\centering
\begin{tikzpicture}
    \node at (0,2.4) {$X_0 \triangleleft s$};
    \draw[edge,red,domain=5:175] plot ({.667*cos(\x)}, {.667*sin(\x)});
    \draw (90:.667) node[above,red] {$sb_2$};
    \draw[edge,red,domain=185:355] plot ({.667*cos(\x)}, {.667*sin(\x)});
    \draw (270:.667) node[above,red] {$sb_3$};
    \node[vertex] (l) at (-2.333,0) {};
    \node[vertex] (cl) at (-.667,0) {};
    \node[vertex] (cr) at (.667,0) {};
    \node[vertex] (r) at (2.333,0) {};
    \draw[edge,blue] (cl) to node[above]{$s b_1$} (l);
    \draw[edge,blue] (cr) to node[above]{$s b_4$} (r);
    \begin{scope}[xshift=6cm]
    \node at (0,2.4) {$X_0 \triangleleft s \triangleleft s b_2$};
    \draw[edge,red,domain=95:175] plot ({.667*cos(\x)}, {.667*sin(\x)});
    \draw (135:.667) node[above left,red] {$sb_2r_2$};
    \draw[edge,red,domain=185:355] plot ({.667*cos(\x)}, {.667*sin(\x)});
    \draw (270:.667) node[above,red] {$sb_3$};
    \draw[edge,red,domain=5:85] plot ({.667*cos(\x)}, {.667*sin(\x)});
    \draw (45:.667) node[above right,red] {$sb_2r_1$};
    \node[vertex] (l) at (-2.333,0) {};
    \node[vertex] (cl) at (-.667,0) {};
    \node[vertex] (c) at (0,.667) {};
    \node[vertex] (ct) at (0,1.667) {};
    \node[vertex] (cr) at (.667,0) {};
    \node[vertex] (r) at (2.333,0) {};
    \draw[edge,blue] (cl) to node[above]{$s b_1$} (l);
    \draw[edge,blue] (cr) to node[above]{$s b_4$} (r);
    \draw[edge,blue] (c) to node[above left]{$s b_2 r_3$} (ct);
    \end{scope}
\end{tikzpicture}
\caption{Two subsequent edge expansions of the base graph of the airplane edge replacement system.}
\label{fig.airplane.expansions}
\end{figure}

For example, \cref{fig.airplane.expansions} depicts the graph expansion $X_0 \triangleleft s$ of the base graph of the airplane edge replacement system and the graph expansion $X_0 \triangleleft s \triangleleft s b_2$.
In general, edges of graph expansions can be naturally identified with certain finite words in this fashion:
if $w = w_1 \dots w_n$ is a word that denotes a $c$-colored edge, then the word $wa$ is an edge of $X_0 \triangleleft w_1 \triangleleft \dots \triangleleft w_n$ if and only if $a$ is an edge of $X_c$.
Thus, whether a letter can be appended to a word $w$ depends precisely on the color $c$ of the last digit of $w$.
This is precisely what edge shifts are for, as we will be clarified in the following subsection.

Before moving there, observe that, whenever $e$ and $f$ are edges of a graph $\Gamma$, the graph expansions $\Gamma \triangleleft e \triangleleft f$ and $\Gamma \triangleleft f \triangleleft e$ are the same.
Moreover, we identify edges of different graph expansion that correspond to the same word.
There is a similar obvious identification of vertices of different graph expansions:
when performing an expansion of an edge $e$ of a graph expansion $\Gamma$, the vertices of the graph $\Gamma \triangleleft e$ are those of $\Gamma$ together with those added by the expansion.

Finally, we give the following useful definition.

\begin{definition}
\label{def.full.expansion}
The \textbf{full expansion sequence} is the sequence of graph expansions $E_0, E_1, E_2, \dots$ obtained by expanding, at each step, every edge of the previous graph, starting from the base graph $E_0 = \Gamma_0,$.
\end{definition}

\subsection{The Symbol Space}
\label{sub.symbol.space}

We start with a definition that will allow us to define a totally disconnected compact metrizable space (often a Cantor space) that expresses how edges are ``nested'' one inside another in graph expansions.

\begin{definition}
\label{def.color.graph}
The \textbf{color graph} $A_\mathcal{R}$ of the edge replacement system $\mathcal{R}$ is the (non-colored) graph consisting of
\begin{itemize}
    \item a vertex $q(c)$ for each color $c \in \mathrm{C}$ and an additional vertex $q(0)$;
    \item an edge labeled by $e$ from $q(\zeta)$ to $q(c)$ for each edge $e$ of $X_\zeta$ that is colored by $c$, where $\zeta \in \{0\} \cup \mathrm{C}$.
\end{itemize}
\end{definition}

In the color graph, $q(0)$ represents the starting state for the base graph $X_0$.

\begin{figure}
\centering
\begin{tikzpicture}
    \node[vertex] (s) at (0,0) {$q(0)$};
    \node[vertex] (b) at (3,0) {$q(\text{\textcolor{blue}{b}})$};
    \node[vertex] (r) at (6,0) {$q(\text{\textcolor{red}{r}})$};
    \draw[edge] (s) to node[above]{$s$} (b);
    \draw[edge] (b) to[loop,out=135, in=90,looseness=7.5] node[left]{$b_1$} (b);
    \draw[edge] (b) to[out=45,in=135] node[above]{$b_2$} (r);
    \draw[edge] (b) to node[above]{$b_3$} (r);
    \draw[edge] (b) to[loop,out=270,in=225,looseness=7.5] node[left]{$b_4$} (b);
    \draw[edge] (r) to[loop,out=90,in=45,looseness=7.5] node[right]{$r_1$} (r);
    \draw[edge] (r) to[out=225,in=315] node[above]{$r_3$} (b);
    \draw[edge] (r) to[loop,out=-45,in=-90,looseness=7.5] node[right]{$r_2$} (r);
\end{tikzpicture}
\caption{The color graph for the airplane edge replacement system.}
\label{fig.airplane.color.graph}
\end{figure}
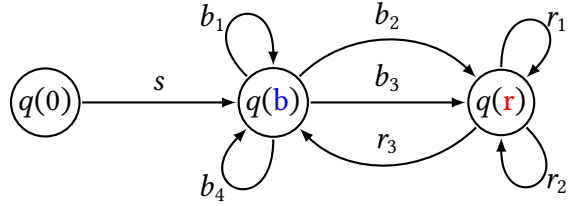

For example, \cref{fig.airplane.color.graph} depicts the color graph for the airplane edge replacement system.
Its edges have been labeled by the edges of the base and replacement graphs (\cref{fig.airplane.replacement}) for ease of reference and in view of the upcoming \cref{prop.language.alphabet}.

Fix a $\zeta \in \{0\} \cup \mathrm{C}$.
The set of edges of $A_\mathcal{R}$ starting from $q(\zeta)$ is precisely the set of edges of the graph $X_\zeta$, which is either the base or a replacement graph.
Each such edge of $A_\mathcal{R}$ terminates at $q(c)$, where $c$ is the color of the edge of $X_\zeta$.
Thus, we have that a word $w = w_1 \dots w_n$ is an edge of some graph expansion of $\mathcal{R}$ if and only if $(w_1, \dots, w_n)$ is a path in $A_\mathcal{R}$ that starts at the vertex $q(0)$.
With \cref{def.edge.shift} in mind and in accordance with \cref{def.alphabet.language}, this can be stated as follows.

\medskip %layout
\begin{proposition}
\label{prop.language.alphabet}
The set of edges that appear among graph expansions of an edge replacement system $\mathcal{R}$ is the language of the edge shift $\Omega(A_\mathcal{R},q(0))$ whose alphabet is the set of edges of the base and replacement graphs.
We denote the language by $\mathbb{L}_\mathcal{R}$ and the alphabet by $\mathbb{A}_\mathcal{R}$.
\end{proposition}

This prompts the following definition.

\begin{definition}
\label{def.symbol.space}
Given an edge replacement system $\mathcal{R}$, its \textbf{symbol space} is the edge shift $\Omega(A_\mathcal{R},q(0))$.
We denote it by $\Omega_\mathcal{R}$.
\end{definition}

For example, the finite words $s b_2 r_1$, $s b_2 r_2$ and $s b_2 r_3$ all belong to the language $\mathbb{L}_\mathcal{A}$ of the airplane edge replacement system, and indeed the respective edges are shown in \cref{fig.airplane.expansions}.
For any infinite sequence $\rho$ in the alphabet $\{r_1, r_2\}$, the infinite sequence $s b_2 \rho$ belongs to the symbol space $\Omega_\mathcal{A}$.
The word $s b_2 b_1$, instead, does not belong to $\mathbb{L}_\mathcal{A}$, since $b_1$ is not an edge of the red replacement graph $X_{\mathrm{c}(s b_2)}$.
Then any infinite sequence that begins with $s b_2 b_1$ cannot belong to $\Omega_\mathcal{A}$.

Note that, for each graph $E_k$ in the full expansion sequence (\cref{def.full.expansion}), the set of its edges is precisely the set of those words in $\mathbb{L}_\mathcal{R}$ that have length $k$.

\begin{remark}
\label{rmk.vertices.as.words}
With the same logic of \cref{def.symbol.space}, we can refer to vertices that appear among graph expansions as finite words $ev$, where the prefix $e \in \mathbb{L}_\mathcal{R}$ represents an edge and $v$ represents some vertex of $\Gamma_{\mathrm{c}(e)}$.

When describing forests (\cref{sec.forest.of.expansions}) and strand diagrams (\cref{sec.SDs}), however, we are not going to be using this convention, since it quickly produces long words that are hard to write in figures and since the names themselves are not important.
We are instead simply going to use a distinct symbol for each vertex:
given a graph with vertices $x$, $y$ and $z$, if an edge expansion produces two new vertices those can be denoted by any two distinct symbols other than $x$, $y$ and $z$.
\end{remark}

\subsection{Expanding edge replacement systems}

Fix an edge replacement system $\mathcal{R}$ and its symbol space $\Omega_\mathcal{R}$.
We will often assume that $\mathcal{R}$ satisfies the mild conditions defined below.

\begin{definition}
\label{def.expanding}
An edge replacement system $(X_0,R,\mathrm{C})$ is \textbf{expanding} if
\begin{enumerate}
    \item neither the base graph nor any replacement graph features isolated vertices (\cref{def.isolated.vertex});
    \item the initial and terminal vertices of the replacement graphs are not connected by an edge;
    \item the replacement graphs have at least three vertices and two edges.
\end{enumerate}
\end{definition}

The airplane edge replacement system $\mathcal{A}$ satisfies all three conditions, so it is expanding.

The first condition ensures that, among all graph expansions, there are never going to be isolated vertices.
Even when dealing with edge replacement systems that are not expanding, however, we are always going to suppose that there are not isolated vertices, as announced in \cref{ass.isolated.vertices}.

By the third condition, since each replacement graph has at least one edge, for each edge $w$ of some graph expansion, there is some element of the symbol space $\Omega_\mathcal{R}$ that starts with $w$.
Moreover, since the third condition states that each replacement graph has at least two edges, given an edge $w$ of some graph expansion there always are at least two edges $wa \neq wb$ for $a \neq b$ edges of $X_{\mathrm{c}(w)}$.
This means that, for each finite path $w$ of $A_\mathcal{R}$, there are always multiple elements of $\Omega_\mathcal{R}$ that have $w$ as a prefix.
This implies that $A_\mathcal{R}$ cannot have inescapable cycles (\cref{def.inescapable.cycle}).
By \cref{prop.shift.is.Cantor}, we then have the following.

\medskip %layout
\begin{proposition}
If an edge replacement system $\mathcal{R}$ is expanding, then its symbol space $\Omega_\mathcal{R}$ is a Cantor space.
\end{proposition}

Together, the three conditions of \cref{def.expanding} are going to be useful in the next subsection, as they imply that a certain \textit{gluing relation} is an equivalence relation.

As will be anticipated in \cref{rmk.generalized.limit.spaces} and developed in \cref{sec.null.expanding}, forgoing the expanding condition sometimes allows us to build limit spaces with isolated points.
There are interesting rearrangement groups arising from non-expanding edge replacement systems, such as the Thompson-like groups $QF$, $QT$ and $QV$ and the Houghton groups $H_n$, which will be discussed is \cref{sub.thompson.like,sub.Houghton}.
These groups can also be built as rearrangement groups of expanding edge replacement systems thanks to \cref{prop.null.expanding.isolated.rarrangement}, even if the resulting limit spaces have less to do with the groups themselves due to the need to replace isolated points with ``rearrangement-rigid'' subspaces.

\section{Limit Spaces}
\label{sec.limit.spaces}

The \textit{limit space} of an edge replacement system is defined as a quotient of the symbol space under a \textit{gluing relation}, as we describe in this subsection.

\begin{definition}
\label{def.glue}
Given an edge replacement system $\mathcal{R}$, the \textbf{gluing relation} of $\mathcal{R}$ is the binary relation on the symbol space $\Omega_\mathcal{R}$ defined by setting $\alpha \sim \beta$ when every finite prefix of $\alpha$ and $\beta$ of equal length represent edges that share at least a vertex.
In symbols this means that, if $\alpha = \alpha_1 \alpha_2 \dots$ and $\beta = \beta_1 \beta_2 \dots$, then
\[ \alpha \sim \beta \iff \forall n \in \mathbb{N},\, \alpha_1 \dots \alpha_n \text{ and } \beta_1 \dots \beta_n \text{ share at least a vertex}. \]
\end{definition}

Intuitively, two sequences of the symbol space are equivalent under the gluing relation when their finite prefixes approach the same vertex.
If two distinct sequences $\alpha$ and $\beta$ of $\Omega_\mathcal{R}$ are glued by $\sim$, then they must share a (possibly trivial) common prefix $x$, after which $x_1 \dots x_k \alpha_{k+1} \dots \alpha_{m}$ and $x_1 \dots x_k \beta_{k+1} \dots \beta_{m}$ must be incident on a common vertex for all $m \geq k$.
For example, in the airplane edge replacement system $\mathcal{A}$ the elements $s b_2 \overline{r_2}$ and $s b_3 \overline{r_1}$ of $\Omega_\mathcal{A}$ are equivalent under the gluing relation, as their prefixes always share the same vertex after the common prefix $s$.
On the other hand, $s b_2 \overline{r_1}$ is not equivalent to the two previous elements of $\Omega_\mathcal{A}$, as it does not share a vertex after their second prefixes.

In general, the gluing relation is not an equivalence relation.
As a simple example, consider a monochromatic edge replacement system where the base graph is a path of three edges $(a,b,c)$ and the replacement graph is a single edge from $\iota$ to $\tau$, which we name $x$.
Note that this is clearly not expanding.
Every graph expansion is a path of three edges $(a x^l, b x^m, c x^n)$ and the symbol space is $\{ a \overline{x}, b \overline{x}, c \overline{x} \}$.
The gluing relation is the following:
$a \overline{x} \sim b \overline{x}$ and $b \overline{x} \sim c \overline{x}$, whereas $a \overline{x} \not\sim c \overline{x}$, so the relation is not transitive.

Fortunately, however, the gluing relation is always an equivalence relation whenever the edge replacement system is expanding (\cref{def.expanding}), as was proved in \cite[Proposition 1.9]{BF19}.
For this very reason, from now on we will work under the following assumption, except when explicitely stated.

\begin{assumption}
Unless explicitly stated, edge replacement systems will be assumed to be expanding.
\end{assumption}

Under this assumption, we can define the limit space as follows.

\begin{definition}
\label{def.limit.space}
Let $\mathcal{R}$ be an expanding edge replacement system with gluing relation $\sim$ and symbol space $\Omega_\mathcal{R}$.
The \textbf{limit space} of $\mathcal{R}$ is the quotient $\Omega_\mathcal{R} / \sim$.
\end{definition}

\begin{figure}\centering
\includegraphics[width=.45\textwidth]{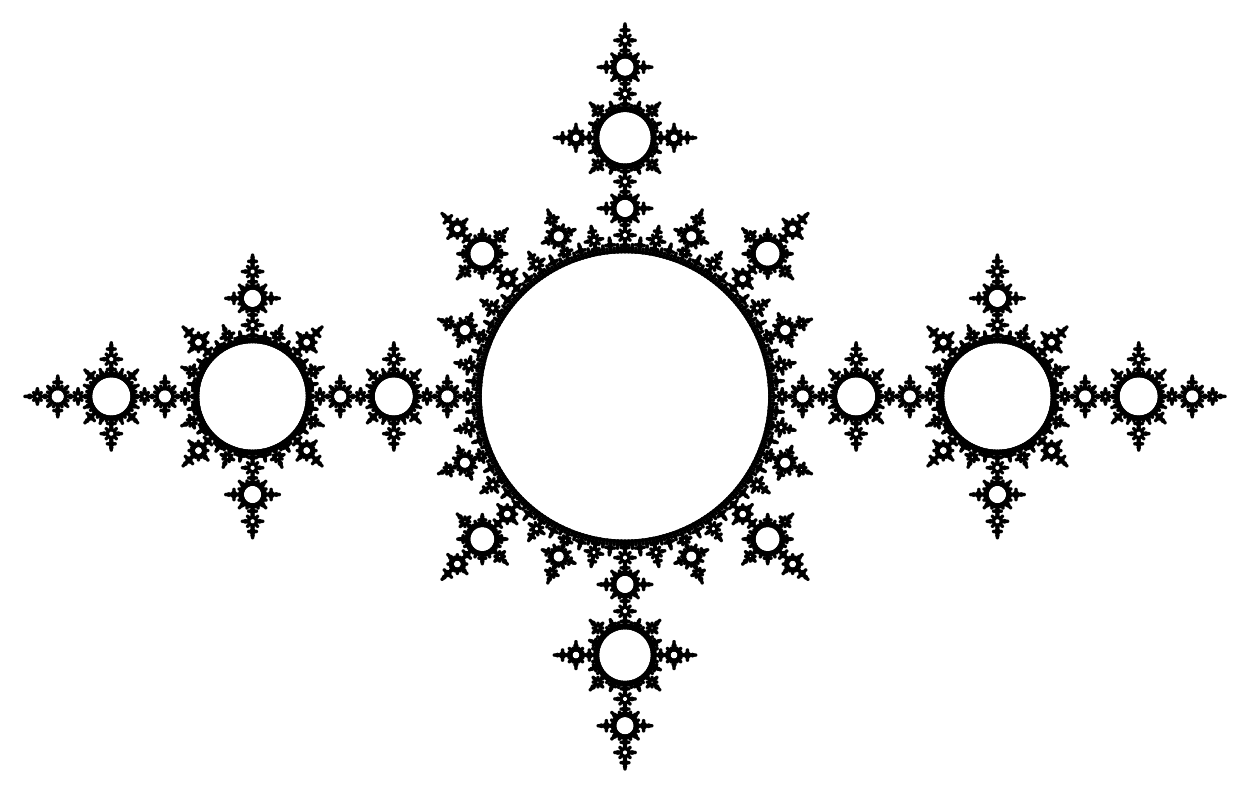}
\caption{The airplane limit space.}
\label{fig.airplane}
\credits{J. Belk and B. Forrest, \cite[Figure 27]{BF19}}
\end{figure}

For example, \cref{fig.airplane} depicts the limit space of the airplane edge replacement system $\mathcal{A}$.
This space is homeomorphic to the boundary of the Julia set for the complex map $f(z) = z^2 + c$, where $c$ is the unique real solution of $c^{3}+2c^{2}+c+1=0$ (approximately $-1.755$).
In short, a Julia set for a complex map $f \colon \mathbb{C} \to \mathbb{C}$ is the subset of $\mathbb{C}$ of those points that have bounded $f$-orbit.
See \cite{JuliaSets} for more details on Julia sets.

When edge replacement systems are expanding, their limit spaces are well-behaved topological spaces, as expressed by the following result from Belk and Forrest.

\medskip %layout
\begin{theorem}[Theorem 1.25 \cite{BF19}]
\label{thm.limit.space}
Limit spaces of expanding edge replacement systems are compact metrizable spaces.
\end{theorem}

\begin{remark}
\label{rmk.generalized.limit.spaces}
Sometimes the limit space may be defined even if the edge replacement system $\mathcal{R}$ is not expanding.
The gluing relation $\sim$ (\cref{def.glue}) is by definition reflexive and symmetric, so the quotient of the edge shift $\Omega_\mathcal{R}$ under the gluing relation $\sim$ can be taken whenever it is transitive.
Even without this assumption, one can always consider the transitive closure $\sim^*$ and generally define the limit space as the quotient $\Omega_\mathcal{R} / \sim^*$, although this is not always a Hausdorff space \cite[Remark 1.23]{BF19}.
This often gives rise to less natural codings for the limit space, but it is going to be useful for example in \cref{sub:orientation}.
\end{remark}

\subsection{Cells of the Limit Space}

\phantomsection\label{txt.point.symbol}
A point of the limit space, which is the equivalence class of a sequence $\alpha \in \Omega_\mathcal{R}$ under the gluing relation $\sim$, will be denoted by $\llbracket \alpha \rrbracket$.
As a natural abuse of notation, given a finite word $w \in \mathbb{L}_\mathcal{R}$ we will write
\[ \llbracket w \rrbracket \coloneq \{ \llbracket w \alpha \rrbracket \mid w \alpha \in \Omega_\mathcal{R} \}, \]
which is the image of the cone $C(w)$ of $\Omega_\mathcal{R}$ (\cref{def.cones}) under the quotient map.

\begin{definition}
\label{def.cell}
Given an edge replacement system $\mathcal{R}$ whose limit space $X$ is defined, a \textbf{cell} of $X$ is a subset $\llbracket w \rrbracket$ of $X$ for a finite word $w \in \mathbb{L}_\mathcal{R}$.
\end{definition}

\begin{figure}\centering
\includegraphics[width=.45\textwidth]{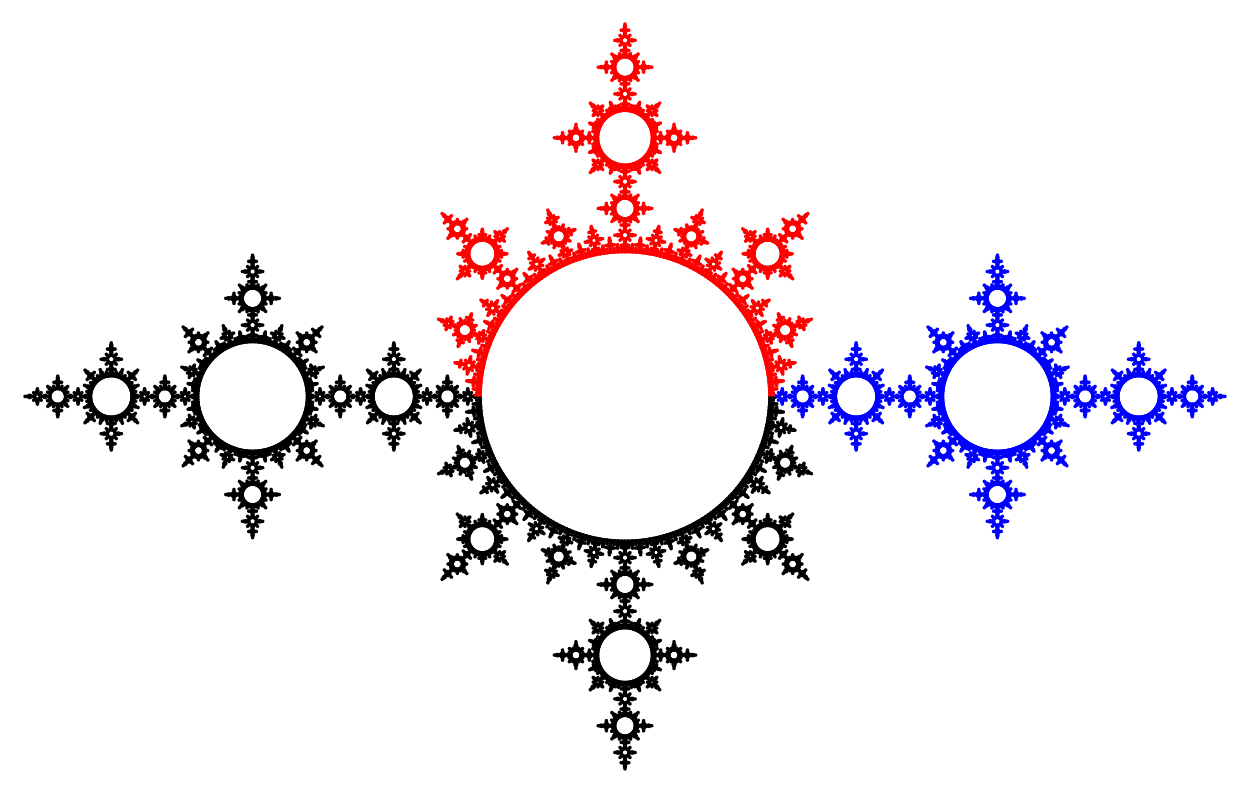}
\caption{Two cells of the airplane limit space.}
\label{fig.airplane.cells}
\end{figure}

For example, \cref{fig.airplane.cells} depicts two cells of the airplane limit space.

With this definition, we have a bijective correspondence between the set of all edges $w$ among all graph expansions and the set of cells $\llbracket w \rrbracket$.
Moreover, it is natural and useful to distinguish cells based on two properties of the originating edge, as expressed by the following definition.

\begin{definition}
\label{def.cell.type}
We say that two cells $\llbracket w_1 \rrbracket$ and $\llbracket w_2 \rrbracket$ have the same \textbf{type} when the edges $w_1$ and $w_2$ are of the same color and they are both loops or neither of them is.
\end{definition}

\phantomsection\label{txt.topological.interior.symbol}
Since cells and their topological interiors will be essential at many points in this dissertation, we will write $\llparenthesis w \rrparenthesis$ for the topological interior of a cell $\llbracket w \rrbracket$.
More on the topological interior is stated later in \cref{rmk.topological.interior}, after vertices of the limit space have been defined.

Cells naturally correspond to self-similar pieces that make up the limit spaces.
They are ``nested'' in the sense expressed right below, which is an excerpt of \cite[Proposition 1.24]{BF19} (we use a different notion of \textit{interior}, but our statement holds for the same arguments given by \cite{BF19}; see \cref{rmk.topological.interior} for more details).

\medskip %layout
\begin{proposition}
Given two distinct cells $\llbracket a \rrbracket$ and $\llbracket b \rrbracket$, either one of the following happens:
\begin{enumerate}
    \item $\llbracket a \rrbracket \subseteq \llbracket b \rrbracket$ and $a$ is a prefix of $b$;
    \item $\llbracket b \rrbracket \subseteq \llbracket a \rrbracket$ and $b$ is a prefix of $a$;
    \item the topological interiors $\llparenthesis a \rrparenthesis$ and $\llparenthesis b \rrparenthesis$ are disjoint and neither $a$ nor $b$ is a prefix of the other.
\end{enumerate}
\end{proposition}

This prompts the following definition.

\begin{definition}
\label{def.cellular.partition}
A \textbf{cellular partition} of the limit space is a collection of cells $\{\llbracket w_1 \rrbracket, \dots, \llbracket w_k \rrbracket \}$ such that $\llbracket w_1 \rrbracket \cup \dots \cup \llbracket w_k \rrbracket$ is the whole limit space and that any pairwise intersection $\llparenthesis w_i \rrparenthesis \cap \llparenthesis w_j \rrparenthesis$ ($i \neq j$) is empty.
\end{definition}

A fundamental feature of cellular partitions is that they correspond in a natural way to graph expansions:
$\{ \llbracket w_1 \rrbracket, \dots, \llbracket w_k \rrbracket \}$ is a cellular partition if and only if $w_1, \dots, w_k$ are the edges of some graph expansion.
Cellular partitions and their correspondence to graph expansions will be essential for the definition of rearrangements in \cref{cha.rearrangements}.

Some cellular partitions are more granular than others in the following sense.

\begin{definition}
\label{def.finer.partitions}
Given two cellular partitions $P_1$ and $P_2$, we say that $P_1$ is \textbf{finer} than $P_2$ if every cell $C \in P_1$ is contained in or is equal to some cell of $P_2$.
\end{definition}

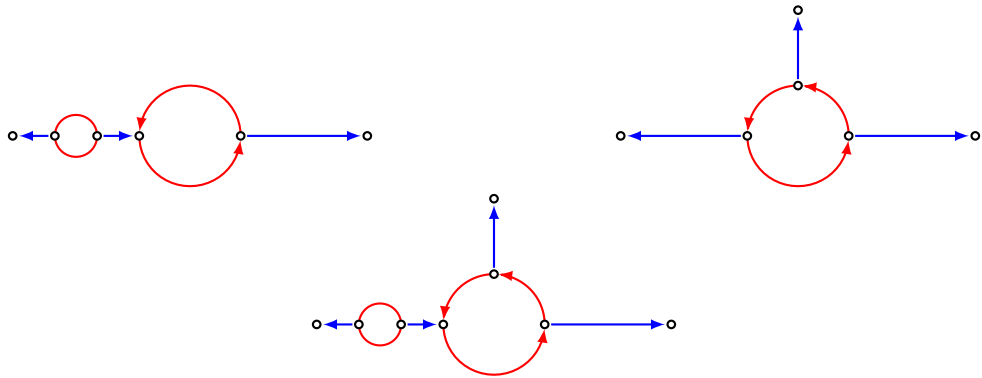
\begin{figure}
\centering
\begin{tikzpicture}
    \begin{scope}[xshift=-4cm]
    \draw[red,domain=10:170] plot ({.278*cos(\x)-1.5}, {.278*sin(\x)});
    \draw[edge,red,domain=5:175] plot ({.667*cos(\x)}, {.667*sin(\x)});
    \draw[edge,red,domain=185:355] plot ({.667*cos(\x)}, {.667*sin(\x)});
    \draw[red,domain=190:350] plot ({.278*cos(\x)-1.5}, {.278*sin(\x)});
    \node[vertex] (l) at (-2.333,0) {};
    \node[vertex] (l1) at (-1.778,0) {};
    \node[vertex] (l2) at (-1.222,0) {};
    \node[vertex] (cl) at (-.667,0) {};
    \node[vertex] (cr) at (.667,0) {};
    \node[vertex] (r) at (2.333,0) {};
    \draw[edge,blue] (l1) to (l);
    \draw[edge,blue] (l2) to (cl);    \draw[edge,blue] (cr) to (r);
    \end{scope}
    \begin{scope}[xshift=4cm]
    \draw[edge,red,domain=5:85] plot ({.667*cos(\x)}, {.667*sin(\x)});
    \draw[edge,red,domain=95:175] plot ({.667*cos(\x)}, {.667*sin(\x)});
    \draw[edge,red,domain=185:355] plot ({.667*cos(\x)}, {.667*sin(\x)});
    \node[vertex] (l) at (-2.333,0) {};
    \node[vertex] (cl) at (-.667,0) {};
    \node[vertex] (c) at (0,.667) {};
    \node[vertex] (ct) at (0,1.667) {};
    \node[vertex] (cr) at (.667,0) {};
    \node[vertex] (r) at (2.333,0) {};
    \draw[edge,blue] (cl) to (l);
    \draw[edge,blue] (cr) to (r);
    \draw[edge,blue] (c) to (ct);
    \end{scope}
    \begin{scope}[yshift=-2.5cm]
    \draw[edge,red,domain=5:85] plot ({.667*cos(\x)}, {.667*sin(\x)});
    \draw[edge,red,domain=95:175] plot ({.667*cos(\x)}, {.667*sin(\x)});
    \draw[edge,red,domain=185:355] plot ({.667*cos(\x)}, {.667*sin(\x)});
    \draw[red,domain=10:170] plot ({.278*cos(\x)-1.5}, {.278*sin(\x)});
    \draw[red,domain=190:350] plot ({.278*cos(\x)-1.5}, {.278*sin(\x)});
    \node[vertex] (l) at (-2.333,0) {};
    \node[vertex] (l1) at (-1.778,0) {};
    \node[vertex] (l2) at (-1.222,0) {};
    \node[vertex] (cl) at (-.667,0) {};
    \node[vertex] (c) at (0,.667) {};
    \node[vertex] (ct) at (0,1.667) {};
    \node[vertex] (cr) at (.667,0) {};
    \node[vertex] (r) at (2.333,0) {};
    \draw[edge,blue] (l1) to (l);
    \draw[edge,blue] (l2) to (cl);
    \draw[edge,blue] (cr) to (r);
    \draw[edge,blue] (c) to (ct);
    \end{scope}
\end{tikzpicture}
\caption{Two cellular partitions of the airplane edge replacement system and their minimal refinement below them, all represented by graph expansions.}
\label{fig.airplane.refinement}
\end{figure}

A simple and intuitive yet useful fact is that two cellular partitions always admit a unique minimal finer common partition, like the ones shown in \cref{fig.airplane.refinement}.

\medskip %layout
\begin{proposition}
\label{prop.common.finer.partition}
Given two cellular partitions $\mathbb{P}_1$ and $\mathbb{P}_2$, there exists a unique least fine cellular partition $\mathbb{P}^*$ that is finer than both $\mathbb{P}_1$ and $\mathbb{P}_2$.
We call it the \textbf{minimal refinement} of $\mathbb{P}_1$ and $\mathbb{P}_2$.
\end{proposition}

The proof is included in that of \cite[Proposition 1.16]{BF19} and essentially consists of exploiting the aforementioned correspondence between cellular partition and graph expansions in order to expand those edges of a cellular partition that are not included in the other until one finds the desired partition.
This proof is constructive, as it describes the actual procedure for computing the minimal refinement.
For example, in \cref{fig.airplane.refinement} one needs to expand the topmost red edge of the first graph expansion and the leftmost blue edge of the second graph expansion to obtain the third one below, which is the minimal refinement.

\subsection{Points of the Limit Space}

In this subsection we introduce and discuss certain special families of points of the limit space.

\subsubsection{Vertices of the Limit Space}

Arguably the more important points of limit spaces are those that descend from vertices of graph expansions.

\begin{definition}
\label{def.vertices}
A \textbf{vertex} of the limit space is a point of the limit space whose representatives are eventually incident on the same vertex of graph expansions.
\end{definition}

\begin{remark}
\label{rmk.glue.implies.vertex}
If a point of the limit space has multiple representatives, then it is always a vertex \cite[Proposition 1.22]{BF19}.
The converse is not true:
for example, $s \overline{b_4}$ is the only representative of the vertex $\llbracket s \overline{b_4} \rrbracket$ in the airplane limit space.
\end{remark}

\phantomsection\label{txt.boundary.vertices}
Depending on whether an edge $w$ is a loop or not, the cell $\llbracket w \rrbracket$ has one or two \textbf{boundary vertices}, which are the points of the limit space corresponding to the initial and terminal vertices of the edge $w$.

\begin{remark}
\label{rmk.topological.interior}
The difference $\llbracket w \rrbracket \setminus \llparenthesis w \rrparenthesis$ between a cell and its topological interior is either empty or it consists of one or both of the boundary vertices of $\llbracket w \rrbracket$.
More precisely, a boundary vertex $p$ of the cell $\llbracket w \rrbracket$ belongs to $\llparenthesis w \rrparenthesis$ if and only if every representative of $p$ has $w$ as a prefix (equivalently, if and only if $p$ has a unique representative).
Indeed, it is easy to see that small enough neighborhoods of $p$ are included in $\llbracket w \rrbracket$ if and only if every representative of $p$ has $w$ as a prefix.
\end{remark}

\subsubsection{Rational and Irrational Points of the Limit Space}

Another distinguishing feature of points is the (ir)rationality of representatives.
The application of this notion to limit spaces of edge replacement systems was initially developed together with Davide Perego for \cite{rationalgluing}, but we ended up not including it in the preprint for the sake of conciseness.
Irrational points of limit spaces had already made an appearance in another joint work with Davide Perego \cite{IG}.

We say that an element of an edge shift is \textbf{rational} if it can be written as $a \overline{b}$ for some finite words $a$ and $b$.
Otherwise, we say that it is \textbf{irrational}.
This translates to limit spaces.

\begin{definition}
\label{def.rational.irrational.points}
A point of the limit space is \textbf{rational} if each of its representatives is a rational element of the symbol space $\Omega_\mathcal{R}$.
A point is instead \textbf{irrational} if it has a unique representative which is an irrational element of $\Omega_\mathcal{R}$.
\end{definition}

Under this definition, there are points that are neither rational nor irrational, as they may have both rational and irrational representatives.
An example will be given shortly (\cref{fig.irrational.vertices}), together with a discussion on this topic (see \cref{cor.either.rational.or.irrational}).

\begin{remark}
Note that the (ir)rationality of a point is not a topological property, i.e., there may exist homeomorphisms of the space that map a rational point to an irrational point.
The simplest example is arguably given by the edge replacement system from \cref{fig.interval.replacement}, which will be considered later to realize Thompson's group $F$;
its limit space is the unit interval $[0,1]$ and rational points are precisely the elements of $\mathbb{Q} \cap [0,1]$.

Once rearrangements will be defined, since they act simply by simply changing prefixes of sequence in $\Omega_{\mathcal{R}}$, it will be clear that the sets of rational and irrational points are instead invariant under rearrangements.
\end{remark}

\subsubsection{Vertices and Rationality}

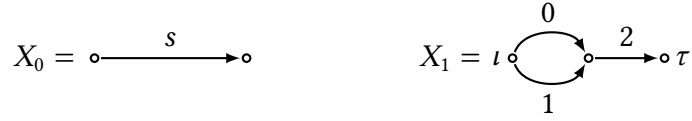
\begin{figure}
\centering
\begin{tikzpicture}
    \node at (-.667,0) {$X_0 =$};
    \node[vertex] (i) at (0,0) {};
    \node[vertex] (t) at (2,0) {};
    \draw[edge] (i) to node[above]{$s$} (t);
    \begin{scope}[xshift=5.5cm]
    \node at (-.825,0) {$X_1 =$};
    \node[vertex] (l) at (0,0) {};
    \draw (l) node[left]{$\iota$};
    \node[vertex] (c) at (1,0) {};
    \node[vertex] (r) at (2,0) {};
    \draw (r) node[right]{$\tau$};
    \draw[edge] (l) to[out=70,in=120] node[above]{$0$} (c);
    \draw[edge] (l) to[out=290, in=240] node[below]{$1$} (c);
    \draw[edge] (c) to node[above]{$2$} (r);
    \end{scope}
\end{tikzpicture}
\caption{An edge replacement system with irrational vertices.}
\label{fig.irrational.vertices}
\end{figure}

In general, rationality does not relate well to vertices.
For example, the edge replacement system depicted in \cref{fig.irrational.vertices} features vertices with irrational representatives, such as the leftmost vertex that is represented by any infinite sequence that belongs to $\{s\} \times \{0,1\}^\omega$, including all of the (uncountably many) irrational sequences.

In most cases that have been considered in the literature, as proved in the following proposition, all vertices are rational because they satisfy the following assumption.

\begin{definition}
\label{def.finite.branching}
We say that an edge replacement system has \textbf{finite branching} if the set of degrees of vertices among graph expansions is bounded.
\end{definition}

\medskip %layout
\begin{proposition}
Consider an edge replacement systems that has finite branching.
Then each of its vertices is rational.
Moreover, the number of representatives of each vertex equals its maximum degree among graph expansions.
\end{proposition}

\begin{proof}
Let $v$ be a vertex of the limit space.
Since the edge replacement system has finite branching, there is a bound on the degrees of $v$ among all graph expansions.
Consider the minimal graph expansion $E$ in which the maximal degree $d$ of $v$ is achieved, meaning that the graph $E$ features precisely $d$ edges that are incident on $v$, each represented by some finite word.
By definition of the gluing relation (\cref{def.glue}), every infinite sequence $\alpha$ representing the point $v$ must have one of these $d$ words as a prefix.

There is precisely one representative of $v$ for each edge $w \in \mathbb{L}_\mathcal{R}$ of $E$ that is incident on $v$.
Indeed, when an edge that is incident on $v$ is expanded, it is replaced by a subgraph that includes at least one edge that is incident on $v$, so there is at least one representative for each edge of $E$ that is incident on $v$.
Moreover, since $d$ is the maximum degree of $v$ among graph expansions, no sequence of edge expansions of $E$ can increase the degree of $v$, so each edge of $E$ that is incident on $v$ is the prefix of at most one representative.

Now, since edge expansions of $E$ do not change the degree of $v$, the expansion of an edge $w \in \mathbb{L}_\mathcal{R}$ of $E$ that is incident on $v$ produces exactly one new edge $wx \in \mathbb{L}_\mathcal{R}$ that is incident on $v$.
The edge $x$ of the replacement graph $X_{\mathrm{c}(w)}$ only depends on $\mathrm{c}(w)$ and on whether $w$ originates from or terminates at $v$.
Since there are finitely many such combinations, after a finite amount of edge expansions one will find some $w x_1 \dots x_k x_{k+1} \dots x_{k+p+1}$ such that $x_{k+1}$ and $x_{k+p+1}$ share the same color and either both originate from or both terminate at $v$.
This ``forces'' the subsequent edges $x_{k+p+i}$ to be the same as $X_{k+i}$.
Thus, if we let $x = w x_1 \dots x_k$ and $y = x_{k+1} \dots x_{k+p+1}$ we have that $v = \llbracket x \overline{y} \rrbracket$, so $v$ is rational.
\end{proof}

Since every point that has multiple representatives is a vertex (see \cref{rmk.glue.implies.vertex}), we have the following consequence.

\medskip %layout
\begin{corollary}
\label{cor.either.rational.or.irrational}
If an edge replacement system has finite branching, then every point is either rational or irrational.
\end{corollary}

Having finite branching is a natural assumption on an edge replacement system:
most edge replacement systems that have been considered in the literature so far all have initial and terminal vertices of replacement graphs that have total degree 1, which easily implies that they have finite branching.
Theorem 4.1 of \cite{BF19}, which provides a sufficient condition for finiteness properties (see \cref{sub.other.results}), requires the edge replacement system to have finite branching.
However, we will not need this assumption in the remainder of this dissertation.

\begin{remark}
\label{rmk.irrational.points.are.dense}
Whether the edge replacement system has finite branching or not, the set of irrational points of a limit space is dense.

Indeed, each cell has countably many points with more than one representative, as they must be vertices of some sequence of edge expansions of that cell, and it has countably many points with eventually periodic representative.
Since each cell is an uncountable set, we can conclude that it contains uncountably many irrational points.
In particular, because of \cref{rmk.density}, the set of irrational points is dense.
\end{remark}

\subsection{More Topological Properties of the Limit Space}

We finally highlight the following fact about cells and balls of limit spaces, which is going to be especially useful throughout \cref{cha.IG}.
Recall that limit spaces are metric spaces (\cref{thm.limit.space}) and note that there is no concern over the choice of metric (see \cite[Notes 1.30 (1)]{BF19}).

\medskip %layout
\begin{lemma}[Lemma 2.2 \cite{IG}]
\label{lem.balls.and.cells}
Each ball of a limit space contains a cell and each cell contains a ball.
\end{lemma}

\begin{proof}
Let $B$ be a ball with radius $r > 0$ centered in $p \in X$.
Note that $B$ must contain some vertex (for example because the representatives of vertices form a dense subset of the symbol space, which is easy to see), so we can assume that the center $p$ of $B$ is itself a vertex up to passing to a subset.

Consider the $n$-th open star of $p$, denoted by $St_n(p)$, which is the union of $\{p\}$ with the interior of the cells corresponding to the edges of the $n$-th full expansion graph (i.e., words of length $n$, recall \cref{def.full.expansion}) that are incident on $p$.
As is proved in \cite[Proposition 1.29]{BF19}, the sequence of diameters of $St_n(p)$ goes to zero as $n$ approaches infinity, so there exists some $N \in \mathbb{N} = \{0,1,\dots\}$ such that $St_N(p) \subseteq B$.
Since each star clearly contains a cell, $B$ must contain a cell too.

The converse is straightforward, since the interior of each cell is open and non-empty.
\end{proof}

Along the lines of this proof, we have actually showed that open stars of vertices form neighborhood bases.
More generally, given a point $p$ of the limit space, one can define the $n$-th open star $St_n(p)$ at $p$ as the union of $\{p\}$ with all the cell interiors $\llparenthesis w \rrparenthesis$ that contain $p$ among the words $w$ of length $n$.
When $p$ is not a vertex, each star is the interior of a single cell.
As a consequence, we have the following fact.

\medskip %layout
\begin{corollary}
\label{cor.basis.of.stars}
For each point $p$ of the limit space, the collection of open stars at $p$ forms a neighborhood basis at $p$.
\end{corollary}

Recall that neighborhood bases can be used to define the topology \cite[Theorem 4.5]{Willard}, so \cref{cor.basis.of.stars} shows that cells really give a complete characterization of the topology of the limit space.

\section{The Forest of Graph Expansions}
\label{sec.forest.of.expansions}

Graph expansions can be represented by certain forests that naturally encode their nested structure.
Essentially, this ``translates'' the data codified by graph expansions into certain finite forests that are colored and labeled on the edges.
This way of representing graph expansions will later lead to forest pair diagrams that represent rearrangements (\cref{sec.forest.pair.diagrams}).

\subsection{Forest Expansions}
\label{sub.forest.expansions}

We are going to fix an arbitrary ordering of the edges of each graph of $\mathcal{R}$.
To encode this data we use the simple notion of rotation systems of graphs, which was defined in \cref{def.rotation.system}.

Recall the definition of the language and alphabet of an edge replacement system (\cref{prop.language.alphabet}):
we think of $e \in \mathbb{A}_\mathcal{R}$ both as letters and as edges of the base or replacement graphs and of $w \in \mathbb{L}_\mathcal{R}$ both as words and as edges of graph expansions.

\phantomsection\label{txt.forest.of.expansions}
Let $\mathcal{R}$ be an edge replacement system and, for each of its graphs, fix an arbitrary ordering of its edges.
We call the \textbf{forest of graph expansions} $\mathbb{F}_\mathcal{R}$ of the edge replacement graph $\mathcal{R}$ the graph defined in the following way, where we refer to its vertices as \textbf{nodes} and its edges as \textbf{branches} to avoid confusion with the other graphs deriving from edge replacement systems.
\begin{itemize}
    \item The nodes are the following:
    there is a root $\varepsilon_s$ for each edge $s$ of the base graph and there is an interior node for each element of $\mathbb{L}_\mathcal{R}$.
    \item There is a branch from $w$ to $we$ for every $w \in \mathbb{L}_\mathcal{R}$ and $e \in \mathbb{A}_\mathcal{R}$ such that $we \in \mathbb{L}_\mathcal{R}$.
    Moreover, for each edge $s$ of the base graph there is a branch from $\varepsilon_s$ to $s$.
    \item A branch from $w$ to $we$ is colored by the same color of $e$ and labeled by the triple $(\iota(we),\tau(we),z)$, where $z \in \mathbb{Z}$ distinguishes between parallel edges (i.e., two distinct edges $e_1$ and $e_2$ such that $\iota(e_1) = \iota(e_2)$ and $\tau(e_1) = \tau(e_2)$).
\end{itemize}
Moreover, we equip $\mathbb{F}_\mathcal{R}$ with an ordering of the roots that corresponds to the chosen ordering of the edges of the base graph and with a rotation system where the sole incoming branch at each interior node $w$ is followed by the branches terminating at each $we$ according to the chosen ordering of the edges $e$ of the replacement graph.

\begin{definition}
\label{def.forest.expansion}
A \textbf{forest expansion} of an edge replacement system $\mathcal{R}$ is a complete subforest (\cref{def.complete.subforest}) of the forest of graph expansions $\mathbb{F}_\mathcal{R}$.
\end{definition}

\begin{remark}
\label{rmk.renaming}
When we refer to vertices, we use the second convention described in \cref{rmk.vertices.as.words}:
it does not matter which symbol we use for vertices, as long as each vertex has its own distinct symbol.

Moreover, we consider two forest expansions to be the same if they only differ by a renaming of symbols (see \cref{fig.airplane.forest.expansions,fig.airplane.forest.renaming} for an example).
More precisely, this renaming of symbols defines an equivalence relation and what we are really interested in are the equivalence classes.
\end{remark}

In most cases, the graph expansions of the edge replacement system do not have parallel edges.
In this case, we can safely omit the third component of the labeling, as we will do from now on.

\cref{fig.airplane.forest.expansions} depicts three forest expansion for the airplane edge replacement system from \cref{fig.airplane.replacement}.
Note that these are trees simply because the base graph consists of a sole edge (the basilica edge replacement system from the upcoming \cref{fig.basilica.replacement} would instead feature two roots, one for each loop of the base graph).
As done in \cref{fig.airplane.forest.expansions}, it is convenient to write each label of a branch vertically on its left and to color the label instead of the branch itself, for the sake of clarity.
We also omit the orientation of the branches, which are instead represented simply by their downward direction.

\begin{figure}
    \centering
    \begin{subfigure}[b]{\textwidth}
    \centering
    \begin{tikzpicture}[font=\small,scale=1.2]
        \node[node] (root) at (0,0) {};
        \node[node] (s) at (0,-1) {};
        \draw (root) to node[blue,left,align=center]{l\\r} (s);
    \end{tikzpicture}
    \hspace{1.75cm}
    \begin{tikzpicture}
        \node at (0,1) {$X_0$};
        \node[vertex] (l) at (-.75,0) {};
        \node[vertex] (r) at (.75,0) {};
        \draw[edge,blue] (l) node[black,above]{l} to (r) node[black,above]{r};
    \end{tikzpicture}
    \caption{The forest expansion that corresponds to the base graph $X_0$.}
    \label{fig.airplane.base.forest}
    \end{subfigure}
    \\
    \vspace{.5cm}
    \begin{subfigure}[b]{\textwidth}
    \centering
    \begin{tikzpicture}[font=\small,scale=1.2]
        \node[node] (root) at (0,0) {};
        \node[node] (s) at (0,-1) {};
        \node[node] (l) at (-1.5,-2) {};
        \node[node] (b) at (-.5,-2) {};
        \node[node] (r) at (.5,-2) {};
        \node[node] (t) at (1.5,-2) {};
        \draw (root) to (s) node[blue,above left,align=center]{l\\r};
        \draw (s) to[out=180,in=90,looseness=1.25] (l) node[blue,above left,align=center]{x\\l};
        \draw (s) to[out=180,in=90,looseness=1.25] (b) node[red,above left,align=center]{x\\y};
        \draw (s) to[out=0,in=90,looseness=1.25] (r) node[blue,above left,align=center]{y\\r};
        \draw (s) to[out=0,in=90,looseness=1.25] (t) node[red,above left,xshift=-.2cm,align=center]{y\\x};
    \end{tikzpicture}
    \hspace{1.5cm}
    \begin{tikzpicture}
        \node at (0,1.5) {$X_0 \triangleleft s$};
        \draw[edge,red,domain=5:175] plot ({.667*cos(\x)}, {.667*sin(\x)});
        \draw (90:.667);
        \draw[edge,red,domain=185:355] plot ({.667*cos(\x)}, {.667*sin(\x)});
        \draw (270:.667);
        \node[vertex] (l) at (-2.333,0) {};
        \node[vertex] (cl) at (-.667,0) {};
        \node[vertex] (cr) at (.667,0) {};
        \node[vertex] (r) at (2.333,0) {};
        \draw[edge,blue] (cl) node[black,above left]{x} to (l) node[black,above]{l};
        \draw[edge,blue] (cr) node[black,above right]{y} to (r) node[black,above]{r};
    \end{tikzpicture}
    \caption{A forest expansion and the corresponding graph expansion $X_0 \triangleleft s$.}
    \label{fig.airplane.forest.expansions.b}
    \end{subfigure}
    \\
    \vspace{.5cm}
    \begin{subfigure}[b]{\textwidth}
    \centering
    \begin{tikzpicture}[font=\small,scale=1.2]
        \node[node] (root) at (0,0) {};
        \node[node] (s) at (0,-1) {};
        \node[node] (l) at (-2.5,-3) {};
        \node[node] (b) at (-1.5,-3) {};
        \node[node] (r) at (-.5,-3) {};
        \node[node] (t) at (1.5,-2) {};
            \node[node] (t1) at (.5,-3) {};
            \node[node] (t2) at (1.5,-3) {};
            \node[node] (t3) at (2.5,-3) {};
        \draw (root) to (s) node[blue,above left,align=center]{l\\r};
        \draw (s) to[out=180,in=90,looseness=1.1] ($(l)+(0,1)$) to (l) node[blue,above left,align=center]{x\\l};
        \draw (s) to[out=180,in=90,looseness=1.25] ($(b)+(0,1)$) to (b) node[red,above left,align=center]{x\\y};
        \draw (s) to[out=180,in=90,looseness=1.25] ($(r)+(0,1)$) to (r) node[blue,above left,align=center]{y\\r};
        \draw (s) to[out=0,in=90,looseness=1.25] (t) node[red,above left,xshift=-.2cm,align=center]{y\\x};
            \draw (t) to[out=180,in=90,looseness=1.25] (t1) node[red,above left,align=center]{y\\z};
            \draw (t) to (t2) node[blue,above left,xshift=-.2cm,align=center]{z\\v};
            \draw (t) to[out=0,in=90,looseness=1.25] (t3) node[red,above left,xshift=-.2cm,align=center]{z\\x};
    \end{tikzpicture}
    \hspace{1.25cm}
    \begin{tikzpicture}
        \node at (0,2.5) {$X_0 \triangleleft s \triangleleft s b_2$};
        \draw[edge,red,domain=95:175] plot ({.667*cos(\x)}, {.667*sin(\x)});
        \draw (135:.667);
        \draw[edge,red,domain=185:355] plot ({.667*cos(\x)}, {.667*sin(\x)});
        \draw (270:.667);
        \draw[edge,red,domain=5:85] plot ({.667*cos(\x)}, {.667*sin(\x)});
        \draw (45:.667);
        \node[vertex] (l) at (-2.333,0) {};
        \node[vertex] (cl) at (-.667,0) {};
        \node[vertex] (c) at (0,.667) {};
        \node[vertex] (ct) at (0,1.667) {};
        \node[vertex] (cr) at (.667,0) {};
        \node[vertex] (r) at (2.333,0) {};
        \draw[edge,blue] (cl) node[black,above left]{x} to (l) node[black,above]{l};
        \draw[edge,blue] (cr) node[black,above right]{y} to (r) node[black,above]{r};
        \draw[edge,blue] (c) node[black,above left]{z} to (ct) node[black,left]{v};
    \end{tikzpicture}
    \caption{A forest expansion and the corresponding graph expansion $X_0 \triangleleft s \triangleleft s b_2$.}
    \end{subfigure}
    \caption{Three forest expansions of the airplane edge replacement system from \cref{fig.airplane.replacement} and the corresponding graph expansions from \cref{fig.airplane.expansions}.}
    \label{fig.airplane.forest.expansions}
\end{figure}
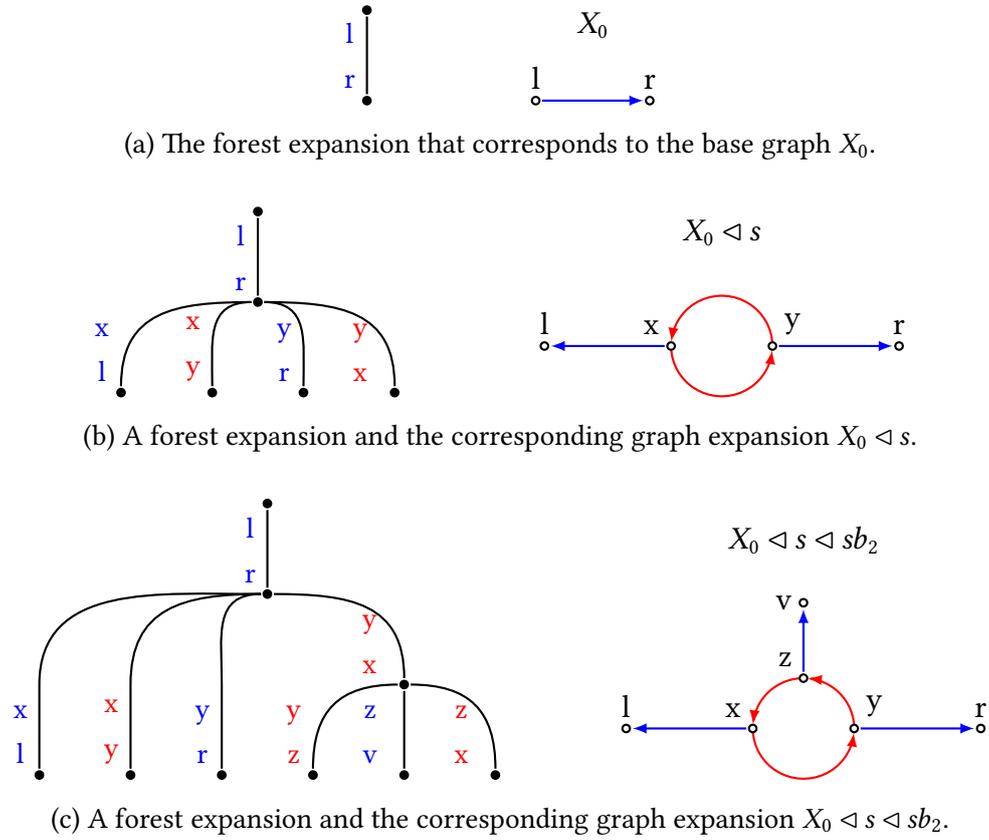

\begin{figure}
    \centering
    \begin{tikzpicture}[font=\small,scale=1.2]
        \node[node] (root) at (0,0) {};
        \node[node] (s) at (0,-1) {};
        \node[node] (l) at (-1.5,-2) {};
        \node[node] (b) at (-.5,-2) {};
        \node[node] (r) at (.5,-2) {};
        \node[node] (t) at (1.5,-2) {};
        \draw (root) to (s) node[blue,above left,align=center]{a\\b};
        \draw (s) to[out=180,in=90,looseness=1.25] (l) node[blue,above left,align=center]{l\\a};
        \draw (s) to[out=180,in=90,looseness=1.25] (b) node[red,above left,align=center]{l\\r};
        \draw (s) to[out=0,in=90,looseness=1.25] (r) node[blue,above left,align=center]{r\\b};
        \draw (s) to[out=0,in=90,looseness=1.25] (t) node[red,above left,xshift=-.2cm,align=center]{r\\l};
    \end{tikzpicture}
    \caption{A forest expansion that is equivalent to that from \cref{fig.airplane.forest.expansions.b}.}
    \label{fig.airplane.forest.renaming}
\end{figure}
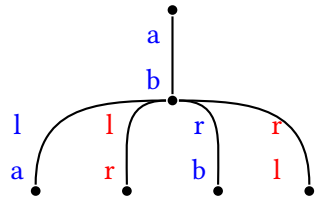

\subsection{Correspondence with Graph Expansions}
\label{sub.forests.and.graphs}

As suggested by \cref{fig.airplane.forest.expansions}, there is a natural bijection between forest expansions and graph expansions which is built inductively as follows.

The \textbf{base forest} $F_0$ of $\mathcal{R}$ is the smallest forest expansion, which only features the roots and, for each root, its unique child.
Since our graphs do not feature isolated vertices (\cref{ass.isolated.vertices}), the labels of the branches that terminate at the sinks of $F_0$ read precisely the base graph $\Gamma_0$ of $\mathcal{R}$, in the sense that each labeled $(v,w,z)$ corresponds to an edge from $v$ to $w$.

Inductively, suppose that we have a forest expansion $F$ that corresponds to a graph expansion $\Gamma$.
Let $e$ be a $c$-colored edge of $\Gamma$ and let $v$ be the corresponding sink of $F$.
We can expand $F$ to $F \triangleleft e$ by attaching below $v$ a child for each edge of the replacement graph $\Gamma_c$, in their previously fixed order.
Doing so, we are adding a branch for each edge $e'$ of $\Gamma_c$;
according to this correspondence, we label a branch by $(\iota(e'),\tau(e'),z)$,
where $z$ distinguishes parallel edges (i.e., if the bottom branches already have labels $(\iota(e'),\tau(e'),z)$, we use a distinct $z$).
Then $F \triangleleft e$ corresponds to $\Gamma \triangleleft e$ in the sense that the labels of the bottom branches of $F \triangleleft e$ read the graph $\Gamma \triangleleft e$.

\begin{figure}
    \centering
    \begin{subfigure}[b]{\textwidth}
    \centering
    \begin{tikzpicture}[font=\small,scale=1.2]
        \node[node] (root) at (0,0) {};
        \node[node] (s) at (0,-1) {};
        \node[node] (l) at (-1.5,-2) {};
        \node[node] (b) at (-.5,-2) {};
        \node[node] (r) at (.5,-2) {};
        \node[node] (t) at (1.5,-2) {};
        \draw (root) to (s) node[blue,above left,align=center]{$\iota$\\$\tau$};
        \draw (s) to[out=180,in=90,looseness=1.25] (l) node[blue,above left,align=center]{x\\$\iota$};
        \draw (s) to[out=180,in=90,looseness=1.25] (b) node[red,above left,align=center]{x\\y};
        \draw (s) to[out=0,in=90,looseness=1.25] (r) node[blue,above left,align=center]{y\\$\tau$};
        \draw (s) to[out=0,in=90,looseness=1.25] (t) node[red,above left,xshift=-.2cm,align=center]{y\\x};
    \end{tikzpicture}
    \hspace{1.5cm}
    \begin{tikzpicture}
        \node at (0,1.5) {$X_{\text{\textcolor{blue}{b}}}$};
        \draw[edge,red,domain=5:175] plot ({.667*cos(\x)}, {.667*sin(\x)});
        \draw (90:.667);
        \draw[edge,red,domain=185:355] plot ({.667*cos(\x)}, {.667*sin(\x)});
        \draw (270:.667);
        \node[vertex] (l) at (-2.333,0) {};
        \node[vertex] (cl) at (-.667,0) {};
        \node[vertex] (cr) at (.667,0) {};
        \node[vertex] (r) at (2.333,0) {};
        \draw[edge,blue] (cl) node[black,above left]{x} to (l) node[black,above]{$\iota$};
        \draw[edge,blue] (cr) node[black,above right]{y} to (r) node[black,above]{$\tau$};
    \end{tikzpicture}
    \caption{The blue replacement tree $T_{\text{\textcolor{blue}{b}}}$ beside the blue replacement graph.}
    \end{subfigure}
    \\
    \vspace{.5cm}
    \begin{subfigure}[b]{\textwidth}
    \centering
    \begin{tikzpicture}[font=\small,scale=1.2]
        \node[node] (root) at (0,0) {};
        \node[node] (s) at (0,-1) {};
        \node[node] (l) at (-1,-2) {};
        \node[node] (c) at (0,-2) {};
        \node[node] (r) at (1,-2) {};
        \draw (root) to (s) node[red,above left,align=center]{$\iota$\\$\tau$};
        \draw (s) to[out=180,in=90,looseness=1.25] (l) node[red,above left,align=center]{$\iota$\\z};
        \draw (s) to (c) node[blue,above left,align=center]{z\\v};
        \draw (s) to[out=0,in=90,looseness=1.25] (r) node[red,above left,xshift=-.2cm,align=center]{z\\$\tau$};
    \end{tikzpicture}
    \hspace{1.25cm}
    \begin{tikzpicture}
        \node at (0,2) {$X_{\text{\textcolor{red}{r}}}$};
        \node[vertex] (l) at (-1.667,0) {}; \draw (-1.667,0) node[above]{$\iota$};
        \node[vertex] (r) at (1.667,0) {}; \draw (1.667,0) node[above]{$\tau$};
        \node[vertex] (c) at (0,0) {};
        \node[vertex] (ct) at (0,1.2) {};
        \draw[edge,red] (l) node[black,above]{$\iota$} to (c);
        \draw[edge,red] (c) to (r) node[black,above]{$\tau$};
        \draw[edge,blue] (c) node[black,above left]{z} to (ct) node[black,left]{v};
    \end{tikzpicture}
    \caption{The red replacement tree $T_{\text{\textcolor{red}{r}}}$ beside the red replacement graph.}
    \end{subfigure}
    \caption{The replacement trees of the airplane edge replacement system.}
    \label{fig.airplane.replacement.trees}
\end{figure}
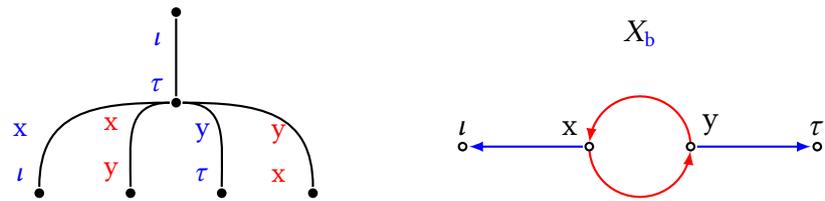
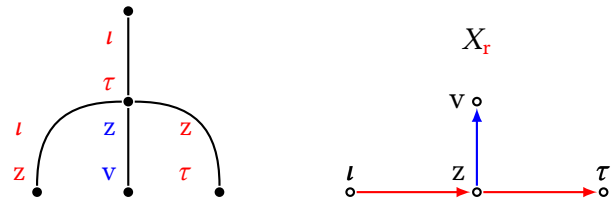

Such expansions are performed locally by replacing a $c$-colored bottom strand with a certain tree, which we call \textbf{replacement tree} $T_c$.
This tree depends on the color $c$ and its labels are adjusted according to the label of the branch being replaced:
we use the symbols $\iota$ and $\tau$ to denote the initial and terminal vertices of the replacement graph $\Gamma_c$, which are to be replaced by the initial and terminal vertices of the edge $e$ being expanded, and every other symbol is temporary and should be replaced with new ones that have not been used before (in accordance with \cref{rmk.renaming}).
Ultimately, the base and the replacement graphs are entirely codified by the base forest and replacement trees.
For example, the base forest of the airplane edge replacement system is depicted in \cref{fig.airplane.base.forest} and the blue and red replacement trees are portrayed in \cref{fig.airplane.replacement.trees}.

Since the graph expansion that corresponds to a forest expansion $F$ can be read on its bottom branches, of which there is precisely one for each leaf of $F$, it is natural to give the following definition.

\begin{definition}
\label{def.leaf.graph}
Given a forest expansion $F$, the graph represented by the labeling of its bottom branches is called the \textbf{leaf graph} of $F$.
\end{definition}

For example, the three graphs on the right in \cref{fig.airplane.forest.expansions} are the leaf graphs of the forest expansions on the right.

%%%%%%%%%%%%%%%%%%%%%%%%%

\chapter{Rearrangement Groups}
\label{cha.rearrangements}

In this chapter we define rearrangements, which are essentially the homeomorphisms of the limit spaces of edge replacement systems that descend from prefix-exchange transformations of the underlying edge shift.

Throughout this chapter, let $\mathcal{R}$ be an edge replacement system with a properly defined limit space $X$ (not necessarily expanding, see \cref{rmk.generalized.limit.spaces}).
Even when the limit space is not defined, rearrangements can still be formally defined as certain gluing-preserving prefix-exchange homeomorphisms of the edge shift $\Omega_\mathcal{R}$;
they would not be regarded as homeomorphisms of a limit space, so we will not be interested in that case.

\section{Rearrangements of Limit Spaces}

\begin{definition}
\label{def.canonical.homeomorphism}
Consider two cells $\llbracket a \rrbracket$ and $\llbracket b \rrbracket$ of the same type (\cref{def.cell.type}).
The \textbf{canonical homeomorphism} between $\llbracket a \rrbracket$ and $\llbracket b \rrbracket$ is the prefix-exchange homeomorphism
\[ \llbracket a \rrbracket \to \llbracket b \rrbracket, \, \llbracket a \alpha \rrbracket \mapsto \llbracket b \alpha \rrbracket. \]
\end{definition}

\begin{definition}
\label{def.rearrangement}
A \textbf{rearrangement} of the limit space is a homeomorphism $\phi$ of the limit space such that there exists a cellular partition (\cref{def.cellular.partition}) on whose cells the restrictions of $\phi$ are canonical homeomorphisms.
\end{definition}

Rearrangements form a group under composition (\cite[Proposition 1.16]{BF19}), so we write $G_\mathcal{R}$ for the group of rearrangements associated to the edge replacement system $\mathcal{R}$.
Explicit ways to compute compositions and inverse of rearrangements will be described with the use of graph pair diagrams (\cref{sec.graph.pair.diagrams}), forest pair diagrams (\cref{sec.forest.pair.diagrams}) and strand diagrams (\cref{sec.SDs}).

\begin{remark}
\label{rmk.rearrangement.action.symbol.space}
Rearrangements come equipped with a natural prefix-exchange action on the symbol space $\Omega_\mathcal{R}$:
if $\phi$ is a rearrangement that maps a cell $\llbracket w \rrbracket$ canonically to a cell $\llbracket v \rrbracket$, this action maps $w \alpha$ to $v \alpha$.
With a slight abuse of notation, thus, we will write $\phi(w) = v$ in this case.

This defines an action by homeomorphisms on the symbol space $\Omega_\mathcal{R}$ and, when the gluing relation is trivial (i.e., the limit space is the edge shift $\Omega_\mathcal{R}$ itself), this is precisely the same action of the rearrangement group.
This special case will be discussed in \cref{sub.topological.full.groups}.
\end{remark}

\section{Graph Pair Diagrams}
\label{sec.graph.pair.diagrams}

In practice, it is very convenient to represent rearrangements by graph isomorphisms between graph expansions, as described in this subsection.

\begin{definition}
\label{def.graph.pair.diagrams}
A \textbf{graph pair diagram} is a graph isomorphism between two graph expansions of the edge replacement system.
If $D$ and $R$ are the graph expansions and $\phi$ is the graph isomorphism, we denote the graph pair diagram by the triple $(D, \phi, R)$.
\end{definition}

\begin{figure}
\centering
\begin{minipage}{.45\textwidth}
\centering
\includegraphics[width=\textwidth]{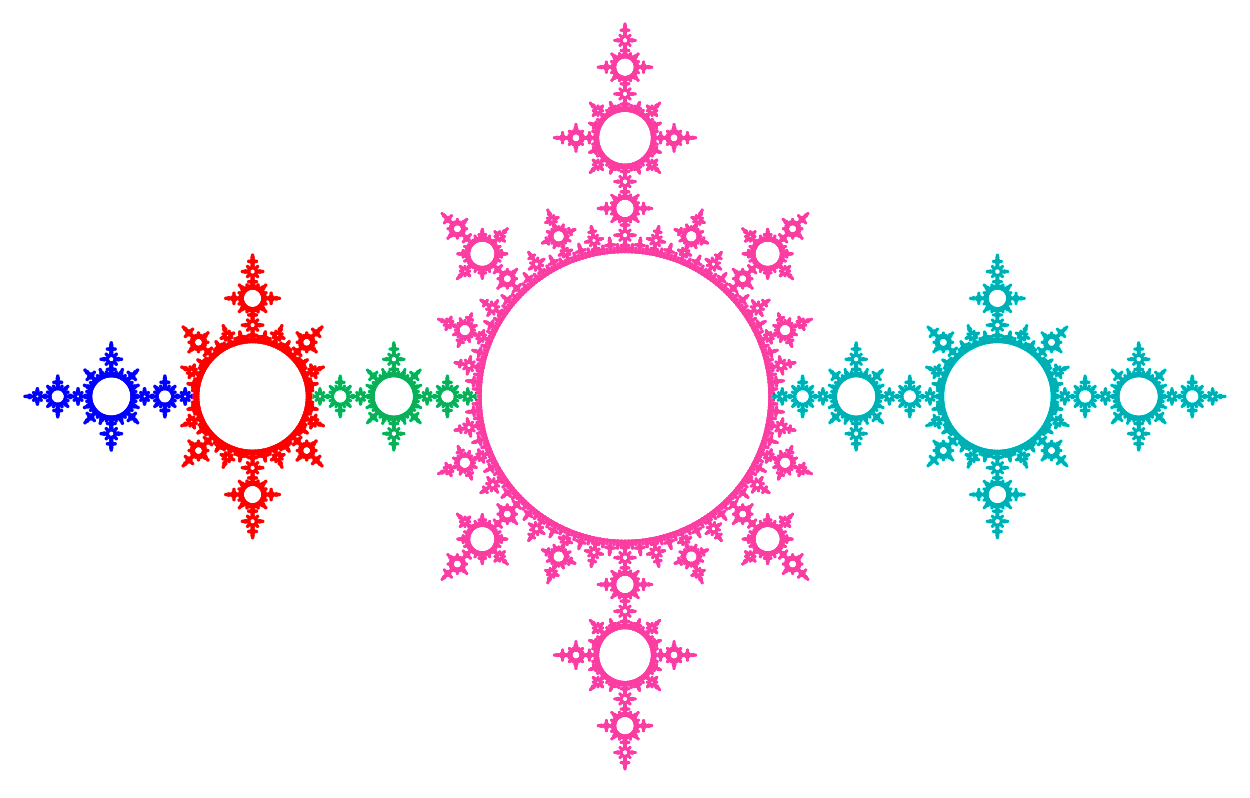}
\end{minipage}%
\begin{minipage}{.1\textwidth}
\centering
\begin{tikzpicture}[scale=1.25]
    \draw[-to] (-1/8,0) -- node[above]{$\alpha$} (1/8,0);
\end{tikzpicture}
\end{minipage}%
\begin{minipage}{.45\textwidth}
\centering
\includegraphics[width=\textwidth]{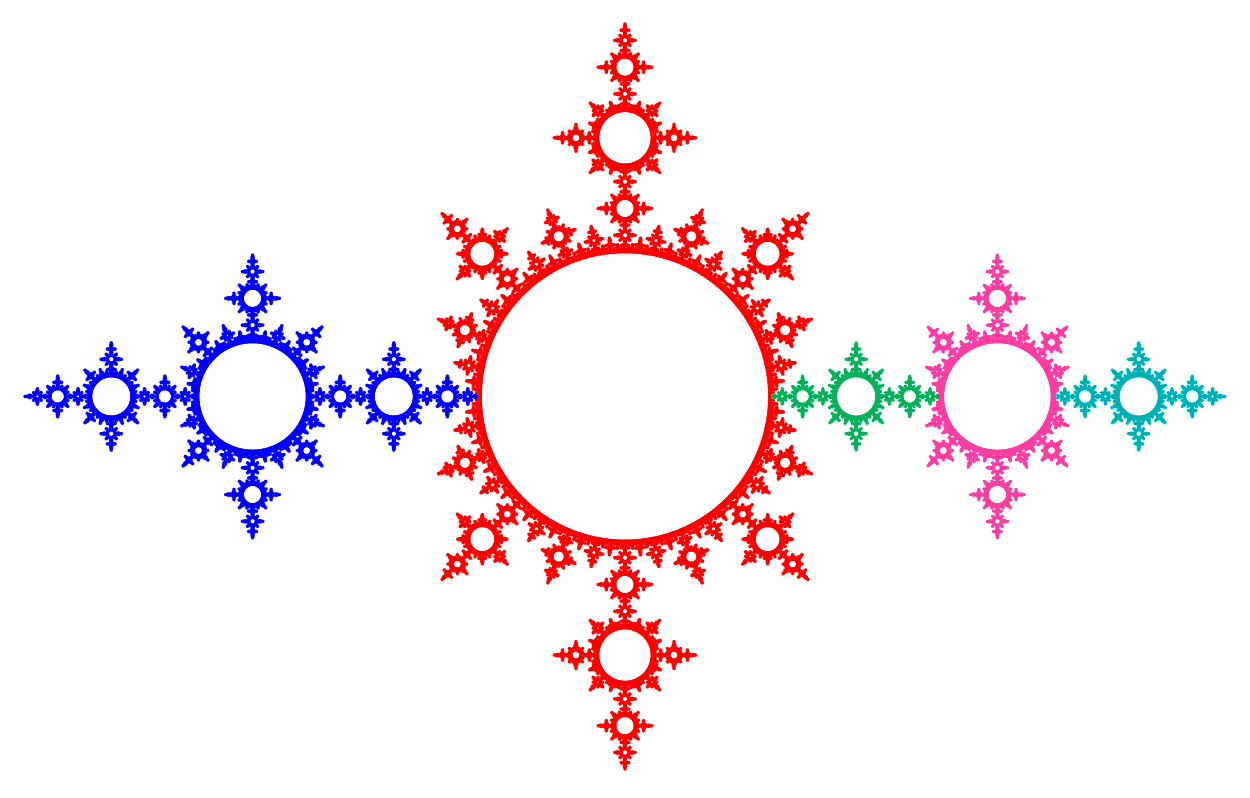}
\end{minipage}
\\
\vspace*{10pt}
\begin{tikzpicture}[scale=1.25]
    \begin{scope}[xshift=-3cm]
    \draw[Magenta,edge,domain=5:175] plot ({.667*cos(\x)}, {.667*sin(\x)});
    \draw[Magenta,edge,domain=185:355] plot ({.667*cos(\x)}, {.667*sin(\x)});
    \draw[red,domain=10:170] plot ({.278*cos(\x)-1.5}, {.278*sin(\x)});
    \draw[red,domain=190:350] plot ({.278*cos(\x)-1.5}, {.278*sin(\x)});
    \node[vertex] (l) at (-2.333,0) {};
    \node[vertex] (l1) at (-1.778,0) {};
    \node[vertex] (l2) at (-1.222,0) {};
    \node[vertex] (cl) at (-.667,0) {};
    \node[vertex] (cr) at (.667,0) {};
    \node[vertex] (r) at (2.333,0) {};
    \draw[edge,blue] (l1) to (l);
    \draw[Green] (l2) to (cl);
    \draw[Green,fill=white] (-.9445,0) circle (.0926);
    \draw[edge,TealBlue] (cr) to (r);
    \end{scope}
    \draw[-to] (-1/8,0) to node[above]{$\alpha$} (1/8,0);
    \begin{scope}[xshift=3cm]
    \draw[Magenta,domain=10:170] plot ({.278*cos(\x)+1.5}, {.278*sin(\x)});
    \draw[Magenta,domain=190:350] plot ({.278*cos(\x)+1.5}, {.278*sin(\x)});
    \draw[red,edge,domain=5:175] plot ({.667*cos(\x)}, {.667*sin(\x)});
    \draw[red,edge,domain=185:355] plot ({.667*cos(\x)}, {.667*sin(\x)});
    \node[vertex] (l) at (-2.333,0) {};
    \node[vertex] (cl) at (-.667,0) {};
    \node[vertex] (cr) at (.667,0) {};
    \node[vertex] (r1) at (1.222,0) {};
    \node[vertex] (r2) at (1.778,0) {};
    \node[vertex] (r) at (2.333,0) {};
    \draw[edge,blue] (cl) to (l);
    \draw[Green] (r1) to (cr);
    \draw[Green,fill=white] (.9445,0) circle (.0926);
    \draw[edge,TealBlue] (r2) to (r);
    \end{scope}
\end{tikzpicture}
\caption{A rearrangement of the airplane limit space, along with a graph pair diagram that represents it.}
\label{fig.airplane.rearrangement}
\end{figure}
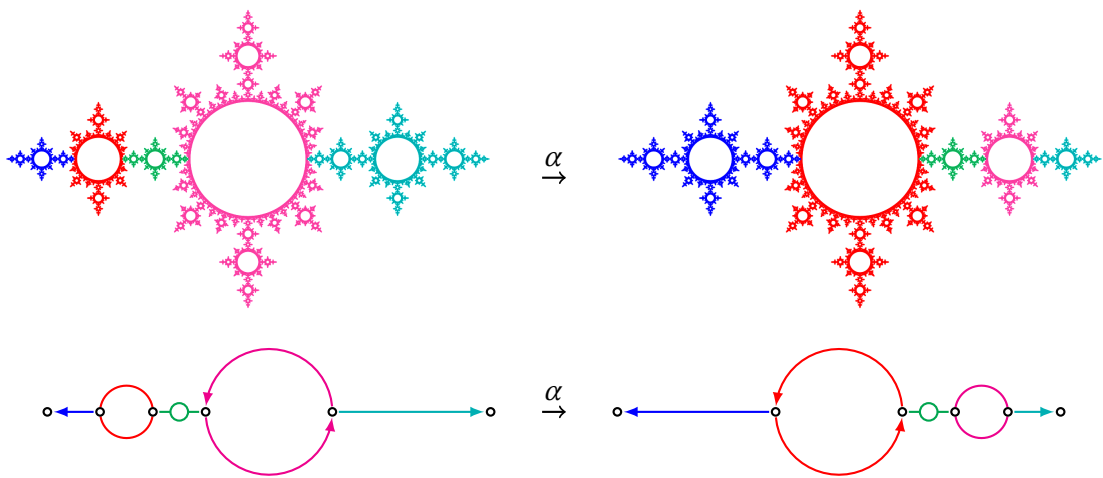

An example for the airplane edge replacement system is depicted in \cref{fig.airplane.rearrangement}, where colors are used to show where the rearrangement maps canonically each cell and where the graph isomorphisms maps each edge.
In \cref{fig.airplane.rearrangement} and in many other figures from here on, certain vertices and the orientation of certain edges is omitted for the sake of clarity.

Recall that we are working under \cref{ass.isolated.vertices}, so the base and replacement graphs of the edge replacement system $\mathcal{R}$ do not have isolated vertices.
In this setting, graph isomorphisms are determined by their bijection of edges.

With this observation in mind, one can see that each graph pair diagram defines a rearrangement.
Indeed, the graph isomorphism $\phi$, thought of as a bijection between the edges of $D$ and those of $R$, induces a map
\[ \llbracket w \alpha \rrbracket \mapsto \llbracket \phi(w) \alpha \rrbracket, \, \forall w \text{ edge of } D, \forall \alpha \text{ such that } w \alpha \in \Omega_\mathcal{R}, \]
which is a rearrangement.

Conversely, below we show how each rearrangement can be represented by a graph pair diagram.
Recall that cellular partition correspond one-to-one to graph expansions and note that, if $\mathbb{P}$ is a cellular partition on whose cells $\phi$ restricts to canonical homeomorphisms, then $\phi(\mathbb{P}) = \{ \phi(C) \mid C \in \mathbb{P} \}$ is also a cellular partition.
Given a rearrangement $\phi$, let $\{ \llbracket w_1 \rrbracket, \dots, \llbracket w_k \rrbracket \}$ be a cellular partition on whose cells $\phi$ restricts to canonical homeomorphisms.
Consider the graph expansions $D$ and $R$ whose edges are $w_1, \dots, w_k$ and $\phi(w_1), \dots, \phi(w_k)$, respectively.
The map $w_i \mapsto \phi(w_i)$ (as in \cref{rmk.rearrangement.action.symbol.space}) is a graph isomorphism, so $(D, \phi, R)$ is a graph pair diagram for $\phi$.

\subsection{Expansion of Graph Pair Diagrams}
\label{sub.reduced.graph.pair.diagrams}

Each rearrangement is represented by infinitely many graph pair diagrams.
Indeed, if a rearrangement $g$ is represented by $(D,\phi,R)$ and $e$ is an edge of $D$, then the graph pair diagram $(D,\phi,R)$ can be \textbf{expanded} to $(D \triangleleft e, \phi^*, R \triangleleft \phi(e))$, where $\phi^*$ agrees with $\phi$ on every edge of $D$ other than $e$ and maps the edges $ex$ resulting from the edge expansion of $e$ in $D$ to $\phi(e)x$.
If a graph pair diagram cannot be obtained from any other using such an expansion, we say that it is \textbf{reduced}.

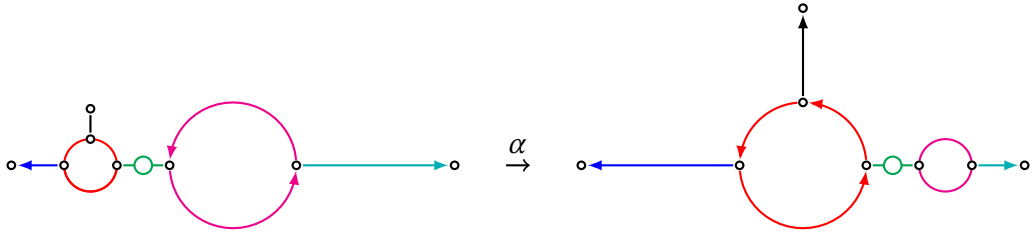
\begin{figure}
\centering
\begin{tikzpicture}[scale=1.25]
    \begin{scope}[xshift=-3cm]
    \draw[Magenta,edge,domain=5:175] plot ({.667*cos(\x)}, {.667*sin(\x)});
    \draw[Magenta,edge,domain=185:355] plot ({.667*cos(\x)}, {.667*sin(\x)});
    \draw[red,domain=10:80] plot ({.278*cos(\x)-1.5}, {.278*sin(\x)});
    \draw[red,domain=100:170] plot ({.278*cos(\x)-1.5}, {.278*sin(\x)});
    \draw[red,domain=190:350] plot ({.278*cos(\x)-1.5}, {.278*sin(\x)});
    \draw[red,domain=190:350] plot ({.278*cos(\x)-1.5}, {.278*sin(\x)});
    \node[vertex] (l) at (-2.333,0) {};
    \node[vertex] (l1) at (-1.778,0) {};
    \node[vertex] (lc) at (-1.5,.278) {};
    \node[vertex] (lct) at (-1.5,.6) {};
    \node[vertex] (l2) at (-1.222,0) {};
    \node[vertex] (cl) at (-.667,0) {};
    \node[vertex] (cr) at (.667,0) {};
    \node[vertex] (r) at (2.333,0) {};
    \draw[edge,blue] (l1) to (l);
    \draw[Green] (l2) to (cl);
    \draw[Green,fill=white] (-.9445,0) circle (.0926);
    \draw[edge,TealBlue] (cr) to (r);
    \draw[black] (lc) to (lct);
    \end{scope}
    \draw[-to] (-1/8,0) to node[above]{$\alpha$} (1/8,0);
    \begin{scope}[xshift=3cm]
    \draw[red,edge,domain=5:85] plot ({.667*cos(\x)}, {.667*sin(\x)});
    \draw[red,edge,domain=95:175] plot ({.667*cos(\x)}, {.667*sin(\x)});
    \draw[red,edge,domain=185:355] plot ({.667*cos(\x)}, {.667*sin(\x)});
    \draw[Magenta,domain=10:170] plot ({.278*cos(\x)+1.5}, {.278*sin(\x)});
    \draw[Magenta,domain=190:350] plot ({.278*cos(\x)+1.5}, {.278*sin(\x)});
    \node[vertex] (l) at (-2.333,0) {};
    \node[vertex] (cl) at (-.667,0) {};
    \node[vertex] (c) at (0,.667) {};
    \node[vertex] (ct) at (0,1.667) {};
    \node[vertex] (cr) at (.667,0) {};
    \node[vertex] (r1) at (1.222,0) {};
    \node[vertex] (r2) at (1.778,0) {};
    \node[vertex] (r) at (2.333,0) {};
    \draw[edge,blue] (cl) to (l);
    \draw[edge,black] (c) to (ct);
    \draw[Green] (r1) to (cr);
    \draw[Green,fill=white] (.9445,0) circle (.0926);
    \draw[edge,TealBlue] (r2) to (r);
    \end{scope}
\end{tikzpicture}
\caption{A graph pair diagram that is an expansion of the one from \cref{fig.airplane.rearrangement}.}
\label{fig.gpd.expansion}
\end{figure}

Expansions of graph pair diagrams are essentially computed by expanding the domain graph by an edge $e$, the range graph by its image $\phi(e)$ and by replacing $\phi$ with a graph isomorphism $\phi^*$ that is the same as $\phi$ except that the subgraph that results from the edge expansions of $D$ and $R$ are mapped ``canonically'' between each other.
For instance, \cref{fig.gpd.expansion} depicts an expansion of the graph pair diagram from \cref{fig.airplane.rearrangement}.
Reductions of graph pair diagrams (the converse of expansions) can be computed with the same ideas.
Since the graph pair diagram from \cref{fig.gpd.expansion} is an expansion of the one from \cref{fig.airplane.rearrangement}, it is not reduced;
instead, the graph pair diagram from \cref{fig.airplane.rearrangement} is reduced.

It is easy to see that the operation of expansion (or reduction) of graph pair diagrams does not change the rearrangement being represented.
In other words, if we are given a graph pair diagrams for a rearrangement $g$, we can perform any expansion (or reduction) on that diagram and we still end up with a diagram that represents the same rearrangement.
In this sense, rearrangements can be thought of as classes of equivalent graph pair diagrams.

A graph pair diagram that is not reduced is ``redundant'' in the same sense in which a fraction whose numerator and denominator share a nontrivial divisor is.
Fortunately, as rational numbers are represented by unique reduced fractions, rearrangements are represented by unique reduced graph pair diagrams.

\medskip %layout
\begin{lemma}[Proposition 1.19 \cite{BF19}]
Every rearrangement is represented by a unique reduced graph pair diagram.
\end{lemma}

The proof is very simple and can be found in \cite{BF19}.

\subsection{Composition of Graph Pair Diagrams}
\label{sub.composition.graph.pair.diagrams}

Computing compositions of graph pair diagrams is reminiscent of how one computes sums of fractions by first finding a common denominator.

Suppose that we are given two rearrangements $g$ and $h$ with graph pair diagrams $(D_g,\phi_g,R_g)$ and $(D_h,\phi_h,R_h)$, respectively.
By \cref{prop.common.finer.partition}, there exists a graph expansion that is a common expansion of $D_g$ and $R_h$.
We can thus always replace the two graph pair diagrams with their expansions in such a way that $D_g = R_h$.
Given such graph pair diagrams, the composition of rearrangements $g h$ is represented by the graph pair diagram $(D_h, \phi_h \circ \phi_g, R_g)$.

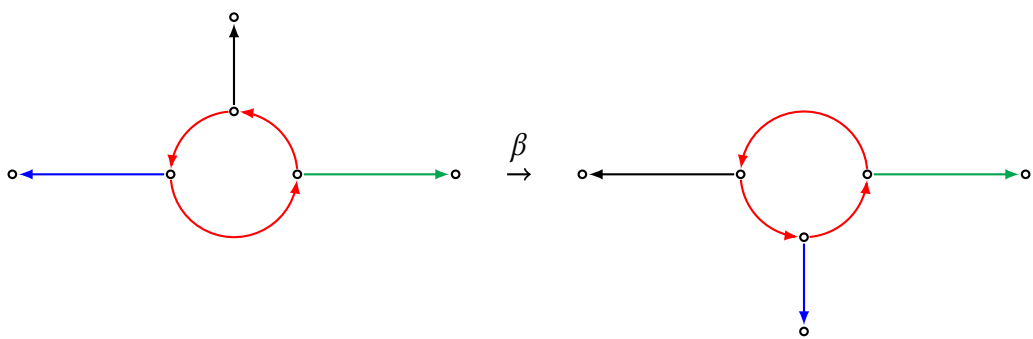
\begin{figure}
\centering
\begin{tikzpicture}[scale=1.25]
    \begin{scope}[xshift=-3cm]
    \draw[edge,red,domain=5:85] plot ({.667*cos(\x)}, {.667*sin(\x)});
    \draw[edge,red,domain=95:175] plot ({.667*cos(\x)}, {.667*sin(\x)});
    \draw[edge,red,domain=185:355] plot ({.667*cos(\x)}, {.667*sin(\x)});
    \node[vertex] (l) at (-2.333,0) {};
    \node[vertex] (cl) at (-.667,0) {};
    \node[vertex] (c) at (0,.667) {};
    \node[vertex] (ct) at (0,1.667) {};
    \node[vertex] (cr) at (.667,0) {};
    \node[vertex] (r) at (2.333,0) {};
    \draw[edge,blue] (cl) to (l);
    \draw[edge,Green] (cr) to (r);
    \draw[edge,black] (c) to (ct);
    \end{scope}
    \draw[-to] (-1/8,0) to node[above]{$\beta$} (1/8,0);
    \begin{scope}[xshift=3cm]
    \draw[edge,red,domain=185:265] plot ({.667*cos(\x)}, {.667*sin(\x)});
    \draw[edge,red,domain=275:355] plot ({.667*cos(\x)}, {.667*sin(\x)});
    \node[vertex] (l) at (-2.333,0) {};
    \node[vertex] (cl) at (-.667,0) {};
    \node[vertex] (c) at (0,-.667) {};
    \node[vertex] (cb) at (0,-1.667) {};
    \node[vertex] (cr) at (.667,0) {};
    \node[vertex] (r) at (2.333,0) {};
    \draw[edge,black] (cl) to (l);
    \draw[edge,Green] (cr) to (r);
    \draw[edge,red,domain=5:175] plot ({.667*cos(\x)}, {.667*sin(\x)});
    \draw[edge,blue] (c) to (cb);
    \end{scope}
\end{tikzpicture}
\caption{Another rearrangement of the airplane edge replacement system.}
\label{fig.gpd.beta}
\end{figure}

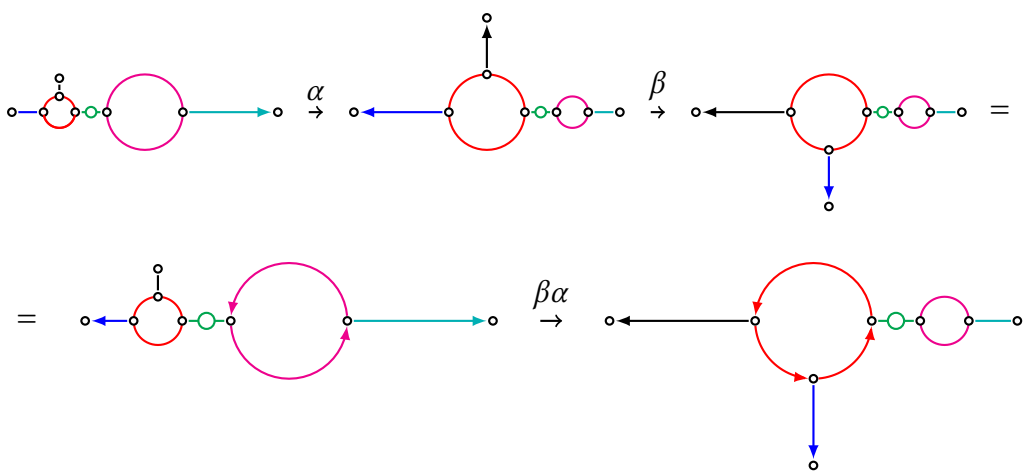
\begin{figure}
\centering
\begin{tikzpicture}[scale=.75]
    \begin{scope}[xshift=-3cm]
    \draw[Magenta,domain=5:175] plot ({.667*cos(\x)}, {.667*sin(\x)});
    \draw[Magenta,domain=185:355] plot ({.667*cos(\x)}, {.667*sin(\x)});
    \draw[red,domain=10:80] plot ({.278*cos(\x)-1.5}, {.278*sin(\x)});
    \draw[red,domain=100:170] plot ({.278*cos(\x)-1.5}, {.278*sin(\x)});
    \draw[red,domain=190:350] plot ({.278*cos(\x)-1.5}, {.278*sin(\x)});
    \node[vertex] (l) at (-2.333,0) {};
    \node[vertex] (l1) at (-1.778,0) {};
    \node[vertex] (lc) at (-1.5,.278) {};
    \node[vertex] (lct) at (-1.5,.6) {};
    \node[vertex] (l2) at (-1.222,0) {};
    \node[vertex] (cl) at (-.667,0) {};
    \node[vertex] (cr) at (.667,0) {};
    \node[vertex] (r) at (2.333,0) {};
    \draw[blue] (l1) to (l);
    \draw[Green] (l2) to (cl);
    \draw[Green,fill=white] (-.9445,0) circle (.0926);
    \draw[edge,TealBlue] (cr) to (r);
    \draw[red,domain=190:350] plot ({.278*cos(\x)-1.5}, {.278*sin(\x)});
    \draw[black] (lc) to (lct);
    \end{scope}
    \draw[-to] (-1/8,0) to node[above]{$\alpha$} (1/8,0);
    \begin{scope}[xshift=3cm]
    \draw[red,domain=5:85] plot ({.667*cos(\x)}, {.667*sin(\x)});
    \draw[red,domain=95:175] plot ({.667*cos(\x)}, {.667*sin(\x)});
    \draw[red,domain=185:355] plot ({.667*cos(\x)}, {.667*sin(\x)});
    \draw[Magenta,domain=10:170] plot ({.278*cos(\x)+1.5}, {.278*sin(\x)});
    \draw[Magenta,domain=190:350] plot ({.278*cos(\x)+1.5}, {.278*sin(\x)});
    \node[vertex] (l) at (-2.333,0) {};
    \node[vertex] (cl) at (-.667,0) {};
    \node[vertex] (c) at (0,.667) {};
    \node[vertex] (ct) at (0,1.667) {};
    \node[vertex] (cr) at (.667,0) {};
    \node[vertex] (r1) at (1.222,0) {};
    \node[vertex] (r2) at (1.778,0) {};
    \node[vertex] (r) at (2.333,0) {};
    \draw[edge,blue] (cl) to (l);
    \draw[edge,black] (c) to (ct);
    \draw[Green] (r1) to (cr);
    \draw[Green,fill=white] (.9445,0) circle (.0926);
    \draw[TealBlue] (r2) to (r);
    \end{scope}
    \begin{scope}[xshift=6cm]
    \draw[-to] (-1/8,0) to node[above]{$\beta$} (1/8,0);
    \end{scope}
    \begin{scope}[xshift=9cm]
    \draw[red,domain=5:175] plot ({.667*cos(\x)}, {.667*sin(\x)});
    \draw[red,domain=185:265] plot ({.667*cos(\x)}, {.667*sin(\x)});
    \draw[red,domain=275:355] plot ({.667*cos(\x)}, {.667*sin(\x)});
    \draw[Magenta,domain=10:170] plot ({.278*cos(\x)+1.5}, {.278*sin(\x)});
    \draw[Magenta,domain=190:350] plot ({.278*cos(\x)+1.5}, {.278*sin(\x)});
    \node[vertex] (l) at (-2.333,0) {};
    \node[vertex] (cl) at (-.667,0) {};
    \node[vertex] (c) at (0,-.667) {};
    \node[vertex] (cb) at (0,-1.667) {};
    \node[vertex] (cr) at (.667,0) {};
    \node[vertex] (r1) at (1.222,0) {};
    \node[vertex] (r2) at (1.778,0) {};
    \node[vertex] (r) at (2.333,0) {};
    \draw[edge,blue] (c) to (cb);
    \draw[edge,black] (cl) to (l);
    \draw[Green] (r1) to (cr);
    \draw[Green,fill=white] (.9445,0) circle (.0926);
    \draw[TealBlue] (r2) to (r);
    \end{scope}
    \begin{scope}[xshift=12cm]
    \draw (0,0) node{$=$};
    \end{scope}
\end{tikzpicture}
\\
\vspace{\baselineskip}
\begin{tikzpicture}[scale=1.15]
    \begin{scope}[xshift=-6cm]
    \draw (0,0) node{$=$};
    \end{scope}
    \begin{scope}[xshift=-3cm]
    \draw[Magenta,edge,domain=5:175] plot ({.667*cos(\x)}, {.667*sin(\x)});
    \draw[Magenta,edge,domain=185:355] plot ({.667*cos(\x)}, {.667*sin(\x)});
    \draw[red,domain=10:80] plot ({.278*cos(\x)-1.5}, {.278*sin(\x)});
    \draw[red,domain=100:170] plot ({.278*cos(\x)-1.5}, {.278*sin(\x)});
    \draw[red,domain=190:350] plot ({.278*cos(\x)-1.5}, {.278*sin(\x)});
    \node[vertex] (l) at (-2.333,0) {};
    \node[vertex] (l1) at (-1.778,0) {};
    \node[vertex] (lc) at (-1.5,.278) {};
    \node[vertex] (lct) at (-1.5,.6) {};
    \node[vertex] (l2) at (-1.222,0) {};
    \node[vertex] (cl) at (-.667,0) {};
    \node[vertex] (cr) at (.667,0) {};
    \node[vertex] (r) at (2.333,0) {};
    \draw[edge,blue] (l1) to (l);
    \draw[Green] (l2) to (cl);
    \draw[Green,fill=white] (-.9445,0) circle (.0926);
    \draw[edge,TealBlue] (cr) to (r);
    \draw[black] (lc) to (lct);
    \end{scope}
    \draw[-to] (-1/8,0) to node[above]{$\beta \alpha$} (1/8,0);
    \begin{scope}[xshift=3cm]
    \draw[red,edge,domain=5:175] plot ({.667*cos(\x)}, {.667*sin(\x)});
    \draw[red,edge,domain=185:265] plot ({.667*cos(\x)}, {.667*sin(\x)});
    \draw[red,edge,domain=275:355] plot ({.667*cos(\x)}, {.667*sin(\x)});
    \draw[Magenta,domain=10:170] plot ({.278*cos(\x)+1.5}, {.278*sin(\x)});
    \draw[Magenta,domain=190:350] plot ({.278*cos(\x)+1.5}, {.278*sin(\x)});
    \node[vertex] (l) at (-2.333,0) {};
    \node[vertex] (cl) at (-.667,0) {};
    \node[vertex] (c) at (0,-.667) {};
    \node[vertex] (cb) at (0,-1.667) {};
    \node[vertex] (cr) at (.667,0) {};
    \node[vertex] (r1) at (1.222,0) {};
    \node[vertex] (r2) at (1.778,0) {};
    \node[vertex] (r) at (2.333,0) {};
    \draw[edge,blue] (c) to (cb);
    \draw[edge,black] (cl) to (l);
    \draw[Green] (r1) to (cr);
    \draw[Green,fill=white] (.9445,0) circle (.0926);
    \draw[TealBlue] (r2) to (r);
    \end{scope}
\end{tikzpicture}
\caption{The composition $\beta \alpha$ of two rearrangements.}
\label{fig.gpd.composition}
\end{figure}

For example, \cref{fig.gpd.composition} depicts the composition $\beta \alpha$, where $\alpha$ is the rearrangement depicted in \cref{fig.airplane.rearrangement,fig.gpd.expansion} and $\beta$ is depicted in \cref{fig.gpd.beta}.

It is easy to verify that, if a rearrangement is represented by a graph pair diagram $(D,\phi,R)$, then its inverse is represented by $(R, \phi^{-1}, D)$.

\section{Undirected Edges}
\label{sub.undirected.edges}

In certain cases, the orientation of an edge (i.e., whether the edge originates from $v$ and terminates at $w$ or originates from $w$ and terminates at $v$) is irrelevant for the purpose of computing its edge expansion.
For example, this is the case for blue edges of the airplane edge replacement system from \cref{fig.airplane.replacement}.
This is formalized by the following definition.

\begin{definition}
\label{def.undirected.edges}
Given a set of edge replacement rules (\cref{def.replacement.rules}), we say that a color $c$ is an \textbf{undirected color} and that $c$-colored edges are \textbf{undirected} if the replacement graph $\Gamma_c$ admits a graph automorphism $\psi$ that switches its initial and terminal vertices $\iota_c$ and $\tau_c$.
\end{definition}

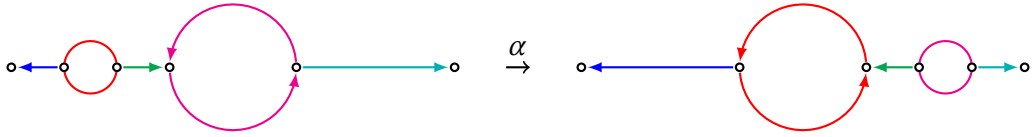
\begin{figure}\centering
\begin{tikzpicture}[scale=1.25]
    \begin{scope}[xshift=-3cm]
    \draw[Magenta,edge,domain=5:175] plot ({.667*cos(\x)}, {.667*sin(\x)});
    \draw[Magenta,edge,domain=185:355] plot ({.667*cos(\x)}, {.667*sin(\x)});
    \draw[red,domain=10:170] plot ({.278*cos(\x)-1.5}, {.278*sin(\x)});
    \draw[red,domain=190:350] plot ({.278*cos(\x)-1.5}, {.278*sin(\x)});
    \node[vertex] (l) at (-2.333,0) {};
    \node[vertex] (l1) at (-1.778,0) {};
    \node[vertex] (l2) at (-1.222,0) {};
    \node[vertex] (cl) at (-.667,0) {};
    \node[vertex] (cr) at (.667,0) {};
    \node[vertex] (r) at (2.333,0) {};
    \draw[edge,blue] (l1) to (l);
    \draw[edge,Green] (l2) to (cl);
    \draw[edge,TealBlue] (cr) to (r);
    \end{scope}
    \draw[-to] (-1/8,0) to node[above]{$\alpha$} (1/8,0);
    \begin{scope}[xshift=3cm]
    \draw[Magenta,domain=10:170] plot ({.278*cos(\x)+1.5}, {.278*sin(\x)});
    \draw[Magenta,domain=190:350] plot ({.278*cos(\x)+1.5}, {.278*sin(\x)});
    \draw[red,edge,domain=5:175] plot ({.667*cos(\x)}, {.667*sin(\x)});
    \draw[red,edge,domain=185:355] plot ({.667*cos(\x)}, {.667*sin(\x)});
    \node[vertex] (l) at (-2.333,0) {};
    \node[vertex] (cl) at (-.667,0) {};
    \node[vertex] (cr) at (.667,0) {};
    \node[vertex] (r1) at (1.222,0) {};
    \node[vertex] (r2) at (1.778,0) {};
    \node[vertex] (r) at (2.333,0) {};
    \draw[edge,blue] (cl) to (l);
    \draw[edge,Green] (r1) to (cr);
    \draw[edge,TealBlue] (r2) to (r);
    \end{scope}
    \end{tikzpicture}
    \caption{The same rearrangement of \cref{fig.airplane.rearrangement} represented using undirected blue edges. Note that the orientation of the edge represented in green is being reversed.}
    \label{fig.airplane.rearrangement.undirected}
\end{figure}

When a color $c$ is undirected, it is convenient to deal with $c$-colored edges as if they were not directed edges, meaning that a $c$-colored edge originating from $v$ and terminating at $w$ can also be treated as if it were originating from $w$ and terminating at $v$ instead.
In particular, graph isomorphisms are allowed to switch the orientation of $c$-colored edges, as happens for example to the edge highlighted in green in \cref{fig.airplane.rearrangement.undirected}.

This can be done under the convention that, when the orientation of an undirected edge $e$ is reversed, we are implicitly applying the automorphism $\psi$ that switches the initial and terminal vertices to the edge expansion of $e$.
In practice, we are ``hiding'' a fixed subgraph automorphism $\psi$ under the inversion of orientation.
In symbols, if the orientation of an edge $e$ of $D$ is reversed by a graph pair diagram $(D,\phi,R)$ then we mean that each $ex$ is mapped to $\phi(e) \psi(x)$.

Note that, if multiple automorphisms switch $\iota$ and $\tau$, one needs to fix one such isomorphism $\psi$.
In this case, if one wants to apply an orientation-inverting automorphism that is not $\psi$, the undirected edge needs to be expanded.

\section{Null-expanding Isolated Colors}
\label{sec.null.expanding}

As noted in \cref{rmk.generalized.limit.spaces}, sometimes it is possible to define limit spaces for non-expanding edge replacement system, but they may lack some nice features such as being Hausdorff.
There are also other reasons why it is often useful to assume that an edge replacement system in expanding (for example, throughout \cref{cha.conjugacy} we will assume this condition for combinatorial reasons).

Among non-expanding edge replacement systems, a particularly well-behaved case in that in which the only obstruction is having a color whose edges are isolated and do not expand in the sense of the following definition.

\begin{definition}
\label{def.null.expanding.isolated}
Given an edge replacement system, we say that a color $c$ is:
\begin{itemize}
    \item \textbf{isolated} if every $c$-colored edge of any graph expansion does not share vertices with any other edge;
    \item \textbf{null-expanding} if $q(c)$ belongs to an inescapable cycle (\cref{def.inescapable.cycle}) of the color graph (\cref{def.color.graph}).
\end{itemize}
\end{definition}

The following fact shows that null-expanding isolated colors are not a real obstruction, as we can create expanding edge replacement systems with the same rearrangement groups.

\medskip %layout
\begin{proposition}
\label{prop.null.expanding.isolated.rarrangement}
Let $\mathcal{R}$ be an edge replacement system and assume that the replacement graphs satisfy the requirements of expanding edge replacement systems (\cref{def.expanding}) with the exception of null-expanding isolated colors $c_1, \dots, c_k$.
Then there exists an expanding edge replacement system $\mathcal{R}^*$ whose rearrangement group is isomorphic to that of $\mathcal{R}$.
\end{proposition}

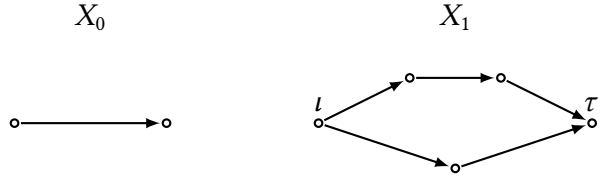
\begin{figure}
\centering
\begin{tikzpicture}
    \node at (1,1.4) {$X_0$};
    \node[vertex] (i) at (0,0) {};
    \node[vertex] (t) at (2,0) {};
    \draw[edge] (i) to (t);
    \begin{scope}[xshift=4cm]
    \node at (1.8,1.4) {$X_1$};
    \node[vertex] (l) at (0,0) {};
    \node[vertex] (t1) at (1.2,.6) {};
    \node[vertex] (t2) at (2.4,.6) {};
    \node[vertex] (b) at (1.8,-.6) {};
    \node[vertex] (r) at (3.6,0) {};
    \draw[edge] (l) node[above]{$\iota$} to (t1);
    \draw[edge] (t1) to (t2);
    \draw[edge] (t2) to (r);
    \draw[edge] (l) to (b);
    \draw[edge] (b) to (r) node[above]{$\tau$};
    \end{scope}
\end{tikzpicture}
\caption{An edge replacement system $\mathcal{N}$ with trivial rearrangement group.}
\label{fig.replacement.rearrangementless}
\end{figure}

\begin{sketch}
Our strategy is based on the edge replacement system $\mathcal{N}$ portrayed in \cref{fig.replacement.rearrangementless}.
Its key property is that, as noted in \cite[Example 2.5]{BF19}, this edge replacement system has trivial rearrangement group.
In this proof, let $\Gamma$ denote the replacement graph of $\mathcal{N}$.

Together, the two notions defined in \cref{def.null.expanding.isolated} imply that, for each null-expanding isolated color $c_i$, its replacement graph must consist of a sole edge colored by some null-expanding isolated color $c_j$ ($i$ and $j$ are colors in the same inescapable cycle).
With this in mind, we build the edge replacement system $\mathcal{R}^*$ by modifying the replacement graphs $X_{c_1}, \dots, X_{c_k}$ of $\mathcal{R}$ in the following way:
if $X_{c_i}$ consists of a $c_j$-colored edge, 
we replace it with a copy $\Gamma_{c_i}$ of $\Gamma$ whose edges are all $c_j$-colored.

Thus, every isolated point $p$ of the symbol space $\Omega_\mathcal{R}$ is replaced by a Cantor space that, under the gluing relation of $\mathcal{R}^*$ (which is an equivalence relation since $\mathcal{R}^*$ is expanding) descends to some cell $\mathcal{K}_p$ that is an isolated homeomorphic copy of the limit space of $\mathcal{N}$.

Now, it is not hard to see that there is an embedding $G_\mathcal{R} \to G_{\mathcal{R}^*}$ obtained by modifying each rearrangement of $\mathcal{R}$ so that each isolated point $p$ is replaced canonically by the cell $\mathcal{K}_p$.
Thus, the image of a rearrangement of $\mathcal{R}$ under this embedding maps canonically between copies of the limit space of $\mathcal{N}$ in the same way in which the original rearrangement of $\mathcal{R}$ maps between the isolated points.
Since all of these cells have trivial rearrangements and can be permuted freely, one can show that this embedding is surjective, so it is actually a group isomorphism, as needed. 
\end{sketch}

\section{Forest Pair Diagrams}
\label{sec.forest.pair.diagrams}

Since rearrangements are essentially isomorphisms between graph expansions (\cref{def.graph.pair.diagrams}) and since graph expansions can be represented by forest expansions (recall that the graph expansion is the leaf graph of forest expansions, \cref{def.leaf.graph}), it is possible to describe rearrangements as pairs of forest expansions together with some bijection between the leaves, as we illustrate in this section.
These diagrams are going to be useful in \cref{sub.IG.non.torsion} to show the existence of so-called \textit{wandering cells} and essential throughout \cref{cha.conjugacy} to define strand diagrams that provide a solution to the conjugacy problem.

Recall that, by \cref{ass.isolated.vertices}, our graphs do not feature isolated vertices.
Under this assumption, a graph isomorphism is entirely described by its action on the edges of the domain and range graphs.

\begin{definition}
\label{def.forest.pair.diagram}
A \textbf{forest pair diagram} is a triple $(F_D, \phi, F_R)$ such that $F_D$ and $F_R$ are forest expansions and $\phi$ is a graph isomorphism between the leaf graphs of $F_D$ and $F_R$.
We refer to $F_D$ and $F_R$ as the \textbf{domain} and \textbf{range} forests of the diagram.
\end{definition}

We will later see in \cref{rmk.forest.pair.diagrams.convention} how the graph isomorphism $\phi$ can be hidden thanks to a choice of labeling for $F_R$ that is induced from $F_D$.

Observe that there is an obvious bijection between the set of graph pair diagrams and that of forest pair diagrams:
\[ (D,\phi,R) \mapsto (F_D, \phi, F_R). \]
Here $\phi$ denotes two distinct maps, but this abuse of notation makes sense:
the first $\phi$ is a graph isomorphism $D \to R$, so it is a bijection between the edges of $D$ and $R$ that preserves adjacency;
the second $\phi$ is a bijection between leaves of $F_D$ and those of $F_R$ that is also a graph isomorphism between the leaf graphs $D$ and $R$.

\begin{figure}
    \centering
    \begin{tikzpicture}[font=\footnotesize,scale=1]
        \node[node] (root) at (0,0) {};
        \node[node] (s) at (0,-1) {};
        \node[node] (l) at (-2.5,-2) {};
            \node[node] (l1) at (-4,-3) {};
                \node[node] (l11) at (-5,-4) {};
                \node[node] (l12) at (-4,-4) {};
                \node[node] (l13) at (-3,-4) {};
                \node[node] (l14) at (-2,-4) {};
            \node[node] (l2) at (-1,-4) {};
            \node[node] (l3) at (0,-4) {};
            \node[node] (l4) at (1,-4) {};
        \node[node] (b) at (2,-4) {};
        \node[node] (r) at (3,-4) {};
        \node[node] (t) at (4,-4) {};
        \draw (root) to (s) node[blue,above left,align=center]{l\\r};
        \draw (s) to[out=180,in=90,looseness=1.1] (l) node[blue,above left,align=center]{x\\l};
            \draw (l) to[out=180,in=90,looseness=1.1] (l1) node[blue,above left,align=center]{v\\x};
                \draw (l1) to[out=180,in=90,looseness=1.25] (l11) node[blue,above left,align=center]{u\\v} node[black,below]{$1$};
                \draw (l1) to (l12) node[red,above left,align=center]{u\\z} node[black,below]{$2$};
                \draw (l1) to[out=0,in=90,looseness=1.25] (l13) node[blue,above left,xshift=-.2cm,align=center]{z\\x} node[black,below]{$3$};
                \draw (l1) to[out=0,in=90,looseness=1.25] (l14) node[red,above left,xshift=-.3cm,align=center]{z\\u} node[black,below]{$4$};
            \draw (l) to[out=0,in=90,looseness=1.25] (l2) node[red,above left,align=center]{v\\w} node[black,below]{$5$};
            \draw (l) to[out=0,in=90,looseness=1.25] (l3) node[blue,above left,align=center]{w\\l} node[black,below]{$6$};
            \draw (l) to[out=0,in=90,looseness=1.25] (l4) node[red,above left,align=center]{w\\v} node[black,below]{$7$};
        \draw (s) to[out=0,in=90,looseness=1.25] ($(b)+(0,1)$) to (b) node[red,above left,align=center]{x\\y} node[black,below]{$8$};
        \draw (s) to[out=0,in=90,looseness=1.25] ($(r)+(0,1)$) to (r) node[blue,above left,align=center]{y\\r} node[black,below]{$9$};
        \draw (s) to[out=0,in=90,looseness=1.25] ($(t)+(0,1)$) to (t) node[red,above left,align=center]{y\\x} node[black,below]{$10$};
        \draw[-to] (0,-4.85) to node[above]{} (0,-5.15);
        \begin{scope}[yshift=-5.8cm]
        \node[node] (root) at (0,0) {};
        \node[node] (s) at (0,-1) {};
        \node[node] (l) at (-4.5,-4) {};
        \node[node] (b) at (-3.5,-4) {};
        \node[node] (r) at (.5,-2) {};
            \node[node] (r1) at (-1,-3) {};
                \node[node] (r11) at (-2.5,-4) {};
                \node[node] (r12) at (-1.5,-4) {};
                \node[node] (r13) at (-.5,-4) {};
                \node[node] (r14) at (.5,-4) {};
            \node[node] (r2) at (1.5,-4) {};
            \node[node] (r3) at (2.5,-4) {};
            \node[node] (r4) at (3.5,-4) {};
        \node[node] (t) at (4.5,-4) {};
        \draw (root) to (s) node[blue,above left,align=center]{l\\r};
        \draw (s) to[out=180,in=90,looseness=1.1] ($(l)+(0,1)$) to (l) node[blue,above left,align=center]{x\\l}  node[black,below]{$\phi(6)$};
        \draw (s) to[out=180,in=90,looseness=1.25] ($(b)+(0,1)$) to (b) node[red,above left,align=center]{x\\y}  node[black,below]{$\phi(7)$};
        \draw (s) to[out=0,in=90,looseness=1.25] (r) node[blue,above left,align=center]{y\\r};
            \draw (r) to[out=180,in=90,looseness=1.25] (r1) node[blue,above left,align=center]{u\\y};
                \draw (r1) to[out=180,in=90,looseness=1.25] (r11) node[blue,above left,align=center]{a\\u}  node[black,below]{$\phi(3)$};
                \draw (r1) to[out=180,in=90,looseness=1.25] (r12) node[red,above left,align=center]{a\\b}  node[black,below]{$\phi(2)$};
                \draw (r1) to[out=0,in=90,looseness=1.25] (r13) node[blue,above left,align=center]{b\\y}  node[black,below]{$\phi(1)$};
                \draw (r1) to[out=0,in=90,looseness=1.25] (r14) node[red,above left,xshift=-.3cm,align=center]{b\\a}  node[black,below]{$\phi(4)$};
            \draw (r) to[out=0,in=90,looseness=1.25] (r2) node[red,above left,align=center]{u\\p}  node[black,below]{$\phi(8)$};
            \draw (r) to[out=0,in=90,looseness=1.25] (r3) node[blue,above left,align=center]{p\\r}  node[black,below]{$\phi(9)$};
            \draw (r) to[out=0,in=90,looseness=1.25] (r4) node[red,above left,align=center]{p\\u}  node[black,below]{$\phi(10)$};
        \draw (s) to[out=0,in=90,looseness=1.1] ($(t)+(0,1)$) to (t) node[red,above left,align=center]{y\\x}  node[black,below]{$\phi(5)$};
        \end{scope}
    \end{tikzpicture}
    \caption{A forest pair diagram of the airplane edge replacement system that is equivalent to the graph pair diagram in \cref{fig.airplane.rearrangement}.}
    \label{fig.airplane.forest.pair.diagram.directed}
\end{figure}
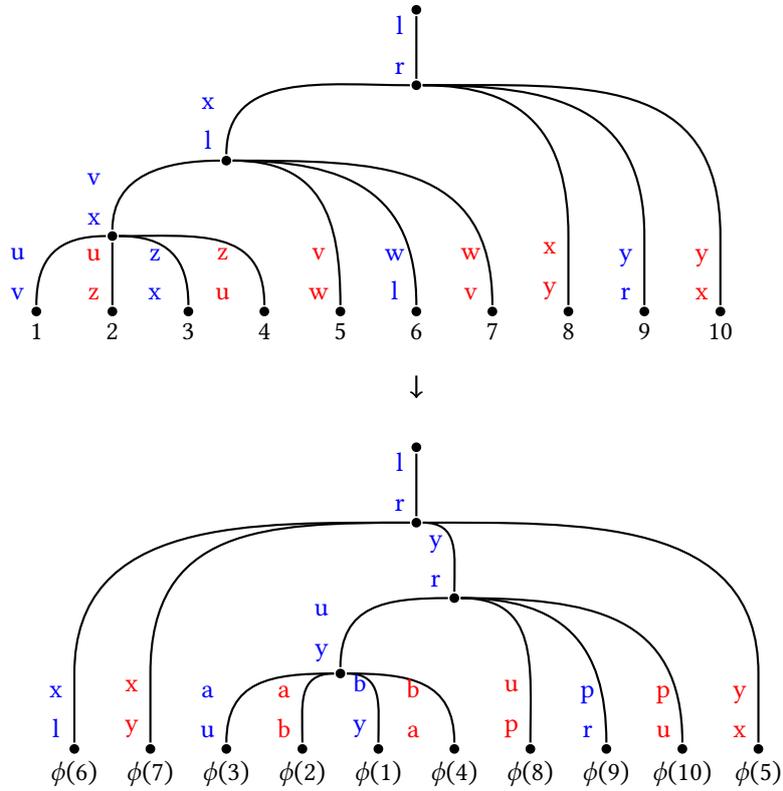

An example of a forest pair diagram for the airplane edge replacement system is displayed in \cref{fig.airplane.forest.pair.diagram.directed}.

\begin{figure}
    \centering
    \begin{tikzpicture}[font=\footnotesize,scale=1]
        \node[node] (root) at (0,0) {};
        \node[node] (s) at (0,-1) {};
        \node[node] (l) at (-1.5,-2) {};
            \node[node] (l1) at (-3,-3) {};
            \node[node] (l2) at (-2,-3) {};
            \node[node] (l3) at (-1,-3) {};
            \node[node] (l4) at (0,-3) {};
        \node[node] (b) at (1,-3) {};
        \node[node] (r) at (2,-3) {};
        \node[node] (t) at (3,-3) {};
        \draw (root) to (s) node[blue,above left,align=center]{l\\r};
        \draw (s) to[out=180,in=90,looseness=1.25] (l) node[blue,above left,align=center]{x\\l};
            \draw (l) to[out=180,in=90,looseness=1.25] (l1) node[blue,above left,align=center]{v\\x} node[black,below]{$1$};
            \draw (l) to[out=180,in=90,looseness=1.25] (l2) node[red,above left,align=center]{v\\w} node[black,below]{$2$};
            \draw (l) to[out=0,in=90,looseness=1.25] (l3) node[blue,above left,align=center]{w\\l} node[black,below]{$3$};
            \draw (l) to[out=0,in=90,looseness=1.25] (l4) node[red,above left,xshift=-.3cm,align=center]{w\\v} node[black,below]{$4$};
        \draw (s) to[out=0,in=90,looseness=1.25] ($(b)+(0,1)$) to (b) node[red,above left,align=center]{x\\y} node[black,below]{$5$};
        \draw (s) to[out=0,in=90,looseness=1.2] ($(r)+(0,1)$) to (r) node[blue,above left,align=center]{y\\r} node[black,below]{$6$};
        \draw (s) to[out=0,in=90,looseness=1] ($(t)+(0,1)$) to (t) node[red,above left,align=center]{y\\x} node[black,below]{$7$};
        \draw[-to] (3.45,-1.75) to node[above]{} (3.75,-1.75);
        \begin{scope}[xshift=7.2cm]
        \node[node] (root) at (0,0) {};
        \node[node] (s) at (0,-1) {};
        \node[node] (l) at (-3,-3) {};
        \node[node] (b) at (-2,-3) {};
        \node[node] (r) at (.5,-2) {};
            \node[node] (r1) at (-1,-3) {};
            \node[node] (r2) at (0,-3) {};
            \node[node] (r3) at (1,-3) {};
            \node[node] (r4) at (2,-3) {};
        \node[node] (t) at (3,-3) {};
        \draw (root) to (s) node[blue,above left,align=center]{l\\r};
        \draw (s) to[out=180,in=90,looseness=1] ($(l)+(0,1)$) to (l) node[blue,above left,align=center]{x\\l}  node[black,below]{$\phi(3)$};
        \draw (s) to[out=180,in=90,looseness=1.25] ($(b)+(0,1)$) to (b) node[red,above left,align=center]{x\\y}  node[black,below]{$\phi(4)$};
        \draw (s) to[out=0,in=90,looseness=1.25] (r) node[blue,above left,align=center]{y\\r};
            \draw (r) to[out=180,in=90,looseness=1.25] (r1) node[blue,above left,align=center]{u\\y}  node[black,below]{$\phi(1)$};
            \draw (r) to[out=180,in=90,looseness=1.25] (r2) node[red,above left,align=center]{u\\p}  node[black,below]{$\phi(5)$};
            \draw (r) to[out=0,in=90,looseness=1.25] (r3) node[blue,above left,align=center]{p\\r}  node[black,below]{$\phi(6)$};
            \draw (r) to[out=0,in=90,looseness=1.25] (r4) node[red,above left,xshift=-.3cm,align=center]{p\\u}  node[black,below]{$\phi(7)$};
        \draw (s) to[out=0,in=90,looseness=1] ($(t)+(0,1)$) to (t) node[red,above left,align=center]{y\\x}  node[black,below]{$\phi(2)$};
        \end{scope}
    \end{tikzpicture}
    \caption{A forest pair diagram of the airplane edge replacement system that is equivalent to the graph pair diagram in \cref{fig.airplane.rearrangement.undirected} using undirected edges. The first edge of the domain leaf graph inverts its orientation, i.e., the edge $(v,x)$ of the domain graph is mapped to the edge $(u,y)$ of the range graph, but the vertex $v$ is mapped to $y$ and $x$ to $u$.}
    \label{fig.airplane.forest.pair.diagram.old}
\end{figure}
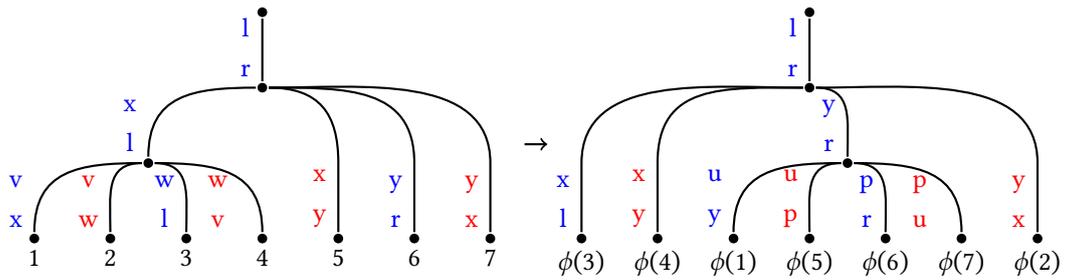

It is convenient to use the convention explained in \cref{sub.undirected.edges}, which allows us to treat any undirected edge of the leaf graph that is originating and terminating at some $v$ and $w$, respectively, as if it were originating from $w$ and terminating at $v$ for the purpose of determining a graph isomorphism of the leaf graphs.
As explained in \cref{sub.undirected.edges}, we are ``hiding'' a fixed orientation-inverting graph isomorphism between the edge expansions of the inverted edge.
As an example, \cref{fig.airplane.forest.pair.diagram.old} depicts the same rearrangement from \cref{fig.airplane.forest.pair.diagram.directed} and the two correspond to the graph pair diagrams portrayed in \cref{fig.airplane.rearrangement} and \cref{fig.airplane.rearrangement.undirected}, respectively.

\begin{figure}
    \centering
    \begin{tikzpicture}[font=\footnotesize,scale=1]
        \node[node] (root) at (0,0) {};
        \node[node] (s) at (0,-1) {};
        \node[node] (l) at (-1.5,-2) {};
            \node[node] (l1) at (-3,-3) {};
            \node[node] (l2) at (-2,-3) {};
            \node[node] (l3) at (-1,-3) {};
            \node[node] (l4) at (0,-3) {};
        \node[node] (b) at (1,-3) {};
        \node[node] (r) at (2,-3) {};
        \node[node] (t) at (3,-3) {};
        \draw (root) to (s) node[blue,above left,align=center]{l\\r};
        \draw (s) to[out=180,in=90,looseness=1.25] (l) node[blue,above left,align=center]{x\\l};
            \draw (l) to[out=180,in=90,looseness=1.25] (l1) node[blue,above left,align=center]{v\\x};
            \draw (l) to[out=180,in=90,looseness=1.25] (l2) node[red,above left,align=center]{v\\w};
            \draw (l) to[out=0,in=90,looseness=1.25] (l3) node[blue,above left,align=center]{w\\l};
            \draw (l) to[out=0,in=90,looseness=1.25] (l4) node[red,above left,xshift=-.3cm,align=center]{w\\v};
        \draw (s) to[out=0,in=90,looseness=1.25] ($(b)+(0,1)$) to (b) node[red,above left,align=center]{x\\y};
        \draw (s) to[out=0,in=90,looseness=1.2] ($(r)+(0,1)$) to (r) node[blue,above left,align=center]{y\\r};
        \draw (s) to[out=0,in=90,looseness=1] ($(t)+(0,1)$) to (t) node[red,above left,align=center]{y\\x};
        \draw[-to] (3.45,-1.75) to node[above]{} (3.75,-1.75);
        \begin{scope}[xshift=7.2cm]
        \node[node] (root) at (0,0) {};
        \node[node] (s) at (0,-1) {};
        \node[node] (l) at (-3,-3) {};
        \node[node] (b) at (-2,-3) {};
        \node[node] (r) at (.5,-2) {};
            \node[node] (r1) at (-1,-3) {};
            \node[node] (r2) at (0,-3) {};
            \node[node] (r3) at (1,-3) {};
            \node[node] (r4) at (2,-3) {};
        \node[node] (t) at (3,-3) {};
        \draw (root) to (s) node[blue,above left,align=center]{l\\r};
        \draw (s) to[out=180,in=90,looseness=1] ($(l)+(0,1)$) to (l) node[blue,above left,align=center]{w\\l};
        \draw (s) to[out=180,in=90,looseness=1.25] ($(b)+(0,1)$) to (b) node[red,above left,align=center]{w\\v}  node[black,below]{};
        \draw (s) to[out=0,in=90,looseness=1.25] (r) node[blue,above left,align=center]{v\\r};
            \draw (r) to[out=180,in=90,looseness=1.25] (r1) node[blue,above left,align=center]{x\\v};
            \draw (r) to[out=180,in=90,looseness=1.25] (r2) node[red,above left,align=center]{x\\y};
            \draw (r) to[out=0,in=90,looseness=1.25] (r3) node[blue,above left,align=center]{y\\r};
            \draw (r) to[out=0,in=90,looseness=1.25] (r4) node[red,above left,xshift=-.3cm,align=center]{y\\x};
        \draw (s) to[out=0,in=90,looseness=1] ($(t)+(0,1)$) to (t) node[red,above left,align=center]{v\\w};
        \end{scope}
    \end{tikzpicture}
    \caption{A forest pair diagram for the same element represented in \cref{fig.airplane.forest.pair.diagram.old} after a renaming of symbols such that the permutation $\phi$ of the bottom branches is determined by the labeling.}
    \label{fig.airplane.forest.pair.diagram}
\end{figure}
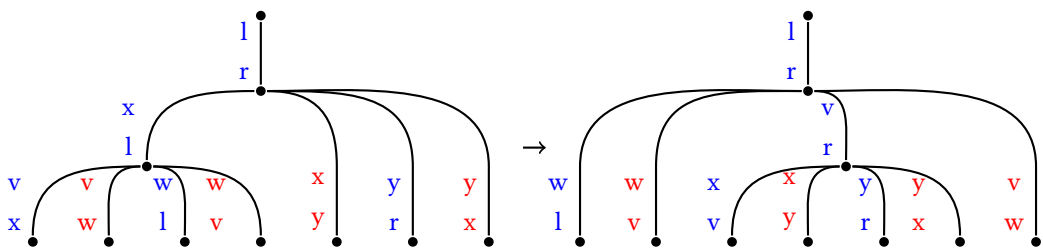

\begin{remark}
\label{rmk.forest.pair.diagrams.convention}
Since we can rename the symbols of a forest expansion (as described in \cref{rmk.renaming}), we can always simply name each vertex in the range forest $F_R$ by the name of its preimage in $F_D$ under $\sigma$, which for example results in \cref{fig.airplane.forest.pair.diagram} for the same element of \cref{fig.airplane.forest.pair.diagram.old}.
Again, in \cref{fig.airplane.forest.pair.diagram} $(v,x)$ is being mapped to $(x,v)$, which is allowed because blue edges here are undirected (\cref{sub.undirected.edges}).
Then the image under $\phi$ of a leaf of $F_D$ labeled by $(v,w,z)$ is precisely the leaf of $F_R$ labeled by the same triple, and so the isomorphism $\phi$ between the leaf graphs is described by the labeling itself. 
With this choice, the forest pair diagram is thus entirely determined by the pair $(F_D, F_R)$.
From now on, we will always be using this simpler notation when referring to forest pair diagrams.

Moreover, in accordance with \cref{rmk.renaming}, given a forest pair diagram $(F_D, F_R)$, we can rename symbols of both $F_D$ and $F_R$ in a coherent way, by which we mean that a symbol $a$ is changed to a symbol $b$ in the domain forest $F_D$ if and only if the symbol $a$ is changed to the symbol $b$ in the range forest $F_R$.
If two forest pair diagrams differ from such a renaming of symbols, they represent the same rearrangement, so we consider them to be the same forest pair diagram.
\end{remark}

\phantomsection\label{txt.reduced.forest.pair.diagrams}
Similarly to how graph pair diagrams can be reduced and expanded (\cref{sub.reduced.graph.pair.diagrams}), forest pair diagrams can also be reduced and expanded by simply removing or adding a replacement tree in $F_D$ and its image in $F_R$ at one of the bottom branches.
When reducing, the permutation of the bottom branches of that instance of the replacement tree must be trivial.
When expanding, one needs to use new symbols that do not appear elsewhere.
This corresponds precisely to expansions and reductions of graph pair diagrams, so each rearrangement is represented by a unique reduced forest pair diagram.

\phantomsection\label{txt.composition.forest.pair.diagrams}
Composition of forest pair diagrams is similar to the composition of graph pair diagrams (which was illustrated in \cref{fig.gpd.composition}):
given two forest pair diagrams for two elements $g$ and $h$, one needs to compute expansions $(F_D^g,F_R^g)$ and $(F_D^h,F_R^h)$ such that $F_D^g = F_R^h$ (using \cref{prop.common.finer.partition}).
Then, carefully renaming the symbols of either the first or the second diagram so that $F_D^g$ and $F_R^h$ have exactly the same labels will produce a forest pair diagram $(F_D^h,F_R^g)$.
Since we will mostly be composing strand diagrams, which follow the same rules of composition of forest pair diagrams, we will give more details in \cref{sub.SDs.composition}, when we describe the composition of strand diagrams.

\section{Examples of Rearrangement Groups}
\label{sec.examples}

\subsection{Thompson Groups}
\label{sub.Thompson}

The celebrated trio of Thompson groups $F$, $T$ and $V$, introduced by mathematician Richard J. Thompson in the '60s, provide the main inspiration for rearrangement groups.
They have various different (yet equivalent) definitions, which highlights their versatility and ubiquity through many topics in mathematics.
Here we will only provide a basic treatment on Thompson groups.
For a standard introduction see \cite{cfp}.

One of the main reasons for the initial interest in Thompson groups concerns their finiteness properties.
A finitely generated or finitely presented group is said to be of type $F_1$ or $F_2$, respectively.
More generally, we say that a group $G$ has finiteness property $F_n$ if it acts freely, faithfully, properly, cellularly, and cocompactly on an $(n-1)$-connected cell complex.
Moreover, a group has finiteness property $F_\infty$ if it is $F_n$ for every $n \in \mathbb{N}$.
We will not deal with finiteness properties except for finite generation, but an interested reader may look at \cite{Brown} for more details on finiteness properties of Thompson-like groups.

\begin{remark}
Throughout this section, a reader that is familiar with Thompson groups will see that the forest pair diagrams from \cref{def.forest.pair.diagram} produce the usual tree pair diagrams that are commonly found in the literature to represent the elements of Thompson groups, except that our diagrams have labels.
Such labels only imply that the bijection between the leaves of the domain and the range must be trivial for $F$, cyclic for $T$ or any for $V$.
This suggests that most rearrangement groups sit between Thompson groups $F$ and $V$.
Indeed, \cite[Proposition 2.8]{BF19} shows that $F$ embeds in many rearrangement groups and a way to embed rearrangement groups in $V$ is explored in \cref{sec.embedding.into.V}.
\end{remark}

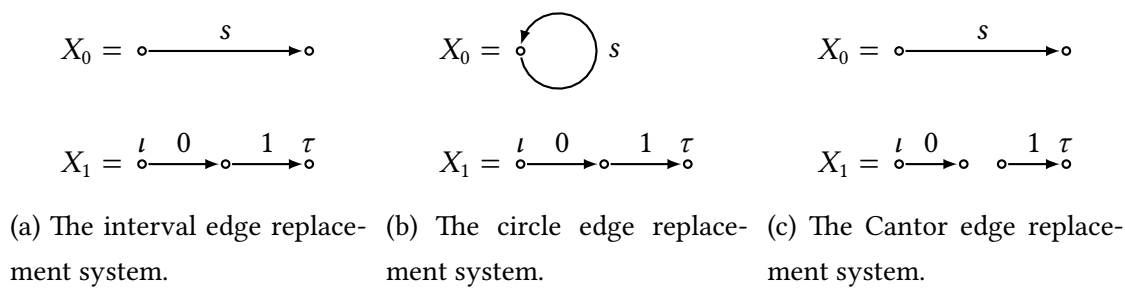
\begin{figure}
\centering
\begin{subfigure}[B]{.32\textwidth}
\centering
\begin{tikzpicture}
    \node at (-.667,0) {$X_0 =$};
    \node[vertex] (i) at (0,0) {};
    \node[vertex] (t) at (2.2,0) {};
    \draw[edge] (i) to node[above]{$s$} (t);
    \begin{scope}[yshift=-1.5cm]
    \node at (-.667,0) {$X_1 =$};
    \node[vertex] (l) at (0,0) {};
    \draw (l) node[above]{$\iota$};
    \node[vertex] (c) at (1.1,0) {};
    \node[vertex] (r) at (2.2,0) {};
    \draw (r) node[above]{$\tau$};
    \draw[edge] (l) to node[above]{$0$} (c);
    \draw[edge] (c) to node[above]{$1$} (r);
    \end{scope}
\end{tikzpicture}
\caption{The interval edge replacement system.}
\label{fig.interval.replacement}
\end{subfigure}
\hfill
\begin{subfigure}[B]{.325\textwidth}
\centering
\begin{tikzpicture}
    \node at (-.667,0) {$X_0 =$};
    \node[vertex] (i) at (0,0) {};
    \draw[edge,domain=190:530] plot ({.5*cos(\x)+.5}, {.5*sin(\x)});
    \draw (1,0) node[right] {$s$};
    \begin{scope}[yshift=-1.5cm]
    \node at (-.667,0) {$X_1 =$};
    \node[vertex] (l) at (0,0) {};
    \draw (l) node[above]{$\iota$};
    \node[vertex] (c) at (1.1,0) {};
    \node[vertex] (r) at (2.2,0) {};
    \draw (r) node[above]{$\tau$};
    \draw[edge] (l) to node[above]{$0$} (c);
    \draw[edge] (c) to node[above]{$1$} (r);
    \end{scope}
\end{tikzpicture}
\caption{The circle edge replacement system.}
\label{fig.circle.replacement}
\end{subfigure}
\hfill
\begin{subfigure}[B]{.32\textwidth}
\centering
\begin{tikzpicture}
    \node at (-.667,0) {$X_0 =$};
    \node[vertex] (i) at (0,0) {};
    \node[vertex] (t) at (2.2,0) {};
    \draw[edge] (i) to node[above]{$s$} (t);
    \begin{scope}[yshift=-1.5cm]
    \node at (-.667,0) {$X_1 =$};
    \node[vertex] (l) at (0,0) {};
    \draw (l) node[above]{$\iota$};
    \node[vertex] (c1) at (.85,0) {};
    \node[vertex] (c2) at (1.35,0) {};
    \node[vertex] (r) at (2.2,0) {};
    \draw (r) node[above]{$\tau$};
    \draw[edge] (l) to node[above]{$0$} (c1);
    \draw[edge] (c2) to node[above]{$1$} (r);
    \end{scope}
\end{tikzpicture}
\caption{The Cantor edge replacement system.}
\label{fig.Cantor.replacement}
\end{subfigure}
\caption{The edge replacement systems for Thompson groups.}
\label{fig.Thompson.replacement}
\end{figure}

\subsubsection{Thompson's group \texorpdfstring{$F$}{F}}
\label{sub.F}

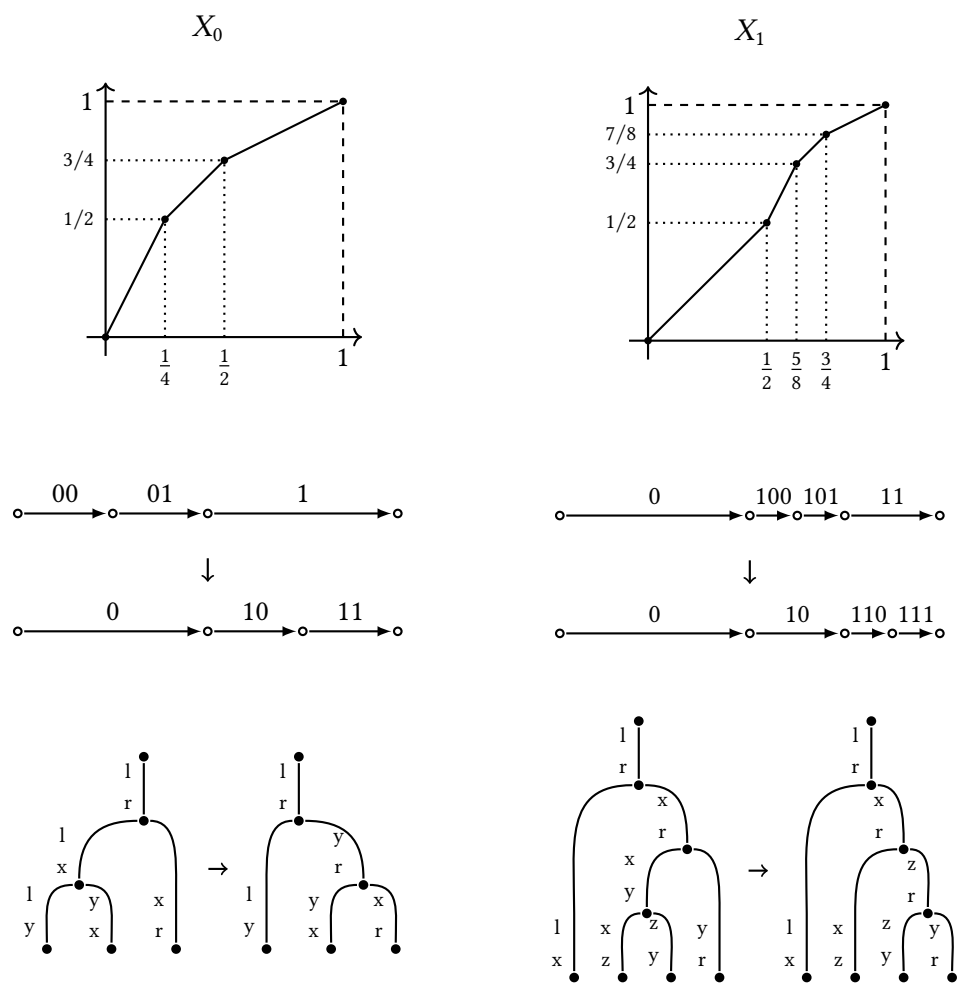
\begin{figure}
\centering
\begin{minipage}{.48\textwidth}
\centering
$X_0$
\\
\vspace{.5cm}
% X_0 as a map
\begin{tikzpicture}[scale=1.25]
    \draw[->] (-0.2, 0) -- (2.7, 0);
    \draw[->] (0, -0.2) -- (0, 2.7);
    \node at (2.5,2.5) [circle,fill,inner sep=1]{};
    \node at (0,0) [circle,fill,inner sep=1]{};
    \draw[dashed] (0,2.5) node[left]{\small$1$} -- (2.5,2.5);
    \draw[dashed] (2.5,0) node[below]{\small$1$} -- (2.5,2.5);
    \draw[scale=2.5, domain=0:0.25, smooth, variable=\x] plot({\x}, {2*\x});
    \node at (2.5*0.25,2.5*0.5) [circle,fill,inner sep=1]{};
    \draw[dotted] (2.5*0.25,0) node[below] {$\frac{1}{4}$} -- (2.5*0.25,2.5*0.5);
    \draw[dotted] (0,2.5*0.5) node[left] {\scriptsize$1/2$} -- (2.5*0.25,2.5*0.5);
    \draw[scale=2.5, domain=0.25:0.5, smooth, variable=\x] plot({\x},{\x+0.25});
    \node at (2.5*0.5,2.5*0.75) [circle,fill,inner sep=1]{};
    \draw[dotted] (2.5*0.5,0) node[below] {$\frac{1}{2}$} -- (2.5*0.5,2.5*0.75);
    \draw[dotted] (0,2.5*0.75) node[left] {\scriptsize$3/4$} -- (2.5*0.5,2.5*0.75);
    \draw[scale=2.5, domain=0.5:1, smooth, variable=\x] plot({\x},{0.5*\x+0.5});
\end{tikzpicture}
\\
\vspace{1cm}
% X_0 as a GPD
\begin{tikzpicture}[scale=1.25]
    \begin{scope}[yshift=.6cm]
    \node[vertex] (0) at (0,0) {};
    \node[vertex] (1/4) at (1,0) {};
    \node[vertex] (1/2) at (2,0) {};
    \node[vertex] (1) at (4,0) {};
    \draw[edge] (0) to node[above]{\small$00$} (1/4);
    \draw[edge] (1/4) to node[above]{\small$01$} (1/2);
    \draw[edge] (1/2) to node[above]{\small$1$} (1);
    \end{scope}
    \draw[-to] (2,1/8) to (2,-1/8);
    \begin{scope}[yshift=-.65cm]
    \node[vertex] (0) at (0,0) {};
    \node[vertex] (1/2) at (2,0) {};
    \node[vertex] (3/4) at (3,0) {};
    \node[vertex] (1) at (4,0) {};
    \draw[edge] (0) to node[above]{\small$0$} (1/2);
    \draw[edge] (1/2) to node[above]{\small$10$} (3/4);
    \draw[edge] (3/4) to node[above]{\small$11$} (1);
    \end{scope}
\end{tikzpicture}
\\
\vspace{1.5cm}
% X_0 as an FPD
    \begin{tikzpicture}[font=\scriptsize,scale=.85]
        \node[node] (root) at (0,0) {};
        \node[node] (s) at (0,-1) {};
        \node[node] (0) at (-1,-2) {};
            \node[node] (00) at (-1.5,-3) {};
            \node[node] (01) at (-.5,-3) {};
        \node[node] (1) at (.5,-3) {};
        \draw (root) to (s) node[above left,align=center]{l\\r};
        \draw (s) to[out=180,in=90,looseness=1.25] (0) node[above left,align=center]{l\\x};
            \draw (0) to[out=180,in=90,looseness=1.25] (00) node[above left,align=center]{l\\y};
            \draw (0) to[out=0,in=90,looseness=1.25] (01) node[above left,align=center]{y\\x};
        \draw (s) to[out=0,in=90,looseness=1] ($(1)+(0,1)$) to (1) node[above left,align=center]{x\\r};
        \draw[-to] (1,-1.75) to node[above]{} (1.3,-1.75);
        \begin{scope}[xshift=2.4cm]
        \node[node] (root) at (0,0) {};
        \node[node] (s) at (0,-1) {};
        \node[node] (0) at (-.5,-3) {};
        \node[node] (1) at (1,-2) {};
            \node[node] (10) at (.5,-3) {};
            \node[node] (11) at (1.5,-3) {};
        \node[node,white] at (0,-3.5) {};
        \draw (root) to (s) node[above left,align=center]{l\\r};
        \draw (s) to[out=180,in=90,looseness=1] ($(0)+(0,1)$) to (0) node[above left,align=center]{l\\y};
        \draw (s) to[out=0,in=90,looseness=1.25] (1) node[above left,xshift=-.1cm,align=center]{y\\r};
            \draw (1) to[out=180,in=90,looseness=1.25] (10) node[above left,align=center]{y\\x};
            \draw (1) to[out=0,in=90,looseness=1.25] (11) node[above left,align=center]{x\\r};
        \end{scope}
    \end{tikzpicture}
\end{minipage}
\begin{minipage}{.48\textwidth}
\centering
$X_1$
\\
\vspace{.5cm}
% X_1 as a map
\begin{tikzpicture}[scale=1.25]
    \draw[->] (-0.2, 0) -- (2.7, 0);
    \draw[->] (0, -0.2) -- (0, 2.7);
    \node at (2.5,2.5) [circle,fill,inner sep=1]{};
    \node at (0,0) [circle,fill,inner sep=1]{};
    \draw[dashed] (0,2.5) node[left]{\small$1$} -- (2.5,2.5);
    \draw[dashed] (2.5,0) node[below]{\small$1$} -- (2.5,2.5);
    \draw[scale=2.5, domain=0:0.5, smooth, variable=\x] plot({\x}, {\x});
    \node at (2.5*0.5,2.5*0.5) [circle,fill,inner sep=1]{};
    \draw[dotted] (2.5*0.5,0) node[below] {$\frac{1}{2}$} -- (2.5*0.5,2.5*0.5);
    \draw[dotted] (0,2.5*0.5) node[left] {\scriptsize$1/2$} -- (2.5*0.5,2.5*0.5);
    \draw[scale=2.5, domain=0.5:0.625, smooth, variable=\x] plot({\x},{2*\x-0.5});
    \node at (2.5*5/8,2.5*0.75) [circle,fill,inner sep=1]{};
    \draw[dotted] (2.5*5/8,0) node[below] {$\frac{5}{8}$} -- (2.5*5/8,2.5*0.75);
    \draw[dotted] (2.5*0.75,0) node[below] {$\frac{3}{4}$} -- (2.5*0.75,2.5*7/8);
    \draw[scale=2.5, domain=0.625:0.75, smooth, variable=\x] plot({\x},{\x+0.125});
    \node at (2.5*0.75,2.5*7/8) [circle,fill,inner sep=1]{};
    \draw[dotted] (0,2.5*0.75) node[left] {\scriptsize$3/4$} -- (2.5*5/8,2.5*0.75);
    \draw[dotted] (0,2.5*7/8) node[left] {\scriptsize$7/8$} -- (2.5*0.75,2.5*7/8);
    \draw[scale=2.5, domain=0.75:1, smooth, variable=\x] plot({\x},{0.5*\x+0.5});
\end{tikzpicture}
\\
\vspace{1cm}
% X_1 as a GPD
\begin{tikzpicture}[scale=1.25]
    \begin{scope}[yshift=.6cm]
    \node[vertex] (0) at (0,0) {};
    \node[vertex] (1/2) at (2,0) {};
    \node[vertex] (5/8) at (2.5,0) {};
    \node[vertex] (3/4) at (3,0) {};
    \node[vertex] (1) at (4,0) {};
    \draw[edge] (0) to node[above]{\footnotesize$0$} (1/2);
    \draw[edge] (1/2) to node[above]{\footnotesize$100$} (5/8);
    \draw[edge] (5/8) to node[above]{\footnotesize$101$} (3/4);
    \draw[edge] (3/4) to node[above]{\footnotesize$11$} (1);
    \end{scope}
    \draw[-to] (2,1/8) to (2,-1/8);
    \begin{scope}[yshift=-.65cm]
    \node[vertex] (0) at (0,0) {};
    \node[vertex] (1/2) at (2,0) {};
    \node[vertex] (3/4) at (3,0) {};
    \node[vertex] (7/8) at (3.5,0) {};
    \node[vertex] (1) at (4,0) {};
    \draw[edge] (0) to node[above]{\footnotesize$0$} (1/2);
    \draw[edge] (1/2) to node[above]{\footnotesize$10$} (3/4);
    \draw[edge] (3/4) to node[above]{\footnotesize$110$} (7/8);
    \draw[edge] (7/8) to node[above]{\footnotesize$111$} (1);
    \end{scope}
\end{tikzpicture}
\\
\vspace{1cm}
% X_1 as a FPD
    \begin{tikzpicture}[font=\scriptsize,scale=.85]
        \node[node] (root) at (0,0) {};
        \node[node] (s) at (0,-1) {};
        \node[node] (0) at (-1,-4) {};
        \node[node] (1) at (.75,-2) {};
            \node[node] (10) at (.125,-3) {};
                \node[node] (100) at (-.25,-4) {};
                \node[node] (101) at (.5,-4) {};
            \node[node] (11) at (1.25,-4) {};
        \draw (root) to (s) node[above left,align=center]{l\\r};
        \draw (s) to[out=180,in=90,looseness=1.25] ($(0)+(0,1)$) to (0) node[above left,align=center]{l\\x};
        \draw (s) to[out=0,in=90,looseness=1.25] (1) node[above left,xshift=-.1cm,align=center]{x\\r};
            \draw (1) to[out=180,in=90,looseness=1.25] (10) node[above left,align=center]{x\\y};
                \draw (10) to[out=180,in=90,looseness=1] (100) node[above left,align=center]{x\\z};
                \draw (10) to[out=0,in=90,looseness=1] (101) node[above left,align=center]{z\\y};
            \draw (1) to[out=0,in=90,looseness=1] ($(11)+(0,1)$) to (11) node[above left,align=center]{y\\r};
        \draw[-to] (1.7,-2.333) to node[above]{} (2,-2.333);
        \begin{scope}[xshift=3.6cm]
        \node[node] (root) at (0,0) {};
        \node[node] (s) at (0,-1) {};
        \node[node] (0) at (-1,-4) {};
        \node[node] (1) at (.5,-2) {};
            \node[node] (10) at (-.25,-4) {};
            \node[node] (11) at (.875,-3) {};
                \node[node] (110) at (.5,-4) {};
                \node[node] (111) at (1.25,-4) {};
        \draw (root) to (s) node[above left,align=center]{l\\r};
        \draw (s) to[out=180,in=90,looseness=1.25] ($(0)+(0,1)$) to (0) node[above left,align=center]{l\\x};
        \draw (s) to[out=0,in=90,looseness=1] (1) node[above left,xshift=-.1cm,align=center]{x\\r};
            \draw (1) to[out=180,in=90,looseness=1.25] ($(10)+(0,1)$) to (10) node[above left,align=center]{x\\z};
            \draw (1) to[out=0,in=90,looseness=1] (11) node[above left,align=center]{z\\r};
                \draw (11) to[out=180,in=90,looseness=1] (110) node[above left,align=center]{z\\y};
                \draw (11) to[out=0,in=90,looseness=1] (111) node[above left,align=center]{y\\r};
        \end{scope}
    \end{tikzpicture}
\end{minipage}
\caption{The generators $X_0$ and $X_1$ of Thompson's group $F$.}
\label{fig.F.generators}
\end{figure}

The smallest of the three Thompson groups, commonly denoted by $F$, can be defined as the group of those orientation-preserving piecewise-linear homeomorphisms of the unit interval $[0,1]$ with finitely many breakpoints, all of which are dyadic (i.e., they belong to $\mathbb{Z}[\frac{1}{2}]$), and such that the slopes are integer powers of $2$.
It admits the following two famous presentations, the first of which is more evocative and the second is finite:
\[ F = \langle X_0, X_1, \ldots \mid X_i^{-1} X_j X_i = X_{n+1} \forall i \le j \rangle, \]
\[ F = \langle X_0, X_1 \mid [X_0 X_1^{-1}, X_0^{-1} X_1 X_0] = [X_0 X_1^{-1}, X_0^{-2} X_1 X_0^2] = 1 \rangle. \]
The elements $X_0$ and $X_1$ from these presentations are represented in \cref{fig.F.generators}.

The same group can be realized as the rearrangement group of the edge replacement system depicted in \cref{fig.interval.replacement}.
The limit space is the unit interval (naturally codified by binary sequences, see \cref{ex.interval}) and, under this identification, the prefix-exchange maps defined by rearrangements coincide with the homeomorphism of Thompson's group $F$ as previously defined.
\cref{fig.F.generators} also depicts $X_0$ and $X_1$ as rearrangements, with both graph pair diagrams and forest pair diagrams.
Note that there is a unique graph isomorphism between each two graph expansions with the same number of edges, so there is no need to specify what the isomorphism is in a graph pair diagram.

Among the many properties of Thompson's group $F$, here we recall some of the most notable.

\medskip %layout
\begin{proposition}
\label{prop.F}
Thompson's group $F$ enjoys the following properties.
\begin{itemize}
    \item $F$ is torsion-free (easy to verify).
    \item $F$ acts transitively on the set of dyadic points of $(0,1)$ \cite[Lemma 4.2]{cfp}.
    \item $F$ is dense in $\mathrm{Homeo}^+([0,1])$ \cite[Corollary A5.8]{PL-homeo}.
    \item $F$ has finiteness type $F_\infty$ \cite[Theorem 4.17]{Brown}.
    \item $F$ has exponential growth \cite[Corollary 4.7]{cfp}.
    \item $F$ is totally ordered \cite[Theorem 4.11]{cfp}.
\end{itemize}
\end{proposition}

The commutator subgroup of $F$ is also of particular interest.

\medskip %layout
\begin{proposition}
\label{prop.FF}
The commutator subgroup $[F,F]$ of Thompson's group $F$ enjoys the following properties:
\begin{itemize}
    \item $[F,F]$ consists exactly of those elements of $F$ that act trivially around $0$ and around $1$ \cite[Theorem 4.1]{cfp}.
    \item $[F,F]$ acts transitively on the set of dyadic points of $(0,1)$ \cite[Lemma 4.2]{cfp}.
    \item $[F,F]$ is simple \cite[Theorem 4.5]{cfp}.
\end{itemize}
\end{proposition}

One of the reasons why Thompson's group $F$ is still very relevant in geometric group theory is that, despite decades of attempts, it is still unknown whether it is amenable or not.
See \cite{amenable} for a recent treatment on this open question.

\subsubsection{Thompson's group \texorpdfstring{$T$}{T}}
\label{sub.T}

Thompson's group $T$ sits intermediately between $F$ and $V$.
Its definition is very similar to that of $F$, but instead of acting on the interval it acts on the circle, which allows cyclic permutations of the dyadic intervals (or, equivalently, of the edges in graph pair diagrams or of the leaves in forest pair diagrams).

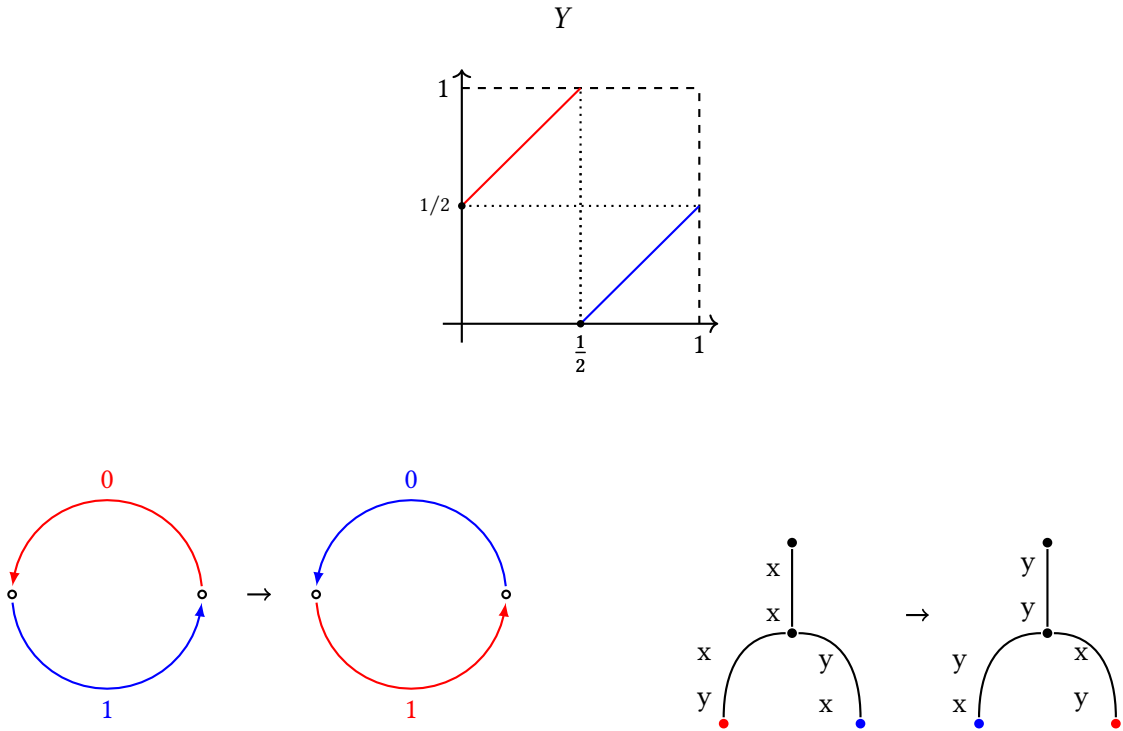
\begin{figure}
\centering
$Y$
\\
\vspace{.5cm}
% Y as a map
\begin{tikzpicture}[scale=1.25]
    \draw[->] (-0.2, 0) -- (2.7, 0);
    \draw[->] (0, -0.2) -- (0, 2.7);
    \draw[dashed] (0,2.5) node[left]{\small$1$} -- (2.5,2.5);
    \draw[dashed] (2.5,0) node[below]{\small$1$} -- (2.5,2.5);
    \draw[red, scale=2.5, domain=0:0.5, smooth, variable=\x] plot({\x}, {\x+.5});
    %\node at (2.5*0.5,2.5*1) [circle,fill,inner sep=1]{};
    \node at (0,2.5*.5) [circle,fill,inner sep=1]{};
    \draw[dotted] (2.5*0.5,0) node[below] {$\frac{1}{2}$} -- (2.5*0.5,2.5*1);
    \draw[blue, scale=2.5, domain=0.5:1, smooth, variable=\x] plot({\x},{\x-.5});
    \node at (2.5*0.5,2.5*0) [circle,fill,inner sep=1]{};
    %\node at (2.5*1,2.5*.5) [circle,fill,inner sep=1]{};
    \draw[dotted] (2.5*0.5,0) node[below] {$\frac{1}{2}$} -- (2.5*0.5,2.5*0.75);
    \draw[dotted] (0,2.5*0.5) node[left] {\scriptsize$1/2$} -- (2.5*1,2.5*0.5);
\end{tikzpicture}
\\
\vspace{1cm}
% Y as a GPD
\begin{tikzpicture}[scale=1.25]
    \begin{scope}[xshift=-1.6cm]
    \draw[edge,red,domain=5:175] plot ({cos(\x)}, {sin(\x)});
    \draw (90:1) node[above,red] {\small$0$};
    \draw (270:1) node[below,blue] {\small$1$};
    \draw[edge,blue,domain=185:355] plot ({cos(\x)}, {sin(\x)});
    \node[vertex] (0) at (0:1) {};
    \node[vertex] (1/2) at (180:1) {};
    \end{scope}
    \draw[-to] (-1/8,0) to (1/8,0);
    \begin{scope}[xshift=1.6cm]
    \draw[edge,blue,domain=5:175] plot ({cos(\x)}, {sin(\x)});
    \draw[edge,red,domain=185:355] plot ({cos(\x)}, {sin(\x)});
    \draw (90:1) node[above,blue] {\small$0$};
    \draw (270:1) node[below,red] {\small$1$};
    \node[vertex] (0) at (0:1) {};
    \node[vertex] (1/2) at (180:1) {};
    \end{scope}
\end{tikzpicture}
\hfill
% Y as a FPD
    \begin{tikzpicture}[font=\small,scale=1.2]
        \node[node] (root) at (0,0) {};
        \node[node] (s) at (0,-1) {};
        \node[node,red] (0) at (-.75,-2) {};
        \node[node,blue] (1) at (.75,-2) {};
        \draw (root) to (s) node[above left,align=center]{x\\x};
        \draw (s) to[out=180,in=90,looseness=1.1] (0) node[above left,align=center]{x\\y};
        \draw (s) to[out=0,in=90,looseness=1.1] (1) node[above left,xshift=-.2cm,align=center]{y\\x};
        \draw[-to] (1.25,-.8) to node[above]{} (1.5,-.8);
        \begin{scope}[xshift=2.8cm]
        \node[node] (root) at (0,0) {};
        \node[node] (s) at (0,-1) {};
        \node[node,blue] (0) at (-.75,-2) {};
        \node[node,red] (1) at (.75,-2) {};
        \draw (root) to (s) node[above left,align=center]{y\\y};
        \draw (s) to[out=180,in=90,looseness=1.1] (0) node[above left,align=center]{y\\x};
        \draw (s) to[out=0,in=90,looseness=1.1] (1) node[above left,xshift=-.2cm,align=center]{x\\y};
        \end{scope}
    \end{tikzpicture}
\caption{The element $Y$ of Thompson's group $T$}
\label{fig.T.generator}
\end{figure}

More explicitly, Thompson's group $T$ can be defined as the group of those orientation-preserving piecewise-linear homeomorphisms of the unit circle $S^1 = [0,1] / \{0,1\}$ with finitely many breakpoints, all of which are dyadic (i.e., they belong to $\mathbb{Z}[\frac{1}{2}]$), and such that the slopes are integer powers of $2$.
The most classical set of generators consists of the $X_0$ and $X_1$ (from Thompson's group $F$, see \cref{fig.F.generators}) along with the element $Y$ depicted \cref{fig.T.generator}.
Differently from Thompson's group $F$, we may now have multiple graph isomorphisms between two graph expansions with the same number of edges, so we use colors in the graph pair diagram to describe the graph isomorphism (the blue edge of the domain graph is mapped to the blue edge of the range graph and so on).

$T$ can be realized as the rearrangement group of the edge replacement system depicted in \cref{fig.circle.replacement}.
The limit space is the unit circle (codified by binary sequences).
\cref{fig.T.generator} also depicts $Y$ as a rearrangement, both with graph pair diagrams and forest pair diagrams.

Here are some notable properties of Thompson's group $T$.

\medskip %layout
\begin{proposition}
\label{prop.T}
Thompson's group $T$ enjoys the following properties:
\begin{itemize}
    \item $T$ contains any finite cyclic group (easy to verify).
    \item The stabilizer of any dyadic point of $S^1$ under the action of $T$ is isomorphic to $F$ (easy to verify).
    \item $T$ acts 2-transitively on the set of dyadic points of $S^1$ (follows from the previous point).
    \item $T$ is dense in $\mathrm{Homeo}^+(S^1)$ \cite[Proposition 4.3]{Tdense}.
    \item $T$ has type $F_\infty$ \cite[Theorem 4.17]{Brown} and in particular it is finitely presented (a small presentation was produced in \cite{Tpresentation1}, which contains a typo; see \cite[Proposition 1.3]{Tpresentation2} for the correct presentation).
    \item $T$ has exponential growth (follows from the second point).
    \item $T$ is simple \cite[Theorem 5.8]{cfp}.
    \item $T$ contains non-abelian free groups (can be showed with a standard ping-pong argument on 4 disjoint intervals).
\end{itemize}
\end{proposition}

Together with $V$, Thompson's group $T$ is the first known example of finitely presented infinite simple group.
As of today, with the exception of the Burger-Mozes groups \cite{BurgerMozes}, most finitely presented infinite simple groups are ``Thompson-like'' in the sense that they arise from generalizations of $T$ and $V$.

\subsubsection{Thompson's group \texorpdfstring{$V$}{V}}
\label{sub.V}

The larger of the three Thompson groups, $V$ does not act on $[0,1]$ by homeomorphisms, but it can still be seen as a group of certain right-continuous piecewise-linear self-bijections of $[0,1)$ that are called \textit{affine interval exchange transformations} (see for example \cite{AIETdistortion}).
In practice, moving from $F$ to $T$ we allowed for cyclic permutations of the dyadic intervals;
moving to $V$, we now allow for \textit{any} permutation of the dyadic intervals.

However, it is easier to forget about $[0,1]$ and $S^1$ and instead define $V$ as a group of homeomorphisms of the Cantor space $\mathfrak{C} = \{0,1\}^\omega$:
Thompson's group $V$ is the group of every prefix-exchange homeomorphism of $\mathfrak{C}$, i.e., any self-bijection $g$ of $\mathfrak{C}$ for which there exist words $u_1, \dots, u_m$ and $v_1, \dots, v_m$ in $\{0,1\}^*$ such that every point of $\mathfrak{C}$ has some $u_i$ and $v_i$ as a prefix and that the action of $g$ is $g(u_i \alpha) = v_i \alpha$ for all $i = 1, \dots, m$.

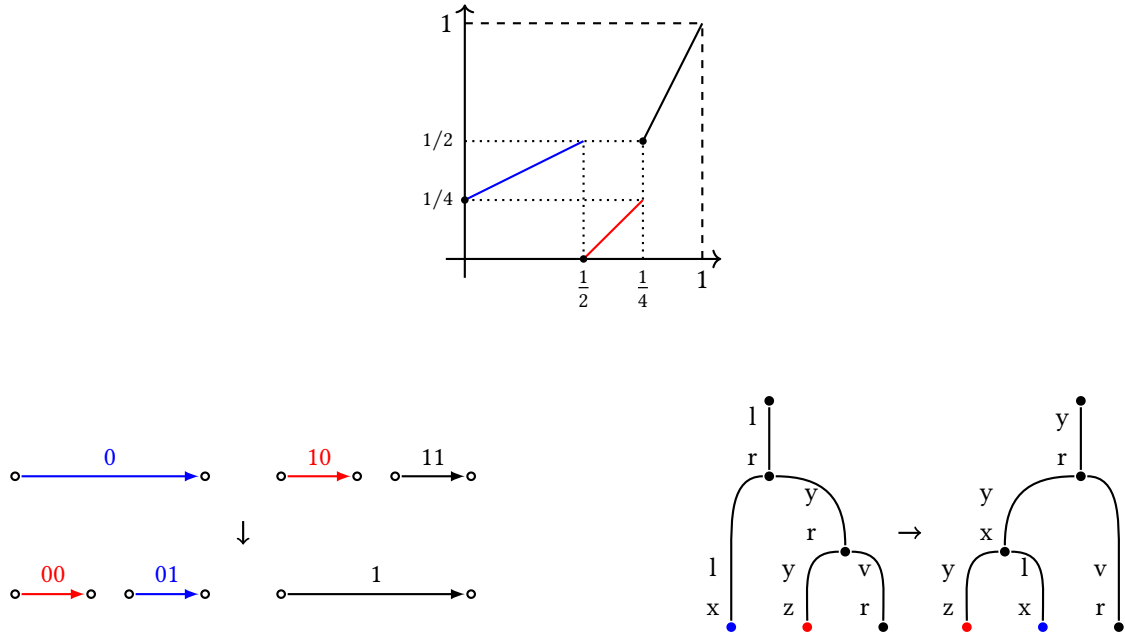
\begin{figure}
\centering
% The element as a map
\begin{tikzpicture}[scale=1.25]
    \draw[->] (-0.2, 0) -- (2.7, 0);
    \draw[->] (0, -0.2) -- (0, 2.7);
    \draw[dashed] (0,2.5) node[left]{\small$1$} -- (2.5,2.5);
    \draw[dashed] (2.5,0) node[below]{\small$1$} -- (2.5,2.5);
    \draw[blue, scale=2.5, domain=0:0.5, smooth, variable=\x] plot({\x}, {.5*\x+.25});
    \node at (2.5*0,2.5*.25) [circle,fill,inner sep=1]{};
    %\node at (2.5*0.5,2.5*.5) [circle,fill,inner sep=1]{};
    \draw[red, scale=2.5, domain=0.5:0.75, smooth, variable=\x] plot({\x},{\x-.5});
    \node at (2.5*0.5,2.5*0) [circle,fill,inner sep=1]{};
    %\node at (2.5*0.75,2.5*.25) [circle,fill,inner sep=1]{};
    \draw[scale=2.5, domain=0.75:1, smooth, variable=\x] plot({\x},{2*\x-1});
    \node at (2.5*0.75,2.5*0.5) [circle,fill,inner sep=1]{};
    %\node at (2.5,2.5) [circle,fill,inner sep=1]{};
    \draw[dotted] (2.5*0.5,0) node[below] {$\frac{1}{2}$} -- (2.5*0.5,2.5*0.5);
    \draw[dotted] (2.5*0.75,0) node[below] {$\frac{1}{4}$} -- (2.5*0.75,2.5*0.5);
    \draw[dotted] (0,2.5*0.5) node[left] {\scriptsize$1/2$} -- (2.5*0.75,2.5*0.5);
    \draw[dotted] (0,2.5*0.25) node[left] {\scriptsize$1/4$} -- (2.5*0.75,2.5*0.25);
\end{tikzpicture}
\\
\vspace{1cm}
% The element as a GPD
\begin{tikzpicture}[scale=1.25]
    \begin{scope}[yshift=.6cm]
    \node[vertex] (0l) at (0,0) {};
    \node[vertex] (0r) at (2,0) {};
    \node[vertex] (10l) at (2.8,0) {};
    \node[vertex] (10r) at (3.6,0) {};
    \node[vertex] (11l) at (4,0) {};
    \node[vertex] (11r) at (4.8,0) {};
    \draw[edge,blue] (0l) to node[above]{\footnotesize$0$} (0r);
    \draw[edge,red] (10l) to node[above]{\footnotesize$10$} (10r);
    \draw[edge,black] (11l) to node[above]{\footnotesize$11$} (11r);
    \end{scope}
    \draw[-to] (2.4,1/8) to (2.4,-1/8);
    \begin{scope}[yshift=-.65cm]
    \node[vertex] (00l) at (0,0) {};
    \node[vertex] (00r) at (.8,0) {};
    \node[vertex] (01l) at (1.2,0) {};
    \node[vertex] (01r) at (2,0) {};
    \node[vertex] (1l) at (2.8,0) {};
    \node[vertex] (1r) at (4.8,0) {};
    \draw[edge,red] (00l) to node[above]{\footnotesize$0$0} (00r);
    \draw[edge,blue] (01l) to node[above]{\footnotesize$01$} (01r);
    \draw[edge,black] (1l) to node[above]{\footnotesize$1$} (1r);
    \end{scope}
    \node[node,white] at (0,-1) {};
\end{tikzpicture}
\hfill
% The element as a FPD
    \begin{tikzpicture}[font=\footnotesize,scale=1]
        \node[node] (root) at (0,0) {};
        \node[node] (s) at (0,-1) {};
        \node[node,blue] (0) at (-.5,-3) {};
        \node[node] (1) at (1,-2) {};
            \node[node,red] (10) at (.5,-3) {};
            \node[node] (11) at (1.5,-3) {};
        \draw (root) to (s) node[above left,align=center]{l\\r};
        \draw (s) to[out=180,in=90,looseness=1] ($(0)+(0,1)$) to (0) node[above left,align=center]{l\\x};
        \draw (s) to[out=0,in=90,looseness=1.25] (1) node[above left,xshift=-.2cm,align=center]{y\\r};
            \draw (1) to[out=180,in=90,looseness=1.25] (10) node[above left,align=center]{y\\z};
            \draw (1) to[out=0,in=90,looseness=1.25] (11) node[above left,align=center]{v\\r};
        \draw[-to] (1.7,-1.75) to node[above]{} (2,-1.75);
        \begin{scope}[xshift=4.1cm]
        \node[node] (root) at (0,0) {};
        \node[node] (s) at (0,-1) {};
        \node[node] (0) at (-1,-2) {};
            \node[node,red] (00) at (-1.5,-3) {};
            \node[node,blue] (01) at (-.5,-3) {};
        \node[node] (1) at (.5,-3) {};
        \draw (root) to (s) node[above left,align=center]{y\\r};
        \draw (s) to[out=180,in=90,looseness=1.25] (0) node[above left,align=center]{y\\x};
            \draw (0) to[out=180,in=90,looseness=1.25] (00) node[above left,align=center]{y\\z};
            \draw (0) to[out=0,in=90,looseness=1.25] (01) node[above left,align=center]{l\\x};
        \draw (s) to[out=0,in=90,looseness=1] ($(1)+(0,1)$) to (1) node[above left,align=center]{v\\r};
        \end{scope}
    \end{tikzpicture}
\caption{An element of Thompson's group $V$}
\label{fig.V.element}
\end{figure}

Thompson's group $V$ can be realized as the rearrangement group of the edge replacement system depicted in \cref{fig.Cantor.replacement}, whose limit space is the standard binary Cantor space $\mathfrak{C} = \{0,1\}^\omega$.
Its action on $\mathfrak{C}$ is precisely the prefix-exchange action described in the previous paragraph.
\cref{fig.V.element} depicts an element of $V$ as an example, both as a bijection of $[0,1]$ and as a rearrangement, with a graph pair diagram and a forest pair diagram.
As we did with Thompson's group $T$, we use colors in the graph pair diagrams to describe the graph isomorphism between the domain and range graphs.

\medskip %layout
\begin{proposition}
\label{prop.V}
Thompson's group $V$ enjoys the following properties:
\begin{itemize}
    \item $V$ contains any finite group (easy to verify).
    \item $V$ contains copies of Thompson groups $F$ and $T$ (easy to verify).
    \item for every $n \in \mathbb{N}$, $V$ acts $n$-transitively on the sets $\{w\overline{0} \mid w \in \{0,1\}^*\}$ and $\{w\overline{1} \mid w \in \{0,1\}^*\}$ (easy to verify).
    \item $V$ is dense in $\mathrm{Homeo}(\mathfrak{C})$, see for example \cref{rmk:V:dns}.
    \item $V$ has type $F_\infty$ \cite[Theorem 4.17]{Brown} and in particular it is finitely presented (see \cite{Vpresentation} for an explicit and small presentation).
    \item $V$ has exponential growth (follows from the second point).
    \item $V$ is simple \cite[Theorem 6.9]{cfp}.
    \item $V$ contains non-abelian free groups (follows from the second point).
\end{itemize}
\end{proposition}

Aside from being, together with $T$, the first known infinite finitely presented simple group, the importance of Thompson's group $V$ is also highlighted by the current state of the Lehnert's conjecture about co-context free groups.
A group is \textit{co-context free} if the set of words in one of its generating sets that are not the identity is a context-free language.
For the sake of conciseness, we will not delve deeper in this definition, but we refer an interested reader to \cite{coCF}.
A group $G$ is \textit{universal co-context free} if it is co-context-free and contains an embedded copy of every co-context-free group.
In his dissertation \cite{QV}, J. Lehnert conjectured that a certain Thompson-like group $QV$ (also known as $\mathrm{QAut}(\mathcal{T}_{2,c})$ and which happens to be a rearrangement group, see \cref{sub.thompson.like}) is universal co-context free.
Thanks to advancements by Bleak, Matucci and Neunh\"{o}ffer \cite{LehnertConjecture}, Lenhert's conjecture has been shown to be equivalent to Thompson's group $V$ being universal co-context free.

Thompson's group $V$ can also be defined as the topological full group of the full shift on the alphabet $\{0,1\}$ (\cref{ex.full.shift}).
These groups will be introduced later in \cref{sub.topological.full.groups}.

\subsubsection{The Higman-Thompson Groups}
\label{sub.higman.thompson.groups}

A natural generalization of the main trio of Thompson groups is the family of so-called Higman-Thompson groups.
Commonly referred to as $F_{n,r}$, $T_{n,r}$ and $V_{n,r}$ with $n \geq 1, r \geq 2$ (but the notation often varies in the literature), and those that we called $V_{n,r}$ were introduced in \cite{higman1974finitely} in the context of automorphism groups of certain algebras.
More generally, the Higman-Thompson groups can be defined as groups of homeomorphisms of an interval, a circle and a Cantor space in the same way in which $F$, $T$ and $V$ are.
This definition can be easily given in terms of rearrangement groups of (monochromatic) edge replacement systems built as follows.
\begin{itemize}
    \item For $F_{n,r}$, the base graph is a path of $r$ edges and the replacement graph is a path of $n$ edges, where $\iota$ and $\tau$ are the initial and terminal vertices of the path, respectively.
    \item For $T_{n,r}$, the base graph is a cycle of $r$ edges and the replacement graph is a path of $n$ edges, where $\iota$ and $\tau$ are the initial and terminal vertices of the path, respectively.
    \item For $V_{n,r}$, the base graph consists of $r$ disjoint edges and the replacement graph consists of $n$ disjoint edges, where $\iota$ and $\tau$ are the initial and terminal vertices of two distinct edges, respectively.
\end{itemize}
Note that setting $n=2$ and $r=1$ produces the edge replacement systems for Thompson groups $F$, $T$ and $V$ depicted in \cref{fig.Thompson.replacement}, i.e., $F=F_{2,1}$, $T=T_{2,1}$ and $V=V_{2,1}$.
Moreover, observe that the limit spaces of the $F_{n,r}$'s are intervals, those of the $T_{n,r}$'s are circles and those of the $V_{n,r}$'s are Cantor spaces.

Many of the properties of $F$, $T$ and $V$ generalize to the Higman-Thompson groups $F_{n,r}$, $T_{n,r}$ and $V_{n,r}$, respectively.
For instance, the previously mentioned \cite[Theorem 4.17]{Brown} is not only about the three main Thompson groups, but also states that the Higman-Thompson groups are all of type $F_\infty$.
Determining which of these groups (or which of their subgroups) are simple is more complicated:
Higman showed that $V_{n,r}$ is simple when $n$ is even and has an index-2 simple normal subgroup when $n$ is odd \cite[Section 5]{higman1974finitely};
Brown studied simple subgroups of $F_{n,r}$ and $T_{n,r}$ \cite[Theorems 4.13, 4.15]{Brown}.

The $V$-like Higman-Thompson groups $V_{n,r}$ can also be defined as the topological full groups of certain full shifts.
These groups will be introduced later in \cref{sub.topological.full.groups}.

\subsection{Rearrangements of Notable Fractals}
\label{sub.fractal.rearrangement.groups}

One of the main appeals of the notion of rearrangement groups is that they have manageable prefix-exchange actions on complicated fractal spaces.
For the main three Thompson groups and the Higman-Thompson groups (and also the groups described in the next \cref{sub.topological.full.groups}), these spaces are simply an interval, a circle or a Cantor space.
Many groups that turn out to be rearrangement groups (such as those discussed in \cref{sub.thompson.like,,sub.Houghton}) do not act on particularly interesting spaces.
This subsection is instead about rearrangement groups of limit spaces that are homeomorphic to notable fractals, which motivated the introduction of these groups in the first place.

\subsubsection{The airplane}
\label{sub.airplane}

The airplane edge replacement system was our guiding example throughout the introductory part of this section and it was depicted in \cref{fig.airplane.replacement}.
It was first mentioned as an example of non-monochromatic edge replacement system in \cite{BF19}, which is instead mostly concerned with monochromatic ones, and its rearrangement group was later studied by the author of this dissertation in \cite{airplane}, adapted from the author's master thesis \cite{airplanethesis}.
We already mentioned that the airplane limit space, portrayed in \cref{fig.airplane}, is homeomorphic to the boundary of the Julia set for the complex map $f(z) = z^2 + c$, where $c$ is approximately $-1.755$.

Denote by $T_A$ the rearrangement group of the airplane edge replacement system.
The following statement is a collection of most known results about $T_A$.

\medskip %layout
\begin{proposition}
\label{prop.TA}
The airplane rearrangement group $T_A$ enjoys the following properties:
\begin{itemize}
    \item $T_A$ is generated by a copy of Thompson's group $F$ together with a copy of Thompson's group $T$ and thus it is finitely generated \cite[Theorem 6.2]{airplane}.
    \item $T_A$ embeds into Thompson's group $T$ \cite[Theorem 8.1]{airplane}.
    \item The abelianization $T_A / [T_A,T_A]$ is $\mathbb{Z}$ \cite[Corollary 7.4]{airplane}.
    \item The commutator subgroup $[T_A,T_A]$ is simple \cite[Theorem 7.8]{airplane}.
    \item The commutator subgroup $[T_A,T_A]$ is finitely generated \cite[Theorem 7.12]{airplane}.
    \item $T_A$ and its commutator subgroup $[T_A,T_A]$ act 2-transitively on the set of interior regions of the airplane limit space \cite[Corollary 9.4]{airplane}.
    \item $T_A$ contains copies of the basilica rearrangement group $T_B$ discussed below \cite[Theorem 10.1]{airplane}.
\end{itemize}
\end{proposition}

\subsubsection{The basilica}
\label{sub.basilica}

The basilica rearrangement group $T_B$ was the first new rearrangement group to be studied.
It was defined and studied in \cite{Belk_2015}, before edge replacement systems and limit spaces were introduced by the same authors in \cite{BF19}.
For this reason, the original treatment of $T_B$ frames it as a group of piecewise-linear homeomorphisms of the circle that preserve the lamination of the basilica Julia set.

\begin{figure}
\centering
\begin{tikzpicture}
    \useasboundingbox (-2.25,-.9) rectangle (8.5,2);
    %\draw[help lines,step=.1cm] (-2.25,-.9) grid (8.5,2);
    %
    \node at (0,1.6) {$X_0$};
    \node[vertex] (s) at (0,0) {};
    \draw[edge] (s) to[loop,out=135,in=225,min distance=2.75cm,looseness=10] node[above left]{$L$} (s);
    \draw[edge] (s) to[loop,out=-45,in=45,min distance=2.75cm,looseness=10] node[above right]{$R$} (s);
    \begin{scope}[xshift=6.5cm,yshift=-.5cm]
    \node at (0,2.1) {$X_1$};
    \node[vertex] (l) at (-1.5,0) {};
    \draw (l) node[above]{$\iota$};
    \node[vertex] (c) at (0,0) {};
    \node[vertex] (r) at (1.5,0) {};
    \draw (r) node[above]{$\tau$};
    \draw[edge] (l) to node[above]{$0$} (c);
    \draw[edge] (c) to[loop,out=130,in=50,min distance=1.8cm,looseness=10] node[above]{$1$} (c);
    \draw[edge] (c) to node[above]{$2$} (r);
    \end{scope}
\end{tikzpicture}
\caption{The basilica edge replacement system.}
\label{fig.basilica.replacement}
\end{figure}
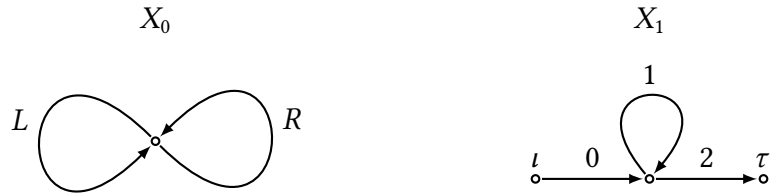

\begin{figure}
\centering
\includegraphics[width=.45\textwidth]{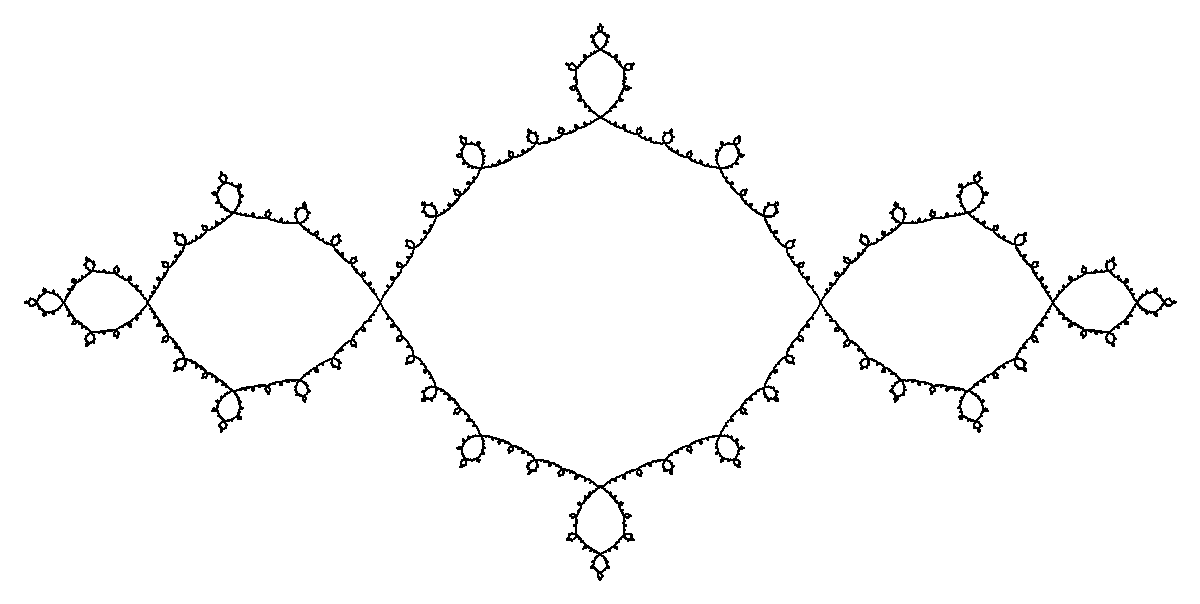}
\caption{The basilica limit space.}
\label{fig.basilica}
\credits{J. Belk and B. Forrest, \cite[Figure 2(b)]{Belk_2015}}
\end{figure}

It can be equivalently defined as the rearrangement group of the edge replacement system depicted in \cref{fig.basilica.replacement} (the author's master thesis \cite{airplanethesis} features a ``translation'' of the arguments of \cite{Belk_2015} in the language of edge replacement systems).
Its limit space, which is portrayed in \cref{fig.basilica}, is a homeomorphic copy of the boundary of the basilica Julia set, which is the Julia set for the complex map $z \mapsto z^2-1$.
The following statement is a collection of most known results about $T_B$.

\medskip %layout
\begin{proposition}
\label{prop.TB}
The basilica rearrangement group $T_B$ enjoys the following properties:
\begin{itemize}
    \item $T_B$ is generated by a copy of Thompson's group $T$ together with a copy of $\mathbb{Z}$ and thus it is finitely generated \cite[Theorem 7.1]{Belk_2015}.
    \item $T_B$ is not finitely presented \cite{WZ19}.
    \item $T_B$ embeds into Thompson's group $T$ (follows, for example, from \cref{prop.TA}).
    \item The abelianization $T_B / [T_B,T_B]$ is $\mathbb{Z} / 2\mathbb{Z}$ \cite[Corollary 8.2]{Belk_2015}.
    \item $T_B$ acts transitively on the set of interior regions of the basilica limit space \cite[Corollary 7.2]{Belk_2015} whereas its commutator subgroup $[T_B,T_B]$ has two orbits \cite[Theorem 8.1]{Belk_2015}.
    \item The commutator subgroup $[T_B,T_B]$ is simple \cite[Theorem 8.4]{Belk_2015}.
    \item The commutator subgroup $[T_B,T_B]$ is finitely generated \cite[Theorem 8.1]{Belk_2015}.
\end{itemize}
\end{proposition}

\subsubsection{Other Fractals}
\label{sub.other.fractals}

Many other edge replacement systems have been considered or briefly mentioned in the literature.
Here we describe three of them.

\begin{figure}
\centering
\begin{tikzpicture}
    \useasboundingbox (-2.2,-1.333) rectangle (8,2);
    %\draw[help lines,step=.1cm] (-2.2,-1.333) grid (8,2);
    %
    \node at (0,1.6) {$X_0$};
    \node[vertex] (s) at (0,0) {};
    \draw[edge] (s) to[loop,out=15,in=85,min distance=2cm,looseness=10] node[right]{$A$} (s);
    \draw[edge] (s) to[loop,out=135,in=205,min distance=2cm,looseness=10] node[above left]{$B$} (s);
    \draw[edge] (s) to[loop,out=255,in=325,min distance=2cm,looseness=10] node[right=.25cm]{$C$} (s);
    \begin{scope}[xshift=6cm,yshift=-.5cm]
    \node at (0,2.1) {$X_1$};
    \node[vertex] (l) at (-1.5,0) {};
    \draw (l) node[above]{$\iota$};
    \node[vertex] (c) at (0,0) {};
    \node[vertex] (r) at (1.5,0) {};
    \draw (r) node[above]{$\tau$};
    \draw[edge] (l) to node[above left]{$0$} (c);
    \draw[edge] (c) to[loop,out=150,in=100,min distance=1.8cm,looseness=10] node[above]{$1$} (c);
    \draw[edge] (c) to[loop,out=80,in=30,min distance=1.8cm,looseness=10] node[above]{$2$} (c);
    \draw[edge] (c) to node[above right]{$3$} (r);
    \end{scope}
\end{tikzpicture}
\caption{The Douady rabbit edge replacement systems.}
\label{fig.rabbit.replacement}
\end{figure}
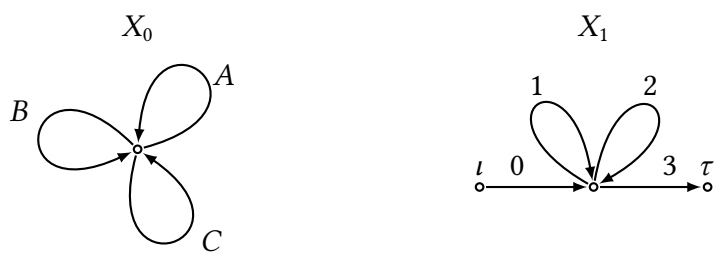

\begin{figure}
\centering
\includegraphics[width=.45\textwidth]{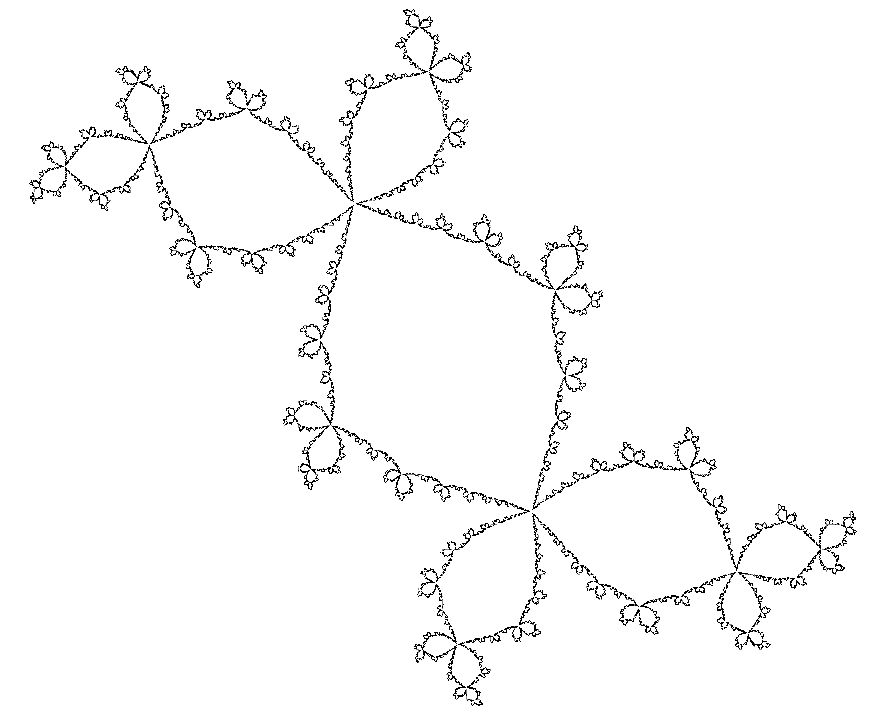}
\caption{The Douady rabbit limit space.}
\label{fig.Douady.rabbit}
\credits{Prokofiev, \href{https://commons.wikimedia.org/wiki/File:Douady_rabbit.png}{wikimedia} (licensed under \href{https://creativecommons.org/licenses/by/3.0}{CC-BY-SA-3.0})}
\end{figure}

\cref{fig.rabbit.replacement} depicts an edge replacement system whose limit space, portrayed in \cref{fig.Douady.rabbit}, is a homeomorphic copy of the boundary of the Douady rabbit, which is the Julia set for the complex map $z \mapsto z^2+c$ where $c$ is approximately $0.123 - 0.745i$.
This can be seen as a generalization of the basilica in the sense that, wherever the basilica edge replacement system has two loops, the Douady rabbit edge replacement system has one more based on the same vertex.
Further generalizations can be constructed by including additional loops, which allows us to produce homeomorphic copies of the boundaries of Julia sets such as the one for $z \mapsto z^2+c$ for $c \approx 0.28 + 0.53i$.

\begin{figure}
\centering
\begin{tikzpicture}
    \node at (0,1.6) {$X_0$};
    \node[vertex] (t) at (0,.8) {};
    \node[vertex] (b) at (0,-.8) {};
    \draw[edge] (b) to[out=180,in=180,looseness=2] node[left]{$l$} (t);
    \draw[edge] (b) to[out=90,in=270] node[left]{$c$} (t);
    \draw[edge] (b) to[out=0,in=0,looseness=2] node[left]{$r$} (t);
    \begin{scope}[xshift=5.5cm]
    \node at (0,1.6) {$X_1$};
    \node[vertex] (l) at (-1.75,0) {}; \draw (-1.75,0) node[above]{$\iota$};
    \node[vertex] (cl) at (-.5,0) {};
    \node[vertex] (cr) at (.5,0) {};
    \node[vertex] (r) at (1.75,0) {}; \draw (1.75,0) node[above]{$\tau$};
    \draw[edge] (cl) to node[above]{1} (l);
    \draw[edge] (cr) to node[above]{4} (r);
    \draw[edge] (cr) to[out=90,in=90,looseness=1.4] node[above]{2} (cl);
    \draw[edge] (cl) to[out=270,in=270,looseness=1.4] node[above]{3} (cr);
    \end{scope}
\end{tikzpicture}
\caption{The bubble bath edge replacement systems.}
\label{fig.bubble.bath.replacement}
\end{figure}
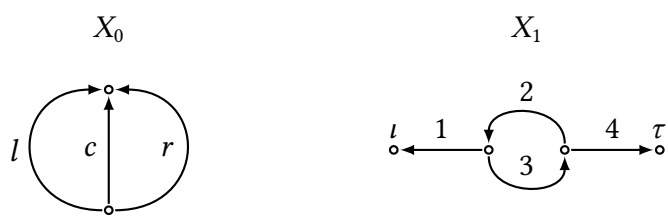

\begin{figure}
\centering
\includegraphics[width=.45\textwidth]{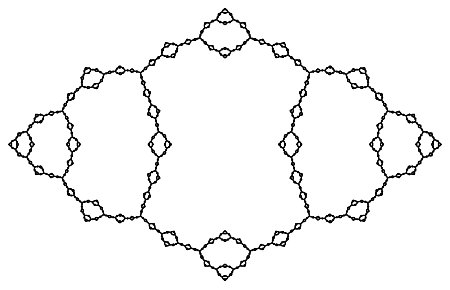}
\caption{The bubble bath limit space.}
\label{fig.bubble.bath}
\credits{J. Belk and B. Forrest, \cite[Figure 19(A)]{BF19}}
\end{figure}

\cref{fig.bubble.bath.replacement} depicts an edge replacement system whose limit space, portrayed in \cref{fig.bubble.bath}, is a homeomorphic copy of the so-called bubble bath Julia set, which is the Julia set for the complex map $z \mapsto \frac{1-z^2}{z^2}$.
A group that is similar to the rearrangement group of this edge replacement system was studied in the senior project \cite{BubbleBath}.

\cref{cha.dendrites} of this dissertation is about the family of dendrite rearrangement groups, which are the rearrangement groups for edge replacement systems depicted schematically in \cref{fig.dendrite.replacement}.
Their limit spaces are Wa\.zewski dendrites, such as the one portrayed in \cref{fig.dendrite.3}, which is also homeomorphic to the Julia set for the complex map $z \mapsto z^2-i$.
We postpone further discussions on dendrites to \cref{cha.dendrites}.

\subsection{Topological Full Groups of Edge Shifts}
\label{sub.topological.full.groups}

Recall that Thompson's group $V$ is the group of those homeomorphisms of the binary Cantor space $\{0,1\}^\omega$ that act by a finite permutation of prefixes (see \cref{sub.V}).
The $V$-like Higman-Thompson groups $V_{n,r}$ (\cref{sub.higman.thompson.groups}) are the same, except that the binary Cantor space is replaced by a different (yet homeomorphic, recall \hyperref[thm.brouwer]{Brouwer's Theorem}) copy of the Cantor space: $\{0,1,\dots,r-1\} \times \{0,1,\dots,n-1\}^\omega$.
This sort of groups can be defined more generally for edge shifts, as defined below.
They are commonly referred to as \textit{topological full groups of edge shifts} (or \textit{of (one-sided) subshifts of finite type}).

These groups were introduced in \cite{Matsumoto} and studied extensively in \cite{Matui}.
We will however follow the introduction to these groups given in the Nuclear Type Systems subsection of \cite[Section 6]{types}, as its notation is much more similar to that of this dissertation.
Most sources assume that the edge shift is \textit{irreducible}, meaning that the graph generating the edge shift (as in \cref{def.edge.shift}) is strongly connected (i.e., given any vertices $v$ and $w$, there is a walk from $v$ to $w$).
Note that, in particular, irreducibility implies that the edge shift is either a singleton (when the graph is a cycle) or a Cantor space by \cref{prop.shift.is.Cantor} (because a strongly connected graph cannot feature an inescapable cycle unless it is itself a cycle).
The topological full groups can be defined without this assumption, as we will do below, but most known results about them assume that the edge shift is a Cantor space.

Consider an edge shift $\Omega(\Gamma,v_0)$.
If $x$ and $y$ are elements of the language of the edge shift (\cref{def.alphabet.language}) that end with the same symbol, then there is a \textbf{prefix-exchange homeomorphism} between the cones of $x$ and $y$ (\cref{def.cones}), which is the following:
\[ C(x) \to C(y), \; x \alpha \mapsto y \alpha. \]
The \textbf{topological full group of the edge shift} $\Omega(\Gamma,v_0)$ is the group $V(\Gamma,v_0)$ of those self-homeomorphisms $\phi \colon \Omega(\Gamma,v_0) \to \Omega(\Gamma,v_0)$ such that, for every $\alpha \in \Omega(\Gamma,v_0)$, there exists a prefix $x$ of $\alpha$ and a prefix $y$ of $\phi(\alpha)$ such that the restriction $\phi|_{C(x)} \colon C(x) \to C(y)$ is a prefix-exchange homeomorphism.

Since $\Omega(\Gamma,v_0)$ is compact, a self-homeomorphism $\phi \colon \Omega(\Gamma,v_0) \to \Omega(\Gamma,v_0)$ lies in $V(\Gamma,v_0)$ if and only if there exist two partitions $C(x_1), \dots, C(x_m)$ and $C(y_1), \dots, C(y_m)$ of the edge shift into cones such that for every $i \in \{1, \dots, m\}$ the restriction $\phi|_{C(x_i)}$ is a prefix-exchange homeomorphism from $C(x_i)$ to $C(y_i)$.

Matui has proven the following results about topological full groups of irreducible edge shifts.

\medskip %layout
\begin{theorem}[\cite{Matui}]
\label{thm.Matui}
Suppose that the edge shift $\Omega(\Gamma,v_0)$ is a Cantor space and consider its topological full group $V(\Gamma,v_0)$.
\begin{enumerate}
    \item $V(\Gamma,v_0)$ has type $F_\infty$ and so, in particular, it is finitely presented (\cite[Theorem 6.21]{Matui}).
    \item The commutator subgroup $[V(\Gamma,v_0),V(\Gamma,v_0)]$ is simple and finitely generated (\cite[Theorem 4.16 and Corollary 6.25]{Matui}).
    \item The abelianization of $V(\Gamma,v_0)$ has been computed (\cite[Corollary 6.24]{Matui}).
\end{enumerate}
\end{theorem}

\subsubsection{Edge replacement systems for Topological Full Groups of Edge Shifts}

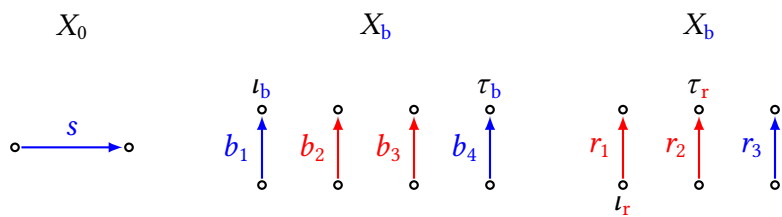
\begin{figure}
\centering
\begin{tikzpicture}
    \node at (0,1.6) {$X_0$};
    \node[vertex] (l) at (-.75,0) {};
    \node[vertex] (r) at (.75,0) {};
    \draw[edge,blue] (l) to node[above]{$s$} (r);
    \begin{scope}[xshift=4cm]
    \node at (0,1.6) {$X_{\text{\textcolor{blue}{b}}}$};
    \node[vertex] (1b) at (-1.5,-.5) {};
    \node[vertex] (1t) at (-1.5,.5) {};
    \node[vertex] (2b) at (-.5,-.5) {};
    \node[vertex] (2t) at (-.5,.5) {};
    \node[vertex] (3b) at (.5,-.5) {};
    \node[vertex] (3t) at (.5,.5) {};
    \node[vertex] (4b) at (1.5,-.5) {};
    \node[vertex] (4t) at (1.5,.5) {};
    \draw[edge,blue] (1b) to node[left]{$b_1$} (1t);
    \draw[edge,red] (2b) to node[left]{$b_2$} (2t);
    \draw[edge,red] (3b) to node[left]{$b_3$} (3t);
    \draw[edge,blue] (4b) to node[left]{$b_4$} (4t);
    \draw (-1.5,.5) node[above]{$\iota_{\text{\textcolor{blue}{b}}}$};
    \draw (1.5,.5) node[above]{$\tau_{\text{\textcolor{blue}{b}}}$};
    \end{scope}
    \begin{scope}[xshift=8.25cm]
    \node at (0,1.6) {$X_{\text{\textcolor{blue}{b}}}$};
    \node[vertex] (1b) at (-1,-.5) {};
    \node[vertex] (1t) at (-1,.5) {};
    \node[vertex] (2b) at (0,-.5) {};
    \node[vertex] (2t) at (0,.5) {};
    \node[vertex] (3b) at (1,-.5) {};
    \node[vertex] (3t) at (1,.5) {};
    \draw[edge,red] (1b) to node[left]{$r_1$} (1t);
    \draw[edge,red] (2b) to node[left]{$r_2$} (2t);
    \draw[edge,blue] (3b) to node[left]{$r_3$} (3t);
    \draw (-1,-.5) node[below]{$\iota_{\text{\textcolor{red}{r}}}$};
    \draw (0,.5) node[above]{$\tau_{\text{\textcolor{red}{r}}}$};
    \end{scope}
\end{tikzpicture}
\caption{The edge replacement system for the topological full group of the edge shift on the airplane color graph (\cref{fig.airplane.color.graph}).}
\label{fig.top.full.group.airplane}
\end{figure}

As first noted in \cite{conjugacy}, given an edge shift $\Omega(\Gamma,v_0)$, its topological full group $V(\Gamma,v_0)$ is a rearrangement group for the edge replacement system having $\Gamma$ as its color graph (except for the starting vertex $q(0)$, see \cref{def.color.graph}), $\Omega(\Gamma,v_0)$ as its underlying symbol space and a trivial gluing relation (i.e., $\alpha \sim \beta \iff \alpha = \beta$).
More precisely, this edge replacement system can be built as follows.
\begin{itemize}
    \item The set of colors is the set of vertices of $\Gamma$.
    \item The base graph is a single edge colored by $v_0$.
    \item The replacement graph for the $v$-colored edge is a disjoint union of edges, one for each edge of $\Gamma$ originating from $v$.
    If such an edge of $\Gamma$ terminates at $w$, the corresponding edge of the replacement graph is colored by $w$.
\end{itemize}
For example, \cref{fig.top.full.group.airplane} depicts the edge replacement system for $V(\Gamma,v_0)$ where $\Gamma$ is the color graph shown in \cref{fig.airplane.color.graph} and $v_0=q(0)$.

Since all edges are disjoint in every graph expansion, the gluing relation $\sim$ is trivial, so each cell $\llbracket x \rrbracket$ is a cone $C(x)$ and the limit space $\Omega(\Gamma,v_0) / \sim$ is the edge shift $\Omega(\Gamma,v_0)$ itself.
Now, the graph isomorphisms between graph expansions are precisely the color-preserving bijections between partitions of $\Omega(\Gamma,v_0)$ into cones, and every such permutation can be realized.
Since the canonical homeomorphisms of the limit space defined in \cref{def.canonical.homeomorphism} are precisely prefix-exchange homeomorphisms in this case, the rearrangement group really is $V(\Gamma,v_0)$.

\medskip %layout
\begin{proposition}
Every topological full group of an edge shift is a rearrangement group.
\end{proposition}

Topological full groups of edge shifts are the ``largest'' rearrangement groups one can build, as illustrated in \cref{prop.embedding.into.topological.full.groups} below.

\subsubsection{Rearrangement Groups Embed into Topological Full Groups of Edge Shifts}
\label{sub.embedding.into.TFGoES}

From its definition (\cref{def.rearrangement}), one can see that each rearrangement group consists of certain prefix-exchange homeomorphisms of an edge shift (\cref{def.edge.shift}) that are also isomorphisms of graph expansions.
The group of those homeomorphisms of an edge shift that act as prefix-exchange homeomorphisms on finitely many cones is precisely the topological full group of the edge shift, so the following fact is not surprising.

\medskip %layout
\begin{proposition}
\label{prop.embedding.into.topological.full.groups}
Given an edge replacement system $\mathcal{R}$ and its symbol space $\Omega_\mathcal{R}$, the corresponding rearrangement group $G_\mathcal{R}$ embeds into the topological full group of $\Omega_\mathcal{R}$.
\end{proposition}

\begin{proof}
Consider the edge replacement system $\mathcal{R}^*$ obtained from $\mathcal{R}$ by ``ungluing'' all the edges of both the base and replacement graphs, so that each graph $X_i^*$ of $\mathcal{R}^*$ is a set of pairwise disjoint edges, one for each edge of $X_i$ (where $X_i$ and $X_i^*$ are the graphs of $\mathcal{R}$ and $\mathcal{R}^*$, respectively).
For example, this produces the edge replacement system for the Higman-Thompson group $V_{n,r}$ (in particular, Thompson's group $V$) from those for $F_{n,r}$ or $T_{n,r}$ ($F$ and $T$).
As another example, starting from the airplane edge replacement system (\cref{fig.airplane.replacement}) one obtains the edge replacement system shown in \cref{fig.top.full.group.airplane}.

The edge replacement systems obtained in this way correspond to those constructed in \cref{sub.topological.full.groups}, so their rearrangement groups are the topological full groups of edge shifts.
The edge shift is $\Omega_\mathcal{R}$, since $\mathcal{R}$ and $\mathcal{R}^*$ share the same symbol space.

Finally, it is easy to see that the rearrangement group of $\mathcal{R}^*$ embeds in that of $\mathcal{R}$.
Indeed, the rearrangements of $\mathcal{R}^*$ are all the prefix-exchange homeomorphisms of $\Omega_\mathcal{R}^*$, whereas the rearrangements of $\mathcal{R}$ are those prefix-exchange homeomorphisms of $\Omega_\mathcal{R}$ that preserve the edge adjacency of the graph expansions.
Since $\Omega_\mathcal{R}^* = \Omega_\mathcal{R}$, each rearrangement of $\mathcal{R}$ is clearly a rearrangement of $\mathcal{R}^*$, as needed.
\end{proof}

This result is useful in \cref{sec.embedding.into.V}, where we show that every rearrangement group embeds in Thompson's group $V$.

\subsection{Thompson-like Groups \texorpdfstring{$QF$, $QT$, $QV$}{QF, QT, QV}}
\label{sub.thompson.like}

In \cref{sub.V} we briefly mentioned a Thompson-like group usually denoted by $QV$ or $\mathrm{QAut}(\mathcal{T}_{2,c})$, which was first introduced in Lehnert's dissertation \cite{QV}.
Here we define this group, together with its two smaller siblings $QF$ and $QT$, and we show that they are rearrangement groups.
This was first noted in \cite{conjugacy}.

This group is defined as follows.
Let $\mathcal{T}_{2,c}$ be the infinite rooted binary tree together with a coloring of the edges with colors $0$ and $1$ such that, for each vertex of the tree, the ``left'' edge is colored by $0$ and the ``right'' edge is colored by $1$.
The set of vertices of $\mathcal{T}_{2,c}$ can thus be seen as the set of paths (including the ``empty path'') from the root, which clearly corresponds to the set $\{0,1\}^*$ of finite words over the alphabet $\{0,1\}$.
An element of $QV$ is a self-bijection of the set of vertices $\{0,1\}^*$ that respects the edge and color relation except for finitely many vertices.
The groups $QF$ and $QT$ are defined in the same way, with the additional requirement that elements of $QF$ must preserve, respectively, the linear and cyclic lexicographic orders of the ends of $\mathcal{T}_{2,c}$.
Essentially, elements of $QF$, $QT$ and $QV$ are allowed to freely permute the vertices between two finite complete subtrees of $\mathcal{T}_{2,c}$ and must act by prefix-exchange permutations below the subtrees;
this permutation must preserve the linear or cyclic lexicographic order in the case of $QF$ and $QT$.

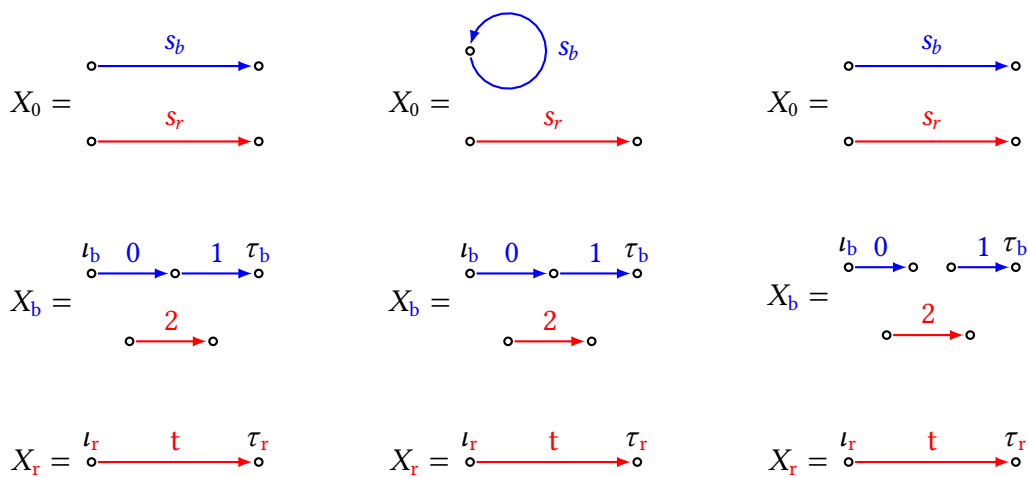
\begin{figure}
\centering
\begin{subfigure}[B]{.32\textwidth}
\centering
\begin{tikzpicture}
    \begin{scope}[yshift=1cm]
    \node at (-.667,-.5) {$X_0 =$};
    \node[vertex] (i) at (0,0) {};
    \node[vertex] (t) at (2.2,0) {};
    \node[vertex] (ri) at (0,-1) {};
    \node[vertex] (rt) at (2.2,-1) {};
    \draw[edge,blue] (i) to node[above]{$s_b$} (t);
    \draw[edge,red] (ri) to node[above]{$s_r$} (rt);
    \end{scope}
    \begin{scope}[yshift=-1.75cm]
    \node at (-.667,-.4) {$X_{\text{\textcolor{blue}{b}}} =$};
    \node[vertex] (l) at (0,0) {};
    \draw (l) node[above]{$\iota_{\text{\textcolor{blue}{b}}}$};
    \node[vertex] (c) at (1.1,0) {};
    \node[vertex] (r) at (2.2,0) {};
    \draw (r) node[above]{$\tau_{\text{\textcolor{blue}{b}}}$};
    \draw[edge,blue] (l) to node[above]{$0$} (c);
    \draw[edge,blue] (c) to node[above]{$1$} (r);
    \node[vertex] (lb) at (.5,-.9) {};
    \node[vertex] (rb) at (1.6,-.9) {};
    \draw[edge,red] (lb) to node[above]{2} (rb);
    \end{scope}
    \begin{scope}[yshift=-4.25cm]
    \node at (-.667,0) {$X_{\text{\textcolor{red}{r}}} =$};
    \node[vertex] (l) at (0,0) {};
    \draw (l) node[above]{$\iota_{\text{\textcolor{red}{r}}}$};
    \node[vertex] (r) at (2.2,0) {};
    \draw (r) node[above]{$\tau_{\text{\textcolor{red}{r}}}$};
    \draw[edge,red] (l) to node[above]{t} (r);
    \end{scope}
\end{tikzpicture}
\caption{The edge replacement system for $QF$.}
\end{subfigure}
\hfill
\begin{subfigure}[B]{.325\textwidth}
\centering
\begin{tikzpicture}
    \begin{scope}[yshift=1cm]
    \node at (-.667,-.5) {$X_0 =$};
    \draw[edge,blue,domain=190:530] plot ({.5*cos(\x)+.5}, {.5*sin(\x)+.2});
    \draw (1,.2) node[right,blue] {$s_b$};
    \node[vertex] (i) at (0,.2) {};
    \node[vertex] (ri) at (0,-1) {};
    \node[vertex] (rt) at (2.2,-1) {};
    \draw[edge,red] (ri) to node[above]{$s_r$} (rt);
    \end{scope}
    \begin{scope}[yshift=-1.75cm]
    \node at (-.667,-.4) {$X_{\text{\textcolor{blue}{b}}} =$};
    \node[vertex] (l) at (0,0) {};
    \draw (l) node[above]{$\iota_{\text{\textcolor{blue}{b}}}$};
    \node[vertex] (c) at (1.1,0) {};
    \node[vertex] (r) at (2.2,0) {};
    \draw (r) node[above]{$\tau_{\text{\textcolor{blue}{b}}}$};
    \draw[edge,blue] (l) to node[above]{$0$} (c);
    \draw[edge,blue] (c) to node[above]{$1$} (r);
    \node[vertex] (lb) at (.5,-.9) {};
    \node[vertex] (rb) at (1.6,-.9) {};
    \draw[edge,red] (lb) to node[above]{2} (rb);
    \end{scope}
    \begin{scope}[yshift=-4.25cm]
    \node at (-.667,0) {$X_{\text{\textcolor{red}{r}}} =$};
    \node[vertex] (l) at (0,0) {};
    \draw (l) node[above]{$\iota_{\text{\textcolor{red}{r}}}$};
    \node[vertex] (r) at (2.2,0) {};
    \draw (r) node[above]{$\tau_{\text{\textcolor{red}{r}}}$};
    \draw[edge,red] (l) to node[above]{t} (r);
    \end{scope}
\end{tikzpicture}
\caption{The edge replacement system for $QT$.}
\end{subfigure}
\hfill
\begin{subfigure}[B]{.32\textwidth}
\centering
\begin{tikzpicture}
    \begin{scope}[yshift=1cm]
    \node at (-.667,-.5) {$X_0 =$};
    \node[vertex] (i) at (0,0) {};
    \node[vertex] (t) at (2.2,0) {};
    \node[vertex] (ri) at (0,-1) {};
    \node[vertex] (rt) at (2.2,-1) {};
    \draw[edge,blue] (i) to node[above]{$s_b$} (t);
    \draw[edge,red] (ri) to node[above]{$s_r$} (rt);
    \end{scope}
    \begin{scope}[yshift=-1.667cm]
    \node at (-.667,-.4) {$X_{\text{\textcolor{blue}{b}}} =$};
    \node[vertex] (l) at (0,0) {};
    \draw (l) node[above]{$\iota_{\text{\textcolor{blue}{b}}}$};
    \node[vertex] (c1) at (.85,0) {};
    \node[vertex] (c2) at (1.35,0) {};
    \node[vertex] (r) at (2.2,0) {};
    \draw (r) node[above]{$\tau_{\text{\textcolor{blue}{b}}}$};
    \draw[edge,blue] (l) to node[above]{$0$} (c1);
    \draw[edge,blue] (c2) to node[above]{$1$} (r);
    \node[vertex] (lb) at (.5,-.9) {};
    \node[vertex] (rb) at (1.6,-.9) {};
    \draw[edge,red] (lb) to node[above]{2} (rb);
    \end{scope}
    \begin{scope}[yshift=-4.25cm]
    \node at (-.667,0) {$X_{\text{\textcolor{red}{r}}} =$};
    \node[vertex] (l) at (0,0) {};
    \draw (l) node[above]{$\iota_{\text{\textcolor{red}{r}}}$};
    \node[vertex] (r) at (2.2,0) {};
    \draw (r) node[above]{$\tau_{\text{\textcolor{red}{r}}}$};
    \draw[edge,red] (l) to node[above]{t} (r);
    \end{scope}
\end{tikzpicture}
\caption{The edge replacement system for $QV$.}
\end{subfigure}
\caption{The edge replacement systems for the Thompson-like groups $QF$, $QT$ and $QV$.}
\label{fig.Q.replacement}
\end{figure}

Consider the edge replacement systems depicted in \cref{fig.Q.replacement}.
A graph expansion consists of some blue path, cycle or set of disjoint edges together with a red edge for each edge expansion that is performed.
This corresponds to a finite complete subtree of $\mathcal{T}_{2,c}$ in the sense that each red edge $wt^k$ ($k \geq 0$ and $w$ any word ending with either $s_r$ or $2$) corresponds to an internal vertex $w$ of the subtree and the blue edges correspond to its leaves.
Note that a rearrangement of one of these three edge replacement systems consists of any permutation of the red edges together with the trivial, a cyclic or any permutation of the blue edges, respectively.
With this in mind, it is clear that the rearrangement groups of the edge replacement systems are exactly $QF$, $QT$ and $QV$.

Observe that the edge replacement systems depicted in \cref{fig.Q.replacement} are not expanding (\cref{def.expanding}).
However, it is easy to see that their gluing relation is an equivalence relation, and the limit spaces are, respectively, $[0,1] \cup N$, $S^1 \cup N$ and $\mathfrak{C} \cup N$ (disjoint unions), where $N$ is the countable set $\{0,1\}^*$ equipped with the discrete topology.
Using \cref{prop.null.expanding.isolated.rarrangement} one could also produce an expanding edge replacement system for these groups by replacing the red replacement graphs from \cref{fig.Q.replacement} with the replacement graph of \cref{fig.replacement.rearrangementless}.

\subsection{The Houghton Groups}
\label{sub.Houghton}

For each $n \in \mathbb{N}$, the $n$-th Houghton group $H_n$ is the group of those permutations of $X_n = \{1, \dots, n\} \times \mathbb{N}$ that, except for finitely many points, act as translations on each of the $n$ copies of $\mathbb{N}$.
First introduced by Houghton in \cite{Houghton1978TheFC}, Brown later showed that each $H_n$ is of type $F_{n-1}$ but not $F_n$.
Here we show how they can be realized as rearrangement groups, which was first noted in \cite{conjugacy}.

A finite subset $S$ of $X_n$ is \textit{complete} if, whenever $(i,n)$ belongs to $S$, for each $j \le i$ the element $(j,n)$ also belongs to $S$.
A \textit{co-complete} subset of $X_n$ is $(\{i\} \times \mathbb{N}) \setminus S$ for some finite complete subset $S$ of $X_n$.
Similarly to how in the previous subsection red edges represented the vertices of the tree $\mathcal{T}_{2,c}$ and blue edges together represented the complement of complete subtrees, here black edges will represent the points of $X_n$ and edges of different colors together will represent the co-complete subsets.

\begin{figure}
\centering
\begin{tikzpicture}
    \node at (-2.5,0) {$X_0 =$};
    \node[vertex] (1i) at (0:.6) {};
    \node[vertex] (1o) at (0:1.6) {};
    \node[vertex] (2i) at (45:.6) {};
    \node[vertex] (2o) at (45:1.6) {};
    \node[vertex] (3i) at (90:.6) {};
    \node[vertex] (3o) at (90:1.6) {};
    \node[vertex] (4i) at (135:.6) {};
    \node[vertex] (4o) at (135:1.6) {};
    \node[vertex] (5i) at (180:.6) {};
    \node[vertex] (5o) at (180:1.6) {};
    \draw[edge,blue] (1i) to node[above]{\small$s_1$} (1o);
    \draw[edge,Orange] (2i) to node[above left]{\small$s_2$} (2o);
    \draw[edge,Green] (3i) to node[left]{\small$s_3$} (3o);
    \draw[edge,red] (4i) to node[below left]{\small$s_4$} (4o);
    \draw[edge,Turquoise] (5i) to node[below]{\small$s_5$} (5o);
    \draw[dotted,gray] (-40:.6) to (-40:1.6);
    \draw[dotted,gray] (220:.6) to (220:1.6);
    \node[gray,align=center] at (270:1.1) {total\\of $n$ edges};
    \begin{scope}[xshift=4.6cm,yshift=.8cm]
    \node at (-1.1,0) {$X_{\text{\textcolor{Plum}{color}}} =$};
    \node at (4.3,0) {$\forall \text{\textcolor{Plum}{color}} \neq \text{black}$};
    \node[vertex] (ci) at (0,0) {};
    \node[vertex] (c1) at (1,0) {};
    \node[vertex] (c2) at (1.5,0) {};
    \node[vertex] (ct) at (2.5,0) {};
    \draw (ci) node[above]{$\iota$};
    \draw (ct) node[above]{$\tau$};
    \draw[edge,black] (ci) to node[above]{\small$b$} (c1);
    \draw[edge,Plum] (c2) to node[above]{\scriptsize{color}} (ct);
    \end{scope}
    \begin{scope}[xshift=4.6cm,yshift=-.8cm]
    \node at (-1.1,0) {$X_{\text{black}} =$};
    \node[vertex] (bi) at (0,0) {};
    \node[vertex] (bt) at (2.5,0) {};
    \draw (bi) node[above]{$\iota$};
    \draw (bt) node[above]{$\tau$};
    \draw[edge,black] (bi) to node[above]{$t$} (bt);
    \end{scope}
\end{tikzpicture}
\caption{Edge replacement systems for the Houghton groups $H_n$. There are $n+1$ distinguished colors, one of which is black and the others have similar replacement graphs.}
\label{fig.houghton.replacement}
\end{figure}
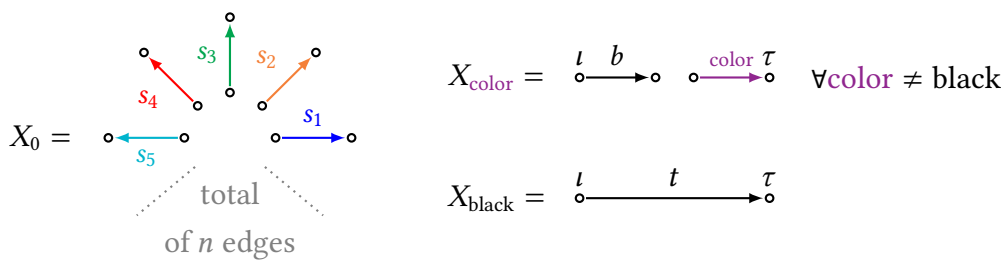

Consider the edge replacement systems schematically depicted in \cref{fig.houghton.replacement} and fix a natural number $n$.
We denote by $\{ 1, \dots, n, \text{black} \}$ the set of colors.
There is a correspondence between graph expansions of this edge replacement system and finite complete subsets of $X_n$ given by mapping each black edge $s_i^k b t^r$ to $(i, k)$ (the edges of other colors represent the $n$ co-complete subsets determined by this finite complete set).
Note that a rearrangement simply consists of any permutation between two finite prefix-closed sets of black edges.
Such sets correspond to complete subsets of $X_n$, so it is clear that rearrangements correspond to the elements of the Houghton group $H_n$.

The edge replacement systems depicted in \cref{fig.houghton.replacement} are not expanding, but the gluing relation is an equivalence relation.
The limit space is thus defined and it is homeomorphic to $\{1, \dots, n\} \times (\{ \frac{k-1}{k} \}_{k \in \mathbb{N}} \cup \{1\})$.
Using \cref{prop.null.expanding.isolated.rarrangement} one could also produce an expanding edge replacement system for the Houghton groups by replacing the black replacement graph from \cref{fig.houghton.replacement} with the replacement graph of \cref{fig.replacement.rearrangementless}.

\section{Other Known Results}
\label{sub.other.results}

In this section we briefly mention three additional facts about rearrangement groups, both of which are due to Belk and Forrest \cite{BF19}.
We will not give many details about them, but we refer an interested reader to the original source for a complete treatment.

\subsection{Finite Subgroups}

Subsection 2.3 of \cite{BF19} is about finite subgroups of rearrangement groups.
It features a proof of the following useful fact and an example of its applications.

\medskip %layout
\begin{theorem}[Theorem 2.9 \cite{BF19}]
\label{thm.finite.subgroups}
Every finite subgroup of a rearrangement group is a subgroup of the automorphism group of some graph expansion of the edge replacement system.
\end{theorem}

This will be useful to prove \cref{thm:non:isomorphic}, where we will show that the rearrangement groups of dendrites of different orders are all distinct because they are distinguished by their finite subgroups.

\subsection{Action on a CAT(0) cubical complex}

Every rearrangement group acts properly by isometries on a locally finite CAT(0) cubical complex.

These complexes are generalizations of those built by Farley for Thompson groups $F$, $T$ and $V$ and, more generally, for diagram and picture groups \cite{FarleyDiagramGroups,FarleyPictureGroups}.
The construction is based on a groupoid of generalized rearrangements that are based on the data given by the edge replacement rules but not on the base graph;
this very groupoid will be discussed in \cref{cha.conjugacy} of this dissertation, as it will be essential for describing a technique to tackle the conjugacy problem of rearrangement groups.
The 1-skeleton of the complex is built from a poset on certain equivalence classes of generalized rearrangements.

A fully detailed construction is given in \cite[Section 3]{BF19}

\subsection{A Condition for Finiteness Properties}

Section 4 of \cite{BF19} makes use of the action of rearrangement groups on the aforementioned CAT(0) cubical complexes to obtain a sufficient condition for the group to be of type $F_\infty$.
Their techniques generalize the classical argument applied on Thompson groups $F$, $T$ and $V$, which combines Bestvina-Brady Morse theory \cite{BestvinaBrady} and Brown's criterion \cite{Brown} to deduce that a group enjoys a finiteness property $F_n$ if it acts properly by isometries on a contractible cube complex with certain properties on sublevels cut by Morse functions.
Belk and Forrest use these ideas to show the following result.

\medskip %layout
\begin{theorem}[Theorem 4.1 \cite{BF19}]
\label{thm.finiteness}
Let $\mathcal{R}$ be an edge replacement system with finite branching (\cref{def.finite.branching}).
Let $\Gamma(\mathcal{R})$ be collection of all graphs that are obtained by finite sequences of edge expansions and reductions from the base graph of $\mathcal{R}$.
Assume that, for every $m \geq 1$, all but finitely many graphs of $\Gamma(\mathcal{R})$ admit at least $m$ distinct reductions.
Then the rearrangement group of $\mathcal{R}$ is of type $F_\infty$.
\end{theorem}

This immediately show that a family of rearrangement groups defined in \cite[Example 1.12]{BF19}, which are called \textit{Vicsek rearrangement groups}, are of type $F_\infty$.
Even if the Vicsek rearrangement groups act on dendrites, it looks like these are not the same groups as the dendrite rearrangement groups discussed in \cref{cha.dendrites}.

Unfortunately, the hypothesis of \cref{thm.finiteness} does not apply to the airplane, basilica, dendrite nor bubble bath rearrangement groups (all of which were introduced or mentioned in \cref{sub.fractal.rearrangement.groups}).
Figure 35 of \cite{BF19} describes why that is the case with the basilica edge replacement system, and a similar issue affects the other edge replacement systems mentioned above.

In fact, the basilica rearrangement group has been proved to not even be of finitely presented \cite{WZ19}.
The airplane rearrangement group is instead of type $F_\infty$ (it is stated but not proved in \cite[Example 2.13]{BF19}) and \cite{dendrite} asks whether dendrite rearrangement groups are of type $F_\infty$.
An upcoming work by Davide Perego and the author of this dissertation concerning a larger class of groups will tackle this question.

%%%%%%%%%%%%%%%%%%%%%%%%%

\chapter{Rearrangement Groups of Dendrites}
\label{cha.dendrites}

Among the many fractals that can be built as limit spaces of edge replacement systems, dendrites (defined later in \cref{def:dendrite}) have already inspired multiple works from group theorists, arguably because of their tree-like structure (the sense in which they are tree-like can actually be formalized: see Definition 2.12 and Theorem 10.32 of \cite{continua}).
For example, \cite{UniversalDendrites} is about the group of homeomorphisms of universal dendrites (i.e., dendrites that contain a homeomorphic copy of every dendrite), \cite{AmenableDendrites} is about actions of amenable countable groups on dendrites and \cite{DM18} is about the dynamics of group actions on dendrites.
The recent work \cite{DM19}, which is going to be a fundamental reference for this chapter, provides a rich study of the group of homeomorphisms of a dendrite.
Moreover, \cite{Kaleid} uses Wa\.zewski dendrites to build simple permutation groups with interesting dynamical and topological properties, and \cite{Duc20} provides a deep study of the group of homeomorphisms of the infinite-order Wa\.zewski dendrite $D_\infty$.

\begin{figure}
\centering
\includegraphics[width=.45\textwidth]{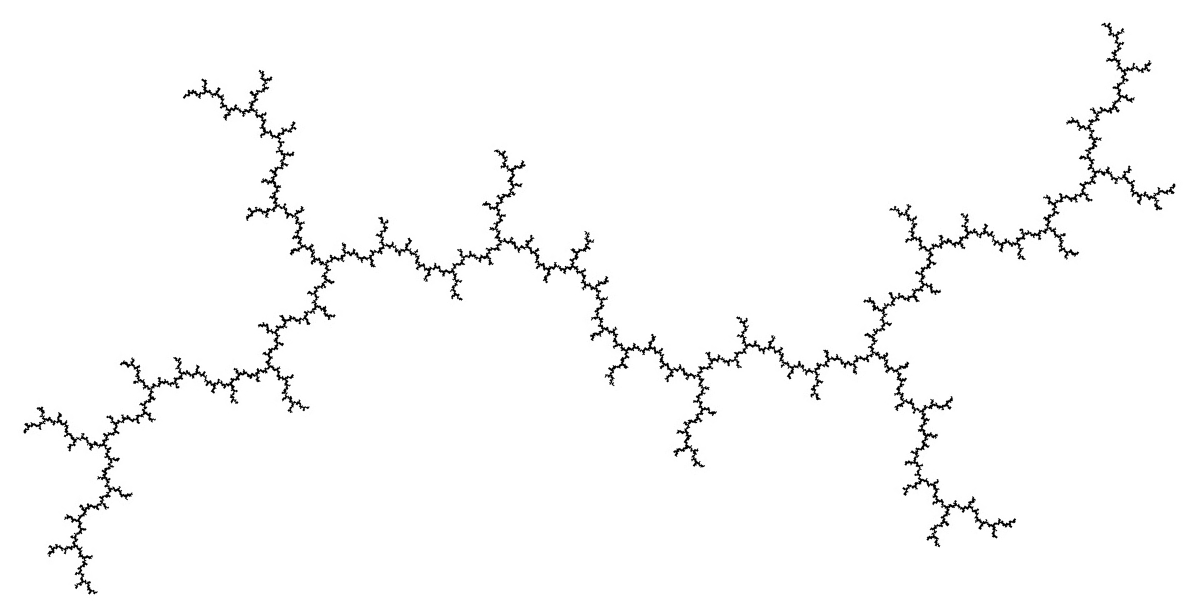}
\caption{The Wa\.zewski dendrite $D_3$.}
\label{fig.dendrite.3}
\credits{Adam Majewski, \href{https://commons.wikimedia.org/wiki/File:Julia_IIM_1.jpg}{wikimedia} (licensed under \href{https://creativecommons.org/licenses/by/3.0}{CC-BY-3.0})}
\end{figure}

This chapter, which presents the content of the work \cite{dendrite}, belongs to the intersection of these two fields of study in group theory, that of Thompson-like groups and that of groups acting on dendrites.
We will make use of strategies from both worlds, without assuming that the reader is familiar with any of them.
We will consider Wa\.zewski dendrites $D_n$ ($n \geq 3$, \cref{fig.dendrite.3} depicts $D_3$), which can be thought as ``dense $n$-regular trees'', and we will study rearrangement groups of Thompson-like homeomorphisms of them.
After breifly providing formal definitions and the necessary background on Wa\.zewski dendrites in \cref{sec.dendrites.background}, in \cref{sec:dendrite:rearrangement:groups} we will exhibit the edge replacement systems $\mathcal{D}_n$ that realize Wa\.zewski dendrites as limit spaces, and their rearrangement groups $G_n$ are the subjects of the subsequent \cref{sec:gen,,sec:dns,,sec:comm,,sec:properties}, where we prove the results collected in the following statement:

\medskip %layout
\begin{maintheorem*}
For each $n \geq 3$, the dendrite rearrangement group $G_n$ is a countable group satisfying the following statements.
\begin{description}
    \item[\cref{thm:gen}:] $G_n$ is finitely generated by a copy of Thompson's group $F$ and a copy of the permutation group $\mathrm{Sym}(n)$.
    \item[\cref{thm:dense}:] $G_n$ is dense in the full group of homeomorphisms of the Wa\.zewski dendrite $D_n$.
    \item[\cref{thm:commutator}:] The abelianization $G_n / [G_n, G_n]$ is isomorphic to $\mathbb{Z}_2 \oplus \mathbb{Z}$.
    \item[\cref{thm:comm:simple}:] The commutator subgroup $[G_n, G_n]$ of $G_n$ is simple if $n \geq 4$.
    \item[\cref{thm:comm:fg}:] The commutator subgroup $[G_n,G_n]$ is finitely generated.
    \item[\cref{thm:non:isomorphic}:] $G_n$ and $G_m$ are isomorphic if and only if $n = m$.
    \item[\cref{sec:properties}:] $G_n$ embeds into $V$, but it does not embed into $T$ nor does $T$ embed into $G_n$.
    Moreover, $G_n$ is not invariably generated, its conjugacy problem is solvable, it contains non-abelian free subgroups and the cardinality of its finite subgroups has been classified.
\end{description}
\end{maintheorem*}

It is not yet known whether or not the commutator subgroup of $G_3$ is simple, as will be discussed in \cref{rmk:3:comm:simple}.

Finally, \cref{sec:generalizations} explores the relation between dendrite rearrangement groups and the airplane rearrangement group $T_A$ studied in \cite{airplane}, proving that $T_A$ is dense in the group of all orientation-preserving homeomorphisms of the Wa\.zewski dendrite $D_\infty$, and hints at possible generalizations of these groups, such as the rearrangement groups of generalized Wa\.zewski dendrites $D_S$ for finite subsets $S \subset \mathbb{N}_{\geq3}$.

We believe that results about density of rearrangement groups in their ``ambient'' overgroup (which are not uncommon, as discussed in \cref{sub:dns}) could produce interesting connections with computer science.
These results could provide a tool for ``approximating'' uncountable homeomorphism groups of fractals by Thompson-like groups, which, being countable and showing nice behaviour with regards to decision problems (such as often having solvable conjugacy problem as shown in \cref{cha.conjugacy}), are arguably quite suitable for computer processing.

Overall, dendrite rearrangement groups seem to share similar behaviors to both Thompson's group $F$ (they are generated by a copy of $F$ along with a finite group; they inherit many transitive properties from those of $F$ and the commutator subgroup is partly characterized by the behaviour at the endpoints of the space on which it acts) and Thompson's group $V$ (such as featuring plenty of torsion elements, not being invariably generated and including copies of non-abelian free subgroups), while resembling neither $F$ nor $V$.

\section{Background on Dendrites}
\label{sec.dendrites.background}

This short section contains the necessary backgrounds about dendrites and Wa\.zewski dendrites.
Chapter X of \cite{continua} gives further general details on these spaces.

We begin with the general definition of dendrites, then we introduce a concept of order of a point which will be useful to define what Wa\.zewski dendrites are.

\begin{definition}
\label{def:dendrite}
A \textbf{dendrite} is a non-empty, locally connected, compact, connected metric space that contains no simple closed curves.
A dendrite is \textbf{degenerate} if it is a singleton.
A \textbf{subdendrite} is
a dendrite that is the subspace 
of a dendrite.
\end{definition}

There are many equivalent definitions of dendrites.
For example, we could also define a dendrite to be a non-empty compact connected metric space such that the intersection of any two connected subsets is a connected subset (\cite[Theorem 10.10]{continua}).
See also (1.1) and (1.2) in \cite[Section V.1]{AnalyticTopology} for further characterizations of dendrites among continua and locally connected continua.
Non-degenerate dendrites are also the underlying topological spaces behind the so-called \textit{metric trees} (\cite[Proposition 2.2]{MetricTrees}).
Examples of dendrites are the compactification of any simplicial tree (as noted in \cite[Proposition 12.2]{DM18}) and certain Julia sets (for instance, for the complex map $z \to z^2 + i$).

\phantomsection\label{txt.order.of.points}
The \textbf{order} of a point $p$ of a dendrite $X$ is the cardinality of the set of connected components of $X \setminus \{p\}$.
This is equivalent to the more general notion of Menger–Urysohn order, which is discussed in \cite[7, \S 51, p. 274-307]{Top2}; see also \cite[Definition 9.3]{continua}.
A distinction of points of a dendrite based on their order will be given later in \cref{sub:points:arcs}.

For the purpose of this work, we are going to consider the most regular non-trivial dendrites, which are called Wa\.zewski dendrites and are defined right below.
We would like to mention that a general dendrite might be so irregular that its homeomorphism group is trivial, so it is natural that the most regular dendrites would produce larger and often more interesting groups of homeomorphisms.

\begin{definition}
\label{def:wazewski:dendrite}
Given a natural number $n \geq 3$, a \textbf{Wa\.zewski dendrite} of order $n$, denoted by $D_n$, is a dendrite each of whose points is either of order $1$, $2$ or $n$ (no intermediate values), and such that every arc of $D_n$ (which is any subspace of $D_n$ that is homeomorphic to the closed interval $[0,1]$) contains a point of order $n$.
\end{definition}

For example, \cref{fig.dendrite.3} depicts the Wa\.zewski dendrite of order $3$.

By Theorem 6.2 of \cite{selfhomeomorphisms}, $D_n \simeq D_m$ if and only if $n = m$, so we can actually talk about \textit{the} Wa\.zewski dendrite of order $n$.
In the next subsections we will build these dendrites as limit spaces of edge replacement systems.

\section{Rearrangement Groups of Wa\.zewski Dendrites}
\label{sec:dendrite:rearrangement:groups}

In this section we introduce the main topic of this chapter:
the Thompson-like groups $G_n$, which will be defined in \cref{sub:den} as rearrangement groups of edge replacement systems $\mathcal{D}_n$ for the Wa\.zewski dendrites $D_n$.
In \cref{sub:points:arcs} we define branch points, endpoints and (rational) arcs of $D_n$, which will be important in the remainder of this work.
Finally, in \cref{sub:subgroups} we start discussing two key subgroups of the $G_n$'s:
finite permutation groups around a branch point and Thompson groups acting on certain linear portions of the dendrite.

\subsection{Edge replacement systems for Dendrites}
\label{sub:den}

\begin{figure}
\centering
\begin{tikzpicture}
    \node at (-2,0) {$\Gamma =$};
    \node[vertex] (c) at (0,0) {};
    \node[vertex] (1) at (0:1.5) {};
    \node[vertex] (2) at (72:1.5) {};
    \node[vertex] (3) at (144:1.5) {};
    \node[vertex] (n-1) at (216:1.5) {};
    \node[vertex] (n) at (288:1.5) {};
    \draw[edge] (c) to node[above]{$1$} (1);
    \draw[edge] (c) to node[above left]{$2$} (2);
    \draw[edge,dotted] (c) to (3);
    \draw[edge,dotted] (c) to (n-1);
    \draw[edge] (c) to node[right]{$n$} (n);
    \begin{scope}[xshift=7cm]
    \node at (-2.5,0) {$R =$};
    \node[vertex] (b) at (0,0) {};
    \node[vertex] (i) at (180:1.5) {}; \draw (0:-1.5) node[above]{$\iota$};
    \node[vertex] (r2) at (120:1.5) {};
    \node[vertex] (rn-1) at (60:1.5) {};
    \node[vertex] (t) at (0:1.5) {}; \draw (0:1.5) node[above]{$\tau$};
    \draw[edge] (b) to node[below]{$1$} (i);
    \draw[edge,dotted] (b) to (r2);
    \draw[edge,dotted] (b) to (rn-1);
    \draw[edge] (b) to node[below]{$n$} (t);
    \end{scope}
\end{tikzpicture}
\caption{A schematic depiction of the dendrite edge replacement system $\mathcal{D}_n$.}
\label{fig.dendrite.replacement}
\end{figure}
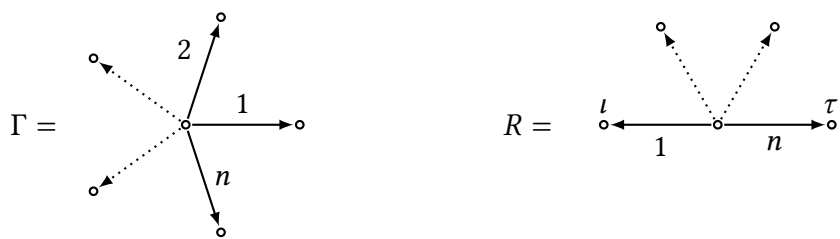

In this subsection we define the edge replacement systems $\mathcal{D}_n$ that produce the Wa\.zewski dendrites $D_n$ as limit spaces, for any natural number $n \geq 3$.
\cref{fig.dendrite.replacement} schematically depicts the generic edge replacement system $\mathcal{D}_n$, which is defined with more precision as follows.
The \textbf{dendrite edge replacement system} $\mathcal{D}_n$ is the monochromatic edge replacement system whose base and replacement graphs are the same tree consisting of a vertex of degree $n$ that is the origin of $n$ edges terminating at $n$ distinguished leaves;
the initial and terminal vertices $\iota$ and $\tau$ of the replacement graph are two distinct leaves.

\phantomsection\label{txt:rep:sys:den}
The edge replacement system $\mathcal{D}_n$ is expanding (\cref{def.expanding}), so the gluing relation is an equivalence relation and the limit space of $\mathcal{D}_n$ is defined as in \cref{def.limit.space}.
Recall that limit spaces are (non-empty) compact metrizable space by \cref{thm.limit.space}.
Moreover, the limit space $D_n$ of $\mathcal{D}_n$ is connected because both the base graph and the replacement graphs are connected, see \cite[Remark 1.27]{BF19}) and it is locally connected (proven in \cref{lem:loc:conn} below).
It contains no simple closed curve, as none of the graph expansions of $\mathcal{D}_n$ contain simple closed paths.
Thus, the limit space of $\mathcal{D}_n$ is a dendrite (\cref{def:dendrite}), and since its points have order $1$, $2$ or $n$, it is homeomorphic to the Wa\.zewski dendrite $D_n$ (\cref{def:wazewski:dendrite}).
We will thus use the same symbol $D_n$ to denote both the Wa\.zewski dendrite and the limit space of $\mathcal{D}_n$.

\begin{remark}
\label{rmk.dendrite.undirected}
Note that the edges of $\mathcal{D}_n$ are \textbf{undirected}, as explained in \cref{def.undirected.edges}.
This means that we can forget about their edge orientation without losing any meaningful information.
In practice, this simply means that we will be able to forget about the orientation of every edge altogether when representing rearrangements, because if an edge $x$ is mapped to an edge $y$ and the orientation is swapped, then what is really happening is $xi \mapsto x \tilde{i}$, where $\tilde{i} = n+1-i$, for all $i \in \{1, \dots, n\}$.
The orientation is still important to correctly describe the gluing relation and codify cells.
\end{remark}

\subsection{Branch Points, Branches, Endpoints and Arcs}
\label{sub:points:arcs}

In this subsection we give a few definition that will be useful throughout this chapter.

\begin{definition}
\label{def:branch}
A \textbf{branch point} of $D_n$ is a point $p \in D_n$ such that $D_n \setminus \{p\}$ has exactly $n$ connected components.
Each of these connected components is called a \textbf{branch} of $D_n$ at the branch point $p$.
\end{definition}

Observe that each branch at $p$ is the topological interior of the union of finitely many cells, one of which has boundary $\{ p \}$.
However, note that not every topological interior of the connected union of finitely many cells is a branch.

We will often refer to the following specific branch point:
\[ \llbracket 11\overline{n} \rrbracket =  \llbracket 21\overline{n} \rrbracket = \dots = \llbracket n1\overline{n} \rrbracket, \]
where the bar on top of $n$ denotes periodicity.
We call this point the \textbf{central branch point} and we denote by $p_0$.
There is nothing topologically special about this branch point, but its nice combinatorial description will make it handy in the remainder of this work.

\begin{definition}
\label{def:ext}
An \textbf{endpoint} of $D_n$ is a point $q \in D_n$ such that $D_n \setminus \{q\}$ has exactly one connected component.
\end{definition}

Given a subset $X$ of $D_n$, we denote by $\mathrm{Br}(X)$ and $\mathrm{En}(X)$ the sets of all branch points and endpoints, respectively, that are contained in $X$.
When $X = D_n$ itself, we simply write $\mathrm{Br}$ and $\mathrm{En}$.

\begin{remark}
\label{rmk:br:en:dns}
A subset of a limit space is dense if and only if it intersects non-trivially every cell (this is shown later in \cref{rmk.density}).
Since every cell contains branch points and endpoints, we can conclude that both $\mathrm{Br}(X)$ and $\mathrm{En}(X)$ are dense in $X$ whenever $X$ is a cell or the entire limit space $D_n$.
It is also useful to note that $\mathrm{Br}$ is countable, whereas $\mathrm{En}$ is uncountable.
\end{remark}

Note that every vertex of a graph expansion of $\mathcal{D}_n$ either corresponds to a branch point or an endpoint, depending on whether its degree is $n$ or $1$.
Moreover, every branch point is a vertex of some graph expansion, but not every endpoint is a vertex, since $\mathrm{En}$ is uncountable.
For concrete examples of endpoints that are not vertices, consider any sequence that does not contain $1$ and that is not eventually repeating.

\phantomsection\label{txt.rationa.endpoints}
We say that an endpoint is a \textbf{rational endpoint} if it is a vertex of some graph expansion.
We denote by $\mathrm{REn}(X)$ the set of those rational endpoints that are contained in $X \subseteq D_n$, and we simply write $\mathrm{REn}$ when $X$ is the entire dendrite $D_n$.
Note that, by the same \cref{rmk.density} mentioned in \cref{rmk:br:en:dns}, $\mathrm{REn}$ is dense in $D_n$ because every cell contains rational endpoints.

Note that this notion of rationality of endpoints does not correspond exactly to the general rationality of points of limit spaces described in \cref{def.rational.irrational.points}, since for example the point $\llbracket \overline{23} \rrbracket$ is a rational point but not a rational endpoint.

Moreover, being a rational endpoint is not a topological property:
the set $\mathrm{REn}$ depends on the chosen homeomorphism between the limit space and the Wa\.zewski dendrite $D_n$, and in general a homeomorphism need not preserve the rationality of an endpoint.
Since we have already fixed a homeomorphism between the limit space and the dendrite to begin with, there will be no need to distinguish between endpoints that are rational under a homeomorphism but not under another.
Differently from generic homeomorphisms, it is easy to see that rearrangements instead act by permutation on $\mathrm{REn}$.

On the other hand, every branch point of $D_n$ is rational, since it is a vertex of some graph expansion.
In this sense, branch points are more useful when dealing with rearrangements, since the set of branch points is invariant under any homeomorphism.

\begin{remark}
\label{rmk:br:ex:pts}
Branch points and rational endpoints have nice and useful combinatorial characterizations based on the sequences that represent them in $\mathcal{D}_n$.
\begin{itemize}
    \item Each branch point of $D_n$ corresponds exactly to the $n$ sequences
    \[ x 1 1 \overline{n}, \: x 2 1 \overline{n}, \: \dots, \: x n 1 \overline{n} \]
    for a unique finite word $x$ (possibly empty).
    In particular, there are natural bijections between the set $\mathrm{Br}$ of branch points of $D_n$, the set of finite words $x$ in the alphabet $\{1, \dots, n\}$ (including the empty word) and the set that consists of the cells $\llbracket x \rrbracket$ of $\mathcal{D}_n$ together with the entire $D_n$ (which can be thought of as the cell associated to the empty word).
    \item Each rational endpoint of $D_n$ correspond uniquely to a sequence
    \[ i \overline{n} \text{ or } x \overline{n} \]
    for a unique $i \in \{1, \dots, n\}$, or for a unique finite non-empty word $x$ in the alphabet $\{1, \dots, n\}$ that does not end with $1$ nor with $n$.
\end{itemize}
\end{remark}

\medskip %layout
\begin{lemma}
\label{lem:loc:conn}
For all $n \geq 3$, the limit space $D_n$ is locally connected.
\end{lemma}

\begin{proof}
Given a point $p \in D_n$, we want to build a basis of connected open sets at $p$.
We distinguish the case in which $p$ corresponds to a vertex and the case in which it does not, and we will use the fact that each cell is connected (because the replacement graph is connected).

If $p$ corresponds to a vertex of some graph expansion, then it is codified by a finite amount of sequences $x \alpha^{(1)}, \dots, x \alpha^{(m)}$ with a common prefix $x$ (possibly empty) and $m$ infinite sequences $\alpha^{(1)}, \dots, \alpha^{(m)}$.
More precisely, by \cref{rmk:br:ex:pts} if $p$ is a branch point then $m=n$ and if $p$ is a rational endpoint then $m=1$.
In both cases, by \cref{cor.basis.of.stars} a basis of connected open sets is given by the cells corresponding to the edges that are adjacent to $p$ in the $k$-th full expansion of $\mathcal{D}_n$ (\cref{def.full.expansion}) in the following way:
\[ \left\{ \{p\} \cup \llparenthesis x \alpha^{(1)}_1 \dots \alpha^{(1)}_k \rrparenthesis \cup \dots \cup \llparenthesis x \alpha^{(m)}_1 \dots \alpha^{(m)}_k \rrparenthesis \mid k \in \mathbb{N} \right\}. \]
Recall that $\llparenthesis w \rrparenthesis$ denotes the topological interior of a cell $\llbracket w \rrbracket$.

Now consider a point $p$ that is not a vertex.
It corresponds to a unique infinite sequence $\alpha = a_1 a_2 \dots$, and must belong to the interior $\llparenthesis a_1 \dots a_k \rrparenthesis$ of every cell $\llbracket a_1 \dots a_k \rrbracket$, so a basis of connected open sets is simply given by
\[ \{ \llparenthesis a_1 \dots a_k \rrparenthesis \mid k \in \mathbb{N} \}. \]
\end{proof}

\begin{definition}
\label{def:arc}
Since dendrites do not contain simple closed curves, any two points are joined by a unique topological arc (i.e., an embedded copy of the closed unit interval), which we simply call an \textbf{arc} of $D_n$.
Given two points $p$ and $q$, we denote by $[p,q]$ the unique arc joining them.
\end{definition}

For dendrites, the existence and uniqueness of arcs joining each pair of distinct points (by \cref{def:wazewski:dendrite}) plays a significant role in their topological properties and in their groups of homeomorphisms
Indeed, this notion will be important for identifying the Thompson subgroups described later in \cref{sub:thomp}, where we will consider the following special types of arcs.

\begin{definition}
\label{def:arc:type}
An arc $A = [p,q]$ is \textbf{rational} if $p$ and $q$ are vertices of some graph expansion of $\mathcal{D}_n$.
A rational arc $[p,q]$ is of one of the following three types:
\begin{itemize}
    \item it is a \textbf{BB-arc} (\textit{branch point to branch point}) if $p, q \in \mathrm{Br}$;
    \item it is an \textbf{EE-arc} (\textit{rational endpoint to rational endpoint}) if $p, q \in \mathrm{REn}$;
    \item it is an \textbf{BE-arc} (\textit{branch point to rational endpoint}) if $p \in \mathrm{Br}$ and $q \in \mathrm{REn}$ or if $p \in \mathrm{REn}$ and $q \in \mathrm{Br}$.
\end{itemize}
\end{definition}

Similarly to the rationality of an endpoint, the rationality of an arc is not a topological property, meaning that homeomorphisms need not preserve it.
It is true, however, that a rearrangement must map a rational arc to a rational arc and a non-rational arc to a non-rational arc.
In practice, we will only be working with arcs that are rational throughout the rest of this chapter.

We will frequently use the following correspondence between paths in a graph expansion and rational arcs.

\begin{remark}
\label{rmk:paths}
There is a natural surjective correspondence from the set of paths in graph expansions to the set of all rational arcs, which is the following.
If $e = x_1 \dots x_k$ is an edge of some graph expansion, we associate to it the arc $[p,q]$ for $p = \llbracket x_1 \dots x_k 1 \overline{n} \rrbracket$ and $q = \llbracket x_1 \dots x_k \overline{n} \rrbracket$.
This arc is precisely the set of points of $D_n$ that are represented by some sequence $x_1 \dots x_k \alpha$ for $\alpha$ in the alphabet $\{1,n\}$.
Inductively extending this construction gives the desired correspondence from paths of a graph expansions to rational arcs, which is surjective because clearly we can find a path between any two points of $D_n$ that are vertices.
In short, if $P$ is the path consisting of the edges $e_1, \dots, e_l$ (all of which are finite words), then the associated arc is
\[ A = \big\{ \llbracket e_i \alpha \rrbracket \mid i=1, \dots, n \text{ and } \alpha \text{ is a sequence in } \{1,n\} \big\}. \]
This correspondence is not injective, but the preimage of each rational arc contains a unique path with a minimal number of edges, which is obtained from any other edge of the same preimage by performing as many graph reductions as possible.
When restricted to the set of these minimal paths, this correspondence is a bijection.
\end{remark}

\begin{remark}
\label{rmk:arcs:dyadic}
For each EE-arc $A$, there is a canonical homeomorphism between $A$ and the unit interval $[0,1]$ that induces a bijection between the branch points of $A$ and the dyadic points of $(0,1)$.
The homeomorphism is built as described below using the path $P$ from the previous \cref{rmk:paths}.

Given a single edge $e = x_1 \dots x_l$, we consider the homeomorphism
\[ \Phi_e: \big\{ \llbracket e \alpha \rrbracket \mid \alpha \text{ is a sequence in } \{1,n\} \big\} \to \big[ 0, 2^{-l} \big] \]
that is built inductively by letting $e 1$ and $e n$ correspond to $2^{-l} [0,1/2]$ and $2^{-l} [1/2,1]$ in the natural order-preserving fashion, and so forth for further edge expansions.

Given a path $P$ of edges $e_1, \dots, e_k$, we build $\Phi_P$ as the concatenation of the $\Phi_{e_i}$'s, by which we mean that the first $\Phi_{e_1}$ maps to $[0, 2^{-l_1}]$, the second $\Phi_{e_2}$ to the adjacent interval of length $2^{-l_2}$ (which is $[2^{-l_1}, 2^{-l_1}+2^{-l_2}]$) and so forth.
For a generic path $P$, this covers some dyadic interval $[0, z] \subseteq [0,1]$, for some $z \in \mathbb{Z}[1/2]$, and it is easy to see that $\Phi_P (A) = [0,1]$ if and only if $P$ corresponds to an EE-arc.

One may note that this construction is very similar to the usual homeomorphism between $[0,1]$ and the limit space of the edge replacement system for Thompson group $F$ (depicted in \cref{fig.interval.replacement}), which is simply given by identifying each number in $[0,1]$ with its binary expansion.
For dendrite edge replacement systems (\cref{fig.dendrite.replacement}), this is less straightforward due to the fact that the edge $1$ is directed towards the initial vertex.
\end{remark}

\subsection{Key Subgroups of Dendrite Rearrangement Groups}
\label{sub:subgroups}

From here onwards, given an $n \geq 3$, we will denote by $G_n$ the rearrangement group of $\mathcal{D}_n$, and we now begin its study.

In general, even if an edge replacement system produces an appealing limit space, its rearrangement group might be uninteresting or even trivial (see for example \cref{fig.replacement.rearrangementless}, which is \cite[Example 2.5]{BF19}).
Here, however, we prove that every dendrite edge replacement system yields a rearrangement group that contains infinitely many copies of the finite permutation group $\mathrm{Sym}(n)$ and infinitely many copies of Thompson's group $F$ (\cref{sub:thomp}), immediately hinting at the fact that these rearrangement groups are, in fact, quite large.
How large they are will be made clear with \cref{thm:dense}, stating that each dendrite rearrangement group $G_n$ is dense in the full group of homeomorphisms of the Wa\.zewski dendrite $D_n$.

\subsubsection{Permutation Subgroups}
\label{sub:perm}

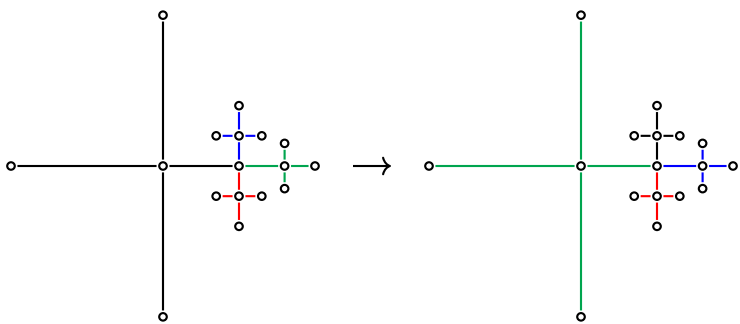
\begin{figure}
\centering
\begin{tikzpicture}
    \node[vertex] (c) at (0,0) {};
    \node[vertex] (l) at (-2,0) {};
    \node[vertex] (t) at (0,2) {};
    \node[vertex] (b) at (0,-2) {};
    \node[vertex] (r) at (2,0) {};
            \node[vertex] (rrc) at (1.6,0) {};
            \node[vertex] (rrt) at (1.6,.3) {};
            \node[vertex] (rrb) at (1.6,-.3) {};
        \node[vertex] (rc) at (1,0) {};
        \node[vertex] (rt) at (1,.8) {};
            \node[vertex] (rtc) at (1,.4) {};
            \node[vertex] (rtl) at (.7,.4) {};
            \node[vertex] (rtr) at (1.3,.4) {};
        \node[vertex] (rb) at (1,-.8) {};
            \node[vertex] (rbc) at (1,-.4) {};
            \node[vertex] (rbl) at (.7,-.4) {};
            \node[vertex] (rbr) at (1.3,-.4) {};
    \draw[black] (c) to (l);
    \draw[black] (c) to (t);
    \draw[black] (c) to (b);
    \draw[black] (rc) to (c);
    \draw[Green] (rrc) to (rc);
    \draw[Green] (rrc) to (rrt);
    \draw[Green] (rrc) to (rrb);
    \draw[Green] (rrc) to (r);
    \draw[blue] (rtc) to (rtl);
    \draw[blue] (rtc) to (rtr);
    \draw[blue] (rtc) to (rc);
    \draw[blue] (rtc) to (rt);
    \draw[red] (rbc) to (rbl);
    \draw[red] (rbc) to (rbr);
    \draw[red] (rbc) to (rc);
    \draw[red] (rbc) to (rb);
    \draw[->] (2.5,0) to (3,0);
    \begin{scope}[xshift=5.5cm]
    \node[vertex] (c) at (0,0) {};
    \node[vertex] (l) at (-2,0) {};
    \node[vertex] (t) at (0,2) {};
    \node[vertex] (b) at (0,-2) {};
    \node[vertex] (r) at (2,0) {};
            \node[vertex] (rrc) at (1.6,0) {};
            \node[vertex] (rrt) at (1.6,.3) {};
            \node[vertex] (rrb) at (1.6,-.3) {};
        \node[vertex] (rc) at (1,0) {};
        \node[vertex] (rt) at (1,.8) {};
            \node[vertex] (rtc) at (1,.4) {};
            \node[vertex] (rtl) at (.7,.4) {};
            \node[vertex] (rtr) at (1.3,.4) {};
        \node[vertex] (rb) at (1,-.8) {};
            \node[vertex] (rbc) at (1,-.4) {};
            \node[vertex] (rbl) at (.7,-.4) {};
            \node[vertex] (rbr) at (1.3,-.4) {};
    \draw[Green] (c) to (l);
    \draw[Green] (c) to (t);
    \draw[Green] (c) to (b);
    \draw[Green] (rc) to (c);
    \draw[blue] (rrc) to (rc);
    \draw[blue] (rrc) to (rrt);
    \draw[blue] (rrc) to (rrb);
    \draw[blue] (rrc) to (r);
    \draw[black] (rtc) to (rtl);
    \draw[black] (rtc) to (rtr);
    \draw[black] (rtc) to (rc);
    \draw[black] (rtc) to (rt);
    \draw[red] (rbc) to (rbl);
    \draw[red] (rbc) to (rbr);
    \draw[red] (rbc) to (rc);
    \draw[red] (rbc) to (rb);
    \end{scope}
\end{tikzpicture}
\caption{An element of the permutation subgroup $K_p$ in $G_4$.}
\label{fig:perm}
\end{figure}

The first fundamental family of subgroups, discussed in more detail in the Lemma below, consists of those groups that fix a branch point $p$ and act by ``rigid'' permutations on the branches at $p$.
We have one such subgroup for each of the countably many branch points $p$ of $D_n$, and an example of an element of one such subgroup is depicted in \cref{fig:perm}.
This depiction is ambiguous, as each of the three external cells of each color could themselves be permuted.
However, we are implying that a natural rotation system, which is defined in the paragraph below, is being preserved at each vertex, possibly except for $p$ itself.

\phantomsection\label{txt:rot:sys}
Recall that a rotation system on a graph is an assignment of a circular order to the edges incident on each vertex (\cref{def.rotation.system}).
By \cref{rmk:br:ex:pts}, each branch point $p$ corresponds to the sequences $xi1\overline{n}$ for all $i=1, \dots, n$.
The $n$ edges that are adjacent to $p$ in a graph expansion correspond to sequences starting with $xi$, and each has a distinct $i$ among $1, \dots, n$.
The circular order at $p$ is defined by the order of these $i$'s.
This rotation system is particularly nice because it is preserved by edge expansions, so it actually induces an orientation on $D_n$.

\medskip %layout
\begin{lemma}
\label{lem:perm}
For every $n \geq 3$ and for each $p \in \mathrm{Br}$, the rearrangement group $G_n$ of $\mathcal{D}_n$ contains an isomorphic copy $K_p$ of the symmetric group $\mathrm{Sym}(n)$ (the group of permutations of $n$ elements) that acts on $D_n$ by permuting the branches at $p$ and preserving the rotation systems at every point with the possible exception of $p$.
\end{lemma}

\begin{proof}
Consider a branch point $p$ and let $E$ be the minimal graph expansion of $\mathcal{D}_n$ where $p$ appears.
Let $\Lambda$ be the subgraph of $E$ corresponding to the connected component of $E \setminus \{p\}$ that contains the central branch point $p_0 = \llbracket i 1 \overline{n} \rrbracket$ (for any choice of $i \in \{1, \dots, n\}$).
The subgraphs corresponding to the other $n-1$ connected components of $E \setminus \{p\}$ can be expanded in such a way that they become isomorphic to $\Lambda$, where the isomorphism fixes $p$ and preserves the rotation system, so we let $E^*$ be the graph expansion of $\mathcal{D}_n$ where the subgraphs corresponding to all of the connected components $\Lambda_i$ of $E^* \setminus \{p\}$ are all isomorphic as trees rooted at $p$ equipped with their rotation systems.
It is easy to see that there is a unique rotation system-preserving rooted tree isomorphism between each of the expanded components.
The group $K_p$ is the subgroup of $G_n$ consisting of those elements whose graph pair diagram has domain and range equal to $E^*$, where the graph isomorphism is a permutation of the $\Lambda_i$'s as trees rooted at $p$ equipped with their rotation systems as discussed above.
These elements correspond one-to-one to the permutations of the $n$ connected components $\Lambda_i$ of $E^* \setminus \{p\}$, and with this identification the group $K_p$ is isomorphic to $\mathrm{Sym}(n)$.
\end{proof}

\subsubsection{Thompson Subgroups}
\label{sub:thomp}

The second fundamental family of subgroups of dendrite rearrangement groups consists of copies of Thompson's group $F$ that act on rational arcs of $D_n$ (which were defined in \cref{def:arc:type}).
Differently from the permutation subgroups $K_p$ described in the previous subsection, here we will start by defining a main copy $H$ of this subgroup and we will later prove in \cref{lem:thomps} that there is an isomorphic copy of the same subgroup acting on every rational arc of $D_n$.

We denote by $A_0$ the EE-arc $[q_1, q_n]$ between endpoints $q_1 = \llbracket 1 \overline{n} \rrbracket$ and $q_n = \llbracket \overline{n} \rrbracket$.
Note that, by \cref{rmk:paths}, $A_0$ corresponds to the set of those points of $D_n$ that are represented by some infinite sequence in the alphabet $\{ 1,n \}$.
The canonical homeomorphism from $A_0$ to $[0,1]$ described in \cref{rmk:arcs:dyadic} prompts the following Lemma.

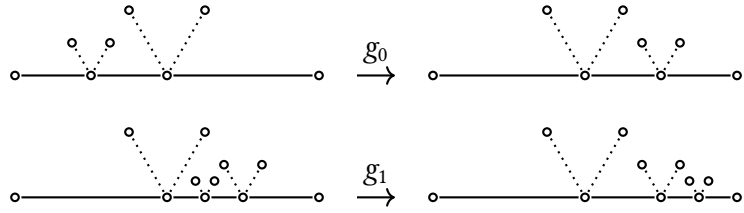
\begin{figure}
\centering
\begin{tikzpicture}
    \node[vertex] (0) at (0,0) {};
    \node[vertex] (1/4) at (1,0) {};
        \node[vertex] (1/4 1) at ($(1/4)+(120:.5)$) {};
        \node[vertex] (1/4 2) at ($(1/4)+(60:.5)$) {};
    \node[vertex] (1/2) at (2,0) {};
        \node[vertex] (1/2 1) at ($(1/2)+(120:1)$) {};
        \node[vertex] (1/2 2) at ($(1/2)+(60:1)$) {};
    \node[vertex] (1) at (4,0) {};
    \draw (0) to (1/4);
    \draw (1/4) to (1/2);
        \draw[dotted] (1/4) to (1/4 1);
        \draw[dotted] (1/4) to (1/4 2);
    \draw (1/2) to (1);
        \draw[dotted] (1/2) to (1/2 1);
        \draw[dotted] (1/2) to (1/2 2);
    \draw[->] (4.5,0) to node[above]{$g_0$} (5,0);
    \begin{scope}[xshift=5.5cm]
    \node[vertex] (0) at (0,0) {};
    \node[vertex] (1/2) at (2,0) {};
        \node[vertex] (1/2 1) at ($(1/2)+(120:1)$) {};
        \node[vertex] (1/2 2) at ($(1/2)+(60:1)$) {};
    \node[vertex] (3/4) at (3,0) {};
        \node[vertex] (3/4 1) at ($(3/4)+(120:.5)$) {};
        \node[vertex] (3/4 2) at ($(3/4)+(60:.5)$) {};
    \node[vertex] (1) at (4,0) {};
    \draw (0) to (1/2);
    \draw (1/2) to (3/4);
        \draw[dotted] (1/2) to (1/2 1);
        \draw[dotted] (1/2) to (1/2 2);
    \draw (3/4) to (1);
        \draw[dotted] (3/4) to (3/4 1);
        \draw[dotted] (3/4) to (3/4 2);
    \end{scope}
\end{tikzpicture}
\\
\vspace{15pt}
\begin{tikzpicture}
    \node[vertex] (0) at (0,0) {};
    \node[vertex] (1/2) at (2,0) {};
        \node[vertex] (1/2 1) at ($(1/2)+(120:1)$) {};
        \node[vertex] (1/2 2) at ($(1/2)+(60:1)$) {};
    \node[vertex] (5/8) at (2.5,0) {};
        \node[vertex] (5/8 1) at ($(5/8)+(120:.25)$) {};
        \node[vertex] (5/8 2) at ($(5/8)+(60:.25)$) {};
    \node[vertex] (3/4) at (3,0) {};
        \node[vertex] (3/4 1) at ($(3/4)+(120:.5)$) {};
        \node[vertex] (3/4 2) at ($(3/4)+(60:.5)$) {};
    \node[vertex] (1) at (4,0) {};
    \draw (0) to (1/2);
    \draw (1/2) to (5/8);
        \draw[dotted] (1/2) to (1/2 1);
        \draw[dotted] (1/2) to (1/2 2);
    \draw (5/8) to (3/4);
        \draw[dotted] (5/8) to (5/8 1);
        \draw[dotted] (5/8) to (5/8 2);
    \draw (3/4) to (1);
        \draw[dotted] (3/4) to (3/4 1);
        \draw[dotted] (3/4) to (3/4 2);
    \draw[->] (4.5,0) to node[above]{$g_1$} (5,0);
    \begin{scope}[xshift=5.5cm]
    \node[vertex] (0) at (0,0) {};
    \node[vertex] (1/2) at (2,0) {};
        \node[vertex] (1/2 1) at ($(1/2)+(120:1)$) {};
        \node[vertex] (1/2 2) at ($(1/2)+(60:1)$) {};
    \node[vertex] (3/4) at (3,0) {};
        \node[vertex] (3/4 1) at ($(3/4)+(120:.5)$) {};
        \node[vertex] (3/4 2) at ($(3/4)+(60:.5)$) {};
    \node[vertex] (7/8) at (3.5,0) {};
        \node[vertex] (7/8 1) at ($(7/8)+(120:.25)$) {};
        \node[vertex] (7/8 2) at ($(7/8)+(60:.25)$) {};
    \node[vertex] (1) at (4,0) {};
    \draw (0) to (1/2);
    \draw (1/2) to (3/4);
        \draw[dotted] (1/2) to (1/2 1);
        \draw[dotted] (1/2) to (1/2 2);
    \draw (3/4) to (7/8);
        \draw[dotted] (3/4) to (3/4 1);
        \draw[dotted] (3/4) to (3/4 2);
    \draw (7/8) to (1);
        \draw[dotted] (7/8) to (7/8 1);
        \draw[dotted] (7/8) to (7/8 2);
    \end{scope}
\end{tikzpicture}
\caption{The two generators $g_0$ and $g_1$ of the copy $H$ of Thompson's group $F$ (drawn with undirected edges, see \cref{rmk.dendrite.undirected}).}
\label{fig:H:generators}
\end{figure}

\medskip %layout
\begin{lemma}
\label{lem:thomp}
For every $n \geq 3$, the dendrite rearrangement group $G_n$ contains an isomorphic copy $H$ of Thompson's group $F$ that acts on $A_0$ as $F$ does on the unit interval.
In particular, this subgroup $H$ is generated by the two elements portrayed in \cref{fig:H:generators}.
\end{lemma}

\begin{proof}
Consider the elements $g_0$ and $g_1$ depicted \cref{fig:H:generators} (a precise definition in terms of sequences is given below) and let $H = \langle g_0, g_1 \rangle$.
A quick comparison with the standard generators $X_0$ and $X_1$ of Thompson's group $F$ (\cref{fig.F.generators}) immediately shows that the subgroup $H$ is isomorphic to Thompson's group $F$.
\end{proof}

It may be noted that this copy $H$ of Thompson's group $F$ can also be produced by an application of \cite[Proposition 2.8]{BF19}.
In practice, up to ignoring a 180-degree flip (and, in general, the elements that switch the endpoints $\llbracket 1\overline{n} \rrbracket$ and $\llbracket \overline{n} \rrbracket$), the subgroup $H$ can be thought of as the rearrangement group obtained by forgetting about the edges $2, 3, \dots, n-1$ in both the base and the replacement graphs.

More detailed definitions of the generators $g_0$ and $g_1$ are the following (here we do not use undirected edges for the sake of completeness, so $g_0( \llbracket 11 \rrbracket) = \llbracket n1 \rrbracket$ and $g_1( \llbracket n1n \rrbracket) = \llbracket nn1 \rrbracket$ are expanded and follow the rules described in \cref{rmk.dendrite.undirected}):
\begin{center}
\begin{tabular}{ l r c l }
    $g_0:$ & $1 n$ & $\mapsto$ & $1$ \\
           & $1 i$ & $\mapsto$ & $\tilde{i}$ \\
           & $1 1 1$ & $\mapsto$ & $n 1 n = n 1 \tilde{1}$ \\
           & $1 1 i$ & $\mapsto$ & $n 1 \tilde{i}$ \\
           & $1 1 n$ & $\mapsto$ & $n 1 1 = n 1 \tilde{n}$ \\
           & $i$  & $\mapsto$ & $ni$ \\
           & $n$ & $\mapsto$ & $n n$ \\
           & \\
           &
\end{tabular}
\hspace{20pt}
\begin{tabular}{ l r c l }
    $g_1:$ & $1$ & $\mapsto$ & $1$ \\
           & $i$ & $\mapsto$ & $i$ \\
           & $n 1 n$ & $\mapsto$ & $n 1 = n \tilde{n}$ \\
           & $n 1 i$ & $\mapsto$ & $n \tilde{i}$ \\
           & $n 1 1 1$ & $\mapsto$ & $n n 1 n = n n 1 \tilde{1}$ \\
           & $n 1 1 i$ & $\mapsto$ & $n n 1 \tilde{i}$ \\
           & $n 1 1 n$ & $\mapsto$ & $n n 1 1 = n n 1 \tilde{n}$ \\
           & $n i$  & $\mapsto$ & $n n i$ \\
           & $n n$ & $\mapsto$ & $n n n$
\end{tabular}
\end{center}
where $\tilde{j}$ denotes $n+1-j$ and the $i$'s always represent all possible choices among $2, \dots, n-1$.

In \cref{sec:gen}, when more about the transitive properties of dendrite rearrangement groups will have been discussed, \cref{lem:thomps} will show that there are many more copies of Thompson's group $F$ nested inside $G_n$, and these will help us shed light on more aspects of these rearrangement groups, including generation and density.

Finally, we point out that actions of Thompson's group $F$ on dendrites have been considered in the literature.
For example, in \cite[Corollary 1.5]{DM18} it is proved that such actions have an orbit that consists of one or two points, and in \cite[Corollary 4.5]{FiniteOrbits} it is shown that such actions are equicontinuous on their minimal sets.
The dissertation \cite{SmithDendrite} is also about a Thompson-like group acting on a dendrite (namely, the Julia set for the complex map $z \to z^2 +i$, which is homeomorphic to the Wa\.zewski dendrite $D_3$):
differently from our dendrite rearrangement groups, the group studied in \cite{SmithDendrite} is much more $T$-like, as it is defined as a group of those homeomorphic of the unit circle that preserve a lamination induced by the Julia set.
It seems likely that this group belongs to the family of orientation-preserving dendrite rearrangement groups defined in \cref{sub:orientation}.

%%%%%%%%%%%%%%%%%%%%%%%%%

\section{Generation of Dendrite Rearrangement Groups}
\label{sec:gen}

In this section we prove that dendrite rearrangement groups are generated by the Thompson subgroup $H$ from \cref{lem:thomp} along with the permutation subgroup $K_{p_0}$ on the central branch point $p_0 = \llbracket 11\overline{n} \rrbracket = \dots = \llbracket n1\overline{n} \rrbracket$ from \cref{lem:perm}.

\subsection{Transitive Properties of the Thompson Subgroups}
\label{sub.transitive.properties.thompson.subgroups}

For any $k \in \mathbb{N}$, let $\mathrm{Br}^{(k)}(X)$ be the set of all subsets of $\mathrm{Br}(X)$ consisting of exactly $k$ elements.
From \cref{lem:thomp}, it is clear that the action of $H$ on the branch points of $A_0^\mathrm{o}$ (the interior of the arc $A_0$) corresponds to the action of $F$ on the dyadic points of $(0,1)$.
In particular, since it is known that for every $k \in \mathbb{N}$ Thompson's group $F$ acts transitively on the set of $k$-tuples of dyadic points of $(0,1)$, we have the following.

\medskip %layout
\begin{lemma}
\label{lem:thomp:trans}
For every $k \in \mathbb{N}$, the subgroup $H$ acts transitively on $\mathrm{Br}^{(k)}(A_0^\mathrm{o})$.
\end{lemma}

Using this Lemma, we can prove the following useful Proposition.

\medskip %layout
\begin{proposition}
\label{prop:gen:trans}
The subgroup $\langle H, K_p \mid p \in \mathrm{Br} \rangle$ is transitive on $\mathrm{Br}$.
\end{proposition}

\begin{proof}
Let $p$ be a branch point and consider the path $P$ from $p_0$ to $p$ in the minimal graph expansion where $p$ appears.
We want to prove that $p$ can be mapped to $p_0$ by an element of $\langle H, K_p \mid p \in \mathrm{Br} \rangle$.
We will do this by decreasing the length $k$ of $P$ in order to conclude by induction.

If $k=0$ then $p_0 = p$ and there is nothing to prove, so we can immediately suppose that $k \geq 1$.
Let $p_1$ be the vertex of $P$ that is adjacent to $p_0$.
We can suppose that $p_1$ belongs to the branch $\llbracket 1 \rrbracket$, as otherwise applying an element of $K_{p_0}$ will fix that.
Then both $p_0$ and $p_1$ are branch points of the arc $A_0$, so \cref{lem:thomp:trans} allows us to find an element of $H$ that maps $p_1$ to $p_0$ and that reduces the length of the path as needed thanks to well-known transitive properties of Thompson's group $F$ (for a recent example of a similar application of such properties, one can see Case 2 of the proof of \cite[Theorem 5.2]{airplane}).
\end{proof}

Now let $A$ be any rational arc (\cref{def:arc:type}) and consider the set $\mathrm{Br}(A^\mathrm{o})$ of branch points that are contained in the interior of $A$.
The next Lemma shows that $G_n$ contains a copy of Thompson's group $F$ for every arc $A = [p,q]$, where $p, q \in \mathrm{Br} \cup \mathrm{REn}$.
The overall strategy of the proof is similar that of \cite[Theorem 5.2]{airplane}, which is about the generating set of the airplane rearrangement group $T_A$.
The similarity is no coincidence:
\cref{lem:thomps} will be useful shortly, when proving that $G_n$ is generated by $H$ and $K_{p_0}$ (\cref{thm:gen}), and the close relation between the airplane rearrangement group and dendrite rearrangement groups will be discussed in \cref{sub:Air}.

\medskip %layout
\begin{lemma}
\label{lem:thomps}
Given a dendrite edge replacement system $\mathcal{D}_n$, suppose that $A$ is a rational arc.
Then $G_n$ contain an isomorphic copy $H_A$ of Thompson's group $F$ that acts on $A$ as Thompson's group $F$ does on $[0,1]$ and on $\mathrm{Br}(A^\mathrm{o})$ as $F$ does on the set of dyadic points of $(0,1)$.
In particular, $H_A$ acts transitively on $\mathrm{Br}^{(k)}(A^\mathrm{o})$.
\end{lemma}

\begin{proof}
It suffices to prove the statement for an EE-arc $A$.
Indeed, each BE-arc or BB-arc is contained in an EE-arc and corresponds to a dyadic interval under the identification of $A$ with $[0,1]$ described in \cref{rmk:arcs:dyadic}, and it is known that $F$ contains a copy of itself on every dyadic interval that is contained in $[0,1]$.

Let $A = [p,q]$ be an EE-arc.
As explained in \cref{rmk:paths}, consider the minimal graph expansion $E$ of $\mathcal{D}_n$ that contains the rational endpoints $p$ and $q$ as vertices.
In this graph, the arc $A$ corresponds to a path $P$ in the sense that, if $e_1, \dots, e_k$ are the edges of $P$,
then
\[ A = \big\{ \llbracket e_i \alpha \rrbracket \in D_n \mid i = 1, \dots, n \text{ and } \alpha \text{ is a sequence in } \{ 1, n \} \big\}. \]
By induction on the length $k$ of the path $P$ in the minimal graph where both $p$ and $q$ appear, we will find a $g \in \langle K_p \mid p \in \mathrm{Br} \rangle$ such that $g(P)$ corresponds to $A_0$ in the sense described above, so that $g(A) = A_0$.
The goal will be to find an element that modifies $P$ so that a reduction occurs, in order to decrease $k$.

In order to make this easier, we will suppose that $P$ contains the central branch point $p_0$.
We can do this without loss of generality up to conjugating by an element of $\langle H, K_p \mid p \in \mathrm{Br} \rangle$ that maps some branch point of $P$ to $p_0$, which exists thanks to \cref{prop:gen:trans}.

The base case is $k=2$, meaning that $P$ consists of two edges of the base graph.
In this case, the two edges correspond simply to distinct 1-letter words $i$ and $j$ chosen among $1, \dots, n$.
Then $g$ is the unique element of $K_{p_0}$ that switches $i$ with $1$ and $j$ with $n$, leaving everything else unchanged.

For $k \geq 3$, either $p$ or $q$ is not among the endpoints $\llbracket i \overline{n} \rrbracket$, and without loss of generality we will assume that it is $p$.
Denote by $r$ the vertex of $P$ that is adjacent to $p$.
Observe that $r$ must be a branch point, so $r = \llbracket x i 1 \overline{n} \rrbracket$ for some $x$ that does not end with $1$ nor with $n$ (\cref{rmk:br:ex:pts}).
Then it is not possible that $p$ is $\llbracket x \overline{n} \rrbracket$, otherwise $E$ would not be the minimal graph expansion where $p$ and $q$ appear, because the cell $\llbracket x \rrbracket$ could be reduced.
\cref{fig:lem:thomps} shows why this is the case with an example.
We now find an element $g \in K_r$ that causes precisely this reduction, i.e., $g$ is such that $g(P)$ has one edge less than $P$ in the minimal graph expansion where both $g(p)$ and $g(q) = q$ appear (note that $g(q) = q$ because, by an assumption made at the beginning of the proof, $P$ must include the central branch point $p_0$, so $q$ is not included in any of the two cells permuted by $g$).
An example is depicted in \cref{fig:path:reduction}.
Since $r = \llbracket x i 1 \overline{n} \rrbracket$, we have that $p = \llbracket x j \overline{n} \rrbracket$ for some $j$ that is not $n$ (as noted above) nor $1$ (or $p$ would not be an endpoint).
Then we choose $g$ to be the unique involution that switches the edges $xj$ and $xn$ without changing anything else.
Now, $g(P)$ admits the aforementioned reduction and we can apply our induction hypothesis on $g(p)$ and $g(q) = q$, finding an element $h \in \langle K_p \mid p \in \mathrm{Br} \rangle$ such that $hg(A) = A_0$ and $H^{hg}$ is the desired copy of Thompson's group $F$ acting on $A$.
\end{proof}

\begin{figure}
\centering
\begin{tikzpicture}[scale=1.333]
    \node[vertex] (c) at (0,0) {};
    \node[vertex] (l) at (-2,0) {};
    \node[vertex] (t) at (0,2) {};
    \node[vertex] (b) at (0,-2) {};
    \node[vertex] (r) at (2,0) {};
        \node[vertex] (rc) at (1,0) {};
        \node[vertex,blue] (rt) at (1,.8) {}; \draw (1,.8) node[above right,blue]{$p$};
            \node[vertex,Green] (rtc) at (1,.4) {}; \draw (1,.4) node[above right,Green]{$r$};
            \node[vertex] (rtl) at (.7,.4) {};
            \node[vertex] (rtr) at (1.3,.4) {};
        \node[vertex] (rb) at (1,-.8) {};
    \draw[black] (c) to (l);
    \draw[black] (c) to (t);
    \draw[black] (c) to (b);
    \draw[black] (rc) to (c);
    \draw[black] (rc) to (r);
    \draw[red] (rtc) to (rtl);
    \draw[red] (rtc) to (rtr);
    \draw[red] (rtc) to (rc);
    \draw[red] (rtc) to (rt);
    \draw[black] (rc) to (rb);
\end{tikzpicture}
\caption{In the proof of \cref{lem:thomps}, \textcolor{blue}{$p$} and \textcolor{Green}{$r$} cannot be as in this example or the \textcolor{red}{red edges} could be reduced.}
\label{fig:lem:thomps}
\end{figure}
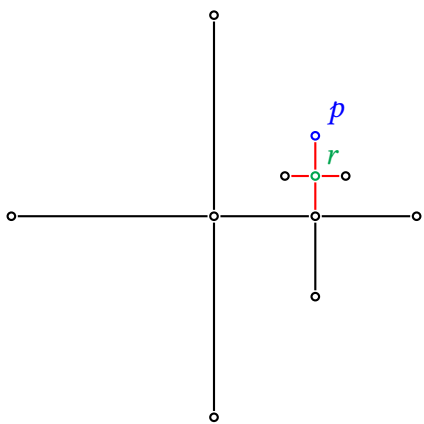

\begin{figure}
\centering
\begin{tikzpicture}[scale=1.333]
    \node[vertex] (c) at (0,0) {};
    \node[vertex] (l) at (-2,0) {};
    \node[vertex] (t) at (0,2) {};
    \node[vertex] (b) at (0,-2) {};
    \node[vertex] (r) at (2,0) {};
        \node[vertex] (rc) at (1,0) {};
        \node[vertex] (rt) at (1,.8) {};
            \node[vertex,Green] (rtc) at (1,.4) {}; \draw (1,.4) node[above left,Green]{$r$};
            \node[vertex] (rtl) at (.7,.4) {};
            \node[vertex,blue] (rtr) at (1.3,.4) {}; \draw (1.3,.4) node[right,blue]{$p$};
        \node[vertex] (rb) at (1,-.8) {};
    \draw[black] (c) to (l);
    \draw[black] (c) to (t);
    \draw[black] (c) to (b);
    \draw[black] (rc) to (c);
    \draw[black] (rc) to (r);
    \draw[red] (rtc) to (rtl);
    \draw[red] (rtc) to (rtr);
    \draw[red] (rtc) to (rc);
    \draw[red] (rtc) to (rt);
    \draw[black] (rc) to (rb);
    \draw[<->,gray] (rtr) to[out=90,in=0,looseness=1.25] node[above right]{$g$} (rt);
\end{tikzpicture}
\caption{An example of the element $g \in K_r$ found in the proof of \cref{lem:thomps} that causes a reduction of the \textcolor{red}{red edges}.}
\label{fig:path:reduction}
\end{figure}
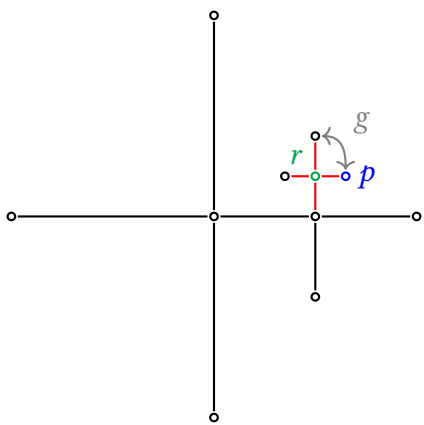

Reading along the lines of the previous proof, one may note that we have actually proved that the $H_A$'s are included in the subgroup $\langle H, K_p \mid p \in \mathrm{Br} \rangle$.
However, as stated in the next Proposition, the sole permutation subgroup $K_{p_0}$, together with $H$, generates every other permutation subgroup $K_p$, meaning that
\[ \langle H, K_p \mid p \in \mathrm{Br} \rangle = \langle H, K_{p_0} \rangle. \]
Thus, the other copies $H_A$ of Thompson's group $F$ discussed in the previous Lemma are all included in $\langle H, K_{p_0} \rangle$.

\medskip %layout
\begin{proposition}
\label{prop:thomps:in:gen}
For every $p \in \mathrm{Br}$, the permutation subgroup $K_p$ is included in $\langle H, K_{p_0} \rangle$.
It follows, in particular, that for every rational arc $A$ the Thompson subgroup $H_A$ is included in $\langle H, K_{p_0} \rangle$.
\end{proposition}

The proof of this Proposition is omitted, as it is not too different from those already seen throughout this section.
One essentially needs to build an element $g \in \langle H, K_{p_0} \rangle$ such that $g(p) = p_0$ and that is ``rigid'' everywhere else, so that $g$ conjugates $K_{p_0}$ to $K_p$.
In practice, one can use strategies similar to those used to prove \cref{prop:gen:trans} and \cref{lem:thomps}, using the $k$-transitivity of $H$ from \cref{lem:thomp:trans} to make sure that the element that is being built is ``rigid'' as required.

\begin{remark}
From the previous discussion we can deduce that, given any two rational arcs $A_1$ and $A_2$, the subgroups $H_{A_1}$ and $H_{A_2}$ are conjugate in $G_n$ if and only if $A_1$ and $A_2$ are of the same type (in the sense of \cref{def:arc:type}).
The permutation subgroups $K_p$ are instead all conjugate copies of one another.
\end{remark}

\subsection{A Finite Generating Set}
\label{sub:gen}

Using similar ideas to those of \cref{lem:thomps} and the fact that Thompson's group $F$ is finitely generated, we can produce a finite generating set for $G_n$.

\medskip %layout
\begin{theorem}
\label{thm:gen}
$G_n = \langle H, K_{p_0} \rangle$. In particular, $G_n$ is finitely generated.
\end{theorem}

\begin{proof}
Let $g \in G_n$.
We need to show that $g \in \langle H, K_{p_0} \rangle$, knowing by the discussions in the previous subsection that $\langle H, K_{p_0} \rangle$ includes $H_A$ and $K_p$ for every rational arc $A$ and every branch point $p$.
Thanks to \cref{prop:thomps:in:gen} and the transitive properties of each $H_A$ on the branch points of $A^\mathrm{o}$, up to composing $g$ with a suitable element of $\langle H, K_{p_0} \rangle$ we can suppose that $g$ fixes the central branch point $p_0$.
We consider the reduced graph pair diagram $(D, R, f)$ for $g$, and we proceed by induction on the number $k$ of branch points that appear as vertices in $D$.
The base case $k=1$ simply means that $g \in K_{p_0}$, so we can immediately suppose that $k \geq 2$.

The simplest case is when the $k$ branch points of $D$ are scattered among at least two distinct branches at $p_0$.
In this case, there is an element $h$ of $K_{p_0}$ such that $h g$ fixes setwise each of the $n$ branches at $p_0$.
The element $h g$ can then be decomposed as the product of its restrictions in each branch at $p_0$, and since there are at least two of these branches containing one of the $k$ branch points, each of these restrictions must contain less than $k$ branch points.
This allows conclude by applying our induction hypothesis on each of the restrictions.

Suppose instead that the $k$ branch points of $D$ all belong to the same branch at $p_0$.
We can suppose that the branch is $\llbracket 1 \rrbracket$ up to composing $g$ with a suitable element of $K_{p_0}$.
A similar strategy to that described at the end of the proof of \cref{prop:gen:trans} (which again makes heavy use of the transitive properties of $H$ on $\mathrm{Br}(A_0)$ stated in \cref{lem:thomp:trans} and is not very different from Case 2 of the proof of \cite[Theorem 5.2]{airplane}) shows that it is possible to compose $g$ on its right and on its left with suitably chosen elements of $H$ in order to cause a reduction in the graph pair diagram, essentially by mapping the vertex in $\llbracket 1 \rrbracket$ that is closest to $p_0$ to $p_0$ and $p_0$ to the branch point of the cell $\llbracket 1 \rrbracket$, i.e., $\llbracket 1 i 1 \overline{n} \rrbracket$.
This produces an element whose reduced graph pair diagram has less than $k$ branch points in its domain graph, allowing us to conclude by our induction hypothesis.
\end{proof}

The reader may have noted that in this proof we are not truly using the result of \cref{lem:thomps} to its full extent, and we are instead partially replicating the strategy of its proof.
Indeed, one could prove \cref{thm:gen} by using only \cref{lem:thomp:trans} alongside with a weaker version of \cref{lem:thomps} that just states that $\langle H, K_{p_0} \rangle$ is transitive on $\mathrm{Br}$ (this is how it was done in \cite[Lemma 5.1]{airplane} for the airplane rearrangement group).
However, we decided to be slightly more redundant in these proofs in order to produce a more in-depth description of $G_n$, because we found that the many copies of Thompson's group $F$ found by \cref{lem:thomps} give a broad understanding of the kind of transformations of $D_n$ that $G_n$ is capable of.
Additionally, \cref{lem:thomps} will be useful when proving \cref{prop:trans}.

Now that we know that each $G_n$ is finitely generated, it is natural to ask the following question, which we do not investigate here.

\begin{question}
What further finiteness properties do dendrite rearrangement groups have?
Are they finitely presented?
Are they $F_\infty$?
\end{question}

%%%%%%%%%%%%%%%%%%%%%%%%%

\section{Density of Dendrite Rearrangement Groups}
\label{sec:dns}

In this section we find useful transitive properties of the action of dendrite rearrangement groups on the set of branch points that allow us to ``approximate'' any homeomorphism of $D_n$ by rearrangements.
We will denote by $\mathbb{H}_n$ the full group of homeomorphisms of the Wa\.zewski dendrite $D_n$.
Since the topology of $\mathbb{H}_n$ is determined by the action on $\mathrm{Br}$, we will use this transitive property to prove that dendrite rearrangement groups are dense in the respective full homeomorphism groups.
This argument is similar to that used in \cite{BasilicaDense} to show that the rearrangement group $T_B$ of the basilica is dense in the full group of homeomorphisms of the basilica.

\subsection{Transitive Properties of Dendrite Rearrangement Groups}

The next Proposition is an analog of \cite[Proposition 6.1]{DM19}, which goes to show that $G_n$ has many of the transitive properties of $\mathbb{H}_n$ despite being only a countable group.
The construction of the tree $T(\mathcal{F})$ described below is heavily inspired by the one given at the beginning of section 6 of \cite{DM19}, which in turn is not much different from the construction used in \cite[Theorem 4.14]{BF19} for studying the rearrangement group of the Vicsek fractal (which is homeomorphic to a dendrite, see \cref{sub:Vic}).
Similar tree approximations for dendrites are common and have appeared elsewhere in the literature (see, for example, section X.3 of \cite{continua}).

\begin{figure}
\centering
\begin{tikzpicture}
    \node[anchor=south west,inner sep=0] (image) at (0,0) {\includegraphics[width=0.448\textwidth]{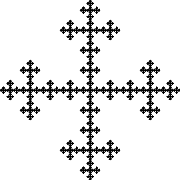}};
    \begin{scope}[x={(image.south east)},y={(image.north west)}]
        \draw[VioletRed] (.4966,.7227) node[circle,fill,inner sep=2pt]{};
        \draw[Green] (.8259,.5034) node[circle,fill,inner sep=2pt]{};
        \draw[BurntOrange] (.6062,.174) node[circle,fill,inner sep=2pt]{};
        \draw[blue] (.3868,.1744) node[circle,fill,inner sep=2pt]{};
    \end{scope}
\end{tikzpicture}
\hspace{15pt}
\begin{tikzpicture}[scale=2.8]
    \draw[white] (-1,-1) -- (-1,1) -- (1,1) -- (1,-1) -- (-1,-1);
    \node[VioletRed,circle,fill,inner sep=2pt] (purple) at (0,.45) {};
    \node[Green,circle,fill,inner sep=2pt] (green) at (.667,0) {};
    \node[BurntOrange,circle,fill,inner sep=2pt] (red) at (.25,-.667) {};
    \node[blue,circle,fill,inner sep=2pt] (blue) at (-.25,-.667) {};
    \node[black,circle,fill,inner sep=2pt] (center) at (0,0) {};
    \node[black,circle,fill,inner sep=2pt] (bottom) at (0,-.667) {};
    \draw (center) to (green);
    \draw (center) to (purple);
    \draw (center) to (bottom);
    \draw (bottom) to (red);
    \draw (bottom) to (blue);
\end{tikzpicture}
\caption{An example of the tree $T(\mathcal{F})$ depicted beside the finite subset $F \subset \mathrm{Br}$, for $n=4$.}
\label{fig:T(F)}
\end{figure}
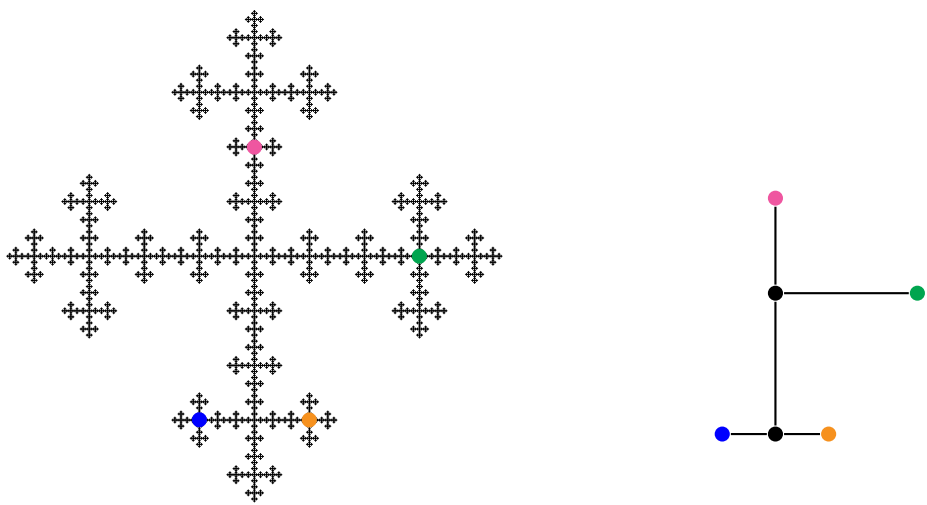

\phantomsection\label{txt:trees}
Let $\mathcal{F}$ be a finite subset of $\mathrm{Br}$ and consider the unique subdendrite $[\mathcal{F}]$ of $D_n$ that contains $\mathcal{F}$ and is minimal with respect to set inclusion.
It is homeomorphic to the geometric realization of a finite simplicial tree that has a vertex for each element of $\mathcal{F}$ (and possibly more vertices where the tree needs to branch out).
This simplicial tree is not unique, since we can choose to add vertices in the middle of any edge, so we choose $T(\mathcal{F})$ to be the minimal among these trees, and since $\mathcal{F} \subset \mathrm{Br}$ we have that each vertex of $T(\mathcal{F})$ is a branch point.
A simple example is depicted in \cref{fig:T(F)}.

\medskip %layout
\begin{proposition}
\label{prop:trans}
Given two finite subsets $\mathcal{F}_1$ and $\mathcal{F}_2$ of $\mathrm{Br}$, any graph isomorphism between $T(\mathcal{F}_1)$ and $T(\mathcal{F}_2)$ can be extended to an element of $G_n$.
\end{proposition}

\begin{proof}
Given a graph isomorphism $f \colon T(\mathcal{F}_1) \to T(\mathcal{F}_2)$, we proceed by induction on the number $m$ of vertices of $T(\mathcal{F}_1)$ and $T(\mathcal{F}_2)$.
If $m = 1$, then $\mathcal{F}_1 = \{p\}$, $\mathcal{F}_2 = \{q\}$ and \cref{prop:gen:trans} provides an element of $G$ that maps $p$ to $q$.

Suppose that $m \geq 2$.
Consider the trees $T(\mathcal{F}_1)^*$ and $T(\mathcal{F}_2)^*$
obtained by removing any leaf $l$ from $T(\mathcal{F}_1)$ and its image $f(l)$ from $T(\mathcal{F}_2)$, respectively.
Observe that $T(\mathcal{F}_1)^* = T(V_1 \setminus \{l\})$ and $T(\mathcal{F}_2)^* = T(V_2 \setminus \{f(l)\})$, where each $V_i$ is the set of vertices of the respective $T(\mathcal{F}_i)$, each corresponding to a branch point.
It is worth noting that $T(V_1 \setminus \{l\})$ and $T(V_2 \setminus \{f(l)\})$ are not always the same as $T(\mathcal{F}_1 \setminus \{l\})$ and $T(\mathcal{F}_2 \setminus \{f(l)\})$, since the two former trees might need to include a vertex that would be unnecessary in the latter.
In \cref{fig:T(F)} the removal from $\mathcal{F}$ of any of its points also causes the removal of one of the two additional vertices of $T(\mathcal{F})$;
for example removing the green vertex also cause the removal of the black vertex that is adjacent to it.
By considering $T(V_1 \setminus \{l\})$ and $T(V_2 \setminus \{f(l)\})$ we are sure that only one vertex is removed.

We can then apply our induction hypothesis, which produces an element $g^* \in G$ that extends the graph isomorphism $f^* \colon T(\mathcal{F}_1)^* \to T(\mathcal{F}_2)^*$ induced (uniquely) from $f$.
It only remains to adjust the element $g^*$ so that it also maps $l$ to $f(l)$.
Let $v$ be unique vertex of $T(\mathcal{F}_2)$ that is adjacent to the leaf $f(l)$.
Note that there is a branch $B_1$ (respectively, $B_2$) at $v$ that includes $v$ and $f(l)$ (respectively, $g^*(l)$) but does not include any other vertex of $T(\mathcal{F}_2)$;
it is possible that $B_1$ and $B_2$ are the same branch.
By \cref{lem:perm} we can choose an element $h_1$ of $K_v$ that switches the branches $B_1$ and $B_2$ and fixes every other point (in case $B_1 = B_2$, we choose $h_1$ to be the identity).
Then we have that $h_1$ fixes every vertex of $T(\mathcal{F}_2)^*$ and that $\tilde{l} = h_1 g^* (l)$ belongs to the same branch $B_1$ that contains $f(l)$.
Now that we have two branch points $\tilde{l} = h_1 g^* (l)$ and $f(l)$ belonging to the same branch $B_1$, we can use \cref{lem:thomps} to find an element $h_2 \in H_{A}$ that maps $h_1 g^* (l)$ to $f(l)$, where $A$ is any rational arc included in $B_1$ whose interior includes $f(l)$ and $\tilde{l}$.
Clearly $h_2$ must fix $v \in B_1$, and since every other vertex of $T(\mathcal{F}_2)^*$ does not belong to $B_1$ we have that the element $g = h_2 h_1 g^*$ maps each vertex of $T(\mathcal{F}_1)$ to its image under the action of $f$, as required.
\end{proof}

This transitive property of $G_n$ on the set of branch points is the same that the full homeomorphism group $\mathbb{H}_n$ has by \cite[Proposition 6.1]{DM19}, so \cref{prop:trans} is essentially saying that every dendrite rearrangement group is as transitive as a subgroup of $\mathbb{H}_n$ can be on the set of branch points.

As an immediate consequence of \cref{prop:trans}, we have the two following facts (which are analogous to Corollaries 6.5 and 6.7 of \cite{DM19}, respectively).

\medskip %layout
\begin{corollary}
The stabilizer of a branch point is a maximal proper subgroup of $G_n$.
\end{corollary}

\begin{corollary}
The action of $G_n$ on the set $\mathrm{Br}$ is oligomorphic.
\end{corollary}

\phantomsection\label{txt.oligomorphic}
Recall that a $G$-action on a set $X$ is \textbf{oligomorphic} if, for every $m \in \mathbb{N}$, the diagonal $G$-action on the set of $m$-tuples of $X$ has finitely many orbits.
For more information about oligomorphic group actions, see \cite{Oligomorphic}.

\subsection{Density in the Full Group of Homeomorphisms}
\label{sub:dns}

Here we prove that \cref{prop:trans} has the following consequence.

\medskip %layout
\begin{theorem}
\label{thm:dense}
For every $n \geq 3$, the rearrangement group $G_n$ for the dendrite edge replacement system $\mathcal{D}_n$ is dense in $\mathbb{H}_n = \mathrm{Homeo}(D_n)$.
\end{theorem}

As done in \cite{DM19}, here we endow the group $\mathbb{H}_n$ with its topology of uniform convergence, which makes it a Polish group (i.e., a separable and completely metrizable topological group).
By \cite[Proposition 2.4]{DM19}, this topology of $\mathbb{H}_n$ is inherited by the pointwise convergence topology of $\mathrm{Sym}\big(\mathrm{Br}\big)$, so in order to prove \cref{thm:dense} it suffices to show that the permutations of $\mathrm{Br}$ induced by elements of $\mathbb{H}_n$ can be ``approximated'' by those of $G_n$.
Thus, \cref{thm:dense} is equivalent to the following.

\medskip %layout
\begin{claim}
\label{clm:dns}
Let $\phi \in \mathbb{H}_n$.
For every $k \geq 1$ and for each $p_1, \dots, p_k \in \mathrm{Br}$ there exists a rearrangement $g_k \in G_n$ such that $g_k(p_i) = \phi(p_i)$ for all $i=1,\dots,k$.
\end{claim}

\begin{proof}
Consider a homeomorphism $\phi$ of $D_n$ such that $\phi(p_i) = q_i$ for branch points $p_1, \dots, p_k, q_1, \dots, q_k$.
Let $\mathcal{F}_1 = \{p_1, \dots, p_k\}$ and $\mathcal{F}_2 = \{q_1, \dots, q_k\}$ and consider the trees $T(\mathcal{F}_1)$ and $T(\mathcal{F}_2)$ obtained as described at \cpageref{txt:trees}.
Consider the unique minimal subdendrites $[\mathcal{F}_1]$ and $[\mathcal{F}_2]$ of $D_n$ that contain $\mathcal{F}_1$ and $\mathcal{F}_2$, respectively;
these are naturally homeomorphic to the geometric realizations of $T(\mathcal{F}_1)$ and $T(\mathcal{F}_2)$.
The restriction of $\phi$ to $[\mathcal{F}_1]$ is a homeomorphism $\phi|_{[\mathcal{F}_1]} \colon [\mathcal{F}_1] \to [\mathcal{F}_2]$ that maps each vertex $p_i$ of $T(\mathcal{F}_1)$ to the vertex $q_i$ of $T(\mathcal{F}_2)$, so it restricts to a graph isomorphism $T(\mathcal{F}_1) \to T(\mathcal{F}_2)$ defined by the mapping $p_i \mapsto q_i$.
By \cref{prop:trans}, there exists a rearrangement $g \in G_n$ that extends the graph isomorphism $T(\mathcal{F}_1) \to T(\mathcal{F}_2)$, meaning that it maps each $p_i$ to $q_i$, as needed.
\end{proof}

In the next section, with \cref{prop:comm:trans,cor:comm:dns}, using these very same arguments, we will see that the commutator subgroup $[G_n, G_n]$ shares strong transitive properties with $G_n$ and thus it is also dense in $\mathbb{H}_n$ despite being significantly ``smaller'' than $G_n$ (it has infinite index by \cref{thm:commutator}).

Being dense in an overgroup is not new behavior in the world of Thompson groups.
Indeed, \cite[Corollary A5.8]{PL-homeo} and \cite[Proposition 4.3]{Tdense} show that Thompson groups $F$ and $T$ are dense in the groups $\mathrm{Homeo}^+([0,1])$ and $\mathrm{Homeo}^+(S^1)$ of orientation-preserving homeomorphisms of the unit interval and the unit circle, respectively.
Also, Thompson's group $V$ is dense in both the group $\mathrm{AAut}(\mathcal{T}_2)$ of almost automorphisms of the binary tree (noted in \cite{Vdense}) and also in the group $\mathrm{Homeo}(\mathfrak{C})$ of homeomorphisms of the Cantor space with its Polish compact-open topology (see \cref{rmk:V:dns} below).
Notably, the very recent work \cite{BasilicaDense} shows that the rearrangement group of the basilica is dense in the group of all orientation-preserving homeomorphisms of the basilica.
Later in \cref{sub:Air} we will see that the airplane rearrangement group is also dense in the group of orientation-preserving homeomorphisms of a dendrite.
All of these results naturally prompt the following question, which was also asked in \cite[Final remarks (b)]{BasilicaDense}.

\begin{question}
When is a rearrangement group of a topological space $X$ dense in the full group of homeomorphisms of $X$?
\end{question}

This question might be related to recent findings and questions about quasi-symmetries of fractals. Namely, \cite{BasilicaQuasiSymmetries} shows that the basilica rearrangement group $T_B$ is in a certain sense dense in the group of planar quasi-symmetries of the basilica Julia set, and \cite{QuasiSymmetries} asks if it is always possible to find a finitely generated dense subgroup of the orientation-preserving group of homeomorphisms.

We believe that this question is interesting for practical reasons.
Finding a Thompson-like group $G$ that is dense in an interesting uncountable group $\mathbb{H}$ essentially allows us to approximate the ``non-computable'' group $\mathbb{H}$ with the simpler group $G$ that can probably be easily handled by a computer.

\begin{remark}
\label{rmk:V:dns}
To the best of the author's knowledge, there is no proof in the literature that $V$ is dense in the group $\mathrm{Homeo}(\mathfrak{C})$ of homeomorphisms of the Cantor set with its compact-open topology.
Thus, although we do not claim the novelty of this result, we provide a quick sketch of proof here for the sake of completeness.

First, observe that $\mathrm{Homeo}(\mathfrak{C})$ with its compact-open topology is the same as the group of automorphisms of the Boolean algebra of clopen subsets of the Cantor space $\mathfrak{C}$ with its pointwise convergence topology (this was noted, for example, at the beginning of section 2 of \cite{TDLCCantor}).
Then it suffices to show that, for every homeomorphism $\phi$ of $\mathfrak{C}$ and for each finite set $\mathcal{S}$ of clopen subsets of $\mathfrak{C}$, there exists an element $g$ of Thompson group $V$ such that $g(A) = \phi(A)$ for all $A \in \mathcal{S}$.
Showing this is straightforward:
if we view the Cantor space $\mathfrak{C}$ as the set of infinite sequences in the alphabet $\{0,1\}$, then any clopen set is the disjoint union of finitely many cones (where by \textit{cone} we mean a subset of $\mathfrak{C}$ consisting of all of those sequences that have a certain common prefix);
if we refine the set $\mathcal{S}$ of clopen subsets of $\mathfrak{C}$ to a set $\mathcal{S}^*$ of cones, then it is immediate to build an element of $V$ that maps each cone of $\mathcal{S}^*$ canonically to its image under the homeomorphism $\phi$.
\end{remark}

%%%%%%%%%%%%%%%%%%%%%%%%%

\section{The Commutator Subgroup}
\label{sec:comm}

In this section we first show that the abelianization of every dendrite rearrangement group is $\mathbb{Z}_2 \oplus \mathbb{Z}$ (\cref{sub:abel}), then we prove that the commutator subgroup $[G_n, G_n]$ is dense in the full homeomorphism group $\mathbb{H}_n$ and that, with the possible exception of $n=3$, it is simple.

\subsection{The Abelianization of Dendrite Rearrangement Groups}
\label{sub:abel}

In this subsection we describe two families of maps:
the local parity maps $\pi_p$ and the local endpoint derivatives $\partial_q$, both satisfying a chain rule that will allow us to define global versions of these maps.
Together, the parity map and the endpoint derivative will allow us to describe the elements of the commutator subgroup $[G_n, G_n]$ of a dendrite rearrangement group and to compute the abelianization of $G_n$ by applying Schreier's Lemma (which was also used in \cite{Belk_2015} for the computing abelianization of the basilica rearrangement group $T_B$).

\subsubsection{The Parity Map}
\label{sub:parity}

As already discussed at \cpageref{txt:rot:sys} when describing the rotation system on the graph expansions of dendrite replacement systems, the branches at a given branch point $p$ are naturally enumerated by $1, \dots, n$.
More precisely, there is a unique finite word $x$ for which $p = \llbracket x i 1 \overline{n} \rrbracket$;
each branch at $p$ intersects precisely one of the cells $\llbracket x i \rrbracket$ (and thus every cell whose address starts with $xi$):
this induces an enumeration of the branches at $p$ (i.e., a bijection with $\{1, \dots, n\}$).
Note that this enumeration is ``stronger'' than the rotation system introduced at \cpageref{txt:rot:sys}, by which we mean that the enumeration induces this rotation system, but a rotation system does not induce an enumeration.

A homeomorphism $g$ of $D_n$ induces a bijection between the branches at $p$ and those at $g(p)$.
This in turn induces permutation $\sigma_{p,g} \in \mathrm{Sym}(n)$ according to the enumeration of the branches at $p$ and $g(p)$.
We define the \textbf{local parity map} at the branch point $p \in \mathrm{Br}$ as the mapping
\[ \pi_p \colon G_n \to \mathbb{Z}_2 \]
such that $\pi_p(g) = 0$ if $\sigma_{p,g}$ is an even permutation and $\pi_p(g) = 1$ if $\sigma_{p,g}$ is an odd permutation.

Note that the local parity maps satisfy the following chain rule:
\[ \pi_p (g h) = \pi_{h(p)} (g) + \pi_p (h). \]

\phantomsection\label{txt.global.parity.map}
Observe that, whenever a branch point $p$ does not appear in a graph pair diagram of a rearrangement $g$, one has $\pi_p(g)=0$ because canonical homeomorphisms of a cell $\llbracket x \rrbracket$ (which were defined in \cref{def.cell}) preserve the rotation system of the branch point $\llbracket x i 1 \overline{n} \rrbracket$.
In particular, this implies that the definition does not depend on the choice of graph pair diagram chosen to represent $g$.
This, together with the previously described chain rule and the fact that every rearrangement induces a permutation of $\mathrm{Br}$, allow us to define the (global) \textbf{parity map $\Pi$} as
\[ \Pi \colon G_n \to \mathbb{Z}_2, \; g \mapsto \sum_{p \in \mathrm{Br}} \pi_p(g). \]
Because of the chain rule identity, it is easy to see that $\Pi$ is a group morphism.

\subsubsection{The Endpoint Derivative}
\label{sub.endpoint.derivative}

We define the \textbf{local endpoint derivative} at the rational endpoint $q \in \mathrm{REn}$ as the mapping
\[ \partial_q \colon G_n \to \mathbb{Z} \]
defined as follows.
If $(D,R,f)$ is a graph pair diagram for $g$, there is a unique edge $x$ of $D$ that is incident on $q$, and a unique edge $y = f(x)$ of $R$ that is incident on $f(q)$.
Denote by $L_D(x)$ and $L_R(y)$ the amount of $n$'s at the end of $x$ and $y$, without considering the first letter in case the word is $n n \dots n$ (the first letter is special in the sense that it denotes an edge of the base graph instead of the replacement graph).
Then
\[ \partial_q(g) = L_D(x) - L_R(y). \]
This definition does not depend on the choice of graph pair diagram chosen to represent $g$, since expanding $x$ in $D$ always causes the edge expansion of $y$ in $R$, and both the resulting lengths $L_D$ and $L_R$ increase by $1$, leaving $\partial_q$ ultimately unchanged.

The intuition behind this definition is that $L_D(x)$ (or $L_R(y)$) represents how small of a portion of a ``straight branch'' $x$ (or $y$) is, when drawing $z1$ and $zn$ straight and aligned with $z$ when expanding an edge $z$.
In this sense, $\partial_q$ looks like a logarithmic derivative in a small enough neighborhood of the endpoint $q$.

Note that, as the local parity maps do, the endpoint derivatives also satisfy the following chain rule:
\[ \partial_q (g h) = \partial_{h(q)} (g) + \partial_q (h). \]

\phantomsection\label{txt.global.endpoint.derivative}
Observe that, whenever a rational endpoint $q$ does not appear in a graph pair diagram of a rearrangement $g$, one has $\partial_q(g)=0$.
Then, for the same reasons used for the parity map, we can define the (global) \textbf{endpoint derivative $\Delta$} as
\[ \Delta \colon G_n \to \mathbb{Z}, \; g \mapsto \sum_{q \in \mathrm{REn}} \partial_q(g). \]
Because of the chain rule identity, $\Delta$ is a group morphism.

Before proceeding with the computation of the abelianization of $G_n$, we would like to mention that this idea of an endpoint derivative was previously used to characterize the commutator subgroup $[T_A, T_A]$ of the airplane rearrangement group.
This in turn is related to the fact that the commutator subgroup of Thompson's group $F$ consists of the elements of $F$ that act trivially around the neighborhoods, i.e., that preserve the numbers of $0$'s (respectively $1$'s) of the leftmost (respectively rightmost) edge of a graph pair diagram.
There was no parity map for $T_A$ essentially because the subgroup of $T_A$ that corresponds to $K_p \leq D_n$ is Thompson group $T$ (see \cref{sub:Air}), which is simple.

\subsubsection{Computing the Abelianization}

Consider the following map:
\[ \Phi = \Pi \times \Delta \colon G_n \to \mathbb{Z}_2 \oplus \mathbb{Z}, \; g \mapsto \big( \Pi(g), \Delta(g) \big). \]
Clearly $\Phi$ is a group morphism, and it is not hard to see how it behaves on the generators $\{g_0, g_1, \tau_2, \dots, \tau_n \}$ of $G_n$, where each $\tau_i$ is the element of $K_{p_0}$ that corresponds to the permutation $(1 i)$ and $g_0$ and $g_1$ are the generators of $H$ and were represented in \cref{fig:H:generators} and described after \cref{lem:thomp}
\begin{align*}
    \Phi(g_0) = (0,0), \; \Phi(g_1) = (0,-1),\\
    \Phi(\tau_i) = (1,0), \; \forall i = 2, \dots, n.
\end{align*}
The parity map on $g_0$ and $g_1$ is $0$ because both of these elements feature exactly two branch points with the same non-trivial permutation (for $g_0$ these are $\llbracket i \overline{n} \rrbracket$ and $\llbracket 1 1 i \overline{n} \rrbracket$ while the permutation is trivial at $\llbracket 1 i \overline{n} \rrbracket$, while for $g_1$ they are $\llbracket n i \overline{n} \rrbracket$ and $\llbracket n 1 1 i \overline{n} \rrbracket$ while the permutation is trivial at $\llbracket i \overline{n} \rrbracket$ and $\llbracket n i \overline{n} \rrbracket$).
In particular, from these computations we deduce that $\Phi$ is surjective, because for all $\zeta \in \{0,1\}$ and $z \in \mathbb{Z}$ one has $\Phi(\tau_2^\zeta g_1^z) = (\zeta, -z)$.
If we let $K = \mathrm{Ker}(\Phi)$, then $G_n / K \simeq \mathbb{Z}_2 \oplus \mathbb{Z}$, so $[G_n, G_n] \leq K$, so our goal is to show that $K \leq [G_n, G_n]$.

Using the following set of transversals for the quotient $G_n / K$
\[ \left\{ \tau_2^{\zeta} g_1^{-z} \mid \zeta \in \{0,1\}, z \in \mathbb{Z} \right\} \]
and the generating set $\{g_0, g_1, \tau_2 \dots, \tau_n \}$ for $G_n$ (\cref{thm:gen}), together with the application of the identities $\tau_2 g_1 = g_1 \tau_2$ and $\tau_i = \tau_i^{-1}$, an application of Schreier's Lemma (see, for example, \cite[Lemma 4.2.1]{MR1970241}) produces the following generating set for $K$:
\begin{equation}
\label{eq:commutator}
    \left\{
    g_1^{-z} \tau_2 \tau_i g_1^z, \;
    \tau_2^\zeta g_1^{-z} g_0 g_1^z \tau_2^\zeta \,
    \mid \, \zeta \in \{0,1\}, z \in \mathbb{Z}, i = 3, \dots, n
    \right\}.
\end{equation}

This generating set is included in the normal closure of $\{ \tau_2 \tau_i, g_0 \mid i = 3, \dots, n \}$ in $G_n$ and the commutator subgroup is normal in $G_n$, so it suffices to show that $[G_n, G_n]$ contains $\tau_2 \tau_3, \dots, \tau_2 \tau_n$ and $g_0$ in order to be able to conclude that $K \leq [G_n, G_n]$.
For each $i \in \{3, \dots, n\}$ the element $\tau_2 \tau_i$ is an even permutation, so it belongs to $[K_{p_0}, K_{p_0}] \leq [G_n, G_n]$.
As for $g_0$, using the fact that $[H,H] \simeq [F,F]$ is the subgroup of the elements of $H$ that act trivially around the endpoints of $\llbracket 1 \rrbracket$ and $\llbracket n \rrbracket$, a direct computation shows that
\begin{equation}
\label{eq:g0}
g_0 = h [\tau_n, g_1] \text{ for some } h \in [H,H],
\end{equation}
which proves that $g_0 \in [G_n, G_n]$ and ultimately that $[G_n, G_n]$ is the kernel of $\Phi$.

In conclusion, the results shown so far about the commutator subgroup are collected in the following theorem:

\medskip %layout
\begin{theorem}
\label{thm:commutator}
The commutator subgroup $[G_n, G_n]$ of a dendrite rearrangement group $G_n$ is the kernel of the surjective morphism $\Phi = \Pi \times \Delta: G_n \to \mathbb{Z}_2 \oplus \mathbb{Z}$.
In particular, the abelianization of $G_n$ is $\mathbb{Z}_2 \oplus \mathbb{Z}$.
Moreover, $[G_n, G_n]$ is generated by the (infinite) set exhibited above in \cref{eq:commutator}.
\end{theorem}

\begin{remark}
\label{rmk:comm:gen}
Using the identities $\tau_i g_1 = g_1 \tau_i$ for $i = 2, \dots, n-1$ and $[\tau_n, g_1^z] = [\tau_n, g_1]^z$ (which can both be verified with quick computations), the part $g_1^{-z} \tau_2 \tau_i g_1^z$ of the generating set in \cref{eq:commutator} can be replaced by the finite subset
\[ \{ \tau_2 \tau_i, [\tau_n, g_1] \mid i = 3, \dots, n \} \]
This produces the following generating set for $[G_n, G_n]$, which will be useful at the end of the proof of \cref{thm:comm:simple}:
\[ \{ \tau_2 \tau_i, \; [\tau_n, g_1], \; \tau_2^\zeta g_1^{-z} g_0 g_1^z \tau_2^\zeta \, \mid \, \zeta \in \{0,1\}, z \in \mathbb{Z}, i = 3, \dots, n \}. \]
This infinite generating set will be refined to a finite generating set in \cref{thm:comm:fg} at the end of this section.
\end{remark}

\subsection{Simplicity of the Commutator Subgroup}
\label{sub.simplicity.commutator}

It turns out that an analogue of \cref{prop:trans} also applies to the commutator subgroup $[G_n, G_n]$.
In order to prove this, we recall that if a group $G$ acts on a space $X$ and $U \subseteq X$ then the \textbf{rigid stabilizer} of $U$ in $G$ is the subgroup of those elements that are supported in $U$, i.e.,
\[ \mathrm{Rist}_G(U) \coloneq \{ g \in G \mid g(x) = x, \, \forall x \in X \setminus U \} \leq G. \]
The notion of rigid stabilizers will also be useful in \cref{thm:comm:simple}.

\medskip %layout
\begin{proposition}
\label{prop:comm:trans}
Given two finite subsets $\mathcal{F}_1$ and $\mathcal{F}_2$ of $\mathrm{Br}$, any graph isomorphism between $T(\mathcal{F}_1)$ and $T(\mathcal{F}_2)$ can be extended to an element of $[G_n,G_n]$.
In particular, $[G_n, G_n]$ acts transitively on $\mathrm{Br}(D_n)$.
\end{proposition}

\begin{proof}
There is no need to repeat the same argument that was used in \cref{prop:trans}.
Instead, one can reason in the following manner.
Because of the aforementioned \cref{prop:trans}, given $\mathcal{F}_1, \mathcal{F}_2$ and a graph isomorphism $f \colon T(\mathcal{F}_1) \to T(\mathcal{F}_2)$ as in the statement, there exists a $g \in G_n$ that extends $f$.
Now, suppose that $\Phi(g) = (\zeta, z) \in \mathbb{Z}_2 \oplus \mathbb{Z}$.
Since $[G_n, G_n]$ is none other than the kernel of $\Phi$ (\cref{thm:commutator}), our goal is to produce an ``adjustment'' of the element $g$ that nullifies $\Phi$ and is the same as $g$ when restricted to $\mathcal{F}_1$.

Since the minimal subdendrite $[\mathcal{F}_2]$ of $D_n$ that contains $\mathcal{F}_2$ cannot be the entire $D_n$, there exists a branch $B^*$ that is fully contained in the complement of $[\mathcal{F}_2]$.
It is not hard to see that the restriction of $\Phi$ on the rigid stabilizer of any branch of $D_n$ is surjective, so in particular one can find an element $h \in \mathrm{Rist}_{G_n}(B^*)$ such that $\Phi(h) = (\zeta, -z)$.
It follows that the element $h g$ belongs to the commutator subgroup of $G_n$ and is the same as $g$ when restricted to $\mathcal{F}_1$, meaning that it extends the graph isomorphism $f: T(\mathcal{F}_1) \to T(\mathcal{F}_2)$ as desired.
\end{proof}

Using the exact same argument of \cref{clm:dns}, this Proposition allows us to immediately prove that the commutator subgroup $[G_n, G_n]$ is also dense in the full group of homeomorphisms of $D_n$.

\medskip %layout
\begin{corollary}
\label{cor:comm:dns}
The commutator subgroup of the dendrite rearrangement group of $D_n$ is dense in the full group of homeomorphisms of $D_n$.
\end{corollary}

As another consequence of \cref{prop:comm:trans}, we have:

\medskip %layout
\begin{corollary}
\label{cor:comm:trans}
$[G_n, G_n]$ acts transitively on the set of branches of $D_n$.
\end{corollary}

\begin{proof}
Consider two branches $B_1$ and $B_2$.
If $p_1$ and $p_2$ are the branch points at which they branch out, respectively, then by \cref{prop:comm:trans} there exists an element $h \in [G_n, G_n]$ that maps $p_1$ to $p_2$.
Since both $h(B_1)$ and $B_2$ are branches at $p_2$, there is an element $\sigma$ of the alternating subgroup $[K_{p_2}, K_{p_2}] \leq [G_n, G_n]$ (which is known to be transitive) such that $\sigma h (B_1) = B_2$.
Thus, the element $\sigma h$ maps $B_1$ to $B_2$ and belongs to $[G_n, G_n]$, as needed.
\end{proof}

This last Corollary allows us to prove the simplicity of $[G_n, G_n]$ for all $n \geq 4$.

\medskip %layout
\begin{theorem}
\label{thm:comm:simple}
For each $n\geq4$, the commutator subgroup $[G_n, G_n]$ is simple.
\end{theorem}

\begin{proof}
This proof shares its structure with the proofs of the simplicity of the commutator subgroups of the basilica and airplane limit spaces (from \cite{Belk_2015} and \cite{airplane}, respectively), which in turn were inspired by the so-called \textit{Epstein's double commutator trick} from \cite{Epstein}.

Given a non-trivial normal subgroup $N$ of $[G_n, G_n]$, we wish to prove that $N = [G_n, G_n]$.
Let $g \in N \setminus \{1\}$.
Since $g$ is non-trivial, there must be an open subset $U$ of $D_n$ such that $U \cap g(U) = \emptyset$.
Every open subset of $D_n$ contains a cell (see \cref{lem.balls.and.cells}) and every cell contains a branch;
thus, there exists a branch $B$ based at some branch point $p$ that is included in $U$, and so $B \cap g(B) = \emptyset$.

Let $h_1, h_2 \in \mathrm{Rist}_{[G_n, G_n]}(B)$.
Note that $g^{-1} h_1 g \in \mathrm{Rist}_{[G_n, G_n]} \big( g^{-1}(B) \big)$, so
\[ [h_1, g] = h_1^{-1} g^{-1} h_1 g \in \mathrm{Rist}_{[G_n, G_n]} \big( B \cup g^{-1}(B) \big) \]
and $[h_1, g] |_B = h_1^{-1} |_B$.
Thus, since $h_2 \in \mathrm{Rist}_{[G_n, G_n]}(B)$, direct computations of $\big[ [h_1, g]^{-1}, h_2 \big]$ restricted on $B$, $g^{-1}(B)$ and elsewhere show that globally
\[ \big[ [h_1, g]^{-1}, h_2 \big] =[h_1,g] h_2^{-1} [h_1,g]^{-1} h_2 = h_1^{-1} h_2^{-1} h_1 h_2 = [h_1, h_2]. \]
Now, since $[h_1, g]^{-1} \in N$ (because $g \in N \trianglelefteq [G_n, G_n] \ni h_1$) and $h_2 \in [G_n, G_n]$, the double commutator of the previous identity belongs to $N$, and so we have that
\[ [h_1, h_2] \in N \text{ for all } h_1, h_2 \in \mathrm{Rist}_{[G_n, G_n]}(B). \]

We do not have much control over the choice of $B$, but this can be overcome with the aid of \cref{cor:comm:trans}:
for any branch $B^*$, there exists some $h \in [G_n, G_n]$ such that $h(B^*) = B$, and conjugating the previous identity by $h$ yields
\begin{equation}
\label{eq:comm:simple}
[h_1, h_2] \in N \text{ for all } h_1, h_2 \in \mathrm{Rist}_{[G_n, G_n]}(B^*).
\end{equation}
As is shown below, this identity allows us to show that $N$ contains the generating set for $[G_n, G_n]$ exhibited in \cref{rmk:comm:gen}, which will conclude this proof.

First, note that $\tau_2 \tau_i \in N$ for all $i = 3, \dots, n$ because $N$ contains the entire subgroup $[K_{p_0}, K_{p_0}]$.
Indeed, each $3$-cycle $(i j k)$ can be written as $[(ij),(jk)]$.
Since $(ij)$ and $(jk)$ do not belong to $[G_n,G_n]$, we need to replace them by some $h_1 = (ij)h$ and $h_2=(jk)h$ in $[G_n,G_n]$ such that $(ijk) = [(ij), (jk)] = [h_1,h_2]$.
Such an $h$ can be found as follows.
Since $n \geq 4$, there is a fourth branch $B$ at $p_0$ that is distinct from $i$, $j$ and $k$;
choose any branch point $p'$ in that branch and let $h$ be any odd permutation at $p'$ whose support is contained in $B$.
Then $h_1 = (ij)h$ and $h_2 = (jk)h$ are compositions of odd permutations, thus $h_1, h_2 \in [G_n,G_n]$ and they act trivially on some branch that is different from $i, j$ and $k$, so one concludes by \cref{eq:comm:simple}.
(The argument to find $h$ does not work when $n=3$; see \cref{rmk:3:comm:simple,qst:3:comm:simple} below.)

The fact that $[\tau_n, g_1]$ belongs to $N$ can be proved using the idea, previously employed in \cref{prop:comm:trans}, of ``adjusting'' the elements $\tau_n$ and $g_1$ so that they end up in the commutator subgroup, using the fact that $[G_n, G_n] = \mathrm{Ker}(\Phi)$.
More specifically, $[\tau_n, g_1]$ can be written as $[\tau_n \rho, g_1 f]$ for suitable elements $\rho$ and $f$:
the element $\rho$ is a $2$-cycle in a permutation subgroup $K_p$ around some ``remote'' branch point $p$, while $f$ acts on some ``remote'' cell $\llbracket xn \rrbracket$ as $g_1$ does on $\llbracket n \rrbracket$, with $x$ a finite word ending, for example, with $2$, in such a way that the supports of $\rho$ and $f$ are disjoint and each of them has trivial intersection with the supports of $\tau_n$ and $g_1$, and the union $U$ of the supports of the four elements is not the entire $D_n$ (by ``remote'' here we mean that we choose such elements so that their supports are disjoint form those of $\tau_n$ and $g_1$).
Since the union $U$ of the supports of $\rho$, $f$, $\tau_n$ and $g_1$ is a proper closed subset of $D_n$, it is easy to see that it is contained in some branch $B$.
It follows that both $\tau_n \rho$ and $g_1 f$ belong to $\mathrm{Rist}_{[G_n, G_n]}(B)$, so we can conclude that $[\tau_n, g_1] = [\tau_n \rho, g_1 f] \in N$ by \cref{eq:comm:simple}.

Now that we know that $N$ contains $[\tau_n, g_1]$, because of \cref{eq:g0} it suffices to prove that $[H,H] \leq N$ in order to show that $g_0 \in N$.
This is done by recalling that the commutator subgroup of Thompson's group $F$ is the subgroup of $F$ that acts trivially around the endpoints of the unit interval, and that it is simple and non-abelian, so the commutator of $[H,H]$ is $[H,H]$ itself.
Then an arbitrary element of $[H,H]$ can be written as the product of commutators of elements of $[H,H]$, each of which must be acting trivially on some neighborhood of the endpoints of $\llbracket 1 \rrbracket$ and $\llbracket n \rrbracket$.
Then these elements act trivially in some small enough branch $B$ and an application of \cref{eq:comm:simple} allows us to conclude that $[H,H] \leq N$, and so that $g_0 \in N$.

Finally, \cref{eq:g0} tells us that $g_0 = h [\tau_n, g_1] = [l_1, r_1] \dots [l_m,r_m] [\tau_n, g_1]$ for some $l_i, r_i \in [H,H]$ (using again the fact that the commutator subgroup of $[H,H]$ is $[H,H]$ itself).
With this identity, the remaining elements of the generating set from \cref{rmk:comm:gen} all belong to $N$ because they are
\[ \tau_2^\zeta g_1^{-z} g_0 g_1^z \tau_2^\zeta = g_0^{g_1^z\tau_2^{\zeta}} = \big[ l_1^{g_1^z\tau_2^{\zeta}}, r_1^{g_1^z\tau_2^{\zeta}} \big] \dots \big[ l_m^{g_1^z\tau_2^{\zeta}}, r_m^{g_1^z\tau_2^{\zeta}} \big] \big[ \tau_n^{g_1^z\tau_2^{\zeta}}, g_1^{g_1^z\tau_2^{\zeta}} \big], \]
where each of the commutators satisfy \cref{eq:comm:simple}.
\end{proof}

\begin{remark}
\label{rmk:3:comm:simple}
As noted in the proof of \cref{thm:comm:simple}, our argument to prove that $\tau_2 \tau_i$ belongs to $N$ for $i = 3, \dots, n$ fails for $n=3$.
The reason is simply that, with only $3$ branches at each branch point, the elements $h_i = (ij)h$ and $h_j = (jk)h$ (with $\{i, j, k\} = \{1, 2, 3\}$ and $h$ as described in the proof) do not belong to a common rigid stabilizer of any branch, so one cannot apply \cref{eq:comm:simple}.
Every other step of the proof works, so it would suffice to find a different way to show that $\tau_2 \tau_3 \in N$ in order to prove that $[G_3, G_3]$ is simple.
However, if the commutator subgroup is not simple, the rest of the proof of \cref{thm:comm:simple} imply that a non-trivial proper normal subgroup $N$ of $[G_3, G_3]$ must include the elements $[\tau_3, g_1]$ and $\tau_2^\zeta g_1^{-z} g_0 g_1^z \tau_2^\zeta$ along with the subgroup $[H, H]$, and so it can be argued that $N$ must also contain $[H_A, H_A]$ for every rational arc $A$, by the normality of $N$ and the transitive properties of the commutator subgroup.
\end{remark}

This question is left unanswered:

\begin{question}
\label{qst:3:comm:simple}
Is $[G_3, G_3]$ simple?
\end{question}

\subsection{A Finite Generating Set for the Commutator Subgroup}

As a final result about the commutator subgroup of $G_n$, we compute a finite generating set.
The proof follows the same strategy of \cite[Theorem 6.12]{airplane} and also makes use of the fact that the commutator subgroup of a Thompson subgroup $H_A$ is the set of elements of $H_A$ that act trivially on some neighborhood of each of the two endpoints of the arc $A$.

Let $q_0 = \llbracket \overline{n} \rrbracket$ and $A_\star = [p_0,q_0]$ (which is the ``right arc'' attached to $p_0$) and denote by $q_1$ the point $\llbracket n i 1 \overline{n} \rrbracket$ (which is the ``middle point'' of $A_\star$).
We start from the generating set from \cref{rmk:comm:gen}, and we will show that those elements all belong to the subgroup
\[ G_\star = \langle \tau_2 \tau_3, \dots, \tau_2 \tau_n, [\tau_n, g_1], g_0, g_0^{\tau_2}, g_\star \rangle, \]
where $g_\star$ is the element of $[H,H]$ that acts on the arc $[p_0,q_1]$ (the left half of the right arc $A_\star = [p_0, q]$) as $g_0$ does on the entire arc $A_0$.
By definition $G_\star$ contains $\tau_2 \tau_i$ for each $i = 1, \dots, n$ and $[\tau_n, g_1]$, so it suffices to show that the elements $\tau_2^\zeta g_1^{-z} g_0 g_1^z \tau_2^\zeta$ all belong to $G_\star$.

\begin{figure}
\centering
\centering
\begin{tikzpicture}
    \node[vertex] (c) at (0,0) {};
    \node[vertex] (l) at (-2,0) {};
        \node[vertex] (lc) at (-1,0) {};
        \node[vertex] (lt) at (-1,.8) {};
    \node[vertex] (t) at (0,2) {};
    \node[vertex] (r) at (2,0) {};
        \node[vertex] (rc) at (1,0) {};
        \node[vertex] (rt) at (1,.8) {};
    \draw[curvy] (rc) to (c);
    \draw[curvy] (c) to (lc);
    \draw[red] (lc) to (l);
    \draw[Green] (rc) to (r);
    \draw[dotted] (lc) to (lt);
    \draw[dotted] (rc) to (rt);
    \draw[dotted] (c) to (t);
    \draw[->] (2.5,0) to (3,0);
    \begin{scope}[xshift=5.5cm]
    \node[vertex] (c) at (0,0) {};
    \node[vertex] (l) at (-2,0) {};
    \node[vertex] (t) at (0,2) {};
    \node[vertex] (r) at (2,0) {};
            \node[vertex] (rrc) at (1.6,0) {};
            \node[vertex] (rrt) at (1.6,.3) {};
        \node[vertex] (rc) at (1,0) {};
        \node[vertex] (rt) at (1,.8) {};
    \draw[curvy] (rrc) to (rc);
    \draw[curvy] (rc) to (c);
    \draw[red] (c) to (l);
    \draw[Green] (rrc) to (r);
    \draw[dotted] (c) to (t);
    \draw[dotted] (rc) to (rt);
    \draw[dotted] (rrc) to (rrt);
    \end{scope}
\end{tikzpicture}
\caption{The general shape of the elements $g_1^{-z} g_0 g_1^z$, with an action on the curvy lines that depends on $z$.}
\label{fig:comm:fg}
\end{figure}
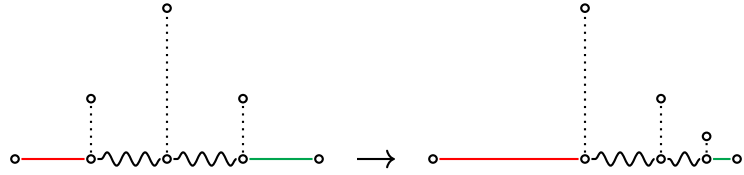

For $\zeta = 0$, observe that every element $g_1^{-z} g_0 g_1^z$ has the shape sketched in \cref{fig:comm:fg}, from which we deduce that it belongs to $[H,H] g_0$.
Now, we can prove that $[H,H] \leq G_\star$ with the following strategy:
if $h \in [H,H]$, then $h^{g_0^k}$ belongs to the commutator subgroup of $H_{A_\star}$ for a negative integer $k$ that is small enough.
But $[H_{A_\star}, H_{A_\star}]$ is included in the commutator subgroup of $\langle [\tau_n, g_1], g_\star \rangle$, as the two elements act on $A_\star$ as the two generators of $F$ do on $[0,1]$ (this step is discussed in more detail in \cite[Theorem 6.12]{airplane}).
Since $[\tau_n,g_1], g_\star \in G_\star$, we have that $[H_{A_\star}, H_{A_\star}]$ is included in $G_\star$.
Thus, $[H,H] \leq G_\star$, so $g_1^{-z} g_0 g_1^z \in G_\star$ for all $z \in \mathbb{Z}$.

As for the case $\zeta = 1$, consider the arc $A_0' = [t,q]$, where $t = \llbracket 2 \overline{n} \rrbracket$ and, as before, $q = \llbracket ni1\overline{n} \rrbracket$.
One can prove that the elements $\tau_2 g_1^{-z} g_0 g_1^z \tau_2$ (similar to \cref{fig:comm:fg}, but replacing the left central branch with the top one) all belong to $G_\star$ essentially by replacing $g_0$ with $g_0^{\tau_2}$ (that acts on $A_0'$ as $g_0$ does on $A_0$) and $H$ with $H_{A_0'}$ in the previous paragraph.
This simply means switching the first and second branches at $p_0$, which preserves the arc $A_\star$ that we used as a pivot, so one ``pushes'' the action of $[H_{A_0'}, H_{A_0'}]$ inside $A_0$ by conjugating by negative powers of $g_0^{\tau_2}$ and obtains that $[H_{A_0'}, H_{A_0'}] \leq G_\star$.

Ultimately, this shows that $G_\star \geq [G_n, G_n]$.
The converse is clear, so we ultimately have the following.

\medskip %layout
\begin{theorem}
\label{thm:comm:fg}
The commutator subgroup of a dendrite rearrangement group $G_n$ is finitely generated, for every $n \geq 3$.
\end{theorem}

%%%%%%%%%%%%%%%%%%%%%%%%%

\section[Further Properties of \texorpdfstring{$G_n$}{Gn}]{Further Properties of Dendrite Rearrangement Groups}
\label{sec:properties}

This section is a collection of additional results about dendrite rearrangement groups.
Here we provide information about their finite subgroups (\cref{sub:fin:subgps}), we prove that they are pairwise non-isomorphic (\cref{sub:distinct}), we show how they relate to the other Thompson groups (\cref{sub:relations:Thompson}), we see that their conjugacy problem is solvable (\cref{sub:conj}), and finally we prove that they are not invariably generated (\cref{sub:IG}).

\subsection{Finite Subgroups}
\label{sub:fin:subgps}

Generalizing the ideas used in \cite[Proposition 2.10]{BF19} (which is about the Vicsek rearrangement group, see \cref{sub:Vic}), we obtain the following result.

\medskip %layout
\begin{proposition}
\label{prop:finite:subgroups}
Each finite subgroup of $G_n$ is a subgroup of the group of automorphisms of some finite $n$-regular tree, and all such subgroups for all finite $n$-regular trees are achieved.

In particular, each finite subgroup of $G_n$ has order
\[ p_0^{\varepsilon} p_1^{k_1} \dots p_m^{k_m}, \]
where $p_1, \dots, p_m$ are the prime divisors of $(n-1)!$, $k_i \in \mathbb{N}$ (possibly zero), $\varepsilon \in \{0,1\}$ and $p_0 = n$ if it is prime and it is $1$ otherwise.
Every such order is achieved by some abelian subgroup.
\end{proposition}

\begin{proof}
First, let $G$ be a finite subgroup of $G_n$.
By \cref{thm.finite.subgroups}, every finite subgroup of a rearrangement group is the subgroup of the automorphism group of some graph expansion of the edge replacement system, so $G \leq \mathrm{Aut}(T)$ for some finite $n$-regular tree $T$.
It is easy to see that every finite $n$-regular tree can be realized as a graph expansion of $\mathcal{D}_n$.

The decomposition of $\mathrm{Aut}(T)$ described in this paragraph is well-known in the literature and is often called \textit{wreath recursion} (see for example \cite[Section 1.4]{NekSSG}).
It is known that every finite tree has a unique \textit{center},
which is either a vertex or an edge that is fixed by every automorphism (see for example \cite[Theorem 3.46]{GroupsGraphsTrees}).
In any case, $\mathrm{Aut}(T)$ can be decomposed as a finite sequence of semidirect products (where the action is by permutations of the components of the direct product) in the following way:
\[ \mathrm{Aut}(T) = A_0 \ltimes (A_1 \times \dots \times A_n), \]
where $A_0 \leq \mathrm{Sym}(n)$ or $A_0 \leq \mathrm{Sym}(2)$ (depending on whether the center is a vertex or an edge) and each $A_i$ is
\[ A_i = A_{i,0} \ltimes (A_{i,1} \times \dots \times A_{i,n-1}), \]
where $A_{i,0} \leq \mathrm{Sym}(n-1)$
and each $A_{i,j}$ is built in this very same way, eventually ending with the trivial group.
In particular, this implies that the order of $\mathrm{Aut}(T)$ (and of each of its subgroups) is of the type exhibited in the statement.

Now, given a number of the form exhibited in the statement, we can write it as a product $q_1 \dots q_k$ of primes $q_i$'s that need not be pairwise distinct.
Without loss of generality, we order them so that $q_1 = n = p_0$ if $n$ is prime and it is one of the factors.
Consider a rooted tree $\mathfrak{T}$ whose root has $q_1$ children whose $i$-th level vertices all have $q_{i+1}$ children (note in particular that each vertex has degree that is at most $n$, by the choice $q_1 = n$ if it is a prime factor).
Denote the root by the empty word and each $i$-th child of $w$ by $wi$.
The abelian group $C_{q_1} \times \dots \times C_{q_k}$ acts faithfully on this tree with the diagonal action defined as follows:
an element $(z_1, \dots, z_k)$ maps each vertex $x_1 \dots x_k$ to $(x_1 + q_1) \dots (x_k + q_k)$.
Since each vertex of $\mathfrak{T}$ has degree that is at most $n$, this action can be extended to a faithful action on the natural completion of $\mathfrak{T}$ to an $n$-regular tree.
In particular, this group embeds into the automorphism group of an $n$-regular tree, so we are done.
\end{proof}

\subsection{Pairwise Distinction}
\label{sub:distinct}

Using \cref{prop:finite:subgroups} from the previous subsection, one can immediately distinguish two dendrite rearrangement groups $G_n$ and $G_m$ using the orders of their finite subgroups whenever $n$ or $m$ is at most $4$ (they can actually be distinguished in this way for all $n, m \leq 9$).
For $n \geq 5$, we can use the same argument that is mentioned in subsection 2.3 of \cite{BF19} for distinguishing the Vicsek rearrangement groups.
For the sake of completeness, here we explain this argument in full detail for dendrite rearrangement groups.

\medskip %layout
\begin{proposition}
\label{prop:alt}
For every $n \geq 3$, the dendrite rearrangement group $G_n$ contains isomorphic copies of the alternating group $A_n$ but not of $A_{n+1}$.
\end{proposition}

\begin{proof}
Thanks to \cref{prop:finite:subgroups} we know what the finite subgroups of $G_n$ are, from which it is immediate to find copies of $A_n$.
More explicitly, it is clear that $A_n$ acts faithfully on the base graph of $\mathcal{D}_n$ itself;
or, alternatively, one can recover this copy of $A_n$ as a subgroup of the permutation subgroup $K_{p_0}$ from \cref{lem:perm}.

If we suppose that $A_{n+1}$ embeds into $G_n$, then by \cref{prop:finite:subgroups} we must have that $A_{n+1}$ acts faithfully on a finite $n$-regular tree.
This is a contradiction because of the following:

\medskip %layout
\begin{claim}
For $m \geq 5$, if the alternating group $A_m$ acts faithfully on a tree $\mathcal{T}$ then there exists a vertex of $\mathcal{T}$ whose degree is at least $m$.
\end{claim}

To prove this Claim, assume that $A_m \simeq H \leq \mathrm{Aut}(\mathcal{T})$.
As noted in the previous \cref{sub:fin:subgps}, $H$ must fix a vertex or an edge of $\mathcal{T}$ (see for example \cite[Theorem 3.46]{GroupsGraphsTrees}).
Up to a barycentric subdivision of the fixed edge, we can assume that $H$ is fixing a vertex $r$, so we can think of $\mathcal{T}$ as if it were a rooted tree.
Then observe that, for each $k \in \mathbb{N}$, the subgroup consisting of the elements of $H$ that act trivially on the first $k$ levels of $\mathcal{T}$ is normal in $H$.
Since $H$ is simple, this implies that for every $k$ either $H$ acts trivially or faithfully on the first $k$ levels.
Thus, since $H$ is not trivial, there exists a $k$ such that the action of $H$ on the first $k-1$ levels is trivial whereas the action on level $k$ is faithful.
Note that then $H$ must embed into a direct product of symmetric groups, i.e.,
\[ H \leq \mathrm{Sym}(m_1) \times \dots \times \mathrm{Sym}(m_l), \]
where the $m_i$'s are the number of children of the $l$ vertices at level $k$.
We now prove that this cannot happen unless there exists a $j$ such that $m_j \geq m$.
This will conclude the proof, since each $m_i$ is either the degree of a vertex of $\mathcal{T}$ or the degree minus $1$, depending on whether $k=0$ or $k>0$.

Denote by $\mathrm{pr}_i$ the $i$-th projection in the aforementioned direct product.
Since $H$ is simple, each $\mathrm{pr}_i(H)$ is either trivial or isomorphic to $H \simeq A_m$.
These cannot all be trivial or $H$ would be trivial itself, so there is a $j$ such that $A_m \simeq \mathrm{pr}_j(H) \leq \mathrm{Sym}(m_j)$.
It is known that $\mathrm{Sym}(m_j)$ cannot contain a copy of $A_m$ unless $m \leq m_j$, so we are done.
\end{proof}

As an immediate consequence of \cref{prop:alt}, we have the following.

\medskip %layout
\begin{theorem}
\label{thm:non:isomorphic}
$G_n \simeq G_m$ if and only if $n = m$.
\end{theorem}

\begin{remark}
\label{rmk:rubin}
Very often results of this kind can be proved as applications of theorems of reconstructions of topological spaces from groups acting on them (i.e., given a group $G$ acting ``nicely'' on two ``nice'' topological spaces $X$ and $Y$, one finds an equivariant homeomorphism between $X$ and $Y$).
The most prominent example of this is the celebrated paper \cite{Rubin}, which features 4 families of group actions that allow a reconstruction of the underlying space.
Of these, the result named \textit{Theorem 2} in the abstract of \cite{Rubin}, which is often just called \textbf{Rubin's Theorem} and is also discussed and proved more succinctly in \cite{ShortRubin}, applies to many Thompson groups and other ``rich'' subgroups of homeomorphism groups.

However, we remark here that neither Rubin's Theorem (nor any of the other 3 from \cite{Rubin}, to the author's knowledge) applies to $G_n$, nor to any group acting by homeomorphisms on a Wa\.zewski dendrite $D_n$.
Indeed, let $U$ be the open set $\llparenthesis 11 \rrparenthesis \cup \dots \cup \llparenthesis n1 \rrparenthesis \cup \{p_0\}$ (which is essentially obtained by taking the inner ``halves'' of the $n$ branches at $p_0$).
Every element of $\mathrm{Rist}_{G_n}(U)$ must fix setwise each of the cells $\llbracket 11 \rrbracket, \dots, \rrbracket n1 \llbracket$, so it must also fix $p_0$, which prevents $\mathrm{Rist}_{G_n}(U)$ from satisfying the conditions of Rubin's Theorem.

It is possible to verify for $G_n$ the hypotheses of \cite[Theorem 0.2]{Rubin} (as was done in \cite[Proposition 7.6]{DM19} for $\mathbb{H}_n$), but these hypotheses alone do not seem to yield any reconstruction theorem, as far as the author knows, since it is not clear whether in this case the structures that Rubin denotes by $HR(G,X)$ are homeomorphic to $X$.
\end{remark}

Before moving on to a comparison with other Thompson groups, we highlight the following fact.

\medskip %layout
\begin{proposition}
\label{prop:inclusion}
$G_n < G_{n+1}$, for all $n \geq 3$.
\end{proposition}

\begin{proof}
Simply consider the subgroup of $G_{n+1}$ generated by its Thompson subgroup $H$ along with the subgroup of $K_{p_0}$ that fixes the cell $\llbracket n+1 \rrbracket$, which is the subgroup of those rearrangements that ``rigidly drag'' every $(n+1)$-th branch along with its branch point.
More precisely, each cell $\llbracket x \rrbracket$ corresponding to a word $x$ that ends with $n+1$ and does not feature $n+1$ anywhere else is mapped canonically to a cell $\llbracket y \rrbracket$ with this same property.
This is an isomorphic copy of $G_n$ inside $G_{n+1}$.
\end{proof}

\subsection{Relations with other Thompson Groups}
\label{sub:relations:Thompson}

As seen in \cref{sub:thomp}, dendrite rearrangement groups feature infinitely many isomorphic copies of Thompson's group $F$.
Here we show how the $G_n$'s compare to the bigger siblings of $F$.

\medskip %layout
\begin{proposition}
\label{prop:dendrite:Thompson:comparison}
For every $n \geq 3$, each dendrite rearrangement group $G_n$ embeds in Thompson's group $V$, but $G_n$ does not embed into Thompson's group $T$ nor $T$ embeds into $G_n$ (so, in particular, $V$ also does not embed into $G_n$).
\end{proposition}

\begin{proof}
An embedding of $G_n$ into $V$ is described in \cref{thm.embedding.into.V} for every rearrangement group.
To achieve it in the specific case of dendrite rearrangement groups, one can simply consider a graph pair diagram for $G_n$ and ``unglue'' all of the edges, considering domain and range graphs consisting of edges each of which is not adjacent to any other.
The result is a valid graph pair diagram for the Higman-Thompson group $V_{n,1}$ (defined in \cref{sub.higman.thompson.groups}), since every permutation of the edges is allowed in $V_{n,1}$.
This prompts a natural embedding of $G_n$ inside $V_{n,1}$, which in turn embeds in $V$ by \cite[Theorem 7.2]{higman1974finitely}.

As for Thompson's group $T$, it is known that its finite subgroups are cyclic and that every finite cyclic group embeds in $T$.
Then \cref{prop:finite:subgroups} implies that $G_n$ cannot embed into $T$ (as clearly, for every $n \geq 3$, there are automorphism groups of finite $n$-regular trees that are not cyclic) and also that $T$ cannot embed into $G_n$ (otherwise $G_n$ would have to feature cyclic subgroups of every order, but it does not by the aforementioned \cref{prop:finite:subgroups}).
\end{proof}

Since both the rearrangement groups of the basilica and of the airplane limit spaces (respectively $T_B$ from \cite{Belk_2015} and $T_A$ from \cite{BF19} and \cite{airplane}) contain copies of $T$ and embed into $T$, we can also conclude the following.

\medskip %layout
\begin{corollary}
For every $n \geq 3$, $T_B$ and $T_A$ do not embed into $G_n$ nor does $G_n$ into $T_B$ or $T_A$. 
\end{corollary}

It appears that the groups $T_B$ and $T_A$, along with Thompson's group $T$, have more in common with the orientation-preserving dendrite rearrangement groups mentioned in \cref{sub:orientation} below.
Instead, dendrite rearrangement groups $G_n$ do not preserve an orientation of the limit space $D_n$ on which they act, which is arguably what separates these groups from $T_B$ and $T_A$.

The only other rearrangement group mentioned in the literature that is likely to contain a copy of or embed into some dendrite rearrangement group is the Vicsek rearrangement group, along with its generalizations (\cite[Example 2.1]{BF19}, see also \cref{sub:Vic}).
See \cref{sub:Vic} below for more about this comparison.

\subsection{A Solution to the Conjugacy Problem}
\label{sub:conj}

We say that the \textbf{conjugacy problem} of a group $G$ is solvable if there exists an algorithm that, given two elements of $G$, infallibly decides in finite time whether those elements belong to the same conjugacy class.

A sufficient condition for the conjugacy problem to be solvable, together with an actual method for solving it, will be described in detail in \cref{cha.conjugacy}.
This sufficient condition is essentially the following:
graph reductions (which are the opposite of graph expansions of the edge replacement system) need to be \textbf{confluent}, by which we mean that if $\Gamma \longrightarrow \Gamma_1$ and $\Gamma \longrightarrow \Gamma_2$ are distinct finite sequences of graph reductions, then there exist $\Gamma^*$ and finite sequences of reductions $\Gamma_1 \longrightarrow \Gamma^*$ and $\Gamma_2 \longrightarrow \Gamma^*$.
In essence, one needs to prove that, when two distinct sequences of graph reductions of $\Gamma$ produce different graphs $\Gamma_1$ and $\Gamma_2$, then these two graphs can be reduced to a common $\Gamma^*$.

Dendrite edge replacement systems $\mathcal{D}_n$ have confluent graph reductions.
More precisely, in this case distinct graph reductions are all disjoint.
Indeed, keeping in mind that each graph reduction is uniquely identified by a vertex of degree $n$ and its $n$ incident edges (which is simply the replacement graph), it is easy to see that no two distinct graph reductions are allowed to involve common edges.
Thus, one has the following.

\medskip %layout
\begin{proposition}
Dendrite rearrangement groups have solvable conjugacy problem.
\end{proposition}

We note here that conjugacy in the (topological) group $\mathrm{AAut}(\mathcal{T}_{d,k})$ of almost automorphism of trees has been studied recently in \cite{Goffer2019ConjugacyAD}.
That paper shows that whether two hyperbolic elements of $\mathrm{AAut}(\mathcal{T}_{d,k})$ are conjugate can be established by studying elements of the Higman-Thompson group $V_{d,k}$ that ``approximate'' the two almost automorphisms.
Since we have seen that dendrite rearrangement groups ``approximate'' the full groups of homeomorphisms (\cref{thm:dense}), inspired by the use that \cite{Goffer2019ConjugacyAD} made of the strand diagrams developed in \cite{BM14}, it is natural to ask the following.

\begin{question}
Can the technology of strand diagrams for rearrangement groups developed in \cref{cha.conjugacy} be used to study conjugacy in $\mathrm{Homeo}(D_n)$?
\end{question}

\subsection{Absence of Invariable Generation}
\label{sub:IG}

We say that a group $G$ is \textbf{invariably generated} if it admits a generating set $S$ that still generates $G$ even after one has conjugated each of the elements of $S$ by elements of $G$.
More precisely:
\[ G = \langle s^{g_s} \mid s \in S \rangle \text{ for every choice of } \{g_s \mid s \in S \} \subseteq G. \]
Thompson's group $F$ is invariably generated, whereas $T$ and $V$ are not, as shown in \cite{Gelander2016InvariableGO}.
See \cref{cha.IG} for more information on this property.

A rearrangement group is \textbf{weakly cell-transitive} if, given an arbitrary cell $C$ and an arbitrary proper union $K$ of finitely many cells, there exists a rearrangement that maps $K$ inside $C$ (this property is discussed further in \cref{cha.IG}).
The main result of \cref{cha.IG} is that, if a rearrangement group is weakly cell-transitive, then it is not invariably generated, nor is its commutator subgroup.

\medskip %layout
\begin{proposition}
\label{prop:wct}
Dendrite rearrangement groups are weakly cell-transitive.
\end{proposition}

\begin{proof}
Given a cell $C$ of $D_n$ and a union $K$ of finitely many cells we need to show that $G_n$ can map $K$ inside $C$.
Note that every cell must contain a branch (in fact, it contains infinitely many) and that the union of any finite amount of cells is fully included in some branch.
Then it suffices to know that, given two arbitrary branches $B_1$ and $B_2$, there exists an element $g$ of $G_n$ such that $g(B_1) \subseteq B_2$.
But then, thanks to \cref{cor:comm:trans}, we are already done.
\end{proof}

This prompts the following application of the main result from \cref{cha.IG}, and gives an example of behaviour of $G_n$ that makes it appear more similar to Thompson groups $T$ and $V$ than to $F$ (see \cref{sub:amen} below for another example of this).

\medskip %layout
\begin{proposition}
\label{prop:IG}
Dendrite rearrangement groups and their commutator subgroups are not invariably generated.
\end{proposition}

\subsection{Existence of Free Subgroups}
\label{sub:amen}

As the result in the previous subsection, the following is another example of a behaviour of dendrite rearrangement groups that resemble that of Thompson groups $T$ and $V$ rather than $F$.

\medskip %layout
\begin{proposition}
Dendrite rearrangement groups contain non-abelian free groups.
In particular, they are not amenable.
\end{proposition}

Below are two easy proofs of this statement.
The first makes use of tools about rearrangement groups, while the second uses tools that are about groups acting on dendrites.
However, both proofs are ultimately based on the same ping-pong argument.

\begin{enumerate}
    \item
    In the previous section we have proved that dendrite rearrangement groups are weakly cell-transitive (\cref{prop:wct}).
    As will be noted later in \cref{cor.wct.free.subgroups}, weakly cell-transitive rearrangement groups contain non-abelian free subgroups, so we are done.
    \item
    A group action on a topological space is \textbf{minimal} if each of its orbits is dense.
    The main result of \cite{Shi12} states that, if a group $G$ acts minimally on a non-degenerate dendrite, then it contains a non-abelian free subgroup.
    This applies to $G_n$ because its action on $D_n$ is weakly cell-transitive (\cref{prop:wct}), and such actions of rearrangement groups are minimal, as will be shown in \cref{prop.orbit.dense}.
    (One can also prove that the action is minimal by showing that the action of a dendrite rearrangement group does not have finite orbits and applying \cite[Theorem 5.2]{GM19}.)
\end{enumerate}

%%%%%%%%%%%%%%%%%%%%%%%%%

\section{Other Dendrite Rearrangement Groups}
\label{sec:generalizations}

The main result of this section is the density of the airplane rearrangement group $T_A$ in the orientation-preserving subgroup $\mathbb{H}_\infty^+$ of the full homeomorphism group of the infinite-order Wa\.zewski dendrite $D_\infty$ (\cref{thm:dns:A}).
This will be easily proved with the same strategy used in \cref{sec:dns}, as an application of to the results of \cite{airplane}.

The rest of the section is devoted to introducing multiple possible generalizations of dendrite rearrangement groups, of which $T_A$ is one example.
We will first briefly discuss the Vicsek rearrangement groups from \cite{BF19} (\cref{sub:Vic}), then we will introduce an orientation-preserving version of dendrite rearrangement groups (\cref{sub:orientation}) which will prompt the aforementioned result about the airplane rearrangement group $T_A$ (\cref{sub:Air}), and finally we will exhibit edge replacement systems for the generalized Wa\.zewski dendrites $D_S$ (\cref{sub:gen:den}).
With the exception of the results about $T_A$, this section mostly has a conversational and less formal tone.

\subsection{The Vicsek rearrangement groups}
\label{sub:Vic}

\begin{figure}
\centering
\begin{tikzpicture}
    \node at (-2.2,0) {$\Gamma =$};
    \node[vertex] (C) at (0,0) {};
    \node[vertex] (L) at (-1.5,0) {};
    \node[vertex] (T) at (0,1.5) {};
    \node[vertex] (R) at (1.5,0) {};
    \node[vertex] (B) at (0,-1.5) {};
    \draw[edge] (C) to (L);
    \draw[edge] (C) to (T);
    \draw[edge] (C) to (R);
    \draw[edge] (C) to (B);
    \begin{scope}[xshift=7.5cm]
    \node at (-3.2,0) {$R =$};
    \node[vertex] (C) at (0,0) {};
    \node[vertex] (I) at (-2.5,0) {}; \draw (-2.5,0) node[above]{$\iota$};
    \node[vertex] (L) at (-1.25,0) {};
    \node[vertex] (T) at (0,1.25) {};
    \node[vertex] (R) at (1.25,0) {}; \draw (1.25,0) node[above]{$\tau$};
    \node[vertex] (B) at (0,-1.25) {};
    \draw[edge] (C) to (L);
    \draw[edge] (C) to (T);
    \draw[edge] (C) to (R);
    \draw[edge] (C) to (B);
    \draw[edge] (I) to (L);
    \end{scope}
\end{tikzpicture}
\caption{The Vicsek edge replacement system for $n=4$.}
\label{fig:vicsek:rep:sys}
\end{figure}
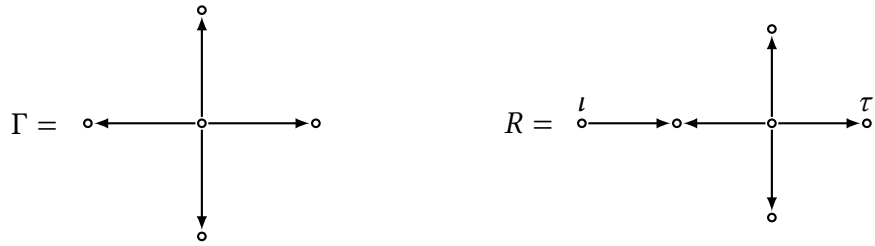

Introduced in \cite[Example 2.1]{BF19}, the Vicsek edge replacement systems are the same of the dendrite edge replacement systems $\mathcal{D}_n$ except for having an additional edge originating from the initial vertex.
For example, \cref{fig:vicsek:rep:sys} depicts the Vicsek edge replacement system for $n=4$.

It is clear that each Vicsek edge replacement system has the same Wa\.zewski dendrite $D_n$ as a limit space (although dendrite edge replacement systems $\mathcal{D}_n$ are arguably more natural edge replacement systems for these limit spaces, while the Vicsek edge replacement system reflects some of the metric aspects of the Vicsek fractals).
It seems likely that the $n$-th Vicsek rearrangement group is generated by a copy of the Higman-Thompson group $F_3$ along with the permutation group $\mathrm{Sym}(n)$ in a similar way in which $G_n$ is generated by $F$ and $\mathrm{Sym}(n)$ by \cref{thm:gen}.
It is known that $F$ and $F_3$ are non-isomorphic groups, as they have different abelianization, which was proved in \cite{Brown}.
However, it does not seem immediately clear whether or not the $n$-th Vicsek rearrangement group is isomorphic to $G_n$, since Rubin's theorem cannot be immediately applied (as discussed in \cref{rmk:rubin}) and the finite subgroups (\cref{sub:fin:subgps}) are likely the same.

Additionally, one may be inspired by this example to build edge replacement systems with even more edges between the initial vertex $\iota$ and the branch point of the replacement graph, which may result in a rearrangement group that has some Higman-Thompson group $F_m$ acting on each EE-arc, in place of $F$ or $F_3$.

We did not investigate these questions any further.

\subsection{Orientation-Preserving Dendrite Rearrangement Groups}
\label{sub:orientation}

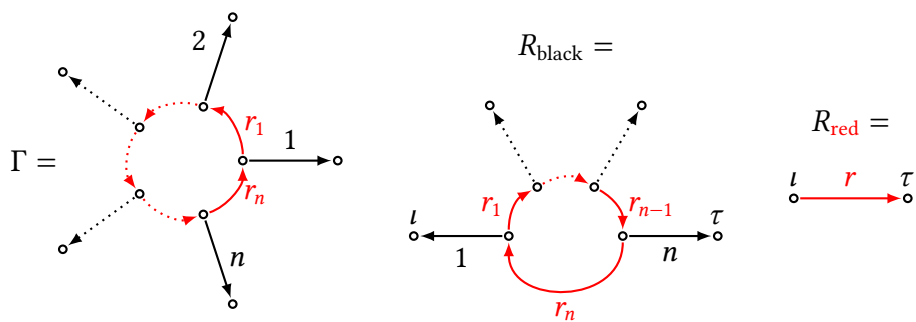
\begin{figure}
\centering
\begin{tikzpicture}
    \node at (-2,0) {$\Gamma =$};
    \node[vertex] (c1) at (0:.75) {};
    \node[vertex] (c2) at (72:.75) {};
    \node[vertex] (c3) at (144:.75) {};
    \node[vertex] (cn-1) at (216:.75) {};
    \node[vertex] (cn) at (288:.75) {};
    \node[vertex] (1) at (0:2) {};
    \node[vertex] (2) at (72:2) {};
    \node[vertex] (3) at (144:2) {};
    \node[vertex] (n-1) at (216:2) {};
    \node[vertex] (n) at (288:2) {};
    \draw[edge] (c1) to node[above]{$1$} (1);
    \draw[edge] (c2) to node[above left]{$2$} (2);
    \draw[edge,dotted] (c3) to (3);
    \draw[edge,dotted] (cn-1) to (n-1);
    \draw[edge] (cn) to node[right]{$n$} (n);
    \draw[edge,red] (c1) to[out=90,in=-12,looseness=.85] node[right]{$r_1$} (c2);
    \draw[edge,red,dotted] (c2) to[out=162,in=54,looseness=.85] (c3);
    \draw[edge,red,dotted] (c3) to[out=234,in=126,looseness=.85] (cn-1);
    \draw[edge,red,dotted] (cn-1) to[out=306,in=198,looseness=.85] (cn);
    \draw[edge,red] (cn) to[out=18,in=270,looseness=.85] node[right]{$r_n$} (c1);
    \begin{scope}[xshift=5cm,yshift=-1cm]
    \node at (0,2.5) {$R_{\text{black}} =$};
    \node[vertex] (b1) at (180:.75) {};
    \node[vertex] (b2) at (120:.75) {};
    \node[vertex] (bn-1) at (60:.75) {};
    \node[vertex] (bn) at (0:.75) {};
    \node[vertex] (r1) at (180:2) {}; \draw (180:2) node[above]{$\iota$};
    \node[vertex] (r2) at (120:2) {};
    \node[vertex] (rn-1) at (60:2) {};
    \node[vertex] (rn) at (0:2) {}; \draw (0:2) node[above]{$\tau$};
    \draw[edge] (b1) to node[below]{$1$} (r1);
    \draw[edge,dotted] (b2) to (r2);
    \draw[edge,dotted] (bn-1) to (rn-1);
    \draw[edge] (bn) to node[below]{$n$} (rn);
    \draw[edge,red] (bn-1) to[out=-30,in=90,looseness=1] node[right]{$r_{n-1}$} (bn);
    \draw[edge,red,dotted] (b2) to[out=30,in=150,looseness=1] (bn-1);
    \draw[edge,red] (b1) to[out=90,in=210,looseness=1] node[left]{$r_1$} (b2);
    \draw[edge,red] (bn) to[out=270,in=270,looseness=1.45] node[below]{$r_n$} (b1);
    \end{scope}
    \begin{scope}[xshift=8.75cm,yshift=-.5cm]
    \node at (0,1) {$R_{\text{\textcolor{red}{red}}} =$};
    \node[vertex] (ir) at (-.75,0) {}; \draw (-.75,0) node[above]{$\iota$};
    \node[vertex] (tr) at (.75,0) {}; \draw (.75,0) node[above]{$\tau$};
    \draw[edge,red] (ir) to node[above]{$r$} (tr);
    \end{scope}
\end{tikzpicture}
\caption{A schematic depiction of the oriented dendrite edge replacement system $\mathcal{D}_n^+$.}
\label{fig:rep:sys:+}
\end{figure}

The \textbf{oriented dendrite edge replacement system} $\mathcal{D}_n^+$, depicted in \cref{fig:rep:sys:+}, is the 2-colored edge replacement system whose base graph and black replacement graph are the same graph consisting of a directed closed path of $n$ red edges and $n$ vertices, along with a black edge originating from each of these $n$ vertices, each terminating at one of $n$ other vertices.
The initial and terminal vertices $\iota$ and $\tau$ of the black replacement graph are two distinct vertices of in-degree $1$ and out-degree $0$.
The red replacement graph is a single red edge.

$\mathcal{D}_n^+$ is not an expanding edge replacement system (\cref{def.expanding}), since the red replacement graph only has two vertices linked by an edge.
Even if this does not pose an obstacle to the definition of a rearrangement group $G_n^+$ based on $\mathcal{D}_n^+$ (it would not be defined as a group of homeomorphisms, but only as a subgroup of the topological full group of an edge shift), the gluing relation (\cref{def.glue}) might not be (and in fact is not) an equivalence relation, so the limit space cannot be defined in the usual way (\cref{def.limit.space}).
However, since the gluing relation of $\mathcal{D}_n^+$ is reflexive and symmetric, it suffices to take the quotient of the symbol space of $\mathcal{D}_n^+$ by the transitive closure of the gluing relation in order to obtain an equivalence relation and a well defined limit space.

As noted in \cite[Remark 1.23]{BF19}, this limit space might not be as well-behaved as in the expanding case.
In our case, since the red replacement graph consists of a sole red edge, it is null-expanding (\cref{def.null.expanding.isolated}), and it is easy to see that each sequence $x \overline{r}$, where $xr$ is a null-expanding edge, ends up collapsing to a single point in the limit space, so the $n$ red cells $\llbracket y r_1 \rrbracket, \dots, \llbracket y r_n \rrbracket$ (for any finite word $y$ ending with a black edge) all degenerate to the same single point (these will correspond exactly to the branch points).
It is easy then to check that the limit spaces of $\mathcal{D}_n$ and $\mathcal{D}_n^+$ are the same, so the possible issues of the limit space not being well-behaved do not show up here.
More precisely, there is a canonical homeomorphism between the two limit spaces given by mapping each sequence $\alpha$ of $\mathcal{D}_n^+$ that does not feature null-expanding edges to the same sequence of $\mathcal{D}_n$.
Every other sequences of $\mathcal{D}_n^+$ must be of the form $x r_i \overline{r}$, so it represents the same branch point as $x i 1 \overline{n}$ and thus are mapped to this same sequence of $\mathcal{D}_n$.

\phantomsection\label{txt.oriented.dendrite.rearrangements}
One can then define the rearrangement groups $G_n^+$ for the edge replacement systems $\mathcal{D}_n^+$ by associating to each graph pair diagram a homeomorphism in the usual way.
The addition of the red edges will simply force every rearrangement to preserve the orientation of the dendrite, by which we mean that they preserve a circular ordering of the set of branches at each branch point.
Moreover, an expanding edge replacement system could arguably be produced by applying the strategy used in the proof of \cref{prop.null.expanding.isolated.rarrangement}, although the limit space would not be a dendrite.
We leave the study of these orientation-preserving dendrite groups for future research, but we list here a few remarks about them:
\begin{enumerate}
    \item The group $G_3^+$ is likely isomorphic to the finitely generated Thompson-like group studied in the dissertation \cite{SmithDendrite}.
    That is the group of those homeomorphisms of $S^1$ that preserve a lamination induced from the Julia set for the complex map $z \to z^2 + i$ (which is homeomorphic to $D_3$).
    We suspect that every $G_n^+$ can be built as a group of homeomorphisms of the circle by using laminations in this way.
    \item The edges of $\mathcal{D}_n^+$ are not undirected, unlike those of $\mathcal{D}_n$ (see \cref{rmk.dendrite.undirected}).
    \item The permutation subgroups from \cref{lem:perm} translate to subgroups of $G_n^+$ that are cyclic groups of order $n$ instead of copies of the symmetric group on $n$ elements.
    \item The Thompson subgroups from \cref{lem:thomp} do not seem to immediately translate to $G_n^+$, but there are probably copies of Thompson's group $F$ inside $G_n^+$ that act on certain arcs of $D_n$.
    \item The transitive property described in \cref{prop:trans} for $G_n$ may translate to an extension of isomorphisms between trees with a rotation system.
    If this happens, it might be used to prove that $G_n^+$ is dense in the group $\mathbb{H}_n^+$ of all orientation-preserving homeomorphisms of $D_n$ (i.e., the group of those homeomorphisms that at each branch point preserve a circular order of the branches).
    \item The parity map from \cref{sub:parity} would be trivial in $G_n^+$, so the abelianization may simply be $\mathbb{Z}$ (which is what happens in the airplane rearrangement group $T_A$, see the next subsection for the relation between $T_A$ and dendrites.)
    Note that in \cite{SmithDendrite} it is conjectured that the abelianization of this dendrite-Thompson group is isomorphic to $\mathbb{Z} \oplus \mathbb{Z}_2 \oplus \mathbb{Z}_3$, picturing a much more complex situation.
    \item If one chooses a different terminal vertex $\tau$ in the black replacement graph (\cref{fig:rep:sys:+}), then some of these remarks do not hold anymore (for example, if $n$ is even and $\tau$ is the terminal vertex of the $n/2$-th edge, then the black edges are actually undirected).
\end{enumerate}

Moreover, it looks like the Kaleidoscopic group $\mathcal{K}(A_3)$ studied in \cite{OrientationD3} (and their generalizations $\mathcal{K}(C_n)$) include the orientation-preserving dendrite rearrangement group $G_3^+$ (and $G_n^+$, respectively) as subgroups.
Indeed, the groups $\mathcal{K}(C_n)$ seem to consist of homeomorphisms of $D_n$ that preserve an orientation at each branch point.
It may then be possible that $G_n^+$ is dense in $\mathcal{K}(C_n)$, but we did not investigate this here.

\subsection{Infinite Order and the Airplane Rearrangement Group}
\label{sub:Air}

The author could not find a natural generalization of the dendrite edge replacement systems $\mathcal{D}_n$ that produces the infinite-order Wa\.zewski dendrite $D_\infty$ (where each branch point has countably infinite order; see \cite{Duc20} for an in-depth study of the topological properties of the group $\mathrm{Homeo}(D_\infty)$).
Instead, if one tries to generalize the ideas described in the previous subsection for the construction of $\mathcal{D}_n^+$ to the infinite-order case, the resulting edge replacement system is expanding and has a limit space that is not a dendrite, and its rearrangement group is the airplane edge replacement system, depicted in \cref{fig.airplane.replacement}, and the rearrangement group $T_A$ was studied in \cite{airplane}.

Differently from $\mathcal{D}_n^+$, the airplane edge replacement system is expanding, and in particular red cells (cells were defined at \cref{def.cell}) do not collapse to branch points, but instead they end up forming infinitely many circles in the limit space, meaning that the limit space is certainly not a dendrite.
However, if one takes the quotient of the airplane limit space by the equivalence relation that relates two points whenever they belong to the same circle of the airplane (i.e., if they are represented by sequences $x \alpha$ and $x \beta$, where $x$ is any finite word ending with $b_2$ or $b_3$ and $\alpha$ and $\beta$ are infinite sequences in the alphabet $\{ r_1, r_2 \}$), then the resulting space is the infinite-order Wa\.zewski dendrite $D_\infty$.
From this perspective, one can think of the airplane rearrangement group $T_A$ as if it were some sort of a orientation-preserving dendrite rearrangement group $G_\infty^+$, because the elements of $T_A$ act by self-permutations on the set of these circles, so a faithful action of $T_A$ on $D_\infty$ is naturally defined.
We will be considering this action in the remainder of this section, but it is important to keep in mind that this is not the canonical action of a rearrangement group on its limit space.

Under this identification, the branch points of $D_\infty$ correspond precisely to the circles of the airplane, and each branch at a fixed branch point corresponds to a dyadic point on the circle (see \cite{airplane} for more details).
The main distinction between a dendrite rearrangement group $G_n$ and $T_A$ thus consists in the fact that the group of permutations of branches around a branch point is Thompson's group $T$ instead of the finite permutation group $\mathrm{Sym}(n)$.

Thanks to the transitive properties of $T_A$ proved in \cite{airplane}, it is not hard to show that an analog of \cref{prop:trans} holds for $T_A$, as sketched below.
The construction of the trees $T(\mathcal{F})$ described at \cpageref{txt:trees} clearly works even if countably many branches are attached to each branch point.
We equip each tree $T(\mathcal{F})$ with the rotation system naturally prompted by the circular ordering of branches of $D_\infty$ induced by the circular ordering of the dyadic points on each circle of the airplane.

\medskip %layout
\begin{proposition}
\label{prop:airplane:trans}
Given two finite subsets $\mathcal{F}_1$ and $\mathcal{F}_2$ of $\mathrm{Br}$, any graph isomorphism between $T(\mathcal{F}_1)$ and $T(\mathcal{F}_2)$ that is compatible with their rotation systems can be extended to an element of $T_A$.
\end{proposition}

\begin{sketch}
By induction on the number $m$ of vertices of $T(\mathcal{F}_1)$ and $T(\mathcal{F}_2)$, the base case is single transitivity, which descends from \cite[Lemma 5.1]{airplane}.
The inductive step is proved as in \cref{prop:trans}:
for $m+1$, exclude a leaf of $T(\mathcal{F}_1)$ and use the induction hypothesis to map the other $m$ vertices where they need to be mapped, then use the fact that Thompson's group $T$ acts $2$-transitively on the set of dyadic points of $S^1$ to find an element that fixes the aforementioned $m$ vertices and moves the previously excluded leaf to wherever it needs to be moved.
\end{sketch}

Consider the group $\mathbb{H}_\infty^+$ of all orientation-preserving homeomorphisms of $D_\infty$, by which we mean the group of those homeomorphisms of $D_\infty$ that preserve a circular order of branches at each branch point.
(This group can be more formally described as a Kaleidoscopic group \cite{Kaleid}, but for this discussion there is no need to introduce such machinery.)
If we equip $\mathbb{H}_\infty^+$ with the subspace topology inherited from $\mathbb{H}_\infty$ then, as was done in \cref{sub:dns} for the density of $G_n$ in $\mathbb{H}_n$, we can prove that $T_A$ is dense in $\mathbb{H}_\infty^+$ by showing that every orientation-preserving homeomorphism can be ``approximated'' by elements of $T_A$, i.e.:

\medskip %layout
\begin{claim}
Let $\phi \in \mathbb{H}_\infty^+$.
For every $k \geq 1$ and for each $p_1, \dots, p_k \in \mathrm{Br}$ there exists a rearrangement $g_k \in T_A$ such that $g_k(p_i) = \phi(p_i)$ for all $i=1,\dots,k$.
\end{claim}

Using the previous \cref{prop:airplane:trans}, the proof of this Claim is precisely the same as that of \cref{clm:dns}, so we have:

\medskip %layout
\begin{theorem}
\label{thm:dns:A}
The airplane rearrangement group $T_A$ is dense in the group $\mathbb{H}_\infty^+$ of all orientation-preserving homeomorphisms of $D_\infty$.
\end{theorem}

\begin{remark}
Theorem 1 of \cite{BasilicaDense} shows that the basilica rearrangement group $T_B$ is dense in the group $\mathrm{Aut}_+(T(\mathcal{B}))$ of orientation-preserving automorphisms of the (countably) infinite-degree regular tree.
Building on this result, the final remarks of \cite{BasilicaDense} mention that the airplane rearrangement group $T_A$ is likely dense in the group of the orientation-preserving (in the same sense that we used elsewhere) homeomorphisms of the Tits-Bruhat $\mathbb{R}$-tree of the field of formal Laurent series $\mathbb{Q}\llbracket t, t^{1/2}, t^{1/4}, \dots \rrbracket$, which we will denote here by $\mathbb{T}$.
This real tree, which is not compact, is homeomorphic to the infinite-degree Wa\.zewski dendrite minus its endpoints.
Indeed, one can define a natural totally bounded metric on $\mathbb{T}$ in order to build its completion $\overline{\mathbb{T}}$, which is compact and is thus homeomorphic to the dendrite $D_\infty$.
Since the action of a homeomorphism of $\mathbb{T}$ entirely determines the action on the set $\overline{\mathbb{T}} \setminus \mathbb{T}$ of endpoints, the groups of (orientation-preserving) homeomorphisms of $D_\infty$ and $\overline{\mathbb{T}}$ are the same.
Hence, \cref{thm:dns:A} supports the claim that $T_A$ is dense in the group of the orientation-preserving homeomorphisms of the Tits-Bruhat $\mathbb{R}$-tree of the field of formal Laurent series $\mathbb{Q} \llbracket t, t^{1/2}, t^{1/4}, \dots \rrbracket$.
\end{remark}

\begin{remark}
As a final note about the infinite-order Wa\.zewski dendrite, we mention that it might be possible to build a (non orientation-preserving) rearrangement group for $D_\infty$ by either allowing the base and replacement graphs of an edge replacement system to be infinite, or maybe by building a direct limit from the inclusions $G_n \leq G_{n+1}$ (\cref{prop:inclusion}).
If one allows any graph isomorphism, the group built with the first method is not countable, as it contains a copy of $\mathrm{Sym}(\infty)$;
the second method instead seems to construct a ``finitary'' countable subgroup of the previous one that embeds into Thompson's group $V$.
We leave the exploration of these generalizations for future research.
\end{remark}

\subsection{Wa\.zewski Dendrite with Multiple Orders}
\label{sub:gen:den}

The last possible generalization that we mention is about the so-called \textbf{generalized Wa\.zewski dendrites}:
for each finite subset $S \subset \mathbb{N}_{\geq3}$, there exists a unique dendrite $D_S$ whose points have orders that belongs to $S \cup \{1, 2\}$ and such that every arc of $D_S$ contains points of every order in $S$.
Section 6 of \cite{DM18} (the entirety of which has been an important reference throughout this work) is about the full homeomorphism groups of generalized Wa\.zewski dendrites.

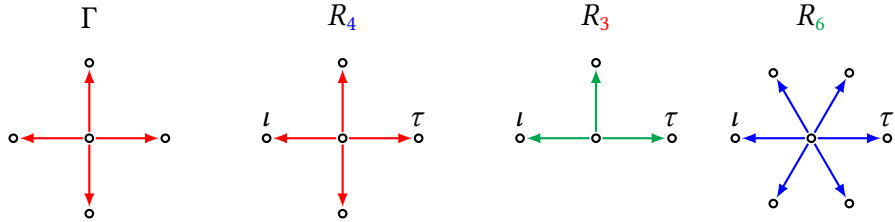
\begin{figure}
\centering
\begin{tikzpicture}
    \node at (0,1.6) {$\Gamma$};
    \node[vertex] (l) at (-1,0) {};
    \node[vertex] (r) at (1,0) {};
    \node[vertex] (c) at (0,0) {};
    \node[vertex] (ct) at (0,1) {};
    \node[vertex] (cb) at (0,-1) {};
    \draw[edge,red] (c) to (l);
    \draw[edge,red] (c) to (r);
    \draw[edge,red] (c) to (ct);
    \draw[edge,red] (c) to (cb);
    \begin{scope}[xshift=3.333cm]
    \node at (0,1.6) {$R_{\textcolor{blue}{4}}$};
    \node[vertex] (l) at (-1,0) {}; \draw (-1,0) node[above]{$\iota$};
    \node[vertex] (r) at (1,0) {}; \draw (1,0) node[above]{$\tau$};
    \node[vertex] (c) at (0,0) {};
    \node[vertex] (ct) at (0,1) {};
    \node[vertex] (cb) at (0,-1) {};
    \draw[edge,red] (c) to (l);
    \draw[edge,red] (c) to (r);
    \draw[edge,red] (c) to (ct);
    \draw[edge,red] (c) to (cb);
    \end{scope}
    \begin{scope}[xshift=6.667cm]
    \node at (0,1.6) {$R_{\textcolor{red}{3}}$};
    \node[vertex] (l) at (-1,0) {}; \draw (-1,0) node[above]{$\iota$};
    \node[vertex] (r) at (1,0) {}; \draw (1,0) node[above]{$\tau$};
    \node[vertex] (c) at (0,0) {};
    \node[vertex] (ct) at (0,1) {};
    \draw[edge,Green] (c) to (l);
    \draw[edge,Green] (c) to (r);
    \draw[edge,Green] (c) to (ct);
    \end{scope}
    \begin{scope}[xshift=9.5cm]
    \node at (0,1.6) {$R_{\textcolor{Green}{6}}$};
    \node[vertex] (l) at (-1,0) {}; \draw (-1,0) node[above]{$\iota$};
    \node[vertex] (r) at (1,0) {}; \draw (1,0) node[above]{$\tau$};
    \node[vertex] (c) at (0,0) {};
    \node[vertex] (v1) at (60:1) {};
    \node[vertex] (v2) at (120:1) {};
    \node[vertex] (v3) at (240:1) {};
    \node[vertex] (v4) at (300:1) {};
    \draw[edge,blue] (c) to (l);
    \draw[edge,blue] (c) to (r);
    \draw[edge,blue] (c) to (v1);
    \draw[edge,blue] (c) to (v2);
    \draw[edge,blue] (c) to (v3);
    \draw[edge,blue] (c) to (v4);
    \end{scope}
\end{tikzpicture}
\caption{An edge replacement system for the generalized Wa\.zewski dendrite $D_{\{3,4,6\}}$ with the ordering $4 \to 3 \to 6$ on $S$.}
\label{fig:rep:sys:col}
\end{figure}

Given a finite $S \subset \mathbb{N}_{\geq3}$, we can build a (polychormatic) edge replacement system in the following way.
First, fix an ordering $s_1, \dots, s_k$ of the elements of $S$.
\begin{itemize}
    \item The set of colors is $\{1, \dots, k\}$.
    \item The $i$-th replacement graph is a tree consisting of a vertex of degree $s_i$ that is the origin of $s_i$ edges colored by $s_{i+1}$ terminating at $s_i$ distinct leaves;
    the initial and terminal vertices $\iota$ and $\tau$ are two distinct leaves.
    \item The base graph is the same as the first replacement graph (in truth, any replacement graph would work).
\end{itemize}
\cref{fig:rep:sys:col} depicts an example of such an edge replacement system, where $S = \{3, 4, 6\}$ and the chosen ordering is $4 \to 3 \to 6$.

The fact that the limit space is $D_S$ can be proved immediately as done at \cpageref{txt:rep:sys:den} for $\mathcal{D}_n$, so in particular the chosen ordering on $S$ does not change the limit space by \cite[Theorem 6.2]{selfhomeomorphisms}.
It is not clear, however, whether the rearrangement group is affected by this choice.
Consequently, in order to keep track of this order, we denote by $G_{s_1, \dots, s_k}$ the corresponding rearrangement group.

It is plausible that $G_{s_1, \dots, s_k}$ is generated by a copy for each of the finite groups $\mathrm{Sym}(s_i)$ along with the subgroup of Thompson's group $F$ that preserves the coloration of the edges.
If this subgroup is transitive enough, then $G_{s_1, \dots, s_k}$ is likely to be dense in $\mathrm{Homeo}(D_S)$.
The study of these groups is left for future research.

Finally, one could ``mix'' the various generalizations proposed throughout this section.
For example, there may be orientation-preserving rearrangement groups of generalized Wa\.zewski dendrites where the Thompson's group acting on each EE-arc is a color-preserving subgroup of some Higman-Thompson group $F_m$.
This too is left for future research.

%%%%%%%%%%%%%%%%%%%%%%%%%

\chapter{The Conjugacy Problem}
\label{cha.conjugacy}

The \textbf{conjugacy problem} is the decision problem of determining whether two given elements of a group $G$ are conjugate in $G$.
The \textbf{conjugacy search problem} is the problem of producing a conjugating element between two given elements of a group $G$ that are conjugate in $G$.
Both problems have been shown to be solvable in the three original Thompson groups:
Guba and Sapir solved it for the whole class of diagram groups, which includes $F$ \cite{guba1997diagram}; it was solved for $V$ by Higman \cite{higman1974finitely} and more generally for Higman-Thompson groups in \cite{HigmanConj}.
In \cite{BM14} Belk and Matucci produced a unified solution for the conjugacy problem in Thompson groups $F$, $T$ and $V$ with a technique that involves the use of special graphs called \textit{strand diagrams} as a way to represent the elements of the groups and, after being ``closed'', their conjugacy classes.
Generalizations of strand diagrams were also used to study conjugacy in almost automorphism groups of trees by Goffer and Lederle \cite{Goffer2019ConjugacyAD} and to solve the conjugacy problem in symmetric Thompson groups (a generalization of $F$, $T$ and $V$) by Aroca \cite{Aroca2018TheCP}.
Diagrams of this kind have also been used by Bux and Santos Rego to solve the conjugacy problem in the braided Thompson's group $V_{br}$ in a forthcoming work \cite{brVconj}, which has been announced in \cite{braidedreport} and is currently being written up.

In this chapter we present the content of the author's work \cite{conjugacy}.
Mainly inspired by \cite{BM14}, the aim of \cite{conjugacy} is to introduce a version of strand diagrams that represent the elements of a rearrangement group, and to use them to solve the conjugacy problem under a certain condition on the edge replacement system that generates the fractal on which the group acts.
More precisely, we prove the following, which is a collection of \cref{thm.conjugacy.solvable,thm.conjugacy.known.groups}:

\medskip %layout
\begin{maintheorem*}
Given an expanding edge replacement system whose edge replacement rules are reduction-confluent, the conjugacy problem and the conjugacy search problem are solvable in the associated rearrangement group. In particular, this method can be used to solve the conjugacy and the conjugacy search problems in the following groups: 
\begin{itemize}
    \item Thompson groups $F$, $T$ and $V$ (\cref{sub.Thompson}) and the Higman-Thompson groups $F_{n,r}$, $T_{n,r}$ and $V_{n,r}$ (\cref{sub.higman.thompson.groups});
    \item the basilica and rabbit rearrangement groups (\cref{sub.basilica});
    \item the airplane rearrangement group $T_A$ (\cref{sub.airplane}) through an adapted method, see below;
    \item the dendrite rearrangement groups $G_n$ (\cref{sub.other.fractals,cha.dendrites});
    \item the Vicsek rearrangement group and its generalizations and the bubble bath rearrangement group (\cite{BF19});
    \item the Houghton groups $H_n$ (\cref{sub.Houghton});
    \item the Thompson-like groups $QV$, $QT$ and $QF$ (\cref{sub.thompson.like});
    \item certain topological full groups of edge shifts (\cref{sub.topological.full.groups}).
\end{itemize}
This solves the conjugacy problem for all rearrangement groups that have been studied so far.
\end{maintheorem*}

We recall that, for the Houghton groups, the conjugacy problem was already solved with different methods in \cite{HoughtonConjugacy}, and even the twisted conjugacy problem was solved in \cite{COX2017390}.
However, the observation that such groups can be realized as rearrangement groups is new and leads to an independent solution to the conjugacy problem.
Moreover, the solution to the conjugacy problem for the groups $QV$, $QT$ and $QF$ and for the rearrangement groups of the airplane, the basilica, the rabbits, the Vicsek and the bubble bath are new, as far as the author knows.
We also remark here that, despite $QV$ being an extension of $V$, \cite[Theorem 3.1]{10.2307/40590898} does not apply to solve the conjugacy problem in $QV$ because the third condition is not met.
Finally, it is interesting to note that $T_B$ is a finitely generated group with solvable conjugacy problem that is not finitely presented (which was proved in \cite{WZ19}).

We employ strand diagrams, but unlike the solution of the conjugacy problem for $F,T,V$ in \cite{BM14} we need to simultaneously make use of two different yet related rewriting systems, one of which is an instance of the so-called \textit{graph rewriting systems}.
We then need to find conditions to find uniquely reduced diagrams, for which we introduce the \textit{reduction-confluence} condition to achieve confluence in our instance of graph rewriting system.
The study of this graph rewriting system did not appear in the unified solution to the conjugacy problem for $F,T,V$ in \cite{BM14} as in those cases it is trivially confluent.
We remark that establishing whether a graph rewriting system is confluent is in general an unsolvable problem and the search for sufficient conditions is actively studied in computer science (see \cite{Plump2005}).

All the groups in the Main Theorem above satisfy the reduction-confluence hypothesis with the exception of $T_A$, but this issue can be circumvented by adding one ``virtual reduction'' to the airplane edge replacement system, as described in \cref{sub.non.confluent}.
More generally, as discussed at the end of this chapter, it seems that strand diagrams can be used to solve the conjugacy problem with this method whenever the reduction system can be made confluent by adding finitely many reduction rules that preserve termination of the rewriting system and the equivalence generated by the rewriting.
Determining when this can be achieved is beyond the scope of this dissertation, since it is related to problems in the more general setting of graph rewriting systems, as will be briefly discussed in \cref{sub.DPO}.

Before delving into strand diagrams, we remark that many generalizations of the conjugacy problem have been studied for Thompson groups:
for example, the simultaneous conjugacy problem for $F$ is solved in \cite{simultaneous} and the twisted conjugacy problem and property $R_\infty$ have been studied in $F$ and $T$ in \cite{Bleak2007TwistedCC, Burillo2013TheCP, Gonalves2018TwistedCI}.
Also, \cite{Belk2021ConjugatorLI} makes use of the technology of strand diagrams to study the conjugator lengths in Thompson groups.
All of these questions can be asked about rearrangement groups as well, and maybe they can be tackled with the aid of strand diagrams.
Additionally, rearrangement groups can be seen as a special case of the recently defined \textit{graph diagram groups} \cite{graphdiagramgroups} and it would be interesting to develop a concept of strand diagram in that setting and investigate conjugacy.

Finally, we observe that when strand diagrams were introduced in \cite{BM14}, they were also used to study dynamics in Thompson groups.
The first version of the forest pair diagrams for rearrangements that is fully developed in \cref{sec.forest.pair.diagrams} was initially used to understand some of the dynamics of rearrangements in \cite{IG} (adapted here into \cref{cha.IG}) and used to show dynamical properties of rearrangements (the existence of wandering cells, see \cref{prop.wandering.cells}).
The strand diagrams for rearrangement groups introduced in this chapter are likely usable to shed further light on the dynamics of elements as done in \cite{BM14} for $F$, $T$ and $V$, but we leave such questions for future study.

\section{Strand Diagrams}
\label{sec.SDs}

Following the ideas of \cite{BM14}, in this section we introduce strand diagrams that represent elements of a rearrangement group (\cref{sub.SDs}).
While doing so, it will be natural and useful to introduce a groupoid consisting of generalized rearrangements obtained by allowing different base graph for the domain and the range graphs, while still keeping the same edge replacement rules (\cref{SUB replacement groupoid}).
Although this is defined here in terms of strand diagrams, this is essentially the same groupoid that was introduced in \cite[Subsection 3.1]{BF19}.

We will then introduce reduction rules to find a unique minimal reduced diagram for each element (\cref{sub.SDs.reduction}), we will see how strand diagrams relate to rearrangements (\cref{sub.SDs.are.rearrangements}) and we will describe how to compose two strand diagrams (\cref{sub.SDs.composition}).

\begin{figure}
    \centering
    \begin{tikzpicture}
        \node[node] (r1) at (0,1) {};
        \node[node] (r2) at (1,1) {};
        \node[node] (r3) at (2,1) {};
        \node[node] (r4) at (3,1) {};
        \node[node] (s1) at (0,0) {};
        \node[node] (s2) at (1,0) {};
        \node[node] (s3) at (2,0) {};
        \node[node] (s4) at (3,0) {};
        \draw (r1) to node[red,left,align=center]{c\\b} (s1);
        \draw (r2) to node[blue,left,align=center]{b\\a} (s2);
        \draw (r3) to node[red,left,align=center]{b\\c} (s3);
        \draw (r4) to node[blue,left,align=center]{c\\d} (s4);
    \end{tikzpicture}
    \caption{The base forest for the airplane edge replacement system modified for \cref{cha.conjugacy} as described in \cref{ass.airplane.diagrams}.}
    \label{fig.airplane.base.forest.conjugacy.chapter}
\end{figure}
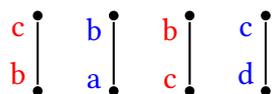

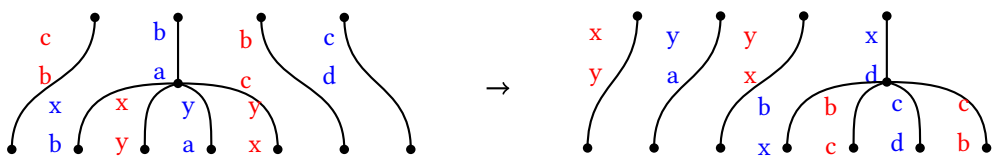
\begin{figure}
    \centering
    \begin{minipage}{.485\textwidth}\centering
    \begin{tikzpicture}[font=\footnotesize,scale=.875]
        \draw (0,0) node[black,circle,fill,inner sep=1.25]{} node[red,xshift=-.65cm,yshift=-.55cm,align=center]{c\\b} to[out=270,in=90,looseness=1.2] (-1.25,-2) node[black,circle,fill,inner sep=1.25]{};
        \draw (1.25,0) node[black,circle,fill,inner sep=1.25]{} -- node[blue,left,align=center]{b\\a} (1.25,-1) node[black,circle,fill,inner sep=1.25]{};
        \draw (2.5,0) node[black,circle,fill,inner sep=1.25]{} node[red,xshift=-.2cm,yshift=-.55cm,align=center]{b\\c} to[out=270,in=90,looseness=1.2] (3.75,-2) node[black,circle,fill,inner sep=1.25]{};
        \draw (3.75,0) node[black,circle,fill,inner sep=1.25]{} node[blue,xshift=-.2cm,yshift=-.55cm,align=center]{c\\d} to[out=270,in=90,looseness=1.2] (4.75,-2) node[black,circle,fill,inner sep=1.25]{};
        \begin{scope}[xshift=1.25cm]
        \draw (0,-1) to[out=180,in=90,looseness=1.2] (-1.5,-2) node[black,circle,fill,inner sep=1.25]{} node[blue,xshift=-.3cm,yshift=.3cm,align=center]{x\\b};
        \draw (0,-1) to[out=195,in=90,looseness=1.2] (-.5,-2) node[black,circle,fill,inner sep=1.25]{} node[red,xshift=-.3cm,yshift=.3cm,align=center]{x\\y};
        \draw (0,-1) to[out=345,in=90,looseness=1.2] (.5,-2) node[black,circle,fill,inner sep=1.25]{} node[blue,xshift=-.3cm,yshift=.3cm,align=center]{y\\a};
        \draw (0,-1) to[out=0,in=90,looseness=1.2] (1.5,-2) node[black,circle,fill,inner sep=1.25]{} node[red,xshift=-.3cm,yshift=.3cm,align=center]{y\\x};
        \end{scope}
    \end{tikzpicture}
    \end{minipage}%
    \begin{minipage}{.03\textwidth}\centering
    \begin{tikzpicture}[scale=1]
        \draw[-to] (0,0) -- (.3,0);
    \end{tikzpicture}
    \end{minipage}%
    \begin{minipage}{.485\textwidth}\centering
    \begin{tikzpicture}[font=\footnotesize,scale=.875]
        \draw (0,0) node[black,circle,fill,inner sep=1.25]{} node[red,xshift=-.55cm,yshift=-.55cm,align=center]{x\\y} to[out=270,in=90,looseness=1.2] (-.75,-2) node[black,circle,fill,inner sep=1.25]{};
        \draw (1.25,0) node[black,circle,fill,inner sep=1.25]{} node[blue,xshift=-.625cm,yshift=-.55cm,align=center]{y\\a} to[out=270,in=90,looseness=1.2] (.25,-2) node[black,circle,fill,inner sep=1.25]{};
        \draw (2.5,0) node[black,circle,fill,inner sep=1.25]{} node[red,xshift=-.7cm,yshift=-.55cm,align=center]{y\\x} to[out=270,in=90,looseness=1.2] (1.25,-2) node[black,circle,fill,inner sep=1.25]{};
        \draw (3.75,0) node[black,circle,fill,inner sep=1.25]{} node[blue,xshift=-.2cm,yshift=-.55cm,align=center]{x\\d} -- (3.75,-1) node[black,circle,fill,inner sep=1.25]{};
        \begin{scope}[xshift=3.75cm]
        \draw (0,-1) to[out=180,in=90,looseness=1.2] (-1.5,-2) node[black,circle,fill,inner sep=1.25]{} node[blue,xshift=-.3cm,yshift=.3cm,align=center]{b\\x};
        \draw (0,-1) to[out=195,in=90,looseness=1.2] (-.5,-2) node[black,circle,fill,inner sep=1.25]{} node[red,xshift=-.3cm,yshift=.3cm,align=center]{b\\c};
        \draw (0,-1) to[out=345,in=90,looseness=1.2] (.5,-2) node[black,circle,fill,inner sep=1.25]{} node[blue,xshift=-.3cm,yshift=.3cm,align=center]{c\\d};
        \draw (0,-1) to[out=0,in=90,looseness=1.2] (1.5,-2) node[black,circle,fill,inner sep=1.25]{} node[red,xshift=-.3cm,yshift=.3cm,align=center]{c\\b};
        \end{scope}
    \end{tikzpicture}
    \end{minipage}
    \caption{A forest pair diagram for the same element represented in \cref{fig.airplane.forest.pair.diagram}, using the convention described in \cref{ass.airplane.diagrams}.}
    \label{fig.airplane.forest.pair.diagram.conjugacy.chapter}
\end{figure}

\begin{assumption}
\label{ass.airplane.diagrams}
In this chapter the airplane edge replacement system (\cref{fig.airplane.replacement}) will be our guiding example.
In order to be able to draw more concise strand diagrams later on and to be consistent with the notation and figures of \cite{conjugacy}, we change the base graph to its first graph expansion.
The base graph will then be the same as the blue replacement graph and thus the base forest will consist of 4 roots, as depicted in \cref{fig.airplane.base.forest.conjugacy.chapter}.
With this change, the rearrangement represented by the forest pair diagram depicted in \cref{fig.airplane.forest.pair.diagram} is instead represented by \cref{fig.airplane.forest.pair.diagram.conjugacy.chapter}.
\end{assumption}

\subsection{Generic Strand Diagrams}
\label{sub.SDs}

\begin{figure}\centering
    \begin{subfigure}[b]{.475\textwidth}\centering
    \begin{tikzpicture}[font=\footnotesize,scale=.875]
        \draw (0,0) node[black,circle,fill,inner sep=1.25]{} node[red,xshift=-.65cm,yshift=-.55cm,align=center]{c\\b} to[out=270,in=90,looseness=1.2] (-1.25,-2);
        \draw (1.25,0) node[black,circle,fill,inner sep=1.25]{} -- node[blue,left,align=center]{b\\a} (1.25,-1) node[black,circle,fill,inner sep=1.25]{};
        \draw (2.5,0) node[black,circle,fill,inner sep=1.25]{} node[red,xshift=-.2cm,yshift=-.55cm,align=center]{b\\c} to[out=270,in=90,looseness=1.2] (3.75,-2);
        \draw (3.75,0) node[black,circle,fill,inner sep=1.25]{} node[blue,xshift=-.2cm,yshift=-.55cm,align=center]{c\\d} to[out=270,in=90,looseness=1.2] (4.75,-2);
        \begin{scope}[xshift=1.25cm]
        \draw (0,-1) to[out=180,in=90,looseness=1.2] (-1.5,-2);
        \draw (0,-1) to[out=195,in=90,looseness=1.2] (-.5,-2);
        \draw (0,-1) to[out=345,in=90,looseness=1.2] (.5,-2);
        \draw (0,-1) to[out=0,in=90,looseness=1.2] (1.5,-2);
        \end{scope}
        \draw (-1.25,-2) to[out=270,in=90,looseness=.8] (4.75,-4);
        \draw (-.25,-2) to[out=270,in=90,looseness=1.2] (1.75,-4);
        \draw (.75,-2) to[out=270,in=90,looseness=1] (-1.25,-4);
        \draw (1.75,-2) to[out=270,in=90,looseness=1] (-.25,-4);
        \draw (2.75,-2) to[out=270,in=90,looseness=1] (.75,-4);
        \draw (3.75,-2) to[out=270,in=90,looseness=1.2] (2.75,-4);
        \draw (4.75,-2) to[out=270,in=90,looseness=1.2] (3.75,-4);
        \begin{scope}[yshift=-4cm]
        \draw (-1.25,0) to[out=270,in=90,looseness=1.2] (0,-2) node[red,xshift=-.45cm,yshift=.3cm,align=center]{x\\y} node[black,circle,fill,inner sep=1.25]{};
        \draw (-.25,0) to[out=270,in=90,looseness=1.2] (1.25,-2) node[blue,xshift=-.45cm,yshift=.3cm,align=center]{y\\a} node[black,circle,fill,inner sep=1.25]{};
        \draw (.75,0) to[out=270,in=90,looseness=1.2] (2.5,-2) node[red,xshift=-.45cm,yshift=.3cm,align=center]{y\\x} node[black,circle,fill,inner sep=1.25]{};
        \draw (3.75,-1) node[black,circle,fill,inner sep=1.25]{} node[blue,xshift=-.2cm,yshift=-.45cm,align=center]{x\\d} -- (3.75,-2) node[black,circle,fill,inner sep=1.25]{};
        \begin{scope}[xshift=3.75cm]
        \draw (0,-1) to[out=180,in=270,looseness=1.2] (-2,0) node[blue,xshift=-.25cm,yshift=0cm,align=center]{x\\b};
        \draw (0,-1) to[out=180,in=270,looseness=1.2] (-1,0);
        \draw (0,-1) to[out=90,in=270,looseness=1.2] (0,0);
        \draw (0,-1) to[out=0,in=270,looseness=1.2] (1,0);
        \end{scope}
        \end{scope}
    \end{tikzpicture}
    \caption{}
    \label{fig_strand_alpha}
    \end{subfigure}
    \begin{subfigure}[b]{.475\textwidth}\centering
    \begin{tikzpicture}[font=\footnotesize,scale=.875]
        \draw (0,0) node[black,circle,fill,inner sep=1.25]{} node[red,xshift=-.65cm,yshift=-.55cm,align=center]{p\\q} to[out=270,in=90,looseness=1.2] (-1.25,-2);
        \draw (1.25,0) node[black,circle,fill,inner sep=1.25]{} -- node[blue,left,align=center]{q\\e} (1.25,-1) node[black,circle,fill,inner sep=1.25]{};
        \draw (2.5,0) node[black,circle,fill,inner sep=1.25]{} node[red,xshift=-.2cm,yshift=-.55cm,align=center]{q\\p} to[out=270,in=90,looseness=1.2] (3.75,-2);
        \draw (3.75,0) node[black,circle,fill,inner sep=1.25]{} node[blue,xshift=-.2cm,yshift=-.55cm,align=center]{p\\f} to[out=270,in=90,looseness=1.2] (4.75,-2);
        \begin{scope}[xshift=1.25cm]
        \draw (0,-1) to[out=180,in=90,looseness=1.2] (-1.5,-2);
        \draw (0,-1) to[out=195,in=90,looseness=1.2] (-.5,-2);
        \draw (0,-1) to[out=345,in=90,looseness=1.2] (.5,-2);
        \draw (0,-1) to[out=0,in=90,looseness=1.2] (1.5,-2);
        \end{scope}
        \draw (-1.25,-2) to[out=270,in=90,looseness=.8] (4.75,-4);
        \draw (-.25,-2) to[out=270,in=90,looseness=1.2] (1.75,-4);
        \draw (.75,-2) to[out=270,in=90,looseness=1] (-1.25,-4);
        \draw (1.75,-2) to[out=270,in=90,looseness=1] (-.25,-4);
        \draw (2.75,-2) to[out=270,in=90,looseness=1] (.75,-4);
        \draw (3.75,-2) to[out=270,in=90,looseness=1.2] (2.75,-4);
        \draw (4.75,-2) to[out=270,in=90,looseness=1.2] (3.75,-4);
        \begin{scope}[yshift=-4cm]
        \draw (-1.25,0) to[out=270,in=90,looseness=1.2] (0,-2) node[red,xshift=-.45cm,yshift=.3cm,align=center]{v\\w} node[black,circle,fill,inner sep=1.25]{};
        \draw (-.25,0) to[out=270,in=90,looseness=1.2] (1.25,-2) node[blue,xshift=-.45cm,yshift=.3cm,align=center]{w\\e} node[black,circle,fill,inner sep=1.25]{};
        \draw (.75,0) to[out=270,in=90,looseness=1.2] (2.5,-2) node[red,xshift=-.45cm,yshift=.3cm,align=center]{w\\v} node[black,circle,fill,inner sep=1.25]{};
        \draw (3.75,-1) node[black,circle,fill,inner sep=1.25]{} node[blue,xshift=-.2cm,yshift=-.45cm,align=center]{v\\f} -- (3.75,-2) node[black,circle,fill,inner sep=1.25]{};
        \begin{scope}[xshift=3.75cm]
        \draw (0,-1) to[out=180,in=270,looseness=1.2] (-2,0) node[blue,xshift=-.25cm,yshift=0cm,align=center]{v\\b};
        \draw (0,-1) to[out=180,in=270,looseness=1.2] (-1,0);
        \draw (0,-1) to[out=90,in=270,looseness=1.2] (0,0);
        \draw (0,-1) to[out=0,in=270,looseness=1.2] (1,0);
        \end{scope}
        \end{scope}
    \end{tikzpicture}
    \caption{}
    \label{fig_strand_renaming}
    \end{subfigure}
    \caption{The strand diagram for the forest pair diagram of \cref{fig.airplane.forest.pair.diagram.conjugacy.chapter} (a) and a renaming of its symbols (b).}
    \label{fig_strand}
\end{figure}
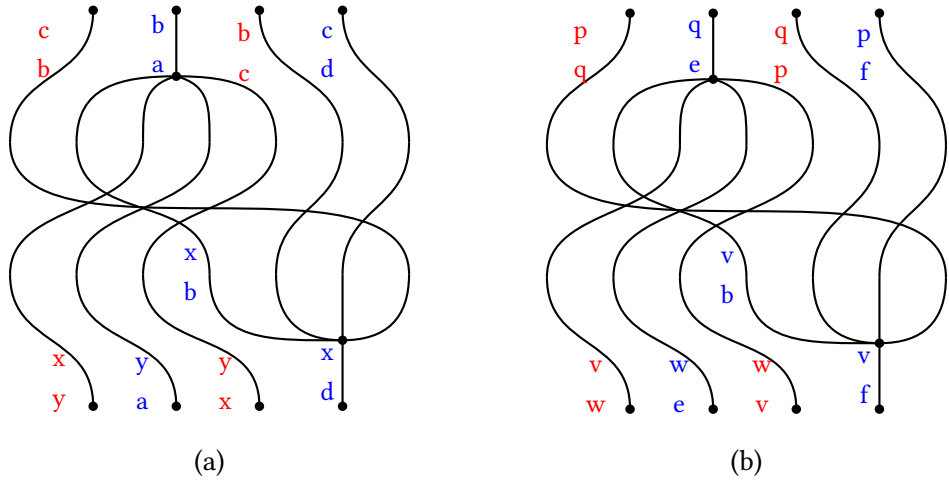

Let $(F_D, F_R)$ be a forest pair diagram (\cref{def.forest.pair.diagram}).
As done in \cref{fig_strand_alpha}, draw $F_D$ and $F_R$ with the range forest turned upside down below the domain forest and join each leaf of $F_D$ to its image in $F_R$, which is the unique leaf of $F_R$ with the same label.
The result is the \textbf{strand diagram} corresponding to $(F_D, F_R)$, consisting of ``strands'' starting at the top, merging and splitting in the middle and hanging at the bottom, along with labels such as those of forest expansions.
Observe that we can recover $F_D$ and $F_R$ (along with the permutation between their leaves implicitly specified by the labels, see \cref{rmk.forest.pair.diagrams.convention}) by ``cutting'' the strand diagram in the unique way that separates the merges from the splits.

Now, these only cover those strand diagrams that are obtained by gluing forest pair diagrams.
In general, the definition of a strand diagram is the one given below.
However, we will see in \cref{sub.SDs.are.rearrangements} that, up to reductions (described in \cref{sub.SDs.reduction}), these are really all of the strand diagrams that are needed to describe rearrangements.

\begin{definition}\label{def.SDs}
A \textbf{strand diagram} is a finite acyclic directed graph whose edges are colored and labeled by ordered triples of symbols, together with a rotation system, such that every vertex is either a source of out-degree $1$, a sink of in-degree $1$, a split or a merge, and with a given ordering of both the sources and the sinks.

By \textbf{split} we mean a vertex that has in-degree $1$ and out-degree at least $2$.
Conversely, by \textbf{merge} we mean a vertex that has out-degree $1$ and in-degree at least $2$.
\end{definition}

In order to avoid confusion between the edges of graph expansions, those of forest expansions (which we called branches) and those of strand diagrams, we will refer to the latter by the term \textbf{strands}.

We consider two strand diagrams equal if there exists a graph isomorphism between them that is compatible with the corresponding rotation systems, and if the two strand diagrams differ by a renaming of the labeling symbols.
An example of such a renaming is depicted in \cref{fig_strand_renaming}.

As done in \cref{fig_strand}, it is convenient to depict sources at the top and sinks at the bottom, both aligned and ordered from left to right in their given order;
this allows us to hide the orientation of the strands, which is always implied to descend from the sources to the sinks.
Moreover, for the sake of clarity we color the label associated to the strand instead of the strand itself, as was done with forest pair diagrams.

The set of \textit{all} strand diagrams is far too large and varied to wield meaningful information about rearrangements (in fact, this general definition does not even take into account any information of the edge replacement system).
This is why we turn our attention to replacement groupoids: classes of strand diagrams determined by the edge replacement rules, as described in the following subsection.

\subsection{Replacement Groupoids}
\label{SUB replacement groupoid}

If $X$ is a replacement graph, we say that a split (merge) is an $X$\textbf{-split} ($X$\textbf{-merge}) if, up to renaming symbols, it is a copy of the replacement tree associated to $X$ (\cref{sub.forests.and.graphs}).
More precisely, an $X$-split consists of a top strand that splits into as many bottom strands as there are edges in $X$, in their given ordering;
if the top strand is labeled by $(\iota,\tau,\epsilon)$, then each bottom strand is labeled by $(v,w,z)$ where $v$ and $w$, respectively, are the initial and terminal vertices of the corresponding edge of $X$, and $z$ is an index that distinguishes parallel edges, up to renaming vertices and indices, with $\iota$ and $\tau$ the initial and terminal vertices of $X$, respectively.
An $X$-merge is the same, with inverted direction of strands.

\phantomsection\label{txt.branching.strands}
By \textbf{branching strand} we refer to the unique top strand of a split or the unique bottom strand of a merge.
Additionally, we say that a symbol is \textbf{generated} by a split (merge) if it appears among the symbols included in the split (merge) excluding the branching strand, i.e., if it represents a new vertex in the graph expansion.

\begin{definition}
\label{def.r.branching}
Let $(R, \mathrm{C})$ be a set of edge replacement rules (\cref{def.replacement.rules}), with $R = \{ X_i \mid i \in \mathrm{C} \}$.
We say that a strand diagram is \textbf{$R$-branching} if:
\begin{enumerate}
    \item every split and every merge is an $X_i$-split or an $X_i$-merge for some $i \in \mathrm{C}$;
    \item whenever a sequence of $k$ (possibly $k=1$) merges is immediately followed by a sequence of $k$ splits that mirror each other, the strands of both sequences are labeled exactly in the same way, in the same order (more precisely, by sequences that ``mirror'' each other we mean that the sequence of merges is followed by a sequence of splits that is precisely its inverse if read without labels);
    \item the same symbol cannot be generated by different splits whose branching strands have different labels;
    likewise, the same symbol cannot be generated by different merges whose branching strands have different labels.
\end{enumerate}
\end{definition}

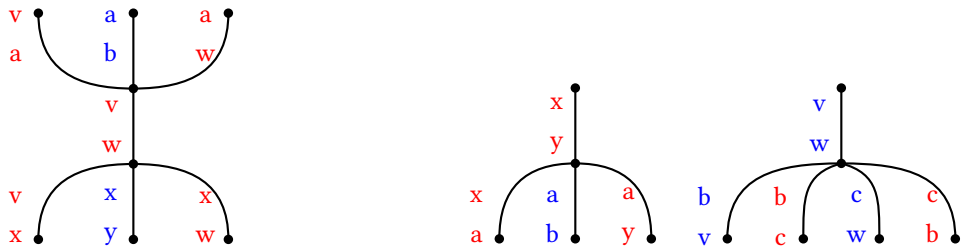
\begin{figure}\centering
    \begin{subfigure}[t]{.333\textwidth}\centering
    \begin{tikzpicture}[font=\footnotesize,scale=1]
        \draw (0,0) to[out=180,in=270,looseness=1.2] (-1.25,1) node[black,circle,fill,inner sep=1.25]{} node[red,xshift=-.3cm,yshift=-.3cm,align=center]{v\\a};
        \draw (0,0) -- (0,1) node[black,circle,fill,inner sep=1.25]{} node[blue,xshift=-.3cm,yshift=-.3cm,align=center]{a\\b};
        \draw (0,0) to[out=0,in=270,looseness=1.2] (1.25,1) node[black,circle,fill,inner sep=1.25]{} node[red,xshift=-.3cm,yshift=-.3cm,align=center]{a\\w};
        \draw (0,0) node[black,circle,fill,inner sep=1.25]{} -- node[red,left,align=center]{v\\w} (0,-1) node[black,circle,fill,inner sep=1.25]{};
        \draw (0,-1) to[out=180,in=90,looseness=1.2] (-1.25,-2) node[black,circle,fill,inner sep=1.25]{} node[red,xshift=-.3cm,yshift=.3cm,align=center]{v\\x};
        \draw (0,-1) -- (0,-2) node[black,circle,fill,inner sep=1.25]{} node[blue,xshift=-.3cm,yshift=.3cm,align=center]{x\\y};
        \draw (0,-1) to[out=0,in=90,looseness=1.2] (1.25,-2) node[black,circle,fill,inner sep=1.25]{} node[red,xshift=-.3cm,yshift=.3cm,align=center]{x\\w};
    \end{tikzpicture}
    \caption{Condition (2) does not hold.}
    \label{fig_notRbranching_2}
    \end{subfigure}
    \hfill
    \begin{subfigure}[t]{.6\textwidth}\centering
    \begin{tikzpicture}[font=\footnotesize,scale=1]
        \draw (0,0) node[black,circle,fill,inner sep=1.25]{} -- node[red,left,align=center]{x\\y} (0,-1) node[black,circle,fill,inner sep=1.25]{};
        \draw (0,-1) to[out=180,in=90,looseness=1.2] (-1,-2) node[black,circle,fill,inner sep=1.25]{} node[red,xshift=-.3cm,yshift=.3cm,align=center]{x\\a};
        \draw (0,-1) -- (0,-2) node[black,circle,fill,inner sep=1.25]{} node[blue,xshift=-.3cm,yshift=.3cm,align=center]{a\\b};
        \draw (0,-1) to[out=0,in=90,looseness=1.2] (1,-2) node[black,circle,fill,inner sep=1.25]{} node[red,xshift=-.3cm,yshift=.3cm,align=center]{a\\y};
        \begin{scope}[xshift=3.5cm]
            \draw (0,0) node[black,circle,fill,inner sep=1.25]{} -- node[blue,left,align=center]{v\\w} (0,-1) node[black,circle,fill,inner sep=1.25]{};
            \draw (0,-1) to[out=180,in=90,looseness=1.2] (-1.5,-2) node[black,circle,fill,inner sep=1.25]{} node[blue,xshift=-.3cm,yshift=.3cm,align=center]{b\\v};
            \draw (0,-1) to[out=195,in=90,looseness=1.2] (-.5,-2) node[black,circle,fill,inner sep=1.25]{} node[red,xshift=-.3cm,yshift=.3cm,align=center]{b\\c};
            \draw (0,-1) to[out=345,in=90,looseness=1.2] (.5,-2) node[black,circle,fill,inner sep=1.25]{} node[blue,xshift=-.3cm,yshift=.3cm,align=center]{c\\w};
            \draw (0,-1) to[out=0,in=90,looseness=1.2] (1.5,-2) node[black,circle,fill,inner sep=1.25]{} node[red,xshift=-.3cm,yshift=.3cm,align=center]{c\\b};
        \end{scope}
    \end{tikzpicture}
    \caption{Condition (3) does not hold because the symbol b appears below two different splits.}
    \label{fig_notRbranching_3}
    \end{subfigure}
    \caption{Strand diagrams that are not $R$-branching (\cref{def.r.branching}.}
    \label{fig_notRbranching}
\end{figure}

These conditions are invariant under renaming symbols, so this definition makes sense.
\cref{fig_notRbranching} depicts examples of strand diagrams that are not $R$-branching, using the edge replacement rules of the airplane edge replacement system.

\begin{remark}
Although these conditions may seem artificial at first, they arise from the following observations, while keeping in mind that splits and merges correspond to edge expansions in the domain and in the range, respectively.
\begin{enumerate}
    \item The first condition essentially states that splits and merges are shaped and labeled like the branchings of forest expansions based on the edge replacement rules $R$;
    \item The second condition ensures that Type 2 reductions, which will be described in the following \cref{sub.SDs.reduction}, can be performed whenever merges are followed by splits (this situation does not arise when pasting the domain and the range of a forest pair diagram as explained at the beginning of \cref{sub.SDs});
    \item The third condition tells us that, in both the underlying domain and range graph expansions, vertices generated by an edge expansion have distinct names that are not used for other vertices.
\end{enumerate}
These conditions are verified in strand diagrams obtained by gluing together the domain and the range forests of a forest pair diagrams as explained at the beginning of \cref{sub.SDs}.
Thus, strand diagrams obtained by gluing the two forests of a forest pair diagram based on the edge replacement rules $(R, \mathrm{C})$ produces an $R$-branching strand diagram.
\end{remark}

In \cref{LEM cut SDs} and the discussion that follows it we will see that $R$-branching strand diagrams, once they have been reduced (as explained later in the next \cref{sub.SDs.reduction}), can be ``cut'' in a unique way that returns a forest pair diagram.
This means that the family of $R$-branching strand diagrams is really the one that we should be looking into in order to study rearrangement groups.

\begin{definition}
\label{def.replacement.groupoid}
The \textbf{replacement groupoid} associated to the set of edge replacement rules $(R, \mathrm{C})$ is the set of all $R$-branching strand diagrams.
\end{definition}

The reasons as to why this is a groupoid are given later in \cref{sub.SDs.composition}, where we define the composition of $R$-branching strand diagrams.

\subsection{Reductions of Strand Diagrams}
\label{sub.SDs.reduction}

An $X$-merge being followed or preceded by a properly aligned $X$-split produce a reduction of the strand diagram in the following way.

\begin{figure}\centering
    \begin{subfigure}[t]{.45\textwidth}\centering
    \begin{tikzpicture}[font=\footnotesize]
        \draw (0,0) -- node[Green,left,align=center]{i\\t} (0,-1) node[black,circle,fill,inner sep=1.25]{};
        \draw (0,-1) to[out=180,in=90,looseness=1.2] (-.85,-2) node[Orange,xshift=-.3cm,yshift=0cm,align=center]{v\textsubscript{1}\\w\textsubscript{1}};
        \draw[dotted] (0,-1) -- (0,-1.5);
        \draw (0,-1) to[out=0,in=90,looseness=1.2] (.85,-2) node[Plum,xshift=-.3cm,yshift=0cm,align=center]{v\textsubscript{n}\\w\textsubscript{n}};
        \draw (0,-2) node{\Large$\dots$};
        \draw (0,-3) to[out=180,in=270,looseness=1.2] (-.85,-2);
        \draw[dotted] (0,-3) -- (0,-2.5);
        \draw (0,-3) to[out=0,in=270,looseness=1.2] (.85,-2);
        \draw (0,-4) -- node[Green,left,align=center]{i\\t} (0,-3) node[black,circle,fill,inner sep=1.25]{};
        \draw[thick,-stealth] (1.3,-2) -- (1.85,-2);
        \draw (2.55,-.5) -- node[Green,left,align=center]{i\\t} (2.55,-3.5);
    \end{tikzpicture}
    \label{fig_SD_red1}
    \caption{Type 1 reduction: an $X$-split on top of an $X$-merge produces a single strand.}
    \end{subfigure}
    \hfill
    \begin{subfigure}[t]{.5\textwidth}\centering
    \begin{tikzpicture}[font=\footnotesize]
        \draw (0,-1) to[out=180,in=270,looseness=1.2] (-.85,0)  node[Orange,xshift=-.3cm,yshift=-.25cm,align=center]{v\textsubscript{1}\\w\textsubscript{1}};
        \draw[dotted] (0,-1) -- (0,-.5);
        \draw (0,-1) to[out=0,in=270,looseness=1.2] (.85,0)  node[Plum,xshift=-.3cm,yshift=-.25cm,align=center]{v\textsubscript{n}\\w\textsubscript{n}};
        \draw (0,-.15) node{\Large$\dots$};
        \draw (0,-2.5) node[black,circle,fill,inner sep=1.25]{} -- node[Green,left,align=center]{i\\t} (0,-1) node[black,circle,fill,inner sep=1.25]{};
        \draw (0,-2.5) to[out=180,in=90,looseness=1.2] (-.85,-3.5) node[Orange,xshift=-.3cm,yshift=.3cm,align=center]{v\textsubscript{1}\\w\textsubscript{1}};
        \draw[dotted] (0,-2.5) -- (0,-3);
        \draw (0,-2.5) to[out=0,in=90,looseness=1.2] (.85,-3.5) node[Plum,xshift=-.25cm,yshift=.25cm,align=center]{v\textsubscript{n}\\w\textsubscript{n}};
        \draw (0,-3.35) node{\Large$\dots$};
        \draw[thick,-stealth] (1,-1.75) -- (1.55,-1.75);
        \draw (2.5,-.4) -- node[Orange,left,align=center]{v\textsubscript{1}\\w\textsubscript{1}} (2.5,-3.1);
        \draw[dotted] (3.35,-.5) -- (3.35,-1.25);
        \draw (3.35,-1.75) node{\Large$\dots$};
        \draw[dotted] (3.35,-3) -- (3.35,-2.25);
        \draw (4.2,-.4) -- node[Plum,left,align=center]{v\textsubscript{n}\\w\textsubscript{n}} (4.2,-3.1);
    \end{tikzpicture}
    \caption{Type 2 reduction: an $X$-merge on top of an $X$-split produces multiple strands.}
    \label{fig_SD_red2}
    \end{subfigure}
    \caption{A schematic depiction of reductions of $R$-branching strand diagrams. Each strand should be labeled and colored according to the edge replacement rules.}
    \label{fig_SD_reductions}
\end{figure}
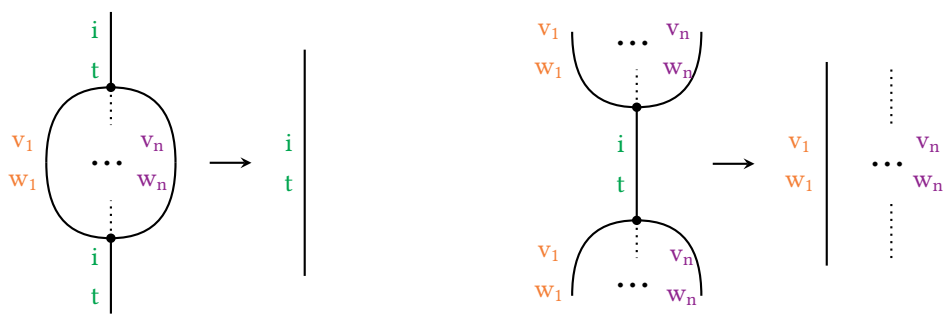

\begin{definition}
\label{def.SD.reductions}
A \textbf{reduction} of an $R$-branching strand diagram is either of the two types of moves shown in \cref{fig_SD_reductions}.
An $R$-branching strand diagram is \textbf{reduced} if no reduction can be performed on it.
Two $R$-branching strand diagrams are \textbf{equivalent} if one can be obtained from the other by a sequence of reductions and inverse reductions.
\end{definition}

Observe that the reduction of an $R$-branching strand diagram is an $R$-branching strand diagram.
Also note that reductions of a forest pair diagram (discussed in \cref{sec.forest.pair.diagrams}) correspond to Type 1 reductions of the associated strand diagram.
Instead, Type 2 reductions cannot appear when gluing together forests from a forest pair diagram (as seen at the beginning of \cref{sub.SDs}), but they often emerge in the composition of $R$-branching strand diagrams, which is described later in \cref{sub.SDs.composition}.

\medskip %layout
\begin{lemma}
\label{LEM reduced SD}
Every $R$-branching strand diagram is equivalent to a unique reduced $R$-branching strand diagram.
\end{lemma}

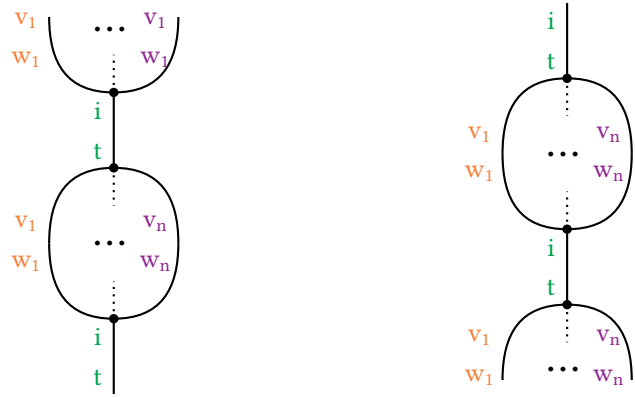
\begin{figure}\centering
    \begin{subfigure}[t]{.4\textwidth}\centering
    \begin{tikzpicture}[font=\footnotesize,scale=1]
        \draw (0,0) to[out=180,in=270,looseness=1.2] (-.85,1) node[Orange,xshift=-.3cm,yshift=-.3cm,align=center]{v\textsubscript{1}\\w\textsubscript{1}};
        \draw[dotted] (0,0) -- (0,.5);
        \draw (0,.85) node{\Large$\dots$};
        \draw (0,0) to[out=0,in=270,looseness=1.2] (.85,1) node[Plum,xshift=-.3cm,yshift=-.3cm,align=center]{v\textsubscript{1}\\w\textsubscript{1}};
        \draw (0,0) node[black,circle,fill,inner sep=1.25]{} -- node[Green,left,align=center]{i\\t} (0,-1) node[black,circle,fill,inner sep=1.25]{};
        \draw (0,-1) to[out=180,in=90,looseness=1.2] (-.85,-2) node[Orange,xshift=-.3cm,yshift=0cm,align=center]{v\textsubscript{1}\\w\textsubscript{1}};
        \draw[dotted] (0,-1) -- (0,-1.5);
        \draw (0,-1) to[out=0,in=90,looseness=1.2] (.85,-2) node[Plum,xshift=-.3cm,yshift=0cm,align=center]{v\textsubscript{n}\\w\textsubscript{n}};
        \draw (0,-2) node{\Large$\dots$};
        \draw (0,-3) to[out=180,in=270,looseness=1.2] (-.85,-2);
        \draw[dotted] (0,-3) -- (0,-2.5);
        \draw (0,-3) to[out=0,in=270,looseness=1.2] (.85,-2);
        \draw (0,-4) -- node[Green,left,align=center]{i\\t} (0,-3) node[black,circle,fill,inner sep=1.25]{};
    \end{tikzpicture}
    \end{subfigure}
    \begin{subfigure}[t]{.4\textwidth}\centering
    \begin{tikzpicture}[font=\footnotesize,scale=1]
        \draw (0,0) -- node[Green,left,align=center]{i\\t} (0,-1) node[black,circle,fill,inner sep=1.25]{};
        \draw (0,-1) to[out=180,in=90,looseness=1.2] (-.85,-2) node[Orange,xshift=-.3cm,yshift=0cm,align=center]{v\textsubscript{1}\\w\textsubscript{1}};
        \draw[dotted] (0,-1) -- (0,-1.5);
        \draw (0,-1) to[out=0,in=90,looseness=1.2] (.85,-2) node[Plum,xshift=-.3cm,yshift=0cm,align=center]{v\textsubscript{n}\\w\textsubscript{n}};
        \draw (0,-2) node{\Large$\dots$};
        \draw (0,-3) to[out=180,in=270,looseness=1.2] (-.85,-2);
        \draw[dotted] (0,-3) -- (0,-2.5);
        \draw (0,-3) to[out=0,in=270,looseness=1.2] (.85,-2);
        \draw (0,-4) node[black,circle,fill,inner sep=1.25]{} -- node[Green,left,align=center]{i\\t} (0,-3) node[black,circle,fill,inner sep=1.25]{};
        \draw (0,-4) to[out=180,in=90,looseness=1.2] (-.85,-5) node[Orange,xshift=-.3cm,yshift=.3cm,align=center]{v\textsubscript{1}\\w\textsubscript{1}};
        \draw[dotted] (0,-4) -- (0,-4.5);
        \draw (0,-4.85) node{\Large$\dots$};
        \draw (0,-4) to[out=0,in=90,looseness=1.2] (.85,-5) node[Plum,xshift=-.3cm,yshift=.3cm,align=center]{v\textsubscript{n}\\w\textsubscript{n}};
    \end{tikzpicture}
    \end{subfigure}
    \caption{Strand diagrams where two non-disjoint reductions of Type 1 and 2 are possible.}
    \label{fig_strand_diamond}
\end{figure}

\begin{proof}
Using the standard argument, by \hyperref[lem.diamond]{Newman's Diamond Lemme} it suffices to prove that the directed graph whose set of vertices is the set of $R$-branching strand diagrams and whose edges are reductions produce a terminating and locally confluent rewriting system (recall \cref{def.locally.confluent,def.terminating}).

Since each reduction strictly decreases the number of strands in a diagram, it is clear that the directed graph is terminating, so we only need to prove the local confluence.

Suppose $f, g$ and $h$ are $R$-branching strand diagrams such that $f \overset{A}{\longrightarrow} g$ and $f \overset{B}{\longrightarrow} h$ are distinct reductions.
If both $A$ and $B$ are reductions of the same type, then they must be disjoint, by which we mean that they concern different splits and merges of $f$.
In this case it is clear that the reduction $B$ can be applied to $g$ and the reduction $A$ can be applied to $h$, producing the same $R$-branching strand diagram.
If instead the reductions $A$ and $B$ are of different type, say that $A$ is of Type 1 and $B$ is of Type 2, then they are either disjoint, in which case we are done for the same reason we just discussed, or we are in one of the situations described in \cref{fig_strand_diamond}.
Observe that, in either of the two situations, whether we decide to perform the reduction $A$ or $B$, the resulting strand diagram is the same, i.e., $g = h$, and so we are done.
\end{proof}

\subsection{Rearrangements as Strand Diagrams}
\label{sub.SDs.are.rearrangements}

In this subsection we will see that the replacement groupoid associated to $R$ contains every rearrangement of any edge replacement system based on the same edge replacement rules $(R, \mathrm{C})$, independent of the base graph.
After describing the composition of $R$-branching strand diagrams in \cref{sub.SDs.composition}, this will imply that every rearrangement group is a subgroupoid of the replacement groupoid associated to its edge replacement rules.

\medskip %layout
\begin{lemma}\label{LEM cut SDs}
Each reduced $R$-branching strand diagram can be cut in a unique way such that the two resulting parts of the diagram are two forest expansions of the edge replacement systems $(A_0, R, \mathrm{C})$ and $(B_0, R, \mathrm{C})$, where $A_0$ and $B_0$ are the base graphs corresponding to the labelings of the sources and sinks, respectively.
Moreover, there is a unique way of gluing them back together into the original strand diagram, and this is described by a unique bijection between the leaves of the two forest expansions that preserves labels, meaning that this permutation is a graph isomorphism between the leaf graph of the domain forest and the one of the range forest.
\end{lemma}

\begin{proof}
Let $f$ be a reduced $R$-branching strand diagram.
Consider the family of all directed paths starting from a source and ending at a sink of the diagram.
We claim that each of these paths contains a unique strand $s$ such that the strands preceding $s$ in the path all belong to splits and not merges, whereas the strands following $s$ all belong to merges and not splits.
Indeed, suppose that some such path does not contain any such strand $s$.
Then there is some strand that is preceded by a merge and followed by a split.
Because of the second requirement of the definition of $R$-branching (\cref{def.r.branching}), these merge and split must be labeled in the same way.
This is a contradiction, because it gives place to a reduction of type 2, but we required $f$ to be reduced.

The uniqueness of such $s$ is evident: the ``cut'' consists of severing each of the unique strands $s$ we just found.
This results in two forest expansions as in the statement of the Lemma: each of the two forests is shaped and labeled like a forest expansion because of the first and third requirement of the definition of $R$-branching (\cref{def.r.branching}).
\end{proof}

This Lemma essentially tells us that each $R$-branching strand diagram corresponds to a generalized rearrangement between graph expansions obtained using same edge replacement rules $(R, \mathrm{C})$, but possibly two distinct base graphs, as we do not require the top and bottom strands to correspond to the same graph.
After being reduced, such a diagram can be cut in a unique way that produces a generalized forest pair diagram where the roots of the two forests need not represent the same base graph.

\phantomsection\label{TXT X-SDs}
Now, consider an $R$-branching strand diagram with the same number of sources and sinks and such that sources and sinks have the same labels (up to renaming symbols) and colors, in the same order.
Then sources and sinks both identify the same graph $X_0$, along with the same ordering of its edges given by the ordering of sources and sinks.
Such an $R$-branching strand diagram is called an \textbf{$\mathcal{X}$-strand diagram}, where $\mathcal{X}$ is the edge replacement system $(X_0, R, \mathrm{C})$.
By the discussion in the previous paragraph, it is clear that each reduced $\mathcal{X}$-strand diagram corresponds to a unique rearrangement of $\mathcal{X}$, as in this case the two graph expansions are realized starting from the same base graph.

Conversely, we previously noted that rearrangements can be represented by forest pair diagrams and we have already seen at the beginning of \cref{sub.SDs} that we can glue the domain and the range forests of an $\mathcal{X}$-forest pair diagram to obtain a strand diagram.
This results in an $\mathcal{X}$-strand diagram.
Then, given an edge replacement system $\mathcal{X} = (X_0, R, \mathrm{C})$, we have that, up to reductions, rearrangements of $\mathcal{X}$ correspond uniquely to $\mathcal{X}$-strand diagrams.

With the aid of \cref{LEM cut SDs}, we are now finally ready to define the composition of $R$-branching strand diagrams.

\subsection{Composition of Strand Diagrams}
\label{sub.SDs.composition}

Let $f$ and $g$ be $R$-branching strand diagrams.
We can compute their composition $f \circ g$ when the following requirements are met:
\begin{enumerate}[label=(\Alph*)]
    \item the number of sources of $f$ equals the number of sinks of $g$;
    \item the sources of $f$ and the sinks of $g$, in their given order, share the same color and, up to renaming symbols, the same labels;
\end{enumerate}
Observe that these conditions mean that the graphs represented by the sources of $f$ and the sinks of $g$ are the same, up to renaming their vertices, and the ordering of their edges is the same.

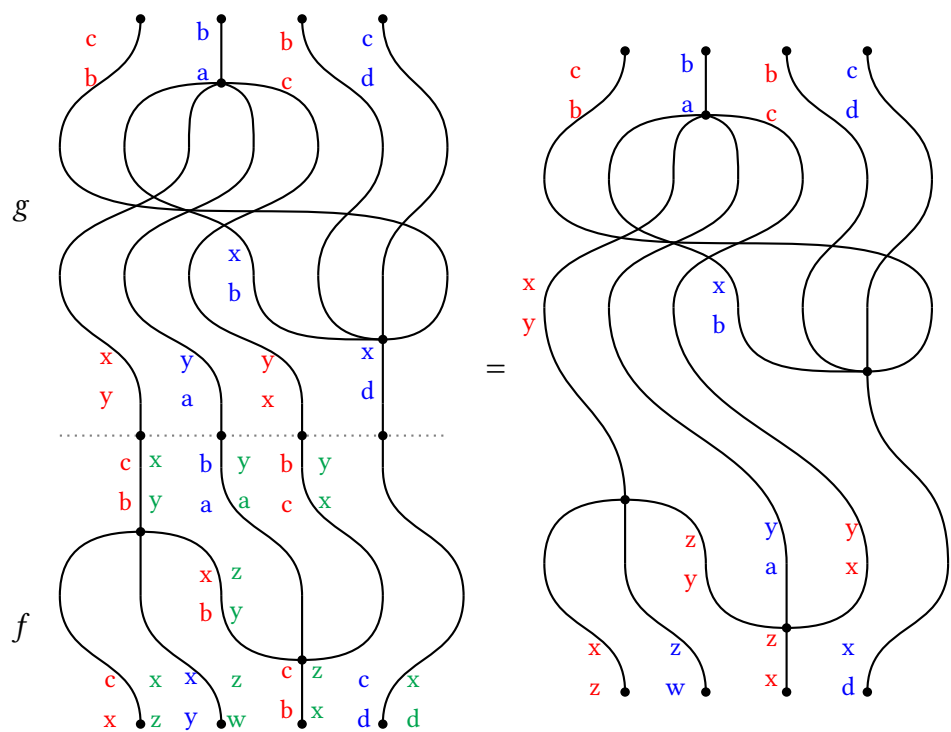
\begin{figure}\centering
\begin{tikzpicture}[font=\footnotesize,scale=.85]
\draw[dotted,gray] (-1.25,-6.5) -- (4.75,-6.5);
\begin{scope} % diagram g
    \draw (-1.85,-3) node{\normalsize$g$};
    \draw (0,0) node[black,circle,fill,inner sep=1.25]{} node[red,xshift=-.65cm,yshift=-.55cm,align=center]{c\\b} to[out=270,in=90,looseness=1.2] (-1.25,-2);
    \draw (1.25,0) node[black,circle,fill,inner sep=1.25]{} -- node[blue,left,align=center]{b\\a} (1.25,-1) node[black,circle,fill,inner sep=1.25]{};
    \draw (2.5,0) node[black,circle,fill,inner sep=1.25]{} node[red,xshift=-.2cm,yshift=-.55cm,align=center]{b\\c} to[out=270,in=90,looseness=1.2] (3.75,-2);
    \draw (3.75,0) node[black,circle,fill,inner sep=1.25]{} node[blue,xshift=-.2cm,yshift=-.55cm,align=center]{c\\d} to[out=270,in=90,looseness=1.2] (4.75,-2);
    \begin{scope}[xshift=1.25cm]
    \draw (0,-1) to[out=180,in=90,looseness=1.2] (-1.5,-2);
    \draw (0,-1) to[out=195,in=90,looseness=1.2] (-.5,-2);
    \draw (0,-1) to[out=345,in=90,looseness=1.2] (.5,-2);
    \draw (0,-1) to[out=0,in=90,looseness=1.2] (1.5,-2);
    \end{scope}
    \draw (-1.25,-2) to[out=270,in=90,looseness=.8] (4.75,-4);
    \draw (-.25,-2) to[out=270,in=90,looseness=1.2] (1.75,-4);
    \draw (.75,-2) to[out=270,in=90,looseness=1] (-1.25,-4);
    \draw (1.75,-2) to[out=270,in=90,looseness=1] (-.25,-4);
    \draw (2.75,-2) to[out=270,in=90,looseness=1] (.75,-4);
    \draw (3.75,-2) to[out=270,in=90,looseness=1.2] (2.75,-4);
    \draw (4.75,-2) to[out=270,in=90,looseness=1.2] (3.75,-4);
    \begin{scope}[yshift=-4cm]
    \draw (-1.25,0) to[out=270,in=90,looseness=1.2] (0,-2) node[red,xshift=-.45cm,yshift=.3cm,align=center]{x\\y};
    \draw (-.25,0) to[out=270,in=90,looseness=1.2] (1.25,-2) node[blue,xshift=-.45cm,yshift=.3cm,align=center]{y\\a};
    \draw (.75,0) to[out=270,in=90,looseness=1.2] (2.5,-2) node[red,xshift=-.45cm,yshift=.3cm,align=center]{y\\x};
    \draw (3.75,-1) node[black,circle,fill,inner sep=1.25]{} node[blue,xshift=-.2cm,yshift=-.45cm,align=center]{x\\d} -- (3.75,-2);
    \begin{scope}[xshift=3.75cm]
    \draw (0,-1) to[out=180,in=270,looseness=1.2] (-2,0) node[blue,xshift=-.25cm,yshift=0cm,align=center]{x\\b};
    \draw (0,-1) to[out=180,in=270,looseness=1.2] (-1,0);
    \draw (0,-1) to[out=90,in=270,looseness=1.2] (0,0);
    \draw (0,-1) to[out=0,in=270,looseness=1.2] (1,0);
    \end{scope}
    \end{scope}
\end{scope}
    \draw (0,-6) -- node[black,circle,fill,inner sep=1.25]{} (0,-7);
    \draw (1.25,-6) -- node[black,circle,fill,inner sep=1.25]{} (1.25,-7);
    \draw (2.5,-6) -- node[black,circle,fill,inner sep=1.25]{} (2.5,-7);
    \draw (3.75,-6) -- node[black,circle,fill,inner sep=1.25]{} (3.75,-7);
\begin{scope}[yshift=-7cm] % diagram f
    \draw (-1.85,-2.5) node{\normalsize$f$};
    \draw (0,0) -- node[red,xshift=-.2cm,yshift=.2cm,align=center]{c\\b} node[Green,xshift=.2cm,yshift=.2cm,align=center]{x\\y} (0,-1) node[black,circle,fill,inner sep=1.25]{};
    \draw (1.25,0) node[blue,xshift=-.2cm,yshift=-.2cm,align=center]{b\\a} node[Green,xshift=.3cm,yshift=-.2cm,align=center]{y\\a} to[out=270,in=90,looseness=1.2] (2.5,-2);
    \draw (2.5,0) node[red,xshift=-.2cm,yshift=-.2cm,align=center]{b\\c} node[Green,xshift=.3cm,yshift=-.2cm,align=center]{y\\x} to[out=270,in=90,looseness=1.2] (3.75,-2);
    \draw (3.75,0) to[out=270,in=90,looseness=1.2] (5,-2);
    \draw (0,-1) to[out=180,in=90,looseness=1.2] (-1.25,-2);
    \draw (0,-1) -- (0,-2);
    \draw (0,-1) to[out=0,in=90,looseness=1.2] (1.25,-2) node[red,xshift=-.2cm,yshift=0cm,align=center]{x\\b} node[Green,xshift=.2cm,yshift=0cm,align=center]{z\\y};
    \begin{scope}[yshift=-2cm]
    \draw (0,-2) node[red,xshift=-.4cm,yshift=.3cm,align=center]{c\\x} node[Green,xshift=.2cm,yshift=.3cm,align=center]{x\\z} node[black,circle,fill,inner sep=1.25]{} to[out=90,in=270,looseness=1.2] (-1.25,0);
    \draw (1.25,-2) node[blue,xshift=-.4cm,yshift=.3cm,align=center]{x\\y} node[Green,xshift=.2cm,yshift=.3cm,align=center]{z\\w} node[black,circle,fill,inner sep=1.25]{} to[out=90,in=270,looseness=1.2] (0,0);
    \draw (2.5,-2) node[black,circle,fill,inner sep=1.25]{} -- node[red,xshift=-.2cm,align=center]{c\\b} node[Green,xshift=.2cm,align=center]{z\\x} (2.5,-1) node[black,circle,fill,inner sep=1.25]{};
    \draw (3.75,-2) node[black,circle,fill,inner sep=1.25]{} node[blue,xshift=-.25cm,yshift=.3cm,align=center]{c\\d} node[Green,xshift=.4cm,yshift=.3cm,align=center]{x\\d} to[out=90,in=270,looseness=1.2] (5,0);
    \draw (2.5,-1) to[out=180,in=270,looseness=1.2] (1.25,0);
    \draw (2.5,-1) -- (2.5,0);
    \draw (2.5,-1) to[out=0,in=270,looseness=1.2] (3.75,0);
    \end{scope}
\end{scope}
    \draw (5.5,-5.5) node{\normalsize$=$};
\begin{scope}[yshift=-.5cm,xshift=7.5cm] % composition
    \draw (0,0) node[black,circle,fill,inner sep=1.25]{} node[red,xshift=-.65cm,yshift=-.55cm,align=center]{c\\b} to[out=270,in=90,looseness=1.2] (-1.25,-2);
    \draw (1.25,0) node[black,circle,fill,inner sep=1.25]{} -- node[blue,left,align=center]{b\\a} (1.25,-1) node[black,circle,fill,inner sep=1.25]{};
    \draw (2.5,0) node[black,circle,fill,inner sep=1.25]{} node[red,xshift=-.2cm,yshift=-.55cm,align=center]{b\\c} to[out=270,in=90,looseness=1.2] (3.75,-2);
    \draw (3.75,0) node[black,circle,fill,inner sep=1.25]{} node[blue,xshift=-.2cm,yshift=-.55cm,align=center]{c\\d} to[out=270,in=90,looseness=1.2] (4.75,-2);
    \begin{scope}[xshift=1.25cm]
    \draw (0,-1) to[out=180,in=90,looseness=1.2] (-1.5,-2);
    \draw (0,-1) to[out=195,in=90,looseness=1.2] (-.5,-2);
    \draw (0,-1) to[out=345,in=90,looseness=1.2] (.5,-2);
    \draw (0,-1) to[out=0,in=90,looseness=1.2] (1.5,-2);
    \end{scope}
    \draw (-1.25,-2) to[out=270,in=90,looseness=.8] (4.75,-4);
    \draw (-.25,-2) to[out=270,in=90,looseness=1.2] (1.75,-4);
    \draw (.75,-2) to[out=270,in=90,looseness=1] (-1.25,-4) node[red,xshift=-.2cm,yshift=0cm,align=center]{x\\y};
    \draw (1.75,-2) to[out=270,in=90,looseness=1] (-.25,-4);
    \draw (2.75,-2) to[out=270,in=90,looseness=1] (.75,-4);
    \draw (3.75,-2) to[out=270,in=90,looseness=1.2] (2.75,-4);
    \draw (4.75,-2) to[out=270,in=90,looseness=1.2] (3.75,-4);
    \begin{scope}[yshift=-4cm]
    \draw (-1.25,0) to[out=270,in=90,looseness=1.2] (0,-3) node[black,circle,fill,inner sep=1.25]{};
    \draw (-.25,0) to[out=270,in=90,looseness=1.1] (2.5,-4) node[blue,xshift=-.2cm,yshift=.2cm,align=center]{y\\a};
    \draw (.75,0) to[out=270,in=90,looseness=1] (3.75,-4) node[red,xshift=-.2cm,yshift=.2cm,align=center]{y\\x};
    \draw (3.75,-1) node[black,circle,fill,inner sep=1.25]{} to[out=270,in=90,looseness=1.4] (5,-4);
    \begin{scope}[xshift=3.75cm]
    \draw (0,-1) to[out=180,in=270,looseness=1.2] (-2,0) node[blue,xshift=-.25cm,yshift=0cm,align=center]{x\\b};
    \draw (0,-1) to[out=180,in=270,looseness=1.2] (-1,0);
    \draw (0,-1) to[out=90,in=270,looseness=1.2] (0,0);
    \draw (0,-1) to[out=0,in=270,looseness=1.2] (1,0);
    \end{scope}
    \end{scope}
    \begin{scope}[yshift=-6cm]
    \draw (0,-1) to[out=180,in=90,looseness=1.2] (-1.25,-2);
    \draw (0,-1) -- (0,-2);
    \draw (0,-1) to[out=0,in=90,looseness=1.2] (1.25,-2) node[red,xshift=-.2cm,yshift=0cm,align=center]{z\\y};
    \begin{scope}[yshift=-2cm]
    \draw (0,-2) node[red,xshift=-.4cm,yshift=.3cm,align=center]{x\\z} node[black,circle,fill,inner sep=1.25]{} to[out=90,in=270,looseness=1.2] (-1.25,0);
    \draw (1.25,-2) node[blue,xshift=-.4cm,yshift=.3cm,align=center]{z\\w} node[black,circle,fill,inner sep=1.25]{} to[out=90,in=270,looseness=1.2] (0,0);
    \draw (2.5,-2) node[black,circle,fill,inner sep=1.25]{} -- node[red,xshift=-.2cm,align=center]{z\\x} (2.5,-1) node[black,circle,fill,inner sep=1.25]{};
    \draw (3.75,-2) node[black,circle,fill,inner sep=1.25]{} node[blue,xshift=-.25cm,yshift=.3cm,align=center]{x\\d} to[out=90,in=270,looseness=1.2] (5,0);
    \draw (2.5,-1) to[out=180,in=270,looseness=1.2] (1.25,0);
    \draw (2.5,-1) -- (2.5,0);
    \draw (2.5,-1) to[out=0,in=270,looseness=1.2] (3.75,0);
    \end{scope}
    \end{scope}
\end{scope}
\end{tikzpicture}
\caption{A composition $f \circ g$ of strand diagrams.}
\label{fig_strand_composition}
\end{figure}

When these requirements hold, we can always find a renaming of the symbols of $f$ (or $g$) such that the diagram obtained by gluing the sinks of $g$ with the sources of $f$ in their given order is an $R$-branching strand diagram (\cref{def.r.branching}).
Indeed, because of requirements A and B we can rename the symbols labeling the sources of $f$ so that they match the labeling of the sinks of $g$, and then we can continue renaming the diagram for $f$ downwards from its sources in the following way.
The branching strand of a merge is automatically determined by its top strands, so there is no ambiguity when renaming merges.
On the other hand, when renaming a split whose branching strand has label $(v,w,z)$, we need to distinguish between two cases:
\begin{itemize}
    \item if $(v,w,z)$ also labels some branching strand of a merge of $g$, then rename the bottom strands of the split of $f$ using the same symbols that also appear in the top strands of the merge of $g$;
    \item if $(v,w,z)$ does not label any branching strand of any merge of $g$, then rename the bottom strands of the split of $f$ using new symbols that do not appear elsewhere.
\end{itemize}
Then the composition $f \circ g$ is the strand diagram obtained by renaming $f$ as described above (or renaming $g$ in the opposite direction, from bottom to top), gluing the sinks of $g$ with the sources of $f$ in their given order and then reducing as explained in \cref{sub.SDs.reduction}.
Indeed, the strand diagram obtained by gluing $g$ and $f$ is $R$-branching (\cref{def.r.branching}), as the first condition is trivially satisfied and the second and third hold because of the steps described above.
An example is given in \cref{fig_strand_composition}, where the green labels announce the renaming of symbols that is needed on $f$.

With this definition of composition between $R$-branching strand diagrams, it is easy to check that the inverse of a diagram is the diagram obtained by reversing the direction of each strand, which makes sinks into sources and sources into sinks.
Essentially, the inverse can be computed just by drawing the original strand diagram ``upside-down''.
It is also clear that this composition is associative, thus the replacement groupoid corresponding to $R$ is really a groupoid.

\subsection{Generators of the Replacement Groupoid}
\label{SUB groupoid generators}

It is worth mentioning that, as a straightforward application of \cref{LEM cut SDs}, we can describe a generating set for the replacement groupoid.
Given a set of replacement graphs $R$, by \textbf{split diagram} (\textbf{merge diagram}) we refer to the $R$-branching strand diagram consisting of a sole split (merge) surrounded by any amount of straight strands.
By \textbf{permutation diagram} we refer to a strand diagram without splits or merges (the reason for this name is that the action of such a diagram is  that of a bijection between the set of sources and the set of sinks).

\medskip %layout
\begin{proposition}
The replacement groupoid associated to $R$ is generated by the infinite set consisting of every permutation diagram and every split diagram (or equivalently every merge diagram).
Moreover, each reduced $R$-branching strand diagram can be written uniquely as a product $M \circ P \circ S$, where $M$ is a product of merge diagrams, $P$ is a permutation diagram and $S$ is a product of split diagrams.
\end{proposition}

A similar fact was first noted in \cite{BF19}:
where simple expansion morphisms, simple contraction morphisms and base isomorphisms correspond to split, merge and permutation strand diagrams, respectively.

\section{Closed Strand Diagrams}
\label{SEC closed strand diagrams}

As was done in \cite{BM14}, we now proceed to ``close'' our strand diagrams (\cref{SUB closed strand diagrams}) and we describe new transformations that can be performed on these closed diagrams (\cref{SUB CSD transformations}).
Then we study the possible uniqueness of reduced closed strand diagrams (\cref{SUB confluent reductions}), which will be crucial when stating a general result about the conjugacy problem in the following \cref{SEC conjugacy problem}.
Finally, in \cref{SUB stable and vanishing} we describe what happens when certain transformations move the ``base line'' of a closed diagram back to where it started and possibly produce new configurations of labels.
This lays the groundwork for the algorithm to solve the conjugacy problem in \cref{SUB algorithm}.

\subsection{Definition of Closed Diagrams}
\label{SUB closed strand diagrams}

\begin{figure}\centering
\begin{tikzpicture}[font=\footnotesize,scale=.875]
    \useasboundingbox (-8,-5.5) rectangle (3,5.5);
    %\draw[help lines,step=.1cm] (-8,-5.5) grid (3,5.5);
    %
    \begin{scope}[xshift=-6cm,yshift=3cm]
    \draw (0,0) node[black,circle,fill,inner sep=1.25]{} node[red,xshift=-.65cm,yshift=-.55cm,align=center]{c\\b} to[out=270,in=90,looseness=1.2] (-1.25,-2);
    \draw (1.25,0) node[black,circle,fill,inner sep=1.25]{} -- node[blue,left,align=center]{b\\a} (1.25,-1) node[black,circle,fill,inner sep=1.25]{};
    \draw (2.5,0) node[black,circle,fill,inner sep=1.25]{} node[red,xshift=-.2cm,yshift=-.55cm,align=center]{b\\c} to[out=270,in=90,looseness=1.2] (3.75,-2);
    \draw (3.75,0) node[black,circle,fill,inner sep=1.25]{} node[blue,xshift=-.2cm,yshift=-.55cm,align=center]{c\\d} to[out=270,in=90,looseness=1.2] (4.75,-2);
    \begin{scope}[xshift=1.25cm]
    \draw (0,-1) to[out=180,in=90,looseness=1.2] (-1.5,-2);
    \draw (0,-1) to[out=195,in=90,looseness=1.2] (-.5,-2);
    \draw (0,-1) to[out=345,in=90,looseness=1.2] (.5,-2);
    \draw (0,-1) to[out=0,in=90,looseness=1.2] (1.5,-2);
    \end{scope}
    \draw (-1.25,-2) to[out=270,in=90,looseness=.8] (4.75,-4);
    \draw (-.25,-2) to[out=270,in=90,looseness=1.2] (1.75,-4);
    \draw (.75,-2) to[out=270,in=90,looseness=1] (-1.25,-4);
    \draw (1.75,-2) to[out=270,in=90,looseness=1] (-.25,-4);
    \draw (2.75,-2) to[out=270,in=90,looseness=1] (.75,-4);
    \draw (3.75,-2) to[out=270,in=90,looseness=1.2] (2.75,-4);
    \draw (4.75,-2) to[out=270,in=90,looseness=1.2] (3.75,-4);
    \begin{scope}[yshift=-4cm]
    \draw (-1.25,0) to[out=270,in=90,looseness=1.2] (0,-2) node[red,xshift=-.45cm,yshift=.3cm,align=center]{x\\y};
    \draw (-.25,0) to[out=270,in=90,looseness=1.2] (1.25,-2) node[blue,xshift=-.45cm,yshift=.3cm,align=center]{y\\a};
    \draw (.75,0) to[out=270,in=90,looseness=1.2] (2.5,-2) node[red,xshift=-.45cm,yshift=.3cm,align=center]{y\\x};
    \draw (3.75,-1) node[black,circle,fill,inner sep=1.25]{} node[blue,xshift=-.2cm,yshift=-.45cm,align=center]{x\\d} -- (3.75,-2);
    \begin{scope}[xshift=3.75cm]
    \draw (0,-1) to[out=180,in=270,looseness=1.2] (-2,0) node[blue,xshift=-.25cm,yshift=0cm,align=center]{x\\b};
    \draw (0,-1) to[out=180,in=270,looseness=1.2] (-1,0);
    \draw (0,-1) to[out=90,in=270,looseness=1.2] (0,0);
    \draw (0,-1) to[out=0,in=270,looseness=1.2] (1,0);
    \end{scope}
    \end{scope}
    \end{scope}
    \draw[dashed] (0,0) to[out=90,in=0] (-2,3) -- (-8,3);
    \draw[fill=gray] (0,0) circle (.25);
    \draw (-2.25,-3) to[out=270,in=270,looseness=1.2] (1.5,-1) -- (1.5,1) to[out=90,in=90,looseness=1.2] (-2.25,3);
    \draw (-3.5,-3) to[out=270,in=270,looseness=1.2] (2,-1.167) -- (2,1.167) to[out=90,in=90,looseness=1.2] (-3.5,3);
    \draw (-4.75,-3) to[out=270,in=270,looseness=1.2] (2.5,-1.333) -- (2.5,1.333) to[out=90,in=90,looseness=1.2] (-4.75,3);
    \draw (-6,-3) to[out=270,in=270,looseness=1.2] (3,-1.5) -- (3,1.5) to[out=90,in=90,looseness=1.2] (-6,3);
\end{tikzpicture}
\caption{The closed strand diagram obtained from the strand diagram depicted in \cref{fig_strand}. The base line is represented by the dashed line.}
\label{fig_CSD}
\end{figure}
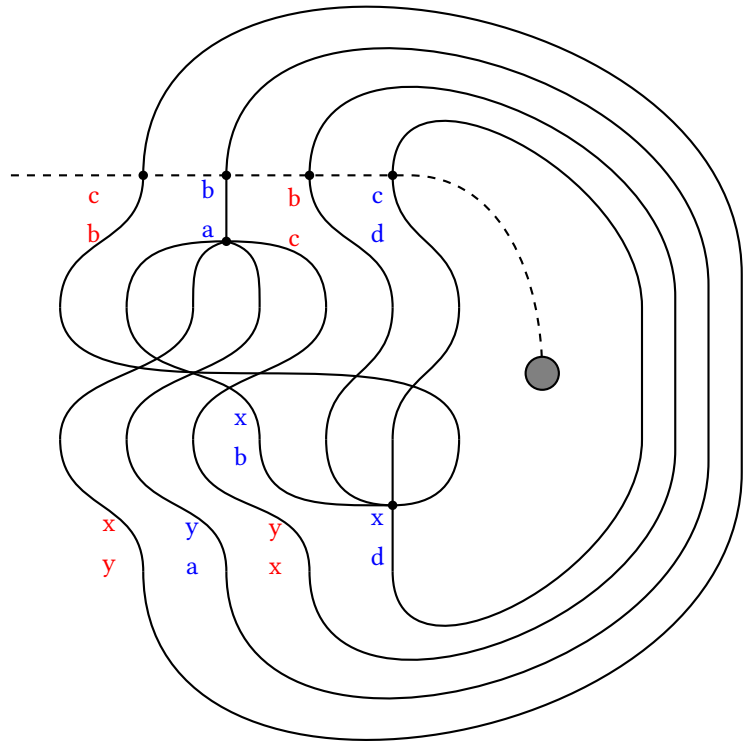

Fix an edge replacement system $\mathcal{X} = (X_0, R, \mathrm{C})$ and its rearrangement group $G_\mathcal{X}$.
As we have seen in the previous section at \cpageref{TXT X-SDs}, we can represent rearrangements as $\mathcal{X}$-strand diagrams, which are $R$-branching strand diagrams (\cref{def.r.branching}) whose sources and sinks both represent the base graph $X_0$ of $\mathcal{X}$ (along with its given ordering of edges).
Precisely because of this last condition, we can essentially ``close'' an $\mathcal{X}$-strand diagram around a unique ``hole'' by attaching each sink to a source in their given order as done in \cref{fig_CSD}.
Labels above and below need not be the same, but they will be the same up to some renaming of symbols.
A formal definition would essentially be the same as that of a strand diagram (\cref{def.SDs}), except that we do not require acyclicity, we do not allow univalent vertices (sources and sinks of out or in-degree $1$) and instead we allow the existence of special vertices of in-degree and out-degree equal to 1 originating from the gluing.

\phantomsection\label{txt.base.things}
More generally, we can close in this manner any $R$-branching strand diagram whose sources and sinks both represent isomorphic graphs, the isomorphism being given by the ordering of sources and sinks.
Diagrams obtained in this manner are called \textbf{$R$-branching closed strand diagrams} (or simply \textit{closed strand diagrams} or \textit{closed diagrams}, if there is no ambiguity in omitting $R$).
Those vertices that were glued (originally sources and sinks) are called \textbf{base points} and the ordered tuple of gluing points is called the \textbf{base line}.
Since sinks and sources represent the same graph, it is natural to refer to that graph as the \textbf{base graph} of the closed diagram.
If $f$ is an $R$-branching strand diagram, we denote by $\llbracket f \rrbracket$ its closure.

\phantomsection\label{txt.opening.closed.strand.diagrams}
A closed strand diagrams $\eta$ can be cut at the base line, producing a unique ``open'' strand diagram $o(\eta)$. In this sense, a closed strand diagram contains the same data as the original $R$-branching strand diagram.
However, as we will see in the next \cref{SEC conjugacy problem}, conjugacy has a natural way of being represented by reductions of parallel strands, shifts of the base line and permutations of the base points.

\subsection{Transformations of Closed Diagrams}
\label{SUB CSD transformations}

Given an $R$-branching closed strand diagram, we can apply three different kinds of transformations: permutations, shifts and reductions.
Essentially, the first two allow permuting the base points and moving the base line, respectively, and they will be called \textit{similarities}.
Reductions, on the other hand, include the same Type 1 and Type 2 reductions of strand diagrams defined in \cref{sub.SDs.reduction} and also a new Type 3 reduction that corresponds to reducing the base graph of the rearrangement.

Each of these transformations encode some kind of (possibly trivial) modification of the base graph that corresponds some (possibly trivial) conjugation of the original rearrangement by some element of the replacement groupoid, as we will see soon.

\subsubsection{Similarities: Permutations and Shifts}
\label{sub.similarities}

\begin{figure}\centering
\begin{tikzpicture}[font=\footnotesize,scale=.875]
    \useasboundingbox (-8,-5.5) rectangle (3,5.5);
    %\draw[help lines,step=.1cm] (-8,-5.5) grid (3,5.5);
    %
    \begin{scope}[xshift=-6cm,yshift=3cm]
    \draw (0,0) node[black,circle,fill,inner sep=1.25]{} to[out=270,in=90,looseness=1.5] node[blue,xshift=-.3cm,align=center]{b\\a} (-.25,-1) node[black,circle,fill,inner sep=1.25]{};
    \draw (1.25,0) node[black,circle,fill,inner sep=1.25]{} node[red,xshift=-.2cm,yshift=-.55cm,align=center]{b\\c} to[out=270,in=90,looseness=1.2] (2.75,-2);
    \draw (2.5,0) node[black,circle,fill,inner sep=1.25]{} node[red,xshift=-.2cm,yshift=-.55cm,align=center]{c\\b} to[out=270,in=90,looseness=1.2] (3.75,-2);
    \draw (3.75,0) node[black,circle,fill,inner sep=1.25]{} node[blue,xshift=-.2cm,yshift=-.55cm,align=center]{c\\d} to[out=270,in=90,looseness=1.2] (4.75,-2);
    \draw (-.25,-1) to[out=180,in=90,looseness=1.2] (-1.25,-2);
    \draw (-.25,-1) -- (-.25,-2);
    \draw (-.25,-1) to[out=345,in=90,looseness=1.2] (.75,-2);
    \draw (-.25,-1) to[out=0,in=90,looseness=1.2] (1.75,-2);
    \draw (-1.25,-2) to[out=270,in=90,looseness=.8] (1.75,-4);
    \draw (-.25,-2) to[out=270,in=90,looseness=1.2] (.75,-4);
    \draw (.75,-2) to[out=270,in=90,looseness=1] (-1.25,-4);
    \draw (1.75,-2) to[out=270,in=90,looseness=1] (-.25,-4);
    \draw (2.75,-2) to[out=270,in=90,looseness=1] (2.75,-4);
    \draw (3.75,-2) to[out=270,in=90,looseness=1.2] (4.75,-4);
    \draw (4.75,-2) to[out=270,in=90,looseness=1.2] (3.75,-4);
    \begin{scope}[yshift=-4cm]
    \draw (-1.25,0) to[out=270,in=90,looseness=1.2] (0,-2) node[blue,xshift=-.45cm,yshift=.3cm,align=center]{y\\a};
    \draw (-.25,0) to[out=270,in=90,looseness=1.2] (1.25,-2) node[red,xshift=-.45cm,yshift=.3cm,align=center]{y\\x};
    \draw (.75,0) to[out=270,in=90,looseness=1.2] (2.5,-2) node[red,xshift=-.45cm,yshift=.3cm,align=center]{x\\y};
    \draw (3.75,-1) node[black,circle,fill,inner sep=1.25]{} node[blue,xshift=-.2cm,yshift=-.45cm,align=center]{x\\d} -- (3.75,-2);
    \begin{scope}[xshift=3.75cm]
    \draw (0,-1) to[out=180,in=270,looseness=1.2] (-2,0) node[blue,xshift=-.25cm,yshift=0cm,align=center]{x\\b};
    \draw (0,-1) to[out=180,in=270,looseness=1.2] (-1,0);
    \draw (0,-1) to[out=90,in=270,looseness=1.2] (0,0);
    \draw (0,-1) to[out=0,in=270,looseness=1.2] (1,0);
    \end{scope}
    \end{scope}
    \end{scope}
    \draw[dashed] (0,0) to[out=90,in=0] (-2,3) -- (-8,3);
    \draw[fill=gray] (0,0) circle (.25);
    \draw (-2.25,-3) to[out=270,in=270,looseness=1.2] (1.5,-1) -- (1.5,1) to[out=90,in=90,looseness=1.2] (-2.25,3);
    \draw (-3.5,-3) to[out=270,in=270,looseness=1.2] (2,-1.167) -- (2,1.167) to[out=90,in=90,looseness=1.2] (-3.5,3);
    \draw (-4.75,-3) to[out=270,in=270,looseness=1.2] (2.5,-1.333) -- (2.5,1.333) to[out=90,in=90,looseness=1.2] (-4.75,3);
    \draw (-6,-3) to[out=270,in=270,looseness=1.2] (3,-1.5) -- (3,1.5) to[out=90,in=90,looseness=1.2] (-6,3);
\end{tikzpicture}
\caption{The diagram obtained by performing a permutation on the closed strand diagram represented in \cref{fig_CSD}.}
\label{fig_CSD_perm}
\end{figure}
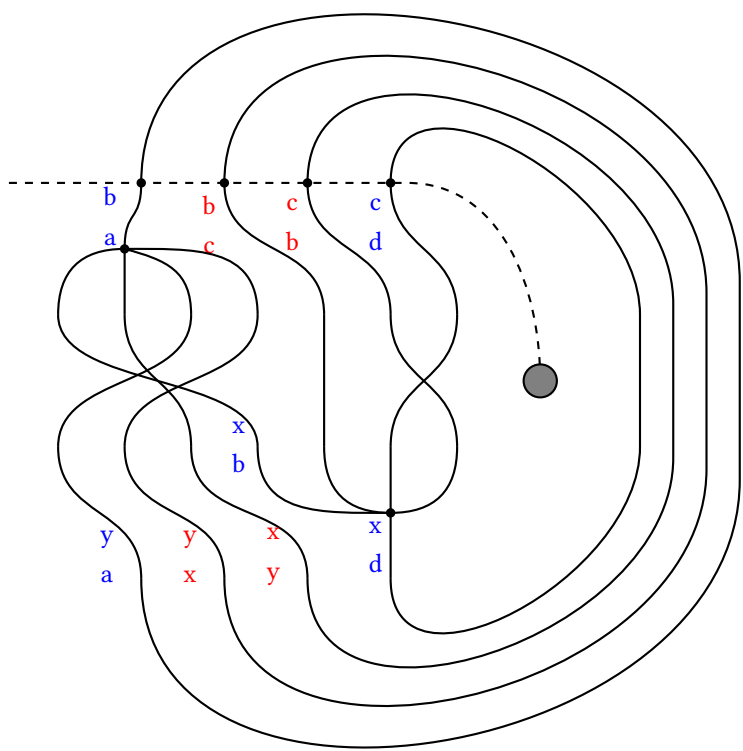

\phantomsection\label{TXT permutations}
The simplest transformation is the \textbf{permutation} of the base line, which consists simply of changing the order of the base points.
For example, \cref{fig_CSD_perm} depict the result of a permutation of the diagram portrayed in \cref{fig_CSD} (the explicit permutation is $(1 3 2) \in S_4$, where each number from $1$ to $4$ denotes a base point in their given ordering).

Clearly permutations transform $R$-branching closed diagrams into $R$-branching closed diagrams.
Transformations of this kind do not change the base graph of the diagram, but they change the order of its edges.

As we will see in \cref{PROP conjugator}, we need this kind of transformations in order to represent conjugations by permutation diagrams (which were defined in \cref{SUB groupoid generators}).

\medskip

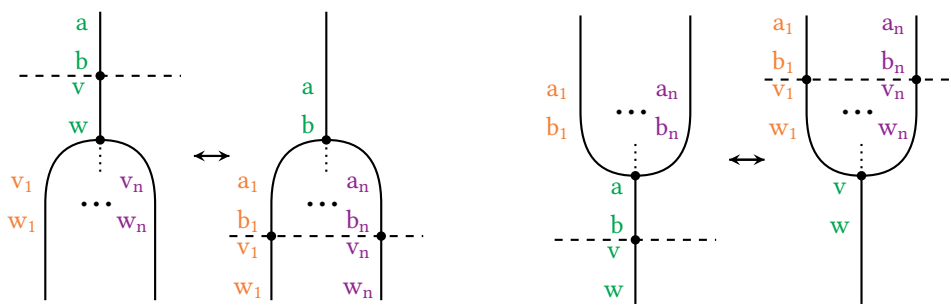
\begin{figure}\centering
\begin{subfigure}{.475\textwidth}\centering
    \begin{tikzpicture}[font=\footnotesize,scale=.85]
        \draw[dashed] (-1.25,0) -- (1.25,0);
        \draw (0,1) -- node[Green,left,align=center]{a\\b} (0,0) node[black,circle,fill,inner sep=1.25]{} -- node[Green,left,align=center]{v\\w} (0,-1) node[black,circle,fill,inner sep=1.25]{};
        \draw (0,-1) to[out=180,in=90,looseness=1.2] (-.85,-2) node[Orange,xshift=-.3cm,yshift=0cm,align=center]{v\textsubscript{1}\\w\textsubscript{1}} -- (-.85,-3.5);
        \draw[dotted] (0,-1) -- (0,-1.5);
        \draw (0,-1) to[out=0,in=90,looseness=1.2] (.85,-2) node[Plum,xshift=-.3cm,yshift=0cm,align=center]{v\textsubscript{n}\\w\textsubscript{n}} -- (.85,-3.5);
        \draw (0,-2) node{\Large$\dots$};
        \draw[thick,-stealth] (1.45,-1.25) -- (2,-1.25);
        \draw[thick,-stealth] (2,-1.25) -- (1.45,-1.25);
        \begin{scope}[xshift=3.5cm]
        \draw[dashed] (-1.5,-2.5) -- (1.5,-2.5);
        \draw (0,1) -- (0,0) -- node[Green,left,align=center]{a\\b} (0,-1) node[black,circle,fill,inner sep=1.25]{};
        \draw (0,-1) to[out=180,in=90,looseness=1.2] (-.85,-2) node[Orange,xshift=-.3cm,yshift=0cm,align=center]{a\textsubscript{1}\\b\textsubscript{1}} -- (-.85,-2.5) node[black,circle,fill,inner sep=1.25]{} -- (-.85,-3.5) node[Orange,xshift=-.3cm,yshift=.4cm,align=center]{v\textsubscript{1}\\w\textsubscript{1}};
        \draw[dotted] (0,-1) -- (0,-1.5);
        \draw (0,-1) to[out=0,in=90,looseness=1.2] (.85,-2) node[Plum,xshift=-.3cm,yshift=0cm,align=center]{a\textsubscript{n}\\b\textsubscript{n}} -- (.85,-2.5) node[black,circle,fill,inner sep=1.25]{} -- (.85,-3.5) node[Plum,xshift=-.3cm,yshift=.4cm,align=center]{v\textsubscript{n}\\w\textsubscript{n}};
        \draw (0,-2) node{\Large$\dots$};
        \end{scope}
    \end{tikzpicture}
\end{subfigure}
\begin{subfigure}{.475\textwidth}\centering
    \begin{tikzpicture}[font=\footnotesize,scale=.85,yscale=-1]
        \draw[dashed] (-1.25,0) -- (1.25,0);
        \draw (0,1) -- node[Green,left,align=center]{v\\w} (0,0) node[black,circle,fill,inner sep=1.25]{} -- node[Green,left,align=center]{a\\b} (0,-1) node[black,circle,fill,inner sep=1.25]{};
        \draw (0,-1) to[out=180,in=90,looseness=1.2] (-.85,-2) node[Orange,xshift=-.3cm,yshift=0cm,align=center]{a\textsubscript{1}\\b\textsubscript{1}} -- (-.85,-3.5);
        \draw[dotted] (0,-1) -- (0,-1.5);
        \draw (0,-1) to[out=0,in=90,looseness=1.2] (.85,-2) node[Plum,xshift=-.3cm,yshift=0cm,align=center]{a\textsubscript{n}\\b\textsubscript{n}} -- (.85,-3.5);
        \draw (0,-2) node{\Large$\dots$};
        \draw[thick,-stealth] (1.45,-1.25) -- (2,-1.25);
        \draw[thick,-stealth] (2,-1.25) -- (1.45,-1.25);
        \begin{scope}[xshift=3.5cm]
        \draw[dashed] (-1.5,-2.5) -- (1.5,-2.5);
        \draw (0,1) -- (0,0) -- node[Green,left,align=center]{v\\w} (0,-1) node[black,circle,fill,inner sep=1.25]{};
        \draw (0,-1) to[out=180,in=90,looseness=1.2] (-.85,-2) node[Orange,xshift=-.3cm,yshift=0cm,align=center]{v\textsubscript{1}\\w\textsubscript{1}} -- (-.85,-2.5) node[black,circle,fill,inner sep=1.25]{} -- (-.85,-3.5) node[Orange,xshift=-.3cm,yshift=-.4cm,align=center]{a\textsubscript{1}\\b\textsubscript{1}};
        \draw[dotted] (0,-1) -- (0,-1.5);
        \draw (0,-1) to[out=0,in=90,looseness=1.2] (.85,-2) node[Plum,xshift=-.3cm,yshift=0cm,align=center]{v\textsubscript{n}\\w\textsubscript{n}} -- (.85,-2.5) node[black,circle,fill,inner sep=1.25]{} -- (.85,-3.5) node[Plum,xshift=-.3cm,yshift=-.4cm,align=center]{a\textsubscript{n}\\b\textsubscript{n}};
        \draw (0,-2) node{\Large$\dots$};
        \end{scope}
    \end{tikzpicture}
\end{subfigure}
\caption{From left to right: two expanding shifts. From right to left: two reducing shifts.}
\label{fig_CSD_shift}
\end{figure}

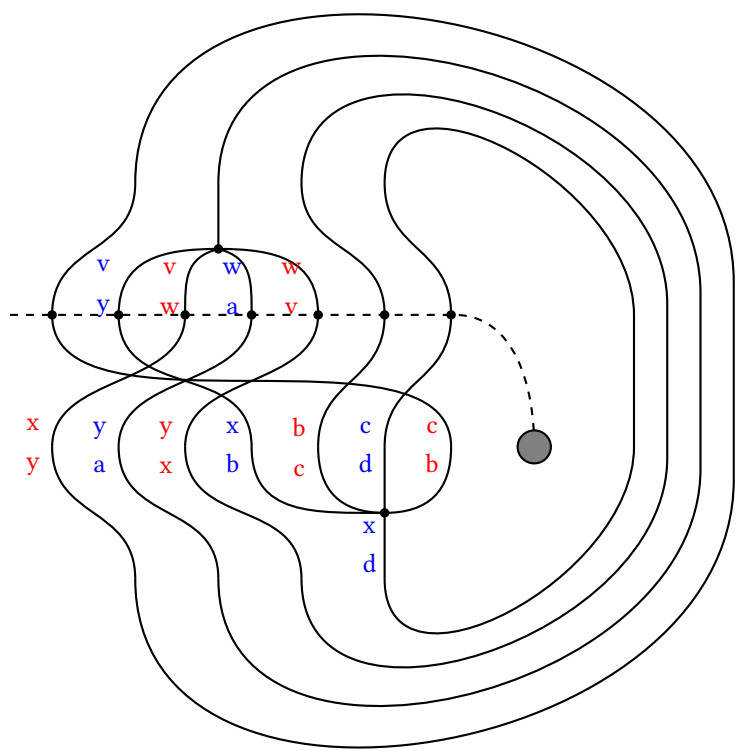
\begin{figure}\centering
\begin{tikzpicture}[font=\footnotesize,scale=.875]
    \useasboundingbox (-8,-4.5) rectangle (3,6.5);
    %\draw[help lines,step=.1cm] (-8,-4.5) grid (3,6.5);
    %
    \begin{scope}[xshift=-6cm,yshift=3cm]
    \draw (0,1) to[out=270,in=90,looseness=1.2] (-1.25,-1) node[black,circle,fill,inner sep=1.25]{};
    \draw (1.25,1) -- (1.25,0) node[black,circle,fill,inner sep=1.25]{};
    \draw (2.5,1) to[out=270,in=90,looseness=1.2] (3.75,-1) node[black,circle,fill,inner sep=1.25]{};
    \draw (3.75,1) to[out=270,in=90,looseness=1.2] (4.75,-1) node[black,circle,fill,inner sep=1.25]{};
    \begin{scope}[xshift=1.25cm]
    \draw (0,0) to[out=180,in=90,looseness=1.2] (-1.5,-1) node[black,circle,fill,inner sep=1.25]{} node[blue,xshift=-.2cm,yshift=.35cm,align=center]{v\\y};
    \draw (0,0) to[out=195,in=90,looseness=1.2] (-.5,-1) node[black,circle,fill,inner sep=1.25]{} node[red,xshift=-.2cm,yshift=.35cm,align=center]{v\\w};
    \draw (0,0) to[out=345,in=90,looseness=1.2] (.5,-1) node[black,circle,fill,inner sep=1.25]{} node[blue,xshift=-.25cm,yshift=.35cm,align=center]{w\\a};
    \draw (0,0) to[out=0,in=90,looseness=1.2] (1.5,-1) node[black,circle,fill,inner sep=1.25]{} node[red,xshift=-.35cm,yshift=.35cm,align=center]{w\\v};
    \end{scope}
    \draw (-1.25,-1) to[out=270,in=90,looseness=.8] (4.75,-3);
    \draw (-.25,-1) to[out=270,in=90,looseness=1.2] (1.75,-3);
    \draw (.75,-1) to[out=270,in=90,looseness=1] (-1.25,-3);
    \draw (1.75,-1) to[out=270,in=90,looseness=1] (-.25,-3);
    \draw (2.75,-1) to[out=270,in=90,looseness=1] (.75,-3);
    \draw (3.75,-1) to[out=270,in=90,looseness=1.2] (2.75,-3);
    \draw (4.75,-1) to[out=270,in=90,looseness=1.2] (3.75,-3);
    \begin{scope}[yshift=-3cm]
    \draw (-1.25,0) node[red,xshift=-.25cm,yshift=0cm,align=center]{x\\y} to[out=270,in=90,looseness=1.2] (0,-2);
    \draw (-.25,0) node[blue,xshift=-.25cm,yshift=0cm,align=center]{y\\a} to[out=270,in=90,looseness=1.2] (1.25,-2);
    \draw (.75,0) node[red,xshift=-.25cm,yshift=0cm,align=center]{y\\x} to[out=270,in=90,looseness=1.2] (2.5,-2);
    \draw (3.75,-1) node[black,circle,fill,inner sep=1.25]{} node[blue,xshift=-.2cm,yshift=-.45cm,align=center]{x\\d} -- (3.75,-2);
    \begin{scope}[xshift=3.75cm]
    \draw (0,-1) to[out=180,in=270,looseness=1.2] (-2,0) node[blue,xshift=-.25cm,yshift=0cm,align=center]{x\\b};
    \draw (0,-1) to[out=180,in=270,looseness=1.2] (-1,0) node[red,xshift=-.25cm,yshift=0cm,align=center]{b\\c};
    \draw (0,-1) to[out=90,in=270,looseness=1.2] (0,0) node[blue,xshift=-.25cm,yshift=-0cm,align=center]{c\\d};
    \draw (0,-1) to[out=0,in=270,looseness=1.2] (1,0) node[red,xshift=-.25cm,yshift=0cm,align=center]{c\\b};
    \end{scope}
    \end{scope}
    \end{scope}
    \draw[dashed] (0,0) to[out=90,in=0] (-1.2,2) -- (-8,2);
    \draw[fill=gray] (0,0) circle (.25);
    \draw (-2.25,-2) to[out=270,in=270,looseness=1.2] (1.5,0) -- (1.5,2) to[out=90,in=90,looseness=1.2] (-2.25,4);
    \draw (-3.5,-2) to[out=270,in=270,looseness=1.2] (2,-.167) -- (2,2.167) to[out=90,in=90,looseness=1.2] (-3.5,4);
    \draw (-4.75,-2) to[out=270,in=270,looseness=1.2] (2.5,-.333) -- (2.5,2.333) to[out=90,in=90,looseness=1.2] (-4.75,4);
    \draw (-6,-2) to[out=270,in=270,looseness=1.2] (3,-.5) -- (3,2.5) to[out=90,in=90,looseness=1.2] (-6,4);
\end{tikzpicture}
\caption{The diagram obtained by performing the (expanding) shift by the split $(b,a)$ on the closed strand diagram represented in \cref{fig_CSD}.}
\label{fig_CSD_shift_ex}
\end{figure}

\phantomsection\label{TXT shifts}
A \textbf{shift} of the base line is the transformation of an $R$-branching closed strand diagram shown in \cref{fig_CSD_shift}, which essentially consists of moving the base line in such a way that it crosses exactly one split or one merge.
An example is depicted in \cref{fig_CSD_shift_ex}.
We say that a shift is \textbf{expanding} and that it \textit{expands} the split (merge) that it crosses if the number of base points is increased by the shift.
If instead the number of base points is decreased, we say that the shift is \textbf{reducing} and that it \textit{reduces} the split (merge) that it crosses.
Observe that, when a shift expands a split (merge), new symbols may need to be generated below the split (above the merge) in order to make sure that the resulting diagram is $R$-branching (\cref{def.r.branching}), as is the case with the letters v and w in \cref{fig_CSD_shift_ex}.

As we will see in \cref{PROP conjugator}, shifts of the base line represent conjugations by split and merge diagrams (which were defined in \cref{SUB groupoid generators}).

\medskip

\phantomsection\label{TXT similarity}
The two types of transformations defined so far are called \textbf{similarities}, and we say that two closed strand diagrams are \textbf{similar} if they only differ by the application of shifts and permutations.
If $\eta$ is a closed diagram, we denote by $\llbracket \eta \rrbracket$ its similarity class.
In practice, two closed diagrams are similar if they are the same up to ``forgetting'' about the base line and the base points.

\begin{remark}\label{RMK cohomology}
    The base line can be thought of as a cocycle of the graph cohomology of the graph obtained from the closed strand diagram by ``forgetting'' about the base line and the base points.
    In this sense, similarities preserve cohomology, and cohomology classes of this graph correspond to similarity classes of closed strand diagrams.
    This is going to be useful for the algorithm in \cref{SUB algorithm}.
\end{remark}

\subsubsection{Reductions of Closed Strand Diagrams}
\label{sub.closed.strand.diagrams.reductions}

A \textbf{Type 1 or 2 reduction} of an $R$-branching closed diagram is either of the two moves shown in \cref{fig_SD_reductions}.
These are the same reductions of strand diagrams, but it is important to keep in mind that these cannot ``cross'' the base line.
In practice, shifts (defined at \cpageref{TXT shifts}) are often needed to move the base line and ``unlock'' these types of reductions.

\begin{figure}\centering
\begin{tikzpicture}[scale=.65,font=\scriptsize]
    \useasboundingbox (-8.25,-13.5) rectangle (7.25,5.5);
    %\draw[help lines,step=.1cm] (-8.25,-13.5) grid (7.25,5.5);
    %
    \draw[dashed] (0,0) -- (-8,0);
    \draw[fill=gray] (0,0) circle (.25);
    \draw (-7,1) -- (-7,0) node[black,circle,fill,inner sep=1.25]{} -- node[Orange,left,align=center]{v\textsubscript{1}\\w\textsubscript{1}} (-7,-1);
    \draw[dotted] (-6,1) -- (-6,-1);
    \draw (-5,1) -- (-5,0) node[black,circle,fill,inner sep=1.25]{} -- node[Plum,left,align=center]{v\textsubscript{n}\\w\textsubscript{n}} (-5,-1);
    \draw (-3.5,1) -- (-3.5,0) node[black,circle,fill,inner sep=1.25]{} -- node[Orange,left,align=center]{v\textsubscript{1}\textsuperscript{$\prime$}\\w\textsubscript{1}\textsuperscript{$\prime$}} (-3.5,-1);
    \draw[dotted] (-2.5,1) -- (-2.5,-1);
    \draw (-1.5,1) -- (-1.5,0) node[black,circle,fill,inner sep=1.25]{} -- node[Plum,left,align=center]{v\textsubscript{n}\textsuperscript{$\prime$}\\w\textsubscript{n}\textsuperscript{$\prime$}} (-1.5,-1);
    \draw (-7,-1) to[out=270,in=270] (7,-1.5);
    \draw[dotted] (-6,-1) to[out=270,in=270] (6,-1.5);
    \draw (-5,-1) to[out=270,in=270] (5,-1.5);
    \draw (-3.5,-1) to[out=270,in=270] (3.5,-1.5);
    \draw[dotted] (-2.5,-1) to[out=270,in=270] (2.5,-1.5);
    \draw (-1.5,-1) to[out=270,in=270] (1.5,-1.5);
    \draw (1.5,-1.5) to[out=90,in=270] (5,1.5);
    \draw[dotted] (2.5,-1.5) to[out=90,in=270] (6,1.5);
    \draw (3.5,-1.5) to[out=90,in=270] (7,1.5);
    \draw (5,-1.5) to[out=90,in=270] (1.5,1.5);
    \draw[dotted] (6,-1.5) to[out=90,in=270] (2.5,1.5);
    \draw (7,-1.5) to[out=90,in=270] (3.5,1.5);
    \draw (7,1.5) to[out=90,in=90] (-7,1);
    \draw[dotted] (6,1.5) to[out=90,in=90] (-6,1);
    \draw (5,1.5) to[out=90,in=90] (-5,1);
    \draw (3.5,1.5) to[out=90,in=90] (-3.5,1);
    \draw[dotted] (2.5,1.5) to[out=90,in=90] (-2.5,1);
    \draw (1.5,1.5) to[out=90,in=90] (-1.5,1);
    \draw[thick,-stealth] (0,-6) -- (0,-7);
    \begin{scope}[yshift=-10.5cm]
    \draw[dashed] (0,0) -- (-4,0);
    \draw[fill=gray] (0,0) circle (.25);
    \draw (-3,0) -- node[Green,left,align=center]{i\\t} (-3,-1) to[out=270,in=270] (3,-1) to[out=90,in=270] (1.5,1) to[out=90,in=90] (-1.5,1) -- (-1.5,0);
    \draw (-1.5,0) node[black,circle,fill,inner sep=1.25]{} -- node[Green,left,align=center]{i\textsuperscript{$\prime$}\\t\textsuperscript{$\prime$}} (-1.5,-1) to[out=270,in=270] (1.5,-1) to[out=90,in=270] (3,1) to[out=90,in=90] (-3,1) -- (-3,0) node[black,circle,fill,inner sep=1.25]{};
    \end{scope}
\end{tikzpicture}
\caption{A schematic depiction of a Type 3 reduction with winding number equal to 2.}
\label{fig_3reduction}
\end{figure}
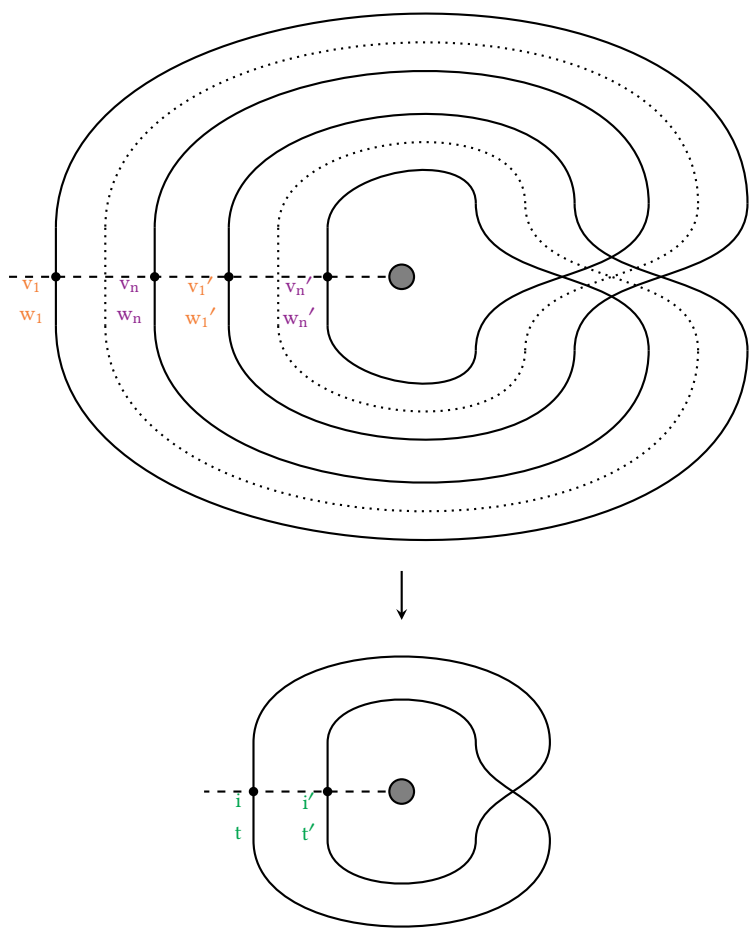

\phantomsection\label{TXT 3 reductions}
A \textbf{looping strand} is a path of strands without splits nor merges.
A \textbf{Type 3 reduction} is a move obtained by considering a sequence of looping strands that are consecutive in their given order at the intersection with the base line and whose labels correspond to the bottom labels of some replacement tree $T$ up to some renaming;
then replacing those looping strands with a single looping strand labeled as the top strand of $T$, up to the same renaming.
Each loop can wind multiple times around the central ``hole'', in which case there are multiple copies of base points that are labeled as the top strand of $T$;
the strands that are not involved in the reduction must intersect the base line on the left or on the right of those that are involved (this can be achieved by performing a permutation).
For example, \cref{fig_3reduction} depicts the general shape of a Type 3 reduction with winding number equal to 2.

Additionally, a Type 3 reduction is only allowed when the labeling symbols being ``deleted'' by the reduction do not appear among labels of strands not involved in the reduction.
The reason for this is that we would otherwise lose data about edge adjacency, since Type 3 reductions encode ``anti-expansions'' of some portion of the base graph and thus delete certain vertices.
Finally, it is important to keep in mind that the looping strands involved in a Type 3 reduction must be in the right order, which is why in practice permutations of the base points (defined at \cpageref{TXT permutations}) are often needed to ``unlock'' this kind of reductions.

\vfill %layout

\subsection{Confluent Reduction Systems and Reduced Closed Diagrams}
\label{SUB confluent reductions}

We say that two closed strand diagrams are \textbf{equivalent} if we can obtain one from the other by applying a sequence of reductions.
The similarity class of a closed strand diagram (defined at \cpageref{TXT similarity}) is \textbf{reduced} if no reduction can be performed on any diagram that it contains.
We will just say that a closed strand diagram is reduced if it belongs to a reduced similarity class, as these only differ by moving the base line or the base points.

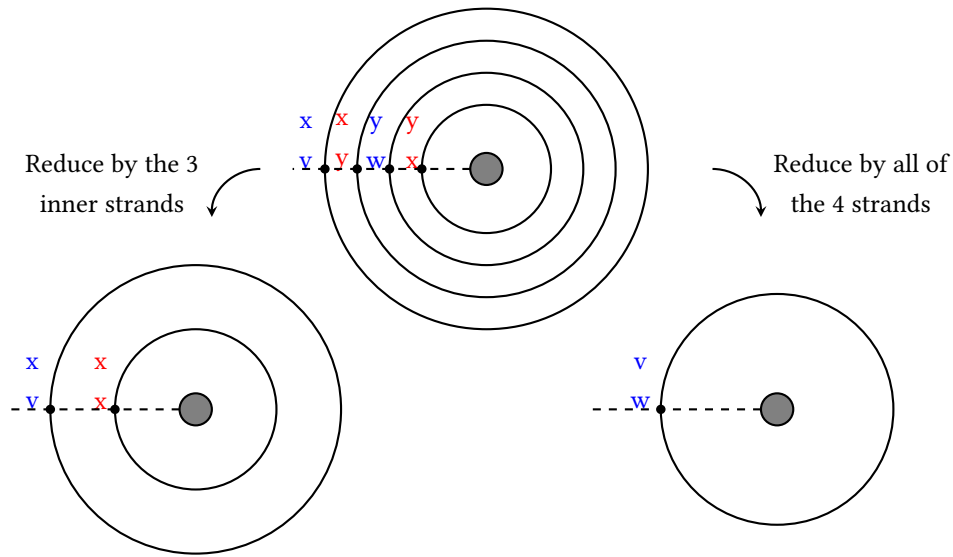
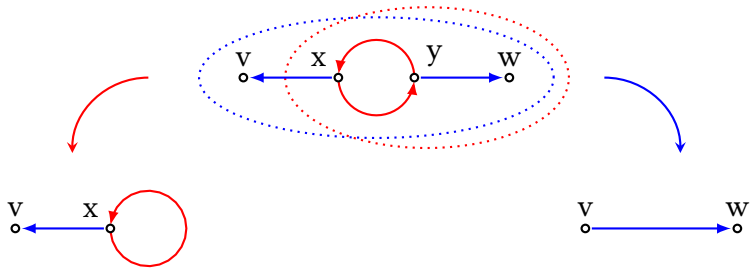
\begin{figure}\centering
\begin{subfigure}[c]{\textwidth}
\begin{tikzpicture}[font=\footnotesize,scale=.85]
    \draw (0,0) circle (1) node[red,xshift=-0.975cm,yshift=.35cm,align=center]{y\\x};
    \draw (0,0) circle (1.5) node[blue,xshift=-1.45cm,yshift=.35cm,align=center]{y\\w};
    \draw (0,0) circle (2) node[red,xshift=-1.9cm,yshift=.35cm,align=center]{x\\y};
    \draw (0,0) circle (2.5) node[blue,xshift=-2.375cm,yshift=.35cm,align=center]{x\\v};
    \draw[dashed] (0,0) -- (-1,0) node[black,circle,fill,inner sep=1.25]{} -- (-1.5,0) node[black,circle,fill,inner sep=1.25]{} -- (-2,0) node[black,circle,fill,inner sep=1.25]{} -- (-2.5,0) node[black,circle,fill,inner sep=1.25]{} -- (-3,0);
    \draw[fill=gray] (0,0) circle (.25);
    \draw[thick,-stealth] (-3.5,0) to[out=180,in=90] node[xshift=-1.5cm,align=center]{Reduce by the 3\\inner strands} (-4.25,-.75);
    \begin{scope}[xshift=-4.5cm,yshift=-3.75cm]
    \draw (0,0) circle (1.25) node[red,xshift=-1.25cm,yshift=.35cm,align=center]{x\\x};
    \draw (0,0) circle (2.25) node[blue,xshift=-2.15cm,yshift=.35cm,align=center]{x\\v};
    \draw[dashed] (0,0) -- (-1.25,0) node[black,circle,fill,inner sep=1.25]{} -- (-2.25,0) node[black,circle,fill,inner sep=1.25]{} -- (-2.9,0);
    \draw[fill=gray] (0,0) circle (.25);
    \end{scope}
    \draw[thick,-stealth] (3.5,0) to[out=0,in=90] node[xshift=1.5cm,align=center]{Reduce by all of\\the 4 strands}  (4.25,-.75);
    \begin{scope}[xshift=4.5cm,yshift=-3.75cm]
    \draw (0,0) circle (1.8) node[blue,xshift=-1.8cm,yshift=.35cm,align=center]{v\\w};
    \draw[dashed] (0,0) -- (-1.8,0) node[black,circle,fill,inner sep=1.25]{} -- (-2.9,0);
    \draw[fill=gray] (0,0) circle (.25);
    \end{scope}
\end{tikzpicture}
\caption{Two reduced closed strand diagrams obtained from the identity of $T_A$.}
\end{subfigure}
\begin{subfigure}[c]{\textwidth}
\begin{tikzpicture}
    \draw[blue,dotted] (0,0) ellipse (2.333cm and .8cm);
    \draw[red,dotted] (.667,0) ellipse (1.8667cm and .9333cm);
    \draw[edge,red,domain=5:175] plot ({.5*cos(\x)}, {.5*sin(\x)});
    \draw (90:.5);
    \draw[edge,red,domain=185:355] plot ({.5*cos(\x)}, {.5*sin(\x)});
    \draw (270:.5);
    \node[vertex] (l) at (-1.75,0) {}; \draw (-1.75,0) node[above]{v};
    \node[vertex] (cl) at (-.5,0) {}; \draw (-.5,0) node[above left]{x};
    \node[vertex] (cr) at (.5,0) {}; \draw (.5,0) node[above right]{y};
    \node[vertex] (r) at (1.75,0) {}; \draw (1.75,0) node[above]{w};
    \draw[edge,blue] (cl) to (l);
    \draw[edge,blue] (cr) to (r);
    \draw[-stealth,red] (-3,0) to[out=180,in=90] (-4,-1);
    \begin{scope}[xshift=-3cm,yshift=-2cm]
        \draw[edge,red,domain=-175:175] plot ({.5*cos(\x)}, {.5*sin(\x)});
        \node[vertex] (l) at (-1.75,0) {}; \draw (-1.75,0) node[above]{v};
        \node[vertex] (cl) at (-.5,0) {}; \draw (-.5,0) node[above left]{x};
        \draw[edge,blue] (cl) to (l);
    \end{scope}
    \draw[thick,-stealth,blue] (3,0) to[out=0,in=90] (4,-1);
    \begin{scope}[xshift=4.75cm,yshift=-2cm]
        \node[vertex] (vblue) at (-2,0) {};
        \node[vertex] (wblue) at (0,0) {};
        \draw[edge,blue] (vblue) node[above,black]{v} -- (wblue) node[above,black]{w};
    \end{scope}
    \draw[white] (0,1.25);
    \draw[white] (-6,0);
    \draw[white] (6,0);
\end{tikzpicture}
\caption{The graph reductions corresponding to the strand diagram reductions shown above.}
\end{subfigure}
\caption{The reason why the airplane edge replacement system (\cref{fig.airplane.replacement}) is not reduction-confluent.}
\label{fig_TA_reduced}
\end{figure}

\phantomsection\label{TXT airplane non confluent}
Observe that, in general, reduced similarity classes of closed strand diagrams may not be unique.
For example, \cref{fig_TA_reduced} depicts two distinct reduced closed diagrams in the airplane edge replacement system.
However, under the hypothesis of \textit{reduction-confluence} of the edge replacement system, we can immediately circumvent this issue, as is shown below.

\subsubsection{Reduction-confluent Edge Replacement Rules}
\label{SUB reduction systems}

Given a set of edge replacement rules $(R, \mathrm{C})$, we define its \textbf{reduction system} as the directed graph whose set of vertices is the set of all graphs with edges colored by $\mathrm{C}$ and whose edges are described by the anti-expansions of $(R, \mathrm{C})$.
More explicitly, given a graph $\Gamma$ and a color $c \in \mathrm{C}$, suppose that $\Gamma$ contains a subgraph $\Delta$ that is isomorphic to $R_c$ with the possible exception of having $\iota$ and $\tau$ glued together, and such that each edge of $\Gamma$ that is adjacent to some vertex of $\Delta$ except for $\iota$ and $\tau$ belongs to $\Delta$.
In this case, there is an edge ($\Gamma$, $\Gamma^*$), where $\Gamma^*$ is the graph obtained from $\Gamma$ by replacing $\Delta$ with an edge starting at $\iota$ and terminating at $\tau$ colored by $c$.

\begin{definition}\label{def.reduction.confluent}
    We say that a set of edge replacement rules $(R, \mathrm{C})$ is \textbf{reduction-confluent} if its reduction system is confluent.
    This means that, whenever $a \dashrightarrow b$ and $a \dashrightarrow c$ are two finite sequences of reductions of $a$, there exist a graph expansion $d$ (which is a vertex of the rewriting system) and two finite sequences of reductions $b \dashrightarrow d$ and $c \dashrightarrow d$.
\end{definition}

Now, if the edge replacement rules are both reduction-confluent and expanding (\cref{def.expanding}), then it is clear that the reduction system is terminating.
So, using \hyperref[lem.diamond]{Newman's Diamond Lemme}, we have that each connected component of the rewriting system has a unique reduced graph.

\phantomsection\label{TXT Basilica Vicsek BubbleBath Houghton QV}
The edge replacement systems for Thompson groups $F, T$ and $V$ (and the Higman-Thompson groups) are clearly reduction-confluent.
In truth, most of the edge replacement systems discussed in the literature (see \cref{sec.examples}) are reduction-confluent:
the dendrite, the basilica (and the rabbits), the Vicsek and the bubble bath reduction rules are confluent, as are the those for the Houghton groups $H_n$ and for the groups $QV$, $QT$, $QF$.
All of this can be easily proved using \hyperref[lem.diamond]{Newman's Diamond Lemme}.

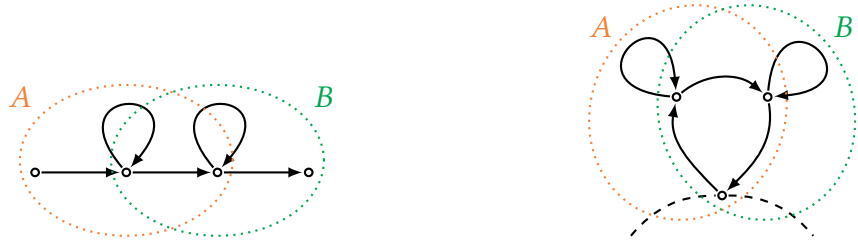
\begin{figure}\centering
\begin{subfigure}[b]{.575\textwidth}\centering
\begin{tikzpicture}[scale=.8]
    \node[vertex] (1) at (-1.5,0) {};
    \node[vertex] (2) at (0,0) {};
    \node[vertex] (3) at (1.5,0) {};
    \node[vertex] (4) at (3,0) {};
    \draw[edge] (1) to (2);
    \draw[edge] (2) to[loop,out=130,in=50,min distance=1.8cm,looseness=10] (2);
    \draw[edge] (2) to (3);
    \draw[edge] (3) to[loop,out=130,in=50,min distance=1.8cm,looseness=10] (3);
    \draw[edge] (3) to (4);
    \draw[dotted, Orange] (0,.2) ellipse (1.75cm and 1.25cm);
    \draw[Orange] (-1.75,1.25) node{$A$};
    \draw[dotted, Green] (1.5,.2) ellipse (1.75cm and 1.25cm);
    \draw[Green] (3.25,1.25) node{$B$};
\end{tikzpicture}
\end{subfigure}
\begin{subfigure}[b]{.4\textwidth}\centering
\begin{tikzpicture}[scale=.8]
    \node[vertex] (b) at (0,0) {};
    \node[vertex] (l) at (-.75,1.6333) {};
    \node[vertex] (r) at (.75,1.6333) {};
    \draw[edge] (b) to[out=135,in=260] (l);
    \draw[edge] (l) to[loop,out=175,in=95,min distance=1.8cm,looseness=10] (l);
    \draw[edge] (l) to[out=40,in=140] (r);
    \draw[edge] (r) to[loop,out=85,in=5,min distance=1.8cm,looseness=10] (r);
    \draw[edge] (r) to[out=280,in=45] (b);
    \draw[dashed] (-1.5,-.667) to[out=50,in=180] (b);
    \draw[dashed] (b) to[out=0,in=130] (1.5,-.667);
    \draw[dotted, Orange, rotate=-20] (-1,1.1) ellipse (1.6cm and 1.8cm);
    \draw[Orange] (-2,2.8) node{$A$};
    \draw[dotted, Green, rotate=20] (1,1.1) ellipse (1.6cm and 1.8cm);
    \draw[Green] (2,2.8) node{$B$};
\end{tikzpicture}
\end{subfigure}
\caption{The only two possible distinct non-disjoint graph reductions of the basilica edge replacement rules.}
\label{fig_B_reduction_confluent}
\end{figure}

For example, consider the basilica edge replacement system (depicted in \cref{fig.basilica.replacement}).
We say that two graph reductions of the same graph are disjoint if they involve subgraphs that share no edges (they can share vertices, provided they are the initial or terminal vertices of the replacement graph).
If two reductions are disjoint, then each can be applied after the other and there is nothing to prove.
Let $A$ and $B$ be two distinct graph reductions of the same graph that are not disjoint.
Then, since each reduction must identify a copy of the basilica edge replacement graph, $A$ and $B$ cannot involve the same loop, otherwise they would be the same reduction.
Thus there are two cases, both portrayed in \cref{fig_B_reduction_confluent}.
In each case, both reductions produce the same graph, so we are done.

The airplane edge replacement system instead is not reduction-confluent, as the base graph itself can be reduced in three ways that produce two distinct reduced graphs (these correspond to the ones shown in \cref{fig_TA_reduced}).

\begin{remark}
Among the edge replacement systems discussed here, the only one with more than one color happens to be exactly the one that is not reduction-confluent.
This is just a coincidence, as it is not hard to build edge replacement systems with one color that are not reduction-confluent and edge replacement systems that are reduction-confluent despite having more than one color.
\end{remark}

\subsubsection{Uniqueness of Reduced Closed Diagrams}

The following Lemma tells us that, if we are dealing with an edge replacement system that is based on reduction-confluent edge replacement rules, then each closed strand diagram is equivalent to a unique reduced closed strand diagram up to similarity.

\medskip %layout
\begin{lemma}
\label{LEM reduced CSD}
Suppose that the edge replacement rules $(R, \mathrm{C})$ are reduction-confluent.
Then, for each similarity class $\llbracket \eta \rrbracket$ of $R$-branching closed diagrams, there is a unique reduced similarity class of $R$-branching closed diagrams that is equivalent to $\llbracket \eta \rrbracket$.
\end{lemma}

\begin{proof}
As we did earlier for \cref{LEM reduced SD}, we will use \hyperref[lem.diamond]{Newman's Diamond Lemma}.
Consider the directed graph defined as follows: the set of vertices is the set of similarity classes of $R$-branching closed diagrams and we have an edge $\llbracket \eta \rrbracket \longrightarrow \llbracket \zeta \rrbracket$ for each distinct reduction from a closed diagram in the same similarity class as $\eta$ to one in the same similarity class as $\zeta$.
It then suffices to prove that this graph is terminating and locally confluent.

It is clear that similar closed diagrams have the same number of splits and merges.
Since each reduction strictly decreases the number of splits and merges of a closed diagram, the directed graph is terminating, so we only need to prove the local confluence.

Suppose $\eta, \zeta$ and $\kappa$ are closed strand diagrams such that $\llbracket \eta \rrbracket \overset{A}{\longrightarrow} \llbracket \zeta \rrbracket$ and $\llbracket \eta \rrbracket \overset{B}{\longrightarrow} \llbracket \kappa \rrbracket$ are distinct reductions.
Without loss of generality we may assume that the reductions of closed diagrams are $\eta_A \overset{A}{\longrightarrow} \zeta$ and $\eta_B \overset{B}{\longrightarrow} \kappa$, where $\eta_A$ and $\eta_B$ are similar.
This means that the reductions $A$ and $B$ are performed without shifting the base line nor permutating the order of the base points, whereas we can transform $\eta_A$ into $\eta_B$ with a finite sequence of permutations and shifts.
Note that then $\eta_A$ and $\eta_B$ have the same strands, splits and merges, and the only possible distinctions must involve the position of the base line and the order of the base points.

First, note that if none of the reductions $A$ and $B$ are of Type 3 then, with two sole exceptions, they must be disjoint, and if they are disjoint clearly the order in which $A$ and $B$ are performed does not affect the resulting diagram (as in the proof of \cref{LEM reduced SD}).
The first exception is that portrayed in \cref{fig_strand_diamond}:
as already discussed in \cref{LEM reduced SD}, the order in which one performs the two reductions does not matter.
The other exception is represented in \cref{fig_onion}:
in this case, after a Type 2 reduction on $\llbracket \eta \rrbracket$ we can apply a Type 3 reduction to get the same diagram obtained by performing a Type 1 reduction on $\llbracket \eta \rrbracket$, so local confluence is verified in this case.
Therefore, we can suppose that $B$ is of Type 3.
But then, since Type 3 reductions can only involve strands devoid of splits and merges, it is clear that if $A$ is of Type 1 or 2 then $A$ and $B$ are disjoint (by which we mean that they involve different strands), so we are done in this case.

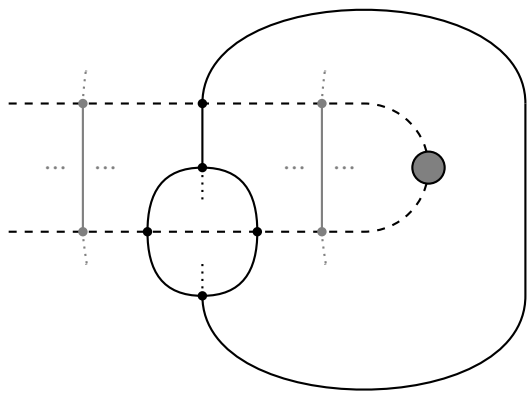
\begin{figure}\centering
\begin{tikzpicture}[font=\small,scale=.85]
    \useasboundingbox (-3,-4.5) rectangle (5,1.5);
    %\draw[help lines,step=.1cm] (-3,-4.5) grid (5,1.5);
    %
    \draw (0,-1) -- (0,0) node[black,circle,fill,inner sep=1.25]{} to[out=90,in=90] (5,0);
    \draw (0,-1) node[black,circle,fill,inner sep=1.25]{} to[out=180,in=90,looseness=1.2] (-.85,-2);
    \draw[dotted] (0,-1) -- (0,-1.5);
    \draw (0,-1) to[out=0,in=90,looseness=1.2] (.85,-2) node[black,circle,fill,inner sep=1.25]{};
    \draw (0,-3) node[black,circle,fill,inner sep=1.25]{} to[out=180,in=270,looseness=1.2] (-.85,-2) node[black,circle,fill,inner sep=1.25]{};
    \draw[dotted] (0,-3) -- (0,-2.5);
    \draw (0,-3) to[out=0,in=270,looseness=1.2] (.85,-2);
    \draw (0,-3) to[out=270,in=270] (5,-3) -- (5,0);
    \begin{scope}[xshift=3.5cm,yshift=-1cm]
    \draw[dashed] (0,0) to[out=90,in=0] (-1,1) -- (-6.5,1);
    \draw[dashed] (0,0) to[out=270,in=0] (-1,-1) -- (-6.5,-1);
    \draw[fill=gray] (0,0) circle (.25);
    \end{scope}
    \draw[gray,dotted] (-1.85,-2) to[out=270,in=285] (-1.8,-2.5);
    \draw[gray] (-1.85,-2) node[circle,fill,inner sep=1.25]{} -- node[left]{$\dots$} node[right]{$\dots$} (-1.85,0) node[circle,fill,inner sep=1.25]{};
    \draw[gray,dotted] (-1.85,0) to[out=90,in=85] (-1.8,.5);
    \draw[gray,dotted] (1.85,-2) to[out=270,in=285] (1.9,-2.5);
    \draw[gray] (1.85,-2) node[circle,fill,inner sep=1.25]{} -- node[left]{$\dots$} node[right]{$\dots$} (1.85,0) node[circle,fill,inner sep=1.25]{};
    \draw[gray,dotted] (1.85,0) to[out=90,in=85] (1.9,.5);
\end{tikzpicture}
\caption{Non-disjoint Type 1 and 2 reductions from the proof of \cref{LEM reduced CSD}. Labels are omitted for the sake of clarity, but they should allow a Type 1 reduction and a Type 2 reduction, and the two base lines represented differ by a shift that ``unlocks'' these reductions.}
\label{fig_onion}
\end{figure}

Finally, suppose that $A$ and $B$ are both of Type 3.
Consider the base graph $\Gamma$ of the diagram (which is the graph represented by the base points, as seen in \cref{SUB closed strand diagrams}).
Observe that Type 3 reductions (defined at \cpageref{TXT 3 reductions}) correspond to reductions of the base graph $\Gamma$ of the closed strand diagram, and it is clear that distinct Type 3 reductions describe distinct reductions of the base graph.
Now, if the reductions $A$ and $B$ are disjoint, by which we mean that they involve different strands of $\llbracket \eta \rrbracket$, then $B$ can be applied to $\zeta$ and $A$ to $\kappa$, so we are done.
If instead they are not disjoint, then consider the subgraph $\Delta$ of $\Gamma$ that is involved in both reductions.
More precisely, $\Delta$ is the union of the two subgraphs isomorphic to replacement graphs which determine the reductions induced by $A$ and $B$ on $\Gamma$.
Denote by $\Delta^A$ and $\Delta^B$ the graphs obtained from $\Delta$ by applying these reductions.
Now, because the edge replacement rules $(R, \mathrm{C})$ are reduction-confluent, there is a graph $\Delta^*$ that is a common reduction of both $\Delta^A$ and $\Delta^B$.
The graph reductions that one needs to perform to obtain $\Delta^*$ from $\Delta^A$ and $\Delta^B$ can clearly be performed on the copy of $\Delta$ contained as a subgraph in $\Gamma$, and these correspond to Type 3 reductions of $\llbracket \eta \rrbracket$ that make these reductions locally confluent, so we are done.
\end{proof}

It is worth noting that the reduction-confluent hypothesis is only needed to prove confluence of Type 3 reductions.
Type 1 and 2 reductions by themselves are always confluent.

\subsection{Stable and Vanishing Symbols}
\label{SUB stable and vanishing}

In this subsection we investigate what happens when a sequence of shifts moves the base line back to where it started.
This will be important when developing the algorithm for solving the conjugacy problem in \cref{SUB algorithm}

\begin{figure}\centering
\begin{tikzpicture}[font=\footnotesize,xscale=1.2]
    \begin{scope}[yshift=3cm]
    \draw[dotted] (0,.333) -- (0,0);
    \draw[dotted] (1,.333) -- (1,0);
    \draw[dotted] (2,.333) -- (2,0);
    \draw[dotted] (3,.333) -- (3,0);
    \draw[dashed,gray] (-1,0) -- (4,0);
    \draw (0,0) -- node[left,align=center]{c\\d} (0,-1) to[out=270,in=90] (1,-2) -- node[left,align=center]{c\\d} (1,-3);
    \draw (1,0) -- node[left,align=center]{a\\b} (1,-1) to[out=270,in=90] (0,-2) -- node[left,align=center]{a\\b} (0,-3);
    \draw (2,0) -- node[left,align=center]{e\\x\textsubscript{-1}} (2,-1) -- (2,-2) -- node[left,align=center]{e\\x\textsubscript{0}} (2,-3);
    \draw (3,0) -- node[left,align=center]{y\textsubscript{-1}\\f} (3,-1) -- (3,-2) -- node[left,align=center]{y\textsubscript{0}\\f} (3,-3);
    \draw (2,-1) to[out=0,in=90] (2.5,-1.5) node[left,align=center,xshift=1.25pt]{y\textsubscript{0}\\x\textsubscript{-1}} to[out=270,in=180] (3,-2);
    \end{scope}
    \draw[dashed,Green] (-1,0) -- (4,0);
    \draw (0,0) -- node[left,align=center]{a\\b} (0,-1) to[out=270,in=90] (1,-2) -- node[left,align=center]{a\\b} (1,-3);
    \draw (1,0) -- node[left,align=center]{c\\d} (1,-1) to[out=270,in=90] (0,-2) -- node[left,align=center]{c\\d} (0,-3);
    \draw (2,0) -- node[left,align=center]{e\\x\textsubscript{0}} (2,-1) -- (2,-2) -- node[left,align=center]{e\\x\textsubscript{1}} (2,-3);
    \draw (3,0) -- node[left,align=center]{y\textsubscript{0}\\f} (3,-1) -- (3,-2) -- node[left,align=center]{y\textsubscript{1}\\f} (3,-3);
    \draw (2,-1) to[out=0,in=90] (2.5,-1.5) node[left,align=center,xshift=1.25pt]{y\textsubscript{1}\\x\textsubscript{0}} to[out=270,in=180] (3,-2);
    \draw[dashed,Green] (-1,-3) -- (4,-3);
    \draw (-2.5,-1.5) node{Original copy of $Z \:$ \scalebox{1.5}[2.5]{$\Biggl\{$}};
    \begin{scope}[yshift=-3cm]
    \draw (0,0) -- node[left,align=center]{c\\d} (0,-1) to[out=270,in=90] (1,-2) -- node[left,align=center]{c\\d} (1,-3);
    \draw (1,0) -- node[left,align=center]{a\\b} (1,-1) to[out=270,in=90] (0,-2) -- node[left,align=center]{a\\b} (0,-3);
    \draw (2,0) -- node[left,align=center]{e\\x\textsubscript{1}} (2,-1) -- (2,-2) -- node[left,align=center]{e\\x\textsubscript{2}} (2,-3);
    \draw (3,0) -- node[left,align=center]{y\textsubscript{1}\\f} (3,-1) -- (3,-2) -- node[left,align=center]{y\textsubscript{2}\\f} (3,-3);
    \draw (2,-1) to[out=0,in=90] (2.5,-1.5) node[left,align=center,xshift=1.25pt]{y\textsubscript{2}\\x\textsubscript{1}} to[out=270,in=180] (3,-2);
    \draw[dashed,gray] (-1,-3) -- (4,-3);
    \draw[dotted] (0,-3.333) -- (0,-3);
    \draw[dotted] (1,-3.333) -- (1,-3);
    \draw[dotted] (2,-3.333) -- (2,-3);
    \draw[dotted] (3,-3.333) -- (3,-3);
    \end{scope}
\end{tikzpicture}
\caption{An example of $Z^\infty$ for an element $Z$ based on the edge replacement rules for Thompson's group $V$, depicted in \cref{fig.Cantor.replacement}.}
\label{fig_infinite_power}
\end{figure}
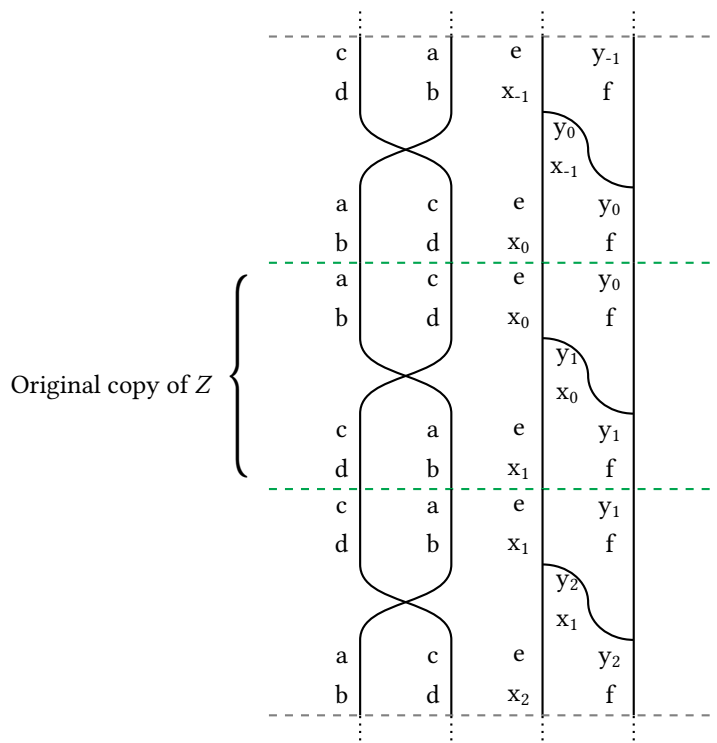

Let $Z$ be an $\mathcal{X}$-strand diagram for some edge replacement system $\mathcal{X}$ (as defined at \cpageref{TXT X-SDs}).
Denote by $Z^\infty$ the infinite power of $Z$, which is the infinite strand diagram obtained by inductively taking higher and higher powers of $Z$ and keeping track of the lines under which the copies of $Z$ are attached.
More precisely, in order to build $Z^\infty$, start from $Z$ with two ``base lines'', one at the top and one at the bottom;
these will be called \textbf{main base lines}.
Start attaching copies of $Z$ at both sides, drawing new base lines where the elements are glued.
Recall that labels need to be adjusted when gluing together the copies of $Z$, as described in \cref{sub.SDs.composition};
keep unchanged the labels of the original copy of $Z$, and adjust the labels above and below the main base lines.
The result of this infinite procedure is $Z^\infty$.
An example is portrayed in \cref{fig_infinite_power}, where the main base lines are depicted in green.

\begin{definition}
\label{def.stable.vanishing}
    Each symbol that appears in $Z$ is either:
    \begin{itemize}
        \item \textbf{stable} if it appears infinitely many times in $Z^\infty$,
        \item \textbf{vanishing} if it only appears finitely many times in $Z^\infty$.
    \end{itemize}
    In other words, stable symbols appear throughout all of $Z^\infty$, whereas the instances of a vanishing symbol are limited to a finite portion of $Z^\infty$ (i.e., a composition of finitely many copies of $Z$ inside $Z^\infty$).
    Thus, all instances of a vanishing symbol are contained in some minimal portion $Z^I \subset Z^\infty$, where $I$ is an interval of integers, and outside of $Z^I$ such a symbol is replaced by other vanishing symbols.
\end{definition}

In the example portrayed in \cref{fig_infinite_power}, the symbols a, b, c, d, e and f are stable, while the symbols x\textsubscript{$i$} and y\textsubscript{$i$} ($i \in \mathbb{Z}$) are vanishing.

For the curious reader, keeping in mind that each symbol represents a vertex in a graph expansion, one can see that each stable symbol corresponds to a point of the limit space that is periodic under the action of the rearrangement represented by $Z$.

\begin{remark}
\label{RMK vanishing symbols}
    Because of point (3) of \cref{def.r.branching}, if a symbol is vanishing, then it can only appear in one connected component of the closed diagram $\zeta$ obtained from $Z$.
\end{remark}

\begin{remark}
\label{RMK stable symbols}
    Observe that there is a finite amount of distinct stable symbols, even if (by definition) there are infinitely many instances in which each stable symbol appears.
    Indeed, since each copy of $Z$ has finitely many symbols, the trajectory of a stable symbol is periodic (i.e., navigating through $Z^\infty$, a stable symbol will eventually be found twice in the same position).
\end{remark}

The next Proposition tells us that moving the base line does not change the labeling, up to a permutation of the stable symbols (for example, in \cref{fig_infinite_power} the pairs of symbols (a,b) and (c,d) may be switched), and the results of this permutation can be listed by an algorithm, as stated in \cref{cor.algorithm.labels} below.

\medskip %layout
\begin{proposition}
    Let $\zeta$ be the closed diagram obtained from $Z$.
    If a sequence of shifts moves the base line of $\zeta$ back to the original position, then the resulting labeling of the diagram remains the same up to renaming symbols.
    This renaming induces a permutation of the stable symbols in each connected component of $\zeta$.
\end{proposition}

\begin{proof}
    Observe that each shift of the base line of $\zeta$ corresponds to applying the same shift to the base lines of $Z^\infty$.
    Thus, moving the base line of $\zeta$ to the original position is equivalent to moving the ``window'' between the two main base lines to some other portion of $Z^\infty$ that is the same up to renaming symbols.
    More precisely, this results in moving every base line of $Z^\infty$ up or down by $n_i$ steps (where $n_i \in \mathbb{Z}$) through each connected component $Z_i$ of $Z^\infty$.
    Now, as noted in \cref{RMK vanishing symbols}, vanishing symbols are distinct in each connected component, so one can always rename them to get them back to the original configuration.
    Finally, by \cref{RMK stable symbols} there is a finite amount of distinct stable symbols, so they must be permuted independently in each connected component.
\end{proof}

Since there are only finitely many distinct stable symbols, the permutation associated by cycling once through a connected component of a closed strand diagram has finite order.
Thus, if one cycles through the base line in this way multiple times, at some point one will find the original configuration of labels of the connected component (up to renaming the vanishing symbols).
Going through this process for every connected component yields an algorithm that, by exhaustion, lists every possible such configuration of the stable symbols, so we have the following.

\medskip %layout
\begin{corollary}
\label{cor.algorithm.labels}
There is an algorithm that lists the finitely many possible configurations of labels of a closed diagram $\zeta$ (up to renaming) obtained by shifting the base line from a fixed position back to itself.
\end{corollary}

\section{Solving the Conjugacy Problem}
\label{SEC conjugacy problem}

In this last section we examine conjugacy in the case of reduction-confluent edge replacement rules (\cref{SUB conf}) and we provide a pseudo-algorithm to solve the conjugacy problem under this assumption (\cref{SUB algorithm}).
Then we discuss the problem without assumption of reduction-confluence, solving the conjugacy problem in the airplane rearrangement group (\cref{sub.non.confluent}).
Finally, in \cref{sub.DPO} we link this problematic to the general setting of DPO graph rewriting systems studied in computer science.

\subsection{Reduction-confluent Edge Replacement Rules}
\label{SUB conf}

Let $\mathcal{X} = (X_0, R, \mathrm{C})$ be a fixed an edge replacement system.

\medskip %layout
\begin{proposition}
Suppose that $(R, \mathrm{C})$ are reduction-confluent edge replacement rules.
Then, if two rearrangements $f$ and $g$ of the same edge replacement system $\mathcal{X}$ are conjugate in $G_\mathcal{X}$, their reduced closed diagrams are similar.
\end{proposition}

\begin{figure}\centering
\begin{tikzpicture}[scale=.4]
    \useasboundingbox (-6,-26) rectangle (12.25,4.5);
    %\draw[help lines,step=.5cm] (-6,-26) grid (12.25,4.5);
    %
    \begin{scope}[draw=Orange]
    \draw (-1.75,-.5) node{\textcolor{Orange}{$\cdots$}};
    \draw (1.75,-.5) node{\textcolor{Orange}{$\cdots$}};
    \draw (-3.5,0) -- (-3.5,-2);
    \draw (-3.5,-5) -- (-3.5,-7);
    \draw[dotted] (-4,-2.5) to[out=90,in=180] (-3.5,-2) to[out=0,in=90] (-3,-2.5);
    \draw[dotted] (-4,-4.5) to[out=270,in=180] (-3.5,-5) to[out=0,in=270] (-3,-4.5);
    \draw (3.5,0) -- (3.5,-2);
    \draw (3.5,-5) -- (3.5,-7);
    \draw[dotted] (4,-2.5) to[out=90,in=0] (3.5,-2) to[out=180,in=90] (3,-2.5);
    \draw[dotted] (4,-4.5) to[out=270,in=0] (3.5,-5) to[out=180,in=270] (3,-4.5);
    \draw (0,0) -- (0,-1);
    \draw[dotted] (0,-1) -- (0,-2);
    \draw (0,-1) to[out=180,in=90] (-1.5,-2);
    \draw (0,-1) to[out=0,in=90] (1.5,-2);
    \draw[dotted] (-2,-2.5) to[out=90,in=180] (-1.5,-2) to[out=0,in=90] (-1,-2.5);
    \draw[dotted] (2,-2.5) to[out=90,in=0] (1.5,-2) to[out=180,in=90] (1,-2.5);
    \draw (-5.3,-3.5) node{\textcolor{Orange}{\Large$h^{-1}$}};
    \draw (0,-3.5) node{\textcolor{Orange}{\LARGE$\cdots$}};
    \draw (0,-7) -- (0,-6);
    \draw[dotted] (0,-6) -- (0,-5);
    \draw (0,-6) to[out=0,in=270] (1.5,-5);
    \draw (0,-6) to[out=180,in=270] (-1.5,-5);
    \draw[dotted] (-2,-4.5) to[out=270,in=180] (-1.5,-5) to[out=0,in=270] (-1,-4.5);
    \draw[dotted] (2,-4.5) to[out=270,in=0] (1.5,-5) to[out=180,in=270] (1,-4.5);
    \end{scope}
    \begin{scope}[draw=Plum,yshift=-7cm]
    \draw (-3.5,0) -- (-3.5,-2);
    \draw (-3.5,-5) -- (-3.5,-7);
    \draw[dotted] (-4,-2.5) to[out=90,in=180] (-3.5,-2) to[out=0,in=90] (-3,-2.5);
    \draw[dotted] (-4,-4.5) to[out=270,in=180] (-3.5,-5) to[out=0,in=270] (-3,-4.5);
    \draw (3.5,0) -- (3.5,-2);
    \draw (3.5,-5) -- (3.5,-7);
    \draw[dotted] (4,-2.5) to[out=90,in=0] (3.5,-2) to[out=180,in=90] (3,-2.5);
    \draw[dotted] (4,-4.5) to[out=270,in=0] (3.5,-5) to[out=180,in=270] (3,-4.5);
    \draw (0,0) -- (0,-1);
    \draw[dotted] (0,-1) -- (0,-2);
    \draw (0,-1) to[out=180,in=90] (-1.5,-2);
    \draw (0,-1) to[out=0,in=90] (1.5,-2);
    \draw[dotted] (-2,-2.5) to[out=90,in=180] (-1.5,-2) to[out=0,in=90] (-1,-2.5);
    \draw[dotted] (2,-2.5) to[out=90,in=0] (1.5,-2) to[out=180,in=90] (1,-2.5);
    \draw (-5.3,-3.5) node{\textcolor{Plum}{\Large$g$}};
    \draw (0,-3.5) node{\textcolor{Plum}{\LARGE$\cdots$}};
    \draw (0,-7) -- (0,-6);
    \draw[dotted] (0,-6) -- (0,-5);
    \draw (0,-6) to[out=0,in=270] (1.5,-5);
    \draw (0,-6) to[out=180,in=270] (-1.5,-5);
    \draw[dotted] (-2,-4.5) to[out=270,in=180] (-1.5,-5) to[out=0,in=270] (-1,-4.5);
    \draw[dotted] (2,-4.5) to[out=270,in=0] (1.5,-5) to[out=180,in=270] (1,-4.5);
    \end{scope}
    \begin{scope}[draw=Orange,yshift=-14cm]
    \draw (-3.5,0) -- (-3.5,-2);
    \draw (-3.5,-5) -- (-3.5,-7);
    \draw[dotted] (-4,-2.5) to[out=90,in=180] (-3.5,-2) to[out=0,in=90] (-3,-2.5);
    \draw[dotted] (-4,-4.5) to[out=270,in=180] (-3.5,-5) to[out=0,in=270] (-3,-4.5);
    \draw (3.5,0) -- (3.5,-2);
    \draw (3.5,-5) -- (3.5,-7);
    \draw[dotted] (4,-2.5) to[out=90,in=0] (3.5,-2) to[out=180,in=90] (3,-2.5);
    \draw[dotted] (4,-4.5) to[out=270,in=0] (3.5,-5) to[out=180,in=270] (3,-4.5);
    \draw (0,0) -- (0,-1);
    \draw[dotted] (0,-1) -- (0,-2);
    \draw (0,-1) to[out=180,in=90] (-1.5,-2);
    \draw (0,-1) to[out=0,in=90] (1.5,-2);
    \draw[dotted] (-2,-2.5) to[out=90,in=180] (-1.5,-2) to[out=0,in=90] (-1,-2.5);
    \draw[dotted] (2,-2.5) to[out=90,in=0] (1.5,-2) to[out=180,in=90] (1,-2.5);
    \draw (-5.3,-3.5) node{\textcolor{Orange}{\Large$h$}};
    \draw (0,-3.5) node{\textcolor{Orange}{\LARGE$\cdots$}};
    \draw (0,-7) -- (0,-6);
    \draw[dotted] (0,-6) -- (0,-5);
    \draw (0,-6) to[out=0,in=270] (1.5,-5);
    \draw (0,-6) to[out=180,in=270] (-1.5,-5);
    \draw[dotted] (-2,-4.5) to[out=270,in=180] (-1.5,-5) to[out=0,in=270] (-1,-4.5);
    \draw[dotted] (2,-4.5) to[out=270,in=0] (1.5,-5) to[out=180,in=270] (1,-4.5);
    \draw (-1.75,-6.5) node{\textcolor{Orange}{$\cdots$}};
    \draw (1.75,-6.5) node{\textcolor{Orange}{$\cdots$}};
    \end{scope}
    \draw (-3.5,-21) to[out=270,in=270,looseness=1] (12.25,-22) -- (12.25,0) to[out=90,in=90,looseness=1] (-3.5,0);
    \draw (0,-21) to[out=270,in=270,looseness=1] (10.5,-21.5) -- (10.5,0) to[out=90,in=90,looseness=1] (0,0);
    \draw (3.5,-21) to[out=270,in=270,looseness=1] (8.75,-21) -- (8.75,0) to[out=90,in=90,looseness=1] (3.5,0);
    \draw[dashed] (7,-10.5) to[out=90,in=0,looseness=.75] (3.5,0) -- (-5,0);
    \draw[dashed,ForestGreen] (7,-10.5) to[out=90,in=0,looseness=1] (3.5,-7) -- (-5,-7);
    \draw[fill=gray] (7,-10.5) circle (.25);
\end{tikzpicture}
\caption{A closed strand diagram obtained from $f = g^h$ and, in \textcolor{ForestGreen}{green}, a different position of the base line.}
\label{fig_conj}
\end{figure}
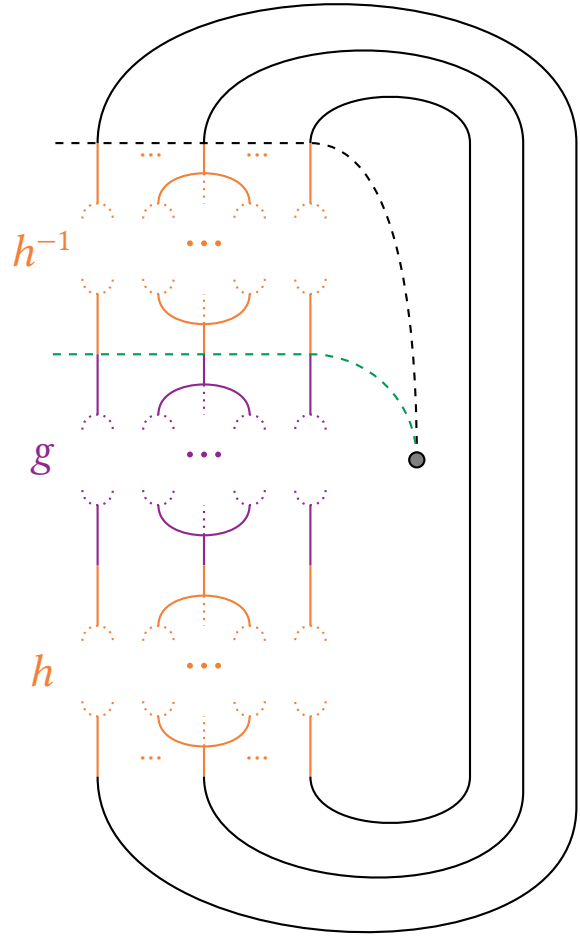

\begin{proof}
Suppose that $f = g^h$, for $f$, $g$ and $h$ in $G_\mathcal{X}$.
Then you can represent the strand diagram for $f$ as shown in \cref{fig_conj}.
This is usually not a reduced strand diagram, but we do not need it to be.
If we now close the diagram, we can perform shifts to move the base line counterclockwise until it is right between $h$ and $g$, as represented in green in \cref{fig_conj}.
In this configuration, $h$ and $h^{-1}$ cancel out completely, leaving the same closed strand diagram that one would obtain when closing $g$.
Thus the closures of $f$ and $g$ are equivalent to a common closed strand diagram, and so because of \cref{LEM reduced CSD} they must have the same reduced closed strand diagram, up to similarity.
\end{proof}

The next result is the bridge that links reductions to conjugacy.
Note that it does not require the hypothesis of reduction-confluence. 

\medskip %layout
\begin{proposition}
\label{PROP conjugator}
Given $f, g \in G_\mathcal{X}$, if up to similarity the closures of any of their strand diagrams can be rewritten as a common reduced closed diagram using reductions and inverses of reductions, then $f$ and $g$ are conjugate in $G_\mathcal{X}$.
\end{proposition}

\begin{proof}
First, observe that it does not matter which strand diagrams we choose to represent $f$ and $g$.
Indeed, two distinct strand diagrams for the same elements would only differ by a sequence of Type 1 and Type 2 reductions, which can also be performed on the closure of such diagrams.

We now note that, if two closed strand diagrams $\eta$ and $\zeta$ differ by a transformation, then the respective ``open'' strand diagrams $o(\eta)$ and $o(\zeta)$ are conjugate by a specific element of the replacement groupoid which depends on the transformation, as explained in the points below.

\begin{itemize}
\phantomsection\label{TXT transformations to conjugacy}
    \item If the transformation is a shift, then the open strand diagrams $o(\eta)$ and $o(\zeta)$ differ from conjugating by a split or a merge diagram.
    \item If the transformation is a permutation, then the open strand diagrams $o(\eta)$ and $o(\zeta)$ differ from conjugating by a permutation diagram.
    \item If the transformation is a Type 1 or 2 reduction, then the open strand diagrams $o(\eta)$ and $o(\zeta)$ represent the same rearrangement.
    \item If the transformation is a Type 3 reduction, then the open strand diagrams $o(\eta)$ and $o(\zeta)$ differ from conjugating by an amount of split diagrams equal to the winding number of the loops involved in the reduction.
\end{itemize}

Suppose that two closed diagrams $\eta$ and $\zeta$ obtained from $f$ and $g$ have the same reduced closed diagram $\rho$.
Then there exist sequences of transformations from $\eta$ to $\rho$ and from $\zeta$ to $\rho$.
Since each transformation corresponds to a (possibly trivial) conjugation of the open diagrams, it follows that both $o(\eta)^{h_\eta} = o(\rho)$ and $o(\zeta)^{h_\eta} = o(\rho)$, where $h_\eta$  and $h_\zeta$ are compositions of split diagrams, merge diagrams and permutation diagrams, which are elements of the replacement groupoid.
Thus, $o(\eta) = o(\zeta)^{h_\zeta \circ h_\eta^{-1}}$.

Finally, we need to show that $h_\zeta \circ h_\eta^{-1}$ is a strand diagram that represents an element of $G_\mathcal{X}$.
The composition $o(\zeta)^{h_\zeta \circ h_\eta}$ exists and results in $o(\eta)$, so $h_\zeta \circ h_\eta^{-1}$ is an element of the replacement groupoid whose source and sink graphs represent the base graph of $\mathcal{X}$, i.e., $h_\zeta \circ h_\eta^{-1}$ is an $\mathcal{X}$-strand diagram (as defined at \cpageref{TXT X-SDs}).
As discussed at \cpageref{TXT X-SDs}, these are exactly the elements of $\mathcal{X}$, so we are done.
\end{proof}

Together these two Propositions yield the following result.

\medskip %layout
\begin{proposition}
\label{PROP conjugate}
Suppose that $\mathcal{X} = (X_0, R, \mathrm{C})$ is an expanding edge replacement system whose set of edge replacement rules $(R, \mathrm{C})$ is reduction-confluent.
Two elements $f, g \in G_\mathcal{X}$ are conjugate in $G_\mathcal{X}$ if and only if their reduced closed diagrams are similar.
\end{proposition}

\subsection{The Algorithm}
\label{SUB algorithm}

Using the results from the last subsection, under the assumption of reduction-confluence we can algorithmically check whether two rearrangements $f$ and $g$ are conjugate in the rearrangement group:
it suffices to consider the closed strand diagrams $\eta$ and $\zeta$ for $f$ and $g$, respectively, compute their reduced diagrams $\eta^*$ and $\zeta^*$ and finally check whether they are the same up to similarity.
In more detail, the algorithm is the following.

\begin{enumerate}
    \setcounter{enumi}{-1}
    \item Consider the closed strand diagrams $\eta$ and $\zeta$ obtained from $f$ and $g$.
    \item Compute their reduced closed diagrams $\eta^*$ and $\zeta^*$ in the following way.
    \begin{itemize}
        \item Check whether the diagram contains a Type 1 or 2 reduction up to permutations and shifts (which can be done by ``forgetting'' about the base line and the base points). If it does, perform the reduction, possibly performing a permutation of the base points or moving the base line with a shift if necessary to ``unlock'' the reduction. Continue doing this until no Type 1 and 2 reductions are possible. This procedure must eventually stop, since Type 1 and 2 reductions decrease the amount of strands.
        \item Check whether the diagram contains a Type 3 reduction, up to permutations. If it does, perform the reduction, possibly performing a permutation of the base points if necessary to ``unlock'' the reduction.
    \end{itemize}
    \item In order to check similarity, consider the graphs obtained by ``forgetting'' about the labels, the base line and the base points of the two diagrams and proceed as follows.
    \begin{itemize}
        \item Check whether the graphs are isomorphic (this terminates because the graphs are finite). If they are not isomorphic, then the strand diagrams cannot be similar. Otherwise, execute the next steps multiple times, separately identifying the two graphs under each isomorphism.
        \item Check whether the two positions of the base lines are the same up to similarity.
        This is equivalent to asking whether they describe cohomologous cocycles by \cref{RMK cohomology}, so it can be efficiently computed.
        Explicit examples of these computations for Thompson's group $V$ are contained in the dissertation \cite{Algorithm}, whose content can be adapted to suit the setting of rearrangement groups. If the answer is negative, then the diagrams are not similar.
        \item Once established that the positions of the base lines are the same up to similarity, it only remains to check that labels are the same up to renaming. To do so, find the same position of the base line by computing, step by step and in parallel, all possible similarities of the second diagram $\zeta^*$ until the position of the base line matches the one in $\eta^*$ (we know that this process terminates because in the previous step we have established that the positions of the base lines are the same up to similarity).
        Once the same position of the base lines has been found, check whether labels are the same up to renaming and up to the finitely many permutations of stable symbols (see \cref{cor.algorithm.labels}).
    \end{itemize}
\end{enumerate}

Note that the order in which the reductions are performed in step (1) does not matter, as there is a unique reduced similarity class thanks to \cref{LEM reduced CSD}

Then we have the following:

\medskip %layout
\begin{theorem}
\label{thm.conjugacy.solvable}
Suppose that $\mathcal{X} = (X_0, R, \mathrm{C})$ is an expanding edge replacement system whose set of edge replacement rules $(R, \mathrm{C})$ is reduction-confluent.
Then the conjugacy problem of $G_\mathcal{X}$ is solvable.
\end{theorem}

Using this method, not only we can say whether two given rearrangements are conjugate, but when they are we can also explicitly compute a conjugating element, as seen in the proof of \cref{PROP conjugator}.
Thus, we also have:

\medskip %layout
\begin{corollary}
Suppose that $\mathcal{X} = (X_0, R, \mathrm{C})$ is an expanding edge replacement system whose set of edge replacement rules $(R, \mathrm{C})$ is reduction-confluent.
Then the conjugacy search problem of $G_\mathcal{X}$ is solvable.
\end{corollary}

As a consequence of \cref{thm.conjugacy.solvable}, since the Higman-Thompson groups, the Houghton groups, the groups $QV$, $QT$ and $QF$, the basilica, the rabbit, the Vicsek and the bubble bath edge replacement rules are reduction-confluent (as noted at \cpageref{TXT Basilica Vicsek BubbleBath Houghton QV}), the conjugacy problem and the conjugacy search problem in their rearrangement groups can be solved using closed strand diagrams.

\begin{remark}
    This algorithm is arguably not the most efficient that it can be.
    For example, an algorithm that checks whether the two base lines describe cohomologous cycles could also provide an explicit sequence of shifts that bring the first base line to the second, allowing to partially skip the subsequent step of the algorithm.
    There might be more places where the algorithm may be optimized, which might be of interest for possible applications to group-based cryptography (see for example section 3 of \cite{crypto}), but a full analysis of the efficiency of this method is beyond the scope of this work.
\end{remark}

\subsection{Non Reduction-confluent Edge Replacement Rules}
\label{sub.non.confluent}

\cref{thm.conjugacy.solvable} shows that reduction-confluence of the edge replacement system (\cref{def.reduction.confluent}) is a sufficient condition for solving the conjugacy problem with closed strand diagrams.
However, we can still adapt this method to solve the conjugacy problem in certain rearrangement groups that do not have this property.

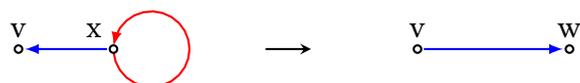
\begin{figure}\centering
\begin{tikzpicture}[scale=1]
    \begin{scope}
    \draw[edge,red,domain=-175:175] plot ({.5*cos(\x)}, {.5*sin(\x)});
    \node[vertex] (l) at (-1.75,0) {}; \draw (-1.75,0) node[above]{v};
    \node[vertex] (cl) at (-.5,0) {}; \draw (-.5,0) node[above left]{x};
    \draw[edge,blue] (cl) to (l);
    \end{scope}
    \draw[thick, -stealth] (1.5,0) -- (2.1,0);
    \node[vertex] (vblue) at (3.5,0) {};
    \node[vertex] (wblue) at (5.5,0) {};
    \draw[edge,blue] (vblue) node[above,black]{v} -- (wblue) node[above,black]{w};
\end{tikzpicture}
\caption{The additional graph reduction that allows us to solve the conjugacy problem in the airplane rearrangement group.}
\label{fig_airplane_graph_reduction}
\end{figure}

For example, at \cpageref{TXT airplane non confluent} we showed that the airplane edge replacement system $\mathcal{A}$ is not reduction-confluent.
However, this issue is only caused by the fact that, whenever we have a subgraph that is isomorphic to the base graph of $\mathcal{A}$ and that features at least one vertex of degree 1, this subgraph can be reduced in two different manners, as shown in \cref{fig_TA_reduced} (w is a new symbol).
It is intuitive that such a small issue can be overcome, and indeed one can modify the reduction system (defined in \cref{SUB reduction systems}) by adding the reduction described in \cref{fig_airplane_graph_reduction}.
This additional rule means that a subgraph isomorphic to the one on the left in the figure, provided no other edge is incident on x, can be replaced with the one on the right, maintaining the adjacencies of the vertex v.
This can be stated more precisely and generally using DPO graph rewriting systems, as will be discussed in the following \cref{sub.DPO}, but the idea should be clear without having to introduce that machinery.

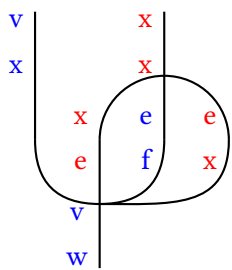
\begin{figure}\centering
\begin{tikzpicture}[font=\small,scale=.85]
    \draw (0,0) -- node[left,blue,align=center,looseness=1.2]{v\\x} (0,-1) -- (0,-2);
    \draw (2,0) -- node[left,red,align=center,looseness=1.2]{x\\x} (2,-1);
    \draw (2,-1) to[out=180,in=90] (1,-2) node[left,red,align=center]{x\\e};
    \draw (2,-1) to (2,-2) node[left,blue,align=center]{e\\f};
    \draw (2,-1) to[out=0,in=90] (3,-2) node[left,red,align=center]{e\\x};
    \draw (0,-2) to[out=270,in=180,looseness=1.2] (1,-3);
    \draw (1,-2) to (1,-3);
    \draw (2,-2) to[out=270,in=0,looseness=1.2] (1,-3);
    \draw (3,-2) to[out=270,in=0,looseness=1.25] (1,-3);
    \draw (1,-3) -- node[left,blue,align=center]{v\\w} (1,-4);
\end{tikzpicture}
\caption{The conjugator corresponding to a Type 3 reduction associated with the graph reduction shown in \cref{fig_airplane_graph_reduction}.}
\label{fig_airplane_conj}
\end{figure}

It is not hard to prove that this newly built reduction system (comprised of the two original reductions alongside the new one) is locally confluent and terminating, thus by the usual argument we find that each connected component has a unique reduced graph.
The reason why this translates to conjugacy is that the newly added reduction of graphs is precisely the same as a red edge expansion followed by a blue reduction, thus we must allow a special kind of Type 3 reduction that, from looping strands as in the bottom left of \cref{fig_TA_reduced}, produces a looping strand as in the bottom right of that same figure.
Then one can prove an analogous result to \cref{PROP conjugator} by adding the following rule to the list at \cpageref{TXT transformations to conjugacy}:
\begin{itemize}
    \item If the transformation is a Type 3 reduction of the newly added kind, then the open strand diagrams differ from conjugation by the diagram in \cref{fig_airplane_conj} (a number of these diagrams equal to the winding number of the loops).
\end{itemize}
It is easy to check that the new rewriting system is confluent, by the same reasoning used to prove \cref{LEM reduced CSD} along with the fact that the reduction system comprised of the two original reductions alongside the new one is confluent.

Thus, collectively we have the following result.

\medskip %layout
\begin{theorem}\label{thm.conjugacy.known.groups}
Strand diagrams can be used to solve the conjugacy problem and the conjugacy search problem in the Higman-Thompson groups, the airplane rearrangement group $T_A$, the rabbit rearrangement groups (including the basilica rearrangement group $T_B$) the Vicsek and bubble bath rearrangement groups, the Houghton groups $H_n$, the groups $QV$, $QT$, $QF$ and those topological full groups of edge shifts whose edge replacement system is reduction-confluent.
\end{theorem}

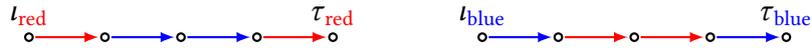
\begin{figure}\centering
\begin{subfigure}{.4\textwidth}\centering
    \begin{tikzpicture}[scale=1]
        \node[vertex] (1) at (1,0) {};
        \node[vertex] (2) at (2,0) {};
        \node[vertex] (3) at (3,0) {};
        \node[vertex] (4) at (4,0) {};
        \node[vertex] (5) at (5,0) {};
        \draw[edge,red] (1) node[above,black]{$\iota_{\text{\textcolor{red}{red}}}$} to (2);
        \draw[edge,blue] (2) to (3);
        \draw[edge,blue] (3) to (4);
        \draw[edge,red] (4) to (5) node[above,black]{$\tau_{\text{\textcolor{red}{red}}}$};
    \end{tikzpicture}
\end{subfigure}
\begin{subfigure}{.4\textwidth}\centering
    \begin{tikzpicture}[scale=1]
        \node[vertex] (1) at (1,0) {};
        \node[vertex] (2) at (2,0) {};
        \node[vertex] (3) at (3,0) {};
        \node[vertex] (4) at (4,0) {};
        \node[vertex] (5) at (5,0) {};
        \draw[edge,blue] (1) node[above,black]{$\iota_{\text{\textcolor{blue}{blue}}}$} to (2);
        \draw[edge,red] (2) to (3);
        \draw[edge,red] (3) to (4);
        \draw[edge,blue] (4) to (5) node[above,black]{$\tau_{\text{\textcolor{blue}{blue}}}$};
    \end{tikzpicture}
\end{subfigure}
\caption{The \textcolor{red}{red} and \textcolor{blue}{blue} replacement graphs for a non confluent-reduction set of edge replacement rules where closed strand diagrams do not seem to easily solve the conjugacy problem.}
\label{fig replacement Bad}
\end{figure}

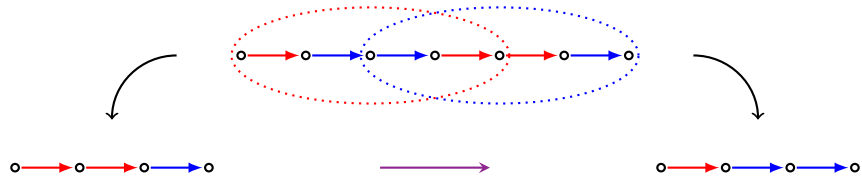
\begin{figure}\centering
\begin{tikzpicture}[scale=.85]
    \begin{scope}[xshift=-3cm]
    \node[vertex] (1) at (0,0) {};
    \node[vertex] (2) at (1,0) {};
    \node[vertex] (3) at (2,0) {};
    \node[vertex] (4) at (3,0) {};
    \node[vertex] (5) at (4,0) {};
    \node[vertex] (6) at (5,0) {};
    \node[vertex] (7) at (6,0) {};
    \draw[edge,red] (1) to (2);
    \draw[edge,blue] (2) to (3);
    \draw[edge,blue] (3) to (4);
    \draw[edge,red] (4) to (5);
    \draw[edge,red] (5) to (6);
    \draw[edge,blue] (6) to (7);
    \draw[red,dotted] (2,0) ellipse (2.15cm and .75cm);
    \draw[thick,-to] (-1,0) to[out=180,in=90] (-2,-1);
    \draw[blue,dotted] (4,0) ellipse (2.15cm and .75cm);
    \draw[thick,-to] (7,0) to[out=0,in=90] (8,-1);
    \end{scope}
    \begin{scope}[xshift=-9.5cm,yshift=-1.75cm]
    \node[vertex] (4) at (3,0) {};
    \node[vertex] (5) at (4,0) {};
    \node[vertex] (6) at (5,0) {};
    \node[vertex] (7) at (6,0) {};
    \draw[edge,red] (4) to (5);
    \draw[edge,red] (5) to (6);
    \draw[edge,blue] (6) to (7);
    \end{scope}
    \begin{scope}[xshift=3.5cm,yshift=-1.75cm]
    \node[vertex] (1) at (0,0) {};
    \node[vertex] (2) at (1,0) {};
    \node[vertex] (3) at (2,0) {};
    \node[vertex] (4) at (3,0) {};
    \draw[edge,red] (1) to (2);
    \draw[edge,blue] (2) to (3);
    \draw[edge,blue] (3) to (4);
    \end{scope}
    \draw[Plum,thick,-stealth] (-.85,-1.75) -- (.85,-1.75);
\end{tikzpicture}
\caption{Adding \textcolor{Plum}{this reduction} to the reduction system associated to the edge replacement rules of \cref{fig replacement Bad} does not make the reduction system confluent, differently from how the new reduction of \cref{fig_airplane_graph_reduction} worked for the airplane reduction system.}
\label{fig replacement Bad_conf}
\end{figure}

This does not seem to immediately apply to edge replacement rules such as the one portrayed in \cref{fig replacement Bad}, as adding the obvious reduction depicted in \cref{fig replacement Bad_conf} creates new non-confluent paths of reductions.
It may still be possible that adding a large (but finite) amount of new reductions might solve this.
It should also be noted that the rearrangement group for this specific edge replacement system is a diagram group, thus the conjugacy problem is solvable if and only if the semigroup $\langle r,b \mid r = r b b r, b = b r r b \rangle$ ($r$ and $b$ standing for red and blue, respectively) has solvable word problem, by a result in \cite{guba1997diagram}.

\subsection{Graph Rewriting Systems and Confluence}
\label{sub.DPO}

Edge replacement rules (\cref{def.replacement.rules}), reduction systems (\cref{SUB reduction systems}) and the modified reduction system discussed for the airplane from the previous \cref{sub.non.confluent} are all instances of \textbf{graph rewriting systems}, which are specific instances of abstract rewriting systems whose objects are graphs.
The modern DPO (double pushout) approach to graph rewriting systems was first introduced in \cite{4569741}, but the precise definition may vary, as explained and investigated in \cite{habel_müller_plump_2001}.

An application of this general theory is the Critical Pair Lemma, which is Lemma 10 from \cite{Plump2005}.
It states that it suffices to check specific pairs of reductions (called \textit{critical pairs}) in order to establish that a graph rewriting system is locally confluent (which implies being confluent if it is terminating).
Our proof that the basilica edge replacement system is reduction-confluent (\cpageref{TXT Basilica Vicsek BubbleBath Houghton QV}) can actually be seen as an instance of this Lemma.
However, the Critical Pair Lemma does not provide a necessary condition, and in general the problem of deciding whether a graph rewriting system is confluent is undecidable (Theorem 5 of \cite{Plump2005}).
Nevertheless, since reductions of expanding edge replacement systems are a very specific and arguably quite simple family of graph rewriting systems, it may very well be that checking for confluence in this context is decidable.
Checking whether this is the case or not is an interesting question that we will not tackle here.

Finally, what we did in the previous \cref{sub.non.confluent} to solve the conjugacy problem for the airplane edge replacement system can be called a \textit{confluent completion} of the reduction system: we have added new reductions that make the reduction system confluent, while preserving termination and the equivalence generated by the rewriting.
In the different context of \textit{term rewriting systems}, the Knuth-Bendix completion algorithm (\cite{Knuth1983}) does exactly what we just described.
However, such an algorithm in the larger context of graph rewriting systems has not been developed yet, and is currently of interest to the community of graph rewriting systems (as stated in section 5 of \cite{Plump2005}).
So it is possible that some future application of the theory of graph rewriting systems will apply to reduction systems of edge replacement systems, possibly providing a wider sufficient condition for the feasibility of the method presented here to solve the conjugacy problem in rearrangement groups.
If a completion algorithm can be produced for the reduction system of any edge replacement system, then the conjugacy (search) problem in every rearrangement group is solved by building confluent completions that allow algorithmic comparisons of closed strand diagrams.

%%%%%%%%%%%%%%%%%%%%%%%%%

\chapter{Invariable Generation}
\label{cha.IG}

A group $G$ is \textbf{invariably generated} (\textbf{IG} for short) if there exists a subset $S \subseteq G$ such that, for every choice $g_s \in G$ for $s \in S$, the group $G$ is generated by $\{ s^{g_s} \mid s \in S \}$ (with $g^h$ we mean $h^{-1} g h$).
In other words, $G$ is invariably generated if there exists a generating set $S$ such that, if one modifies $S$ by conjugating each $s \in S$ by some $g_s \in G$ (possibly a distinct $g_s$ for each $s$), then they still end up with a generating set for $G$.

This notion has been present in the literature for a long time, for example in \cite{Wiegold1976TransitiveGW, Wiegold1976TransitiveGW2} under a different but equivalent definition.
However, the term \textit{Invariable Generation} was first used in \cite{DIXON199225} by Dixon in his study about computational Galois theory.
Finite groups are invariably generated, hence interesting questions in that setting generally involve sizes of invariable generating sets \cite{KANTOR2011302, GARZONI2020218}.
On the other hand, the context of infinite groups presents both examples of IG and non-IG groups: Houghton groups (defined in \cref{sub.Houghton}) are IG for all $n \geq 2$ \cite{COX2022120}, whereas certain convergence groups (in particular hyperbolic groups) and certain arithmetic groups are non-IG \cite{8200290, 10.1093/imrn/rnw137}.
Moreover, Gelander, Golan and Juschenko proved that Thompson group $F$ is (finitely) invariably generated (later generalized in \cite{FFIG}), while Thompson groups $T$ and $V$ are not invariably generated \cite{Gelander2016InvariableGO}.
Since Thompson groups have so many generalizations, it is natural to ask whether invariable generation can be studied for those generalizations using similar methods.

This chapter presents the content of Davide Perego's and the author's work \cite{IG}, which generalize the results from \cite{Gelander2016InvariableGO} and their proofs of non-invariable generation for Thompson groups $T$ and $V$ to the setting of rearrangement groups.
The main result can be summarized as follows.

\medskip %layout
\begin{maintheorem*}
\label{THM main intro}
Every CO-transitive subgroup $G$ of a rearrangement group $G_\mathcal{X}$ is not invariably generated.
\end{maintheorem*}

\phantomsection\label{txt.CO.transitivity}
By \textbf{CO-transitive} (compact-open transitive \cite{CO-compact}) we mean that $G$ acts on a space $X$ in such a way that, for each proper compact $K$ and each non-empty open $U$ of $X$, there is an element of $G$ that maps $K$ inside $U$.
We will work with a condition (\textit{weak cell-transitivity}, introduced in \cref{sec.wct}) that is equivalent to CO-transitivity and easier to verify in the language of edge replacement systems.

The Main Theorem has itself a generalization which is worth mentioning:

\medskip %layout
\begin{corollary*}
If $1 \neq N \trianglelefteq G \leq G_\mathcal{X}$, where $G_\mathcal{X}$ is a rearrangement group and $G$ is weakly cell-transitive, then $N$ is not invariably generated.
\end{corollary*}

In particular, many known rearrangement groups are CO-transitive, so we have the following new results:

\medskip %layout
\begin{corollary*}
The following groups (and their commutator subgroups) are not invariably generated:
\begin{itemize}
    \item the Higman-Thompson groups $T_{n,r}$ and $V_{n,r}$ (\cref{sub.higman.thompson.groups});
    \item the basilica and rabbit rearrangement groups (\cref{sub.basilica});
    \item the airplane rearrangement group $T_A$ (\cref{sub.airplane});
    \item the dendrite rearrangement groups (\cref{sub.other.fractals,cha.dendrites});
    \item the Vicsek rearrangement group and its generalizations (\cite{BF19});
    \item topological full groups of edge shifts (\cref{sub.topological.full.groups}).
\end{itemize}
\end{corollary*}

We notice that there are currently no known necessary or sufficient conditions to establish whether a subgroup of a rearrangement group is itself a rearrangement group (except for stabilizers of rational points, see \cref{sec.rational.stabilizers}), but that our results give sufficient conditions for non-invariable generation among subgroups of a rearrangement group.

\section{Weak Cell-Transitivity}
\label{sec.wct}

Let $G_\mathcal{X}$ be the rearrangement group associated to an expanding edge replacement system $\mathcal{X} = (X_0, R, C)$ (\cref{def.expanding}), where $X_0$ is the base graph, $R = \{ X_i \mid i \in C \}$ is the set of replacement graphs and $C$ is the set of colors.
We denote by $X$ the limit space associated to $\mathcal{X}$.

\begin{definition}
\label{def.wct}
A subgroup $G$ of a rearrangement group $G_\mathcal{X}$ is \textbf{weakly cell-transitive} when, for each cell $C$ and each proper union of finitely many cells $A$, there exists a $g \in G$ such that $g(A) \subseteq C$.
\end{definition}

As far as the writer know, this class contains most of the rearrangement groups that have been studied so far, along with some notable subgroups.
Indeed, the basilica rearrangement group $T_B$ and its commutator subgroup $[T_B, T_B]$ (\cref{sub.basilica}), the airplane rearrangement group $T_A$ and $[T_A, T_A]$ (\cref{sub.airplane}), the dendrite rearrangement groups (\cref{cha.dendrites}) and the rearrangement group of the Vicsek fractal (Example 2.1 of \cite{BF19}, see also \cref{sub:Vic}) are all examples of weakly cell-transitive rearrangement groups.
Indeed, taking $T_A$ as a model, it is easy to argue that for every finite union of cells there exists a rearrangement that maps it into one of the four cells of the base graph.
It can also be seen that there are three orbits of cells (red cells, ``external'' blue cells and ``internal'' blue cells) and that every cell contains an element of each orbit.
Combining the two, one can immediately conclude that $T_A$ is weakly cell-transitive, and a similar argument can be applied to the commutator subgroup $[T_A, T_A]$, the basilica and dendrite rearrangement groups $T_B$ and $G_n$, their commutator subgroups $[T_B, T_B]$ and $[G_n, G_n]$ and the Vicsek rearrangement group.

Thompson groups $T$ and $V$ (\cref{sub.Thompson}) are weakly cell-transitive too, and in truth their actions are actually transitive on the set of cells.
The same holds for Higman-Thompson groups $T_{n,r}$ and $V_{n,r}$ (\cref{sub.higman.thompson.groups}).

As a counterexample, Thompson group $F$ (\cref{sub.Thompson}) is instead not weakly cell-transitive, since a union of cells containing an endpoint of $[0,1]$ cannot be mapped inside a cell that does not contain said endpoint.
Indeed, in \cite{Gelander2016InvariableGO} it has been shown that $F$ is invariably generated.
For the same reasons, Higman-Thompson groups $F_{n,r}$ (\cref{sub.higman.thompson.groups}) are not weakly cell-transitive, although whether these groups are invariably generated is yet to be investigated.

\subsection{Other classes of actions by homeomorphisms}

In this subsection we explore the relationship between weakly cell-transitive actions of rearrangement groups and other classes of group actions on topological spaces.

\phantomsection\label{txt.minimal}
First, recall that an action $G \curvearrowright X$ is said to be \textbf{minimal} if the orbit of each point of $X$ is dense in $X$.

\begin{remark}
\label{rmk.density}
A subset $D$ of a limit space $X$ is dense in $X$ if and only if it has non-trivial intersection with each cell of $X$.
Indeed, the topological interior of each cell is non-empty, so if $D$ is dense then it must intersect each cell non-trivially.
Conversely, if $A$ is an open subset of $X$, then it must contain some ball, which in turn must contain some cell by \cref{lem.balls.and.cells};
therefore, if $D$ intersects each cell non-trivially then it intersects each open set non-trivially.
\end{remark}

\medskip %layout
\begin{proposition}
\label{prop.orbit.dense}
For a subgroup $G$ of a rearrangement group $G_\mathcal{X}$, the following statements are equivalent:
\begin{enumerate}
    \item $G$ is weakly cell-transitive (\cref{def.wct});
    \item $G$ is \textbf{flexible} (in the sense of \cite{flexible}, although the term was first used in \cite{Belk2019OnTA}), i.e., the group acts on a compact Hausdorff $X$ in such a way that, for any two proper closed sets of $X$ with non-empty topological interior, there is an element that maps one inside the other;
    \item $G$ is \textbf{CO-transitive} (compact-open transitive, \cite{CO-compact}), i.e., for each proper compact $K$ and each non-empty open $U$, there is an element that maps $K$ inside $U$.
\end{enumerate}
Moreover, if $G_\mathcal{X}$ is weakly cell-transitive, then the action of $G_\mathcal{X}$ on the limit space $X$ is minimal.
\end{proposition}

\begin{proof}
We will prove that (1) is equivalent to (2). One can prove the equivalence between (1) and (3) essentially in the same way.
To prove this we will use the fact that each cell contains balls of radius greater than zero (\cref{lem.balls.and.cells}).
First, recall that the limit space of an expanding edge replacement system is compact and Hausdorff (\cref{thm.limit.space}).
If $G$ is flexible, $C$ is a cell and $A$ is a proper union of finitely many cells, then $C$ and $A$ are proper closed subsets of the limit space $X$ with non-empty topological interior, so we are done.
Conversely, consider two proper closed subsets $E_1$ and $E_2$ of $X$ with non-empty topological interior.
Since the complement $E_1^\mathsf{C}$ is proper, open and nonempty, it must contain some proper topological interior $\mathring{C_0}$ of a cell $C_{0}$, so $E_1$ is contained in $\mathring{C_0}^\mathsf{C}$.
This is the union of finitely many cells since, for every cellular partition (\cref{def.cellular.partition}) containing $C_0$, by \cref{rmk.topological.interior} the complement of $\mathring{C_0}$ is the union of every other cell of the partition.
Also, $E_2$ contains $\mathring{E_2}$, which is non-empty, so it contains some cell $C$.
Then, since $G$ is weakly cell-transitive, there exists some $g \in G$ such that $g(\mathring{C_0}^\mathsf{C}) \subseteq C$, thus $g(E_1) \subseteq g(\mathring{C_0}^\mathsf{C}) \subseteq C \subseteq E_2$.

Finally, suppose that $G$ is weakly cell-transitive.
In order to show that the action is minimal, let $p$ be a point of the limit space $X$ and consider a cell $A$ that contains $p$.
By hypothesis, for any cell $C$ there exists a $g \in G$ such that $g(A) \subseteq C$.
Then $g(p) \in C$, so the orbit of $p$ intersects $C$ non-trivially, and by \cref{rmk.density} we are done.
\end{proof}

\phantomsection\label{txt.vigorous}
Additionally, since the action of a rearrangement group is defined on cells, it induces an action on the symbol space (which is the action that allows us to embed rearrangement groups in topological full groups of edge shifts, see \cref{sub.embedding.into.TFGoES}).
Recall that the symbol space is often a Cantor space (\cref{prop.shift.is.Cantor}).
A group acting on a Cantor space is \textbf{vigorous} if, for any clopen subset $A$ of the Cantor set, for any two clopen proper subsets $B$ and $C$ of $A$, there is a $g$ in the pointwise stabilizer of the complement of $A$ that maps $B$ inside $C$.
If the symbol space is a Cantor space and if a subgroup of a rearrangement group is vigorous with its action on the symbol space, then it is weakly cell-transitive.
Indeed, take $A$ as the entire Cantor set: cells correspond to ``glued'' cones, which are clopen.

A necessary condition for weak cell-transitivity is that there are edges of every type among the expansions of each replacement graph (recall that the type of an edge depends on its color and on whether it is a loop or not, \cref{def.cell.type}).
In truth, this condition is very broad, as it is also a necessary condition for the minimality of the rearrangement group.
Indeed, if the replacement graph $R_c$ does not feature edges of a certain type among its expansions, then any point included in some edge of that type cannot be mapped inside a $c$-colored cell.
This shows that, for example, the Houghton groups (\cref{sub.Houghton}) and the Thompson-like groups $QF$, $QT$ and $QV$ (\cref{sub.thompson.like}) are not weakly cell-transitive.
In general, finding a fully combinatorial necessary or sufficient condition for weak cell-transivitiy seems unlikely, as even edge replacement systems that do not seem to have special properties may end up producing trivial rearrangement groups (see \cref{fig.replacement.rearrangementless}, which is \cite[Example 2.5]{BF19}).

\begin{figure}
\centering
\begin{tikzpicture}[scale=1.25]
    \node at (-.667,0) {$X_0 =$};
    \node[vertex] (L) at (0,0) {};
    \node[vertex] (R) at (1.25,0) {};
    \draw[edge,domain=190:530] plot ({.5*cos(\x)+.5}, {.5*sin(\x)});
    \draw[edge,domain=190:530] plot ({.5*cos(\x)+1.75}, {.5*sin(\x)});
    \draw (.5,.5) node[above] {$L$};
    \draw (1.75,.5) node[above] {$R$};
    \begin{scope}[xshift=4.5cm]
    \node at (-.667,0) {$X_1 = $};
    \node[vertex] (l) at (0,0) {};
    \draw (l) node[above]{$\iota_1$};
    \node[vertex] (c) at (1.1,0) {};
    \node[vertex] (r) at (2.2,0) {};
    \draw (r) node[above]{$\tau_1$};
    \draw[edge] (l) to node[above]{$0$} (c);
    \draw[edge] (c) to node[above]{$1$} (r);
    \end{scope}
\end{tikzpicture}
\caption{An edge replacement system whose rearrangement group has minimal action but is not weakly cell-transitive.}
\label{fig.double.T}
\end{figure}
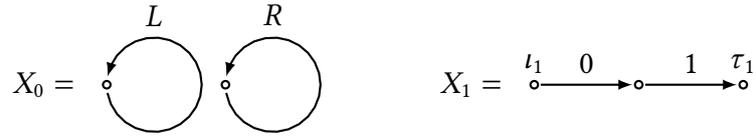

An example of a rearrangement group whose action is minimal but not weakly cell-transitive is given by the edge replacement system of \cref{fig.double.T}, whose rearrangement group is $T^2 \rtimes \mathrm{Sym}(2)$.
This sort of rearrangement groups that are obtained from ``products'' of edge replacement systems (as will be described in \cref{sub.products}) seem to be often minimal but not weakly cell-transitive, so it would be interesting to find an example of minimal and not weakly cell-transitive rearrangement group with connected base graph.

\section{Wandering Sets}

As in \cite{Gelander2016InvariableGO}, we introduce the notion of wandering sets, which will be essential in the proof of the \hyperref[THM main]{main Theorem} of this chapter.

\begin{definition}
\label{def.wandering}
Let $G \curvearrowright X$ and $g \in G$. We say that a subset $A \subseteq X$ is:
\begin{itemize}
    \item \textbf{$g$-wandering} if, for every $n \in \mathbb{Z}$ such that $g^n \neq 1_G$, we have that $g^n(A) \cap A = \emptyset$;
    \item \textbf{weakly $g$-wandering} if, for every $n \in \mathbb{Z}$, either $g^n$ fixes $A$ pointwise or $g^n(A) \cap A = \emptyset$.
\end{itemize} 
\end{definition}

The remainder of this section is dedicated to proving the following result:

\medskip %layout
\begin{proposition}
\label{prop.wandering.cells}
Let $G_\mathcal{X}$ be rearrangement group. Then each element $g \in G_\mathcal{X}$ admits a weakly wandering cell.
\end{proposition}

This is easy to see for rearrangement groups of monochromatic edge replacement systems, since they embed nicely inside Higman-Thompson groups $V_{n,r}$, for which the results of \cite{Gelander2016InvariableGO} apply with barely any modification.
In order to prove this result in the general setting of polychromatic edge replacement systems, we will need to generalize the results of \cite{Gelander2016InvariableGO}.
We will distinguish between two cases depending on whether the element of $G_\mathcal{X}$ is torsion or not.
Lemmas \ref{lem.wandering.sets.torsion} and \ref{lem.wandering.non.torsion} contained in the two following subsections will thus prove \cref{prop.wandering.cells}.

Before we prove \cref{prop.wandering.cells}, we point out the following consequence:

\medskip %layout
\begin{corollary}
\label{cor.wct.free.subgroups}
Each weakly cell-transitive subgroup of a rearrangement group contains a copy of the free group on two generators, and in particular it is not amenable.
\end{corollary}

This can be easily proved using \cref{cor.wandering.conjugacy} and the ping-pong Lemma.

\subsection{Torsion Elements}

The Lemma below is a generalization of a property of Thompson groups $F$, $T$ and $V$: if an element $g$ fixes some irrational point $\alpha$ of $[0,1]$, then it must fix some small enough neighborhood of $\alpha$.
Recall that, in the general context of limit spaces of edge replacement systems, a point $p$ is irrational if it has a unique representative that is not eventually periodic (\cref{def.rational.irrational.points}).
Recall that irrational points always exist, and in truth they densely populate the limit space (\cref{rmk.irrational.points.are.dense}).

\medskip %layout
\begin{lemma}
\label{lem.irrational.stabilizer}
If a rearrangement $g$ fixes an irrational point $p$ of $X$, then there exists a small enough cell $C$ containing $p$ that is fixed by $g$ pointwise.
\end{lemma}

\begin{proof}
Suppose $p$ is an irrational point, and let $e_1 e_2 \dots$ be its unique representative, which must be an irrational element of the symbol space $\Omega_\mathcal{R}$.
Let $(D,\sigma,R)$ be a graph pair diagram for $g$ (\cref{def.graph.pair.diagrams}).
Since the point $p$ is irrational, it is included in a unique cell $C_D = \llbracket e_1 \dots e_n \rrbracket$ of $D$, which $g$ maps canonically to some cell $C_R$ of $R$.
Since $p$ has a unique representative and is contained in $C_R$, we must have $C_R = \llbracket e_1 \dots e_m \rrbracket$ for some $m$.
Now, $(D,\sigma,R)$ is a graph pair diagram for $g$, so the action on the cells of $D$ is canonical, meaning that the effect of $g$ on the points represented by some $e_1 \dots e_n \alpha$ is the replacement of the prefix $e_1 \dots e_n$ with $e_1 \dots e_m$ (\cref{def.canonical.homeomorphism}).
In particular, $p = g(p)$ is represented by both $e_1 e_2 \dots e_n e_{n+1} \dots$ and $e_1 \dots e_m e_{n+1} \dots$.
Since $p$ has unique representative, this implies that $e_1 \dots e_n e_{n+1} \dots$ is the same as $e_1 \dots e_m e_{n+1} \dots$.
This is only possible if either $e_1 \dots e_n = e_1 \dots e_m$ or $e_{n+1} e_{n+2} \dots$ is periodic, but the latter is impossible since $p$ is irrational.
\end{proof}

\medskip %layout
\begin{lemma}
\label{lem.wandering.sets.torsion}
If $g \in G_\mathcal{X}$ is torsion, then there exists a weakly $g$-wandering cell.
\end{lemma}

\begin{proof}
Suppose $g$ is a non-trivial torsion element of $G_\mathcal{X}$.
Then, because of \cref{rmk.irrational.points.are.dense}, there must exist an irrational point $p$ that is not fixed by $g$.
Let $n$ be the period of $p$ under the action of $g$.
Since $g, g^2, \dots, g^{n-1}$ are homeomorphisms of $X$ and $X$ is a Hausdorff topological space (\cref{thm.limit.space}) there exist a small enough ball $B$ of $p$ such that $g(B), \dots, g^{n-1}(B)$ are disjoint from $B$.
By \cref{lem.irrational.stabilizer}, we can assume that $g^n$ fixes $B$ pointwise.
Then, by \cref{lem.balls.and.cells}, there exists a cell contained in $B$, and so that cell is weakly $g$-wandering.
\end{proof}

\subsection{Non-torsion Elements}
\label{sub.IG.non.torsion}

Following the strategy of \cite{Gelander2016InvariableGO}, we need to gather some results about the dynamics of rearrangements in order to prove the existence of wandering cells for non-torsion elements.
The following Lemmas are generalizations of Lemmas 10.2, 10.3 and 10.5 of \cite{Brin2004HigherDT}, along with the general structure of its section 10, which is about \textit{revealing pairs} of trees for elements of $V$.

\phantomsection\label{txt.caret}
We are going to be using the forest pair diagrams from \cref{sec.forest.pair.diagrams}.
We do not really need to worry about the labelings of such diagrams, here.
What matters is that graph pair diagrams are in bijection with forest pair diagrams, so that each rearrangement can be represented by a bijection between the leaves of a pair of forest and that such diagrams can be expanded by attaching a tree to a bottom branch.
In the remainder of this section, by \textbf{caret} of a forest expansion we will mean a branch $b$ that terminates at an interior node $v$ together with all of the branches that originate from $v$, the amount of which depends on the color of the upper branch $b$.
Recall that carets are instances of replacement trees and represent edge expansions, as was described in \cref{sub.forests.and.graphs}.

\phantomsection\label{txt.forest.difference}
Given a forest pair diagram $(F_D,F_R)$, we can compute $F_D - F_R$ (and $F_R - F_D$), by which we mean the set of carets that belong to $F_D$ but not to $F_R$.
These correspond one-to-one to edge expansions performed in $D$ that have not been performed in $R$.

\begin{definition}
\label{def.domain.imbalance}
Given a forest pair diagram $(F_D, F_R)$, its \textbf{domain imbalance} is the number of carets of $F_D - F_R$.
\end{definition}

\begin{remark}
\label{rmk.imbalances}
For polychromatic edge replacement systems, this need not be the same as the number of carets of $F_R - F_D$, which we could call \textit{range imbalance}.
This is no issue since, among all of the representatives for the same rearrangement, domain imbalance and range imbalance always differ by the same number.

Indeed, representatives for the same rearrangements differ by sequences of edge expansions or reductions of domain and range, so it suffices to see how imbalances change when going from $(D,\sigma,R)$ to an edge expansion $(D \triangleleft e, \sigma', R \triangleleft \sigma(e))$, which consists of appending a caret $U$ at $e$ in $F_D$ and the same caret $U$ at $\sigma(e)$ in $F_R$.
If $e = \sigma(e)$, clearly the imbalances remain the same, as the added carets cancel out.
If $e \neq \sigma(e)$, we distinguish cases based on whether $e$ is an internal node of $F_R$ and whether $\sigma(e)$ is an internal node of $F_D$:
\begin{enumerate}
    \item if both happen, then the edge expansion decreases domain imbalance and range imbalance by one;
    \item if only one happens, then both the domain and the range imbalance are left unchanged by the edge expansion;
    \item if neither of the two happen, then the edge expansion increases both domain and range imbalance by one.
\end{enumerate}
In all of these cases, domain and range imbalance change by the same amount.
\end{remark}

\phantomsection\label{txt.components}
When computing $F_D - F_R$, we can also consider its \textbf{components}, by which we mean the pairwise disjoint ``trees'' in which $F_D$ breaks up after carets of $F_R$ have been removed.
Each component represents a maximal set of edge expansions that are ``not independent'', by which we mean expansions of edges $e_1, \dots, e_k$ where $D = E \triangleleft e_1 \triangleleft \dots \triangleleft e_k$ such that, for $i \geq 2$, each $e_i$ is not an edge of $E$ but is instead created by subsequent edge expansions.

From now on, we will make the following assumption about the graph pair diagrams that we consider.

\begin{assumption}
\label{ass.minimality.choice}
Assume that $(D,\sigma,R)$ is a graph pair diagram with the least domain imbalance of all of the representatives of the same rearrangement (which is the same as assuming it has least range imbalance, by \cref{rmk.imbalances}).
Among the representatives with minimal domain imbalance, we choose one that has the smallest number of components of $F_D-F_R$.
Among all such representatives, we choose one that has the smallest number of components of $F_R-F_D$.
In the remainder of this subsection we will prove facts about graph pair diagrams $(D,\sigma,R)$ chosen in this manner.
\end{assumption}

\phantomsection\label{txt.expandable.sequence}
Before we start with the first Lemma, consider a sequence $u_1, \dots, u_n$ of edges of $D$ such that $u_1, \dots, u_n$ and $\sigma(u_n)$ are pairwise distinct and such that $u_{i+1} = \sigma(u_i)$ for all $1 \leq i < n$.
We call this an \textbf{expandable sequence} of edges of $D$ (leaves of $F_D$).
Observe that, in particular, these edges must have the same color, and that $u_2, \dots, u_n, \sigma(u_n)$ must be edges of $R$.
We can then perform the same edge expansion (or sequence of edge expansions) on each of $u_1, \dots, u_n$ in $D$ and each of $\sigma(u_1), \dots, \sigma(u_n)$ in $R$, obtaining a new graph pair diagram that represents the same rearrangement.
We call this operation an \textbf{iterated expansion}.
Observe that, when looking at the pair of forests $F_D$ and $F_R$, this operation corresponds exactly to attaching the same subtree (representing the edge expansion) to the leaves $u_1, \dots, u_n$ in $D$ and to the leaves $\sigma(u_1), \dots, \sigma(u_n)$ in $R$.

\medskip %layout
\begin{lemma}
\label{LEM Brin 1}
There is no expandable sequence $u_1, \dots, u_n$ such that $u_1$ is an interior node of $F_R$ and $\sigma(u_n)$ is an interior node of $F_D$.
\end{lemma}

\begin{proof}
Suppose by contradiction that the statement does not hold and consider the component $U$ of $F_D - F_R$ with root at $\sigma(u_n)$.
We can perform an iterated expansion by $U$ along $u_1, \dots, u_n$, which causes the following modifications on $F_D - F_R$:
\begin{itemize}
    \item the copy of $U$ appended at $\sigma(u_n)$ is removed;
    \item the copies of $U$ appended to both $F_D$ and $F_R$ at each $u_i$ for $2 \leq i \leq n$ do not contribute to $F_D - F_R$, as they cancel out;
    \item since $u_1$ is an interior node of $F_R$, the number of carets added to $F_D$ at $u_1$ from the copy of $U$ is strictly less than the number of carets of $U$.
\end{itemize}
This would lower the domain imbalance, which is impossible by \cref{ass.minimality.choice}.
\end{proof}

\medskip %layout
\begin{lemma}
\label{LEM Brin 2}
There is no expandable sequence $u_1, \dots, u_n$ such that all of the following hold:
\begin{enumerate}
    \item $u_1$ is not a node of $F_R$,
    \item $\sigma(u_n)$ is an interior node of $F_D$,
    \item the component of $F_D - F_R$ rooted at $\sigma(u_n)$ and the one containing $u_1$ are not the same component.
\end{enumerate}
\end{lemma}

\begin{proof}
Suppose by contradiction that there is such an expandable sequence $u_1, \dots, u_n$.
Let $U$ and $V$ respectively be the components of $F_D - F_R$ rooted at $\sigma(u_n)$ and containing $u_1$, and consider an iterated expansion by $U$ along the sequence $u_1, \dots, u_n$.
This iterated expansion does not change the domain imbalance, but it removes the component $U$ from $F_D - F_R$, while introducing no new component to $F_D - F_R$.
Indeed, the subtrees appended in both $F_D$ and $F_R$ to $u_i$ for $2 \leq i \leq n$ cancel out in $F_D - F_R$ and a copy of $U$ is appended at the already existing component $V$.
This is not possible by the assumption of minimal number of components of $F_D - F_R$ among those with minimal domain imbalance (\cref{ass.minimality.choice}).
\end{proof}

\medskip %layout
\begin{lemma}
\label{LEM Brin 3}
There is no expandable sequence $u_1, \dots, u_n$ such that all of the following hold:
\begin{enumerate}
    \item $\sigma(u_n)$ is not a node of $F_D$,
    \item $u_1$ is an interior node of $F_R$,
    \item the component of $F_R - F_D$ rooted at $u_1$ and the one containing $\sigma(u_n)$ are not the same component.
\end{enumerate}
\end{lemma}

\begin{proof}
This proof is dual to that of the previous Lemma.
Suppose by contradiction that there is such an expandable sequence $u_1, \dots, u_n$.
Let $U$ and $V$ respectively be the components of $F_R - F_D$ containing $\sigma(u_n)$ and rooted at $u_1$ and consider an iterated expansion by $V$ along the sequence $u_1, \dots, u_n$.
This iterated expansion does not change the domain imbalance, but it removes the component $V$ from $F_R - F_D$, while introducing no new component to $F_R - F_D$.
Indeed, the subtrees appended in both $F_D$ and $F_R$ to $u_i$ for $2 \leq i \leq n$ cancel out in $F_R - F_D$ and a copy of $V$ is appended at the already existing component $U$.
This is not possible by our assumption of minimal number of components of $F_R - F_D$ among those with minimal number of components of $F_D - F_R$, which in turn were chosen among those with minimal domain imbalance (\cref{ass.minimality.choice}).
Note that the argument works because the number of components of $F_D - F_R$ has not increased, since $\sigma(u_n)$ is not a node of $F_D$ by assumption.
\end{proof}

\medskip %layout
\begin{lemma}
\label{LEM Brin 4}
\begin{enumerate}
    \item For each non-trivial component $U$ of $F_D - F_R$, there is a unique leaf $l = l(U)$ of $U$ such that there is an expandable sequence $u_1, \dots, u_n$ starting at $u_1 = l$ and ending at $\sigma(u_n) = r$, where $r = r(U)$ denotes the root of $U$.
    \item For each non-trivial component $V$ of $F_R - F_D$, there is a unique leaf $l = l(V)$ of $V$ such that there is an expandable sequence $u_1, \dots, u_n$ starting at $u_1 = r$ and ending at $\sigma(u_n) = l$, where $r = r(V)$ denotes the root of $V$.
\end{enumerate}
\end{lemma}

\begin{proof}
We will only prove the second point of the Lemma, as the first point is proven in a very similar manner, using \cref{LEM Brin 2} in place of \cref{LEM Brin 3}.

Observe that $r$ is a leaf of $F_D$, so let $u_1 = r$.
This can be regarded as a trivial expandable sequence.
Since $u_1 = r$ is an interior node of $F_R$, by \cref{LEM Brin 1} we have that $\sigma(u_1)$ cannot be an interior node of $F_D$.
For the same reason, by \cref{LEM Brin 3} we have that either $\sigma(u_1)$ is a node of $F_D$ or it belongs to the component $V$ of $F_R - F_D$.
Thus, $\sigma(u_1)$ is either a leaf of $F_D$ or it is a leaf of $V$, and only one of these is possible.

Now, assume that $u_1, \dots, u_n$ is an expandable sequence with $r = u_1$ and with $\sigma(u_n)$ either a leaf of $V$ or a leaf of $F_D$.
There is at least one such expandable sequence, as shown in the previous paragraph.
If $\sigma(u_n)$ is a leaf of $V$, we are done, so assume $\sigma(u_n)$ is a leaf of $F_D$ instead.
We can then create a longer expandable sequence by adding the element $u_{n+1} = \sigma(u_n)$ at the end of the sequence.
This new sequence satisfies the same property of the previous one: since $u_1 = r$ is not an interior node of $F_R$, by \cref{LEM Brin 3} we have that $\sigma(u_n) = u_{n+1}$ can either be a common leaf of $F_D$ and $F_R$ or a leaf of $V$.

This expandable sequence cannot extend infinitely.
Indeed, if this were not the case, since there is a finite amount of leaves of $F_D$ among which we can choose and $g$ is an invertible map, for some $n$ we would have $u_1 = \sigma(u_n)$.
But $u_1$ would be a leaf of $F_R$, whereas it should be an interior node of $F_R$.
This means that we must reach a leaf $l$ with the property we are looking for, so we are done.
\end{proof}

It follows from the second point of \cref{LEM Brin 4} that, for each non-torsion rearrangement $g$, if we choose a graph pair diagram $(D,\sigma,R)$ that satisfies \cref{ass.minimality.choice}, then for each leaf $r$ of $F_D$ that is an interior node in $F_R$, there is a leaf $l$ of $F_R$ located below $r$ and an $n \in \mathbb{N}^* = \{1, 2, \dots \}$ such that $r, g(r), \dots, g^{n-1}(r)$ are distinct leaves of $F_D$ and $g^n(r) = l$.
In terms of graph pair diagrams, this means that for each edge $e$ of $D$ that has been expanded in $R$, there is an edge $e^*$ of $R$ generated by an edge expansion of $e$ and an $n \in \mathbb{N}^*$ such that the topological interiors of the cells generated by the edges $e, g(e), \dots, g^{n-1}(e)$ of $D$ are such that $g^n(e) = e^*$ and they are pairwise disjoint (observe that boundary vertices do not interfere with this because of \cref{rmk.topological.interior}).
Essentially, given the topological interior $\mathring{C}$ of a cell whose generating edge is in $D$ and has been expanded in $R$, iterating $g$ maps $\mathring{C}$ to disjoint topological interiors of cells for $n-1$ steps, and then to a topological interior contained inside $\mathring{C}$ itself at the $n$-th step.

\medskip %layout
\begin{lemma}
\label{lem.wandering.non.torsion}
If $g \in G_\mathcal{X}$ is not torsion, then there exists a $g$-wandering cell.
\end{lemma}

\begin{proof}
Let $(D,\sigma,R)$ be a graph pair diagram for $g$ as discussed above.
Since $g$ is not torsion, $D \neq R$, so there exists some edge $e$ of $D$ that has been expanded in $R$.
Consider the edge $e^*$ and the number $n \in \mathbb{N}^*$ given by the discussion above.
Let $f$ be an edge of $R$ that is distinct from $e^*$ and that has been generated by an edge expansion of $e$, so $\llbracket f\rrbracket \subseteq \llbracket e \rrbracket$.
We claim that the topological interior $\llparenthesis f \rrparenthesis$ is $g$-wandering, which would conclude the proof, as it implies that every cell contained in $\llparenthesis f \rrparenthesis$ is also $g$-wandering.

Indeed, let $0 \neq m \in \mathbb{Z}$.
For every $U \subseteq X$, if $g^m(U) \cap U \neq \emptyset$ then $U \cap g^{-m}(U) \neq \emptyset$, so we can assume that $m > 0$.
Let $m = qn + r$, where $q \geq 0$ and $r \in \{0, 1, \dots, n-1\}$.
Since $g^n(\llparenthesis e \rrparenthesis) = \llparenthesis e^* \rrparenthesis \subseteq \llparenthesis e \rrparenthesis$, we have $g^{qn}( \llparenthesis e \rrparenthesis) \subseteq \llparenthesis e^* \rrparenthesis$, and so $g^{qn}( \llparenthesis f \rrparenthesis) \subseteq \llparenthesis e^* \rrparenthesis$.

Now, if $r=0$ we are done, as $f$ and $e^*$ being distinct edges of $D$ implies $\llparenthesis f \rrparenthesis \cap \llparenthesis e^* \rrparenthesis = \emptyset$.
If instead $r \neq 0$, then $g^m( \llparenthesis f \rrparenthesis) = g^{qn + r} ( \llparenthesis f \rrparenthesis) \subseteq g^r( \llparenthesis e \rrparenthesis)$, which is disjoint from $\llparenthesis e \rrparenthesis$ because $r \in \{1, \dots, n-1\}$.
Thus, since $\llparenthesis f \rrparenthesis \subseteq \llparenthesis e \rrparenthesis$ we have $g^m( \llparenthesis f \rrparenthesis) \cap \llparenthesis f \rrparenthesis = \emptyset$, so we are done.
\end{proof}

\section{Proof of the Theorem}

The results of the previous section give us the following Corollary:

\medskip %layout
\begin{corollary}
\label{cor.wandering.conjugacy}
Suppose that $G$ is a weakly cell-transitive subgroup of some rearrangement group $G_\mathcal{X}$.
Let $g \in G$ and let $A \subset X$ be a union of finitely many cells of $\mathcal{X}$. Then there exists some $h \in G$ such that $A$ is weakly $g^h$-wandering.
\end{corollary}

\begin{proof}
Because of \cref{prop.wandering.cells}, there exists a weakly $g$-wandering cell $I \subset X$.
Since $G$ is weakly cell-transitive, there exists an $h \in G$ such that $A \subseteq h^{-1} (I)$.
Now, it is clear that $h^{-1} (I)$ is weakly $g^h$-wandering, so $A$ is too and we are done.
\end{proof}

\medskip %layout
\begin{theorem}
\label{THM main}
Every weakly cell-transitive subgroup $G$ of a rearrangement group $G_\mathcal{X}$ is not invariably generated.
\end{theorem}

Note that this is equivalent to the \hyperref[THM main intro]{Main Theorem} stated in the introduction of this chapter because of \cref{prop.orbit.dense}.

\begin{proof}
Observe that there are at most countably many non-trivial conjugacy classes of $G$, since every rearrangement group $G_\mathcal{X}$ is itself at most countably infinite
(indeed, each rearrangement can be represented as a graph isomorphism of some full expansion graph $E_n$, as defined in \cref{def.full.expansion}, and for each $n \in \mathbb{N}$ there are finitely many such isomorphisms).
First, assume that there are exactly countably many conjugacy classes, so let $\{Cl_n\}_{n \in \mathbb{N}^*}$ be the non-trivial conjugacy classes of $G$, where $\mathbb{N}^* = \{1, 2, \dots \}$.

Consider a point $p$ of the limit space $X$.
Recall that each point of the limit space is the equivalence class under the gluing relation (\cref{def.glue}) of some element of the symbol space (\cref{def.symbol.space});
each element of the symbol space corresponds to an infinite sequence of cells decreasing with respect to set inclusion.
Let $C_0 C_1 C_2 \dots$ be such a sequence for the point $p$, i.e., $p = \llbracket w_0 w_1 w_2 \dots \rrbracket$ and $C_i = \llbracket w_0 \dots w_i \rrbracket$.
For each $n \in \mathbb{N}^* = \{ 1, 2, \dots \}$, let $I_n$ be the set difference of the topological interior of $C_n$ with the cell $C_{n+1}$.
We intentionally do not consider $n = 0$ in order to make sure that there is at least some cell that is not contained in the union of these $I_n$, as the cell $C_0$ is the whole limit space whenever the base graph consists solely of an edge.
Note that these $I_n$'s are pairwise disjoint and converge to the point $p$.

Now, for each $n \in \mathbb{N}^*$, the complement $I_n^{\mathsf{C}}$ is a non-trivial union of finitely many cells.
More explicitly, if $E$ is a graph expansion containing the edge $e$ corresponding to $C_n$, then consider $E \triangleleft e$: the set $I_n^{\mathsf{C}}$ is the union of all cells of $E \triangleleft e$ that also appeared in $E$ along with $C_{n+1}$.
Hence, by \cref{cor.wandering.conjugacy} there exists a $\gamma_n$ in the $n$-th conjugacy class $Cl_n$ such that $I_n^{\mathsf{C}}$ is weakly $\gamma_n$-wandering.

We claim that the orbit of $p$ under the action of the subgroup $H = \langle \gamma_{n} \mid n \in \mathbb{N}^* \rangle$ of $G$ is contained in $I = \bigcup_{n \in \mathbb{N}^*} I_n$.
Then, since there is at least a cell that is not contained in $I$, the action of $H$ is not minimal.
So, by virtue of \cref{prop.orbit.dense}, we would get $H \neq G$ as required.
Thus, \cref{LEM rearrangement ping-pong} below completes the proof in the case of countably many conjugacy classes of $G$.

If there are only finitely many conjugacy classes, then we can do everything in exactly the same way, except that we will have a finitely generated $H = \langle \gamma_1, \dots, \gamma_k \rangle$ under whose action the orbit of $p$ is contained in the union of finitely many $I_n$'s, leading to the same conclusion for the same \cref{LEM rearrangement ping-pong}.
\end{proof}

The following Lemma and its proof are essentially the same as Lemma 18 of \cite{Gelander2016InvariableGO}.
The only difference is that our Lemma applies to cells of a limit space instead of intervals of $S^1$.
Even so, we include its proof here for the sake of completeness.

\medskip %layout
\begin{lemma}
\label{LEM rearrangement ping-pong}
Let $q$ be an element of the $H$-orbit of $p$ that is different from $p$.
Let $g \in H$ be an element of minimal word-length over the alphabet $A \coloneq \{ \gamma_n \mid n \in J \}$ such that $g(p) = q$ (where $J$ is either $\mathbb{N}$ or $\{1, \dots,, k \}$).
Assume that $g = \gamma_{i_1}^{k_1} \dots \gamma_{i_m}^{k_m}$, where:
\begin{itemize}
    \item $k_1, \dots, k_m \neq 0$,
    \item $i_{j+1} \neq i_j$ for each $j=1, \dots, m-1$,
    \item $|k_1| + \dots + |k_m|$ is the word-length of $g$ over $A$.
\end{itemize}
Then $q \in  I_{i_m}$.
\end{lemma}

\begin{proof}
The proof is by induction on $m$.
If $m=1$, then $g = \gamma_{i_1}^{k_1}$, and the set $I_{i_1}^{\mathsf{C}}$ is weakly $\gamma_{i_1}$-wandering.
Then, since $\gamma_{i_1}^{k_1} (p) = q \neq p$ and $p \in I_{i_1}^\mathsf{C}$, we have that $\gamma_{i_1}^{k_1}$ does not fix $I_{i_1}^{\mathsf{C}}$ pointwise.
Therefore, the set $\gamma_{i_1}^{k_1} (I_{i_1}^{\mathsf{C}})$ is disjoint from $I_{i_1}^{\mathsf{C}}$, meaning that $\gamma_{i_1}^{k_1} (I_{i_1}^{\mathsf{C}}) \subseteq I_{i_1}$.
Thus, since $\gamma_{i_1}^{k_1} (p) = q \notin I_{i_1}$, we must have $q \in I_{i_1}$ as needed.

Assume the lemma holds for $m$ and let $g = \gamma_{i_1}^{k_1} \dots \gamma_{i_m}^{k_m} \gamma_{i_{m+1}}^{k_{m+1}}$.
Let $g_1 = \gamma_{i_1}^{k_1} \dots \gamma_{i_m}^{k_m}$ and $r = g_1 (p)$.
By induction hypothesis, $r \in I_{i_m}$.
The facts that the $I_n$'s are pairwise disjoint and $i_m \neq i_{m+1}$ imply that $r \notin I_{i_{m+1}}$.
By assumption of minimality of the word-length of $g$, we have $q = g(p) = \gamma_{i_{m+1}}^{k_{m+1}} (r) \neq r$, so $\gamma_{i_{m+1}}^{k_{m+1}}$ does not fix $I_{i_{m+1}}^{\mathsf{C}}$ pointwise.
Therefore, since $I_{i_{m+1}}^{\mathsf{C}}$ is weakly $\gamma_{i_{m+1}}$-wandering, we must have $\gamma_{i_{m+1}}^{k_{m+1}} (I_{i_{m+1}}^{\mathsf{C}}) \subseteq I_{i_{m+1}}$.
Now, $r \notin I_{i_{m+1}}$, so $q = \gamma_{i_{m+1}}^{k_{m+1}} (r) \in I_{i_{m+1}}$.
\end{proof}

As a final remark, it is easy to see that the proofs of the previous \hyperref[cor.wandering.conjugacy]{Corollary} and \hyperref[THM main]{Theorem} can also be applied to normal subgroups of a weakly cell-transitive subgroup of a rearrangement group, providing the following more general statement:

\medskip %layout
\begin{proposition}
\label{PROP main}
If $1 \neq N \trianglelefteq G \leq G_\mathcal{X}$, where $G_\mathcal{X}$ is a rearrangement group and $G$ is weakly cell-transitive, then $N$ is not invariably generated.
\end{proposition}

In particular, we have the following notable result:

\medskip %layout
\begin{corollary}
The commutator subgroup of a weakly cell-transitive rearrangement group is not invariably generated.
\end{corollary}

%%%%%%%%%%%%%%%%%%%%%%%%%

\chapter{Rationality of the Gluing Relation}
\label{cha.rational.gluing}

In this chapter, which covers the content of Davide Perego's and the author's work \cite{rationalgluing}, we construct automata that recognize pairs of infinite sequences precisely when they are equivalent under the gluing relation that defines the limit space, meaning that the gluing relation is \textit{rational}.

Different ways to define rationality have been employed in the literature, especially in symbolic dynamics with notions like edge shifts (see \cref{sub.edge.shifts}) and sofic shifts (see \cite[Chapter 3]{LindMarcus}).
We have already described how the well-known \hyperref[thm.Alexandroff.Hausdorff]{Alexandroff-Hausdorff Theorem} allows us to define a \textit{coding} any compact metrizable space $X$:
every point of $X$ is associated to an equivalence class of infinite sequences of a Cantor space $\mathfrak{C}$ (often an edge shift) that are ``glued'' by an equivalence relation.
For instance, we recall the natural and simple \cref{ex.interval}, which is the binary coding of the unit interval.
While a general coding can be quite unwieldy, interesting and useful codings generally arise when there is an explicit rule for gluing two infinite sequences together.
This is what motivates our interest in the \textit{rationality} of gluing relations:
this property essentially expresses that there is indeed a rule that described when two sequences glue together, and that rule is ``computable''.

There are several examples of this in the literature, though they may be addressed with different names or be stronger or weaker versions of what we are describing here.
The Gromov boundary of a hyperbolic group has rational gluing with respect to different codings.
In \cite{CP,P}, it is proved that it is \textit{semi-Markovian} with respect to the boundary of the tree of geodesics and other trees introduced by Gromov himself in \cite{G}.
Being semi-Markovian is generally a stronger notion than rationality (see again \cite{P}).
In \cite{Perego}, it is shown that the gluing associated to the \textit{horofunction boundary} is rational. Limit spaces of contracting self similar groups have rational gluing with respect to the regular tree on which the groups act on (see \cite{N2}).
In \cite{Bandt}, Bandt explores the opposite point of view, building self-similar spaces starting from automata that define gluing relations and our construction fits in Bandt's formalism.

Here we prove that the gluing relation produced by edge replacement systems is rational,
as conjectured in \cite[Example 7.14]{Perego}.
We provide an explicit construction that associates to any expanding edge replacement system a finite-state automaton that recognizes pairs of equivalent infinite sequences.
We will work with expanding edge replacement systems (\cref{def.expanding}), but we have seen in \cref{prop.null.expanding.isolated.rarrangement} that this is not a strong assumption.

\section{The Gluing Automaton}

In this section, after establishing a useful assumption about colors and loops and expressing how automata (\cref{def.automaton}) can be used to describe a relation between infinite sequences (\cref{def.rational.relation}), we construct gluing automata for the gluing relation defined by edge replacement systems (\cref{def.glue}).

\subsection{Colors and Loops}
\label{sub.loops}

Throughout this chapter, we make the following assumption on the graph expansions of our edge replacement systems, as it will make it easier to build our gluing automaton.

\begin{assumption}
\label{ass.loops}
We assume that edges of the same color are either all loops or all not loops.
\end{assumption}

Under this assumption, it is natural to identify the two boundary vertices for colors that represent loops.
Each replacement graph $\Gamma_c$ is then equipped with either one or two boundary vertices:
a single vertex $\lambda_c$ for loops and two distinct vertices $\iota_c$ and $\tau_c$ for non-loops.

\begin{figure}
\centering
\begin{tikzpicture}
    \node at (-.8,0) {$\Gamma_0 =$};
    \node[vertex] (i) at (1,0) {};
    \draw[edge,red] (i) to[out=150, in=-150, min distance=1.2cm] node[above left]{$L$} (i);
    \draw[edge,red] (i) to[out=-30, in=30, min distance=1.2cm] node[above right]{$R$} (i);
    \begin{scope}[xshift=4.25cm]
    \node at (-.825,0) {$\Gamma_{\textcolor{blue}{1}} =$};
    \node[vertex] (l) at (0,0) {};
    \draw (l) node[above]{$\iota_{\textcolor{blue}{1}}$};
    \node[vertex] (c) at (1,0) {};
    \node[vertex] (r) at (2,0) {};
    \draw (r) node[above]{$\tau_{\textcolor{blue}{1}}$};
    \draw[edge,blue] (l) to node[above]{$1$} (c);
    \draw[edge,red] (c) to[out=60, in=120, min distance=1.2cm] node[above right]{$2$} (c);
    \draw[edge,blue] (c) to node[above]{$3$} (r);
    \end{scope}
    \begin{scope}[xshift=8.5cm]
    \node at (-1,0) {$\Gamma_{\textcolor{red}{2}} =$};
    \node[vertex] (l) at (0,0) {};
    \draw (l) node[left]{$\lambda_{\textcolor{red}{2}}$};
    \node[vertex] (c) at (1,0) {};
    \draw[edge,blue] (l) to[out=290, in=240] node[below]{$4$} (c);
    \draw[edge,red] (c) to[out=-30, in=30, min distance=1.2cm] node[above right]{$5$} (c);
    \draw[edge,blue] (c) to[out=110,in=60] node[above]{$6$} (l);
    \end{scope}
\end{tikzpicture}
\caption{The basilica edge replacement system $\mathcal{B}$ under \cref{ass.loops}, which can be obtained as a modification of \cref{fig.basilica.replacement}.}
\label{fig.basilica.replacement.loops}
\end{figure}
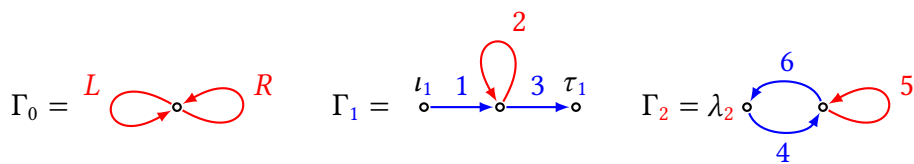

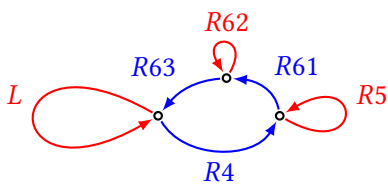
\begin{figure}
    \centering
    \begin{tikzpicture}[scale=2]
        \node[vertex] (i) at (1,0) {};
        \node[vertex] (r) at (1.8,0) {};
        \node[vertex] (rt) at (1.45,.25) {};
        \draw[edge,red] (i) to[out=150, in=-150, min distance=1.2cm] node[above left]{\small$L$} (i);
        \draw[edge,blue] (i) to[out=300, in=240] node[below]{\small$R4$} (r);
        \draw[edge,red] (r) to[out=-30, in=30, min distance=.6cm] node[above right]{\small$R5$} (r);
        \draw[edge,blue] (r) to[out=100,in=-5] node[above right]{\small$R61$} (rt);
        \draw[edge,red] (rt) to[out=60, in=120, min distance=.3cm] node[above]{\small$R62$} (rt);
        \draw[edge,blue] (rt) to[out=185,in=50] node[above left]{\small$R63$} (i);
    \end{tikzpicture}
    \vspace{-12pt}
    \caption{A graph expansion of the edge replacement system $\mathcal{B}$.}
    \label{fig.basilica.expansion.loops}
\end{figure}

\begin{remark}
This assumption is not really a restriction on the limit spaces that we can produce, nor on their rearrangement groups.
Indeed, given a generic edge replacement systems, we can always split a color $c$ into two colors $c_\text{loop}$ and $c_\text{non-loop}$ and create a new edge replacement system that produces the same limit space and the same rearrangement group.
For example, \cref{fig.basilica.replacement.loops} depicts the basilica edge replacement system (originally portrayed in \cref{fig.basilica.replacement}) with this modification and \cref{fig.basilica.expansion.loops} shows a graph expansion.

More precisely, the new edge replacement system $\mathcal{R}^*$ is obtained by coloring by $c_\text{loop}$ or $c_\text{non-loop}$ the $c$-colored edges of the graphs of the original edge replacement system $\mathcal{R}$ depending on whether they are loops or not, after having identified $\iota_c$ and $\tau_c$ with a unique vertex $\lambda_c$ in the $c_\text{loop}$ replacement graph.
Since rearrangements preserve the type of cells (which is the color together with the property being loops or not, \cref{def.cell.type}), the rearrangement groups of $\mathcal{R}$ and $\mathcal{R}^*$ are truly the same.
\end{remark}

\subsection{Construction of the Gluing Automaton}
\label{sub.gluing.automaton}

From here on, we will be frequently using the definitions about automata that were given in \cref{sub.automata}.
In particular, recall that an infinite sequence is recognized by an automaton if the automaton features a run that reads the sequence on its labels.

\begin{definition}
\label{def.rational.relation}
Let $\Omega$ be a set of infinite sequences.
We say that an equivalence relation on $\Omega$ is \textbf{rational} when there exists an automaton that recognizes a pair $(w, v)$ if and only if the two elements $w, v \in \Omega$ are related.
\end{definition}

Note that an automaton that recognizes a relation on $\Sigma^\omega = \{ x_1 x_2 \ldots \mid x_i \in \Sigma \}$ must have $\Sigma^2$ as alphabet.
With a slight abuse of notation, a sequence $(x_1,y_1) (x_2,y_2) \dots$ here represents a pair of sequences $x_1 x_2 \ldots$ and $y_1 y_2 \ldots$ belonging to $\Omega_\mathcal{R}$.

For the rest of this chapter, consider an expanding edge replacement system $\mathcal{R}$ with colors $\mathrm{C} = \{1, \dots, k\}$, base graph $\Gamma_0$ and replacement graphs $\Gamma_1, \dots, \Gamma_k$.
The gluing automaton $\mathrm{Gl}_\mathcal{R} = \left( \Sigma^2, Q, \rightarrow, q_0(0) \right)$ is built as follows.

\subsubsection{The alphabet}

$\Sigma^2$ consists of all possible pairs of symbols of $\Sigma$, where $\Sigma$ has exactly one symbol for each edge of the base and replacement graphs, i.e.,
\[ \Sigma^2 = \{ (x,y) \mid x,y \in \mathbb{A} \}, \]
where $\mathbb{A}$ denotes the alphabet of the symbol space, which is set of edges of the base and replacement graphs (see \cref{prop.language.alphabet}).
As expressed by \cref{prop.language.alphabet}, when processing one of the two entries of the $l$-th letter of a sequence in this alphabet, it will be natural to refer to that entry as if it were an edge of the $l$-th graph of the full expansion sequence (i.e., the sequence of graphs obtained by expanding, at each step, every edge of the previous graph, starting with the base graph, as in \cref{def.full.expansion}).
In this way, the two entries of the $l$-th letter of a sequence in $\Sigma^2$ represent two edges of the same graph, so in what follows we will simply describe their possible adjacency in the $l$-th graph of the full expansion sequence without explicitly mentioning the graph itself.

\subsubsection{The set of states}
\label{sub.states.gluing.automaton}

$Q$ is the union of the two following sets $Q_0$ and $Q_1$:
\[ Q_0 = \big\{ q_0(i) \mid i = 0, 1, \dots, k \big\}, \]
\[ Q_1 = \big\{ q_1 \big(\begin{smallmatrix} i & \gamma \\ j & \delta \end{smallmatrix}\big) \mid i, j \in \mathrm{C}, \gamma, \delta \in \mathcal{E}, \gamma \neq \mathrm{db}^-, \delta=\mathrm{db}^\pm \iff \gamma=\mathrm{db}^+ \big\}, \]
where $\mathcal{E} = \{\mathrm{in}, \mathrm{out}, \mathrm{lp}, \mathrm{db}^+, \mathrm{db}^- \}$.
The initial state is $q_0(0) \in Q_0$.

The following subsection will describe in more detail what these states mean.
The intuition is that, while processing the edges of our two sequences, the first set $Q_0$ consists of an initial state $q_0(0)$ along with a state $q_0(i)$ for each color $i$ that indicates that the two entries being processed are the same $i$-colored edge, while the second set $Q_1$ has a state for each quadruple $\big(\begin{smallmatrix} i & \gamma \\ j & \delta \end{smallmatrix}\big)$ where $i, j$ are the colors of the edges and $\gamma, \delta$ distinguish between the types of adjacency between the edges.
The symbols $\mathrm{lp}$ and $\mathrm{db}$ stand for \textit{loop} and \textit{double}, respectively, and their meaning will be made clear when describing the transition functions.

\subsubsection{The transitions}
\label{sub.transitions.gluing.automaton}

The transition function $\rightarrow$ is the partial mapping $Q \times \Sigma \to Q$ defined as follows, where we distinguish between three different types of transitions.
\begin{description}
    \item[Type 0]
        $q_0(i) \xrightarrow{(e,e)} q_0 (\mathrm{c}(e))$ for all $i \in \{0\} \cup \mathrm{C}$ and all edges $e$ of $\Gamma_i$.
    \item[Type 1]
        $q_0(i) \xrightarrow{(a,b)} q_1 \big(\begin{smallmatrix} \mathrm{c}(a) & \alpha \\ \mathrm{c}(b) & \beta \end{smallmatrix}\big)$ for all $i \in \{0\} \cup \mathrm{C}$ and all distinct edges $a$ and $b$ of $\Gamma_i$, where $\alpha$ and $\beta$ depend on the adjacency between $a$ and $b$ in $\Gamma_i$ as shown in \cref{tab:transition:type:1}.
    \item[Type 2]
        $q_1 \big(\begin{smallmatrix} i & \gamma \\ j & \delta \end{smallmatrix}\big) \xrightarrow{(a,b)} q_1 \big(\begin{smallmatrix} \mathrm{c}(a) & \alpha \\ \mathrm{c}(b) & \beta \end{smallmatrix}\big)$ for all $i,j \in \{0\} \cup \mathrm{C}$ whenever the edges $a$ and $b$ of $\Gamma_i$ and $\Gamma_j$, respectively, satisfy the conditions described in \cref{tab:transition:type:2}, which also depend on $\gamma$ and $\delta$ and determine $\alpha$ and $\beta$.
\end{description}

The meaning of these transitions (especially Type 2) will be made clear in \cref{sec:gl:aut}.
The idea is that Type 0 transitions simply induce a copy of the color graph (\cref{def.color.graph}), Type 1 transitions arise when two sequences first diverge and Type 2 transitions describe the adjacency of sequences after their divergence.
Note that, when processing the pair $(a,b)$ in a state $q_0(i)$, there is a Type 1 transition if and only if $a$ and $b$ are adjacent in $\Gamma_i$.

Observe that loops at $\lambda_i$ may arise when the original edge replacement system (before gluing the initial and terminal vertices as explained in \cref{sub.loops}) had loops at $\iota$ or at $\tau$.
Since the original edge replacement system is expanding, this is the only possible occurrence of the column $\iota(a) = \tau(a) = \lambda_i$ (row $\iota(b) = \tau(b) = \lambda_j$) of \cref{tab:transition:type:2}.

Moreover, in Type 2 transitions $\alpha$ and $\beta$ will never be $\mathrm{db}^+$ nor $\mathrm{db}^-$, since the original edge replacement system is expanding.
Furthermore, note that when $\gamma$ (or $\delta$) is $\mathrm{db}^+$ or $\mathrm{db}^-$ one needs to check the incidence of $a$ (or $b$) with both the initial and terminal vertices $\iota_i$ and $\tau_i$ (or $\iota_j$ and $\tau_j$) of the replacement graph.
However, since the edge replacement system is expanding, an edge of a replacement graph cannot be incident on both the initial and terminal vertices;
thus, at most one of these can occur, so there is no ambiguity in the table.

\begin{table}
\centering
\renewcommand{\arraystretch}{1.3333}
\begin{tabularx}{\textwidth}{c|c|c|c}
    \begin{tikzpicture}
        \node[vertex] (A) at (0,0) {};
        \node[vertex] (B) at (1,0) {};
        \node[vertex] (C) at (2,0) {};
        \draw[edge] (A) to node[above]{$a$} (B);
        \draw[edge] (B) to node[above]{$b$} (C);
    \end{tikzpicture}
    &%
    \begin{tikzpicture}
        \node[vertex] (A) at (0,0) {};
        \node[vertex] (B) at (1,0) {};
        \node[vertex] (C) at (2,0) {};
        \draw[edge] (A) to node[above]{$a$} (B);
        \draw[edge] (C) to node[above]{$b$} (B);
    \end{tikzpicture}
    &%
    \begin{tikzpicture}
        \node[vertex] (A) at (0,0) {};
        \node[vertex] (B) at (1,0) {};
        \node[vertex] (C) at (2,0) {};
        \draw[edge] (B) to node[above]{$a$} (A);
        \draw[edge] (B) to node[above]{$b$} (C);
    \end{tikzpicture}
    &%
    \begin{tikzpicture}
        \node[vertex] (A) at (0,0) {};
        \node[vertex] (B) at (1,0) {};
        \node[vertex] (C) at (2,0) {};
        \draw[edge] (B) to node[above]{$a$} (A);
        \draw[edge] (C) to node[above]{$b$} (B);
    \end{tikzpicture}%
    \\
    $q_1 \big(\begin{smallmatrix} \mathrm{c}(a) & \mathrm{in} \\ \mathrm{c}(b) & \mathrm{out} \end{smallmatrix}\big)$%
    &%
    $q_1 \big(\begin{smallmatrix} \mathrm{c}(a) & \mathrm{in} \\ \mathrm{c}(b) & \mathrm{in} \end{smallmatrix}\big)$%
    &%
    $q_1 \big(\begin{smallmatrix} \mathrm{c}(a) & \mathrm{out} \\ \mathrm{c}(b) & \mathrm{out} \end{smallmatrix}\big)$%
    &%
    $q_1 \big(\begin{smallmatrix} \mathrm{c}(a) & \mathrm{out} \\ \mathrm{c}(b) & \mathrm{in} \end{smallmatrix}\big)$%
    \\[.2cm]\hline
    \begin{tikzpicture}
        \node[vertex] (A) at (0,0) {};
        \node[vertex] (B) at (1,0) {};
        \draw[edge] (A) to node[above]{$a$} (B);
        \draw[edge] (B) to[loop right, min distance=1cm] node[right]{$b$} (B);
    \end{tikzpicture}
    &%
    \begin{tikzpicture}
        \node[vertex] (A) at (0,0) {};
        \node[vertex] (B) at (1,0) {};
        \draw[edge] (B) to node[above]{$a$} (A);
        \draw[edge] (B) to[loop right, min distance=1cm] node[right]{$b$} (B);
    \end{tikzpicture}
    &%
    \begin{tikzpicture}
        \node[vertex] (B) at (0,0) {};
        \node[vertex] (C) at (1,0) {};
        \draw[edge] (B) to[loop left, min distance=1cm] node[left]{$a$} (B);
        \draw[edge] (B) to node[above]{$b$} (C);
    \end{tikzpicture}
    &%
    \begin{tikzpicture}
        \node[vertex] (B) at (0,0) {};
        \node[vertex] (C) at (1,0) {};
        \draw[edge] (B) to[loop left, min distance=1cm] node[left]{$a$} (B);
        \draw[edge] (C) to node[above]{$b$} (B);
    \end{tikzpicture}%
    \\
    $q_1 \big(\begin{smallmatrix} \mathrm{c}(a) & \mathrm{in} \\ \mathrm{c}(b) & \mathrm{lp} \end{smallmatrix}\big)$%
    &%
    $q_1 \big(\begin{smallmatrix} \mathrm{c}(a) & \mathrm{out} \\ \mathrm{c}(b) & \mathrm{lp} \end{smallmatrix}\big)$%
    &%
    $q_1 \big(\begin{smallmatrix} \mathrm{c}(a) & \mathrm{lp} \\ \mathrm{c}(b) & \mathrm{out} \end{smallmatrix}\big)$%
    &%
    $q_1 \big(\begin{smallmatrix} \mathrm{c}(a) & \mathrm{lp} \\ \mathrm{c}(b) & \mathrm{in} \end{smallmatrix}\big)$%
    \\[.2cm]\hline
    \begin{tikzpicture}
        \node[vertex] (B) at (0,0) {};
        \draw[edge] (B) to[loop left, min distance=1cm] node[left]{$a$} (B);
        \draw[edge] (B) to[loop right, min distance=1cm] node[right]{$b$} (B);
    \end{tikzpicture}
    &%
    \begin{tikzpicture}
        \node[vertex] (A) at (0,0) {};
        \node[vertex] (C) at (1.5,0) {};
        \draw[edge] (A) to[bend right] node[above]{$a$} (C);
        \draw[edge] (A) to[bend left] node[above]{$b$} (C);
    \end{tikzpicture}
    &%
    \begin{tikzpicture}
        \node[vertex] (A) at (0,0) {};
        \node[vertex] (C) at (1.5,0) {};
        \draw[edge] (A) to[bend left] node[above]{$a$} (C);
        \draw[edge] (C) to[bend left] node[above]{$b$} (A);
    \end{tikzpicture}
    &%
    \\
    $q_1 \big(\begin{smallmatrix} \mathrm{c}(a) & \mathrm{lp} \\ \mathrm{c}(b) & \mathrm{lp} \end{smallmatrix}\big)$%
    &%
    $q_1 \big(\begin{smallmatrix} \mathrm{c}(a) & \mathrm{db}^+ \\ \mathrm{c}(b) & \mathrm{db}^+ \end{smallmatrix}\big)$%
    &%
    $q_1 \big(\begin{smallmatrix} \mathrm{c}(a) & \mathrm{db}^+ \\ \mathrm{c}(b) & \mathrm{db}^- \end{smallmatrix}\big)$%
    &%
\end{tabularx}
\vspace{.2cm}
\caption{Type 1 transitions and the types of adjacency.}
\label{tab:transition:type:1}
\end{table}

\begin{sidewaystable}
%\begin{table}
\renewcommand{\arraystretch}{1.2}
\centering
\scriptsize
%\begin{tabular}{*{28}{|K{.365cm}}|}
\begin{tabular}{*{20}{|c}|}
%\begin{tabular}{|c|c|c|c|c|c|c|c|c|c|c|c|c|c|c|c|c|c|c|c|c|c|c|c|c|c|c|c|}
%\begin{tabularx}{\textwidth}{|C|C|C|C|C|C|C|C|C|C|C|C|C|C|C|C|C|C|C|C|C|C|C|C|C|C|C|C|}
    \cline{1-20}
    \multicolumn{3}{|c|}{} & $\gamma$ &%
    \multicolumn{3}{c|}{$\mathrm{in}$}&%
    \multicolumn{3}{c|}{$\mathrm{out}$}&%
    \multicolumn{3}{c|}{$\mathrm{lp}$}&%
    \multicolumn{6}{c|}{$\mathrm{db}^+$}&%
    \multicolumn{1}{c|}{$\mathrm{db}^-$}\\
    \cline{4-20}
    \multicolumn{3}{|c|}{} & $\iota(a)$ &%
    $\tau_i$ & & $\tau_i$ &%
    $\iota_i$ & & $\iota_i$ &%
    $\lambda_i$ & & $\lambda_i$ &%
    $\tau_i$ & & $\tau_i$ & $\iota_i$ & & $\iota_i$ & any%
    \\
    \cline{4-20}
    \multicolumn{3}{|c|}{} & $\tau(a)$ &%
    & $\tau_i$ & $\tau_i$ &%
    & $\iota_i$ & $\iota_i$ &%
    & $\lambda_i$ & $\lambda_i$ &%
    & $\tau_i$ & $\tau_i$ & & $\iota_i$ & $\iota_i$ & any%
    \\
    \cline{1-20}
    $\delta$ & $\iota(b)$ & $\tau(b)$ & &%
    \multicolumn{3}{c|}{} & \multicolumn{3}{c|}{} & \multicolumn{3}{c|}{} & \multicolumn{6}{c|}{} &\\
    \cline{1-20}
    %
    % 1st row
    \multirow{3}{*}{$\mathrm{in}$} &
    $\tau_j$ & & &%
    $\begin{smallmatrix} \mathrm{out} \\ \mathrm{out} \end{smallmatrix}$ & $\begin{smallmatrix} \mathrm{in} \\ \mathrm{out} \end{smallmatrix}$ & $\begin{smallmatrix} \mathrm{lp} \\ \mathrm{out} \end{smallmatrix}$ &%
    $\begin{smallmatrix} \mathrm{out} \\ \mathrm{out} \end{smallmatrix}$ & $\begin{smallmatrix} \mathrm{in} \\ \mathrm{out} \end{smallmatrix}$ & $\begin{smallmatrix} \mathrm{lp} \\ \mathrm{out} \end{smallmatrix}$ &%
    $\begin{smallmatrix} \mathrm{out} \\ \mathrm{out} \end{smallmatrix}$ & $\begin{smallmatrix} \mathrm{in} \\ \mathrm{out} \end{smallmatrix}$ & $\begin{smallmatrix} \mathrm{lp} \\ \mathrm{out} \end{smallmatrix}$ &%
    & & & & & &%
    \multirow{21}{*}{\rotatebox[origin=c]{270}{empty}}\\
    \cline{2-3}\cline{5-19}
    %
    % 2nd row
    & & $\tau_j$ & &%
    $\begin{smallmatrix} \mathrm{out} \\ \mathrm{in} \end{smallmatrix}$ & $\begin{smallmatrix} \mathrm{in} \\ \mathrm{in} \end{smallmatrix}$ & $\begin{smallmatrix} \mathrm{lp} \\ \mathrm{in} \end{smallmatrix}$ &%
    $\begin{smallmatrix} \mathrm{out} \\ \mathrm{in} \end{smallmatrix}$ & $\begin{smallmatrix} \mathrm{in} \\ \mathrm{in} \end{smallmatrix}$ & $\begin{smallmatrix} \mathrm{lp} \\ \mathrm{in} \end{smallmatrix}$ &%
    $\begin{smallmatrix} \mathrm{out} \\ \mathrm{in} \end{smallmatrix}$ & $\begin{smallmatrix} \mathrm{in} \\ \mathrm{in} \end{smallmatrix}$ & $\begin{smallmatrix} \mathrm{lp} \\ \mathrm{in} \end{smallmatrix}$ &%
    & & & & & &%
    \\
    \cline{2-3}\cline{5-19}
    %
    % 3rd row
    & $\tau_j$ & $\tau_j$ & &%
    $\begin{smallmatrix} \mathrm{out} \\ \mathrm{lp} \end{smallmatrix}$ & $\begin{smallmatrix} \mathrm{in} \\ \mathrm{lp} \end{smallmatrix}$ & $\begin{smallmatrix} \mathrm{lp} \\ \mathrm{lp} \end{smallmatrix}$ &%
    $\begin{smallmatrix} \mathrm{out} \\ \mathrm{lp} \end{smallmatrix}$ & $\begin{smallmatrix} \mathrm{in} \\ \mathrm{lp} \end{smallmatrix}$ & $\begin{smallmatrix} \mathrm{lp} \\ \mathrm{lp} \end{smallmatrix}$ &%
    $\begin{smallmatrix} \mathrm{out} \\ \mathrm{lp} \end{smallmatrix}$ & $\begin{smallmatrix} \mathrm{in} \\ \mathrm{lp} \end{smallmatrix}$ & $\begin{smallmatrix} \mathrm{lp} \\ \mathrm{lp} \end{smallmatrix}$ &%
    & & & & & &%
    \\
    \cline{1-19}
    %
    % 4th row
    \multirow{3}{*}{$\mathrm{out}$} &
    $\iota_j$ & & &%
    $\begin{smallmatrix} \mathrm{out} \\ \mathrm{out} \end{smallmatrix}$ & $\begin{smallmatrix} \mathrm{in} \\ \mathrm{out} \end{smallmatrix}$ & $\begin{smallmatrix} \mathrm{lp} \\ \mathrm{out} \end{smallmatrix}$ &%
    $\begin{smallmatrix} \mathrm{out} \\ \mathrm{out} \end{smallmatrix}$ & $\begin{smallmatrix} \mathrm{in} \\ \mathrm{out} \end{smallmatrix}$ & $\begin{smallmatrix} \mathrm{lp} \\ \mathrm{out} \end{smallmatrix}$ &%
    $\begin{smallmatrix} \mathrm{out} \\ \mathrm{out} \end{smallmatrix}$ & $\begin{smallmatrix} \mathrm{in} \\ \mathrm{out} \end{smallmatrix}$ & $\begin{smallmatrix} \mathrm{lp} \\ \mathrm{out} \end{smallmatrix}$ &%
    & & & & & &%
    \\
    \cline{2-3}\cline{5-19}
    %
    % 5th row
    & & $\iota_j$ & &%
    $\begin{smallmatrix} \mathrm{out} \\ \mathrm{in} \end{smallmatrix}$ & $\begin{smallmatrix} \mathrm{in} \\ \mathrm{in} \end{smallmatrix}$ & $\begin{smallmatrix} \mathrm{lp} \\ \mathrm{in} \end{smallmatrix}$ &%
    $\begin{smallmatrix} \mathrm{out} \\ \mathrm{in} \end{smallmatrix}$ & $\begin{smallmatrix} \mathrm{in} \\ \mathrm{in} \end{smallmatrix}$ & $\begin{smallmatrix} \mathrm{lp} \\ \mathrm{in} \end{smallmatrix}$ &%
    $\begin{smallmatrix} \mathrm{out} \\ \mathrm{in} \end{smallmatrix}$ & $\begin{smallmatrix} \mathrm{in} \\ \mathrm{in} \end{smallmatrix}$ & $\begin{smallmatrix} \mathrm{lp} \\ \mathrm{in} \end{smallmatrix}$ &%
    & & & & & &%
    \\
    \cline{2-3}\cline{5-19}
    %
    % 6th row
    & $\iota_j$ & $\iota_j$ & &%
    $\begin{smallmatrix} \mathrm{out} \\ \mathrm{lp} \end{smallmatrix}$ & $\begin{smallmatrix} \mathrm{in} \\ \mathrm{lp} \end{smallmatrix}$ & $\begin{smallmatrix} \mathrm{lp} \\ \mathrm{lp} \end{smallmatrix}$ &%
    $\begin{smallmatrix} \mathrm{out} \\ \mathrm{lp} \end{smallmatrix}$ & $\begin{smallmatrix} \mathrm{in} \\ \mathrm{lp} \end{smallmatrix}$ & $\begin{smallmatrix} \mathrm{lp} \\ \mathrm{lp} \end{smallmatrix}$ &%
    $\begin{smallmatrix} \mathrm{out} \\ \mathrm{lp} \end{smallmatrix}$ & $\begin{smallmatrix} \mathrm{in} \\ \mathrm{lp} \end{smallmatrix}$ & $\begin{smallmatrix} \mathrm{lp} \\ \mathrm{lp} \end{smallmatrix}$ &%
    & & & & & &%
    \\
    \cline{1-19}
    %
    % 7th row
    \multirow{3}{*}{$\mathrm{lp}$} &
    $\lambda_j$ & & &%
    $\begin{smallmatrix} \mathrm{out} \\ \mathrm{out} \end{smallmatrix}$ & $\begin{smallmatrix} \mathrm{in} \\ \mathrm{out} \end{smallmatrix}$ & $\begin{smallmatrix} \mathrm{lp} \\ \mathrm{out} \end{smallmatrix}$ &%
    $\begin{smallmatrix} \mathrm{out} \\ \mathrm{out} \end{smallmatrix}$ & $\begin{smallmatrix} \mathrm{in} \\ \mathrm{out} \end{smallmatrix}$ & $\begin{smallmatrix} \mathrm{lp} \\ \mathrm{out} \end{smallmatrix}$ &%
    $\begin{smallmatrix} \mathrm{out} \\ \mathrm{out} \end{smallmatrix}$ & $\begin{smallmatrix} \mathrm{in} \\ \mathrm{out} \end{smallmatrix}$ & $\begin{smallmatrix} \mathrm{lp} \\ \mathrm{out} \end{smallmatrix}$ &%
    & & & & & &%
    \\
    \cline{2-3}\cline{5-19}
    %
    % 8th row
    & & $\lambda_j$ & &%
    $\begin{smallmatrix} \mathrm{out} \\ \mathrm{in} \end{smallmatrix}$ & $\begin{smallmatrix} \mathrm{in} \\ \mathrm{in} \end{smallmatrix}$ & $\begin{smallmatrix} \mathrm{lp} \\ \mathrm{in} \end{smallmatrix}$ &%
    $\begin{smallmatrix} \mathrm{out} \\ \mathrm{in} \end{smallmatrix}$ & $\begin{smallmatrix} \mathrm{in} \\ \mathrm{in} \end{smallmatrix}$ & $\begin{smallmatrix} \mathrm{lp} \\ \mathrm{in} \end{smallmatrix}$ &%
    $\begin{smallmatrix} \mathrm{out} \\ \mathrm{in} \end{smallmatrix}$ & $\begin{smallmatrix} \mathrm{in} \\ \mathrm{in} \end{smallmatrix}$ & $\begin{smallmatrix} \mathrm{lp} \\ \mathrm{in} \end{smallmatrix}$ &%
    & & & & & &%
    \\
    \cline{2-3}\cline{5-19}
    %
    % 9th row
    & $\lambda_j$ & $\lambda_j$ & &%
    $\begin{smallmatrix} \mathrm{out} \\ \mathrm{lp} \end{smallmatrix}$ & $\begin{smallmatrix} \mathrm{in} \\ \mathrm{lp} \end{smallmatrix}$ & $\begin{smallmatrix} \mathrm{lp} \\ \mathrm{lp} \end{smallmatrix}$ &%
    $\begin{smallmatrix} \mathrm{out} \\ \mathrm{lp} \end{smallmatrix}$ & $\begin{smallmatrix} \mathrm{in} \\ \mathrm{lp} \end{smallmatrix}$ & $\begin{smallmatrix} \mathrm{lp} \\ \mathrm{lp} \end{smallmatrix}$ &%
    $\begin{smallmatrix} \mathrm{out} \\ \mathrm{lp} \end{smallmatrix}$ & $\begin{smallmatrix} \mathrm{in} \\ \mathrm{lp} \end{smallmatrix}$ & $\begin{smallmatrix} \mathrm{lp} \\ \mathrm{lp} \end{smallmatrix}$ &%
    & & & & & &%
    \\
    \cline{1-19}
    %
    % 10th row
    \multirow{6}{*}{$\mathrm{db}^+$} &
    $\tau_j$ & & &%
    & & &%
    & & &%
    & & &%
    $\begin{smallmatrix} \mathrm{out} \\ \mathrm{out} \end{smallmatrix}$ & $\begin{smallmatrix} \mathrm{in} \\ \mathrm{out} \end{smallmatrix}$ & $\begin{smallmatrix} \mathrm{lp} \\ \mathrm{out} \end{smallmatrix}$ & & & &%
    \\
    \cline{2-3}\cline{5-19}
    %
    % 11th row
    & & $\tau_j$ & &%
    & & &%
    & & &%
    & & &%
    $\begin{smallmatrix} \mathrm{out} \\ \mathrm{in} \end{smallmatrix}$ & $\begin{smallmatrix} \mathrm{in} \\ \mathrm{in} \end{smallmatrix}$ & $\begin{smallmatrix} \mathrm{lp} \\ \mathrm{in} \end{smallmatrix}$ & & & &%
    \\
    \cline{2-3}\cline{5-19}
    %
    % 12th row
    & $\tau_j$ & $\tau_j$ & &%
    & & &%
    & & &%
    & & &%
    $\begin{smallmatrix} \mathrm{out} \\ \mathrm{lp} \end{smallmatrix}$ & $\begin{smallmatrix} \mathrm{in} \\ \mathrm{lp} \end{smallmatrix}$ & $\begin{smallmatrix} \mathrm{lp} \\ \mathrm{lp} \end{smallmatrix}$ & & & &%
    \\
    \cline{2-3}\cline{5-19}
    %
    % 13th row
    & $\iota_j$ & & &%
    & & &%
    & & &%
    & & &%
    & & & $\begin{smallmatrix} \mathrm{out} \\ \mathrm{out} \end{smallmatrix}$ & $\begin{smallmatrix} \mathrm{in} \\ \mathrm{out} \end{smallmatrix}$ & $\begin{smallmatrix} \mathrm{lp} \\ \mathrm{out} \end{smallmatrix}$ &%
    \\
    \cline{2-3}\cline{5-19}
    %
    % 14th row
    & & $\iota_j$ & &%
    & & &%
    & & &%
    & & &%
    & & & $\begin{smallmatrix} \mathrm{out} \\ \mathrm{in} \end{smallmatrix}$ & $\begin{smallmatrix} \mathrm{in} \\ \mathrm{in} \end{smallmatrix}$ & $\begin{smallmatrix} \mathrm{lp} \\ \mathrm{in} \end{smallmatrix}$ &%
    \\
    \cline{2-3}\cline{5-19}
    %
    % 15th row
    & $\iota_j$ & $\iota_j$ & &%
    & & &%
    & & &%
    & & &%
    & & & $\begin{smallmatrix} \mathrm{out} \\ \mathrm{lp} \end{smallmatrix}$ & $\begin{smallmatrix} \mathrm{in} \\ \mathrm{lp} \end{smallmatrix}$ & $\begin{smallmatrix} \mathrm{lp} \\ \mathrm{lp} \end{smallmatrix}$ &%
    \\
    \cline{1-19}
    %
    % 16th row
    \multirow{6}{*}{$\mathrm{db}^-$} &
    $\tau_j$ & & &%
    & & &%
    & & &%
    & & &%
    & & & $\begin{smallmatrix} \mathrm{out} \\ \mathrm{out} \end{smallmatrix}$ & $\begin{smallmatrix} \mathrm{in} \\ \mathrm{out} \end{smallmatrix}$ & $\begin{smallmatrix} \mathrm{lp} \\ \mathrm{out} \end{smallmatrix}$ &%
    \\
    \cline{2-3}\cline{5-19}
    %
    % 17th row
    & & $\tau_j$ & &%
    & & &%
    & & &%
    & & &%
    & & & $\begin{smallmatrix} \mathrm{out} \\ \mathrm{in} \end{smallmatrix}$ & $\begin{smallmatrix} \mathrm{in} \\ \mathrm{in} \end{smallmatrix}$ & $\begin{smallmatrix} \mathrm{lp} \\ \mathrm{in} \end{smallmatrix}$ &%
    \\
    \cline{2-3}\cline{5-19}
    %
    % 18th row
    & $\tau_j$ & $\tau_j$ & &%
    & & &%
    & & &%
    & & &%
    & & & $\begin{smallmatrix} \mathrm{out} \\ \mathrm{lp} \end{smallmatrix}$ & $\begin{smallmatrix} \mathrm{in} \\ \mathrm{lp} \end{smallmatrix}$ & $\begin{smallmatrix} \mathrm{lp} \\ \mathrm{lp} \end{smallmatrix}$ &%
    \\
    \cline{2-3}\cline{5-19}
    %
    % 19th row
    & $\iota_j$ & & &%
    & & &%
    & & &%
    & & &%
    $\begin{smallmatrix} \mathrm{out} \\ \mathrm{out} \end{smallmatrix}$ & $\begin{smallmatrix} \mathrm{in} \\ \mathrm{out} \end{smallmatrix}$ & $\begin{smallmatrix} \mathrm{lp} \\ \mathrm{out} \end{smallmatrix}$ & & & &%
    \\
    \cline{2-3}\cline{5-19}
    %
    % 20th row
    & & $\iota_j$ & &%
    & & &%
    & & &%
    & & &%
    $\begin{smallmatrix} \mathrm{out} \\ \mathrm{in} \end{smallmatrix}$ & $\begin{smallmatrix} \mathrm{in} \\ \mathrm{in} \end{smallmatrix}$ & $\begin{smallmatrix} \mathrm{lp} \\ \mathrm{in} \end{smallmatrix}$ & & & &%
    \\
    \cline{2-3}\cline{5-19}
    %
    % 21st row
    & $\iota_j$ & $\iota_j$ & &%
    & & &%
    & & &%
    & & &%
    $\begin{smallmatrix} \mathrm{out} \\ \mathrm{lp} \end{smallmatrix}$ & $\begin{smallmatrix} \mathrm{in} \\ \mathrm{lp} \end{smallmatrix}$ & $\begin{smallmatrix} \mathrm{lp} \\ \mathrm{lp} \end{smallmatrix}$ & & & &%
    \\
    \cline{1-20}
\end{tabular}
%\end{tabularx}
\caption{Type 2 transitions. Each entry $\begin{smallmatrix} \alpha \\ \beta \end{smallmatrix}$ represents a transition to the state $q_1 \big( \begin{smallmatrix} \mathrm{c}(a) & \alpha \\ \mathrm{c}(b) & \beta \end{smallmatrix} \big)$. An empty entry at $\big( \begin{smallmatrix} \gamma \\ \delta \end{smallmatrix} \big)$ means that the state $q_1 \big( \begin{smallmatrix} i & \gamma \\ j & \delta \end{smallmatrix} \big)$ does not exist.}
\label{tab:transition:type:2}
%\end{table}
\end{sidewaystable}

For instance, \cref{fig:F:gluing:automaton} depicts the gluing automaton $\mathrm{Gl}_\mathcal{I}$ for the interval edge replacement system from \cref{fig.interval.replacement}, which is a classical example.
\cref{fig:dendrite:gluing:automaton} depicts the gluing automaton $\mathrm{Gl}_{\mathcal{D}_n}$ for the dendrite edge replacement systems depicted in \cref{fig.dendrite.replacement}, except that the base graph is replaced by a single edge named $s$ in order to draw a simpler gluing automaton (it is easy to see that the limit spaces and rearrangement groups are the same).
Finally, \cref{tab:Basilica:gluing:automaton} represents the gluing automaton of the basilica edge replacement system:
the automaton $\mathrm{Gl}_\mathcal{B}$ features a transition $q^{(1)} \xrightarrow{(a,b)} q^{(2)}$ if and only if \cref{tab:Basilica:gluing:automaton} features a row that reads $q^{(1)} \mid (a,b) \mid q^{(2)}$.

For example, to deduce the transition $q_1 \big(\begin{smallmatrix} 1 & \mathrm{in} \\ 1 & \mathrm{out} \end{smallmatrix}\big) \xrightarrow{(1,0)} q_1 \big(\begin{smallmatrix} 1 & \mathrm{in} \\ 1 & \mathrm{out} \end{smallmatrix}\big)$ in \cref{fig:F:gluing:automaton} from \cref{tab:transition:type:2} one needs to observe that $\gamma=\mathrm{in}$ and $\tau(1)$ is the terminal vertex $\tau_1$ of the replacement graph, so one needs to consider the second column, and that $\delta=\mathrm{out}$ and $\iota(0)$ is the initial vertex $\iota_1$ of the replacement graph, so one needs to consider the fourth row.
This entry of the table reads $\begin{smallmatrix} \mathrm{in} \\ \mathrm{out} \end{smallmatrix}$.

\begin{figure}
\centering
\begin{tikzpicture}
    \node[state] (start) at (-2.5,0) {$q_0(0)$};
    \node[state] (i) at (0,0) {$q_0(1)$};
    \node[state] (in-out) at (3,1) {$q_1 \big(\begin{smallmatrix} 1 & \mathrm{in} \\ 1 & \mathrm{out} \end{smallmatrix}\big)$};
    \node[state] (out-in) at (3,-1) {$q_1 \big(\begin{smallmatrix} 1 & \mathrm{out} \\ 1 & \mathrm{in} \end{smallmatrix}\big)$};
    \draw[edge] (-3.5,0) to (start);
    \draw[edge] (start) to node[above]{$(s,s)$} (i);
    \draw[edge] (i) to[loop above, min distance=1cm] node[above]{$(0,0)$} (i);
    \draw[edge] (i) to[loop below, min distance=1cm] node[below]{$(1,1)$} (i);
    \draw[edge] (i) to node[above]{$(0,1)$} (in-out);
    \draw[edge] (i) to node[below]{$(1,0)$} (out-in);
    \draw[edge] (in-out) to[loop right, min distance=1cm] node[right]{$(1,0)$} (in-out);
    \draw[edge] (out-in) to[loop right, min distance=1cm] node[right]{$(0,1)$} (out-in);
\end{tikzpicture}
\caption{The gluing automaton $\mathrm{Gl}_\mathcal{I}$ for the interval edge replacement system.}
\label{fig:F:gluing:automaton}
\end{figure}

\begin{figure}
\centering
\begin{tikzpicture}
    \node[state] (start) at (-2.5,0) {$q_0(0)$};
    \node[state] (i) at (0,0) {$q_0(1)$};
    \node[state] (out-out) at (3,0) {$q_1 \big(\begin{smallmatrix} 1 & \mathrm{out} \\ 1 & \mathrm{out} \end{smallmatrix}\big)$};
    \node[state] (in-in) at (6,0) {$q_1 \big(\begin{smallmatrix} 1 & \mathrm{in} \\ 1 & \mathrm{in} \end{smallmatrix}\big)$};
    \draw[edge] (-3.5,0) to (start);
    \draw[edge] (start) to node[above]{$(s,s)$} (i);
    \draw[edge] (i) to[loop above, min distance=1cm] node[above]{$(i,i)$} (i);
    \draw[edge] (i) to node[above]{$(i,j)$} (out-out);
    \draw[edge] (out-out) to node[above]{$(1,1)$} (in-in);
    \draw[edge] (in-in) to[loop above, min distance=1cm] node[above]{$(n,n)$} (in-in);
\end{tikzpicture}
\caption{The gluing automaton $\mathrm{Gl}_{\mathcal{D}_n}$ for the dendrite edge replacement systems.
The label $(i,j)$ indicates that there is a transition for all possible pairs $(i,j)$, where $i$ and $j$ range from $1$ to $n$.}
\label{fig:dendrite:gluing:automaton}
\end{figure}
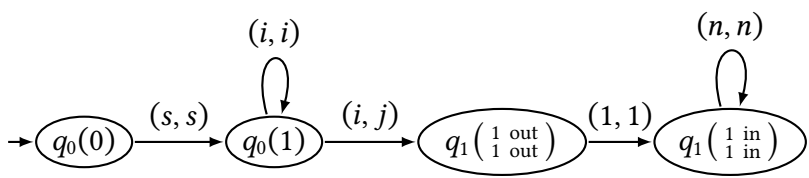

\begin{table}
\centering
\renewcommand{\arraystretch}{1.333}
\small
\begin{subtable}[T]{0.45\textwidth}
\centering
{\normalsize Types 0 and 1}\\
\vspace{.25cm}
\begin{tabular}{r|c|l}
    \multirow{4}{*}{$q_0(0)$} & $(L,L)$ & $q_0(\textcolor{red}{2})$ \\
    & $(R,R)$ & $q_0(\textcolor{red}{2})$ \\
    & $(L,R)$ & $q_1 \big(\begin{smallmatrix} \textcolor{red}{2} & \mathrm{lp} \\ \textcolor{red}{2} & \mathrm{lp} \end{smallmatrix}\big)$ \\
    & $(R,L)$ & $q_1 \big(\begin{smallmatrix} \textcolor{red}{2} & \mathrm{lp} \\ \textcolor{red}{2} & \mathrm{lp} \end{smallmatrix}\big)$ \\
    \cline{1-3}
    \multirow{9}{*}{$q_0(\textcolor{blue}{1})$} & $(1,1)$ & $q_0(\textcolor{blue}{1})$ \\
    & $(1,2)$ & $q_1 \big(\begin{smallmatrix} \textcolor{blue}{1} & \mathrm{in} \\ \textcolor{red}{2} & \mathrm{lp} \end{smallmatrix}\big)$ \\
    & $(1,3)$ & $q_1 \big(\begin{smallmatrix} \textcolor{blue}{1} & \mathrm{in} \\ \textcolor{blue}{1} & \mathrm{out} \end{smallmatrix}\big)$ \\
    & $(2,1)$ & $q_1 \big(\begin{smallmatrix} \textcolor{red}{2} & \mathrm{lp} \\ \textcolor{blue}{1} & \mathrm{in} \end{smallmatrix}\big)$ \\
    & $(2,2)$ & $q_0(\textcolor{red}{2})$ \\
    & $(2,3)$ & $q_1 \big(\begin{smallmatrix} \textcolor{red}{2} & \mathrm{lp} \\ \textcolor{blue}{1} & \mathrm{out} \end{smallmatrix}\big)$ \\
    & $(3,1)$ & $q_1 \big(\begin{smallmatrix} \textcolor{blue}{1} & \mathrm{out} \\ \textcolor{blue}{1} & \mathrm{in} \end{smallmatrix}\big)$ \\
    & $(3,2)$ & $q_1 \big(\begin{smallmatrix} \textcolor{blue}{1} & \mathrm{out} \\ \textcolor{red}{2} & \mathrm{lp} \end{smallmatrix}\big)$ \\
    & $(3,3)$ & $q_0(\textcolor{blue}{1})$ \\
    \cline{1-3}
    \multirow{9}{*}{$q_0(\textcolor{red}{2})$} & $(4,4)$ & $q_0(\textcolor{blue}{1})$ \\
    & $(4,5)$ & $q_1 \big(\begin{smallmatrix} \textcolor{blue}{1} & \mathrm{in} \\ \textcolor{red}{2} & \mathrm{lp} \end{smallmatrix}\big)$ \\
    & $(4,6)$ & $q_1 \big(\begin{smallmatrix} \textcolor{blue}{1} & \mathrm{db}^+ \\ \textcolor{blue}{1} & \mathrm{db}^- \end{smallmatrix}\big)$ \\
    & $(5,4)$ & $q_1 \big(\begin{smallmatrix} \textcolor{red}{2} & \mathrm{lp} \\ \textcolor{blue}{1} & \mathrm{in} \end{smallmatrix}\big)$ \\
    & $(5,5)$ & $q_0(\textcolor{red}{2})$ \\
    & $(5,6)$ & $q_1 \big(\begin{smallmatrix} \textcolor{red}{2} & \mathrm{lp} \\ \textcolor{blue}{1} & \mathrm{out} \end{smallmatrix}\big)$ \\
    & $(6,4)$ & $q_1 \big(\begin{smallmatrix} \textcolor{blue}{1} & \mathrm{db}^+ \\ \textcolor{blue}{1} & \mathrm{db}^- \end{smallmatrix}\big)$ \\
    & $(6,5)$ & $q_1 \big(\begin{smallmatrix} \textcolor{blue}{1} & \mathrm{out} \\ \textcolor{red}{2} & \mathrm{lp} \end{smallmatrix}\big)$ \\
    & $(6,6)$ & $q_0(\textcolor{blue}{1})$ \\
\end{tabular}
\end{subtable}
\begin{subtable}[T]{0.45\textwidth}
\centering
{\normalsize Type 2}\\
\vspace{.25cm}
\begin{tabular}{r|c|l}
    \multirow{4}{*}{$q_1 \big(\begin{smallmatrix} \textcolor{red}{2} & \mathrm{lp} \\ \textcolor{red}{2} & \mathrm{lp} \end{smallmatrix}\big)$} & $(4,4)$ & $q_1 \big(\begin{smallmatrix} \textcolor{blue}{1} & \mathrm{out} \\ \textcolor{blue}{1} & \mathrm{out} \end{smallmatrix}\big)$ \\
    & $(4,6)$ & $q_1 \big(\begin{smallmatrix} \textcolor{blue}{1} & \mathrm{out} \\ \textcolor{blue}{1} & \mathrm{in} \end{smallmatrix}\big)$ \\
    & $(6,4)$ & $q_1 \big(\begin{smallmatrix} \textcolor{blue}{1} & \mathrm{in} \\ \textcolor{blue}{1} & \mathrm{out} \end{smallmatrix}\big)$ \\
    & $(6,6)$ & $q_1 \big(\begin{smallmatrix} \textcolor{blue}{1} & \mathrm{in} \\ \textcolor{blue}{1} & \mathrm{in} \end{smallmatrix}\big)$ \\
    \cline{1-3}
    \multirow{2}{*}{$q_1 \big(\begin{smallmatrix} \textcolor{blue}{1} & \mathrm{in} \\ \textcolor{red}{2} & \mathrm{lp} \end{smallmatrix}\big)$} & $(3,4)$ & $q_1 \big(\begin{smallmatrix} \textcolor{blue}{1} & \mathrm{in} \\ \textcolor{blue}{1} & \mathrm{out} \end{smallmatrix}\big)$ \\
    & $(3,6)$ & $q_1 \big(\begin{smallmatrix} \textcolor{blue}{1} & \mathrm{in} \\ \textcolor{blue}{1} & \mathrm{in} \end{smallmatrix}\big)$ \\
    \cline{1-3}
    \multirow{1}{*}{$q_1 \big(\begin{smallmatrix} \textcolor{blue}{1} & \mathrm{in} \\ \textcolor{blue}{1} & \mathrm{out} \end{smallmatrix}\big)$} & $(3,1)$ & $q_1 \big(\begin{smallmatrix} \textcolor{blue}{1} & \mathrm{in} \\ \textcolor{blue}{1} & \mathrm{out} \end{smallmatrix}\big)$ \\
    \cline{1-3}
    \multirow{2}{*}{$q_1 \big(\begin{smallmatrix} \textcolor{red}{2} & \mathrm{lp} \\ \textcolor{blue}{1} & \mathrm{in} \end{smallmatrix}\big)$} & $(4,3)$ & $q_1 \big(\begin{smallmatrix} \textcolor{blue}{1} & \mathrm{out} \\ \textcolor{blue}{1} & \mathrm{in} \end{smallmatrix}\big)$ \\
    & $(6,3)$ & $q_1 \big(\begin{smallmatrix} \textcolor{blue}{1} & \mathrm{in} \\ \textcolor{blue}{1} & \mathrm{in} \end{smallmatrix}\big)$ \\
    \cline{1-3}
    \multirow{2}{*}{$q_1 \big(\begin{smallmatrix} \textcolor{red}{2} & \mathrm{lp} \\ \textcolor{blue}{1} & \mathrm{out} \end{smallmatrix}\big)$} & $(4,1)$ & $q_1 \big(\begin{smallmatrix} \textcolor{blue}{1} & \mathrm{out} \\ \textcolor{blue}{1} & \mathrm{out} \end{smallmatrix}\big)$ \\
    & $(6,1)$ & $q_1 \big(\begin{smallmatrix} \textcolor{blue}{1} & \mathrm{in} \\ \textcolor{blue}{1} & \mathrm{out} \end{smallmatrix}\big)$ \\
    \cline{1-3}
    \multirow{1}{*}{$q_1 \big(\begin{smallmatrix} \textcolor{blue}{1} & \mathrm{out} \\ \textcolor{blue}{1} & \mathrm{in} \end{smallmatrix}\big)$} & $(1,3)$ & $q_1 \big(\begin{smallmatrix} \textcolor{blue}{1} & \mathrm{out} \\ \textcolor{blue}{1} & \mathrm{in} \end{smallmatrix}\big)$ \\
    \cline{1-3}
    \multirow{2}{*}{$q_1 \big(\begin{smallmatrix} \textcolor{blue}{1} & \mathrm{out} \\ \textcolor{red}{2} & \mathrm{lp} \end{smallmatrix}\big)$} & $(1,4)$ & $q_1 \big(\begin{smallmatrix} \textcolor{blue}{1} & \mathrm{out} \\ \textcolor{blue}{1} & \mathrm{out} \end{smallmatrix}\big)$ \\
    & $(1,6)$ & $q_1 \big(\begin{smallmatrix} \textcolor{blue}{1} & \mathrm{out} \\ \textcolor{blue}{1} & \mathrm{in} \end{smallmatrix}\big)$ \\
    \cline{1-3}
    \multirow{1}{*}{$q_1 \big(\begin{smallmatrix} \textcolor{blue}{1} & \mathrm{out} \\ \textcolor{blue}{1} & \mathrm{out} \end{smallmatrix}\big)$} & $(1,1)$ & $q_1 \big(\begin{smallmatrix} \textcolor{blue}{1} & \mathrm{out} \\ \textcolor{blue}{1} & \mathrm{out} \end{smallmatrix}\big)$ \\
    \cline{1-3}
    \multirow{1}{*}{$q_1 \big(\begin{smallmatrix} \textcolor{blue}{1} & \mathrm{in} \\ \textcolor{blue}{1} & \mathrm{in} \end{smallmatrix}\big)$} & $(3,3)$ & $q_1 \big(\begin{smallmatrix} \textcolor{blue}{1} & \mathrm{in} \\ \textcolor{blue}{1} & \mathrm{in} \end{smallmatrix}\big)$ \\
    \cline{1-3}
    \multirow{2}{*}{$q_1 \big(\begin{smallmatrix} \textcolor{blue}{1} & \mathrm{db}^+ \\ \textcolor{blue}{1} & \mathrm{db}^- \end{smallmatrix}\big)$} & $(1,3)$ & $q_1 \big(\begin{smallmatrix} \textcolor{blue}{1} & \mathrm{out} \\ \textcolor{blue}{1} & \mathrm{in} \end{smallmatrix}\big)$ \\
    & $(3,1)$ & $q_1 \big(\begin{smallmatrix} \textcolor{blue}{1} & \mathrm{in} \\ \textcolor{blue}{1} & \mathrm{out} \end{smallmatrix}\big)$ \\
\end{tabular}
\end{subtable}
\caption{The gluing automaton $\mathrm{Gl}_\mathcal{B}$ for the basilica.}
\label{tab:Basilica:gluing:automaton}
\end{table}

\section{Proof of the Theorem}
\label{sec:gl:aut}

The definition of the gluing automaton $\mathrm{Gl}_\mathcal{R}$ might look complicated, but each of its transition has a simple meaning.
The aim of this subsection is to describe how the gluing automaton works with a series of Lemmas that will be useful when proving the rationality of the gluing relation in \cref{thm:rat:gl}.

\medskip %layout
\begin{lemma}[Type 0 transitions generate the color graph]
\label{lem:0:trans}
The full sub-automaton $\mathrm{T}_\mathcal{R} = \left( \Sigma^2, Q_0, \rightarrow, q_0 \right)$ induced by the states $Q_0$ is naturally isomorphic to the color graph.
More precisely, a finite word $(x_0, y_0) \dots (x_m, y_m)$ in the alphabet $\Sigma^2$ is recognized by $\mathrm{T}_\mathcal{R}$ if and only if $x_i = y_i$ and the word $x_0 \dots x_m$ belongs to $\mathbb{A}_\mathcal{R}$.
\end{lemma}

\begin{proof}
Note that $\mathrm{Gl}_\mathcal{R}$ features no transition from $Q_1$ to $Q_0$ and the transitions from $Q_0$ to itself are precisely the Type 0 transitions.
The map $\mathrm{T}_\mathcal{R} \rightarrow A_\mathcal{R}$ that is the identity on $Q_0$ and sends the transition $q_0(i) \xrightarrow{(e,e)} q_0 (\mathrm{c}(e))$ to $q_0(i) \xrightarrow{e} q_0 (\mathrm{c}(e))$ is an isomorphism.
\end{proof}

Note that if two words $x_0 \dots x_m$ and $y_0 \dots y_m$ in $\mathbb{A}_\mathcal{R}$ are the same except for $x_m \neq y_m$, then they represent adjacent edges in the $m$-the full expansion $E_m$ if and only if $x_m$ and $y_m$ are adjacent in $\Gamma_{\mathrm{C}(x_{m-1})} = \Gamma_{\mathrm{C}(y_{m-1})}$.
Thus, even if by definition Type 1 transitions depend on edge adjacencies in the base or the replacement graphs, they really determine adjacencies in $E_m$.
This is because, using distinguished colors for loops and defining their edge replacement systems to have a sole boundary vertex (as described in \cref{ass.loops}) allows us to see the expanded subgraph as an embedded copy of the replacement graph.

\medskip %layout
\begin{lemma}[Type 1 transition implies adjacency]
\label{lem:1:trans:imply:adj}
Given a sequence $(x_0, y_0) (x_1, y_1) \dots$ of $( \Sigma^2)^\omega$ that is recognized by $\mathrm{Gl}_\mathcal{R}$, one and only one of the following hold.
\begin{enumerate}
    \item The sequence is recognized by the sub-automaton $\mathrm{T}_\mathcal{R}$ and $x_i=y_i$ for all $i$.
    \item The run of $(x_0, y_0) (x_1, y_1) \dots$ in $\mathrm{Gl}_\mathcal{R}$ features exactly one Type 1 transition taking place when processing the first $(x_m, y_m)$ such that $x_m \neq y_m$.
    Moreover, $x_0 \dots x_m$ and $y_0 \dots y_m$ are adjacent edges in $E_m$ and the state of $\mathrm{Gl}_\mathcal{R}$ after processing $(x_m, y_m)$ is $q_1 \big(\begin{smallmatrix} \mathrm{c}(x_m) & \alpha \\ \mathrm{c}(y_m) & \beta \end{smallmatrix}\big)$, which determines the type of adjacency by the entries $\alpha$ and $\beta$, as describe in \cref{tab:transition:type:1}.
\end{enumerate}
\end{lemma}

\begin{proof}
Case (1) is the one discussed in \cref{lem:0:trans}, so let us assume that there is an $m \in \mathbb{N}$ such that $x_m \neq y_m$.
Take $m$ to be the minimum such number.
Then, before the letter $(x_m, y_m)$, only Type 0 transitions have taken place.
Note that $\mathrm{Gl}_\mathcal{R}$ features no transitions from $Q_1$ to $Q_0$.
Thus, since Type 1 transitions are all oriented from $Q_0$ to $Q_1$, each run can feature at most one Type 1 transition.

By \cref{lem:0:trans}, $x_0 \dots x_{m-1} = y_0 \dots y_{m-1}$ belongs to $\mathbb{A}_\mathcal{R}$, so it is an edge of $E_{m-1}$.
By the definition of a Type 1 transition, it follows that $x_0 \dots x_m$ and $y_0 \dots y_m$ are adjacent in $E_m$ and that the state of $\mathrm{Gl}_\mathcal{R}$ after processing $(x_m, y_m)$ is $q_1 \big(\begin{smallmatrix} \mathrm{c}(x_m) & \alpha \\ \mathrm{c}(y_m) & \beta \end{smallmatrix}\big)$, where $\alpha$ and $\beta$ describe the type of adjacency as expressed in \cref{tab:transition:type:1}.
\end{proof}

\medskip %layout
\begin{lemma}[Adjacency implies Type 1 transition]
\label{lem:adj:imply:1:trans}
Given two words $x_0 \dots x_m$ and $y_0 \dots y_m$ in $\mathbb{A}_\mathcal{R}$ that are the same except for $x_m \neq y_m$, if they represent adjacent edges of $E_m$ then the word $(x_0, y_0) (x_1, y_1) \dots (x_m, y_m)$ is recognized by $\mathrm{Gl}_\mathcal{R}$, and its run consists of $m$ Type 0 transitions followed by a Type 1 transition.
\\
Moreover, the type of adjacency between $x_m$ and $y_m$ in $\Gamma_{\mathrm{C}(x_{m-1})} = \Gamma_{\mathrm{C}(y_{m-1})}$ and between $x_0 \dots x_m$ and $y_0 \dots y_m$ in $E_m$ is the same, i.e., the last state of the run $q_1 \big(\begin{smallmatrix} \mathrm{c}(x_m) & \alpha \\ \mathrm{c}(y_m) & \beta \end{smallmatrix}\big)$
determines the edge adjacency of $x_0 \dots x_m$ and $y_0 \dots y_m$ in $E_m$ as described in \cref{tab:transition:type:1}.
\end{lemma}

\begin{proof}
Given two words $x_0 \dots x_m$ and $y_0 \dots y_m$ in $\mathbb{A}_\mathcal{R}$ that only differ in $x_m \neq y_m$ and that represent adjacent edges of $E_m$, by \cref{lem:0:trans} the sub-automaton $\mathrm{T}_\mathcal{R}$ recognizes $(x_0, y_0) (x_1, y_1) \ldots (x_{m-1}, y_{m-1})$.
By definition, a Type 1 transition occurs when processing the letter $(x_m, y_m)$, ending in a state $q_1 \big(\begin{smallmatrix} \mathrm{c}(x_m) & \alpha \\ \mathrm{c}(y_m) & \beta \end{smallmatrix}\big)$ that determines the edge adjacency between $x_m$ and $y_m$ in the replacement graph $\Gamma_{\mathrm{c}(x_{m-1})} = \Gamma_{\mathrm{c}(y_{m-1})}$, and thus in $E_m$, as expressed in \cref{tab:transition:type:1}.
\end{proof}

Note that a Type 2 transition is preceded by a (possibly empty) sequence of Type 0 transitions, one Type 1 transition and possibly a sequence of Type 2 transitions.

\medskip %layout
\begin{lemma}[Type 2 transitions imply adjacency]
\label{lem:2:trans:imply:adj}
Let $(x_0, y_0) (x_1, y_1) \ldots (x_k, y_k)$ be recognized by $\mathrm{Gl}_\mathcal{R}$.
If the transition involving $(x_k, y_k)$ is a Type 2 transition, then $x_0 \dots x_k$ and $y_0 \dots y_k$ are adjacent edges of $E_k$.
\end{lemma}

\begin{proof}
First, suppose that the transition involving $(x_k, y_k)$ is the first Type 2 transition, so the previous transition is of Type 1.
Then, by \cref{lem:1:trans:imply:adj}, the edges $x_1 \dots x_{k-1}$ and $y_1 \dots y_{k-1}$ are adjacent and the state $q_1 \big(\begin{smallmatrix} \mathrm{c}(x_{k-1}) & \alpha \\ \mathrm{c}(y_{k-1}) & \beta \end{smallmatrix}\big)$ after processing $(x_{k-1}, y_{k-1})$ describes the type of adjacency as expressed in \cref{tab:transition:type:1}.
Thoroughly going through \cref{tab:transition:type:2} 
reveals that the column and the row of any non empty entry correspond to a pair of adjacent edges in the graph obtained by gluing the replacement graphs $\Gamma_{\mathrm{c}(x_{k-1})}$ and $\Gamma_{\mathrm{c}(y_{k-1})}$ as dictated by $\alpha$ and $\beta$, which appears as a subgraph of $E_k$ that contains the edges $x_1 \dots x_k$ and $y_1 \dots y_k$.
Moreover, the resulting state $q_1 \big(\begin{smallmatrix} \mathrm{c}(x_k) & \alpha' \\ \mathrm{c}(y_k) & \beta' \end{smallmatrix}\big)$ describes the adjacency between $x_0 \dots x_k$ and $y_0 \dots y_k$.

Now suppose that the transition involving $(x_k, y_k)$ is not the first Type 2 transition.
The argument of the previous paragraph applies as soon as the last state in the run of $(x_1,y_1) \dots (x_{k-1}, y_{k-1})$ described the adjacency of $(x_{k-1}, y_{k-1})$.
Thus, we can iteratively apply the argument to see that $x_0 \dots x_k$ and $y_0 \dots y_k$ are adjacent and that the last state of the run determines their type of adjacency.
\end{proof}

\medskip %layout
\begin{lemma}[Adjacency implies a Type 2 transition]
\label{lem:adj:imply:2:trans}
Let $x_0 \dots x_{m+1}$ and $y_0 \dots y_{m+1}$ in $\mathbb{A}_\mathcal{R}$ represent adjacent edges of $E_{m+1}$.
Suppose that $x_0 \dots x_m \neq y_0 \dots y_m$ and that $(x_0, y_0) \dots (x_m, y_m)$ is recognized by $\mathrm{Gl}_\mathcal{R}$.
Then $(x_0, y_0) \dots (x_{m+1}, y_{m+1})$ is recognized by $\mathrm{Gl}_\mathcal{R}$ and the last transition is a Type 2 transition ending in $q_1 \big(\begin{smallmatrix} \mathrm{c}(x_{m+1}) & \alpha \\ \mathrm{c}(y_{m+1}) & \beta \end{smallmatrix}\big)$ that describes the adjacency between $x_0 \dots x_{m+1}$ and $y_0 \dots y_{m+1}$ in $E_{m+1}$.
\end{lemma}

\begin{proof}
First, suppose that the state after having processed $(x_0, y_0), \dots (x_m, y_m)$ is $q_1 \big(\begin{smallmatrix} \mathrm{c}(x_m) & \alpha \\ \mathrm{c}(y_m) & \beta \end{smallmatrix}\big)$ and that it describes the type of adjacency as expressed in \cref{tab:transition:type:1}.
Then a careful look at \cref{tab:transition:type:2} 
shows that the gluing automaton can process the next letter $(x_{m+1}, y_{m+1})$ with a Type 2 transition that ends in the state $q_1 \big(\begin{smallmatrix} \mathrm{c}(x_{m+1}) & \alpha' \\ \mathrm{c}(y_{m+1}) & \beta' \end{smallmatrix}\big)$ which describes the adjacency between $x_0 \dots x_{m+1}$ and $y_0 \dots y_{m+1}$.

Now, by \cref{lem:adj:imply:1:trans} the common prefix $x_0 \dots x_p = y_0 \dots y_p$ of $x_0 \dots x_m$ and $y_0 \dots y_m$ (where $0 \leq p \leq m-1$) is processed by the sub-automaton $\mathrm{T}_\mathcal{R}$, after which a Type 1 transition brings us to the state $q_1 \big(\begin{smallmatrix} \mathrm{c}(x_{p+1}) & \alpha \\ \mathrm{c}(y_{p+1}) & \beta \end{smallmatrix}\big)$ that describes the adjacency of $x_0 \dots x_{p+1}$ and $y_0 \dots y_{p+1}$ in $E_{p+1}$.
Reasoning as in the previous paragraph shows that this is followed by a Type 2 transition to the state $q_1 \big(\begin{smallmatrix} \mathrm{c}(x_{p+2}) & \alpha' \\ \mathrm{c}(y_{p+2}) & \beta' \end{smallmatrix}\big)$ that describes the adjacency between $x_0 \dots x_{p+2}$ and $y_0 \dots y_{p+2}$.
Since each new transition results in a state that describes edge adjacency, the previous paragraph can be iteratively applied to every subsequent prefix of $x_0 \dots x_{m+1}$ and $y_0 \dots y_{m+1}$, concluding with a Type 2 transition from $(x_0, y_0) \dots (x_m, y_m)$ to $(x_0, y_0) \dots (x_{m+1}, y_{m+1})$ that ends in a state that describes the edge adjacency between $x_0 \dots x_{m+1}$ and $y_0 \dots y_{m+1}$.
\end{proof}

Now let us put the Lemmas together.

\medskip %layout
\begin{theorem}
\label{thm:rat:gl}
The gluing relation of an expanding edge replacement system is rational.
\end{theorem}

\begin{proof}
We have to prove that $(x_0, y_0) (x_1, y_1) \ldots$ is recognized by $\mathrm{Gl}_\mathcal{R}$ if and only if the sequences $x_0 x_1 \ldots$ and $y_0 y_1 \ldots$ are equivalent under the gluing relation $\sim$.

\medskip

Suppose $v = x_0 x_1 \dots$ and $w = y_0 y_1 \dots$ are equivalent.
If $v = w$ then $\mathrm{Gl}_\mathcal{R}$ recognizes the sequence $(v,w)$ by \cref{lem:0:trans}, so we are done.
Suppose that instead $v \neq w$ and let $m \geq 0$ be the smallest natural number such that $x_m \neq y_m$.
Since $v \sim w$, by \cref{def.glue} there exists a unique vertex on which $x_0 \ldots x_k$ and $y_0 \ldots y_k$ are incident as edges of $E_k$, for all $k > m$ (for $k=m$ the edges may be parallel, in which case there may be two such vertices).
In particular, $x_0 \ldots x_m$ and $y_0 \ldots y_m$ are adjacent, so by \cref{lem:adj:imply:1:trans} the prefix $(x_0, y_0) \ldots (x_m, y_m)$ is recognized by $\mathrm{Gl}_\mathcal{R}$ and its run ends with a Type 1 transition to the state $q_1 \big(\begin{smallmatrix} \mathrm{c}(x_m) & \alpha \\ \mathrm{c}(y_m) & \beta \end{smallmatrix}\big)$, which describes the type of adjacency between $x_1 \dots x_m$ and $y_0 \dots y_m$ in $E_m$.

For all of the subsequent prefixes of our sequence $(v,w)$ we are in the case described in \cref{lem:adj:imply:2:trans}, which shows that there is a Type 2 transition from any of those prefixes to the next one, ultimately showing that $(v,w)$ is recognized by $\mathrm{Gl}_\mathcal{R}$.

\medskip

Conversely, suppose that $(v, w) = (x_0, y_0) (x_1, y_1) \ldots$ is recognized by $\mathrm{Gl}_\mathcal{R}$.
We need to show that $v$ and $w$ belong to the symbol space $\Omega_\mathcal{R}$ and that $ v \sim w $.

Consider the morphism $\phi$ from the underlying graph of $\mathrm{Gl}_\mathcal{R}$ to the color graph $A_\mathcal{R}$ defined by mapping the states as follows:
    \[ q_0(i) \text{ is mapped to } q_0(i), \]
    \[ q_1 \big(\begin{smallmatrix} \mathrm{c}(a) & \alpha \\ \mathrm{c}(b) & \beta \end{smallmatrix}\big) \text{ is mapped to } q_0(\mathrm{c}(a)); \]
and by mapping the transitions as follows:
    \[ q_0(i) \xrightarrow{(e,e)} q_0 (\mathrm{c}(e)) \text{ is mapped to } q_0(i) \rightarrow q_0 (\mathrm{c}(e)), \]
    \[ q_0(i) \xrightarrow{(a,b)} q_1 \big(\begin{smallmatrix} \mathrm{c}(a) & \alpha \\ \mathrm{c}(b) & \beta \end{smallmatrix}\big) \text{ is mapped to } q_0(i) \rightarrow q_0 (\mathrm{c}(a)), \]
    \[ q_1 \big(\begin{smallmatrix} i & \gamma \\ j & \delta \end{smallmatrix}\big) \xrightarrow{(a,b)} q_1 \big(\begin{smallmatrix} \mathrm{c}(a) & \alpha \\ \mathrm{c}(b) & \beta \end{smallmatrix}\big) \text{ is mapped to } q_0(i) \rightarrow q_0 (\mathrm{c}(a)). \]
By definition of Type 1 and 2 transitions, in the schemes above $a$ is an edge of the graph $\Gamma_i$.
Thus, $\phi(\mathrm{Gl}_\mathcal{R})$ is a subgraph of the color graph, so $\phi: \mathrm{Gl}_\mathcal{R} \to A_\mathcal{R}$ is a well-defined graph morphism.
In fact, it is easy to see that $\phi(\mathrm{Gl}_\mathcal{R})$ is precisely $A_\mathcal{R}$.
Thus, if $(v,w)$ is recognized $\mathrm{Gl}_\mathcal{R}$ then $v$ is in the edge shift $\Omega_\mathcal{R}$ of the color graph.
The same holds for the entry $w$, since $\mathrm{Gl}_\mathcal{R}$ is symmetric under switching inputs $(a,b)$ with $(b,a)$, states $q_i \big(\begin{smallmatrix} i & \gamma \\ j & \delta \end{smallmatrix}\big)$ with $q_i \big(\begin{smallmatrix} j & \delta \\ i & \gamma \end{smallmatrix}\big)$ and symbols $\mathrm{db}^+$ and $\mathrm{db}^-$ .

We now show that $v \sim w$.
If $v = w$ there is nothing to prove, so let $v \neq w$.
By \cref{lem:1:trans:imply:adj}, the run of $(v, w)$ features finitely many Type 0 transitions followed by a Type 1 transition at $(x_0, y_0), \dots (x_m, y_m)$, for some $m \geq 0$, that brings us to a state of $Q_1$, and $x_0 \dots x_m$ is adjacent to $y_0 \dots y_m$.
This shows that the prefixes of $v$ and $w$ of equal length $k \leq m+1$ represent adjacent edges.
Since the transitions of the gluing automaton that start in $Q_1$ are of Type 2 and end in $Q_1$, every subsequent transition of the run of $(v, w)$ is a Type 2 transition.
By \cref{lem:2:trans:imply:adj}, the prefixes of $v$ and $w$ of equal length $k > m+1$ are adjacent, concluding our proof.
\end{proof}

%%%%%%%%%%%%%%%%%%%%%%%%%

\chapter[Additional Results]{Additional Results about Rearrangement Groups}
\label{cha.additional.results}

In this chapter we exhibit three main additional results, one for each section.
The first two \cref{sec.f.g.abelian,sec.rational.stabilizers} describe some closure properties of the family of rearrangement groups which go to show that the class is very large.
These results are not included in any publication nor preprint at the time of writing this dissertation.
The third \cref{sec.embedding.into.V} describe embeddings of rearrangement groups into Thompson's group $V$, which was first noted in \cite{conjugacy} but is described in greater detail here.

\section[Finitely Generated Abelian Groups]{Finitely Generated Abelian Groups are Rearrangement Groups}
\label{sec.f.g.abelian}

In this section we exhibit some closure properties of the class of rearrangement groups and we show that finitely generated abelian groups can be realized as rearrangement groups.
We will first show that direct products of rearrangement groups are rearrangement groups, then we will see how to construct edge replacement systems for any finite group and finally for the infinite cyclic group $\mathbb{Z}$.

\subsection{Products of Rearrangement Groups}
\label{sub.products}

Suppose that $G_X$ and $G_Y$ are rearrangement groups for two edge replacement systems $\mathcal{X}$ and $\mathcal{S}$, respectively.
Without loss of generality, assume that $\mathcal{R}$ and $\mathcal{S}$ have distinct colors.
Consider the edge replacement system $\mathcal{R} \oplus \mathcal{S}$ whose base graph is the disjoint union of those of $\mathcal{R}$ and $\mathcal{S}$ and whose replacement graphs are all of those from $\mathcal{R}$ and $\mathcal{S}$.
Since graph isomorphisms preserve the coloring of the edges and the colors of $\mathcal{R}$ and $\mathcal{S}$ are distinct, graph pair diagrams for the edge replacement system $\mathcal{R} \oplus \mathcal{S}$ can be decomposed into an element of $G_X$ and one of $G_Y$.
It is clear that the two components of each graph pair diagram do not interact with each other, so we have the following fact.

\medskip %layout
\begin{proposition}
\label{prop.direct.product.is.RG}
If $G_X$ and $G_Y$ are rearrangement groups, then $G_X \times G_Y$ can also be realized as a rearrangement group.
\end{proposition}

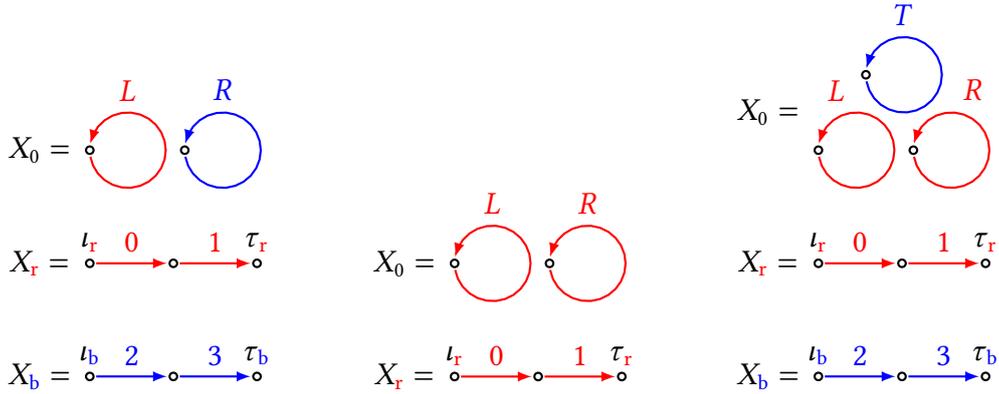
\begin{figure}
\centering
\begin{subfigure}[t]{.32\textwidth}
\centering
\begin{tikzpicture}
    \node at (-.667,0) {$X_0 =$};
    \node[vertex] (L) at (0,0) {};
    \node[vertex] (R) at (1.25,0) {};
    \draw[edge,red,domain=190:530] plot ({.5*cos(\x)+.5}, {.5*sin(\x)});
    \draw[edge,blue,domain=190:530] plot ({.5*cos(\x)+1.75}, {.5*sin(\x)});
    \draw (.5,.5) node[above,red] {$L$};
    \draw (1.75,.5) node[above,blue] {$R$};
    \begin{scope}[yshift=-1.5cm]
    \node at (-.667,0) {$X_{\text{\textcolor{red}{r}}} = $};
    \node[vertex] (l) at (0,0) {};
    \draw (l) node[above]{$\iota_{\text{\textcolor{red}{r}}}$};
    \node[vertex] (c) at (1.1,0) {};
    \node[vertex] (r) at (2.2,0) {};
    \draw (r) node[above]{$\tau_{\text{\textcolor{red}{r}}}$};
    \draw[edge,red] (l) to node[above]{$0$} (c);
    \draw[edge,red] (c) to node[above]{$1$} (r);
    \end{scope}
    \begin{scope}[yshift=-3cm]
    \node at (-.667,0) {$X_{\text{\textcolor{blue}{b}}} = $};
    \node[vertex] (l) at (0,0) {};
    \draw (l) node[above]{$\iota_{\text{\textcolor{blue}{b}}}$};
    \node[vertex] (c) at (1.1,0) {};
    \node[vertex] (r) at (2.2,0) {};
    \draw (r) node[above]{$\tau_{\text{\textcolor{blue}{b}}}$};
    \draw[edge,blue] (l) to node[above]{$2$} (c);
    \draw[edge,blue] (c) to node[above]{$3$} (r);
    \end{scope}
\end{tikzpicture}
\caption{An edge replacement system for $T^2$.}
\label{fig.direct.T}
\end{subfigure}
\begin{subfigure}[t]{.32\textwidth}
\centering
\begin{tikzpicture}
    \node at (-.667,0) {$X_0 =$};
    \node[vertex] (L) at (0,0) {};
    \node[vertex] (R) at (1.25,0) {};
    \draw[edge,red,domain=190:530] plot ({.5*cos(\x)+.5}, {.5*sin(\x)});
    \draw[edge,red,domain=190:530] plot ({.5*cos(\x)+1.75}, {.5*sin(\x)});
    \draw (.5,.5) node[above,red] {$L$};
    \draw (1.75,.5) node[above,red] {$R$};
    \begin{scope}[yshift=-1.5cm]
    \node at (-.667,0) {$X_{\text{\textcolor{red}{r}}} = $};
    \node[vertex] (l) at (0,0) {};
    \draw (l) node[above]{$\iota_{\text{\textcolor{red}{r}}}$};
    \node[vertex] (c) at (1.1,0) {};
    \node[vertex] (r) at (2.2,0) {};
    \draw (r) node[above]{$\tau_{\text{\textcolor{red}{r}}}$};
    \draw[edge,red] (l) to node[above]{$0$} (c);
    \draw[edge,red] (c) to node[above]{$1$} (r);
    \end{scope}
\end{tikzpicture}
\caption{An edge replacement system for $T^2 \rtimes \mathrm{Sym}(2)$.}
\label{fig.complete.wreath.T}
\end{subfigure}
\begin{subfigure}[t]{.32\textwidth}
\centering
\begin{tikzpicture}
    \node at (-.667,.5) {$X_0 =$};
    \node[vertex] (L) at (0,0) {};
    \node[vertex] (R) at (1.25,0) {};
    \node[vertex] (T) at (.625,1) {};
    \draw[edge,red,domain=190:530] plot ({.5*cos(\x)+.5}, {.5*sin(\x)});
    \draw[edge,red,domain=190:530] plot ({.5*cos(\x)+1.75}, {.5*sin(\x)});
    \draw[edge,blue,domain=190:530] plot ({.5*cos(\x)+1.125}, {.5*sin(\x)+1});
    \draw (.5,.5) node[above left,red] {$L$};
    \draw (1.75,.5) node[above right,red] {$R$};
    \draw (1.125,1.5) node[above,blue] {$T$};
    \begin{scope}[yshift=-1.5cm]
    \node at (-.667,0) {$X_{\text{\textcolor{red}{r}}} = $};
    \node[vertex] (l) at (0,0) {};
    \draw (l) node[above]{$\iota_{\text{\textcolor{red}{r}}}$};
    \node[vertex] (c) at (1.1,0) {};
    \node[vertex] (r) at (2.2,0) {};
    \draw (r) node[above]{$\tau_{\text{\textcolor{red}{r}}}$};
    \draw[edge,red] (l) to node[above]{$0$} (c);
    \draw[edge,red] (c) to node[above]{$1$} (r);
    \end{scope}
    \begin{scope}[yshift=-3cm]
    \node at (-.667,0) {$X_{\text{\textcolor{blue}{b}}} = $};
    \node[vertex] (l) at (0,0) {};
    \draw (l) node[above]{$\iota_{\text{\textcolor{blue}{b}}}$};
    \node[vertex] (c) at (1.1,0) {};
    \node[vertex] (r) at (2.2,0) {};
    \draw (r) node[above]{$\tau_{\text{\textcolor{blue}{b}}}$};
    \draw[edge,blue] (l) to node[above]{$2$} (c);
    \draw[edge,blue] (c) to node[above]{$3$} (r);
    \end{scope}
\end{tikzpicture}
\caption{An edge replacement system for  $T^3 \rtimes \mathrm{Sym}(2)$.}
\label{fig.incomplete.wreath.T}
\end{subfigure}
\caption{Three edge replacement systems obtained as sums of other edge replacement systems.}
\label{fig.product.T}
\end{figure}

For example, the edge replacement system depicted in \cref{fig.direct.T} produces the direct product of two copies of Thompson's group $T$, each of which acts solely on one of the cells $\llbracket L \rrbracket$ or $\llbracket R \rrbracket$.

Even if we will not use it, this kind of construction can be generalized to the case in which certain colors and their replacement graphs are common between $\mathcal{R}$ and $\mathcal{S}$, which yields certain wreath products of direct products $G_1 \times \dots \times G_k$ of rearrangement groups with finite permutation groups $S \leq \mathrm{Sym}(k)$.
For example, the edge replacement system depicted in \cref{fig.complete.wreath.T} is given by the disjoint union of the circle edge replacement system (\cref{fig.circle.replacement}) with the same color and its rearrangement group is $(T \times T) \rtimes \mathrm{Sym}(2)$, where each copy of $T$ acts solely on one of the cells $\llbracket L \rrbracket$ or $\llbracket R \rrbracket$ and the group $\mathrm{Sym}(2)$ acts on the limit space by the prefix-exchange map $\llbracket L \alpha \rrbracket \mapsto \llbracket R \alpha \rrbracket$ that permutes the two copies of $T$.
This is the same edge replacement system that was depicted in \cref{fig.double.T} to give an example of a minimal but not weakly cell-transitive rearrangement group.

Sometimes the wreath product is more complicated:
the rearrangement group of the edge replacement system depicted in \cref{fig.incomplete.wreath.T} is $(T \times T \times T) \rtimes \mathrm{Sym}(2)$, where $\mathrm{Sym}(2)$ only permutes two of the copies of $T$.

\subsection{Finite Groups as Rearrangement Groups}

We can show that every finite group is a rearrangement group using the classical Frucht's Theorem from \cite{Frucht}.

\medskip %layout
\begin{theorem}[Frucht's Theorem]
Every finite group is the automorphism group of some finite undirected graph.
\end{theorem}

It is worth mentioning that this result has been proved for infinite groups \cite{FruchtInfinite}.
Due to the fact that graph expansions are always finite graphs, we will only be able to apply Frucht's original statement about finite groups, which allows us to prove the following fact.

\medskip %layout
\begin{proposition}
\label{prop.finite.groups.are.RG}
Every finite group is a rearrangement group.
\end{proposition}

\begin{proof}
It is easy to go from undirected to directed graphs.
Given an undirected graph, one can build a directed barycentric subdivision by replacing each undirected edge $e = \{u,w\}$ with a vertex $v_e$ and two directed edges, one from $v_e$ to $u$ and one from $v_e$ to $w$.
Clearly the automorphism groups of the two graphs are the same, so Frucht's Theorem implies that every finite group is the automorphism group of some finite directed graph.

Now, given a finite group $G$, let $\Gamma$ be a finite directed graph such that $\mathrm{Aut}(\Gamma) = G$.
Consider the monochromatic edge replacement system with trivial replacement graph and $\Gamma$ as its base graph.
This is not expanding and the limit space is likely not defined, but the rearrangement group exists nonetheless and it is clearly $\mathrm{Aut}(\Gamma) = G$.
\end{proof}

The same trick used in \cref{prop.null.expanding.isolated.rarrangement} to produce an expanding edge replacement system with the same rearrangement group as a non-expanding one arguably also applies to this construction.

\subsection{The Infinite Cyclic Group}

Using the last two results we can show that the family of rearrangement groups includes every finitely generated abelian group.

\medskip %layout
\begin{theorem}
Every finitely generated abelian group is a rearrangement group.
\end{theorem}

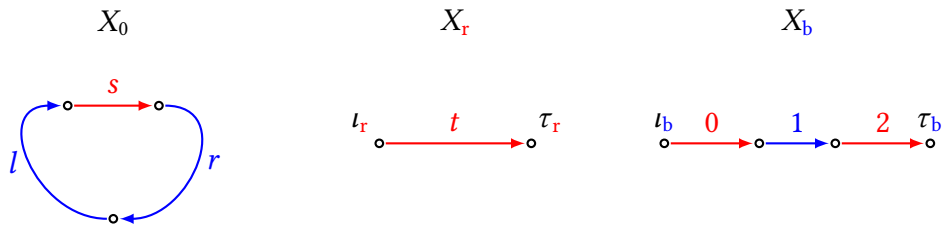
\begin{figure}
\centering
\begin{tikzpicture}
    \node at (0,1.6) {$X_0$};
    \node[vertex] (l) at (-.6,.5) {};
    \node[vertex] (r) at (.6,.5) {};
    \node[vertex] (b) at (0,-1) {};
    \draw[edge,red] (l) to node[above]{$s$} (r);
    \draw[edge,blue] (r) to[out=0,in=0,looseness=1.5] node[right]{$r$} (b);
    \draw[edge,blue] (b) to[out=180,in=180,looseness=1.5] node[left]{$l$} (l);
    \begin{scope}[xshift=4.5cm]
    \node at (0,1.6) {$X_{\text{\textcolor{red}{r}}}$};
    \node[vertex] (i) at (-1,0) {}; \draw (-1.25,0) node[above]{$\iota_{\text{\textcolor{red}{r}}}$};
    \node[vertex] (t) at (1,0) {}; \draw (1.25,0) node[above]{$\tau_{\text{\textcolor{red}{r}}}$};
    \draw[edge,red] (i) to node[above]{$t$} (t);
    \end{scope}
    \begin{scope}[xshift=9cm]
    \node at (0,1.6) {$X_{\text{\textcolor{blue}{b}}}$};
    \node[vertex] (i) at (-1.75,0) {}; \draw (-1.75,0) node[above]{$\iota_{\text{\textcolor{blue}{b}}}$};
    \node[vertex] (l) at (-.5,0) {};
    \node[vertex] (r) at (.5,0) {};
    \node[vertex] (t) at (1.75,0) {}; \draw (1.75,0) node[above]{$\tau_{\text{\textcolor{blue}{b}}}$};
    \draw[edge,red] (i) to node[above]{$0$} (l);
    \draw[edge,blue] (l) to node[above]{$1$} (r);
    \draw[edge,red] (r) to node[above]{$2$} (t);
    \end{scope}
\end{tikzpicture}
\caption{An edge replacement systems for the infinite cyclic group.}
\label{fig.Z.replacement}
\end{figure}

\begin{proof}
It is known that every finitely generated abelian group $G$ can be decomposed as $G_T \times \mathbb{Z}^k$, where $G_T$ is the torsion subgroup of $G$ and $k \in \mathbb{N}$ is the rank of $G$.
$G_T$ is a rearrangement group because of \cref{prop.finite.groups.are.RG}, so by \cref{prop.direct.product.is.RG} we only need to show that $\mathbb{Z}$ is a rearrangement group to conclude the proof.

Consider the edge replacement system depicted in \cref{fig.Z.replacement}.
Every graph expansion is a cycle that has $n+1$ consecutive red edges (those that have $s$ as a prefix) followed by a blue edge, then $n$ red edges and finally another blue edge.
The number $n$ coincides with the total number of expansions of blue edges.
Moreover, the blue edges are $r1^a$ and $l1^b$ for some $a$ and $b$ such that $a+b = n$.
Let us denote by $E(a,b)$ the graph expansions.

Two graphs $E(a,b)$ and $E(c,d)$ are isomorphic if and only if $a+b = c+d$, so each graph pair diagram is of the form $\left( E(a,b), E(c,d) \right)$ (where the graph isomorphism is uniquely determined by the graphs, so it is omitted).
Note that the graph pair diagram can be reduced precisely when both $a$ and $c$ are non-zero or both $b$ and $d$ are non-zero,
and each reduction decreases both $a$ and $c$ or $b$ and $d$ by $1$.
For example, $\left( E(1,2), E(3,0) \right) = \left( E(0,2), E(2,0) \right)$.
Then every reduced pair diagram has either of the forms $\left( E(0,n), E(n,0) \right)$ or $\left( E(n,0), E(0,n) \right)$.

Since rearrangements are uniquely determined by their reduced graph pair diagrams, there is precisely a rearrangement for each integer $z$, where $z = c-a = b-d$.
It is easy to see that composing two rearrangements that correspond to integers $z_1$ and $z_2$ yields the rearrangement that corresponds to their sum $z_1 + z_2$, so $G \simeq \mathbb{Z}$ and the proof is complete.
\end{proof}

\section{Stabilizers of Rational Points}
\label{sec.rational.stabilizers}

Fix an edge replacement system $\mathcal{R}$ and its rearrangement group $G$.
In this section we describe a simple method to build the pointwise stabilizer of any finite set $S$ of rational points (\cref{def.rational.irrational.points}) under the action of $G$ as a rearrangement group $G_S$ for a new edge replacement system $\mathcal{R}^S$ that is a colored variation of the original one.
The new edge replacement system $\mathcal{R}^S$ produces the same limit space as $\mathcal{R}$ and the action of $G_S$ is compatible with that of $G$.

The construction presented here was originally developed together with Davide Perego when working on \cite{rationalgluing} (which is adapted into \cref{cha.rational.gluing}), as it shares some of the same ideas.
The authors then decided not to include this in the final version of \cite{rationalgluing} for the sake of conciseness.

We will only describe how to build an edge replacement system for the stabilizer of a single rational point.
Passing to a finite set $S$ of rational points simply involves a finite iteration of the same construction for each point of $S$.
Ultimately, we prove the following fact.

\medskip %layout
\begin{theorem}
Given a set $S$ of rational points of the limit space of an edge replacement system $\mathcal{R}$, the pointwise stabilizer of $S$ under the action of the rearrangement group $G_\mathcal{R}$ is itself a rearrangement group.
\end{theorem}

Since the edge replacement system that we will construct for the stabilizers has the same graphs of the original one except for being re-colored, when we consider edge replacement systems whose gluing relation is trivial (\cref{sub.topological.full.groups}) we obtain the following consequence.

\medskip %layout
\begin{corollary}
If $G$ is the topological full group of an edge shift $X$ and $S$ is a finite set of rational points of $X$, then the pointwise stabilizer of $S$ is itself a topological full group of an edge shift.
\end{corollary}

For the remainder of this subsection, fix an edge replacement system $\mathcal{R} = (X_0, R, \mathrm{C})$ with $R = \{X_c \mid c \in \mathrm{C}\}$ and denote by $G$ be the corresponding rearrangement group.

\subsection{Stabilizers of Vertices}

Given a vertex $p$, we build the edge replacement system $\mathcal{R}^p = (X_0^p, R^p, \mathrm{C}^p)$ from $\mathcal{R}$, where the set of colors $\mathrm{C}^p$ is $\{ c, c^-, c^+, c^\pm \mid c \in \mathrm{C} \}$ and the base and replacement graphs are defined below.

Let $E$ be the minimal (or any) graph expansion of $\mathcal{R}$ where $p$ appears.
The base graph $X_0^p$ of $\mathcal{R}^p$ is obtained from $E$ with the following modification:
if $e$ is an edge of $E$ that is incident on $p$, then replace its color $c$ with:
\begin{itemize}
    \item $c^-$ if $e$ originates from but does not terminate at $p$,
    \item $c^+$ if $e$ does not originate from but terminates at $p$,
    \item $c^\pm$ if $e$ originates from and terminates at $p$.
\end{itemize}

The replacement graphs for the original colors $c \in \mathrm{C} \subseteq \mathrm{C}^p$ are the same as those of $\mathcal{R}$, while the replacement graphs for the additional colors $c^+$, $c^-$ and $c^\pm$ are obtained in the following way.
\begin{itemize}
    \item $X_{c^-}$ is obtained from $X_c$ with the following modification:
    for every edge $e$ of $X_c$ that is incident at $\iota_c$, replace its color $d$ with:
    \begin{itemize}
        \item $d^-$ if $e$ originates from but does not terminate at $\iota_c$,
        \item $d^+$ if $e$ terminates at but does not originate from $\iota_c$,
        \item $d^\pm$ if $e$ originates from and terminates at $\iota_c$.
    \end{itemize}
    \item $X_{c^+}$ is obtained from $X_c$ with the following modification:
    for every edge $e$ of $X_c$ that is incident at $\tau_c$, replace its color $d$ with
    \begin{itemize}
        \item $d^-$ if $e$ originates from but does not terminate at $\tau_c$,
        \item $d^+$ if $e$ terminates at but does not originates from $\tau_c$,
        \item $d^\pm$ if $e$ originates from and terminates at $\tau_c$.
    \end{itemize}
    \item $X_{c^\pm}$ is obtained from $X_c$ with the following modification:
    for every edge $e$ of $X_c$ that is incident at $\iota_c$ or at $\tau_c$, replace its color $d$ with
    \begin{itemize}
        \item $d^-$ if $e$ originates from $\iota_c$ or $\tau_c$ but does not terminate at $\iota_c$ nor $\tau_c$,
        \item $d^+$ if $e$ terminates at $\iota_c$ or $\tau_c$ but does not originate from $\iota_c$ nor $\tau_c$,
        \item $d^\pm$ if the starting and ending vertices of $e$ are both among $\iota_c$ and $\tau_c$.
    \end{itemize}
\end{itemize}
If $e$ is an edge of a replacement graph $X_c$ of $\mathcal{R}$, then we denote by $e^\varepsilon$ the corresponding edge of $X_{c^\varepsilon}$ (where $\varepsilon$ is either empty or one of the symbols $-$, $+$ or $\pm$).

There is a natural homeomorphism between the limit spaces of $\mathcal{R}$ and $\mathcal{R}^p$ that simply forgets about the marks $-$, $+$ and $\pm$.
More precisely, if we denote by $X$ and $X^p$ the two limit spaces, the homeomorphism is simply the map
\[ \Phi_p \colon X^p \to X,\, x_1^{\varepsilon_1} x_2^{\varepsilon_2} \dots \mapsto x_1 x_2 \dots, \]
where $\varepsilon_i$ is either empty or one of the symbols $-$, $+$ or $\pm$.
This is truly a homeomorphism because the graph expansions (which determine the gluing relation and thus the limit space, recall \cref{def.glue,def.limit.space}) are the same up to forgetting about the marks $-$, $+$ and $\pm$.
In the same fashion, we define a bijection
\[ \phi_p \colon \mathbb{L}_{\mathcal{R}^p} \to \mathbb{L}_\mathcal{R} \]
that forgets about the marks $-$, $+$ and $\pm$ of the edges that appear among the graph expansions of $\mathcal{R}^p$.

\begin{remark}
\label{rmk.vertex.stabilizers}
Carefully analyzing how we decided to append the marks $-$, $+$ and $\pm$ reveals the following property of $\mathcal{R}^p$.
Given any edge $e = x_1 \dots x_k \in \mathbb{L}_{\mathcal{R}^p}$, we have that:
\begin{itemize}
    \item if $e$ originates from but does not terminate at $\Phi_p^{-1}(p)$, then $\mathrm{c}(e)$ has mark $-$;
    \item if $e$ terminates at but does not originate from $\Phi_p^{-1}(p)$, then $\mathrm{c}(e)$ has mark $+$;
    \item if $e$ originates from and terminates at $\Phi_p^{-1}(p)$, then $\mathrm{c}(e)$ has mark $\pm$;
    \item if $e$ is not incident on $\Phi_p^{-1}(p)$ then $\mathrm{c}(e) \in \mathrm{C}$ (i.e., the color of $e$ is a color from the original edge replacement system $\mathcal{R}$).
\end{itemize}
In short, this modification allows us to keep track of edge adjacency to $p$ in every graph expansion.
The idea is in fact reminiscent to that of the gluing automaton in \cref{cha.rational.gluing}.
\end{remark}

\begin{figure}
\centering
\begin{tikzpicture}
    \begin{scope}[xshift=-1cm]
    \node at (0,1.6) {$X_0$};
    \draw[edge,Orange,domain=5:175] plot ({.5*cos(\x)}, {.5*sin(\x)});
    \draw (90:.5) node[above,Orange] {$b_2$};
    \draw[edge,Magenta,domain=185:355] plot ({.5*cos(\x)}, {.5*sin(\x)});
    \draw (270:.5) node[below,Magenta] {$b_3$};
    \node[vertex] (l) at (-1.75,0) {}; \draw (-1.75,0) node[above]{$\iota_{\text{\textcolor{blue}{b}}}$};
    \node[vertex] (cl) at (-.5,0) {};
    \node[vertex] (cr) at (.5,0) {};
    \node[vertex] (r) at (1.75,0) {}; \draw (1.75,0) node[above]{$\tau_{\text{\textcolor{blue}{b}}}$};
    \draw[edge,ForestGreen] (cl) to node[above]{$b_1$} (l);
    \draw[edge,blue] (cr) to node[above]{$b_4$} (r);
    \end{scope}
    \begin{scope}[xshift=4cm]
    \draw[edge,red,domain=5:175] plot ({.5*cos(\x)}, {.5*sin(\x)});
    \draw (90:.5) node[above,red] {$b_2$};
    \draw[edge,red,domain=185:355] plot ({.5*cos(\x)}, {.5*sin(\x)});
    \draw (270:.5) node[below,red] {$b_3$};
    \node at (0,1.6) {$X_{\text{\textcolor{blue}{b}}}$};
    \node[vertex] (l) at (-1.75,0) {}; \draw (-1.75,0) node[above]{$\iota_{\text{\textcolor{blue}{b}}}$};
    \node[vertex] (cl) at (-.5,0) {};
    \node[vertex] (cr) at (.5,0) {};
    \node[vertex] (r) at (1.75,0) {}; \draw (1.75,0) node[above]{$\tau_{\text{\textcolor{blue}{b}}}$};
    \draw[edge,blue] (cl) to node[above]{$b_1$} (l);
    \draw[edge,blue] (cr) to node[above]{$b_4$} (r);
    \end{scope}
    \begin{scope}[xshift=8cm]
    \node at (0,1.6) {$X_{\text{\textcolor{red}{r}}}$};
    \node[vertex] (l) at (-1.25,0) {}; \draw (-1.25,0) node[above]{$\iota_{\text{\textcolor{red}{r}}}$};
    \node[vertex] (r) at (1.25,0) {}; \draw (1.25,0) node[above]{$\tau_{\text{\textcolor{red}{r}}}$};
    \node[vertex] (c) at (0,0) {};
    \node[vertex] (ct) at (0,1) {};
    \draw[edge,red] (l) to node[below]{$r_1$} (c);
    \draw[edge,red] (c) to node[below]{$r_2$} (r);
    \draw[edge,blue] (c) to node[left]{$r_3$} (ct);
    \end{scope}
    \begin{scope}[xshift=4cm,yshift=-3.8cm]
    \draw[edge,red,domain=5:175] plot ({.5*cos(\x)}, {.5*sin(\x)});
    \draw (90:.5) node[above,red] {$b_2^-$};
    \draw[edge,red,domain=185:355] plot ({.5*cos(\x)}, {.5*sin(\x)});
    \draw (270:.5) node[below,red] {$b_3^-$};
    \node at (0,1.6) {$X_{\text{\textcolor{ForestGreen}{b\textsuperscript{$-$}}}}$};
    \node[vertex] (l) at (-1.75,0) {}; \draw (-1.75,0) node[above]{$\iota_{\text{\textcolor{ForestGreen}{b\textsuperscript{$-$}}}}$};
    \node[vertex] (cl) at (-.5,0) {};
    \node[vertex] (cr) at (.5,0) {};
    \node[vertex] (r) at (1.75,0) {}; \draw (1.75,0) node[above]{$\tau_{\text{\textcolor{ForestGreen}{b\textsuperscript{$-$}}}}$};
    \draw[edge,Cyan] (cl) to node[above]{$b_1^-$} (l);
    \draw[edge,blue] (cr) to node[above]{$b_4^-$} (r);
    \end{scope}
    \begin{scope}[xshift=8cm,yshift=-3.8cm]
    \node at (0,1.6) {$X_{\text{\textcolor{Magenta}{r\textsuperscript{$-$}}}}$};
    \node[vertex] (l) at (-1.25,0) {}; \draw (-1.25,0) node[above]{$\iota_{\text{\textcolor{Magenta}{r\textsuperscript{$-$}}}}$};
    \node[vertex] (r) at (1.25,0) {}; \draw (1.25,0) node[above]{$\tau_{\text{\textcolor{Magenta}{r\textsuperscript{$-$}}}}$};
    \node[vertex] (c) at (0,0) {};
    \node[vertex] (ct) at (0,1) {};
    \draw[edge,Magenta] (l) to node[below]{$r_1^-$} (c);
    \draw[edge,red] (c) to node[below]{$r_2^-$} (r);
    \draw[edge,blue] (c) to node[left]{$r_3^-$} (ct);
    \end{scope}
    \begin{scope}[xshift=4cm,yshift=-7.6cm]
    \draw[edge,red,domain=5:175] plot ({.5*cos(\x)}, {.5*sin(\x)});
    \draw (90:.5) node[above,red] {$b_2^+$};
    \draw[edge,red,domain=185:355] plot ({.5*cos(\x)}, {.5*sin(\x)});
    \draw (270:.5) node[below,red] {$b_3^+$};
    \node at (0,1.6) {$X_{\text{\textcolor{Cyan}{b\textsuperscript{$+$}}}}$};
    \node[vertex] (l) at (-1.75,0) {}; \draw (-1.75,0) node[above]{$\iota_{\text{\textcolor{Cyan}{b\textsuperscript{$+$}}}}$};
    \node[vertex] (cl) at (-.5,0) {};
    \node[vertex] (cr) at (.5,0) {};
    \node[vertex] (r) at (1.75,0) {}; \draw (1.75,0) node[above]{$\tau_{\text{\textcolor{Cyan}{b\textsuperscript{$+$}}}}$};
    \draw[edge,blue] (cl) to node[above]{$b_1^+$} (l);
    \draw[edge,Cyan] (cr) to node[above]{$b_4^+$} (r);
    \end{scope}
    \begin{scope}[xshift=8cm,yshift=-7.6cm]
    \node at (0,1.6) {$X_{\text{\textcolor{Orange}{r\textsuperscript{$+$}}}}$};
    \node[vertex] (l) at (-1.25,0) {}; \draw (-1.25,0) node[above]{$\iota_{\text{\textcolor{Orange}{r\textsuperscript{$+$}}}}$};
    \node[vertex] (r) at (1.25,0) {}; \draw (1.25,0) node[above]{$\tau_{\text{\textcolor{Orange}{r\textsuperscript{$+$}}}}$};
    \node[vertex] (c) at (0,0) {};
    \node[vertex] (ct) at (0,1) {};
    \draw[edge,red] (l) to node[below]{$r_1^+$} (c);
    \draw[edge,Orange] (c) to node[below]{$r_2^+$} (r);
    \draw[edge,blue] (c) to node[left]{$r_3^+$} (ct);
    \end{scope}
\end{tikzpicture}
\caption{The edge replacement system for a vertex stabilizer.}
\label{fig.airplane.vertex.stabilizer}
\end{figure}
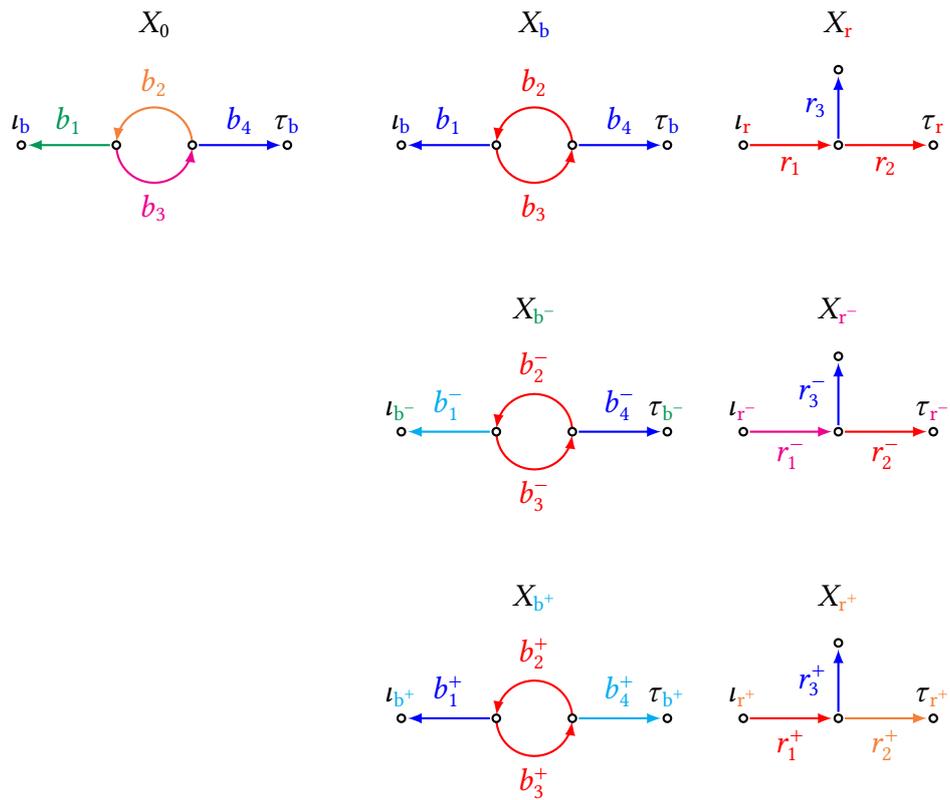

For example, let $p$ be the terminal vertex of the edge $b_2$ in the airplane edge replacement system (\cref{fig.airplane.replacement}).
There are no loops, so we only need to add $-$ and $+$ to the colors blue and red.
\cref{fig.airplane.vertex.stabilizer} depicts the edge replacement system $\mathcal{A}_p$, where {blue}\textsuperscript{$-$} is represented by \textcolor{ForestGreen}{green}, {blue}\textsuperscript{$+$} by \textcolor{Cyan}{cyan}, {red}\textsuperscript{$-$} by \textcolor{Magenta}{magenta} and {red}\textsuperscript{$+$} by \textcolor{Orange}{orange}.
It is clear that every $c$-colored edge of $\mathcal{A}$ that originates from (respectively, terminates at) $p$ is colored by $c^-$ (respectively, $c^+$) in $\mathcal{A}_p$.

Note that, since our graphs do not have isolated vertices by \cref{ass.isolated.vertices}, a graph isomorphism maps a vertex $v$ to a vertex $w$ if and only if it maps an edge originating from (respectively, ending at) $v$ to an edge originating from (respectively, ending at) $w$.
Moreover, it is clear that a rearrangement fixes $p$ if and only if, among its graph pair diagrams that feature $p$ as a vertex, the graph isomorphisms fix $p$.

These observations, together with \cref{rmk.vertex.stabilizers}, yield the desired result.

\medskip %layout
\begin{proposition}
The homeomorphism $\Phi_p \colon X_p \to X$ conjugates the rearrangement group $G_p$ of the edge replacement system $\mathcal{R}^p$ into $G$.
The image of this embedding is precisely the stabilizer of the vertex $p$ under the action of $G$.
\end{proposition}

\subsection{Stabilizers of Other Rational Points}

Recall that the only points of the limit space that may be represented by multiple sequences of the symbol space $\Omega_\mathcal{R}$ are the vertices (\cref{rmk.glue.implies.vertex}).
Thus, given a rational point $q$ that is not a vertex, it must be represented by a unique sequence $\alpha \in \Omega_\mathcal{R}$.
Since $q$ is a rational point, $\alpha$ is a rational sequence, so $\alpha = x \overline{y}$ for some finite words $x = x_1 \dots x_m$ and $y = y_1 \dots y_k$ that, without loss of generality, can be chosen to be of minimal length.
We now describe how to build an edge replacement system $\mathcal{R}^q = (X_0^q, R^q, \mathrm{C}^q)$ from $\mathcal{R}$, where the set of colors $\mathrm{C}^q$ is the same as $\mathrm{C}$ with the addition of $k$ new colors, i.e., $\mathrm{C}^q = \mathrm{C} \cup \{\gamma_i \mid i = 1, \dots, k\}$, and the base and replacement graphs are defined below.

Let $E$ be the minimal graph expansion where the edge $x y_1$ appears.
The base graph $X_0^q$ is obtained from $E$ with a unique modification:
the edge $x y_1$ is colored by $\gamma_1$ instead of the original color of $y_1$.
The replacement graphs for the colors $c \in \mathrm{C} \subset \mathrm{C}^q$ are the same as those of $\mathcal{R}$, while the replacement graphs for the additional colors $\gamma_i$ are obtained with the following modification:
$X_{\gamma_i}$ is the same as $X_{\mathrm{c}(y_i)}$, except that the edge $y_{i+1}$ (where we count cyclically, so that $k+1$ stands for $1$) is colored by $\gamma_{i+1}$.
If $e$ is an edge of a replacement graph $X_c$ of $\mathcal{R}$, then we denote by $e^{\gamma_i}$ the corresponding edge of $X_{\gamma_i}$.

Let $X$ and $X^q$ denote the limit spaces of $\mathcal{R}$ and $\mathcal{R}^q$, respectively.
As with the case of a vertex that we analyzed previously, there is a natural homeomorphism
\[ \Phi_q \colon X^q \to X,\, x_1^{\varepsilon_1} x_2^{\varepsilon_2} \dots \mapsto x_1 x_2 \dots, \]
where $\varepsilon_i$ is either empty or some $\gamma_j$, and there is a natural bijection
\[ \phi_q \colon \mathbb{L}_{\mathcal{R}^q} \to \mathbb{L}_\mathcal{R} \]
that forgets about the marks $\gamma_j$.

\begin{remark}
\label{rmk.rational.stabilizers}
It is easy to see that $\mathcal{R}^q$ enjoys the following property:
a cell $\llbracket e \rrbracket_q$ of $\mathcal{R}^q$ contains $\Phi_q^{-1}(q)$ if and only if the color of the edge $e \in \mathbb{L}_{\mathcal{R}^q}$ is marked by some $\gamma_j$.
Note that, in particular, each graph expansion of $\mathcal{R}^q$ always features precisely one edge whose color is marked by some $\gamma_j$, since points that are not vertices have unique representatives.
\end{remark}

\begin{figure}
\centering
\begin{tikzpicture}
    \begin{scope}[xshift=-2.5cm]
    \node at (0,2) {$X_0$};
    \node[vertex] (b) at (0,0) {};
    \node[vertex] (i) at (180:1.5) {};
    \node[vertex] (u) at (90:1.35) {};
    \node[vertex] (t) at (0:1.5) {};
    \draw[edge] (b) to node[below]{$1$} (i);
    \draw[edge] (b) to node[left]{$2$} (u);
    \draw[edge,Magenta] (b) to node[below]{$3$} (t);
    \end{scope}
    \begin{scope}[xshift=2.5cm]
    \node at (0,2) {$X_1$};
    \node[vertex] (b) at (0,0) {};
    \node[vertex] (i) at (180:1.5) {}; \draw (0:-1.5) node[above]{$\iota_1$};
    \node[vertex] (u) at (90:1.35) {};
    \node[vertex] (t) at (0:1.5) {}; \draw (0:1.5) node[above]{$\tau_1$};
    \draw[edge] (b) to node[below]{$1$} (i);
    \draw[edge] (b) to node[left]{$2$} (u);
    \draw[edge] (b) to node[below]{$3$} (t);
    \end{scope}
    \begin{scope}[xshift=-5cm,yshift=-4cm]
    \node at (0,2) {$X_\text{\textcolor{Magenta}{$\gamma_1$}}$};
    \node[vertex] (b) at (0,0) {};
    \node[vertex] (i) at (180:1.5) {}; \draw (0:-1.5) node[above]{$\iota_\text{\textcolor{Magenta}{$\gamma_1$}}$};
    \node[vertex] (u) at (90:1.35) {};
    \node[vertex] (t) at (0:1.5) {}; \draw (0:1.5) node[above]{$\tau_\text{\textcolor{Magenta}{$\gamma_1$}}$};
    \draw[edge,blue] (b) to node[below]{$1^{\gamma_1}$} (i);
    \draw[edge] (b) to node[left]{$2^{\gamma_1}$} (u);
    \draw[edge] (b) to node[below]{$3^{\gamma_1}$} (t);
    \end{scope}
    \begin{scope}[xshift=0cm,yshift=-4cm]
    \node at (0,2) {$X_\text{\textcolor{blue}{$\gamma_2$}}$};
    \node[vertex] (b) at (0,0) {};
    \node[vertex] (i) at (180:1.5) {}; \draw (0:-1.5) node[above]{$\iota_\text{\textcolor{blue}{$\gamma_2$}}$};
    \node[vertex] (u) at (90:1.35) {};
    \node[vertex] (t) at (0:1.5) {}; \draw (0:1.5) node[above]{$\tau_\text{\textcolor{blue}{$\gamma_2$}}$};
    \draw[edge] (b) to node[below]{$1^\star$} (i);
    \draw[edge,Green] (b) to node[left]{$2^\star$} (u);
    \draw[edge] (b) to node[below]{$3^\star$} (t);
    \end{scope}
    \begin{scope}[xshift=5cm,yshift=-4cm]
    \node at (0,2) {$X_\text{\textcolor{Green}{$\gamma_3$}}$};
    \node[vertex] (b) at (0,0) {};
    \node[vertex] (i) at (180:1.5) {}; \draw (0:-1.5) node[above]{$\iota_\text{\textcolor{green}{$\gamma_3$}}$};
    \node[vertex] (u) at (90:1.35) {};
    \node[vertex] (t) at (0:1.5) {}; \draw (0:1.5) node[above]{$\tau_\text{\textcolor{green}{$\gamma_3$}}$};
    \draw[edge] (b) to node[below]{$1^{\gamma_3}$} (i);
    \draw[edge,Magenta] (b) to node[left]{$2^{\gamma_3}$} (u);
    \draw[edge] (b) to node[below]{$3^{\gamma_3}$} (t);
    \end{scope}
\end{tikzpicture}
\caption{The edge replacement system for the stabilizer of the rational point $q = \llbracket 3 \overline{122} \rrbracket$ of the dendrite $D_3$ that is not a vertex.}
\label{fig.dendrite.rational.stabilizer}
\end{figure}
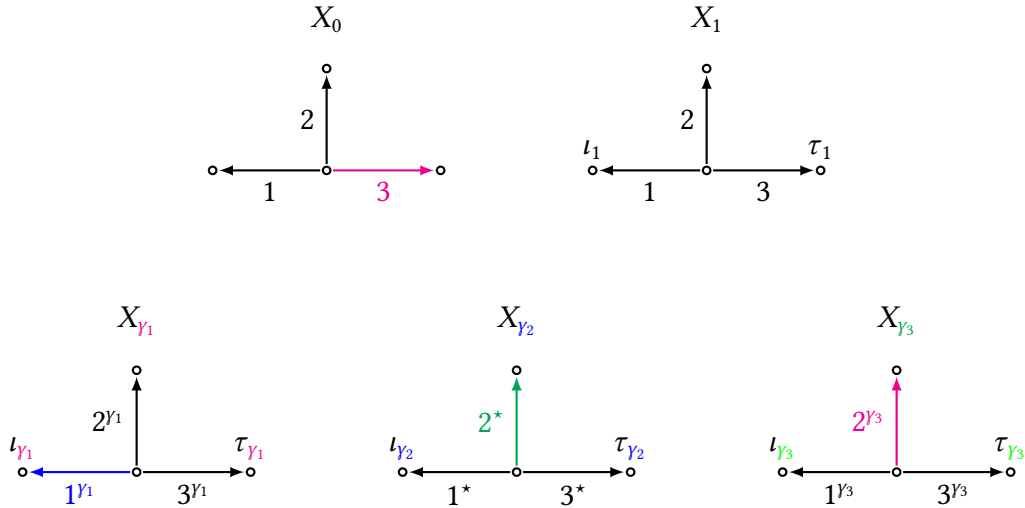

For example, let $q = \llbracket 3 \overline{122} \rrbracket$ in the dendrite edge replacement system $\mathcal{D}_3$ (\cref{fig.dendrite.replacement}).
Since the period is $k=3$, we introduce $3$ new colors $\gamma_1$, $\gamma_2$ and $\gamma_3$, which we represent by \textcolor{Magenta}{magenta}, \textcolor{blue}{blue} and \textcolor{Green}{green}, respectively.
The edge replacement system for the stabilizer of $q$ is depicted in \cref{fig.dendrite.rational.stabilizer}.
It is clear that, for every edge $e$ of $\mathcal{D}_3$, the cell $\llbracket e \rrbracket$ in $\mathcal{D}_3$ contains $q$ if and only if the corresponding edge is colored by red in ${\mathcal{D}_3}_q$.

Now, a rearrangement $g$ of the edge replacement system $\mathcal{R}$ fixes $q = \llbracket x \overline{y} \rrbracket$ if and only if it there exist natural numbers $d_1$, $d_2$ and $l$ (possibly zero) such that $g$ maps canonically (\cref{def.canonical.homeomorphism}) between $\llbracket x y^{d_1} y_1 \dots y_l \rrbracket$ and $\llbracket x y^{d_2} y_1 \dots y_l \rrbracket$.
Together with \cref{rmk.rational.stabilizers}, this shows the desired result.

\medskip %layout
\begin{proposition}
The homeomorphism $\Phi_q \colon X_q \to X$ conjugates the rearrangement group $G_q$ of the edge replacement system $\mathcal{R}^q$ into $G$.
The image of this embedding is precisely the stabilizer of the rational point $q$ under the action of $G$.
\end{proposition}

\section[Rearrangement Groups Embed into \texorpdfstring{$V$}{V}]{Every Rearrangement Group Embeds into Thompson's group \texorpdfstring{$V$}{V}}
\label{sec.embedding.into.V}

In this section we show that all rearrangement groups embed into Thompson's group $V$.
This was first noted in \cite[Remark 1.8]{conjugacy}, but here we develop a more fleshed out and detailed proof that also applies to non-expanding edge replacement systems.
The arguments are not dissimilar to those used in \cite[Theorem 6.14]{types}, which shows that certain stabilizers of Thompson's group $V$ are isomorphic to topological full groups of edge shifts.

In \cref{prop.embedding.into.topological.full.groups} we have showed that each rearrangement group embeds in the topological full group of some edge shift (\cref{sub.topological.full.groups}).
Thus, we restrict our attention to such groups and, in two lemmas, we show that the topological full group of every edge shift embeds in Thompson's group $V$.
This fact is arguably known by the experts, but no proof has been included in the literature so far, as far as the author knows.

\subsection{Topological Full Groups of Binary Edge Shifts}

\begin{definition}
\label{def.binary.edge.shift}
An edge shift is \textbf{binary} if every vertex of the graph that defines it has out-degree $2$.
\end{definition}

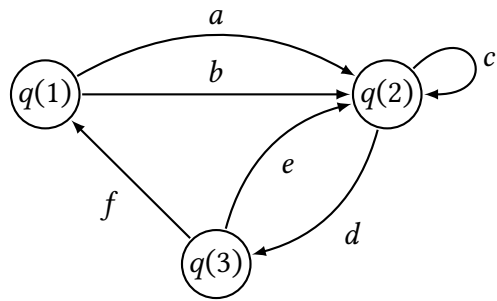
\begin{figure}
\centering
\begin{tikzpicture}
    \node[vertex] (1) at (0,0) {$q(1)$};
    \node[vertex] (2) at (4.5,0) {$q(2)$};
    \node[vertex] (3) at (2.25,-2.25) {$q(3)$};
    \draw[edge] (1) to[out=30,in=150] node[above]{$a$} (2);
    \draw[edge] (1) to node[above]{$b$} (2);
    \draw[edge] (2) to[loop,out=45,in=0,looseness=7.5] node[right]{$c$} (2);
    \draw[edge] (2) to[out=255,in=15] node[below right]{$d$} (3);
    \draw[edge] (3) to[out=75,in=195] node[below right]{$e$} (2);
    \draw[edge] (3) to node[below left]{$f$} (1);
\end{tikzpicture}
\caption{A graph $\Gamma$ that defines binary edge shifts.}
\label{fig.binary.edge.shift}
\end{figure}

Observe that binary edge shifts are precisely those in which, given a valid prefix $w$ in the language of the edge shift, there are exactly two distinct successors $wx$ and $wy$.
For example, \cref{fig.binary.edge.shift} depicts a graph that defines binary edge shifts.

\medskip %layout
\begin{lemma}
\label{lem.binary.shifts.embedding.into.V}
For every binary edge shift, its topological full group embeds into Thompson's group $V$.
\end{lemma}

\begin{proof}
The embedding that we describe below is essentially induced by a map that forgets colors and maps every pair of edges originating from each fixed vertex to $0$ and $1$.

Let $\Gamma$ be the graph that defines the binary edge shift $\Omega$ and let $E$ be its set of edges.
Denote by $\mathfrak{C}$ the usual binary Cantor space (which is the full shift on $\{0,1\}$).
For each vertex $q$ of $\Gamma$, fix an ordering of the two out-going edges, which is simply a bijection $\phi_q \colon \iota^{-1}(v) \to \{0,1\}$.
Together, these bijections induce a surjective map $\phi \colon E \to \{0,1\}$ that restricts to $\phi_q$ on each edge that originates from $q$.
Define the map
\[ \Phi \colon \Omega \to \mathfrak{C},\, x_1 x_2 \dots \mapsto \phi(x_1) \phi(x_2) \dots. \]
For example, if we consider the graph $\Gamma$ of \cref{fig.binary.edge.shift} the binary edge shift $\Omega(\Gamma,q(1))$, we can define $\phi$ so that it maps $a$, $c$ and $e$ to $0$ and $b$, $d$ and $f$ to $1$;
then the word $accdedfb$ is mapped to $00010111$.

We now show that $\Phi$ is a homeomorphism.
If $\Phi(x_1 x_2 \dots) = \Phi(y_1 y_2 \dots)$ then suppose that $x_1 x_2 \dots \neq y_1 y_2 \dots$.
Let $k \in \mathbb{N}$ be the smallest such that $x_{k+1} \neq y_{k+1}$.
Since $x_1 \dots x_k = y_1 \dots y_k$, we only have two choices for $x_{k+1}$ and $y_{k+1}$, as $\iota^{-1}(\tau(x_1 \dots x_k)) = \{x_{k+1}, y_{k+1} \}$.
Then $\phi(x_{k+1}) \neq \phi(y_{k+1})$ because $\phi$ is a bijection when restricted to $\iota^{-1}(\tau(x_1 \dots x_k))$.
But then $\Phi(x_1 x_2 \dots)$ and $\Phi(y_1 y_2 \dots)$ cannot be the same, which is a contradiction.
This shows that $\Phi$ is injective.
By a similar inductive argument, it is easy to show that $\Phi$ is surjective because $\phi$ restricts to a bijection on $\iota^{-1}(\tau(x_1 \dots x_k))$ for each finite prefix $x_1 \dots x_k$, so each ``partial processing'' $\phi(x_1) \dots \phi(x_k)$ can always be continued by finding both an $x_{k+1}$ that maps to $0$ and one that maps to $1$ under $\phi$.
Finally, $\Phi$ maps bijectively between the cones $C(x_1 \dots x_k)$ and $C(\phi(x_1) \dots \phi(x_k))$.
Cones form a basis of clopen sets by \cref{prop.basis.clopen}, so $\Phi$ is a homeomorphism because it is continuous from a compact space to a Hausdorff space.

Let $V_\Gamma$ be a topological full group of the edge shift $\Omega$ and consider the map
\[ \Psi \colon V_\Gamma \to V,\, g \mapsto \Phi g \Phi^{-1}. \]
Since $V_\Gamma$ is comprised of certain prefix-exchange homeomorphisms between finite partitions of $\Omega$ into cones, the image of $\Psi$ is clearly also comprised of prefix-exchange homeomorphisms between finite partitions of $\mathfrak{C}$, which are precisely the elements of Thompson's group $V$.
(Note that the converse would not work, i.e., $\Phi^{-1} V \Phi$ is generally not included in $V_\Gamma$.)
Since the map $\Psi$ is induced by conjugacy, it is evidently an injective group morphism, so we are done.
\end{proof}

We now only need to show that the topological full group of any edge shift embeds into that of some binary edge shift.

\subsection{Embedding into Topological Full Groups of Binary Edge Shifts}

In order to embed the topological full group of an edge shift into that of some binary edge shift, we will modify the graph $\Gamma$ that defines the edge shift so that vertices have out-degree $2$.
For the proof we employ strand diagrams (\cref{sec.SDs}) for ease of visualization, but one can certainly prove the same result without them.

\medskip %layout
\begin{lemma}
\label{lem.embedding.into.binary.shifts}
Suppose that no vertex of $\Gamma$ have out-degree that is less than $2$.
Then the topological full groups of any edge shift that is based on $\Gamma$ embeds into the topological full group of some binary edge shift.
\end{lemma}

\begin{proof}
Consider the topological full group $V(\Gamma,v_0)$ of some edge shift $\Omega(\Gamma,v_0)$.
We have shown in \cref{sub.topological.full.groups} that $V(\Gamma,v_0)$ is a rearrangement group whose edge replacement system has a single edge colored by $v_0$ as its base graph and each replacement graph $X_v$ has a disjoint edge for each edge of $\Gamma$ that originates at a fixed vertex $v$.
Following \cref{SUB groupoid generators} we see that the replacement groupoid (see \cref{SUB replacement groupoid}) associated to these edge replacement rules is generated by split and merge diagrams together with permutation diagrams for every possible color-preserving permutation.
Let $\mathcal{G}$ be this replacement groupoid.

We now modify the split and merge diagrams in order to embed $\mathcal{G}$ in a groupoid $\mathcal{G}^*$ that is based on a binary edge shift.
We first define $\mathcal{G}^*$ as the replacement groupoid of a set of edge replacement rules that are obtained as follows.
Fix a vertex $v$ of $\Gamma$ and let $d$ be the arity of the $v$-colored splits and merges (i.e., $d$ is the out-degree of $v$ in $\Gamma$).
Let $l_1, \dots, l_d$ be the ordering of the leaves of the $v$ replacement tree.
If $d=2$ we do nothing, so suppose that $d \ge 2$.
Then we add to the original edge replacement rules $d-2$ new colors, which we denote by $l_1^*, l_2^*, \dots, l_{d-2}^*$.
The $l_i^*$ replacement graph consists of two disjoint edges:
for all $i \in \{1, \dots, d-3\}$, an edge is colored by $\mathrm{c}(l_i)$ and the other by $l_{i+1}^*$;
for $i=d-3$, one is colored by $\mathrm{c}(l_{d-1})$ and the other by $\mathrm{c}(l_d)$.
We also change the replacement graph $X_v$ to a graph with two disjoint edges, one colored by $\mathrm{c}(l_1)$ and the other by $l_1^*$.
An example of this construction is depicted in \cref{fig.binaryfication}.

\begin{figure}
\centering
\begin{subfigure}[B]{.28\textwidth}
\centering
\begin{tikzpicture}
    \node at (-.75,0) {$X_v$};
    \node[vertex] (a1) at (0,1.2) {};
    \node[vertex] (a2) at (1,1.2) {};
    \node[vertex] (b1) at (0,.4) {};
    \node[vertex] (b2) at (1,.4) {};
    \node[vertex] (c1) at (0,-.4) {};
    \node[vertex] (c2) at (1,-.4) {};
    \node[vertex] (d1) at (0,-1.2) {};
    \node[vertex] (d2) at (1,-1.2) {};
    \draw[edge,blue] (a1) to (a2);
    \draw[edge,red] (b1) to (b2);
    \draw[edge,red] (c1) to (c2);
    \draw[edge,Green] (d1) to (d2);
\end{tikzpicture}
\caption{The original replacement graph $X_v$.}
\end{subfigure}
\hfill
\begin{subfigure}[B]{.7\textwidth}
\centering
\begin{tikzpicture}
    \node at (.5,1.5) {$X_v$};
    \node[vertex] (a1) at (0,.5) {};
    \node[vertex] (a2) at (1,.5) {};
    \node[vertex] (a3) at (0,-.5) {};
    \node[vertex] (a4) at (1,-.5) {};
    \draw[edge,blue] (a1) to (a2);
    \draw[edge,Cyan] (a3) to (a4);
    \begin{scope}[xshift=2.25cm]
    \node at (.5,1.5) {$X_{\textcolor{Cyan}{l_1^*}}$};
    \node[vertex] (b1) at (0,.5) {};
    \node[vertex] (b2) at (1,.5) {};
    \node[vertex] (b3) at (0,-.5) {};
    \node[vertex] (b4) at (1,-.5) {};
    \draw[edge,red] (b1) to (b2);
    \draw[edge,Plum] (b3) to (b4);
    \end{scope}
    \begin{scope}[xshift=4.5cm]
    \node at (.5,1.5) {$X_{\textcolor{Plum}{l_2^*}}$};
    \node[vertex] (c1) at (0,.5) {};
    \node[vertex] (c2) at (1,.5) {};
    \node[vertex] (d1) at (0,-.5) {};
    \node[vertex] (d2) at (1,-.5) {};
    \draw[edge,red] (c1) to (c2);
    \draw[edge,Green] (d1) to (d2);
    \end{scope}
\end{tikzpicture}
\caption{The new replacement graphs $X_v$, $X_{l_1^*}$ and $X_{l_2^*}$.}
\end{subfigure}
\caption{An example of the construction of binary edge shifts described in the proof of \cref{lem.embedding.into.binary.shifts}.}
\label{fig.binaryfication}
\end{figure}
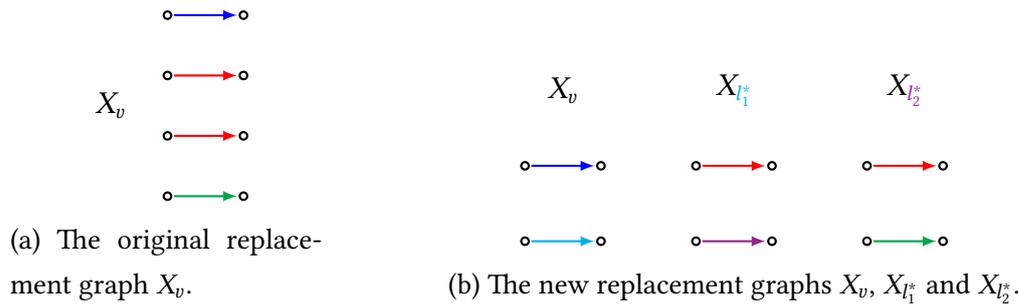

Consider the map $\Phi \colon \mathcal{G} \to \mathcal{G}^*$ that is defined as follows on the generators of $\mathcal{G}$.
A split (merge) diagram with a $v$-colored branching strand is mapped to the composition of $d-1$ split (merge) diagrams as in \cref{fig.binaryfication.trees}.
Permutation diagrams are mapped to themselves (which makes sense because the colors of the original edge replacement rules are also colors of the new one).
One can check that $\Phi(AB) = \Phi(A) \Phi(B)$ for each pair of generators $A$ and $B$ of $\mathcal{G}$, so $\Phi$ extends to a groupoid morphism $\mathcal{G} \to \mathcal{G}^*$.
Similarly, one can see that $\Phi$ is injective.

\begin{figure}
\centering
\begin{tikzpicture}
    \node[node] (root) at (0,1) {};
    \node[node] (s) at (0,0) {};
    \node[node] (1) at (-1.5,-1) {};
    \node[node] (2) at (-.5,-1) {};
    \node[node] (3) at (.5,-1) {};
    \node[node] (4) at (1.5,-1) {};
    \draw (root) to (s);
    \draw[blue] (s) to[out=180,in=90,looseness=1.25] (1);
    \draw[red] (s) to[out=180,in=90,looseness=1.25] (2);
    \draw[red] (s) to[out=0,in=90,looseness=1.25] (3);
    \draw[Green] (s) to[out=0,in=90,looseness=1.25] (4);
    \draw[thick,-stealth] (2.8,-.5) -- (4.35,-.5);
    \begin{scope}[xshift=7cm]
    \node[node] (root) at (0,2) {};
    \node[node] (s) at (0,1) {};
    \node[node] (1) at (-1.5,-2) {};
    \node[node] (1*) at (.5,0) {};
    \node[node] (2) at (-.5,-2) {};
    \node[node] (2*) at (1,-1) {};
    \node[node] (3) at (.5,-2) {};
    \node[node] (4) at (1.5,-2) {};
    \draw (root) to (s);
    \draw[blue] (s) to[out=180,in=90,looseness=1.25] ($(1)+(0,2)$) to (1);
    \draw[Cyan] (s) to[out=0,in=90,looseness=1] (1*);
    \draw[red] (1*) to[out=180,in=90,looseness=1.25] ($(2)+(0,1)$) to (2);
    \draw[Plum] (1*) to[out=0,in=90,looseness=1] (2*);
    \draw[red] (2*) to[out=180,in=90,looseness=1.25] (3);
    \draw[Green] (2*) to[out=0,in=90,looseness=1] (4);
    \end{scope}
\end{tikzpicture}
\caption{The mapping $\Phi$ of the split corresponding to the replacement graph $X_v$ from \cref{fig.binaryfication}.}
\label{fig.binaryfication.trees}
\end{figure}
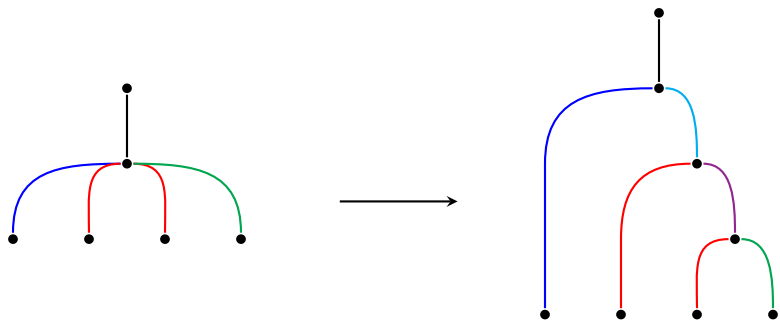

Now, the subgroupoid $V(\Gamma,v_0)$ of $\mathcal{G}$ consists of those elements of $\mathcal{G}$ whose source and sink is labeled by a sole $v_0$-colored edge.
Then its image under $\Phi$ is a subgroupoid of $\mathcal{G}^*$ that consists of certain elements of $\mathcal{G}^*$ whose source and sink is labeled by a sole $v_0$-colored edge.
Thus, $\Phi(V(\Gamma,v_0))$ is a subgroup of the rearrangement group of the modified edge replacement system, which is the topological full group of a binary edge shift, so we are done.
\end{proof}

This proof was inspired by well-known arguments that show how the Higman-Thompson groups $V_{n,r}$ embed into Thompson's group $V$.

\subsection{Dealing with Vertices of Out-Degree \texorpdfstring{$1$}{1}}

The previous lemma works for graphs whose vertices all have degrees that are at least $2$.
Vertices of out-degree $0$ can be ignored, as no infinite walk will go through them, but we still need to take care to the case in which some vertices have out-degree $1$.

\begin{figure}
\centering
\begin{tikzpicture}
    \node[vertex,red] (L) at (0,0) {};
    \node[vertex,red] (R) at (2,0) {};
    \draw[red,edge] (L) node[above]{$v$} to node[above]{$e$} (R) node[above]{$w$};
    \draw[edge,dotted,Green] (-1,.5) to (L);
    \draw[edge,dotted,Green] (-1,-.5) to (L);
    \draw[edge,dotted,blue] (3,.5) to (R);
    \draw[edge,dotted,blue] (R) to (3,-.5);
    \draw[edge,dotted,blue] (R) to[loop,out=270,in=225,min distance=1.2cm] (R);
    \draw[thick,-stealth] (4.3,0) -- (5.7,0);
    \begin{scope}[xshift=8cm]
    \node[vertex,red] (C) at (0,0) {};
    \draw[edge,dotted,Green] (-1,.5) to (C) node[red,yshift=.5cm]{$v = w$};
    \draw[edge,dotted,Green] (-1,-.5) to (C);
    \draw[edge,dotted,blue] (1,.5) to (C);
    \draw[edge,dotted,blue] (C) to (1,-.5);
    \draw[edge,dotted,blue] (C) to[loop,out=247.5,in=292.5,min distance=1.2cm] (C);
    \end{scope}
\end{tikzpicture}
\caption{A schematic depiction of a contraction of $\Gamma$ that does not alter the topology nor the topological full groups of its edge shifts.}
\label{fig.contraction}
\end{figure}
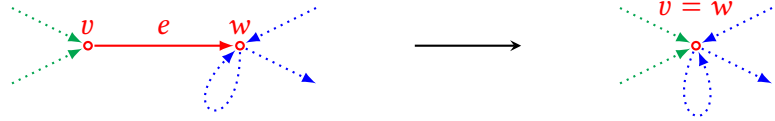

\begin{remark}
Suppose that a graph $\Gamma$ has a vertex $v$ of out-degree $1$ that is not the original and terminal vertex of a loop.
Let $e$ be the unique edge that originates from $v$ and denote by $w$ the terminal vertex of $e$.
Let $\Gamma^*$ be the graph obtained from $\Gamma$ by contracting the edge $e$, meaning that $v$ and $w$ are identified and $e$ is removed.
This move is shown schematically in \cref{fig.contraction}.
Since every walk in $\Gamma$ that goes through $v$ is forced to move to $w$ next, the edge shifts on $\Gamma$ and $\Gamma^*$ are homeomorphic:
the homeomorphism simply removes each letter $e$.
Moreover, conjugating by this homeomorphism provides isomorphisms between the topological full groups $V(\Gamma,u)$ and $V(\Gamma^*,u)$, for each vertex $u$ of $\Gamma$).

This means that we can assume that our group $V(\Gamma,v_0)$ is based on a graph $\Gamma$ where every vertex has out-degree at least $2$ unless its unique out-going edge is a loop.
Note that, in particular, each inescapable cycle (\cref{def.inescapable.cycle}) can be replaced by a loop of an arbitrary color among those that appear in the cycle.

Now, suppose that in $\Gamma$ every vertex has out-degree at least $2$ except for vertices $v_1, \dots, v_k$ whose out-going edges are loops.
Let $V_\Gamma$ be a topological full group of some edge shift that is based on such $\Gamma$.
The edge replacement rules that realize $V_\Gamma$ (as discussed in \cref{sub.topological.full.groups}) can be chosen to be expanding except for the isolated null-expanding edges that correspond to $v_1, \dots, v_k$.
Thus, by \cref{prop.null.expanding.isolated.rarrangement}, $V_\Gamma$ can be realized as the rearrangement group of an expanding edge replacement system, and by \cref{prop.embedding.into.topological.full.groups} it embeds in the topological full group of some edge shift whose graph has vertices of out-degree at least $2$.

Ultimately, we can change the graph $\Gamma$ so that its vertices all have out-degree at least $2$.
\end{remark}

Now \cref{lem.binary.shifts.embedding.into.V,lem.embedding.into.binary.shifts} and the previous remark, together with \cref{prop.embedding.into.topological.full.groups}, immediately yield the desired result.

\medskip %layout
\begin{theorem}
\label{thm.embedding.into.V}
Rearrangement groups embed into Thompson's group $V$.
\end{theorem}

\begin{remark}
Together with this result, the fact that $QV$ is a rearrangement group (\cref{sub.thompson.like}) provides an alternative proof of the fact that it embeds into Thompson's group $V$.
This was the essential result \cite[Theorem 4]{LehnertConjecture} that allowed the restatement of Lehnert's conjecture that is described in \cref{sub.V}.
\end{remark}

%%%%%%%%%%%%%%%%%%%%%%%%%

\emergencystretch=2em
\printbibliography[heading=bibintoc]

\end{document}